\documentclass[]{book}
\usepackage[T1]{fontenc}
\usepackage[utf8]{inputenc}
\usepackage{lmodern}
\usepackage[all]{xy}
\usepackage[a4paper]{geometry}
\usepackage[english]{babel}
\usepackage{fullpage}
\usepackage{eso-pic}
\usepackage{graphicx}
\usepackage{tikz}
\usepackage[tikz]{bclogo}
\usepackage{csquotes}

\usepackage{colortbl}
\usepackage{graphicx}

\usepackage[encapsulated]{CJK}

\usepackage[textsize=tiny]{todonotes}

\usepackage{pdfpages}

\usepackage{multicol}
\usepackage[makeroom]{cancel}
\usepackage{ulem}
\usepackage{amsmath,amsfonts,amssymb}
\usepackage{tikz}

\usepackage{epigraph}
\usepackage{hyperref}
\usepackage[europeanresistors,oldvoltagedirection]{circuitikz}

\usepackage{systeme}
\newcommand{\toto}{\rightrightarrows}

\usepackage{latexsym,amsfonts,amsthm,amsmath,amscd,amssymb,color}

\usepackage[many]{tcolorbox}

\newtheorem{lemma}{Lemma}[section]
\newtheorem{con}[lemma]{Convention}

\newtheorem{proposition}[lemma]{Proposition}
\newtheorem{corollary}[lemma]{Corollary}
\newtheorem{question}[lemma]{Question}

\theoremstyle{definition}
\newtheorem{definition}[lemma]{Definition}
\newtheorem{example}[lemma]{Example}
\newtheorem{remark}[lemma]{Remark}
\newtheorem{exo}[lemma]{Exercise}

\newtcbtheorem
  [use counter*=lemma,number within=section]
  {theorems}
  {Theorem}
  {
		fontupper=\itshape,
		breakable,
		enhanced,
    colback=gray!25,
    colframe=gray!0!black,
    fonttitle=\bfseries,
  }
  {thm}

\newtcbtheorem
  [use counter*=lemma,number within=section]
  {propositions}
  {Proposition}
  {
		fontupper=\itshape,
		breakable,
		enhanced,
    colback=gray!25,
    colframe=gray!0!black,
    fonttitle=\bfseries,
  }
  {thm}

\newtcbtheorem
  [use counter*=lemma,number within=section]
  {conjectures}
  {Conjecture}
  {
		fontupper=\itshape,
		breakable,
		enhanced,
    colback=gray!25,
    colframe=gray!0!black,
    fonttitle=\bfseries,
  }
  {thm}

\newtcbtheorem
  [use counter*=lemma,number within=section]
  {definitions}
  {Definition}
  {
		fontupper=\itshape,
		breakable,
		enhanced,
    colback=gray!25,
    colframe=gray!0!black,
    fonttitle=\bfseries,
  }
  {def}

\newtcbtheorem
  [use counter*=lemma,number within=section]
  {questions}
  {Question}
  {
		fontupper=\itshape,
		breakable,
		enhanced,
    colback=gray!25,
    colframe=gray!0!black,
    fonttitle=\bfseries,
  }
  {thm}

\newtcbtheorem
  [use counter*=lemma,number within=section]
  {corollaries}
  {Corollary}
  {
		fontupper=\itshape,
		breakable,
		enhanced,
    colback=gray!25,
    colframe=gray!0!black,
    fonttitle=\bfseries,
  }
  {thm}

\newtcbtheorem
  [use counter*=lemma,number within=section]
  {notations}
  {Notation}
  {
		fontupper=\itshape,
		breakable,
		enhanced,
    colback=gray!25,
    colframe=gray!0!black,
    fonttitle=\bfseries,
  }
  {thm}

\newcounter{numexof}
\theoremstyle{remark}

\newtheorem{theorem}[lemma]{Theorem}
\newtheorem{exemple}[lemma]{Example}

\newcommand{\dd}{\mathrm{d}}
\newcommand{\E}{\mathcal{E}}
\usetikzlibrary{calc,trees,positioning,arrows,chains,shapes.geometric,
    decorations.pathreplacing,decorations.pathmorphing,shapes,
    matrix,shapes.symbols}

\tikzset{
>=stealth',
  punktchain/.style={
    rectangle, 
    rounded corners, 
    
    draw=black, very thick,
    text width=10em, 
    minimum height=3em, 
    text centered, 
    on chain},
  line/.style={draw, thick, <-},
  element/.style={
    tape,
    top color=white,
    bottom color=blue!50!black!60!,
    minimum width=8em,
    draw=blue!40!black!90, very thick,
    text width=10em, 
    minimum height=3.5em, 
    text centered, 
    on chain},
  every join/.style={->, thick,shorten >=1pt},
  decoration={brace},
  tuborg/.style={decorate},
  tubnode/.style={midway, right=2pt},
}

\newcommand{\HRule}{\rule{\linewidth}{0.5mm}}
\newcommand{\blap}[1]{\vbox to 0pt{#1\vss}}
\newcommand\AtUpperLeftCorner[3]{
  \put(\LenToUnit{#1},\LenToUnit{\dimexpr\paperheight-#2}){\blap{#3}}
}
\newcommand\AtUpperRightCorner[3]{
  \put(\LenToUnit{\dimexpr\paperwidth-#1},\LenToUnit{\dimexpr\paperheight-#2}){\blap{\llap{#3}}}
}

\title{Singular foliations\\
 }
\date{July 18-22, 2022}
 \author{Camille Laurent-Gengoux, Ruben Louis, Leonid Ryvkin}

\makeatletter

\begin{document}

\normalem

\begin{titlepage}
    \enlargethispage{2cm}
 
    \AddToShipoutPicture{
        \AtUpperLeftCorner{1.5cm}{1cm}{\includegraphics[width=4cm]{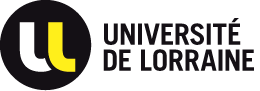}}
        \AtUpperLeftCorner{1.5cm}{3cm}{\includegraphics[width=6cm]{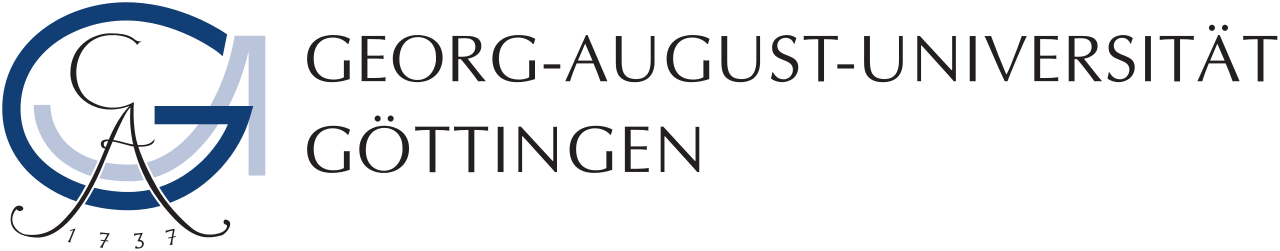}}
        \AtUpperLeftCorner{8.2cm}{1cm}{\includegraphics[width=3cm]{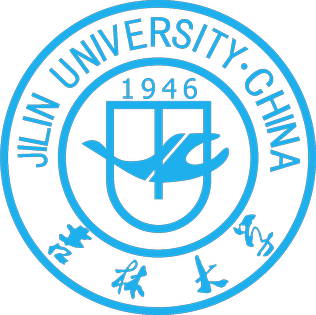}}
        \AtUpperRightCorner{1.5cm}{1.5cm}{\includegraphics[width=7cm]{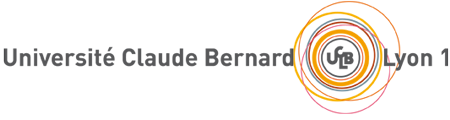}}
   }
 
    \begin{center}
        \vspace*{10cm}
 
        \textsc{\@ An invitation to singular foliations}

        \vspace*{0.5cm}
  \text{CRM Barcelona, Poisson geometry summer school}

\HRule
        \vspace*{0.5cm}
 
        \large{Camille Laurent-Gengoux\footnote{
Universit\'e de Lorraine, CNRS, IECL, F-57000 Metz, France, \texttt{camille.laurent-gengoux@univ-lorraine.fr}}
}, Ruben Louis\footnote{Department of Mathematics, Jilin University, Changchun 130012, Jilin, China and 
Institut f\"ur Mathematik, Georg-August-Universit\"at G\"ottingen, G\"ottingen, Germany, \texttt{ruben.louis@mathematik.uni-goettingen.de}}, Leonid Ryvkin\footnote{Université Claude Bernard Lyon 1, Bâtiment Jean Braconnier 21, avenue Claude Bernard, 69100 Villeurbanne, France, and University of Göttingen, Bunsenstr. 3-5, 37073 Göttingen, Germany, \texttt{ryvkin@math.univ-lyon1.fr}}
    \end{center}
 
    \vspace*{3.2cm}
 
    \begin{center}
        \makebox[\textwidth]{\includegraphics[width=6cm]{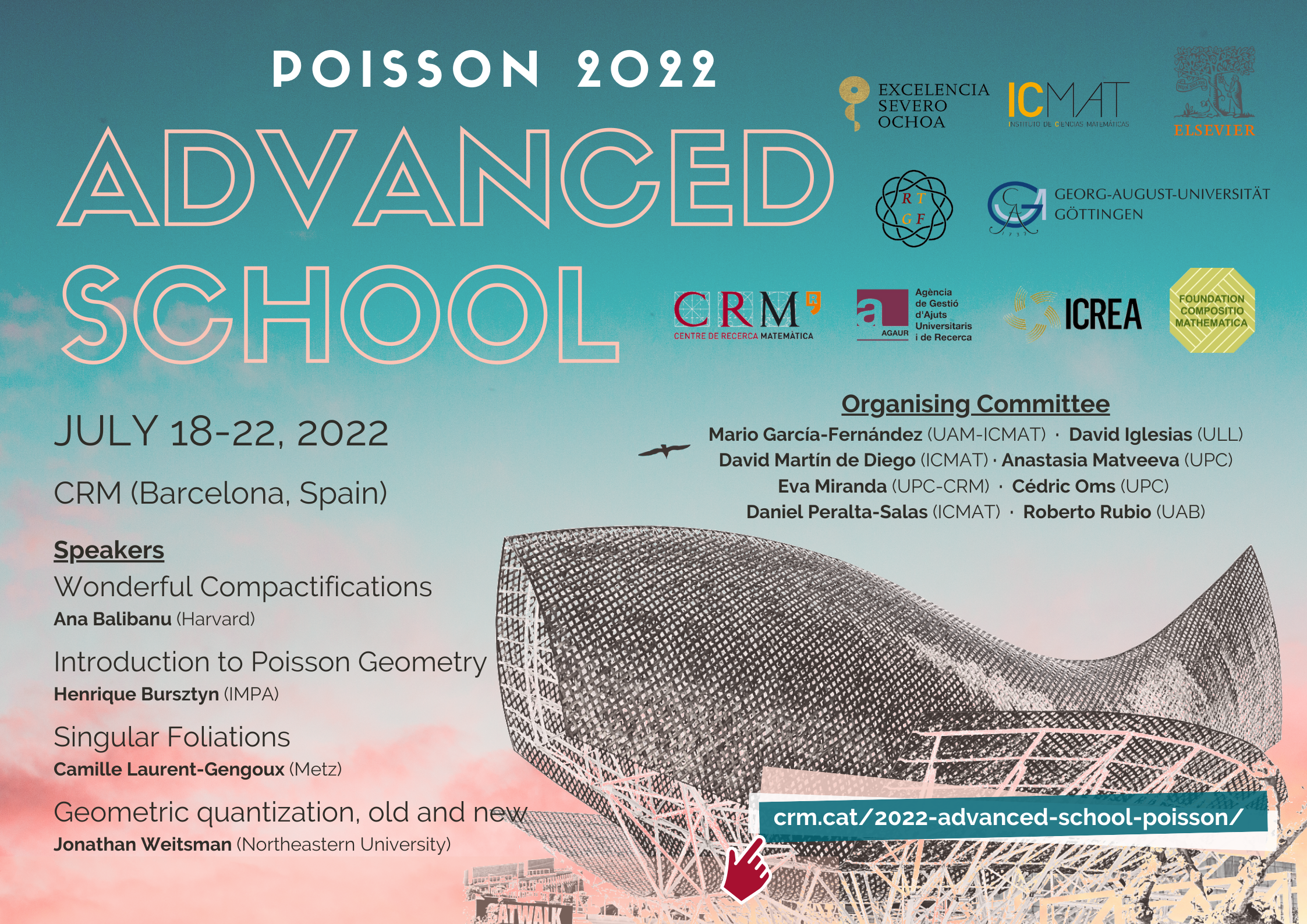}}
         \text{}
    \end{center}
 
\end{titlepage}
\ClearShipoutPicture

	\tableofcontents

\chapter*{Introduction}

\section*{Flying under radar: Singular foliations.}

Singular foliations are so common in mathematics that they often go unnoticed.

Regular foliations have been long studied; The Frobenius theorem \cite{MR1579710,MR515141} is taught quite early in the differential geometry curriculum. Holonomy (or ``first return'') is a very classical notion \cite{Moerdijk}. 
In contrast, singular foliations have never been studied with such an intensity. Still, there is a long story behind foliations that have leaves which are not all the same dimension: 
\begin{enumerate}
    \item  As pointed by Sylvain Lavau  \cite{Lavau18}, the 1960s saw an intense debate about finding a correct definition of a singular foliation. 
    The discussion led to some major discoveries by H. Hermann  (1962), T. Nagano (1966), P. Stefan  (1970),  H. Sussmann  (1973)\footnote{We refer to Sylvain Lavau's excellent article \cite{Lavau18} for the historical aspects. \cite{Lavau18} can also be read as an introduction to the subject.}. See \cite{Hermann}-\cite{Nagano}-\cite{Stepan1}-\cite{Stepan2}-\cite{Sussmann1}-\cite{Sussmann2}. 
    
\item Then the subject seems to have been slightly forgotten, or at least put aside. There were, still, important contributions to the linearization problem by Dominique Cerveau \cite{Cerveau} (where ``singular foliations'' appear under the name of ``involutive distributions'') in 1977 and Pierre Dazord \cite{Dazord} who defined a holonomy map for a singular leaf in 1984. There were other contributions coming from complex geometry, in particular - but not only - about codimension 1 or dimension 1 singular foliations (see the excellent review \cite{MR4454136}) and the theory of residues (see e.g. Paul Baum and Raoul Bott's \cite{zbMATH03423310}, Ali Sinan Sert\"oz's \cite{zbMATH04118703}, Andr{\'e} Belotto da Silva and Daniel Panazzolo's \cite{zbMATH07124409} or Tatsuo Suwa's \cite{zbMATH03931387}). Also, Poisson geometers knew that symplectic leaves of a Poisson manifold, or Lie algebroid leaves, were a sort of ``singular foliation'' \cite{Lichne}-\cite{Weinstein}, but, to our knowledge, rarely developed it as such.
    \item Then, starting in the 2000s, a ``singular foliation's revival'' arose from non-commutative geometry, with pioneering and fundamental works by Iakovos Androulidakis, Claire Debord, Georges Skandalis, and Marco Zambon in particular. It is unfair to summarize their contributions in one sentence, but since we have to do so, let us claim that, from the geometric point of view, a crucial feat is the construction, by Androulidakis and Skandalis \cite{AS}, of a holonomy groupoid of a singular foliation, that extends holonomy groupoids of regular foliations \cite{Moerdijk}, and a smooth groupoid previously constructed by Claire Debord for projective singular foliations \cite{Debord}. A theorem of crucial importance was also established by Claire Debord: although Androulidakis-Skandalis holonomy groupoid is not smooth, it is longitudinally smooth \cite{Debord2}.
   Then, Omar Mohsen \cite{OmarMohsen} introduced a quotient of the holonomy groupoid, now called Mohsen's groupoid.  As an application, the so-called Helffer-Nourrigat conjecture was recently solved by Androulidalis, Mohsen, and Yuncken \cite{AMY}.
    \item[] This holonomy groupoid, or more precisely its natural $C^*$-algebra, is used by this school to define and study elliptic pseudodifferential operators \cite{MR4026451}, analytic indexes \cite{MR4300554}, or to investigate its Baum-Connes conjecture \cite{AS19}, Boutet de Monvel calculus \cite{DS:Boutet} - in one word, to do analysis of singular foliations and their differential operators, to define symbols \cite{MR4239197}. It is now used to solve classical conjectures about hypo-elliptic operators \cite{AMY}. We have no expertise in non-commutative geometry and will not speak much about that subject, although it is certainly the most quickly progressing one among those using singular foliations.
    \item[] Singular foliations, especially those of dimension or codimension $1$ are also currently used in complex dynamics, see, e.g. \cite{zbMATH07372950}-\cite{zbMATH07124409}-\cite{zbMATH06517046}, and  in theoretical physics for perturbation theories of moduli spaces, see e.g., \cite{zbMATH07119821,fischer2024topologicalclassificationsymmetrybreaking}. 
    
The purpose of the present invitation is not to do analysis of singular foliations, although it is certainly the most active topic at the moment. In particular, we are not competent in non-commutative geometry, in index theorem, and pseudo-differential operators. We do not claim to be highly competent neither in holomorphic dynamical systems nor in theoretical physics. 
Our purpose is to simply to introduce the  \underline{geometry} of a singular foliation - maybe we should even say the \underline{differential topology} of  singular foliations, since we will not speak much about metrics. In our opinion, geometry is the easiest entrance to singular foliations. We are also non-competent to include an overview of derived foliations, recently developed by Tony Pantev, Bertrand To\"en and Gabriele Vezzosi (see, e.g. \cite{pantevToen}-\cite{pantevToen2}): there is certainly a link to be explored with what we called the universal Lie $ \infty$-algebroid of a singular foliation in Section \ref{sec:Universal}, while \cite{TV2024} may have some relation with the construction of the holonomy groupoid.

\end{enumerate}
\noindent
\vspace{.5cm}

  Let us go back to the initial debate - in a very anachronistic manner: Should singular foliations be seen:
    \begin{enumerate}
    \item[(0)]  as level sets (called ``leaves'') of (maybe non-independent) functions?
        \item[$\clubsuit $]  as a partition of a manifold into submanifolds?
        \item[$\diamondsuit$] as the data, at each point, of sub-spaces of the tangent space satisfying an involutivity condition?
        \item[$\heartsuit$] as a regular foliation defined on some open subset of the manifold?
        \item[$ \spadesuit$] or as an involutive $\mathcal C^\infty(M) $-module of vector fields (morally thought to be tangent to the leaves)?
    \end{enumerate}
    Definition (0) (i.e., ``level set of non-independent functions'') is opposite to what we intend to study here: leaves would not be manifolds, and even if we work within the context of algebraic geometry (so that these level sets would be affine varieties), there is still a problem: exceptional leaves would be of bigger dimensions than the ``regular'' ones. 
    We do not claim that such a geometry is not interesting by itself, but this is clearly opposite to what we are looking for\footnote{Jokingly, we say that we intend to study lasagna dishes with a few isolated spaghetti, but we do not wish to study isolated lasagnas in a spaghetti dish.}.
    $\heartsuit $ is used in holomorphic geometry, where a singular foliation on a complex manifold $M$ may be defined as a holomorphic regular foliation on a codimension $ \geq 2$ analytic subset of $M$: this definition, however, is essentially equivalent to the holomorphic equivalent of $\clubsuit,\diamondsuit, \spadesuit $ which are essentially equivalent one to the other, when made precise in the right way.

        \noindent
   Now, in the smooth case, the three remaining points of views  $(\clubsuit,\diamondsuit,\spadesuit)$

 have to be made more precise to yield a reasonable definition of a singular foliation. As we shall see in the first chapter, all of them allow counter-examples to properties that we wish to be true. This does not mean that they have to be rejected, but they have to be made precise. 
 
 We may dare to say that after that debate took place in the late 1960s, only two definitions survived to the XXI-st century:
    \begin{enumerate}
        \item[($\spadesuit\star$)] A singular foliation is a  sub-sheaf of the sheaf of vector fields stable under Lie bracket, stable under multiplication by a smooth function,
        and locally finitely generated as a module over smooth functions. \footnote{For those unfamiliar with or hostile to sheaves, this definition can be equivalently stated as: a locally finitely generated involutive sub-$\mathcal C^\infty(M)$-module of the module of compactly supported vector fields.}.
        \item[($\clubsuit\star$)]  A singular foliation is a partition of a manifold into submanifolds called leaves, such that through any vector tangent to a leaf there is at least a vector field tangent to all leaves \cite{Sussmann2,DLPR}.
    \end{enumerate}
    We will work with the first of these definitions, for the following reasons:
    \begin{enumerate}
        \item Definition ($\spadesuit\star $) implies definition ($\clubsuit\star$): Singular foliations in the sense of ($\star$) do admit leaves which are honest submanifolds and partition the manifold\footnote{and this is the least we can require to dare calling an object ``singular foliation'': leaves have to make sense!} and the henceforth obtained partition satisfies ($\clubsuit\star$),
        \item the tangent spaces of these leaves form a (singular) involutive distribution,
        \item it is -according to us- general enough to contain most interesting examples,
        \item but it is restrictive enough to be able to prove strong results,
     while, for instance, singular foliations as in ($\clubsuit\star$) may not admit an AS-holonomy groupoid (at least, not a longitudinally smooth one). 
        \item Last, ($\spadesuit\star$) is used by a now well-established community of non-commutative geometers (Androulidakis, Debord, Mohsen, Skandalis, Yuncken, Skandalis, Zambon - to cite a few) and some theoretical physicists (e.g., Kotov, Strobl), while ($\star\star$) seems to be less commonly used nowadays, although it is not abandoned \cite{Miyamoto,zbMATH07807748,DLPR}.
    \end{enumerate}
For all these reasons, we will present  the theory of singular foliations using Definition ($\spadesuit\star$). Although we had no time or space to present it, we claim that these notes would in fact not present a fundamentally different theory had we decided to use Definition ($\clubsuit\star$).

\section*{Are singular foliations worth studying?}

Is there a point in studying singular foliations? It will depend on where you come from and where you want to go.

First, whoever studies Poisson geometry will encounter a highly non-trivial singular foliation: the symplectic leaves of a Poisson structure. But we claim more: whoever understands classical Poisson geometry has understood objects which are more or less analogous to those used in the geometry of singular foliations. Half of the way is behind you.

 Below, we listed the classical notions of Poisson geometry on the left, and their equivalent objects in the SF-theory on the left\footnote{We use the abbreviations SF= Singular Foliations and AS =Androulidakis-Skandalis}: if you know what the left-hand column is about, understanding the right-hand column should not be overly difficult. Also, we wrote $>,=, <$ to tell which side is, in our subjective opinion, harder to understand.

\vspace{1cm}
\begin{tabular}{|l|c|l|} \hline
  \cellcolor[gray]{0.8} Notion in Poisson geometry & & \cellcolor[gray]{0.8} The equivalent notion in Singular foliation theory \\ \hline
Poisson manifold $(M,\pi) $  & $=$  & Singular foliation $\mathcal F $ on $M$  \\   \hline
Hamiltonian flows are Poisson diffeo. & & Vector fields tangent to $\mathcal F $ are symmetries of $\mathcal F $ \\  (This is almost trivial)& $<< $ & (This is really hard, at least in the smooth case, \\ & & Many existing proofs have gaps…)  \\ \hline 
Weinstein's splitting theorem &$=$&  Singular Foliations' splitting theorem  \\ \hline 
Partition into symplectic leaves & $>$ & Partition into leaves \\ \hline
Transverse Poisson structure (of a leaf) &$=$& Transverse singular foliation (of a leaf)  \\ \hline
Poisson-Dirac reduction & $>$& Induced SF on a transverse submanifold \\ \hline
Lie algebroid structure on $T^* M$      & $>$     & (easy) almost Lie algebroid structures generating $ \mathcal F$ \\  
            & $<<$ & or (harder) the universal Lie $\infty $-algebroid of  $\mathcal F $ \\
      \hline
    Isotropy Lie algebra $\ker \pi_m^\# $  at $m\in M$              & $=$& (easy) isotropy Lie algebra of $\mathcal F $ at $m$  \\ & $ <<$& (harder) isotropy Lie $\infty$-algebra of  $\mathcal F $ at $m$. \\\hline Poisson cohomology & $>$ &Longitudinal cohomology (easy)  \\ & $<<$& Cohomology of the universal Lie $\infty $-algebroid (harder)  \\ \hline
     Symplectic realization & =&  bisubmersions  \\  \hline
     Morita equivalences & = & Equivalences of bisubmersions   \\  \hline
Symplectic Groupoid  & & AS holonomy groupoid   \\ (This is often a smooth groupoid  & $<<<$ & (This is almost never a smooth groupoid,  \\  at worst a stacky groupoid \cite{MR2197411}) & & not even a stacky groupoid) \\\hline
\end{tabular}
\vspace{1cm}

In particular, the AS holonomy groupoid is not like any Lie groupoid Poisson geometry has so far produced. Its non-smoothness is at the origin of the subtle analysis developed by non-commutative geometers. 
Although the AS holonomy groupoid is certainly the most studied aspect of singular foliation at the present time \cite{AS}-\cite{Debord2}-\cite{AS11}-\cite{AS19}-\cite{AZ13}, we will construct it in detail.

\section*{To which area of mathematics do singular foliations belong to?}

As we will see, singular foliations shall be defined as a sub-algebra $\mathcal F$ of vector fields, stable under the Lie bracket and under multiplication by a function, and the leaf through a point $m \in M$ shall be the set of points reachable from $m$ following the flows of vector fields in $\mathcal F$. Those vector fields in $\mathcal F$ are, heuristically, vector fields ``tangent to all leaves''. But there is an additional assumption in the definition that more or less makes consensus nowadays: we should require $ \mathcal F$  to be ``locally finitely generated''. Also, vector fields in $ \mathcal F$ are often supposed to be compactly supported (see Definition \ref{def:consensus}).  Before dealing with those technical points, we have to address a more fundamental question: in which area of mathematics are we?  

The present manuscript is mainly written having in mind  the universe of smooth differential geometry. But singular foliations do make sense in real analytic differential geometry, in complex geometry, and in algebraic geometry as well. And we will try to deal with all three aspects altogether.
For that purpose, we will use the language of sheaves\footnote{The reader interested only in the smooth case may perfectly ignore the word ``sheaf'' and replace it by the corresponding global objects. For technical reasons, it is then better to use ``compactly supported'' objects, or to just replace “the sheaf of” by "the locally defined".}. More precisely:

\begin{enumerate}
    \item In real differential geometry, sheaves can be ignored, and singular foliations on a manifold $M$ will be defined as a locally finitely generated sub-$\mathcal C^\infty(M)$-module of compactly supported vector fields stable under Lie brackets.
    \item In real analytic or holomorphic or algebraic settings, global objects may not exist, or it may be that there are too few of them. One has to work with the sheaf of vector fields, and it does not make sense to consider compactly supported vector fields anymore. 
    Moreover, the “locally finitely generated” condition is equivalent, in this context, to “coherent sheaf”. 
   
    \item In smooth, real analytic or complex settings, singular foliations induce a partition of $M$ into leaves which are smooth, real analytic or complex submanifolds respectively. This is not true anymore in algebraic geometry: the ``leaves'' are not algebraic sub-varieties. This is highly related to the well-known fact that the flow of a polynomial vector field is a real-analytic or holomorphic map, but not a polynomial map in general.
\end{enumerate}

Again, although we will deal with real analytic or holomorphic or algebraic settings, we will mostly take the smooth differential geometry point of view. Also, we will assume that the reader knows everything about differential geometry: classical or less-classical theorems about flows of vector fields will often be admitted, and only those specific to singular foliations shall be detailed.

\subsubsection{Acknowledgments}

The authors would like to thank the organizers of the Poisson 2022 summer school  for offering them the opportunity to present an overview of singular foliations, that eventually led to the present manuscript.
We extend special thanks to Eva Miranda, and to Cédric Oms for the opportunity to publish these notes in the annals of the summer school.

We acknowledge several discussions with several mathematicians, in particular Georges Skandalis, who suggested several improvements, Iakovos Androulidakis, Camilo Angulo, Claire Debord, Simon Raphael Fischer, Anastasios Fotiadis, Noriaki Ikeda, Jun Jiang, Oleksii Kotov, Sylvain Lavau, Hsuan-Yi Liao, Omar Mohsen, Hadi Nahari, Vladimir Salnikov, Yunhe Sheng, Bernd Stratmann, Thomas Strobl, Robert Yuncken, Chenchang Zhu. Last, we would like to thank Cédric Rigaud, who wrote a master degree memoir that we turned into examples.

All three authors have discussed the present manuscript during their stay at IHP in the  summer of 2023. 

C.L.-G.  would like to thank the Tsing Hua University \begin{CJK*}{UTF8}{bsmi}國立清華大學\end{CJK*}  and the NCTS \begin{CJK*}{UTF8}{bsmi}國家理論科學研究中心\end{CJK*}  for their hospitality.

R.L. acknowledges the full financial support provided by Jilin University for his joint postdoctoral position with Göttingen University.

L.R. acknowledges the financial support of the DFG through the grant Higher Lie Theory and the LabEx MiLyon. 

\subsubsection*{Conventions}

Throughout the text: manifolds shall be separated and second countable.  
Vector fields on a manifold $M$ shall be denoted by
$\mathfrak X(M) $, or simply $\mathfrak X $ when there is no ambiguity, but, depending on the context, the previous notation may stand for global vector fields on $M$ or for the sheaf of vector fields. We will always spell out our convention in due place.
Sections of a vector bundle $E$ shall be denoted by $\Gamma(E)$, but, again, we will sometimes make no difference between the notation of the sheaf of local sections or global sections. Sections over an open subset $\mathcal U \subset M $ of a vector bundle $E \to M$ shall be denoted by $\Gamma_{\mathcal U}(E) $. This non-constant notation was chosen to avoid using too heavy symbols coming from sheaf theory all along the text, especially when sheaves are not needed.

We also invite the reader to check the conventions that we use to denote singular foliations, at the end of section \ref{sec:sheafdef}.
When we define them, we need to distinguish two notions of singular foliations, respectively denoted by $ \mathcal F_c$ and $ \mathcal F_\bullet$. This equivalence allows us to simply denote a singular foliation by an $ \mathcal F$, and the rest of the text will simply use that notation. Compactly supported vector fields  or sections shall be denoted by $\mathfrak X_c(M) $ and $\Gamma_c(E)$, respectively. In order to deal with holomorphic, real analytic and smooth settings simultaneously, we will often use the symbol $\mathcal O $ for the relevant sheaf of functions.  Also, for $X$ a vector field on $M$ or $e$ a section of a vector bundle $E \to M$, we denote by $X_{|_{m} }$ and  $e_{|_{m} }$ their values at a point $ m \in M$.

Also, to avoid having to repeat “Let $M$ be  a manifold equipped with a singular foliation”, we will often say "Let $(M,\mathcal F) $" be a foliated manifold.

Restrictions to  an open $\mathcal U \subset M $ or ``any-mathematical-notion-$N$-that-restricts'' will mostly be denoted by $\mathfrak i^*_\mathcal U N $.

\chapter{What is a singular foliation?}

\section{Naive and less naive attempts of a definition of a singular foliation}

In order to understand the geometric ideas behind the consensus definition of a singular foliation, let us make a list of definitions that are natural, but either turned out to be dead ends, or did not yet prevail so far for some reason. 

\vspace{0.2cm}
\noindent
This section is widely inspired by Sylvain Lavau's \cite{Lavau18}, and by Iakovos Androulidakis and Marco Zambon's \cite{MZ16}.

\subsection{Partitionifolds. Is a singular foliation simply a partition by smooth manifolds?}

\vspace{.5cm}
\noindent
Most differential geometers are used to hear the word ``foliation'' as referring to what we will call here ``regular foliation''.
Such a ``regular'' foliation partitions
a manifold into submanifolds, all of the same dimension. As a consequence, the most natural idea that comes to mind when trying to make up a definition of a singular foliation is to try to define them as being a disjoint union of submanifolds called ``leaves'' - now of varying dimension.
This perfectly makes sense, but let us give it another name.

\vspace{0.2cm}
\noindent
Unless otherwise specified, the discussion of this section makes sense in smooth, real analytic, or complex geometry.

\vspace{.5cm}

\begin{definitions}{A first attempt to define singular foliations: partitionifolds}{partitionifolds}
Let $M$ be a manifold.
A \textit{partitionifold}\footnote{We suggest the word \emph{partitionni\'et\'e} in French.} of $M$ is a partition of $M$ into connected immersed submanifolds\footnote{From now on, ``submanifold'' means by default ``immersed submanifolds''.}, called \emph{leaves}.  
\end{definitions}
\vspace{.5cm}

\normalfont
Let us introduce a convenient notation.

\vspace{.5cm}
\begin{notations}{To a point, we associate its leaf}{Not:LeafThroughm}
A partitionifold on a manifold $M$ shall be denoted as a map:
 $$ \begin{array}{llll}  L_\bullet \colon & M & \to & \hbox{ $\{$Submanifolds  of $M\}$}\\ &m &\mapsto & \, L_m \end{array}  $$ that maps a point $ m \in M$ to the submanifold in the partition to which $m$ belongs. Also, for all $m \in M $, $L_m $ shall be called the \emph{leaf through $m$}.
\end{notations}
\vspace{.5cm}

Below are two examples of partitionifolds that we not wish to allow as being decent ``singular foliations''.

\begin{example}
\label{exa:WeirdStraigtLines}
We call \emph{bioriented partition} the partitionifold on $ \mathbb R^2$ given by the following partition
\begin{enumerate}
\item  the leaf $ L$ given by the straight line $ \{y=0\}$, called \emph{central leaf},
\item the  half lines $H_y $ given by $ \{ (x,y) \, |\,  x \in \mathbb  R_-^* \} $ for all $y \neq 0 $,
\item   the half lines $ V^+_x$ given for all $x \geq 0 $  by $ \left\{ (x,y) \, | \, y \in \mathbb  R_+^* \right\} $, and 
\item the half lines $V_x^- $ given for all $x \geq 0 $ by $ \left\{ (x,y) \, | \, y \in \mathbb  R_-^* \right\} $.
\end{enumerate}
\end{example}

\begin{example} \label{exa:magnetic} \textit{The magnetic partition} is an example of partitionifold that behaves badly, although it is ``regular'' in the sense that all its leaves have the same dimension.
It is given as follows:
\begin{center}
\includegraphics[width=5cm]{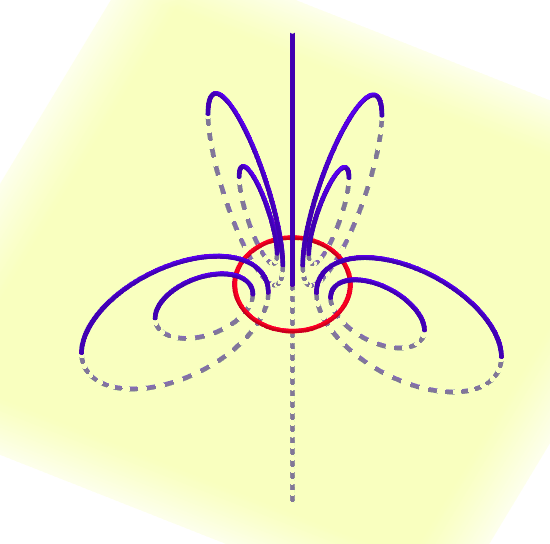}
\end{center}
These can be seen as being the lines of a magnetic field generated by an electric current in the red circle, to which the red circle itself is added. 
\end{example}
\begin{remark} 
The reader used to regular foliations will notice that  with the Bioriented partition or the Magnetic partition are partitionifolds which are
\begin{enumerate}
    \item similar to regular foliations in the sense that all leaves are submanifolds all of the same dimension~$1$,
    \item but they are still  not regular foliations in a neighborhood of the point $ (0,0)$, respectively the red circle, for the bioriented partition, respectively the magnetic partition.
\end{enumerate}
\end{remark}

\begin{example}
\noindent
French speakers may also look at the Agrégation de Mathématiques of 1998, ``Sujet de mathématiques générales'': Its first part is dedicated to the construction of a partitionifold on $\mathbb R^3$ whose leaves are all circles of non-zero radius. 
It is of course not a regular foliation.
\end{example}

\begin{example}\label{exa:lasagn}
\textit{``Isolated lasagna in a dish of spaghettis''.} Consider the partitionifold on $M= \mathbb R^3 $ with coordinates $ (x,y,z)$ whose leaves are defined to be:
\begin{enumerate}
    \item The plane $z=0$ (the ``isolated lasagna'' in red). This is the only leaf of dimension $2$.
    \item The straight lines  parallel to the $x$-axis (the ``spaghettis'' in blue) and not contained in the plane $z=0$. 
\end{enumerate}

\begin{center}
\includegraphics[width=5cm]{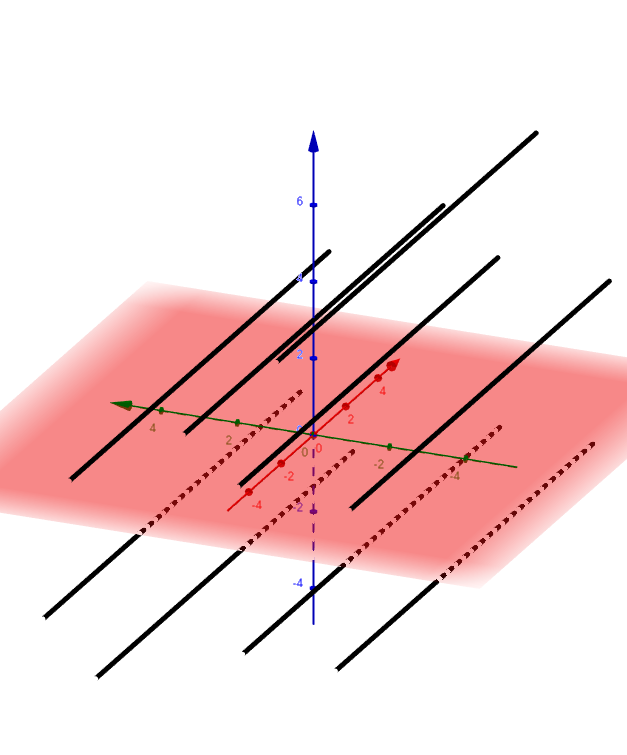}
\end{center}
\flushright{$\square$}
\end{example}

\noindent
To any partitionifold, one can associate 
 a very natural subspace of vector fields, namely those which are tangent to leaves.

\vspace{.5cm}

\begin{notations}{Vector fields tangent to every leaf}{Not:tgtvect}
Let $L_\bullet $ be a partitionifold on $M$.
 We denote by $\mathfrak T (L_\bullet) \subset \mathfrak X(M)$ the sub-sheaf\footnote{The reader unfamiliar or hostile to sheaves can define instead $\mathfrak T (L_\bullet)$ to be the sub-spaces of vector fields satisfying $X_{|_\ell} \in T_\ell L_\ell $ for all $ \ell \in M $. It is sometimes convenient to add the assumption ``compactly supported'' vector fields. Sheaves are only necessary while working within the framework of complex or real analytic geometry.} of vector fields tangent to all leaves, i.e., such that $X_{|_\ell} \in T_\ell L_\ell $ for all $ \ell$ in the open subset on which $X$ is defined. We say that such a vector field is \emph{tangent to the partitionifold $L_\bullet$}.
 
\end{notations}

\vspace{.5cm}
\noindent
It is routine to check that
$\mathfrak T (L_\bullet) \subset \mathfrak X(M)$
is a sub-module over the relevant algebra of functions. 
Also, it is a Lie subalgebra.
\vspace{.2cm}
\normalfont
Consider a curve  $ \gamma: I \to M$ whose  derivative is tangent to the leaf to which it belongs, i.e., such that for all $t \in I$: \begin{equation} 
\label{eq:speedtangent}\frac{d \gamma(t) }{dt} \in T_{\gamma(t)} L_{\gamma(t)}. \end{equation}
Is it true that this curve can not ``jump from a leaf to another leaf''? The answer is no, and the next exercises give a counter examples.

 \begin{exo}
 Consider the smooth partitionifold on $M=\mathbb R$ given by the three subsets
  $$ \mathbb R_-^*, \{0\}, \mathbb R_+^* .$$
 Show that the curve $t\mapsto t^3$ satisfies \eqref{eq:speedtangent} but is not contained in a single leaf.
 \end{exo}

\begin{exo}
Consider the bioriented partition on $\mathbb R^2 $ of Example \ref{exa:WeirdStraigtLines}. Show that the curve 
 $$   t \mapsto  \left\{\begin{array}{rr} \left(e^{-1/t^2},0\right)  & \hbox{ if $ t \leq 0$}\\ \left(0,e^{-1/t^2}\right)  &\hbox{ if $ t \geq 0$}  \end{array}f \right.  $$
 satisfies \eqref{eq:speedtangent} but is not contained in a single leaf.
\end{exo}

\noindent
The previous exercises show that curves may satisfy Equation \eqref{eq:speedtangent} and still jump from leaves to leaves.

But integral curves of vector fields in $\mathfrak T (L_\bullet)$, who automatically satisfy  Equation \eqref{eq:speedtangent}, can not jump from a leaf to a leaf, as we now see. 

\vspace{.5cm}

\begin{propositions}{Not jumping from leaves to leaves}{integralcurves}
Let $L_\bullet $ be a partitionifold on $M$.
 An integral curve $\gamma(t) $ of a vector field $X \in \mathfrak T (L_\bullet) \subset \mathfrak X(M)$ tangent to the partitionifold is always contained in one fixed leaf.
\end{propositions}

\begin{proof}
The statement is not obvious. The difficulty is that a vector field $X$ may be tangent to a submanifold $L\in M $, but its flow may not preserve it\footnote{For instance, consider the Euler vector field  $ X =\sum_{i=1}^nx_i \frac{\partial}{\partial x_i}$ on $M=\mathbb R^n $. The time-$t$ flow is a homothety with multiplying factor $ e^t $. In particular, no open ball centered at $0$ is preserved under the flow, although they are submanifolds.}.
What is true, however, if that if a vector field $ X$ is tangent to a submanifold $L$, for any integral curve starting at $t=t_0$ from a point $\ell\in L $ is ``locally in $L$'', i.e., there is $\epsilon >0 $ such that $\gamma(t) \in  L$ if $ |t-t_0| < \epsilon$.  

In the present situation, since we are given a vector field $X$ is tangent to $L_\bullet$ at all points, this implies that an integral curve $ t \mapsto \gamma(t)$ of $X$, defined on a connected open interval $ I \subset \mathbb R$, ``locally lies in the same leaf'', i.e., for any $t_0 \in I$ there is $\epsilon >0 $ such that $L_{\gamma(t)}= L_{\gamma(t_0)}$ if $ |t-t_0| < \epsilon$. Said otherwise, $\gamma^{-1}(L) $ is an open subset of $ I$ for any leaf $ L$ of $L_\bullet $. 
Since $L_\bullet$ form a partition of $M$, $I$ is the disjoint union of the open sets $(\gamma^{-1}(L))_{L \in {\mathcal L}} $ with $ {\mathcal L}$ the set of leaves of $ L_\bullet$.
Since $I$ is connected, there exists a leaf $L$ such that  
$\gamma^{-1}(L) =I$, i.e.,
the integral curve of $ t \mapsto \gamma(t)$ must be in the same leaf on its full domain.  
\end{proof}

\noindent
\vspace{0.2cm}
For any partitionifold $L_\bullet$ on $M$, and any open subset $\mathcal U \subset M$, a partitionifold on $ \mathcal U$ is obtained by mapping $m \in \mathcal U$ to the connected component of $m $ in $L_m \cap \mathcal U$.
We denote by $\mathfrak i_{ U}^* L^\bullet $ this partitionifold and call it \emph{restriction to $\mathcal U $} of $ L_\bullet$.

\vspace{2mm}
\noindent
Given  partitionifolds $ L_\bullet$ on $M$ and $ L_\bullet'$ on $M'$,
we call \emph{isomorphism  
from  $ L_\bullet$ to $ L_\bullet'$} a diffeomorphism $\phi\colon M \to M' $ such that $ \phi(L_{m} )= \phi(L_{\phi(m)}')$ for all $m \in L $. When $M=M'$ and $ L_\bullet= L_\bullet'$, we shall speak of a \emph{symmetry of $L_\bullet $.}

\vspace{.5cm}
\begin{propositions}{Flows are symmetries}{prop:FlowsAreSymmetries}
Let $M$ be a manifold equipped with a partitionifold $ L_\bullet$.
The flow at time $t$ of a complete vector field
$X \in \mathfrak T(L_\bullet)$
tangent to $L_\bullet $ is a symmetry of~$L_\bullet $.

More generally, for a maybe non-complete vector field $X \in \mathfrak T(L_\bullet)$ tangent to  all leaves, its flow $\phi_t^X $ at time $t$, provided it is well-defined on some open subset $\mathcal U \subset L $,  is an isomorphism from the restriction of $L_\bullet $ to $\mathcal U $ to the restriction of $L_\bullet $ to $\phi_t^X (\mathcal U) $.
\end{propositions}
\begin{proof}
The first part of Proposition \ref{thm:prop:FlowsAreSymmetries} is a consequence of the second one. We therefore only prove the second part. 
Consider two points $m_0,m_1 \in \mathcal U $ that are in the same leaf of $\mathfrak i_{\mathcal  U}^* L_\bullet $, and therefore in the same leaf $L $ of $L_\bullet $. There is a smooth path $m\colon [0,1] \to \mathcal U$ starting from $m_0$ and arriving at $ m_1$ which is entirely contained in $L \cap \mathcal U$. Since integral curves can not jump from one leaf to another one by Proposition \ref{thm:integralcurves}, for every $u \in [0,1], s \in [0,t]$, the map $ s\mapsto \phi_s^X (m(u)) $ is valued in the leaf $L$. In particular, the curve
 $$ u \mapsto \phi_t^X (m(u)) $$
 is entirely contained in $L$. It is also contained in $ \phi_t^X(\mathcal U)$. Hence, $\phi_t^X(m_0) $ and $ \phi_t^X(m_1)$ are in the same leaf of $\mathfrak i_{ \phi_t^X(\mathcal U)}^* L_\bullet$. This proves the claim.
\end{proof}

\noindent
Let $L_\bullet $ be a partitionifold on $M$. Consider $S \subset M $ a submanifold, we can associate to every $s \in S $ the connected component $(L_s \cap S)_0$ of $s$ in the intersection $ L_s \cap S $. The map
 $$ \begin{array}{rcl}
  S &\to  & \hbox{$\{$Connected subsets of $S\}$} \\
 s &\mapsto& \, (L_s \cap S)_0 \end{array} 
 $$
may not be a partitionifold: it is valued in connected subsets, but not in smooth manifolds. However, it is a  classical result of differential geometry that if the intersection is clean, i.e., if for all $s \in S$:
 \begin{equation}
 \label{eq:transverse} T_s S + T_s L_s = T_s M, 
 \end{equation}
then $ L_s \cap S$ is a submanifold for every $s \in S $. It may not be connected, but
 the connected component  $(L_s\cap S)_0 $ of $s \in S$ in $  L_s \cap S$ is now a non-empty connected submanifold of $S$.  In particular, 
 \begin{equation}\begin{array}{rcl}\label{eq:restrictionParitition} 
 S & \to & \left\{ \hbox{ Connected submanifolds of $ S$ }\right\}
 \\
 s &\mapsto&  (L_s\cap S)_0
 \end{array}
 \end{equation} 
 is a partitionifold on $S$.
\vspace{0.5cm}

\begin{notations}{How to denote a restriction?}{not:restriction}
Let $M$ be a manifold equipped with a partitionifold $ L_\bullet$. A manifold $S$ satisfying Condition \eqref{eq:transverse} shall be said to intersect $L_\bullet$ \emph{cleanly}. We denote by $\mathfrak i_S^* L_\bullet $ the partitionifold on $S$ given by Equation \eqref{eq:restrictionParitition}, and call it the
\emph{restriction of $L_\bullet$ to $M$}
\end{notations}

\begin{remark}
Since open subsets of $M$ intersect cleanly any partitionifold $L_\bullet$ on $ M$, the terminology and notations of Notation \ref{thm:not:restriction} match the previous conventions.
\end{remark} 

\noindent\vspace{0.2cm}
Let us conclude this question.

\begin{questions}{Are partitionifolds a good notion of singular foliations?}{ques:parti}
No, it is not!
It is too weak a notion to satisfy any significant theorem, beside the meager ones mentioned above.
\end{questions}

\subsection{Smooth partitionifolds. Is a singular foliation  a smooth partition by submanifolds?}

\noindent
We now suggest a second notion, denoted as $(**)$ in the introduction, that we claim could be the definition of a singular foliation. In fact, it is used as such by several authors \cite{Miyamoto,DLPR,Sussmann2}. It is not the most popular definition, but it is a perfectly workable notion. 

Of course, these authors call this notion ``singular foliation'', but for clarity, we prefer to give it another name. 
Unless otherwise specified, all results of this section are valid on smooth, real analytic and complex manifolds.

\vspace{.5cm}

\begin{definitions}{A more subtle attempt: smooth partitionifolds}{smooth}

A partitionifold $ L_\bullet$ is said to be \emph{smooth} if
for every $\ell \in M$ and every tangent vector $ u \in T_\ell L_\ell $ , there exists a vector field $X$ through $u$, defined in a neighborhood $ \mathcal U$ of $m$, which is tangent to all leaves\footnote{i.e., $X|_m \in T_m L_m $ for all $ m \in \mathcal U$}.
\end{definitions}

\vspace{.5cm}

\vspace{0.2cm} \noindent 
Said differently, a partitionifold is smooth if and only if, for every $\ell \in M $, the evaluation map 
 $$ \begin{array}{rcl}  \mathfrak T(L_\bullet) & \to & T_\ell L_\ell \\  X & \mapsto & X_{|_\ell} \end{array}$$
 is a surjective linear map.

\begin{remark}
    The word ``smooth'' in the expression "smooth partitionifold is extremely confusing. We used it by analogy with smooth distributions \cite{DLPR,Stepan1,Stepan2}. 
But it is confusing because it is defined in the holomorphic or real analytic contexts as well. We kept the name, however, in order to be consistent with existing literature, but also because this confusion happens anyway quite often: a smooth affine variety is not a smooth variety. We hope that it will not cause further confusion.
\end{remark}

\vspace{0.2cm} \noindent 
Let us start by a few non-examples.

\begin{exo}
Show that neither the ``magnetic partition'' (Example \ref{exa:magnetic}) nor the ``isolated lasagna in a spaghetti dish'' (Example \ref{exa:lasagn})  are smooth partitionifolds. \flushright{$\square$}
\end{exo}

\begin{exo} \emph{"Pinched curves".}
Consider the partitionifold on $\mathbb R^2 $, with coordinates $x,y$ whose leaves are the graph of the function $f_\lambda\colon  x\mapsto\lambda \left( {\mathrm{th}}(\sqrt[3]{x}) +1 \right) $ with $\lambda \in \mathbb R$. For each value of $\lambda $, the graph of $f_\lambda $ is a smooth\footnote{Even if $f_\lambda $ is not a smooth function at $x=0$ for $ \lambda \neq 0$, its graph is a smooth submanifold of dimension $1$ in~$\mathbb R^2$.} submanifold of dimension $1$ in~$\mathbb R^2$.

\begin{center}
\includegraphics[width=5cm]{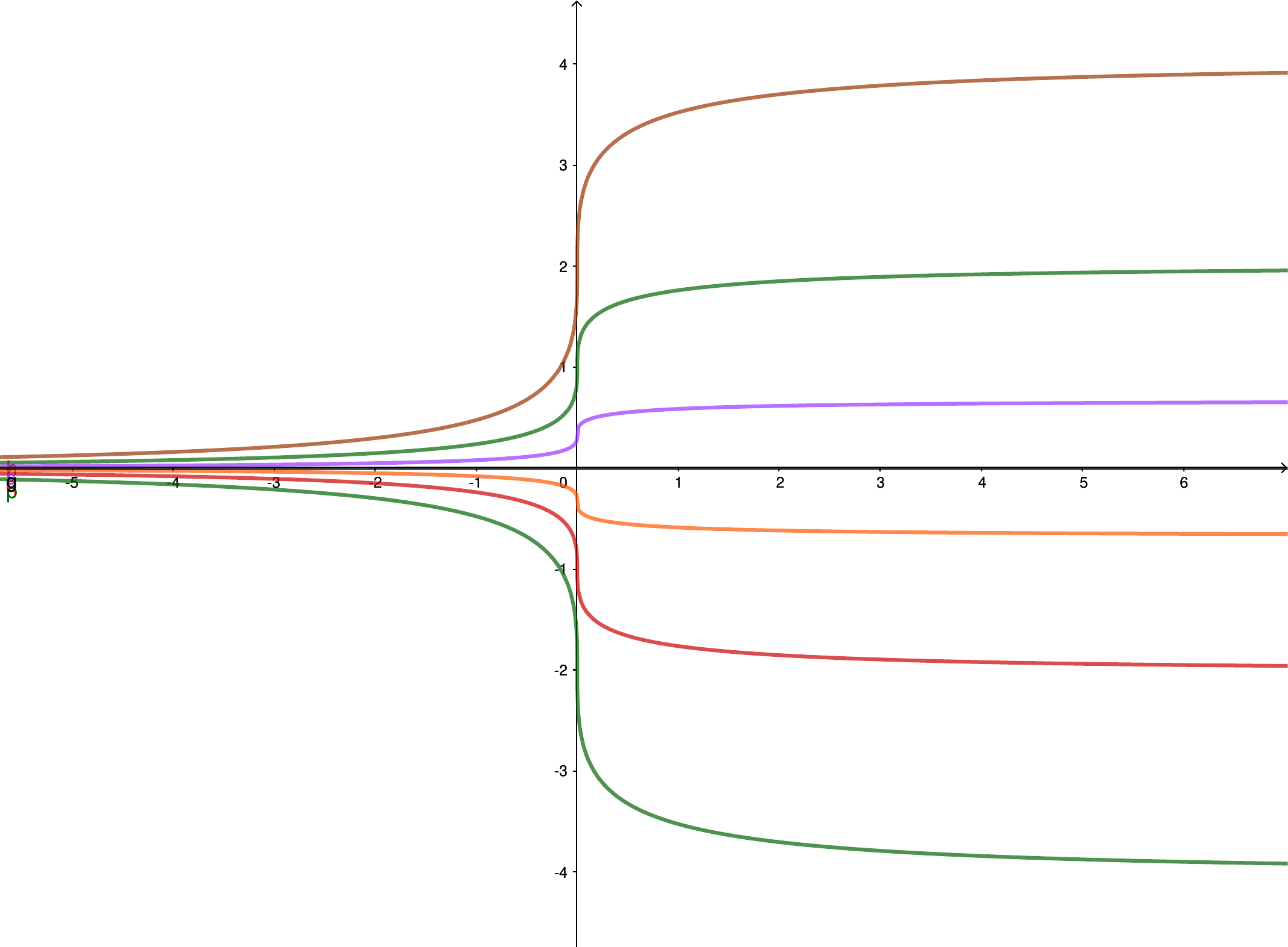}
\end{center}

\noindent
Show that this partitionifold of $\mathbb R^2 $ is not smooth. \emph{Hint:} consider a neighborhood of $ (0,0)$ and use the fact that the tangent space of $L_{(0,y)} $ is for all $y \neq 0 $ a vertical straight line.
\end{exo}

\noindent
For a smooth partitionifold, the following proposition means that singular leaves (not defined yet!) have smaller dimensions than regular ones. This explains why ``isolated lasagnas'' of Example \ref{exa:lasagn} are not smooth partitionifolds.

\begin{proposition}

\label{prop:semicontinu}
Let $M$ be a manifold equipped with a smooth partitionifold $ L_\bullet$, the function: 

$$ \begin{array}{rcl} M &\to& \mathbb N_0 \\ m & \mapsto & \dim(L_m) \end{array} $$
is lower semi-continuous\footnote{I.e., for all $k \in \mathbb N_0 $,  $\{m \in M | \dim(L_m) \geq k\}$ is an open subset in $M$}.
\end{proposition}
\begin{proof}
Let us choose a point $m_0 \in M $, let $r$ be the dimension of the leaf through $L_{m_0}$, and let $(e_1, \dots,e_r)$ be a basis of~$T_{m_0} L_{m_0} $.
By assumption, there exist $r$ vector fields $ X_1, \dots, X_r$ through $(e_1, \dots,e_r)$, defined in a neighborhood $\mathcal U $ of $m_0$ and tangent to all leaves. 
They are therefore independent at each point of a sub-neighborhood $ \mathcal U' \subset \mathcal U$, so that  $\dim(L_m) \geq r$ for all $m \in \mathcal U' $.
\end{proof}

\vspace{0.5cm}
\noindent
Smooth partitionifolds behave much better than partitionifolds, as we will briefly show by giving several reasonable theorems that they satisfy. Their first nice property is that, for a smooth partitionifold $L_\bullet $, along a given leaf $L$, $L_\bullet $ ``always looks locally the same''. This idea that a ``traveler going along a leaf will be bored'' is an important one for singular foliations.

\vspace{.5cm}

\begin{theorems}{Along a leaf, landscape is always identical}{thm:landscape}
Two points on the same leaf of a smooth partitionifold $L_\bullet $ have neighborhoods where the restrictions of $L_\bullet $ are isomorphic.
\end{theorems}

\vspace{.5cm}
\noindent
We start with a lemma.

\begin{lemma}
\label{lem:fond}
Given any two points $\ell_0,\ell_1$ on the same leaf $L$ of a smooth partitionifold $L_\bullet $, there exists a finite number of  vector fields $X_1, \dots, X_n \in \mathfrak T (L_\bullet) $ (i.e., vector fields tangent to all leaves) such that if we apply successively the flows at time $1$ of $X_1, \dots, X_n $ to the point $\ell_0$, we obtain the point $\ell_1$. 
\end{lemma}
\begin{proof}
Let us say that $x \in M$ and $y \in M$ are reachable one from the other when there exists vector fields as in the lemma. This relation clearly defines an equivalence relation, that we call the reachable relation. It follows from Proposition \ref{thm:integralcurves} that each equivalence class of the reachable relation is contained in a given leaf of $ L_\bullet$. It suffices to check that equivalence classes are precisely the leaves of $ L_\bullet$. Since each leaf is a connected set, it suffices to prove that the equivalence classes are open subsets of a leaf. Let $ \ell \in L$ be a point. Let $ X_1, \dots, X_k$ be vector fields in $ \mathcal T_{L_\bullet}$ whose values $e_1:=\left.X_1\right|_{\ell} , \dots, e_k:= \left.X_k\right|_{\ell}$ generate $T_\ell L $. Consider the map
 $$  \begin{array}{rcl}  \mathbb K^k  & \to & L\\ ( \lambda_1, \dots, \lambda_k) & \mapsto  &  \Phi_1^{\sum_{i=1}^k \lambda_i X_i }(\ell) \end{array}  $$
The differential of this map at the point $ (0,\dots, 0)$ is given for  every $(\lambda_1, \dots, \lambda_k) \in \mathbb K^k \simeq  T_{(0,\dots, 0)} \mathbb K^k$ by
$$  (\lambda_1, \dots, \lambda_k) \mapsto  \sum_{i=1}^k \lambda_i e_i $$
In particular, this differential is surjective. By the implicit function theorem, it means that  the function is open in $L$, i.e., that every point in some neighborhood of $ \ell$ is in the same equivalence class as $\ell $ is. This completes the proof. 
\end{proof}

\begin{proof}[Proof of Theorem \ref{thm:thm:landscape}.]
Lemma \ref{lem:fond} implies that if two points $\ell_0,\ell_1 $ are in the same leaves, then there are vector fields $ X_1, \dots, X_r$ tangent to all leaves whose flow at time $1$ maps $\ell_0$ to $ \ell_1$.
Proposition \ref{thm:prop:FlowsAreSymmetries}  then implies that these flows are isomorphisms of partitionifolds. This completes the proof.    
\end{proof}

\noindent
Let $ L$ be a leaf of a smooth partitionifold $L_\bullet $. A pointed submanifold $(\Sigma,\ell) $ that intersect $L$ at $\ell$  (i.e., $ \Sigma \subset M$ is a submanifold and  $\ell \in \Sigma \cap L$) is said to  be a \emph{$L_\bullet$-cut} of $L$ if 
\begin{enumerate}
\item $ \Sigma$ is transverse to $ L$ at $ \ell$, i.e.,
$ T_\ell \Sigma \oplus T_\ell L =T_\ell M, $ and
\item  $ \Sigma$ intersects $ L_\bullet$-cleanly.
\end{enumerate}
\begin{lemma}
 Any pointed submanifold $(\Sigma,\ell) $ transverse to $ L$ at $ \ell$ 
admits a neighborhood $ \Sigma'$ of $\ell $ such that $(\Sigma',\ell)$ is a $L_\bullet$-cut of $L$.
\end{lemma}
\begin{proof}
Let $ X^1, \dots, X^k \in \mathfrak T(L_\bullet)$ be vector fields tangent to all leaves whose evaluations at $\ell $ form a basis of $ T_\ell L$. There exists a neighborhood $\mathcal V$ of $\ell $ in $\Sigma $ such that for all $ \sigma \in \mathcal V$ 
 $$ T_\sigma \Sigma \oplus \langle  X^1_{|_\sigma}, \dots, X^k_{|_\sigma}\rangle = T_\sigma M. $$ This implies $ T_\sigma \Sigma + T_\sigma L_\sigma =T_\sigma M$, which is precisely the definition of intersecting $ L_\bullet$ cleanly.
\end{proof}

\vspace{0.2cm}
\noindent
For any $L_\bullet$-cut $(\Sigma,\ell) $ of a leaf $L$,  $\mathfrak i^*_{ \Sigma } L_\bullet  $ is a smooth partitionifold, that we call a \emph{transverse partitionifold} of the leaf $L$.   
\vspace{.5cm}

\begin{corollaries}{The germ of a slice transverse to a leaf}{cor:germ}
Let $M$ be a manifold equipped with a smooth partitionifold $ L_\bullet$. Any two $ L_\bullet$-cuts  of a given leaf $L$ have neighborhoods  on which their restrictions are isomorphic
\footnote{
More precisely, for any two pointed submanifolds $ (\Sigma_1,\ell_1)$ and $ (\Sigma_2,\ell_2) $ transverse to the same leaf $ L$, there exists neighborhoods $\mathcal U_1 \subset \Sigma_1 ,\mathcal U_2 \subset \Sigma_2$  of $\ell_1, \ell_2 $ and an isomorphism $$\xymatrix{\mathfrak i^*_{ \Sigma_1 \cap \mathcal U_1} L_\bullet    \ar[r]^{\sim}& 
\mathfrak i^*_{ \Sigma_2 \cap \mathcal U_2} L_\bullet }.$$}.
\end{corollaries}
\begin{proof}
This comes from the fact that the vector fields $ X_1, \dots, X_d \in  \mathfrak T (L_\bullet)$ in Lemma \ref{lem:fond} can be chosen such that the composition of their flows (which are symmetries of $ L_\bullet$) not only map $\ell_1 $ to $\ell_2 $ but also $\Sigma_1 $ to $\Sigma_2 $, at least in a neighborhood of $ \ell_1$.
A similar argument is presented in the proof of Theorem \ref{thm:transverseFol}: we refer the reader to that proof for more details.
\end{proof}

\noindent
Corollary \ref{thm:cor:germ} implies that it makes sense to speak of the  transverse model of a leaf. To be more precise, it can be defined as follows.
Let $d$ the dimension of the manifold and $r$ the dimension of the leaf $L$
We call \emph{representative of the transverse partitionifold of $L$} a pair $ (\mathcal U,L^\mathcal U)$ 
with $ \mathcal U$ a neighborhood of $0$ 
such that there exists a $ L$-cut $ (\Sigma,\ell)$ and an isomorphism isomorphic  $ (\Sigma,\mathfrak i_\Sigma^* L_\bullet) \simeq (\mathcal U,L^\mathcal U)$  mapping $ \ell$ to $0$.
More precisely, consider 
pairs $ (\mathcal U,L^\mathcal U)$ 
with $ \mathcal U$ a neighborhood of $0$  $\mathbb K^s $,
and then identify two such pairs
 $$ (\mathcal U,L^\mathcal U) \sim  (\mathcal V,L^\mathcal V) $$
if they have neighborhoods $\mathcal U' $ and $ \mathcal V'$ of $0$ on which the restricted smooth partitionifolds are isomorphic, through an isomorphism that preserves $0$. We call \emph{germs of partitionifolds at $0$ in dimension $s$}
 the equivalence classes of this equivalence relation.
To any leaf $L$ of a partitionifold, one associates a canonical germ of partitionifolds at $0$ in dimension $d-k$

with $d = {\mathrm{dim}}(M)$ and $ k= {\mathrm{dim}}(L)$:
It is by definition the unique class which admits a representative isomorphic, as a smooth partitionifold,  to the smooth partitionifold of one $ L_\bullet$-cut of $ L$.
  We will have very similar theorems for in Section \ref{sec:traveling} for the definition of singular foliation that we will choose.

\vspace{2mm}
\noindent
A classical result of classical theorem in Control Theory called Nagano-Sussmann Theorem, which is a continuation of Rashevsky–Chow Theorem, see Chapter 5 in \cite{zbMATH02105696},  justifies the introduction of smooth partitionifolds.

Let us describe it.
Given a finite family of $d$  vector fields $ X_\bullet := X_1, \dots, X_d$ on a manifold $M$, consider the equivalence relation on $M$ generated by $ m \sim m'$ if $m$ and $m'$ are on the same integral curve of  $X_1, \dots, X_d$.  Two points in the same equivalence class of this relation are said to be \emph{reachable one from the other}, and we call this class the \emph{reachable equivalence relation}.
Let us construct this class differently. Assume that the vector fields
 $X_1, \dots, X_d$ are complete for the sake of simplicity. Consider 
the subgroup $ {\mathrm{Diff}}_{X_\bullet}(M) $ of the group  $ {\mathrm{Diff}}(M) $ of diffeomorphisms
of $M $ generated by the flows at time $t\in \mathbb R $ of the vector fields $X_1 , \dots, X_d$. In equation, a diffeomorphism $\phi  $ of $M$ belongs to $ {\mathrm{Diff}}_{X_\bullet}(M) $ if there exists an integer $n \geq 1$, a $n$-tuple  $ i_1 \dots,i_n \in \{1, \dots, d\}$ and a $n$-tuple  $ t_1 \dots,t_n \in \{1, \dots, d\}$ such that $$ \phi = \Phi_{X_{i_1}}^{i_1} \circ \dots, \circ \Phi_{X_{i_n}}^{i_n} .$$
When the vector fields are not complete, 
$ \mathrm{Diff}_{X_\bullet}(M) $ becomes what is called a pseudo-group, i.e., pairs made of an open subset of $M$ and a diffeomorphism defined on that open subset.
In any case, it makes sense to speak of the orbits of  $ \mathrm{Diff}_{X_\bullet}(M) $: two points $m,m'\in M$ belong to the same orbit of  of $ (X_1, \dots, X_d)$ if there exists $\phi \in  \mathrm{Diff}_{X_\bullet}(M) $ such that $ \phi(m)=m'$. The orbits are precisely the equivalence classes of the previous equivalence relation.

\vspace{0.5cm}

\begin{theorems}{Nagano–Sussmann in terms of smooth partitionifolds}
{thm:anydistributionGivessmooth}
Let $ X_\bullet =(X_1, \dots, X_d)$ be a finite family of vector fields on a manifold $M $. The equivalence classes of the reachable equivalence relation above 
form a smooth partitionifold $L_\bullet$ of $M$. 
\end{theorems}

\vspace{.5cm}

Nagano–Sussmann's theorem \cite{Nagano,Sussmann1} says in fact something more. We say that a subspace $ \mathcal G \subset \mathfrak X(M) $ of vector fields is \emph{completely tangent to a smooth partitionifold $L_\bullet $} if for every $m \in m$, and every $ u \in T_m L_m$ tangent to the leaf through $L$ at $m$, there exists $X \in \mathcal G$ such that $ X(m)=u$. 
Consider the vector fields obtained by considering: 
\begin{enumerate}
\item the Lie subalgebra $ {\mathrm{Lie}}(X_\bullet)$ of $\mathfrak X(M) $ generated by the vector fields in $X_\bullet $, i.e., the space generated by all vector fields of the form
  $$  \left[\cdots \left[X_{i_1},X_{i_2} \right] \dots , X_{i_n}\right]  $$
  for all $n \geq 1 $, $ i_1, \dots,i_n \in \{1, \dots, d\} $
\item then consider all possible push-forwards through all the diffeomorphisms $\phi \in \mathrm{Diff}_{X_\bullet}(M)  $ of vector fields in $ {\mathrm{Lie}}(X_\bullet)$, 
\item then all possible linear combinations, with coefficients in functions, of the previously defined vector fields. 
 \end{enumerate}
 In the real analytic or holomorphic setting, step 2 is not required.
By proceeding as in steps 1,2,3 before, we obtain a sub-sheaf of $ \mathfrak X(M) $ that we denote by $ {\mathrm{Closure}}(X_\bullet) $. It is a module over functions, and it is closed under Lie bracket. It is of course included into $ \mathfrak T(L_\bullet)$. 

\begin{remark}
    \normalfont
    In general, the inclusion ${\mathrm{Closure}}(X_\bullet) \subset \mathfrak T(L_\bullet)$ is strict: take for instance on $ M=\mathbb R$ the vector field $ X_\bullet$ to be the family with one element, namely $X_1= x^2 \frac{\partial}{\partial x}$. $ {\mathrm{Closure}}(X_\bullet) $ is the $ \mathcal C^\infty(M)$-module generated by $X_1$ while  $ \mathfrak T(L_\bullet)$ is made of all vector fields on $\mathbb R $ vanishing at $0$.
\end{remark}

\vspace{.5cm}
\noindent
Nagano-Sussmann theorem is of course an extremely strong result, and would be an excellent reason to work with smooth partitionifolds as a definition of singular foliations.  We will however not it for reasons that we discuss now.
\vspace{0.5cm}

\begin{questions}{Conclusion: Are smooth partitionifolds a good definition of singular foliations?}{ques:smoothparitioniolds}
Our opinion is that it is fine as a definition. They are called ``singular foliations'' by several authors, even recently \cite{Miyamoto}. 
Nagano-Sussmann Theorem above is an excellent reason to use them.
 But it not the definition which is most commonly used at the moment. We will see that they are not practical to deal with differential operators, which explains why non-commutative geometers do not use it. Most importantly, the theory one could develop out of smooth partitionifolds would not be so different from the one we will develop in Section \ref{sec:concensus}: they are different, but parallel.
\end{questions}

\vspace{0.5cm}
\noindent
In addition, here is an oddity that we want our ``singular foliations'' to avoid, presented as an exercise. In Section \ref{sec:concensus}, we will explain why having non-locally finitely generated modules is an issue.

\begin{exo} \label{exo:infinite} \textit{``Vector fields tangent to the leaves are not finitely generated''}.
On $M=\mathbb R $, consider the partitionifold whose $0$-dimensional leaves are $\{1\},\{\frac{1}{2}\}, \{\frac{1}{3}\}, \dots, \{\frac{1}{n}\},\dots  $ and  $ \{0\}$ and whose $1$-dimensional leaves are the open intervals bounded by these points.
Show that vector fields tangent to $L$ are not a finitely generated module over $\mathcal C^\infty(M) $, and that there is no neighborhood $\mathcal U $ of $0$ on which such vectors form a locally finitely generated  $\mathcal C^\infty(\mathcal U) $-module.
\end{exo}

\subsection{Is a singular foliation an involutive distribution?}

\noindent
Before presenting the definition of a singular foliation that we intend to use, let us introduce an alternative manner to define partitionifolds of non-constant dimensions.
It is extremely classical in differential geometry \cite{MR4454136,Moerdijk} that a regular foliation may be defined as being an integrable sub-vector bundle $D \subset TM $. It is therefore tempting to allow the fibers of the vector bundle $D $ to be of non-constant dimension, as long as its sections are closed under the Lie bracket of vector fields:

\vspace{.5cm}

\begin{definitions}{Integrable singular distributions}{IntegrableSingDistr}
A singular distribution on a manifold $M $ is a map  ${\mathfrak D} $ associating to a point $m\in M$ a subspace ${\mathfrak D}_m  \subset T_m M $. 
A singular distribution ${\mathfrak D} $ is said to be:
\begin{enumerate}
    \item \textit{involutive} when $[\Gamma({\mathfrak D}) ,\Gamma({\mathfrak D}) ] \subset \Gamma({\mathfrak D})  $, 
where $\Gamma({\mathfrak D}) \subset \mathfrak X(M) $ is the $ \mathcal C^\infty(M)$-module of vector fields $X$ such that $ X_m \in {\mathfrak D}_m $ for all $m \in M$.
\item \textit{integrable} when there exists a partitionifold $L_\bullet$ such that for all $ m \in M$, $T_m L_m = \mathfrak D_m $.
\end{enumerate}
\end{definitions}

\vspace{.5cm}

\begin{exo}
Show that for any partitionifold $L_\bullet $ on $M$, the map
 $$ \mathfrak D : m \mapsto T_m L_m $$
 
 is an involutive and integrable singular distribution. 
\end{exo}

\begin{exo}\label{exo:verybadinterDistr}
Let $ M$ be a manifold and $m_0 $ a point.
Show that the map
 $$  m \mapsto \left\{ \begin{array}{rr} T_{m_0} M  & \hbox{ if $m=m_0$} \\  0_{T_m M} & \hbox{ if $m\neq m_0$} \end{array} \right.$$
 is an involutive but non-integrable singular distribution. (Here, $0_{E}$ stands for the zero element of a vector space $E$).
\end{exo}

\noindent
For a given singular distribution, what would be the proper definition of a leaf? 
There is a natural manner to define leaves for a singular distribution $\mathfrak D $.
Consider the equivalence relation on $M$ generated by the relation $ x_0 \sim x_1 $ if  there exists a path of class $ C^1$ such that 
 \begin{equation}\label{eq:equiv_relation0}  x(0) =x_0, x(1)=x_1 \hbox{ and } \frac{d}{dt}x(t) \in \mathfrak D_{x(t)} \hbox{ and } \frac{d}{dt}x(t) \neq 0 .\end{equation}
 Equivalently, one could define an equivalence relation as follows: call integral submanifold of $\mathfrak D $ a submanifold $\Sigma$ such that $T_\sigma\Sigma \subset \mathfrak D_\sigma   $ for all $ \sigma \in \Sigma$. We could then consider the equivalence relation generated by the relation $x_0 \sim x_1 $ if there exists an integral submanifold containing both $x_0 $ and $x_1$.  
 The classes of this equivalence relation can not decently be called leaves, because they are not submanifolds, as seen in the following Exercise.
 
 \begin{exo}
\label{exo:trumpet}
Here is an example (the ``trumpet foliation'') of an involutive singular distribution for which one class of the equivalence definition \eqref{eq:equiv_relation0} is not a manifold.
Take $M= \mathbb R^2$ with coordinates $(x,y)$. Let $k(x) = e^{-1/x}$ for $ x>0$ and $k(x)=0 $ for $ x  \leq 0$. 
Divide $\mathbb R^2 $ in three zones:
$$ North := \{ y \geq k(x) \} ,  Middle := \{x > 0 \hbox{ and }  -k(x) < y < k (x) \}, South := \{ y \leq - k(x) \}   $$
Define a singular distribution by:
 $$ \mathfrak D_m = \left\{  \begin{array}{ll} \langle(1,k'(x))\rangle  & \hbox{ for } m \in North \\  T_m \mathbb R^2 & \hbox{ for } m \in Middle  \\  \langle(1,-k'(x))\rangle  & \hbox{ for } m \in South \end{array} \right. $$

\begin{center}
\includegraphics[width=7cm]{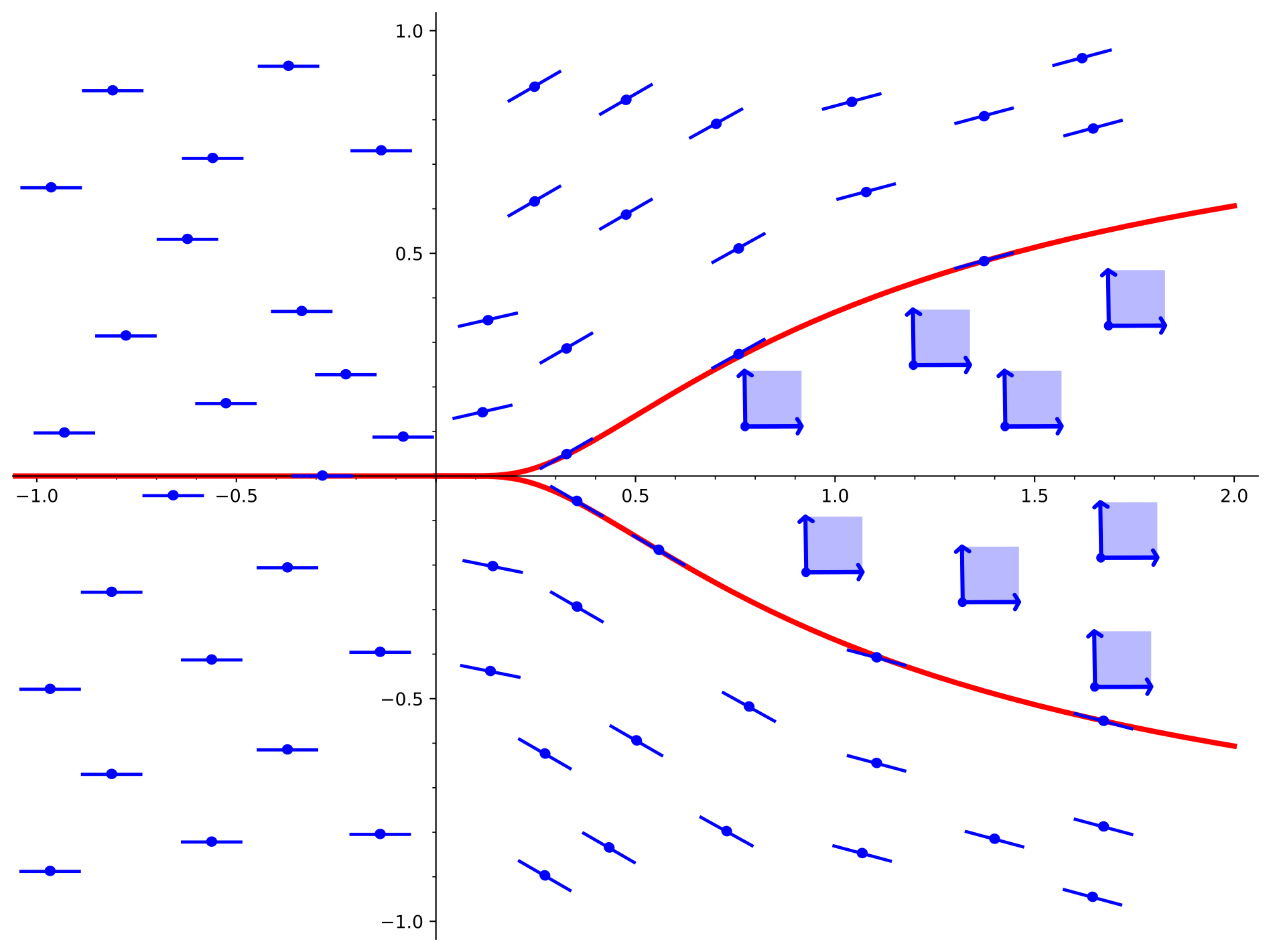}
\end{center}
 
\begin{enumerate}

    \item  Show that $ \mathfrak D$ is involutive
    \item Show that $\mathfrak D$ is not integrable (\emph{Hint}  Show that the equivalence class of $(0,0) $  is  $  \{y=0\} \cup \overline{Middle} $, which is not a manifold).
\end{enumerate}

\end{exo}

\vspace{2mm}
\noindent
It is clear that we have to avoid situations like the one in Exercise \ref{exo:verybadinterDistr}, as well as the tangent spaces of the partitionifolds of Examples \ref{exa:magnetic} and \ref{exa:lasagn} (``magnetic foliation'' or ``isolated lasagnas''). For that purpose, we will impose a second condition, similar to Definition \ref{def:smooth} of smooth partitionifolds.

\vspace{.5cm}

\begin{definitions}{Smooth singular distributions}{Second attempt}
A singular distribution ${\mathfrak D}$ is said to be smooth if for every point $ m \in M$ and $u \in \mathfrak D_m $, there exists a vector field $X \in \Gamma(\mathfrak D)$ through $u$.
\end{definitions}

\begin{exo}
Let $L_\bullet $ be a smooth partitionifold of $M$.
Consider the singular distribution ${\mathfrak D}_L : m \mapsto T_m L_m $.
\begin{enumerate}
    \item Show that it is integrable and involutive,
    \item and smooth
    \item and that the flow of any section in $\Gamma({\mathfrak D}_L ) $ preserves ${\mathfrak D}_L$.
\end{enumerate}
\end{exo}

\begin{exo}
 For an involutive and integrable smooth distribution, show the classes of the equivalence relation \eqref{eq:equiv_relation0} are precisely the leaves of $L_\bullet $
\end{exo}
\vspace{2mm}
\noindent
 The two exercises above seem to indicate that smooth involutive singular distributions are a ``good'' notion. There is however a type of example which is quite annoying:

\begin{exo} 
\label{exo:comb}
Here is a non-integrable distribution, the ``\emph{infinite comb}'', that will be a source of several counter-examples.
Consider on $M = \mathbb R^2 $ with variables $(x,y)$ the singular distribution given by
\begin{equation}  \mathfrak D_{(x,y)} = \left\{ \begin{array}{cll}  \langle \frac{\partial}{\partial x} \rangle & \hbox{if $x \leq 0 $}  & \hbox{i.e., ``Dimension $1$ in the black zone - and horizontal.''} \\  \langle \frac{\partial}{\partial x},
\frac{\partial}{\partial y}\rangle  &\hbox{if $x > 0 $}
 & \hbox{i.e., ``Dimension $2$ in the red zone.''}\\
\end{array}\right. \end{equation}
\begin{center}

\includegraphics[scale=0.2]{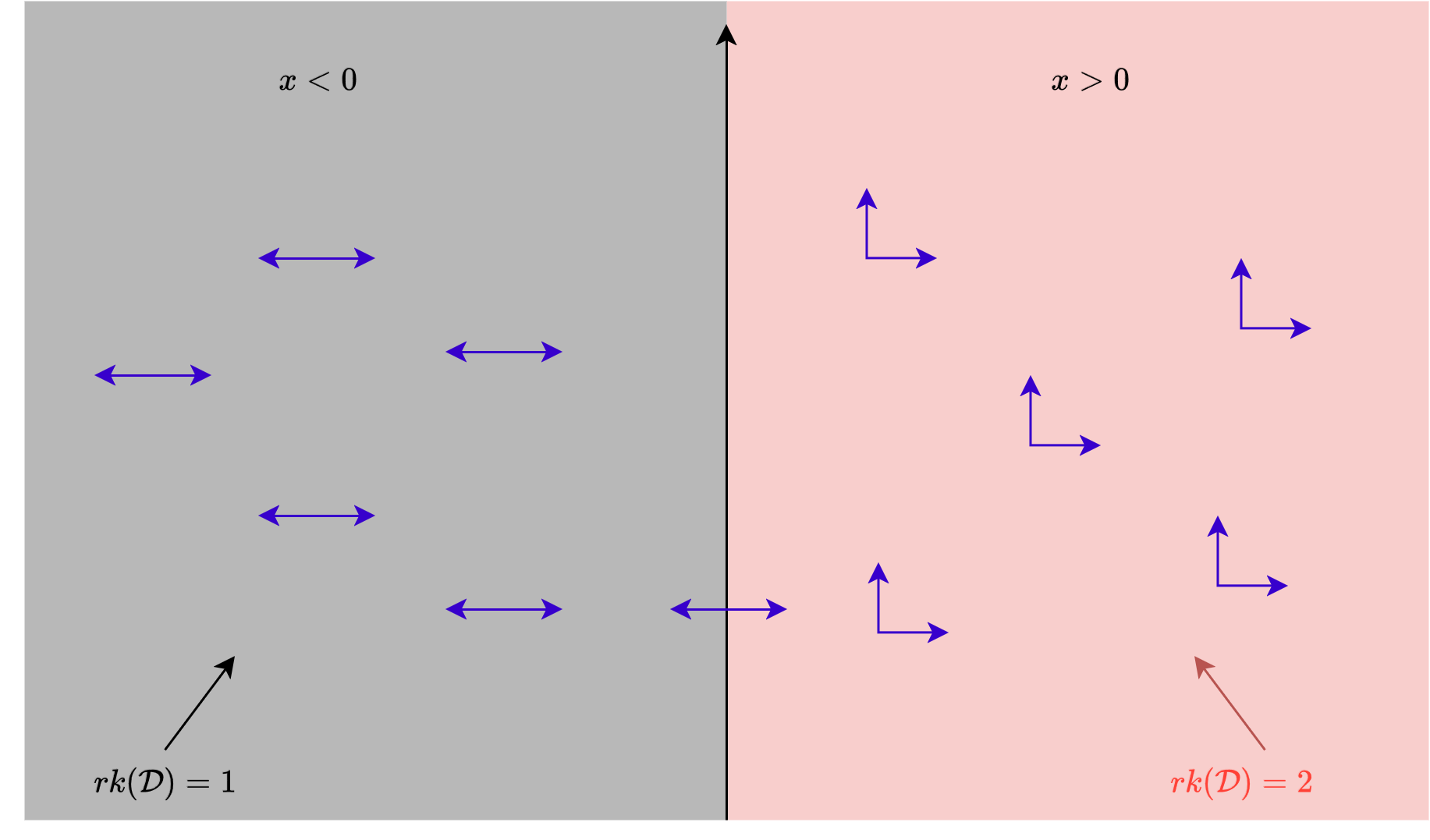}
\end{center}

\begin{enumerate}
    \item Show that the singular distribution $ \mathfrak D$ is smooth.
    \item Show that the singular distribution $ \mathfrak D$ is involutive.

    \item Show that any two points in $\mathbb R^2 $ are in the same equivalence class of the relation \eqref{eq:equiv_relation0}.
    \item Show that it is not integrable.
\end{enumerate}
\end{exo}

\noindent
\vspace{2mm}
The last exercise shows that smooth involutive singular distributions may not be integrable. Another issue of the notion of smooth integrable distribution is that the flow of a complete vector field in $\Gamma(\mathfrak D) $ may not be a symmetry of $\mathfrak D $. Again, the infinite comb is a counter-example:

\begin{exo}
Show that the vector field $-\frac{\partial}{\partial x} $ belongs to $\mathfrak D $ but that its flow does not preserve $\mathfrak D $.
\end{exo}

In the context of smooth differential geometry, Stefan-Sussmann theorems is a way out of these counter-examples.
\vspace{.5cm}
\vspace{0.Cm}

\begin{theorems}{Stefan-Sussmann Theorem (1973) \cite{Stepan1}-\cite{Stepan2}-\cite{Sussmann1}-\cite{Sussmann2}-\cite{DLPR}}{Stefan-Sussmann}
Let $\mathfrak D $ be an involutive smooth distribution on a smooth manifold $M$. The following items are equivalent:
\begin{enumerate}
    \item[(i)] $\mathfrak D$ is integrable.\footnote{I.e., there exists a smooth partitionifold $L_\bullet $ such that $T_m L_m = \mathfrak D_m$ at all points $m \in M$.} 
    \item[(ii)] There exists a $ \mathcal C^\infty(M)$-module $\mathcal F $ of vector fields such that:
    \begin{enumerate}
        \item $\mathcal F $  generates $ \mathfrak D$. \footnote{I.e., for all $m \in M$ and $u \in \mathfrak D_m$, there exists $X \in \mathcal F$ with $ X_{|_m}=u$.} 
        \item The flow $ \phi_t^X$ of any vector field  $X \in \mathcal F $ preserves $\mathfrak D $. \footnote{i.e., $\phi^X_t(\mathfrak D_m)= \mathfrak D_{\phi^X_t(m)} $ for all $m\in M, X \in \mathcal F, t \in \mathbb R $ for which the flow is well-defined in a neighborhood of $m$.}
    \end{enumerate}
\end{enumerate}

\end{theorems}  

\vspace{2mm}
\noindent
In \cite{DLPR}, it is proven that one can add a condition (c) in item \emph{(ii)} of Theorem \ref{thm:Stefan-Sussmann}: one can also impose that $ \mathcal F$ is locally finitely generated, and the equivalence of items \emph{(i)} and \emph{(ii)} still holds true. In general, however, we cannot assume $ [\mathcal F,\mathcal F] \subset \mathcal F$, so $ \mathcal F$ can not be assumed to be what we will call a “singular foliation” in the next section.

Let us conclude this section:

\vspace{.5cm}

\begin{questions}{Are involutive smooth singular distributions a good notion of singular foliations?}{ques:intDistr}
No, it is not!

\vspace{1mm}
\noindent
Foliations should have leaves, and there is an issue with the notion of leaves. See discussion  about infinite combs, which are not integrable.

\vspace{1mm}
\noindent
However, with an additional condition on flows, Stefan-Sussmann Theorem \ref{thm:Stefan-Sussmann} grants integrability. But this condition is hard to check in a concrete manner.
\end{questions}

\section{Singular foliations through vector fields: a consensus definition?}
\label{sec:concensus}

\noindent
We now introduce the definition of singular foliations that we will use effectively. Essentially, the singular foliations will be for us a subspace $\mathcal F $ of the space of vector fields, a subspace supposed to “behave like” vector fields tangent to the leaves of a partition of $M$ by sub-varieties. This is the view that seems to have prevailed for some time. It certainly predominates today if we count the number of publications. It is therefore tempting to say that there is a growing consensus around the definition that we present in the lines below. However, it is not universally used, see e.g. Miyamoto's recent Ph.D.  \cite{Miyamoto,zbMATH07807748}. 
Regardless, this vector field view is not \emph{that} different from the "smooth partitionifold" view: both theories are in some way parallel. Although they are not so easily comparable, in the sense that there is no functorial relationship between them, and notably no inclusion, they are parallel in the sense that both theories have the same shape. There is a natural way to translate a theorem using our definition into a theorem about smooth partitionifolds, and the resulting statement is most likely true.

\vspace{2mm}
\noindent
Note that if what we call a singular foliation is only a subspace of the space of vector fields satisfying several conditions, the leaves are not a priori given in the definition. The existence of leaves will now be a theorem, and in fact a very non-trivial theorem going back to Hermann's Theorem \ref{thm:Hermann}.

\vspace{2mm}
\noindent
We will present the smooth setting apart from the real analytic and holomorphic ones, although there is indeed a manner to merge them all. We will also present a purely algebraic framework, which defines
singular foliations on an affine variety, or even on an arbitrary commutative algebra.

\subsection{Definition of a singular foliation: the smooth case}

\label{sec:defsmooth}

\vspace{2mm}
\noindent
We now present a definition which has become the object of consensus in non-commutative geometry, see e.g. \cite{AS,AS11,AZ2,Debord,OmarMohsen,AMY} by Androulidakis, Debord, Mohsen, Skandalis, Yuncken, and Zambon - to quote a few. This branch of non-commutative geometry works essentially within the context of \underline{smooth} differential geometry, so we start in this context.

\vspace{0.2cm}
Recall that the \emph{support} of a smooth function or a smooth vector field on a  manifold $M$ is the closure of the subset of $M$ of points where it is not equal to zero.
We denote by $ \mathcal C^\infty_c(M)$ and $\mathfrak X_c(M) $ the space of compactly supported smooth functions and vector fields, respectively. Recall that $\mathfrak X_c(M) $ is a module over the algebra $\mathcal C^\infty(M) $.

\vspace{.5cm}

\begin{definitions}{Singular foliations on smooth manifolds}{consensus}
A singular foliation on a smooth manifold $M$ is a subspace
$\mathcal F \subset \mathfrak X_c(M)  $ which 
\begin{enumerate}
    \item[$(\alpha)$] is involutive,
    
    \item[$(\beta)$] is stable under multiplication by an element of $\mathcal C^\infty(M)$, 
    
     \item[$(\gamma)$] is locally finitely generated.
     
\end{enumerate}
\end{definitions}
\vspace{0.5cm}
\begin{remark}
A subspace $\mathcal F \subset \mathfrak X_c(M) $ that satisfies only  conditions $(\alpha)$ and $(\beta)$ is  a \textit{Lie-Rinehart subalgebra}\footnote{ 
In fact, generic Lie-Rinehart subalgebras of vector fields behave very badly in the smooth case and there is not much to say about them. They are more interesting in the complex case. But we need the concept for pedagogical reasons, in particular in order to explain why we impose condition $\gamma$.} of the Lie-Rinehart algebra of compactly supported vector fields. (Please ignore this remark is you never heard of Lie-Rinehart algebras).
\end{remark}

\subsubsection{Definition \ref{def:consensus} explained item by  item.}

\begin{enumerate}
    \item[$(\alpha)$] ``$\mathcal F $ is involutive'' means that there is an inclusion
     $$ [\mathcal F, \mathcal F] \subset \mathcal F$$
     where $[\cdot\,, \cdot]$ stands for the bracket of vector fields. In words, condition $ (\alpha)$ means that $\mathcal F$ is a sub-Lie algebra of the Lie algebra $ \mathfrak X_c (M)$ of compactly supported smooth vector fields on the manifold $M$.
     \item[$(\beta)$] ``$\mathcal F$ is stable under multiplication by an element of $\mathcal C^\infty(M)$'' means that for all $F \in \mathcal C^\infty(M), X \in \mathcal F $, $F X  \in \mathcal F$. In algebraic terminology, it means that  $\mathcal F$ is a $\mathcal C^\infty(M)$-sub-module of the $\mathcal C^\infty(M)$-module $\mathfrak X_c(M) $ of compactly supported vector fields on $M$. 
     \item[$(\gamma)$] The meaning of ``$\mathcal F $ is locally finitely generated'' has to be made very precise. It  means that for any point $ m \in M$, there exists a finite family $X^1, \dots, X^r \in \mathcal F $ and an open neighborhood $\mathcal U$ of $x$ in $M$ such that for every $X \in \mathcal F$, there exists $ f_1, \dots, f_r \in \mathcal C^\infty(M)$ satisfying for all $ x \in \mathcal{U}$:
      $$ X_{|_x} = \sum_{i=1}^r f_i(x) \, \,   X^i_{|_x}   .$$
    \end{enumerate}

\begin{exo}
Show that in Definition \ref{def:consensus}, $ \mathcal C^\infty(M)$ could be replaced by $  \mathcal C^\infty_c(M)$ without making any difference.
\end{exo}

\vspace{1cm}
\noindent
\subsubsection{Definition  \ref{def:consensus} justified item by item.}

Now that we have explained the meaning of the three items, let us justify these items. Recall that the general idea is that $\mathcal F $ shall behave like ``vector fields tangent to the leaves''. 
\begin{enumerate}
    \item[] To start with, why do we use the space $\mathfrak X_c(M) $ of compactly supported vector fields?
    \item[] We would first like to tell the reader that the reasons are above all technical, that there is nothing very deep in this choice. Let $M$ be a non-compact manifold. The following spaces:
\begin{enumerate}
    \item compactly supported vector fields on $M$,
    \item all smooth vector fields on $M$ (compactly supported or not),
\end{enumerate}
are \emph{different} as modules over $\mathcal C^\infty(M) $. But both of them should define the \emph{same} singular foliation: we do not wish to distinguish them.
They obviously have the same leaf: $M$ itself (leaves have not been defined, but it is obvious in this case that $M$ is the only leaf, whatever it means). If we did not impose ``compactly supported'', then we would have to say ``and we will identify $\mathcal F $ with the subspace of compactly supported vector fields in $\mathcal F  $''. To avoid that two different subspaces of $\mathfrak X(M) $  define the same singular foliation, one directly requires that $ \mathcal F \subset \mathfrak X_c(M)$. 

One can get rid of all compactness assumption by using another possible definition that involves sheaves, as we will see in the next section.
\item Why do we assume $(\alpha)$, i.e., that $\mathcal F$ is integrable? If two vector fields $X,Y$ are tangent to a submanifold $L$, so is its bracket. Since $\mathcal F $ must be thought of as being a replacement of vector fields tangent to the leaves, it makes sense to require $[\mathcal F, \mathcal F] \subset \mathcal F $. 
\item Why do we assume $(\beta)$, i.e., that $\mathcal F$ is a $\mathcal C^\infty(M)$-module? Because, if a vector field $X$ is tangent to all leaves, so is $fX$ for all smooth functions $f\in\mathcal C^\infty(M)$.
\item Why do we assume $(\gamma)$, i.e., that $\mathcal F$ is ``locally finitely generated''? The idea is to avoid weird counter-examples as the infinite comb (see Example \ref{exo:comb}). Imposing locally finitely generated guaranties that leaves will make sense. This is the topic of a subsequent section (Section \ref{sec:leavesExist}). 

\end{enumerate}

\begin{exo}
Show that the space of vector fields on $ \mathbb R$ of the form $f(x) \frac{\partial}{\partial x} $ with $f$ a function vanishing with all its derivatives at $0$  satisfy items $(\alpha)$ and $(\beta)$ by not item $(\gamma) $.  
(\emph{Hint:} Look at second item in Exercise \ref{exo:finitelygeneratedex}).
\end{exo}

\begin{exo}
\label{exo:infinitecomb}
    We call\footnote{We owe this example to Thomas Strobl.} \emph{infinite comb}\footnote{This is if course related to exercise \ref{exo:comb}, which describes the singular distribution obtained by taking the distribution associated to these vector fields.} the singular foliation on $ \mathbb R^2$ made of all vector fields of the form 
    $$ F(x,y) \frac{\partial}{\partial x}  + G(x,y)  \frac{\partial}{\partial y} $$
    where $G(x,y)$ is a function with compact support in the half plane $\{(x,y) | x \geq 0 \} $ and $F(x,y)$ is a function with compact support. Here $ (x,y)$ are the coordinates on $\mathbb R^2 $.
Show that these vector fields satisfy items $(\alpha)$ and $(\beta)$ but not item $(\gamma) $.   
    \end{exo}

\vspace{0.2cm}
\noindent
The next exercise is crucial, for quite a few singular foliations are defined as families $ X_1, \dots, X_r$ that satisfy one of the equivalent conditions listed there.

\begin{exo}
\label{exo:finitelygenerated}
Let $M$ be a manifold, and let $X_1, \dots, X_r  \in \mathfrak X(M)$ be vector fields. Show that the following three items are equivalent:
\begin{enumerate}
    \item[(i)] The $\mathcal C_c^\infty(M)$-module\footnote{To be more explicit, this item means that the space of vector fields of the form $\sum_{i=1}^r f_i X_i $, with $ f_1, \dots, f_r$ compactly supported smooth functions on $M$, is a singular foliation.} generated by $ X_1, \dots, X_r$ is a singular foliation,  
    \item[(ii)] 
There exist smooth functions $ c_{i,j}^k \in \mathcal C^\infty(M) $, with $i,j,k \in \{1, \dots, r\}^3$, such that
 \begin{equation}\label{eq:ii} [X_i,X_j] = \sum_{k=1}^n c_{ij}^k X_k  \end{equation}
 for all $i,j \in 1, \dots, r $.
    \item[(iii)] 
There exists smooth functions $ \tilde{c}_{i,j}^k \in \mathcal C^\infty(M) $, with $i,j,k \in \{1, \dots, r\}^3$ satisfying
  \begin{equation}\label{eq:iii} \tilde{c}_{ij}^k = - \tilde{c}_{ji}^k \hbox{ and } [X_i,X_j] = \sum_{k=1}^n \tilde{c}_{ij}^k X_k  \end{equation}
 for all possible indices.
 \item[] {\emph{Hint:}} Show that if a family $ c_{ij}^k$ of functions satisfies (\ref{eq:ii}), then $ \tilde{c}_{ij}^k:=\frac{1}{2} \left(c_{ij}^k - c_{ji}^k\right) $  satisfies (\ref{eq:iii}).
\end{enumerate}
\end{exo}

\begin{definition}
\label{deff:globallyfinitelygenerated}
    A singular foliation on $M$ as in Exercise \ref{exo:finitelygenerated} is said to be  \emph{finitely generated}.
\end{definition}

These singular foliations will be studied in detail in Section \ref{sec:globally}

\subsubsection{Restriction to an open subset}

Let $ \mathcal F$ be a singular foliation on a manifold $M$, with $ \mathcal F$ as in Definition \ref{def:consensus}.
Let $ \mathcal U \subset M$ be an open subset.
Four natural spaces can be associated to it:
\begin{enumerate}
\item[$\clubsuit$] The subspace of $\mathfrak X_c(\mathcal U) $ of vector fields obtained by restriction to $ \mathcal U$ of vector fields in $ \mathcal F$ whose support lies within $ \mathcal U$. 
\item[$\diamondsuit$] The subspace of $\mathfrak X(\mathcal U) $ of vector fields obtained by restriction to $ \mathcal U$ of vector fields in $ \mathcal F$. 
\item[$\heartsuit$] The $ C^\infty(\mathcal U)$-module generated by elements in $\diamondsuit$.
\item[$\spadesuit$] The subspace of all vector fields $ \mathfrak X(\mathcal U) $  that coincide with a vector field in $ \mathcal F$ in a neighborhood of every point $ m \in \mathcal U$.
\end{enumerate}

\begin{remark}
    $ \spadesuit$ can alternatively be defined as the subspace of all vector fields $X \in \mathfrak X (\mathcal U)$ such that   for every compactly supported function $f$ with support in $ \mathcal U$, the vector field $f X$, extended by zero outside $ \mathcal U$, in an element in $ \mathcal F $.
\end{remark}

\begin{lemma}\label{lem:inclusions}
There is a chain of inclusions:
 $$  \clubsuit \, \subset \, \diamondsuit \,\subset\, \heartsuit  \subset \, \spadesuit.  $$
 \end{lemma}
\begin{proof}
The only non-trivial inclusion is $ \heartsuit  \subset \, \spadesuit$.
A vector field $Y \in \mathfrak X(\mathcal U)$ lies in $ \heartsuit$ if it reads 
$$ Y=\sum_{i=1}^n f_i X_i$$
with $X_i \in \mathcal F $ and $ f_i \in \mathcal C^\infty(\mathcal U)$. Choose $m \in \mathcal U$. For each $i=1, \dots,n$, there exists a function $\tilde f_i \in \mathcal C^\infty(M)$ that coincides with  $f_i$ on some neighborhood $ \mathcal V_i$ of $m$ in $ \mathcal U$, but is now defined on the whole manifold $M$.
The vector field  $  \sum_{i=1}^n \tilde{f}_i X_i$ coincides with $Y$ on the neighborhood $ \cap_{i=1}^n \mathcal V_i$ of $m$. By construction, it belongs to $ \mathcal F$. The vector field $Y$ therefore belongs to $ \spadesuit$.

\end{proof}

 \begin{exo}
     Show that for $ M = \mathbb R$, $\mathcal F = \mathfrak X_c(M) $, and $\mathcal U= \mathbb R\backslash \{0\} $,  each one of the inclusions  in Lemma \ref{lem:inclusions}  is a strict inclusion.
     
     {\emph{Hint:}} Show that the support of every vector field in $\heartsuit $ is contained in $ [-R,R]\backslash \{0\} $ for some $R > 0$.
 \end{exo}

\begin{exo}
Show that $ \clubsuit$ is a singular foliation on $\mathcal U$.
\end{exo}

\begin{exo}
Show that a $ \mathcal C^\infty(M)$-module $ \mathcal F \subset \mathfrak X_c(M)$ stable under Lie bracket is a singular foliation (i.e., is finitely generated in the sense of Definition \ref{def:consensus}) if and only if every point admits a neighborhood in which $ \clubsuit$  is a finitely generated singular foliation (see Definition \ref{deff:globallyfinitelygenerated}).
\end{exo}

\begin{exo}
\label{exo:forfinitelygeneratedrestriction}
Let $ \mathcal F$  be a globally finitely generated singular foliation with generators $ X_1, \dots, X_r$ as in Definition \ref{deff:globallyfinitelygenerated}. Show that for every $ \mathcal U \subset M$:
\begin{enumerate}
\item[$ \clubsuit$] is the space of vector fields of the form $ \sum_{i=1}^r f_i X_i$ with $ f_i$ smooth functions on $ \mathcal U$ with compact support in $ \mathcal U$. 
\item[$ \diamondsuit$] is the space of vector fields of the form $ \sum_{i=1}^r f_i X_i$ with $ f_i$ smooth functions on $ \mathcal U$ obtained by restriction of a compactly supported smooth function on $M $.
\item[$ \heartsuit$] is the space of vector fields of the form $ \sum_{i=1}^r f_i X_i$ with $ f_i$ smooth functions on $ \mathcal U$ whose support is included into a compact subset of $M$. 
\item[$ \spadesuit$] is the space of vector fields of the form $ \sum_{i=1}^r f_i X_i$ with $ f_i$ smooth functions on $ \mathcal U$.
\end{enumerate}
{\emph{Warning:}} None of the four items are trivial, except the second one! For the last one, the proof goes as follows. Let $X$ be in $ \spadesuit$. For every $m \in \mathcal U$,
there exists a $r $-tuple $ f_1^r, \dots,f_r^m$ in compactly supported functions on $ M$ such that $X$ and $ \sum_{j=1}^r f_j^r X_j $ coincide on a neighborhood $ \mathcal U_m$ of $m$ in $ \mathcal U$.
One extracts out of the open cover $ (\mathcal U_m)_{m \in M}$ a sub-open-cover $ (\mathcal U_{m_i})_{i \in I}$ such that a partition of unity $ (\chi_i)_{i \in I}$ exists, and such that each $ m\in \mathcal U$ has a neighborhood that belongs to finitely many of the open subsets $ (\mathcal U_{m_i})_{i \in I}$. The vector fields $ \sum_{j=1}^r  (\sum_{i \in I} \chi_i f^i_j) X_j$    and $X $ coincide on the whole open subset $\mathcal U$. 
Since the function  $\sum_{i \in I} \chi_i f^i_j$ is a smooth function on $ \mathcal U$ for $j=1, \dots, r$, this completes the proof.
\end{exo}

\begin{exo}
\label{exo:fromCompactToSheaves}
Is it true that a vector field in $\spadesuit$ which is compactly supported belongs to $ \clubsuit$?

 \end{exo}

It is therefore natural to ask which one of the four spaces $ \clubsuit,\diamondsuit,\heartsuit,\spadesuit$ is the “good” restriction of $ \mathcal F$ to an open subspace  $\mathcal U$.
While the answer depends on the context, the one that we will need soon is $ \spadesuit$. Hence, we will change the notation:

\begin{definition}
Let $ \mathcal F$ be a singular foliation on $M$.
For every $\mathcal U \subset M $, we will denote by $\mathcal F(\mathcal U )$ the subspace of vector fields on $ \mathcal U $ that coincide with a vector field of $ \mathcal F$ is a neighborhood of every point, i.e., the space denoted by $ \spadesuit$ in the above lines.     
\end{definition}

\begin{exo}
\label{exo:isasheaf}
Let $ \mathcal F$ be a singular foliation on $M$, defined as in Definition \ref{def:consensus}.
Show that by associating to an open subset $ \mathcal U \subset M$ the space $ \mathcal F(\mathcal U)  $, one defines 
 a sheaf on $M$.
 Show that it is both a sheaf of Lie algebras, and a module over the sheaf of smooth functions on $M$. 
\end{exo}

We conclude this discussion by a technical but important lemma.
 
\begin{lemma} \label{lem:numbergenerators}
Let $\mathcal F$ be a singular foliation.
For every point $m \in M$, there exist an open neighborhood $\mathcal U $ of $m$ in $M$ and $ X_1, \dots, X_r \in \mathcal F$ such that for any $ \mathcal V \subset \mathcal U $, $\mathcal F({\mathcal V}) $ is generated over $\mathcal C^\infty(\mathcal V) $ by the restriction to $\mathcal V $ of  $ X_1, \dots, X_r$.
\end{lemma}
\begin{proof}
The proof is similar to the proof of item $ \spadesuit$ in Exercise \ref{exo:forfinitelygeneratedrestriction}. 
By definition \ref{def:consensus},
there exists a neighborhood $ \mathcal U$ and $ X_1, \dots, X_r \in \mathcal F$ such that for
every $X \in \mathcal F$, there exists smooth functions
$ f_1, \dots,f_r \in \mathcal C^\infty(M)$ satisfying that 
the vector fields $  X $ and $ \sum_{j=1}^r f_j X_j $ coincide on $ \mathcal U$. 

Now, let $Y$ be a vector field in $ \mathcal F_\mathcal V$. 
By definition, this means that for every $m' \in \mathcal V$, there exists a neighborhood $ \mathcal W_{m'}$ of $ m'$ in $ \mathcal V$ and a vector field $ Y^{m'}$ in $ \mathcal F$ that coincide with $ Y$ on  $ \mathcal W_{m'}$. By assumption on $ \mathcal U$, the restriction to $ \mathcal U$ of $ Y^{m'}$ reads $ \sum_{j=1}^r f_j^{m'} X_j  $ for some smooth functions $ f_1^{m'}, \dots, f_r^{m'} \in \mathcal C^\infty(M)$. Since the open subsets $ (\mathcal W_{m'})_{m' \in \mathcal V}$ cover $ \mathcal V$, one can extract out of this open cover a sub-open cover $(\mathcal W_{m_i'})_{i \in I}$ that admits a partition of unity $ (\chi_i)_{i \in I}$ and satisfies that every point in $ \mathcal V$ has a neighborhood that belongs to finitely many of the open subsets $(\mathcal W_{m_i'})_{i \in I}  $.
By construction, one has on the whole open subset $ \mathcal V$:
 $$  Y = \sum_{j=1}^r \left(  \sum_{i \in I} \chi_i f_j^{m_i'} \right) X_j  .$$
Since  for every $j=1, \dots, r$, the function $\sum_{i \in I} \chi_i f_j^{m_i'} $ is smooth on $ \mathcal V$, this proves the claim.
\end{proof}

\subsection{Smooth singular foliations: a definition using sheaves}
\label{sec:sheafdef}
\noindent
In order to define singular foliations, the use of compactly supported global vector fields is conceptually easy, but some readers may prefer to use sheaves. Sheaves will in any case be essential for the holomorphic case and the real analytical case. Let us therefore use sheaves to give a definition equivalent to Definition \ref{def:consensus}.

\vspace{0.2cm}
\noindent
In this section, we denote by $$\mathfrak X_\bullet : \mathcal U \longrightarrow \mathfrak X(\mathcal U) \hbox{ and } \mathcal C^\infty_\bullet : \mathcal U \longrightarrow \mathcal C^\infty(\mathcal U) $$ the sheaves of vector fields and of smooth functions on the manifold $M$, respectively\footnote{We remind the reader not to confuse, for a given sheaf  $ \mathcal F_\bullet$ on $M$  and a given open subset $ \mathcal U \subset M$, the space $ \mathcal F(\mathcal U)$ or sections of $ \mathcal F_\bullet$ over $ \mathcal U$ with the restriction to $ \mathcal U$ of the sheaf $\mathcal F $, which should be denoted $ \mathcal F_{|_\mathcal U}$.

}.

\vspace{.5cm}

\begin{definitions}{Definition of a smooth singular foliation: version 2, with sheaves}{consensus2}
A singular foliation on a smooth manifold $M$ is a subsheaf
$$\mathcal F_\bullet \colon \mathcal U \mapsto \mathcal F_\mathcal U$$ of the sheaf $\mathfrak X_\bullet  $ of vector fields on $M$ such that 
\begin{enumerate}
    \item[$(\alpha)$] $\mathcal F_\bullet$ is involutive\footnote{I.e., $[\mathcal F(\mathcal U),\mathcal F(\mathcal U)] \subset \mathcal F(\mathcal U)$ for all open subsets $ \mathcal U \subset M$.},
    \item[$(\beta)$] is a sub-sheaf of $\mathcal C^\infty_\bullet$-modules\footnote{I.e., $ \mathcal C^\infty({\mathcal U} )\mathcal F(\mathcal U) \subset \mathcal F(\mathcal U)$ for all $\mathcal U \subset M$.}, 
     \item[$(\gamma)$] is locally finitely generated\footnote{
For sheaves of modules over functions, the meaning of ``locally finitely generated'' needs to be made very precise. We mean that every point admits an open neighborhood $\mathcal U $ on which there exists $ X_1, \dots, X_r \in \mathcal F(\mathcal U)$ such  that for every $\mathcal V \subset \mathcal U $, the restrictions of  $ X_1, \dots, X_r $ to $\mathcal V $ generate $\mathcal{F}(\mathcal V)$ as a $ \mathcal C^\infty(\mathcal V)$-module.}.
\end{enumerate}
\end{definitions}

\vspace{.5cm}
\noindent 
It is of course embarrassing to have two definitions of singular foliations. Fortunately, as pointed out in Alfonso Garmendia's \cite{PhD_Garmendia} and Roy Wangs \cite[Proposition 2.1.9 and Remark 2.1.13]{wangphd}, these two definitions match, as we now see.

\vspace{.5cm}

\begin{propositions}{No difference!}{prop:nodifference}
Let $M$ be a smooth manifold. 
There is a one to one correspondence between:
\begin{enumerate}
    \item[(i)] Singular foliations defined as in Definition \ref{def:consensus}.
     \item[(ii)] Singular foliations defined as in Definition \ref{def:consensus2}.
\end{enumerate}

\end{propositions}
\begin{proof}
The map $(i) \mapsto (ii)$ consists, given $\mathcal F$ as in Definition \ref{def:consensus}, in considering the sheaf $ \mathcal F_\bullet$ of vector fields that coincide locally with an element in $ \mathcal F$, see Exercise \ref{exo:isasheaf}. 
Lemma \ref{lem:numbergenerators} guaranties it is locally finitely generated in the sense of the third item in Definition \ref{def:consensus2}.

The map $(ii) \mapsto (i)$ consists in considering global compactly supported sections $ \mathcal F_c$ of the sheaf $ \mathcal F_\bullet$ in Definition \ref{def:consensus2}. In equation
 $$ \mathcal F_c = \mathcal F (M) \cap \mathfrak X_c(M).  $$
By construction, $\mathcal F_c$ is a $\mathcal C^\infty(M)$-module stable under Lie bracket. Let us check that it is locally finitely generated in the sense of Definition \ref{def:consensus}.
Let $m \in M$ be a point, $ \mathcal U$ a neighborhood of $m$, and $ X_1, \dots, X_r \in \mathcal F_\mathcal U$ be as in item 3 of Definition \ref{def:consensus2}. Let $ \mathcal U' \subset \mathcal U$ be a second neighborhood of $m$ such that there exists $ \chi \in \mathcal C^\infty(M)$ which is equal to $ 1$ on $ \mathcal U'$ and whose support is a compact subset of $\mathcal U $. 
 The vector fields $ \chi X_1, \dots, \chi X_r$ extend to global sections of $ \mathcal F$ that we still denote by $\chi X_1, \dots, \chi X_r$. These sections belongs to $ \mathcal F_c$ since $ \chi$ is compactly supported.
  By definition of "locally finitely generated" in the third item in Definition \ref{def:consensus2}, for every $Y \in \mathcal F_c$, there exists  $f_1, \dots, f_r \in \mathcal C^\infty(\mathcal U) $ such that the restriction of
 $Y $ to $ \mathcal U $ and $ \sum_{i=1}^r f_i X_i $ coincide on $ \mathcal U$. 
 This implies that the restriction to $ \mathcal U'$ of $Y$  and of $  \sum_{i=1}^r \chi f_i \, \chi X_i $ coincide on $ \mathcal U'$. Since $ \chi f_i$ is now a smooth function on $M$, this implies that $ \mathcal U'$ and the vector fields $ \chi X_1, \dots, \chi X_r$ satisfy the condition of "locally finite generated" as defined in the third item of Definition \ref{def:consensus}.

 We have check that both maps above are inverse one to the other.
Let $ \mathcal F$ be a singular foliation as in Definition \ref{def:consensus}. Let $ \mathcal F_\bullet$ be its associated sheaf as in $ (i) \mapsto (ii)$. A global compactly supported section $Y$ of that sheaf 
is a compactly supported vector field $Y$ on $M$ that locally coincides with a vector field in $ \mathcal F$.
We have to check that $Y$ belongs to $ \mathcal F$. This goes as follows: every $m \in M$ has a neighborhood $ \mathcal U_m$ such that there exists $ X^m \in \mathcal F$ that coincides with $X$ on $ \mathcal U_m$. Since the support of $ Y$ is compact, finitely many of the subsets $ \mathcal U_m$ cover it. Let us denote by $  \mathcal U_{m_1}, \dots, \mathcal U_{m_N} $ such a finite family. 
There exists smooth functions $ \chi_1, \dots, \chi_N$, with support in $\mathcal U_{m_1}, \dots, \mathcal U_{m_N}$ respectively, such that $ \sum_{i=1}^r \chi_i =1$ on a neighborhood $\mathcal V $ of the support of $Y$. 
The vector fields $Y $ and $ \sum_{i=1}^N \chi_i X^{m_i} $ coincide on $ \mathcal V$. The vector fields $Y$ and $ \sum_{i=1}^N \chi  \chi_i X^{m_i} $, with $ \chi$ any function which is $1$ on the support of $Y$ and $0$ outside $ \mathcal V$, coincide therefore on the whole manifold $M$. The crucial point is that the sum is now finite, which proves that $ Y$ belongs to $ \mathcal F$.  The composition \emph{(i)} $ \mapsto $  \emph{(ii)} $ \mapsto $  \emph{(i)} is therefore the identity map.

Now, let $ \mathcal F_\bullet$ be a sheaf as in Definition \ref{def:consensus2}. Consider $ \mathcal U \subset $ on open subset and $ Y \in \mathcal F(\mathcal U)$. For every $m \in \mathcal U$, let $ \chi$ be a smooth function on $M$ equal to $ 1$ in a neighborhood of $m$ and whose support is a compact subset of $ \mathcal U$. Then $ \chi Y$ extends to a compactly supported section of $ \mathcal F_\bullet$ on $M$, and coincides with $ Y$ in a neighborhood of $ M$. 
 The composition \emph{(ii)} $ \mapsto $  \emph{(i)} $ \mapsto $  \emph{(ii)} is therefore the identity map. 
This completes the proof.
\end{proof}

\begin{remark}
Most of the considerations above, and the proof of Proposition \ref{thm:prop:nodifference} in particular, are general phenomenons for what are called \emph{fine sheaves}.
\end{remark}

\begin{bclogo}[arrondi = 0.1, logo = \bcdz]{Warning about notations !}

From now on:
\begin{enumerate}
    \item We will call \emph{foliated manifolds} pairs made of a manifold equipped with a singular foliation, and denote them by $ (M,\mathcal F)$.
    \item We will not make any more notation distinction, in the smooth case, between $\mathcal F $ and $\mathcal F_\bullet$ (i.e., between singular foliations seen as sub-modules of compactly supported vector fields or seen as sheaves, i.e., between singular foliations defined as in Definition \ref{def:consensus} or as in Definition \ref{def:consensus2}). We will mostly use the notation $ \mathcal F$.
\end{enumerate}
\end{bclogo}

\subsection{Singular foliations on complex or real analytic manifolds}

 \vspace{0.2cm}
 \noindent
For complex manifolds, singular foliations have to be defined through sheaves. As a matter of fact, most geometric objects have to be defined through sheaves since there are no or few globally defined functions, vector fields - and so on. For the reader unfamiliar with sheaves, or hostile to them, it simply means that they have to be defined locally. In this section, we fix $M$ a complex or real analytic manifold, and we denote by $\mathcal O$ its sheaf of holomorphic or real analytic functions and by $ \mathfrak X$ the sheaf of vector fields on $M$. For any open subset $\mathcal U \subset M $, we denote by $\mathcal O(\mathcal U) $ the $\mathbb C$-algebra of holomorphic or real analytic $\mathbb K$-valued functions on $\mathcal U $, with $ \mathbb K=\mathbb C$ or $ \mathbb K = \mathbb R$ depending on the context.

Let us start with  a question about the limit of what we want to study.
 \vspace{0.2cm}
 \noindent

\begin{question}
Let $M= \mathbb C$.
Consider the sheaf $ \mathcal F^{O_\infty}$of holomorphic vector fields vanishing together with all their derivatives at $0$. Do we want singular foliations to be defined such that this sheaf is one of them? Or do we do not accept such sheaves?
\end{question}

To answer this question, let us describe $ \mathcal F^{O_\infty}$.
For $ \mathcal U \subset \mathbb C$ an open subset that does not contain $ 0$, $ \mathcal F^{O_\infty}(\mathcal U) = \mathfrak X(\mathcal U)$. For $ \mathcal U$ a connected open subset containing $0$, we have
    $ \mathcal F^{O_\infty}(\mathcal U) =0$. 
    For a generic open subset $ \mathcal U$,
    $ \mathcal F^{O_\infty}(\mathcal U) $ is the space of vector fields on $\mathcal U $ which are equal to zero on the connected component $\mathcal U_0$  of $0$.
    
It is therefore tempting to answer “yes” to the question above because $ (\alpha)$ for every $ \mathcal U$, $ \mathfrak X(\mathcal U)$ is stable under Lie bracket,  $ (\beta)$  $ \mathfrak X(\mathcal U)$ is a module over $ \mathcal O(\mathcal U)$, and $ (\gamma)$ it is finitely generated: the generator is the holomorphic vector field 
  $$ f(z) \frac{\partial}{\partial z}  $$
  with $f$ a function which is $ 0$ on $\mathcal U_0$ and $1$ in $\mathcal U \backslash \mathcal U_0 $.
Also, there is a natural candidate for leaves: there are two of them: the point $\{0\}$ and $ \mathbb C \backslash \{0\}$. Two points are in the same leaf, in this sense, if and only if there exists a vector field in $\mathcal F^{O_\infty} $
 whose time $1$ flow maps the first point to the second one.
  
So, is the answer “yes”? Not so fast! This sheaf is \emph{not}
locally finitely generated in the sense of Definition \ref{def:consensus2}. For every $ \mathcal U$ an open subset containing $0$, no matter how “small”, there is an open subset $ \mathcal V$ such that $ \mathcal F^{0_\infty}(\mathcal V)$ is not zero: it suffices to choose $ \mathcal V \subset \mathcal U_0$ an open subset that does not contain $0$. 
    The space $ \mathcal F^{0_\infty}(\mathcal V)$ being not reduced to zero, it can not be obtained by restricting to $ \mathcal V$ generators of $ \mathcal F^{0_\infty}(\mathcal U)$, since those are all zero.

 \vspace{0.2cm}
 \noindent
So, is the answer “no”?  Well, it depends on the context, but we suggest answering “no” to that question, and to define singular foliations as follows:

\vspace{.5cm}

\begin{definitions}{Holomorphic or real analytic singular foliations}{consensus2alg}
A singular foliation on a complex (or real analytic) manifold $M$ is a subsheaf
$\mathcal F $ of the sheaf $\mathfrak X $ of holomorphic (or real analytic) vector fields on $M$ which 
\begin{enumerate}
    \item[$(\alpha)$] is involutive,
    \item[$(\beta)$] is stable under multiplication under $\mathcal O$,
    \item[$(\gamma)$] is locally finitely generated\footnote{In the literature, this property is often referred to as being “of finite type”. For consistency of notations throughout the lecture, we prefer to say “locally finitely generated”}, i.e., for every point $m$, there exists a neighborhood $ \mathcal U$ and vector fields $ X_1, \dots, X_r \in \mathcal F(\mathcal U)$ such that for every open subset $\mathcal V \subset \mathcal U $, $ \mathcal F(\mathcal V)$ is generated as a $ \mathcal O_\mathcal V$-module by the restrictions to $ \mathcal V$ of $ X_1, \dots, X_r$.
\end{enumerate}
\end{definitions}
\vspace{.5cm}

Now, complex geometers may be surprised by the previous definition, because the “locally finitely generated” condition is not commonly used in complex geometry, even under its alternative name “of finite type”. They would prefer to use \emph{coherent sheaves} \cite{zbMATH04118703}. However, it is a classical, but hard, theorem in complex analysis, that locally finitely generated (=of finite type) subsheaves of a coherent sheaf are coherent sheaves and that vector fields form a coherent sheaf\footnote{Look for "Oka's coherence theorem".}. As a consequence, Definition \ref{def:consensus2alg} can be equivalently restated as  in the next proposition.

\begin{proposition} 
\label{prop:coherent}
  A singular foliation on a complex (or real analytic) manifold $M$ is a subsheaf
$\mathcal F_\bullet $ of the sheaf $\mathfrak X $ of vector fields on $M$ which 
\begin{enumerate}
    \item[$(\alpha)$] is involutive,
    \item[$(\beta)$] is stable under multiplication under $\mathcal O$,
    \item[$(\gamma)$] is coherent.
    \end{enumerate}
\end{proposition}

\vspace{0.5cm}
Proposition \ref{prop:coherent} implies that our definition of a singular foliation in the complex case coincides with the one given in  Paul Baum and Raoul Bott's \cite{zbMATH03423310}, Ali Sinan Sert\"oz's \cite{zbMATH04118703},  Andr{\'e} Belotto da Silva and Daniel Panazzolo's \cite{zbMATH07124409}, or Tatsuo Suwa's \cite{zbMATH03931387}.

 \vspace{0.5cm}
 \noindent

We just explained why, in the complex setting, one can re-use mutatis mutandis Definition \ref{def:consensus2}. 
However, any sheaf that satisfies items ($\alpha$) and ($\beta$) in Definition \ref{def:consensus2} in fact satisfies a variation of item ($\gamma$) as well. Let us be precise: unlike the algebra of germs of smooth functions,  the algebra of germs of holomorphic or real analytic functions are Noetherian:

\begin{theorem} \label{thm:Tougeron} \cite{MR240826}
Germs of holomorphic (resp. real analytic) functions near $0 \in \mathbb C^n$ (resp. $\mathbb R^n$) form a Noetherian ring.
\end{theorem}

\begin{remark}
 We warn the reader that Theorem \ref{thm:Tougeron} does \emph{not} imply that  any sub-$ \mathcal O(\mathcal U)$-module $\mathcal F({\mathcal U}) \subset \mathfrak X(\mathcal U)$ is finitely generated for every $ \mathcal U \subset M$.
 Here is a classical counter example: take $M=\mathcal U=\mathbb C$, and consider $ \mathcal F$ to be the sheaf of all vector fields vanishing at all the points $ \{ n \mid n \in \mathbb N\} $, except a finite number of them. The space of global sections, i.e., $ \mathcal F(\mathbb C)$, is not finitely generated.
 \end{remark}

\vspace{0.2cm}

\noindent
Let us explain the consequences of Theorem \ref{thm:Tougeron} in the holomorphic setting: the real analytic setting is similar.
In a chart neighborhood $ \mathcal U$ of a point $m \in M$, with coordinates $z_1, \dots, z_d $, holomorphic vector fields decompose as sums 
 $$ \sum_{i=1}^d f_i( z_1, \dots, z_d)  \frac{\partial}{\partial z_i}$$
 with $ f_1, \dots, f_d \in \mathcal O(\mathcal U)$ being $\mathbb C $-valued holomorphic functions on $\mathcal U $.
In particular, as a module over $ \mathcal O(\mathcal U)$, holomorphic vector fields of $ \mathcal U$ are isomorphic with
$$ \mathfrak X(\mathcal U) \simeq \underbrace{\mathcal O({\mathcal U}) \oplus \cdots \oplus \mathcal O({\mathcal U})}_{d \hbox{ terms }}$$ (with $d$ the dimension of the manifold). 
The isomorphism above makes the germs at $m$ of all elements $ X \in  \mathcal F$  a sub-$\mathcal O_{m}$-module of 
$$ \underbrace{\mathcal O_{m} \oplus \cdots \oplus \mathcal O_{m}}_{d \hbox{ terms }} $$
with $\mathcal O_m $ being the algebra of germs at $m$ of holomorphic functions defined near $m$.
The henceforth obtained sub-module is finitely generated over $\mathcal O_m $
 by Theorem \ref{thm:Tougeron}.
Let $ X_1, \dots,X_r$ be vector fields on $ \mathcal U_1, \dots, \mathcal U_r$ respectively whose germs at $m$ generate this module.
Let $\mathcal U = \cap_{i=1}^r  \mathcal U_i $.  For any connected neighborhood $\mathcal V \subset \mathcal U$ of $m$ in $M$, and any vector field $ X \in \mathcal F(\mathcal V )$, the germ of $X$ at $m$ is in the $ \mathcal O_m$-module generated by $ X_1, \dots, X_r$. This means that there exists  $r$ holomorphic functions $f_1,\dots, f_r$ defined on some neighborhood $ \mathcal W \subset \mathcal V$ such that $X =\sum_{i=1}^r f_i X_i$. But these functions are defined on a smaller subset  $\mathcal W$ of $ \mathcal V$. 
We cannot therefore conclude that $ \mathcal F$ is a locally finitely generated module over the sheaf of holomorphic functions. But the previous results are still a phenomenon that does not appear in the smooth setting. In several contexts (as in \cite{LLS}), it might be enough to work with sheaves  $ \mathcal F \subset \mathfrak X$ such that only conditions $(\alpha) $ and $ (\beta)$ hold in Definition \ref{def:consensus2alg}, leaving aside condition $ (\gamma)$. It is always granted for free that any point $m\in M$ admits a neighborhood $ \mathcal U$ and $ X_1, \dots, X_r \in \mathcal F(\mathcal U)$ such that for any $ \mathcal V \subset \mathcal U$ containing $m$, and any $ Y \in \mathcal F(\mathcal V)$, there exists an open subset $ \mathcal W \subset \mathcal V$ containing $m$ and functions $ f_1, \dots, f_r \in \mathcal O(\mathcal W)$ such that $ Y= \sum_{i=1}^r f_i X_i$ on $ \mathcal W$. 

\subsubsection{Restriction to an open subset}

When smooth singular foliations are defined as in Definition \ref{def:consensus}, then restrictions to open subsets are easy to define, and subtleties like Lemma \ref{lem:inclusions} do not appear.

Let $ \mathcal F$ be a singular foliation on a smooth, complex or real analytic manifold, seen as a sheaf as in Definition \ref{def:consensus2} or \ref{def:consensus2alg}. For every open $ \mathcal U \subset M$, the sheaf $ \mathcal F$ can be restricted to $ \mathcal U$. The restriction defines a singular foliation on $ \mathcal U$ that we denote by $ \mathcal F_\mathcal U$. 

In the smooth case, under the correspondence introduced by Proposition \ref{thm:prop:nodifference}, this singular foliation $\mathcal F_\mathcal U $ corresponds to the subspace of $ \mathfrak X_c(\mathcal U)$ denoted by $ \clubsuit$ in Section \ref{sec:defsmooth}.

\vspace{0.5cm}
We warn the reader not to confuse $ \mathcal F_\mathcal U$, which is a sheaf on $ \mathcal U$, and $ \mathcal F(\mathcal U)$ which stands for its sections on $ \mathcal U$, especially in the complex or real analytic settings.

\vspace{0.5cm}
Now, in the complex case, there is a phenomenon called Hartog's theorem or Riemann's extension theorem that says that a holomorphic object defined outside a codimension $2$ analytic subset extends to the whole space. It is tempting to apply such results to the complement of the singular locus of a singular foliation (see Section \ref{sec:openReg}): this is discussed in Remark \ref{rmk:Hartog}.

\subsection{Singular foliation on an affine complex variety}
\label{sec:affinevarieties}

\vspace{0.2cm}
Let us now define singular foliations on affine varieties over the field $\mathbb C $.
Let us recall that an affine variety (maybe non-irreducible) is a subset  $W \subset \mathbb C^n $ given by polynomial equations. But we will see them in a more algebraic manner.

Denote by $\mathcal O_n $ the algebra of polynomial functions in $ n$-variables, and $\mathcal I_W \subset  \mathcal O_n$ the ideal of functions vanishing on $W$. 
We call \emph{functions on $W$} the quotient ring $ \mathcal O_W:= \frac{\mathcal O_n}{\mathcal I_W}$.
We call \emph{vector fields on $W$} and denote by $\mathfrak X_W $ the $\mathcal O_W $-module of derivations of $\mathcal O_W $. It is equipped with the commutator as a Lie bracket. Since the algebra $\mathcal O_W $ is a Noetherian algebra,  and since $\mathfrak X_W $ is a  $\mathcal O_W $-module of finite rank, any sub-$\mathcal O_W $-module is finitely generated.  The assumption ``locally finitely generated'' is therefore useless in that context, and we suggest the following definition of a singular foliation on an affine variety.

\vspace{.5cm}

\begin{definition}

A singular foliation on an affine variety $W$ is a sub-space
$\mathcal F $ of the $\mathcal O_W$-module of $\mathfrak X_W $ which 
\begin{enumerate}
    \item[$(\alpha)$] is involutive,
    \item[$(\beta)$] is stable under multiplication under $\mathcal O_W $.
\end{enumerate}
\end{definition}

\vspace{0.5cm}
\noindent
Notice that the definition does not make reference to the ``ambient space'': the definition above only makes use of the algebra of functions on $W$. If two affine varieties $W$ and $W'$ in $\mathbb C^n $ and $ \mathbb C^{n'}$ respectively satisfy $\mathcal O_n/\mathcal I_W \simeq \mathcal O_{n'}/\mathcal I_{W'} $, then they are equipped with the same singular foliations.

\vspace{0.2cm}
\noindent
For schemes, or quasi-projective varieties, again, the use of sheaves will be necessary.

\begin{exo}
Write the definition of a singular foliation on an arbitrary scheme.
\end{exo}

\vspace{0.2cm}
\noindent
The above discussion leads to 
 a purely algebraic definition of what a singular foliation is.
Let $\mathcal O $ be a commutative unital algebra
(which may be thought of as being an algebra of ``functions'' – whatever it means). 

\begin{definition}\label{def:algebraic}

A sub-space
$\mathcal  F$
 of ${\mathrm{Der}}(\mathcal O) $ is said to be an \textit{algebraic singular foliation} if:
 \begin{enumerate}
     \item[$ \alpha$] $\mathcal F $ is a stable under the Lie bracket of ${\mathrm{Der}}(\mathcal O) $,
     \item[$ (\beta)$] $\mathcal F $ is a sub-${\mathrm{Der}}(\mathcal O) $-module of  ${\mathrm{Der}}(\mathcal O) $,
     \item[$ (\gamma)$] and it is finitely generated as an $\mathcal O $-module.
 \end{enumerate}

\end{definition}

\vspace{0.2cm}
\noindent
If the algebra $\mathcal O $ is \textit{(i)} finitely generated and \textit{(ii)} Noetherian, then every involutive $\mathcal O$-sub-module of derivations of $\mathcal O$ is an algebraic singular foliation: the condition $ (\gamma)$ is therefore useless in that case.
In particular:

\begin{proposition}
   \label{prop:fromalgtoholo} 
   An algebraic singular foliation on a  smooth\footnote{Smooth here means “that has no singular point” !} affine variety in $W \subset  \mathbb C^N $ is also a holomorphic singular foliation on\footnote{$W$ being now seen as a complex manifold.} $W$.
\end{proposition}

\subsubsection[]{Germification}\label{sec:germification}

Here is an example of an algebraic singular foliation associated to any smooth, real analytic or complex singular foliation $ \mathcal F$ and a manifold $M$ (with sheaf of functions $ \mathcal O$).
Choose a point $m \in M$.

	\begin{definition}
	A \emph{function germ at a point $m\in M$} is an equivalence class of pairs $(\mathcal U,f)$ with $\mathcal U\subset M$ an open subset containing $m$, and $f \in \mathcal O(\mathcal U)$, under the equivalence relation : $(\mathcal U,f)\sim(\mathcal V,g)$ if $f=g$ on an open subset of $\mathcal U\cap \mathcal V$ containing $m$.
	 The quotient comes equipped with a natural algebra structure denoted by $\mathcal O_{m}$. \end{definition}
  
  Notice that $\mathcal O_{m}$ is a local ring. 
  Any vector field $ X \in \mathfrak X(M)$, defined in a neighborhood $\mathcal F $ of $ m$, induces a derivation $ X_m$ of $ \mathcal O_m$. For $ \mathcal F$ a singular foliation on $ M$, denote by $ \mathcal F_m$ the space of derivations of $ \mathcal O_m$ obtained by considering all such derivations for all vector fields in $ \mathcal F(\mathcal U)$  for all
 open neighborhoods  $\mathcal U $ of $m$.
 Axioms $ (\alpha)$, $(\beta)$, and $(\gamma)$ in Definition \ref{def:consensus2} or \ref{def:consensus2alg}
 imply that axioms  $ (\alpha)$, $(\beta)$, and $(\gamma)$ in Definition \ref{def:algebraic} are satisfied.
   The following definition then makes sense.

\begin{definition}
Let $\mathcal{F}$ be an algebraic singular foliation. For a given point $m \in M$, we call \emph{germ of $\mathcal F$ at $m$} the algebraic singular foliation $\mathcal F_m $ over $ \mathcal O_m$.
\end{definition}

We could also work with formal functions (=formal power series) instead of germs, see \cite{Ryvkin2,fischer2024classification} and the discussion following Theorem \ref{thm:LinearCerveau}.

\subsection{Globally finitely generated singular foliations}
\label{sec:globally}

We have  been through quite an extensive discussion about the limits and sense of the ``locally finitely generated'' condition. But quite a few singular foliations are in fact \emph{globally} finitely generated.

Notice that the notion was already introduced in Definition \ref{deff:globallyfinitelygenerated} in the smooth setting using compactly supported vector fields, but we can now define them within the setting of sheaves.

Globally finitely generated can be defined in the same manner in the smooth, complex, or real analytic contexts altogether. We leave  it to the reader to check that the next definition matches Definition\footnote{Notice that the vector fields $ X_1, \dots, X_r$ are global sections of the sheaf $ \mathcal F_\bullet$ if one uses Definition \ref{def:globallyfinitelygenerated} but may not be in $ \mathcal F_c \subset \mathfrak X_c(M)$ if one uses Definition \ref{deff:globallyfinitelygenerated}. This makes no practical difference.} \ref{deff:globallyfinitelygenerated} in the smooth case.

\vspace{.5cm}

\begin{definitions}{A common definition}{globallyfinitelygenerated}
A singular foliation  $\mathcal F $ on a manifold $M$  is said to be \emph{finitely generated} if there exists vector fields\footnote{Called \emph{generators of $ \mathcal F$}.} $X_1, \dots, X_ r$  such that for every open subset $\mathcal U \subset M$, 
 $ \mathcal F({\mathcal U})$ is the $ \mathcal O(\mathcal U)$ module generated by the restrictions to $\mathcal U $ of $X_1, \dots, X_r$.
\end{definitions}

\vspace{.5cm}

\begin{remark}
Let $\mathcal F $ be a singular foliation on $M$. Every point $m \in M$ has a neighborhood on which it is finitely generated.
\end{remark}

Here is an important result.

\begin{lemma}
\label{lem:christoff}
Let $ X_1, \dots,X_r$ be generators of  a finitely generated singular foliation  $\mathcal F $ on a manifold $M$, then there exists a family of functions $ c_{ij}^k \in \mathcal O$, indexed by $ i,j,k =1, \dots, r$, such that for all $i,j=1, \dots, r$:
 $$[X_i,X_j] = \sum_{k=1}^r c_{ij}^k X_k. $$
 \end{lemma}

 Such functions $ c_{ij}^k$ are called \emph{Christoffel symbols of $ \mathcal F$ with respect to $X_1, \cdots,X_r $}.
 Since there are, in general, relations between the generators $X_1, \dots, X_k $, the Christoffel symbols $ c_{ij}^k$ are \underline{not unique}.
 
 \begin{exo}
 \label{exo:christo1}
Let $ X_1, \dots, X_r$ be generators of a finitely generated singular foliation $\mathcal F $, and  $(c_{ij}^k)_{i,j,k=1}^r$ a choice of Christoffel symbols of $\mathcal F $ with respect to $X_1, \dots, X_r$.
\begin{enumerate}
\item 
Show that 
$$ \left(\frac{c_{ij}^k-c_{ji}^k}{2} \right)_{i,j,k=1}^r $$
is again a choice of Christoffel symbols of $\mathcal F $ with respect to $X_1, \dots, X_r$.
\item Show that, without any loss of generality, Christoffel symbols of $\mathcal F $ with respect to $X_1, \dots, X_r$ can be assumed to satisfy $c_{ji}^k = - c_{ij}^k $ for all possible indices.
\end{enumerate}
 \end{exo}

\begin{exo}
\label{exo:infinitestability}
The \emph{``non-finitely-many-generators''} singular foliation -
an example due to Iakovos Androulidakis and Marco Zambon.
On $M= \mathbb R^2 $, call $\mathcal F $ the space of all vector fields $X \in \mathfrak X (\mathbb R^2) $ that vanish at order $n$ at the point of coordinates $(n,0) $. 
i.e., vector fields of the form:
 $$ X = f(x,y) \frac{\partial}{\partial x} + g(x,y) \frac{\partial}{\partial y}  $$
 such that for all $a,b,n \in \mathbb N_0$ with $a + b \leq n $:
  $$ \frac{\partial^{a+b} f}{\partial x^a\partial y^b}{(0,n)}= \frac{\partial^{a+b} g}{\partial x^a\partial y^b}{(0,n)}=0 . $$

\begin{center}
\includegraphics[scale=0.2]{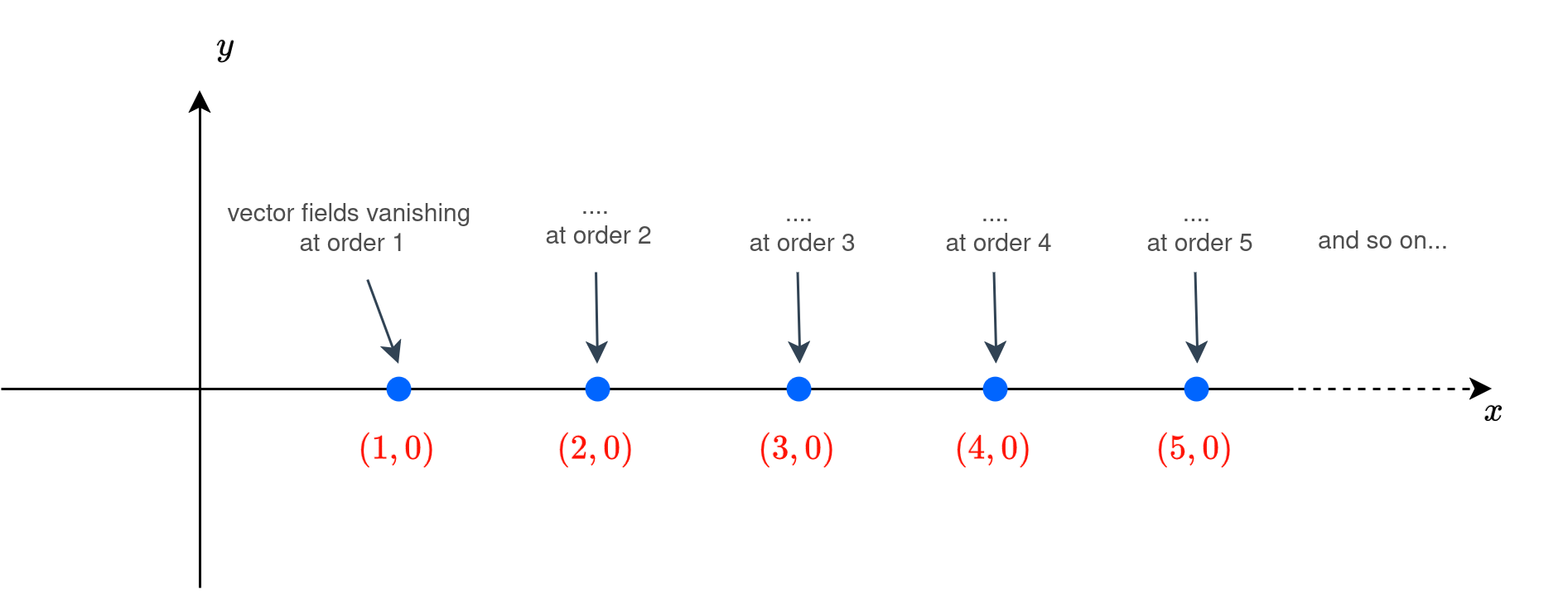}

A representation of the ``non-finitely-many-generators'' singular foliation.
\end{center}
  
Show that
\begin{enumerate}
    \item $\mathcal F $ is an integrable distribution.
    \item $\mathcal F $ is locally finitely generated.
    \item $ \mathcal F$ is not globally finitely generated.
    \item  $ \mathcal F$ is not the image through the anchor map of a Lie algebroid on $\mathbb R^2 $.
\end{enumerate}

\end{exo}

In the smooth setting, if for a singular foliation the number of local generators is bounded, then it is finitely generated.
\vspace{0.5cm}

\begin{proposition}\label{prop:localglobalnumberofgenerators}
Let $M$ be a smooth manifold of dimension $n$ and $\mathcal F$ a singular foliation admitting a global bound $k$ of the number of local generators, i.e., $M$ is covered by open sets $\{U_i\}$ such that $\mathcal F({U_i})$ is generated by $k$ or less vector fields.

Then $\mathcal F$ can be generated by a family of $(n+1)\cdot k$ vector fields. In particular, it is a finitely generated singular foliation\footnote{See Definition \ref{def:globallyfinitelygenerated}}.
\end{proposition}

\begin{proof}
Since the covering dimension of a smooth manifold is equal to its dimension as a manifold, the cover $U$ admits a refinement such that every point is covered by at most $n+1$ points. This is the classical paving principle. Now we can apply Theorem 3 from \cite{ostrand}, stating that there exists a finite open covering $V_1,...,V_{n+1}$ of $M$ such that $V_{j}=\bigsqcup_{\alpha}V_j^\alpha$, where each $V_j^\alpha$ is included in some $U_i$ and $V_j^\alpha$, $V_j^\beta$ have disjoint closures for fixed $j$ and $\alpha\neq\beta$. This latter fact means that there exist functions $\phi_j^\alpha$ such that $\phi_j^\alpha|_{V_{j}^\alpha}=1$ and $\phi_j^\alpha|_{V_{j}^\beta}=0$ for $\alpha\neq \beta$. We pick for each $(j,\alpha)$ a collection of vector fields $Y_j^{\alpha, 1},...,Y_j^{\alpha, k}\in \mathcal F$ whose restrictions generate $\mathcal F_{V_j^\alpha}$. We set 
$$X_j^s= \sum_{\alpha}=\phi_j^\alpha Y_j^{\alpha, s}
$$ 
for $j\in \{1,...,n+1\}$ and $s\in \{1,...,k\}$. By construction, this family of vector fields generates the foliation.
\end{proof}

\begin{remark}One recognizes in the proof of Proposition \ref{prop:localglobalnumberofgenerators} the logic of the argument used to show that any finite-dimensional vector bundle over a smooth manifold is the quotient of (or direct summand in) a trivial vector bundle.
\end{remark}

 \section{Some basic notions about singular foliations (symmetries and regular part)}

For $ \mathcal F$ a singular foliation on $ M$ and $ \mathcal U \subset M$ an open subset, we denote by $ \mathcal F_\mathcal U$ the induced singular foliation on $ \mathcal U$. 

 \subsection{Symmetries and inner symmetries}

\label{sec:symmetry}

In the lines below, we place ourselves in the smooth case.
In the real analytic or complex case, the whole story is similar, but the word ``diffeomorphism'',``smooth'', and ``interval'' are to be replaced by ``biholomorphism'', ``holomorphic'', and ``open ball in $ \mathbb C$'' in the text below.

 Let $ (M, \mathcal F)$ and $(M',\mathcal F') $ be foliated manifolds in the smooth category. We call \emph{isomorphism of singular foliations} 
a diffeomorphism $\phi \colon M \simeq M' $
such that $X \in \mathcal F$ if and only if $ \phi_*(X) \in \mathcal F'$.
When $M'=M $ and $ \mathcal F' = \mathcal F$, we shall speak of a \emph{symmetry of $ \mathcal F$}.  Sometimes, we will speak of a \emph{local symmetry of  a foliated manifold $(M,\mathcal F) $}: it is a triple $ (\mathcal U, \mathcal U', \phi)$ with $\mathcal U , \mathcal U'$ open subsets of $ M$ and $\phi \colon \mathcal U \longrightarrow \mathcal U' $ a diffeomorphism which maps $\mathcal F_{\mathcal U}  $ to  $ \mathcal F'_{\mathcal U'} $. 

\vspace{1mm}
\noindent
Symmetries of $ (M,\mathcal F)$ form a group, and local symmetries form what is called a pseudo-group.

\begin{definition}
\label{symmetries}
We denote by ${\mathrm{Sym}}(\mathcal F) $ the group of symmetries of a singular foliation.
\end{definition}

\vspace{3mm}
\noindent
We now intend to define inner symmetries. This goes through the correct definition of time-dependent vector fields. Again, we work on the smooth category\footnote{In the real-analytic or holomorphic contexts, there is a subtlety  about time-dependent vector fields that appear, e.g., in Definition \ref{innersymmetries}: even in the real analytic or holomorphic cases,  it suffices that this time-dependency be \underline{smooth} in the parameter that we denote by $t$ below.}.

\vspace{1mm}
\noindent
First, let us define time-dependent vector fields valued in a singular foliation.

\begin{definition}
\label{innersymmetries}
Let $I$ be an interval of $\mathbb R $. Let $ (M,\mathcal F)$ be a foliated manifold. 
We say that a family $ (X_t)_{t \in I}$ of vector fields is a \emph{smooth time-dependent vector field\footnote{We sometime also say "vector fields depending smoothly on $t$"} in $\mathcal F $} if, in a neighborhood $\mathcal U$ of every point $m \in M$, there exists local generators $ X_1, \dots, X_r $ of $ \mathcal F$ and smooth functions $ f_1, \dots, f_r \in \mathcal C^\infty(M \times I)$ such that for all $ m \in \mathcal U, t \in I$:
$$ \left. X_t \right._{|_{m}}= \sum_{i=1}^r f_i(m,t) \left.X_i\right._{|_{m}} .$$

\end{definition}

\begin{remark}\label{skandalexample}
\normalfont
According to Definition \ref{innersymmetries}, applied to the case $\mathcal F = \mathfrak X_c (M) $,  a smooth time-dependent vector field is a family $ t \mapsto X_t$ of vector fields that, altogether, form a smooth map from $ M \times I$ to $TM$.  It deserves to be noticed that a smooth time-dependent vector fields valued in $\mathcal F $ may \emph{not} be a smooth time-dependent vector field on $M$ such that $ X_t \in \mathcal F$ for every value of $t$.
Here is a counter example (provided to us by Georges Skandalis). Take the singular foliation on $M=\mathbb R$ generated by the vector field
$$X= e^{-\frac{1}{x^2}} \frac{\partial}{\partial x}$$
Consider the family of vector fields on $M$ given by
$$  X_t = \left\{\begin{array}{cl} \frac{t}{x^2+t^2} X &  \hbox{ 
 for $t \neq 0$} \\ 0 & \hbox{ for $t=0$.} \end{array} \right.  $$
Since $\frac{t}{x^2+t^2}  e^{-\frac{1}{x^2}} $ extends to a smooth function on $ \mathbb R^2$ that vanishes with all its derivatives at the point $ (0,0)$, the previously defined family
is a smooth time-dependent family of vector fields on $M=\mathbb R$. Also, for all $t \in \mathbb R$, $X_t\in \mathcal F$. 
But since the function $\frac{t}{x^2+t^2}$ does not extend to a smooth function at $ (0,0)$, the family $X_t$ is \emph{not} a smooth time-dependent vector field in $\mathcal F$ in the sense of Definition \ref{innersymmetries}. 
\end{remark}

\begin{exo}
Show that in Definition \ref{innersymmetries}, 
``there exists local generators $ X_1, \dots, X_r $ of $ \mathcal F$ and smooth functions $ f_1, \dots, f_r \in \mathcal C^\infty(M \times I)$'' could be replaced by 
``for every local generators $ X_1, \dots, X_r $ of $ \mathcal F$, there exists smooth functions $ f_1, \dots, f_r \in \mathcal C^\infty(M \times I)$''.
\end{exo}

\begin{exo}
\label{exo:linkedwithanchor}
This exercise supposes the notion of anchored bundle, see Section \ref{sec:AnchoredBundle}. 
Let $A \to M$ be a vector bundle. We say that a time-dependent section $ (a_t)_{t \in [0,1]}$ of $A$ is smooth if the map  $(m,t) \mapsto \left. a_t \right._{|_{m}} $ is a smooth map from $M\times [0,1]$ to $A$. 
Show that $ (X_t)_{t \in I}$ is a smooth time-dependent vector field in $\mathcal F $ if and only if for every anchored bundle $(A,\rho) $ over $\mathcal F $, there exists a smooth time-dependent\footnote{A time-dependent section $(a_t)_{t \in I}$ is said to be smooth if  $ (m,t) \mapsto a_t(m)$ is a smooth map from $M\times I $ to $A$.} section $(a_t)_{t \in I} $ of $A \to M$ such that $ \rho(a_t) =X_t$ for every $t \in I$.
\end{exo}

\begin{exo}
 \label{exo:smoothtimedepDirectProduct}
 Let $(M,\mathcal F)$ be a foliated manifold.
 Let $I$ be an open interval of $\mathbb R $.
 
 Let $\hat{\mathcal F} $ be the singular foliation on $M \times I $ which is the direct product of $(M,\mathcal F) $ with $(I,0) $.
 
Let $\tilde{\mathcal F} $ be the singular foliation on $M \times I $ which is the direct product of $(M,\mathcal F) $ with $(I,\mathfrak X(I)) $.

 Consider a family $ (X_t)_{t \in I}$ of vector fields such that $X_t \in \mathcal F$ for all $t \in I$.
 
 \begin{enumerate}
 \item
 Show that the following are equivalent:
\begin{enumerate}
\item[(i)]  $ (X_t)_{t \in I}$ is a smooth time-dependent vector field in $\mathcal F $,
\item[(ii)] the vector field $X_\bullet $ on $M \times I $ whose value at $(m,t) $ is $ (\left. X_t\right._{|_{m}},0)$ belongs to $\hat{\mathcal F} $,
\item[(iii)] the vector field $Y $ on $M \times I $ given by
$$ Y = X_\bullet+ \frac{\partial}{\partial t}$$ 
belongs to $\tilde{\mathcal F} $.
\end{enumerate}
\item Show that the time $ \tau$ flow $ \phi_\tau^Y$ of $Y$ and the time $\tau $-flow $ \phi_\tau^{X_\bullet}$ of $ (X_t)_{t \in I}$ are related by
  $$  \phi_\tau^Y ( m ,t) = \left( \phi_\tau^{X_\bullet}(m), t + \tau\right) $$
  for all $ t,\tau,m$ for which these flows exist.
  \item Conclude that ``any inner symmetry (defined below) of $\mathcal F $ is obtained by restricting to well-chosen transverse submanifolds the flow of a vector field in $\tilde{\mathcal F} $''.
  \end{enumerate}
\end{exo}

\vspace{0.5cm}

\begin{definitions}{Inner symmetries}{InnerSym} 
Let $ (M,\mathcal F)$ be a foliated manifold.
A diffeomorphism of $ (M,\mathcal F)$ is said to be an \emph{inner-symmetry} if it is the time-$1$ flow of a smooth time-dependent vector field $(X_t)_{t \in [0,1]} $ in $\mathcal F $ whose time-$1$ flow exists.
We denote by $ \mathrm{Inner}(\mathcal F)$ the set of inner symmetries of $\mathcal F $.

If the previous diffeomorphism only on some open subset $\mathcal U \subset M $, we speak of a \emph{local inner symmetry}.
\end{definitions}

\vspace{0.5cm}
\noindent
Here is the first important lemma.

\begin{lemma}
\label{lem:innerIs}
Inner symmetries of a singular foliation form a group.
\end{lemma}
\begin{proof}
If a smooth time-dependent vector fields $(X_t)_{t \in [0,1]} $ in $\mathcal F $ that yields the inner-symmetry $ \phi$, then $Y_t= -X_t$ is also a smooth time-dependent vector field and it yields $\phi^{-1} $. Also, upon replacing $ X_t$ by $ f'(t) X_{f(t)}$ with $ f \colon [0,1] \to [0,1] $ any smooth function with $ f=0$ near $0$ and $ f=1$ near $1$, one yields a smooth time-dependent vector field that vanishes near $t=0$ and $t=1$ together with all its derivatives with respect to $t$, whose time-$1$ flow is still $ \phi$. This allows gluing: for any two inner symmetries $\phi,\psi $ given by such smooth time-dependent vector fields $(X_t)_{t \in [0,1]} $ and $(Y_t)_{t \in [0,1]} $ respectively, the vector field 
$$ Z_ t = X_{2t} \hbox{ if $t \in [0, 1/2]$ and }  Z_t = Y_{2t-1} \hbox{ if $t \in [1/2,1]$}.  $$
is a smooth time-dependent vector field in $\mathcal F $ that yields $ \psi \circ \phi$. 
\end{proof}

\vspace{1mm}
\noindent
It is a highly non-trivial result that \textbf{inner symmetries of $(M,\mathcal F) $ are in fact symmetries of $(M,\mathcal F) $}: this will be the point of one of the main theorems in these lectures, namely Theorem \ref{thm:flowofinner}, and more precisely Corollary \ref{coro:flowissymmetry}. Let us admit this theorem and the corollary for the remaining of this section.

\begin{exo}
\label{exo:InnerAndSymmetry}
Let $(M,\mathcal F) $ be a foliated manifold.
\begin{enumerate}
\item 
Show that for any symmetry $ \psi$ 
 and any inner symmetry $ \phi$ of $ \mathcal F$, the composition $\psi \circ \phi \circ \psi^{-1}$ is an inner symmetry.
 \item[] In the process, describe the smooth time-dependent vector field $ Y_t$ whose time-$1$ flow is $ \psi \circ \phi \circ \psi^{-1}$ out of $ \psi$ and out of a smooth time-dependent vector field $ X_t$ whose time-$1$ flow is $ \phi $.
 \end{enumerate}
Note: In view of Corollary \ref{coro:flowissymmetry} (that states that inner symmetries form a normal subgroup of the group of symmetries), this exercise means that inner symmetries form a \emph{normal} subgroups of the group of symmetries.
 \end{exo}

\vspace{1mm}
\noindent
For a symmetry $\phi $ of a singular foliation $ (M,\mathcal F)$, the notion of fixed point is straightforward: it is just a point $ m \in M$ such that $\phi(m)=m $.
There is a notion of ``fixed point'' which is much subtle for inner-symmetry, and that we now introduce.

\begin{definition}
\label{def:veryfixed}
    Let $ (M,\mathcal F)$ be a foliated manifold. 
    A point $m \in M$ is said to be a \emph{very-fixed point} of a (maybe local) inner symmetry $ \phi$ if $ \phi$ is the time $1$-flow of some smooth time-dependent vector field $ (X_t)_{t \in [0,1]}$ in $\mathcal F $ such that
the functions $ f_1, \dots, f_r$ that appear in Definition \ref{innersymmetries} vanish at the point $m $ for all $t \in [0,1]$.

\vspace{1mm}
    \noindent
Equivalently, $m$ is a very-fixed point of $ \phi$ if there exists local generators $ X_1, \dots, X_r $ of $ \mathcal F$ on some neighborhood $ \mathcal U$ of $m$ and smooth functions $ f_1, \dots, f_r \in \mathcal C^\infty(M \times I)$ satisfying $f_i(m,t)=0 $ for all $t \in [0,1] $ such that $ \phi$ is the time-$1$ flow of the vector field given $\forall x \in \mathcal U, t \in [0,1] $ by
$$ \left. X_t \right._{|_{x}}= \sum_{i=1}^r f_i(x,t) \, \left. X_i \right._{|_{x}} .$$ 
\end{definition}

  \begin{exo}
Let $ (M,\mathcal F)$ be a foliated manifold. 
  \begin{enumerate}
\item Show that if $m$ is a very fixed point for an inner symmetry $ \phi$, then it is also a fixed point, i.e., $\phi(m)=m $. 
\item We prove that the converse is not true. Let $ \mathcal F$ be the singular foliation of vector fields vanishing at $0$ on $ \mathbb R^2$. Show that the rotation of center $0$ and of any angle $\theta \neq 0 $ is an inner symmetry. Show that $0$ is a fixed point of that inner symmetry, but is not a very-fixed point.
  \end{enumerate}
  \end{exo}

\begin{exo}
\label{exo:InnerAndSymmetry2}
This exercise is a continuation of Exercise \ref{exo:InnerAndSymmetry}.
 Show that if $m$ is a very-fixed point of an inner symmetry $ \phi$, and is $ \psi$ is a symmetry of $\mathcal F $, then $\psi^{-1}(m) $ is a very-fixed point of the inner symmetry $ \psi^{-1} \circ \phi \circ \psi $.
 \end{exo}

\begin{exo}
 \label{exo:veryfixed} 
 This exercise is a continuation of Exercise \ref{exo:linkedwithanchor}.
Let $ (A,\rho)$ be an anchored  bundle for  $ \mathcal F$ (see Section \ref{sec:AnchoredBundle}). 
Let $m$ be a very fixed point for an inner symmetry obtained as the time-$1$-flow of a vector field $(X_t)_{t \in [0,1]}$ of the form required in Definition \ref{def:veryfixed}.
\begin{enumerate}
\item
Show that the  smooth time-depending section  $(a_t)_{t \in [0,1]}$ such that $ \rho(a_t)=X_t$ for all $ t\in [0,1]$, whose existence is granted by Exercise \ref{exo:linkedwithanchor}, can be chosen to satisfy $ \left.a_t\right._{|_m}=0$ for all $ t\in [0,1]$.
\item Show that any  smooth time-depending section  $(a_t)_{t \in [0,1]}$ such that $ \rho(a_t)=X_t$ satisfies that $ \left.a_{t}\right._{|_m}$ is valued in the strong kernel\footnote{See Section \ref{sec:isotropy} for a definition.} of $ \rho$ at $m$.
\end{enumerate}
\end{exo}

 \begin{remark}
     \normalfont
     For a smooth time-dependent vector field $ X_t$ in $\mathcal F $ of the form described in Definition \ref{def:veryfixed}, notice that for all $t \in [0,1] $, the class $ [X_t]$ of $ X_t$ in the
     isotropy Lie algebra $ \mathfrak g_m(\mathcal F)$ of $ \mathfrak F$ at $m$ is zero (see Section \ref{sec:isotropy} for a definition of this Lie algebra). 
 \end{remark}

\subsection{The rank at a point of a singular foliation}

\noindent
Given a singular foliation on a smooth, real analytic or complex manifold $M$, there are two notions that must not be confused: the rank at that point and the dimension of the tangent space at that point.

The \emph{rank} of an $\mathcal O $-module $\mathcal A $ is the minimal number of its generators.
It is denoted by $rk_{\mathcal O}(\mathcal A)$ and takes values in $\mathbb N \cup \{+\infty\}$.

Let  $m $ be a point in a (smooth, complex, or real analytic) manifold.
We say that a sequence 
$ (\mathcal U_i)_{i\geq 0}$ of open neighborhoods of $m$  \emph{converges to $m$} if for any open neighborhood $\mathcal V $ of $m$, there exists $i_0$ such that for all $ i \geq i_0$, we have $ \mathcal U_i \subset \mathcal V $. The third condition (i.e., “locally finitely generated”) in Definition \ref{def:consensus2}) implies the following result.

\vspace{.5cm}

\begin{propositions}{The rank at a point is well-defined}{rank}
Let $\mathcal F $ be a singular foliation on a smooth, complex, or real analytic manifold $M$.
Let $m \in M$ be a point and $ (\mathcal U_i)_{i \in \mathbb N}$ be a sequence of neighborhoods that converges to $m$.
The sequence
 $$  i \mapsto rk_{\mathcal O_{{\mathcal U}_i}}( \mathcal F_{{\mathcal U}_i}). $$
is finite and constant after a certain rank, and this constant does not depend on the choice of a sequence of open neighborhoods converging to $m$.

It is therefore an integer that depends only on $m$ and $\mathcal F $. It is called the \emph{rank of $\mathcal F $ at $m$}, and denoted by ${\mathrm{rk}}_m (\mathcal F) $.
\end{propositions}

\begin{exo} 
We now work in the smooth case.
Let $\mathcal F $ be a singular foliation of rank less than or equal to $r$ at every point of the manifold $M$. Prove that it is finitely generated.
\emph{Hint:} use Proposition \ref{prop:localglobalnumberofgenerators}.
\end{exo}

\subsection{The tangent space of a singular foliation, and its dimension} 

Let $ \mathcal F$ be a singular foliation on a complex, real analytic or smooth manifold $M$. 

We call \emph{tangent space of $\mathcal F $ at $m \in M$} the subspace of $ T_m M$, denoted by $T_m \mathcal F$, obtained by evaluating at $m$ all vector fields in $ \mathcal F$, defined in any open neighborhood $\mathcal U $ of $m$ in $M$.

\begin{remark}
\label{rmk:tnagentspace}
If a smooth singular foliation $\mathcal F $ is defined through compactly supported vector fields as in Definition \ref{def:consensus}, then $T_m \mathcal F$ is defined by:
 $$ T_m \mathcal F := \{X_{|_m}  \, | \, X \in \mathcal F \}.$$
If a smooth, real analytic or holomorphic singular foliation is defined as a sub-sheaf $\mathcal F_\bullet $ of the sheaf $ \mathfrak X$ of vector fields as in Definition \ref{def:consensus2} or \ref{def:consensus2alg}, then $T_m \mathcal F$ is defined\footnote{Of course, in the smooth setting, both definitions of $ T_m \mathcal F$ coincide.} by:
 $$ T_m \mathcal F := \cup_{\mathcal U \in \mathfrak V_m} \{X_{|_m}  \, | \, X \in \mathcal F_\mathcal U \} $$
 where $\mathfrak V_m$ stands for the set of all open neighborhoods of $m$ in $M$. However, since singular foliation are locally finitely generated, it follows from the axiom that there exists an open neighborhood $\mathcal U $ such that:
  $$  T_m \mathcal F :=  \{X_{|_m}  \, | \, X \in \mathcal F(\mathcal U )\}.  $$
  It suffices to take an neighborhood $\mathcal U $ as in the third item in Definition \ref{def:consensus2}.
\end{remark}

\begin{lemma}
For every point $m \in M$ in a manifold $M$ equipped with a singular foliation $\mathcal F $, the dimension of the tangent space at $m$ is less or equal than the rank of $ \mathcal F$ at $m$. In equation:
 $$  {\mathrm{dim}}(T_m \mathcal F) \leq {\mathrm{rk}}_m (\mathcal F)  $$
\end{lemma}
\begin{proof}
This follows from the discussion in Remark \ref{rmk:tnagentspace}.
\end{proof}

\begin{lemma}
\label{lem:dimTF}
Let $ (M,\mathcal F)$ be a foliated manifold.
The map
$$ \begin{array}{rrcl} {\mathrm{dim}}^{T\mathcal F}: &M &\to  & \mathbb N\\ &m &\mapsto &{\mathrm{dim}}\left(T_m \mathcal F\right)\end{array} $$    
is lower semi-continuous.
\end{lemma}
\begin{proof}

 Let $m \in M$ be a point and  $r =  {\mathrm{dim}}\left(T_m \mathcal F\right)$. 
By definition, there exists $ X_1, \dots, X_r \in \mathcal F$ such that $ X_1|_{m}, \dots, X_r|_{m} $ form a basis of $ T_m \mathcal F$. There exists a neighborhood $\mathcal U $ of $ m$ in $M$ such that for every $ n \in \mathcal U$, 
 $ X_1|_{n}, \dots, X_r|_{n} $ are independent, which implies that $ T_n \mathcal F$ is a vector space of dimension greater or equal to $ r$. 
 This proves that $$\left\{ m '\in M \, | \,  \mathrm{dim}\left(T_{m'} \mathcal F\right) \geq r \right\} $$ is an open subset, which is the content of the lemma. 
\end{proof}

\subsection{The regular part of a singular foliation}

\label{sec:openReg}

Let $\mathcal F $ be a singular foliation on a manifold $M$.

\vspace{.5cm}

\begin{definitions}{Regular point}{regularpoint}
A \emph{regular point} of a singular foliation $ \mathcal F$ on a smooth, complex or real analytic  manifold $M$ is a point $m_0$
in a neighborhood of which the map  ${\mathrm{dim}}^{T\mathcal F} \colon m \mapsto {\mathrm{dim}} \colon(T_m \mathcal F) $ (defined in Lemma \ref{lem:dimTF}) is constant.

A point which is not regular is called a \emph{singular point}.
\end{definitions}
\vspace{0.5cm}

\begin{example}
\label{ex:maxImpliesRegular}
    Any point $m$ where the function $ {\mathrm{dim}}^{T \mathcal F} $ reaches its maximal value is a regular point. 
    This is an obvious consequence of lower semi-continuity (see Lemma \ref{lem:dimTF}).
\end{example}

\begin{remark}
\label{rmk:lowercontinuous}
\noindent
For a smooth singular foliation, there are points which are regular, but are not as in Example \ref{ex:maxImpliesRegular}, i.e., are not points where the map ${\mathrm{dim}}^\mathcal F$
takes its maximal value.

An example is given as follows.
 Let us fix $ \chi(x)$  a smooth function which is $ >0$ on $\mathbb R_+^* $ and equal to $ 0$
on $ \mathbb R_-$.
A singular foliation on $M=\mathbb R $ is given by
 $$ \mathcal F:= \left\{ F(x) \chi(x) \frac{\partial}{\partial x}  \, \middle| \, F(x) \in \mathcal C^\infty(\mathbb R) \right\} .$$
All points of $M=\mathbb R $ are regular except $0$. On $\mathbb R_-^* $, we have 
${\mathrm{dim}}^{T\mathcal F }=0$ while on 
$\mathbb R_+^* $, we have 
${\mathrm{dim}}^{T\mathcal F} =1$.
\end{remark}

\noindent
In the complex or real analytic settings, if $M$ is connected, it is however true that a point $m$ is regular if and only if  $ {\mathrm{dim}}(T_m \mathcal F)$ reaches its maximal value.

By lower semi-continuity of ${\mathrm{dim}}^{T\mathcal F }$, the subset of all regular points of a singular foliation $ \mathcal F$ on a smooth, complex or real analytic  manifold $M$ is an open subset. 
By lower semi-continuity again, it is also a dense subset. 
We call it the \emph{regular part of $\mathcal F$} and denoted by $M_{\mathrm{reg}} $ (at least when there is no ambiguity on the singular foliation that we consider).

\vspace{.5cm}
\begin{propositions}{Regular part}{regularPart}
The regular part $M_{\mathrm{reg}} $  of a singular foliation $ \mathcal F$ is a dense open subset of $M$.

\end{propositions}
\begin{proof}
Let $\mathcal U $ be an arbitrary neighborhood of an arbitrary point $m$. Since the restriction to $\mathcal U $ of the function ${\mathrm{dim}}^{T\mathcal F} $ takes only finitely many values, there exists $m' \in \mathcal U$ where it is maximum. Since it is lower semi-continuous, there is a neighborhood $\mathcal V \subset \mathcal U $ of $m'$ where this function is constant, i.e., a regular point of $ \mathcal F$. 
\end{proof}

\vspace{.5cm}

\begin{exo}
Show that each one of the following statements is wrong:
\begin{enumerate}
\item The set of singular points of a singular foliation is a submanifold of $M$.
\item For every $k \in \mathbb N$, the set $\{m \in M | {\mathrm{dim}}(T_{m} \mathcal F)  =k\} $ is a manifold\footnote{Foliations that satisfy this property have particular features, studied by David Miyamoto's PhD or \cite{zbMATH07807748}.} of $M$. 

\item Singular points are of measure $0$. 
\end{enumerate}

{\textit{Hint}:} For the last question, one can choose $ M = \mathbb R$ and  $ \mathcal F$ the singular foliation generated by $ \chi(x) \frac{\partial}{\partial x}$ where $\chi(x) $ is a smooth function that vanishes on a fat Cantor subset of $ [0,1]$ (or any subset of $ \mathbb R$ of empty interior and non-zero measure).
\end{exo}

By a \emph{Frobenius regular foliation}, we mean the data, on a manifold $N$ of a distribution $ D \subset TN$ which is of constant rank and integrable, i.e., $[\Gamma(D),\Gamma(D)] \subseteq \Gamma(D)$. By the classical Frobenius theorem (cf. e.g.,\cite{zbMATH00041765}), every point has a neighborhood which admits local coordinates on which $D$ is generated by $ \frac{\partial}{\partial x_1}, \dots,  \frac{\partial}{\partial x_k}$.

\vspace{0.5cm}

\begin{proposition}{The regular part of a singular foliation is equipped with a ``good old'' Frobenius regular foliation}{goodold}
Let $(M,\mathcal F) $ be a foliated manifold and $M_{\mathrm{reg}} $ be its regular part. There exists a Frobenius regular foliation $D_{reg}$ on $ M_{\mathrm{reg}}$
such that the sheaves $\mathcal F_{M_{\mathrm{reg}}} $ and $\Gamma(D_{reg}) $ coincide. 
 \end{proposition}
 \begin{proof}
On the regular part, $m \mapsto T_m \mathcal F $ is a smooth distribution $D$ of constant rank, whose sections $\Gamma(D) $ are isomorphic to the restriction  $\mathcal F_{M_{\mathrm{reg}}} $ of the sheaf $ \mathcal F$ to $ M_{\mathrm{reg}}$. Its sections are therefore closed under Lie bracket. Notice that it is also an immediate consequence of the splitting theorems of Section \ref{sec:splitting}.
 \end{proof}

\vspace{0.5cm}
\noindent
In particular, the notions of ``singular foliation such that all points are regular'' and of ``Frobenius regular foliations'' coincide.

\begin{remark}
\label{rmk:Hartog}
Let $ M$ be  a complex manifold and $ \Sigma \subset M$ a subvariety of codimension $ \geq 2$.
By Hartog's principle, a vector field on  $ M \backslash \Sigma$ extends to the whole manifold $M$. Let $D \subset TM$ be a regular integrable distribution on $M \backslash \Sigma $. Its sheaf of sections extends to a sub-sheaf $\mathcal F$ of vector fields on $ \mathfrak X_\bullet$. This sheaf\footnote{Used as a definition of  complex singular foliations, for instance in the on-line book 
“Local dynamics of singular holomorphic foliations” edited by Abate.} satisfies automatically axioms $(\alpha) $ and $ (\beta)$ of Definition \ref{def:consensus2alg}.   
It is not clear to us whether this sheaf automatically satisfies the condition $ (\gamma)$ (but this could be the effect of our limited knowledge). For singular foliations which are of codimension $1$ on their regular part, it is the case in view of Theorem 4.6.2 in \cite{zbMATH07372950}. When $\Sigma$ has codimension $\geq 3$, a theorem of Frisch-Guenot and Siu (cf. e.g., \cite{douadybourbaki} Theorem 2) implies that any coherent sheaf on $M\backslash \Sigma$ extends to a coherent sheaf over $M$, in particular axiom $(\gamma)$ is also satisfied.
\end{remark}

\vspace{0.5cm}

\begin{bclogo}[  arrondi = 0.1, logo = \bcdz]{An alternative definition of regular points!}
Complex geometers may work with an alternative definition of what a regular point is. The issue is that our definition conflicts with a definition of regular points of a coherent sheaf. For instance, for Paul Baum and Raoul Bott's \cite{zbMATH03423310} or Ali Sinan Sert\"oz's \cite{zbMATH04118703}, given $\mathcal F$ a holomorphic singular foliation on a complex manifold $M$, a point $ m\in M$ is regular if it admits an open neighborhood $ \mathcal U \subset M $ such that the restriction $ \mathcal F_{|_\mathcal U}$ is a Debord singular foliation (See Section \ref{sec:Debord}). 
A regular point in our sense is regular in this new sense, but the converse is not true.
 \end{bclogo}

\section{Examples of singular foliations}\label{examples-of-SF}

\vspace{1mm}
\noindent
The purpose of this section is to give an ordered list of examples of singular foliations.

We start with an exercise (which is in fact the content of Section \ref{sec:globally}).

\begin{exo}
\label{exo:compactly_supported}
Let $M$ be a smooth manifold, and $ X_1, \dots, X_r \in \mathfrak X(M)$ be vector fields such that there exists functions $ (c_{ij}^k)_{i,j,k \in \{1, \dots, r} \in \mathcal C^\infty(M)$ with
$$ [X_i,X_j] = \sum_{k=1}^r c_{ij}^k X_k.  $$
\begin{enumerate}
\item
Show that the space $ \mathcal F $ of vector fields of the form $\sum_{i=1}^r f_i X_i $, with $ f_1,\dots,f_r$ compactly supported functions on $ M$, is a singular foliation\footnote{In the sense of Definition \ref{def:consensus}} on $ \mathcal F$.
\item 
State an equivalent result in the real analytic and complex settings, using sheaves.
\end{enumerate}
\end{exo}

We also invite the reader to do the next exercise, in order to find counter-examples.

\begin{exo}
\label{exo:finitelygeneratedex}
Consider the algebra of smooth functions on $ \mathbb R^d$. 
\begin{enumerate}
\item 
Show that smooth functions on $ \mathbb R^d$ vanishing at $0$ together with their $k$ first partial derivatives, i.e., 
     $$ \mathcal I^k := \left\{ F \in \mathcal C^\infty(\mathbb R^d) \, \middle| \, \frac{\partial^{i_1 + \dots + i_d } F}{ \partial x_{1}^{i_1} \cdots \partial x_{d}^{i_d} }(0, \dots, 0)  =0 \hbox{ for all } i_1 , \dots, i_d \in \mathbb N_0 \hbox{ with } i_1+ \cdots + i_d \leq k \right\}  $$
     is an ideal of $C^\infty(\mathbb R^d) $ which is finitely generated.
     Describe explicit generators.
     \item Show that the ideal 
 $$ \mathcal I^\infty := \left\{ F \in \mathcal C^\infty(\mathbb R^d) \, \middle| \, \frac{\partial^{i_1 + \dots + i_d } F}{ \partial x_{1}^{i_1} \cdots \partial x_{d}^{i_d} }  (0, \dots,0)=0 \hbox{ for all } i_1 , \dots, i_d \in \mathbb N_0 \right\}  $$
     is not finitely generated. {\emph{Hint:}} Not easy: one solution is to show that if an ideal $ \mathcal I$ is finitely generated, then the flow  $ \phi_X^t$  of any complete vector field $X$ such that $ X[\mathcal I] \subset \mathcal I$ satisfies $\phi_t^* (\mathcal I)=\mathcal I$. Then show that this property is not satisfied\footnote{To simplify, one can assume that $ X=\frac{\partial}{\partial x_1}$ from the very beginning.} for $X:= \frac{\partial}{\partial x_1}$.
     \end{enumerate}
\end{exo}

\subsection{Regular foliations}

\vspace{1mm}
\noindent
Although it seems grammatically problematic, regular foliations are examples of singular foliations. This point was in fact established in
Section \ref{sec:openReg}.
This is valid in smooth, complex or real analytic settings.

\subsection{Singular foliations and Lie algebroids}

\label{sec:LieAlgebroidsAreSingFoliation}

\vspace{1mm}
\noindent
Let $M$ be a manifold with sheaf of functions $ \mathcal O_M$ (smooth, real-analytic or holomorphic depending on the context).

\vspace{0.2cm}
{\textbf{Convention:}} In this section, $ \mathcal O_\mathcal U$ (rather than $ \mathcal O(\mathcal U)$) will stand for functions over an open subset $ \mathcal U$.
This convention is chosen to be consistent with the convention for sections over  $\mathcal U $ of a vector bundle $A\to M$, namely $ \Gamma_\mathcal U (A)$.
\vspace{0.2cm}

Recall that a \emph{Lie algebroid over $M$} \cite{Mackenzie} is a triple $(A,\rho,[\cdot\,, \cdot])$, with $A$ a vector bundle over $M$, $\rho \colon A \to TM $ a vector bundle morphism over the identity of $M$ called \emph{anchor map}, and $[\cdot\,, \cdot] $ a Lie bracket on the sheaf of sections of $A$ such that  the so-called \emph{Leibniz identity} holds for all $ a,b \in \Gamma(A), f \in \mathcal O_M$:
 $$  [a,fb]=f \, [a,b] + \rho(a)[f] \, b.$$
 Lie algebroids have been extensively reviewed and studied, see e.g., \cite{Mackenzie,CFM,MeinrenkenAlgebroids}. 
We will show that any Lie algebroid induces a singular foliation.
The following Lemma holds true in the smooth, real-analytic or holomorphic categories.

\begin{lemma}
\label{lem:anchorismorphis}
 For any Lie algebroid   $(A,\rho,[\cdot\,, \cdot])$  over $M$, and any open subset $ \mathcal U \subset M$ the anchor map $\Gamma_\mathcal U(A) \to \mathfrak X(\mathcal U) $ is a Lie algebra morphism.
\end{lemma}

\begin{proof}
Recall that $ \Gamma_\mathcal U$ stands for the $ \mathcal O_\mathcal V$-module of sections of $A$ over an arbitrary open subset  $ \mathcal V \subset M$.
The Jacobi identity on $ \Gamma_\mathcal V (A)$ implies that for any $a,b,c \in \Gamma_\mathcal V (A)$ and any $ f \in \mathcal O_\mathcal V$:
 \begin{equation}
 \label{eq:outOfJacobi}
 \left({\mathrm{ad}}_a \circ {\mathrm{ad}}_b  - {\mathrm{ad}}_b \circ {\mathrm{ad}}_a\right) (fc)  = {\mathrm{ad}}_{[a,b]} (fc) .
 \end{equation}
 Now, there are three kinds of terms that will appear in the previous equation \eqref{eq:outOfJacobi} if one uses the Leibniz identity to allow the function $f$ to ``get out'':
\begin{enumerate}
    \item[1] Those for which $f$ will be ``differentiated'' twice: these terms appear on the left-hand side of the equation \eqref{eq:outOfJacobi} as
     $ \rho(a) \circ \rho (b) [f] \, -   \rho(b) \circ \rho (a) [f] \, c $.
    \item[2] Those for which $f$ will be ``differentiated'' only once, there are two types of such terms  
    \begin{enumerate} 
    \item[(2a)] $ \rho([a,b]) [f] \, c$ on the right-hand side of \eqref{eq:outOfJacobi}
    \item[(2b)] $\rho(a)[f] \ [b,c] + \rho(b)[f] \ [a,c] -\rho(a)[f] \ [b,c] - \rho(b)[f] \ [a,c]   $ on the left-hand side of \eqref{eq:outOfJacobi}. These terms cancel out.
    \end{enumerate}
    \item[3] Those where $f$ is not ``differentiated'' at all, those terms are 
    \begin{enumerate} 
    \item[(3a)] $ f \, \left({\mathrm{ad}}_a \circ {\mathrm{ad}}_b  - {\mathrm{ad}}_b \circ {\mathrm{ad}}_a\right) (c)   $ on the left-hand side of \eqref{eq:outOfJacobi}
    \item[(3b)] $f \, \,{\mathrm{ad}}_{[a,b]} (c) $ on the right-hand side of \eqref{eq:outOfJacobi}
    \end{enumerate}
\end{enumerate}
The Jacobi identity being satisfied for triples of sections of $A$, the terms (3a) and (3b) cancel each other out. Also, the terms (2b) add up to zero.
Hence, the terms 1 and (2a) are the only ones remaining. They transform  Equation \eqref{eq:outOfJacobi} into the relation: 
 $$ \left( (\rho([a,b]) - \left[\rho(a),\rho(b)\right]) [f] \right) \, c = 0  . $$
 Now, let us assume that $\mathcal V $ is chosen such that a nowhere vanishing section $c \in \Gamma_\mathcal V  (A) $ exists.
We then have $(\rho([a,b]) - \left[\rho(a),\rho(b)\right]) [f] =0$.
Assume also that  $\mathcal V$ is chosen such that any covector is the differential of at least a function $f \in \mathcal O_V$ (which is always true in the smooth case, but is only true for $ \mathcal V$ ``small'' enough in the complex case). This implies that $ \rho$ is a Lie algebra morphism from $\Gamma_\mathcal V (A) $ to $\mathfrak X (\mathcal V) $. Since every point $m \in M$ admits such a neighborhood $\mathcal V $, however, this implies that $\rho $ is a Lie algebra morphism when restricted to any open subset $\mathcal U \subset M$.
\end{proof}

At this point, it is more convenient to distinguish the smooth case from the complex and real analytic ones.

\subsubsection{The smooth case}

\vspace{1mm}
\noindent
Let us consider that singular foliations on a smooth manifold $M$ are defined as in Definition \ref{def:consensus2}, through compactly supported vector fields.
Let $ (A \to M, [\cdot\,, \cdot], \rho)$ be a smooth  Lie algebroid over $M$.
Let $\mathcal F := \rho\left( \Gamma_c(A)\right) $, where $\Gamma_c $ stands for compactly supported sections. Lemma \ref{lem:anchorismorphis} implies that $\mathcal F $ is closed under Lie bracket. The remaining axioms are obviously satisfied. Hence, the following result holds true.

\vspace{.5cm}

\begin{propositions}{Image through anchor map of Lie algebroids - smooth case}{prop:LAareFoliations}
The image through the anchor map of compactly supported sections of a Lie algebroid over $M$ form a singular foliation on $M$.
\end{propositions}

\vspace{.5cm}
\subsubsection{The complex or real-analytic case}

Proposition \ref{thm:prop:LAareFoliations} can not be extended immediately from the smooth context to the complex or real analytic contexts altogether\footnote{The presentation here also essentially works on affine varieties or schemes.}.
We denote by $ \mathcal O$ the sheaf of real-analytic or holomorphic functions.
Again, for $A\to M $ a vector bundle, we  denote by $\Gamma_{\mathcal U}(A) $ the sections of $A $ over an open subset $\mathcal U$. Of course, $ \Gamma_{\mathcal U}(A) $ is a $ \mathcal O_{{\mathcal U}}$-module, and, assigning to an open subset the sections over it
 $$ {\mathcal U} \mapsto \Gamma_{\mathcal U}(A)  ,$$
one defines  a sheaf of $\mathcal O $-modules over $M$. 
The technical difficulty that appears at this point is that singular foliations are supposed to be sheaves, but 
 $$ {\mathcal U} \mapsto \rho(\Gamma_{\mathcal U}(A))  $$
 is not a sheaf on $M$ (and therefore not a sub-sheaf of the sheaf $\mathfrak X $ of vector fields on $M$). It is only a pre-sheaf\footnote{For instance, if $A$ has no globally defined non-constant functions, e.g.$A$ is trivial and $M$ is compact.}. But the difficulty can be circumvented:
 to turn it into a sheaf, one has to ``sheafify'' it, i.e., to map an open subset $\mathcal U \subset M$ to the sub-$\mathcal O_\mathcal U $-module of vector fields $X \in \mathfrak X(\mathcal U)$ on $\mathcal{U}$ such that  every $m \in \mathcal U$ admits a neighborhood $\mathcal V $ on which there exists $a \in \Gamma_{\mathcal V} (A)$ with  $ \rho(a) =X  $ (on $\mathcal V $).
This defines a sheaf of $\mathcal O $-modules $\underline{\rho(\Gamma(A) )}$ that we call the \emph{image of the Lie algebroid $(A,\rho, [\cdot\,, \cdot])$  through its anchor map}.

\vspace{.5cm}

\begin{propositions}{Image through anchor map of Lie algebroids: complex case}{prop:LAareFoliations2}
Let $(A,\rho, [\cdot\,, \cdot])$ be a Lie algebroid over a complex or real analytic manifold $M$.
The sheaf  $\underline{\rho(\Gamma(A) )}$ 
(=image of the Lie algebroid $(A,\rho, [\cdot\,, \cdot])$  through its anchor map) is a singular foliation on $M$.
\end{propositions}
\begin{proof}
This statement is an immediate consequence of Lemma \ref{lem:anchorismorphis}.
\end{proof}

\begin{exo}
\label{exo:rank1}
Let $X$ be a  vector field on a manifold $M$. Show that $ \mathcal F_X = \{f X | f \in \mathcal C^\infty_c(M) \}$ is a singular foliation on $M$ that comes from a Lie algebroid of rank $1$.
\end{exo}

\begin{exo}(See \cite{CamilleLouis}, Section 3.1.4)
Let $\mathcal F $ be a singular foliation and $\varphi \in \mathcal C^\infty(M) $ a function. 
Check that  $$  \varphi \mathcal F := \{ \varphi X, X \in \mathcal F\} $$
is a singular foliation again. Show that if $ \mathcal F$ is the image through the anchor map of a Lie algebroid, so is $ \varphi \mathcal F$.
\end{exo}

\subsubsection{Lie algebra actions}

\noindent 
Let $\mathfrak g$ be a Lie algebra. We call a Lie algebra morphism $ \mathfrak g \to \mathfrak X (M) $  a \emph{Lie algebra action} of $ \mathfrak g$ on $M$. We denote it by $ x \mapsto \underline{x}$.  In particular, any Lie group action of  a Lie group $G$ on a manifold $M$ induces a Lie algebra action of its Lie algebra $ \mathfrak g$.
The $ \mathcal C^\infty_c(M)$-module\footnote{Recall that $ \mathcal C^\infty_c(M)$ stands for compactly supported smooth functions.} generated by $\{\underline{x}, x \in \mathfrak g\} $ is a singular foliation.
We have therefore proven the first part of the following result:

\vspace{0.3cm}

\begin{propositions}{Lie group actions}{LieALgebraActions}
   The infinitesimal (Lie algebra) action of a Lie group action $G$ on a manifold $M$ induces a singular foliation on $M$.
    It is the image through the anchor map of a Lie algebroid called transformation Lie algebroid.
\end{propositions}
\begin{proof}
The first part of the statement is already proven. For the second part, consider the trivial vector bundle $ \mathfrak g \times M \to M$ equipped with the anchor $ \rho$ and bracket $ [\cdot\,, \cdot]$ defined on constant sections, which are identified with elements of $\mathfrak g $, by
 $$ \rho (x) = \underline{x}  \hbox{ and } [x,y]=[x,y]_{\mathfrak g} $$
 where $ x,y \in \mathfrak g$ are seen as constant sections of $A$, and $[x,y]_{\mathfrak g}$ is the bracket of  $\mathfrak g $.  This describes a Lie algebroid structure \cite{Mackenzie}.
\end{proof}

\vspace{.5cm}
\subsubsection{Projective or ``Debord'' singular foliations}
\label{sec:Debord}
Here is an important class of singular foliations that come from a Lie algebroid. We will state the results in the smooth case, and leave the generalization to the reader for the complex or real analytic settings.

\vspace{.5cm}

\begin{definitions}{\cite{Debord} Generators and no relations}{Debord}

We say that a singular foliation on a smooth manifold $M$ is \emph{Debord} if $\mathcal F $ is a projective $\mathcal C^\infty(M) $-module\footnote{In the complex or real analytic settings, one has to require $\mathcal F $ to be a projective sheaf with respect to the sheaf of functions. Since the rank of $ \mathcal F$ is (locally) finite, it is equivalent to say that every point has a neighborhood $\mathcal U$ where $ \mathcal F_\mathcal U$ is a free $ \mathcal O_\mathcal U$-module.}.
\end{definitions}

\vspace{.5cm}

\noindent
In a concrete manner, Debord foliations are those which admit, in a neighborhood $\mathcal U $ of every point, generators $X_1, \dots, X_r $ between which there is no relation. I.e., if 
 $$ \sum_{i=1}^r f_i X_i =0 ,$$
 then all of the functions $f_1, \dots, f_r \in \mathcal C^\infty(\mathcal U) $ are zero.

\begin{remark}
Equivalently, we could use the language of Definition \ref{def:consensus2} to define Debord foliations in the smooth case. 

\end{remark}

\begin{exo}
Show that the restriction of a Debord singular foliation to an open subset is still a Debord singular foliation.
\end{exo}

\begin{remark}
In addition of honoring the important discoveries of Claire Debord about them in \cite{Debord}, the name ``Debord foliation'' is encouraged by a very practical perspective. Saying ``projective foliations'' would be ambiguous, since  it could mean ``foliations on a projective variety''. 
\end{remark}

\noindent
Debord foliations are dealt with in this section, because they always arise from a Lie algebroid. By the smooth Serre-Swan theorem \cite{zbMATH07201776}, there exists a vector bundle $A \to M $ and a  $\mathcal C^\infty(M) $-module isomorphism
$$\Gamma_c(A) \simeq \mathcal F. $$
Composing this $ \mathcal C^\infty_c(M)$-module isomorphism with the inclusion 
$$\Gamma_c(A) \simeq \mathcal F \, \hookrightarrow \, \mathfrak X_c(M), $$ 
we obtain an inclusion $\Gamma_c(A) \hookrightarrow \mathfrak X_c(M) $. Since it is  a morphism of $\mathcal C_\mathcal C^\infty(M) $-modules, it   has to be given by a vector bundle morphism:
 $$ \rho \colon A \to TM ,$$
 that we call anchor map.
 The anchor $\rho $ is injective at the level of sections. 
  This does \underline{not} imply that the anchor $ \rho$ needs to be injective at all points of $M$, but it certainly has to be injective on a dense open subset of $M$. 
 Lastly, the isomorphism $\Gamma_c(A) \simeq \mathcal F$ extends to non-compactly supported sections, and equips $ \Gamma(A)$ with a Lie bracket, for which $ \rho$ is easily seen to be an anchor map. This proves the following statement: 
\vspace{0.5cm}
\begin{propositions}{Debord algebroids}{Debord} 
A singular foliation on a smooth, real analytic or complex manifold $M$ is Debord if and only if it is the image of a Lie algebroid whose anchor map is injective on a dense open subset of $M$.
\end{propositions}
\vspace{0.5cm}

\begin{exo}
Show that the singular foliation on $ \mathbb R^2$
generated by $ \frac{\partial}{\partial x}$ and    $ \frac{\partial}{\partial y}$ is Debord. 
\end{exo}

\begin{exo}
Show that the singular foliation in Exercise \ref{exo:rank1} is Debord, at least if there is no open subset where the vector field $X$ is identically zero. Show that the anchor is not injective at a point $m$ where $X_{|_m}=0$. 
\end{exo}

\begin{exo}
\label{exo:codim1}
Show that compactly supported vector fields on a manifold $M$ vanishing on a codimension $1$ submanifold form a Debord singular foliation. This theme will be developed further in section \ref{sec:submanif} below.
\end{exo}

\begin{example}
    \normalfont
    An interesting example of a Debord singular foliation will be given in Section \ref{sec:hormander}, see Example \ref{ex:YErp}.
\end{example}

\subsection{Vector fields vanishing at a point at prescribed order}

\label{ex:singFolVanish}

\vspace{1mm}
\noindent
We can also construct singular foliations by playing with the order of vanishing of vector fields at certain points.
We decided to place ourselves within the context of smooth differential geometry in the discussion, but the holomorphic or real-analytic contexts would work as well in more  or less a similar manner, and the conclusions are written in a way that merges all three contexts.
Let us start with an exercise.

\begin{exo}
Let $\mathcal F_1  $ be the space of all compactly supported smooth vector fields on $\mathbb R^n $ vanishing at $0$.
Show that  $\mathcal F_1  $  is a finitely generated singular foliation (see Definition \ref{def:globallyfinitelygenerated}) generated by the finite family of vector fields
 $$  \left\{\,  x_i \frac{\partial}{\partial x_j} \, \middle| \, i,j \in \{1, \dots,n\} \, \right\}. $$
\emph{Hint:} use the so-called ``Hadamard's lemma'', i.e., the fact that any compactly supported smooth function $F$ on $\mathbb R^n $ vanishing at $0$ decomposes as 
  $$ F= \sum_{i=1}^n x_i F_i  $$
  for some compactly supported smooth functions $F_1, \dots, F_n \in \mathcal C^\infty(\mathbb R^n)$.
\end{exo}

This exercise can be easily generalized.
We say that a vector field $X$ on $\mathbb R^n $ vanishes to order $2$ at the origin if its coefficients
 $$ X = \sum_{i=1}^n F_i  (x_1, \dots, x_n) \frac{\partial}{\partial x_i}$$
 satisfy that
  $$  F_i(0, \dots, 0)=0 \hbox{ and } \frac{\partial F_i}{\partial x_j}(0, \dots, 0) =0 \hbox{ for all $ i,j=1, \dots, n$}.$$
  It is a classical result that a compactly supported smooth function $F$ on $\mathbb R^n $ vanishes at $0 \in \mathbb R^n$ together with its differential if and only if it
decomposes as 
  $$ F= \sum_{i,j=1}^n x_i x_j F_{i,j}   $$
  for some compactly supported smooth functions $F_{i,j} \in \mathcal C^\infty(\mathbb R^n)$.
Let $\mathcal F_2 $ be the space of such vector fields. We leave it as an exercise to the reader to check that
  \begin{enumerate}
  \item Show that $\mathcal F_2 $
 is generated, as a $\mathcal C^\infty_c(\mathbb R^n) $-module, by the family 
 $$  \left\{\, x_i x_j \frac{\partial}{\partial k} \, \middle| \, 1 \leq i\leq j \leq n \hbox{ and } k=1, \dots, n \, \right\}.$$
 \item Show that $\mathcal F_2 $ is stable under Lie bracket. \item Conclude that $\mathcal F_2 $ is a finitely generated (see Definition \ref{def:globallyfinitelygenerated}) singular foliation on $\mathbb R^n $.
 \end{enumerate}

\vspace{1mm}
\noindent
Let us generalize the previous two exercises.
For every $k \in \mathbb N$, let $\mathcal F_k $ be the space of compactly supported  smooth vector fields on $\mathbb R^n $ that vanish at $0$ together with their partial derivatives of order $i$ for $i \leq k-1 $.
\begin{enumerate}
    \item It is straightforward that $\mathcal F_k  $ is a $\mathcal C^\infty(\mathbb R^n)$-module stable under Lie bracket.
    \item There is an identification $\mathcal F_k = \mathcal I_0^k \mathfrak X_c(M)$ where $\mathcal I_0$ is the ideal of smooth functions on $\mathbb R^n $ vanishing at the origin.
    \item An explicit family of generators of $\mathcal F_k $ over  $\mathcal C^\infty(\mathbb R^n)$ is therefore given by:
     $$  \left\{\, x_{i_1} \cdots x_{i_k} \frac{\partial}{\partial x_j} \, \middle| \, 1 \leq i_1\leq \cdots \leq i_k \leq n \hbox{ and } j=1, \dots, n \, \right\}.  $$
\item $\mathcal F_k $ is therefore a finitely generated singular foliation in the sense of Definition \ref{def:globallyfinitelygenerated}.
\end{enumerate}

For every $k \geq 1$, the singular foliation $ \mathcal F_k$ can also be seen as a complex, real analytic or algebraic singular foliations on $\mathbb K^n $ with $\mathbb K=\mathbb R $ or $\mathbb C $ depending on the context. In all these realms, it is true that
 $\mathcal F_k  = \mathcal I_0^k \mathfrak X_\bullet $ where $\mathcal I_0$ is the ideal of relevant (sheaf of) functions vanishing at the origin.  Of course, this discussion can be enlarged to any point in a  smooth or complex manifold. In conclusion:

\vspace{.3cm}

\begin{propositions}{Vector fields vanishing at given orders at given points}{eq:orderk}
Let $M$ be a smooth, real analytic or complex manifold. For every point $ m \in M$, and every choice of an integer $ k \geq 1$, the space of vector fields on $M$ vanishing together with their $k-1$ first derivatives at $m$  form a singular foliation on $M$.
\end{propositions}
\vspace{0.3cm}

\begin{exo}
Find all possible relations between the generators of $\mathcal F_k $, for any $ k \geq 1$.
\end{exo}

\begin{exo} This exercise supposes the notion of leaves.
Find its leaves of each one of the singular foliations $\mathcal F_k $ for any $ k \geq 1$. 
\end{exo}

\begin{exo} Can $k$ be replaced by $ + \infty$ in Proposition \ref{thm:eq:orderk}?
\emph{Hint:} Look at the second item in \ref{exo:finitelygeneratedex})
\end{exo}

Proposition \ref{thm:eq:orderk} generalizes to a family of distinct points $ m_1, \dots, m_d $ in an arbitrary manifold $ M$, and \emph{orders} $k_1, \dots, k_d \in \mathbb N$. In the smooth case, one can then consider all compactly supported vector fields on $M$ whose order of vanishing at the point  $ m_i$ is greater or equal to $k_i$ for all $i=1, \dots, d$. One can even choose a sequence $ (m_i)_{i \in \mathbb N}$ of points and $ (k_i)_{i \in \mathbb N}$ orders, provided that it has no accumulation point. In particular, the Androulidakis-Zambon's ``non-finitely-many-generators'' singular foliation of exercise \ref{exo:infinitestability} is of that type.

\subsubsection{More sophisticated examples}

This example (inspired by Grabowska and Grabowski \cite{GG}), appeared in \cite{Ryvkin2}, Example 1.11.
We present it as a real analytic singular foliation on $\mathbb R^n$ (we could see of course also see it as a complex singular foliation on $\mathbb C^n$ or a smooth one).

On $M={\mathbb R}^n $, we attribute to the canonical coordinates $(x_1, \dots, x_n)  $ the strictly positive weights $(i_1, \dots, i_n) $. Equipped with this weight, the ring $\mathcal A $ of real analytic functions on $M$  becomes a graded algebra.
 $$ \mathcal A  = \oplus_{i=0}^\infty   \, \mathcal A_i.  $$
 It is also a filtered algebra, with respect to the filtration:
  $$ \mathcal A^{\geq k}  = \oplus_{i=k}^\infty   \, \mathcal A_i .$$

\begin{example}
Assume $i_1=1, i_2 =2 $ and so on. 
The {weight} of $x_1^3 x^2_3 x_5 $ is $ 1 \times 3+ 3 \times 2 + 1 \times 5 = 14 $, so that $x_1^3 x^2_3 x_5  \in \mathcal A_{14}$.
\end{example}

Let $k$ be a non-negative integer. The space of real analytic vector fields $X$ such that
 $$ X \left[ \mathcal A^{\geq n}\right] \subset \mathcal A^{\geq n+k} \hbox{ for all $ n \in \mathbb N$ }$$
is a module, that we denote by $\mathcal F_{k} $, over real analytic functions. It is stable under Lie bracket.
It is generated by the family
$$   \left\{ ~  x_1^{j_1}  \dots x_n^{j_n} \tfrac{\partial }{ \partial x_a } ~\middle| 
i_1 j_1+  i_2 j_2 + \cdots +  i_n j_n  \geq j_a + k ~ \right\}.$$

If $ (j_1, \dots, j_n)$ satisfies the above condition, then so does $ (j_1', \dots, j_n')$ as long as $ j_i' \geq j_i$ for all indices $i=1, \dots, n$. 
This implies that the generating family of $\mathcal F_k $ can be chosen to be finite.
Therefore, $\mathcal F_k $ is finitely generated, and is a real analytic singular foliation.

\subsection{Singular foliations attached to a sub-variety (I): the algebraic case} 

\label{sec:algebraiccase_subvariety}

We now work within the context of complex algebraic geometry. Let 
$\mathcal O$ be the algebra of polynomial functions on an affine variety $M $. 
Recall that $\mathcal O$  is a quotient of the form:
  \begin{equation}\label{eqdef:affinevariety}  \mathcal O = \frac{\mathbb C[x_1, \dots, x_n]}{\mathcal I_M} \end{equation}
  with $ {\mathcal I}_M$ a prime ideal of $\mathbb C[x_1, \dots, x_n]$.

\begin{example}\normalfont
The reader not familiar with algebraic geometry can assume $M=\mathbb C^n$ so that $\mathcal O= \mathbb C[x_1, \dots,x_n]$ is the algebra of polynomials
in $n$ variables.
  \end{example}

Recall that, by definition, the $ \mathcal O$-module $\mathfrak X (M) $ of vector fields on $M$ is the $\mathcal O $-module $\mathfrak X(M) $  of derivations of $\mathcal O $. 

\begin{example}\normalfont
For $M = \mathbb C^n $, vector fields are simply expressions of the form  
 $$ \sum_{i=1}^n P_i(x_1, \dots, x_n) \frac{\partial}{\partial x_i} ,$$
and are uniquely determined by 
the polynomial functions $ (P_i(x_1, \dots, x_n))_{i=1, \dots,n}$. As an $ \mathcal O $-module therefore, $\mathfrak X(M)= \mathcal O^n $.
\end{example}

We recall the following Lemma.

\begin{lemma}    
\label{lem:isnotherian}
The $ \mathcal O$-module $ \mathfrak X(M)$ of vector fields  on an affine variety is of finite rank. In particular, it is Noetherian. 
\end{lemma}

\begin{proof}
Any vector field $X$ on $M$ is determined by its values on the functions $\bar{x}_1 , \dots, \bar{x}_n$ (the horizontal refers to the quotient in  Equation \eqref{eqdef:affinevariety}).
In particular, vector fields are a finitely generated $ \mathcal O $-module. Now, since $ \mathcal O$ is Noetherian, so is any finite rank $  \mathcal O$-module, which concludes the proof.
\end{proof}

Let $W \subset M$ be an affine sub-variety, i.e., the zero locus\footnote{=the subset of points where all elements in $ \mathcal I_W$ vanish.} of some prime ideal
 $\mathcal I_W \subset \mathcal O$.
Since $\mathcal O $ is Noetherian, this ideal has finitely many generators $\varphi_\bullet = (\varphi_1, \dots, \varphi_k) $. 

\begin{example} \normalfont
On $M=\mathbb C^n $, an affine subvariety $W$ is a subset given by  $\phi_1= \dots = \phi_k =0$ where $\phi_1, \dots, \phi_k $ generate a prime ideal.
\end{example}

In algebraic geometry, geometrical properties have to be translated in a purely algebraic language. For instance, we say that a vector field $X$ \emph{vanishes on $W$} if $ X[\mathcal O] \subset \mathcal I_W$ and is \emph{tangent to $W$} if $X[\mathcal I_W] \subset \mathcal I_W$. A vector field vanishes on $W$ if and only if it belongs to $ \mathcal I_W \mathfrak X(M)$.

\begin{center}
\begin{tabular}{|lcl|} \hline
    $X$ vanishes on $W$& $ \Leftrightarrow  $ & $X(w)=0  $ $ \forall w \in W$,\\
    $X$ is tangent to $W$ & $\Leftrightarrow$  & $X \in T_w W $ for every regular point of $W$. \\
    \hline
\end{tabular}
\end{center}
Here is our main statement:

\vspace{0.3cm}

\begin{propositions}{Two foliations associated to an affine variety}{affine} 
Let $W \subset \mathbb C^n $ be an affine variety. Both vector fields on $ \mathbb C^n$  tangent to $W $ and vector fields on $ \mathbb C^n$ vanishing on $W$ are algebraic singular foliations\footnote{I.e., are finitely generated sub-$\mathcal O$-modules of vector fields on $M$ (=derivations of $\mathcal O = \mathbb C[x_1, \dots,x_n] $) stable under Lie bracket.}.
\end{propositions}

\vspace{0.3cm}

Proposition \ref{thm:affine} is a direct consequence of the following more general result.

\begin{proposition}

For any ideal $\mathcal I \subset \mathcal O $,
the following families are algebraic singular foliations.
\begin{enumerate}
    \item The   $\mathcal O $-module $\mathfrak X_{\mathcal I} $ of all vector fields $X \in  \mathfrak X(M)$ such that $ X[\mathcal I] \subset \mathcal I$.
    \item $ \mathcal I \mathfrak X(M)$.
\end{enumerate}
\end{proposition}
\begin{proof}
By Lemma \ref{lem:isnotherian}, it suffices to check that the previous sets are $\mathcal O$-modules stable under Lie bracket, which is straightforward.
\end{proof}

\begin{exo}
\label{exo:forfrancis}
Let ideal $\mathcal I _W\subset \mathcal O =\mathbb C[z_1, \dots,z_N]$ be the ideal of functions vanishing on an affine variety $W$.
\begin{enumerate}
    \item 
Show that the space of all vector fields $X$ on $ \mathbb C^N$ such that $ X[\mathcal I _W] \subset \mathcal I _W^2$ is an algebraic singular foliation.
\item Is this algebraic singular foliation really different from the algebraic singular foliation of vector fields vanishing on $ W$?
\end{enumerate}
\end{exo}

There are of course more examples.
Let $ M'$ be a second affine variety with algebra of functions $ \mathcal O'$.
Let $\phi \colon M \to M'$ be an affine map, i.e., a map such that $ \phi^* \colon \mathcal O' \to \mathcal O$ is an algebra morphism.
A  vector field $X \in  \mathfrak X(M)$ such that $ X[\varphi^* (F)]=0 $ for all $F \in \mathcal O' $ is called \emph{tangent to the fibers of $ \phi$}. 

\begin{example}
For $M = \mathbb C^n$ and $ M' = \mathbb C^{n'}$, 
$ \phi$ is given by a $n'$-tuple of polynomial functions, i.e.,
 $$ \phi \colon (z_1, \dots, z_n) \mapsto  \left(\phi_1(z_1, \dots, z_n) , \dots, \phi_{n'}(z_1, \dots, z_n)\right)  $$
 and vector fields tangent to fibers are those vector fields $X $ such that $X [ \phi_1]= \cdots =X[\phi_{n'}] =0 $.
\end{example}

Vector fields tangent to the fibers of $ \phi$ form an $ \mathcal O $ -module stable under Lie bracket. Here is an obvious consequence of Lemma \ref{lem:isnotherian}.

\begin{proposition}
    For any map of affine varieties $ \phi \colon M \to M'$, vector fields on $M$ tangent to fibers of $ \phi$ form an algebraic singular foliation on $M$.
\end{proposition}

\vspace{0.5cm}
Let us finish this section with the following remark.

\begin{remark}
When $M$ is a smooth affine variety, i.e., when $M$ has no singular points (e.g., $ M = \mathbb C^n$), then $M$ is also a complex manifold.
All algebraic singular foliations constructed above, and more generally any algebraic singular foliation $\mathcal F$ on $M$, may be seen as a complex singular foliation: it suffices to consider the sheaf of all vector fields which are linear combinations, with coefficients in holomorphic functions, of vector fields in $\mathcal F$. 
In short, it suffices to take the tensor product of $ \mathcal F$ with holomorphic functions.
\end{remark}

\subsection{Singular foliations attached to a submanifold (II): the smooth or complex cases} 
\label{sec:submanif}

This section makes sense in the smooth, real analytic or complex contexts indifferently.

\vspace{0.3cm}

\begin{propositions}{Vector fields tangent to $L$ of vanishing along $L$}{Tgt}

Let $L$ be a submanifold of $M$, and $k \in \mathbb N $ an integer. Both
\begin{enumerate}
    \item the space of vector fields tangent to $L$,
    \item  and the space of vector fields vanishing at order $k$ at all points in $ L$ 

\end{enumerate}
 are singular foliations.

\end{propositions}
\begin{proof}
The proof consists in
\begin{enumerate}
    \item Checking that  the space $\mathfrak X_L(M) $ of all vector fields on $M$ tangent to the sub-manifold $L$, \begin{enumerate}
        \item is a module over functions,
        \item is stable under Lie bracket,
   \item and that,
    in any local coordinates $(x_1, \dots,x_a, y_1, \dots, y_b )$ where $L $ is given by $0=y_1= \dots =y_b $, it is generated by $$\left\{ \frac{\partial}{\partial x_i}, y_j \frac{\partial}{\partial y_k} \middle| 1 \leq i \leq a \hbox{ and } 1 \leq j,k \leq b   \right\} $$
     \end{enumerate}

    \item Then in checking that the second space is algebraically described by $ \mathcal I_L^k \mathfrak X(M) $, where $\mathcal I_L $ stands for the ideal of functions vanishing on $L$. Since the ideal $\mathcal I_L$ is locally finitely generated, this completes the proof. 
\end{enumerate}
\end{proof}

\begin{exo}
Let $L_1, L_2 \subset M$ be submanifolds of $M$ that intersect transversally, i.e., such that:
 $$ T_x L_1 + T_x L_2 = T_x M \hspace{1cm}  \forall x \in L_1 \cap L_2 .$$
Consider the space of all vector fields on $ M $ tangent to both $L_1  $ and $L_2$. Show that it is a singular foliation.
\end{exo}

\begin{exo}
We owe to \cite{francis2023singular}-\cite{BDW} the idea of the following exercise\footnote{Related to Exercise \ref{exo:forfrancis}.}.
Let $M$ be a smooth manifold, and $ \Sigma$ a submanifold defined as the zero locus of a function $ \phi \colon M \to \mathbb R$.
We assume the differential of $\phi $ to be non-zero at each point of $ \Sigma$.
\begin{enumerate}
\item 
Show that the space of all vector fields $X$ such  that
 $ X[\phi] \in \mathcal I_\Sigma^2$, with $ \mathcal I_\Sigma$ the ideal of smooth functions vanishing on $ \Sigma$, is a singular foliation on $ M$.
 \item Show that it contains vector fields vanishing at order $ \geq 2$ along $\Sigma$.
 \item Show that its leaves are the connected components of $ \Sigma$ and $ M \backslash \Sigma$.
 \item Extend the previous results to any sub-manifold $N$  defined as the zero locus of functions $ \phi_1, \dots,\phi_k$ which are independent in a neighborhood of any point of $N$.  
  \end{enumerate}
  In \cite{francis2023singular} and \cite{BDW}, it is explained that such a singular foliation plays the role of "vector fields tangent up to order $2$ to the submanifold $ \Sigma$". 
\end{exo}

\subsection{Hamiltonian vector fields and singular foliations}
\label{sec:hamiltonian}

Let $M$ be a smooth or holomorphic Poisson structure. We denote by $ \mathcal O_\bullet$ the corresponding sheaf of functions. Recall \cite{LPV} that a Poisson structure is a skew-symmetric biderivation $ \{\cdot, \cdot\}$ on the sheaf $\mathcal O_\bullet $ of functions that satisfies the Jacobi identity. More explicitly, for every open subset $ \mathcal U \subset M$, one is given a biderivation
 $$  \begin{array}{rcl} \mathcal O_{\mathcal U} \times \mathcal O_{\mathcal U} & \longrightarrow  & \mathcal O_{\mathcal U} \\ (f,g)&\mapsto  &  \{f,g\}\end{array} ,$$
 compatible with restrictions, and
 which satisfies the Jacobi identity for all $ f,g,h \in \mathcal O_{\mathcal U}$:
  $$ \{f,\{g,h\}\}=  \{\{f,g\},h\}+\{g,\{f,h\}\}   .$$
Since $ \{ \cdot, \cdot \}$ is a biderivation, for any function $h \in \mathcal O_\mathcal U$, 
the map 
 $$  f \mapsto \mathcal \{f,h\}   $$
 is a derivation of the sheaf $\mathcal O_\bullet $. It is therefore given by a vector field $X_h $ which is called the \emph{Hamiltonian vector field of the function $h$}.

 \vspace{0.5cm}

 \begin{propositions}{A Poisson structure induces  a singular foliation}{PoisImpliesFol}
For any Poisson structure on a manifold $M$, the $\mathcal O $-module generated by Hamiltonian vector fields form a singular foliation on $M$.

It is called the \emph{symplectic\footnote{We are not able at this point in the lecture to explain this name: it is justified by the non-trivial observation that the leaves of this singular foliation are naturally equipped with a symplectic structure, see Chapter IV in \cite{LPV}} foliation of $\{ \cdot,\cdot \}$}.
 \end{propositions}
 \begin{proof}
 Let $\mathcal F $ be the sheaf generated by vector fields of the form $ f X_h$ with $f,h$ local functions on $M$.
 It follows from the Jacobi  identity that for any local functions $ h_1,h_2 \in \mathcal O_\mathcal U$
  $ [X_{h_1}, X_{h_2}] = - X_{\{h_1,h_2\}}$ so that for any local functions $ f_1,f_2 \in \mathcal O_\mathcal U$:
   \begin{equation}
   \label{eq:stableHamiltonian}
   \left[ f_1 X_{h_1}, f_2 X_{h_2}\right] = -f_1 f_2 X_{ \{h_1,h_2 \} } + f_1  \{f_2,h_1 \} X_{h_2} - f_2  \{f_1,h_2 \} X_{h_1} . 
   \end{equation}
   This proves that $ \mathcal F$ is closed under Lie bracket. Let us prove that it is locally finitely generated. Let $ m \in M$ be a point and $ (x_1, \dots, x_d)$ be a local chart on a coordinate neighborhood $\mathcal V $. 
   Then for any function $ h(x_1, \dots, x_n ) \in \mathcal O_\mathcal V$, one has:
    $$ X_h =  \sum_{i=1}^n \frac{\partial f}{\partial x_i} X_{x_i}  .$$
    This is \emph{not} obvious, but it follows from the axioms, see Chapter I in \cite{LPV}. In particular, the family of vector fields 
     $$   X_{x_1}, \dots, X_{x_d}   $$
     generates $ \mathcal F$
 on  $\mathcal V $. 
In particular, any point has a neighborhood on which the number of generators is bounded by the dimension of the manifold.   
 \end{proof}

 \vspace{0.5cm}

  \noindent
 It is a classical result that for any Poisson structure $\{\cdot, \cdot\} $ on a manifold $M$ there exists a vector bundle morphism:
  $$  \pi^\# \colon T^* M \longrightarrow TM  $$
  which is skew-symmetric (and therefore comes from a section $ \pi \in \Gamma(\wedge^2 TM)$) such that for any two functions $f$ and $g \in \mathcal O_\mathcal U$:
   $$  \{f,g\}  = \langle \pi^\# (df) , dg \rangle $$
where  $\langle \cdot, \cdot \rangle $ stands for the duality pairing between $ TM$ and $ T^* M$.
It is also a classical result\footnote{See the classical \cite{CDW} (written in French) or the more recent \cite{CFM} for an excellent introduction to the subject.} that $T^* M $ has a Lie algebroid structure, whose bracket $ [\cdot\,, \cdot]$ is  characterized by the two following properties:
\begin{enumerate}
\item its anchor map is $ \pi^\#$,
\item on sections of $ T^*M$, i.e., exact $1$-forms, it is related to the Poisson structure by\footnote{The minus sign can be turned into a "+", it is a matter of convention.} $ [df,dg] = - d \{f,g\}$ for all open $ \mathcal U \subset M$ and $f,g \in \mathcal O_\mathcal U$.
\end{enumerate}
We call it the \emph{cotangent Lie algebroid}, see \cite{CFM}. The following result is obvious, at least if one accepts the existence of the cotangent Lie algebroid.

\begin{lemma}
The symplectic foliation in Proposition    \ref{thm:PoisImpliesFol} is the image of the anchor map of the cotangent Lie algebroid. 
In particular, $ T_m \mathcal F = {\mathrm{Im}}(\pi^\#_m)$ for all $m\in M$.
\end{lemma}

\noindent
Now, there are more singular foliations that are attached to Poisson structures.
Here are some of them. First, for any Poisson subalgebra, i.e., any sheaf of sub-algebra $ \mathcal B_\bullet \subset \mathcal O_\bullet$ such that $ \{\mathcal B_\mathcal U, \mathcal B_\mathcal U\} \subset \mathcal B_\mathcal U$ for every open $ \mathcal U \subset M$, the sheaf of $\mathcal O_\bullet $-sub-modules $\mathcal F_\mathcal B \subset \mathfrak X_\bullet $
generated the Hamiltonian vector fields of functions in $\mathcal B_\bullet $ is stable under Lie bracket in view of Equation \eqref{eq:stableHamiltonian}.
Therefore, as soon as $\mathcal F_{\mathcal B}$ is finitely generated, it becomes a singular foliation on $M$.
It happens in particular in the following two contexts.
\begin{enumerate}
\item Let $M'$ be a manifold equipped with a Poisson structure $ \{\cdot, \cdot\}'$. We say that a map $ \phi \colon M \to M'$ is a \emph{Poisson map} if the pull-back map $\phi^* \colon \mathcal O_{\mathcal U}' \to \mathcal O_{\phi^{-1}(\mathcal U)} $ is a Lie algebra morphism. 
Let $\mathcal B_\bullet := \phi^* \mathcal O_\bullet ' $ be the subsheaf of all functions pulled back from functions on $M' $. This is clearly a Poisson sub-algebra.
Moreover, $\mathcal F_\mathcal B $ is then locally finitely generated. Indeed, for any $m \in M$,  it is generated by the Hamiltonian vector fields
 $$ X_{\phi^* x_1'}, \dots, X_{\phi^* x_{d'}'}  $$
 with $x_1', \dots ,x_{d'}' $ being local coordinates in a neighborhood of $\phi(m) $ in $M'$. 
\item We say that a submanifold $ N \subset M$ is \emph{coisotropic} if the ideal $ \mathcal I_N $ of functions
 vanishing on $N$ is a Poisson subalgebra of $\mathcal O$. This is equivalent to require that 
 the Hamiltonian vector field $ X_h$ is tangent to $N$ for every $ h \in \mathcal I_N$, or to assume that for every $p \in N$, the vector bundle morphism $\pi^\# $ maps the annihilator $T_pN^\perp $ of $ T_p N \subset T_p M$ to $ T_p N$. Again, $\mathcal F_{\mathcal I_N} $ is finitely generated. Near any point not in $N$, $\mathcal F_{\mathcal I_N} $ coincides with the symplectic singular foliation, while near every point in $N$, it is generated by the vector fields $$ \left\{ X_{y_i} , y_j X_{x_c} \mid i,j=1, \dots, k,  \hbox{and } c =k+1, \dots, n \right\} $$ where $ y_1, \dots, y_k, x_{k+1}, \dots, x_n$ are local coordinates into which $ N$ is given by $ y_1=\dots =y_k =0$.
\item The singular foliation $\mathcal F_{\mathcal I_N} $ restricts to a singular foliation on $N$, which is generated by the restrictions to $N$ of the Hamiltonian vector fields $ X_{y_1}, \dots, X_{y_k}$. Equivalently, it is the singular foliation\footnote{The importance of this singular foliation comes from the following fact: when the quotient of $N$ by the leaves of this singular foliation is a manifold, it is automatically a Poisson manifold, and this procedure is called Poisson reduction, see Chapter V in \cite{LPV}} on $N$ generated 
by restrictions to $N$ of Hamiltonian vector fields of functions vanishing on $N$.  

 \end{enumerate}

\begin{exo}
Show that the Hamiltonian vector fields of a Liouville integrable system \cite{LPV} generate a singular foliation.

Show that this still holds for a non-commutative integrable system as defined in \cite{zbMATH06986134}. 
\end{exo}

\subsection{Linear singular foliations}

A faithful finite-dimensional representation of a Lie algebra may be seen as singular foliation: it suffices to consider the singular foliation associated to its transformation Lie algebroid. Let us be more precise.

Notice that for every vector space $V$ of finite dimension, there is a Lie algebra  morphism $ X \mapsto \hat{X}  $ mapping a linear endomorphism of $X \in {\mathrm{End}}(V)$ to the vector field $ \hat{X} $ on $V$ such that $ \hat{X} [\alpha] = X^* (\alpha) $ for any $\alpha \in V^* $ (seen as a function on $V$).

\begin{remark}
Upon choosing a basis $ (e_1, \dots, e_d)$ of $V$, and the corresponding coordinates $(x_1, \dots, x_d) $, this morphism maps a matrix $ (a_{i,j} )_{i=1}^d  $ to the vector field $\sum_{i,j=1}^d a_{i,j} x_i \frac{\partial}{\partial x_j} $.
\end{remark}

Let $\mathfrak g $ be a Lie algebra, and $V$ be a finite-dimensional representation of $\mathfrak g $, described by a Lie algebra morphism
$\eta\colon \mathfrak g \to   {\mathrm{End}}(V)$.
 Consider the $\mathcal O_V $-module\footnote{With $\mathcal O_V $ being smooth, holomorphic, or polynomial functions depending on whether the base field is $\mathbb R $ or $\mathbb C $, and depending on the preferences of the reader.} $\mathcal F^{\eta} $  generated by the vector fields $ \{\widehat{\eta(x)}, x \in \mathfrak g\} $.

\begin{proposition}
\label{LetonV}
Let $(V,\eta) $ be a representation of a Lie algebra $\mathfrak g $. 
Then $\mathcal F^{\eta} $ is a singular foliation on $V$.
 \end{proposition}

 The exercise supposes that the notion of leaves is already familiar to the reader. It also assumes the notion of ``isotropy Lie algebra at a point''. It explains how the initial representation can be deduced from the induced singular foliation in the faithful case.

 \begin{exo}
 Let $(\mathfrak g, V,\eta, \mathcal F^{\eta}) $ be as in Proposition \ref{LetonV}. 
 \begin{enumerate}
 \item Show that the leaves of $\mathcal F^{\eta} $  are the orbits for the Lie group action $G \to {\mathrm{GL}}(V) $ integrating $ \eta$.
 \item  This question supposes that the notion of isotropy Lie algebra at a point is known. Show that the isotropy Lie algebra of $\mathcal F^{\eta} $  at  $0 \in V $ is $\frac{\mathfrak g}{{\mathrm{ker}}(\eta)} $.

 \item Is the following statement correct: ``Two faithful representations $(V,\eta) $ and $V',\eta') $ are isomorphic if and only if their induced singular foliations $\mathcal F^{\eta} $ and $\mathcal F^{\eta'} $ are diffeomorphic''.
    \item Compare the isotropy Lie algebra of $\mathcal F^\eta $ at a point $v \in V$ with the  stabilizer of $v $.

\end{enumerate}
\end{exo}

\begin{example}\label{exa:concentric}
The singular foliation by concentric spheres, i.e., the singular foliation on $\mathbb R^n $ generated by the vector fields 
 $$  \left\{ \,  x_i \frac{\partial}{\partial x_j} - x_j\frac{\partial}{\partial x_i} \, \middle| \, 1 \leq i < j \leq n \,  \right\} $$
 comes from the action of ${\mathfrak{so}}(n) $ on $\mathbb R^n $. 
 Its leaves are by concentric spheres.

\begin{center}
\includegraphics[scale=0.4]{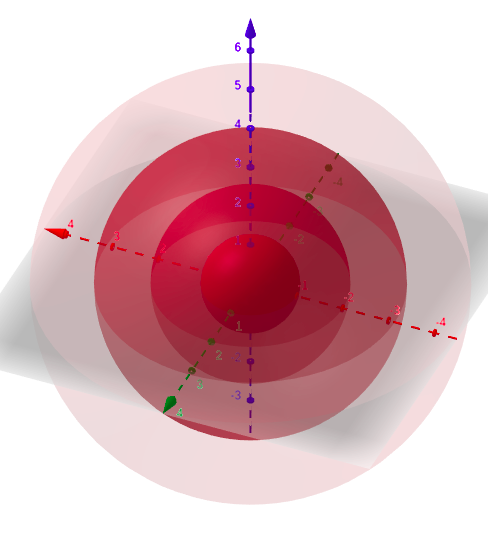}

Concentric spheres in three dimensions.
\end{center}
\end{example}

\subsection{Hörmander singular foliations}

\label{sec:hormander}

Assume that we are given, on a manifold  $M $, a family $(\mathcal G_i )_{i \in \mathbb N} $
of subsheaves of the sheaf $\mathfrak X(M) $ such that for all\footnote{Our $\mathbb N $ contains zero in the present section.} $ i,j \in \mathbb N$:
 $$ [  \mathcal G_i,\mathcal G_j ]  \subset \mathcal G_{i+j} .$$
If, moreover, each of the spaces $\mathcal G_i $ is a locally finitely generated $\mathcal C^\infty(M)$-module and if there exists $N \in \mathbb N$ such that $ \mathcal G_N =\mathcal G_i$ for all $i \geq N$, then we say that the family $(\mathcal G_i)_{i\geq 0} $ is a \emph{filtered subsheaf of $\mathfrak X(M)$}. In that case, $\mathcal{G}_0$ and $\mathcal{G}_N$ are  singular foliations on $M$.

Let us construct a singular foliation on $ \mathbb R \times M$. Denote by $ t$ the parameter on $\mathbb R$. The following Lemma is straightforward. 

\begin{lemma}
Consider a filtered subsheaf $(\mathcal G_i)_{i=0}^N$.
The subsheaf $ \mathcal G[\mathbb R] \subset  \mathfrak X_c(\mathbb R \times M ) $ of vector fields of the form\footnote{Below, $ X_i \in\mathcal G_i \subset \mathfrak X(M)$ is considered as a vector field on $ \mathbb R \times M$ whose value at $ (t,m)$ is $ (0,X_i|_{_m})$.}
 $$ \left\{ \sum_{i=0}^N g_i t^i X_i  \, \middle|  \, g_1, \dots, g_N \in \mathcal C^{\infty}(\mathbb R\times M) , X_1 \in \mathcal G_1, \dots, X_N \in \mathcal G_N \right\}   $$
 is a singular foliation on $\mathbb R \times  M $.

\end{lemma}

  We call $ \mathcal G[\mathbb R]\subset \mathfrak X(\mathbb R\times M)$ the \emph{filtered singular foliation} of $(\mathcal G_i)_{i=0}^N$.

\begin{example}
Let $\alpha $ is a contact $1$-form on a manifold $M$. Then the family:
$$
\mathcal G_i = \left\{ \begin{array}{ll}
0 & \hbox{ for $i=0$ }\\ \Gamma( {\mathrm{Ker}}(\alpha) ) &  \hbox{ for $i=1$ } \\ 
  \mathfrak X(M) & \hbox{ for $i \geq 2$ }\end{array}\right. 
$$
is a filtered subsheaf of $\mathfrak X(M) $. 
Notice that, in this case, the singular foliation  
$$ \mathcal G[\mathbb R] := \left\{f t X+ g t^2 Y \, \middle| \, X \in \Gamma( {\mathrm{Ker}}(\alpha) ), Y \in \mathfrak X(M), f,g \in \mathcal C^\infty(\mathbb R \times M) \right\}$$ 
is not of constant rank, although each one of the $ \mathcal G_i$'s is the section space of a vector bundle of constant rank. 
The leaves of this singular foliations are the points $\{(0,m)\} $ with $m \in M$ and the submanifolds $\{t\} \times M $ with $t \in \mathbb R^*$.
\end{example}

\begin{example}
\label{ex:recursion}
Consider a family $ X_1, \dots, X_r$ of vector fields. Consider the recursively defined\footnote{The brackets $ \langle \cdot \rangle$ below stand for "the module over functions generated by"} family of sub-sheaves of vector fields:
 $$\mathcal G_0=\{0\},\;  \mathcal G_1 := \langle X_1, \dots, X_r \rangle  \hbox{ and }  \mathcal G_{i+1}=  \left\langle \sum_{k=1}^i \left[ \mathcal G_k,\mathcal G_{i+1-k} \right] \, + \, \mathcal G_i \right\rangle .$$

 If the $\mathcal C^\infty(M)$  module generated by $ \mathcal G_i $ is constant after a certain rank $N$, then $ \mathcal G_N$ is a singular foliation on $M$.
 The leaves of the corresponding filtered singular foliation are  the points $(0,m) $ with $m \in M$ and the submanifolds $\{t\} \times L_m $ with $ L_m$ a leaf of the singular foliation $\mathcal G_N $ and $t \in \mathbb R$.
\end{example}

\begin{exo}
Show that the singular foliation constructed in  Example \ref{ex:recursion} can be alternatively defined by using the sequence
 $$  \mathcal G_1 := \langle X_1, \dots, X_r \rangle  \hbox{ and }   \mathcal G_{i+1}=   \left[ \mathcal G_1,\mathcal G_{i} \right] \, + \, \mathcal G_i  .$$
\end{exo}

\vspace{0.3cm}

\begin{definitions}{Hörmander's condition}{Hormander}
Let $M$ be a manifold.
A singular foliation on $\mathbb R\times M$  is said to be a Hörmander singular foliation if it is a singular foliation associated to a filtered sub-sheaves  $(\mathcal G_i)_{i \in \mathbb N} $ of $ \mathfrak X(M)$ such that $ \mathcal G_i = \mathfrak X(M)$ after a certain rank.  
\end{definitions}

\vspace{0.3cm}

For any Hörmander singular foliation such that $ \mathcal G_0=0$, the leaves are the points $\{(0,m)\} $ with $m \in M$ and the submanifolds $ \{t\} \times M$ for $ t \in \mathbb R^*$. In general, leaves are the sets $ \{0\} \times L$ with $L$ a leaf of the singular foliation generated by $\mathcal G_0 $ and the submanifolds $ \{t\} \times M$ for $ t \in \mathbb R^*$.

\begin{exo}
Show that all the isotropy Lie algebras (see Section \ref{sec:isotropy}) of a Hörmander's singular foliation such that $\mathcal G_0=0 $ are nilpotent Lie algebras.
\end{exo}

The name comes from Hörmander's condition in operator theory, whose content we briefly explain. Consider a differential operator on $\mathcal C^\infty(M) $ which is a sum of squares plus one ``linear'' term, i.e., is of the form $D=\sum_{i=1}^r X_i^2+X_{r+1} $ 
with $ X_1, \dots, X_r,X_{r+1} \in \mathfrak X(M)$. 
Then it has been shown by  Hörmander that $D$ is hypo-elliptic if the family $ X_1, \dots, X_r,X_{r+1}$ and its successive Lie brackets  (i.e., the outcome of the construction in Example \ref{ex:recursion}) generate all vector fields on $M$, see \cite{zbMATH03250869,AMY}. 
This condition is known as  Hörmander's condition, hence the name of the previously described singular foliations.

\begin{example}
The following example is given in \cite{AMY}. Let $M = \mathbb R^2 $ with parameters $x,y$ and consider $\partial_x, x\partial_y $ as being $\mathcal G_1 $ and $\mathcal G_i = \mathfrak X(\mathbb R^2)$ for $ i \geq 2$. 
The associated Hörmander singular foliation is the singular foliation on $ \mathbb R \times \mathbb R^2$ generated by the three vector fields $$ t \partial_x, t x\partial_y, t^2\partial_y .$$
Its isotropy Lie algebra (see Section \ref{sec:isotropy}) at the points $ t=x=0$ is the Heisenberg Lie algebra. It is an Abelian Lie algebra of dimension $2$ at all other singular points (i.e., the set $\{ t=0, x \neq 0 \}$). 
 \end{example}

\begin{example}
\label{ex:YErp}
The following example appears in Erik van Erp and Robert Yuncken's \cite{zbMATH06826678}. Let $M$ be a manifold of dimension $d$. Assume that the tangent bundle $TM $ of a manifold $M$ comes equipped with an increasing sequence of sub-bundles $$ 0=E_0 \subsetneq\cdots  \subsetneq E_i \subsetneq E_{i+1} \subsetneq \cdots \subsetneq E_k = TM $$
such that $\mathcal E_i = \Gamma(E_i)  $ is a filtered subsheaf of $\mathfrak X(M) $.
Then its associated singular foliation $  \mathcal E [t] $  on $ \mathbb R \times M $
is Debord. This can be seen as follows.
Let $$ 0=r_0 < r_1  < \cdots < r_i < r_{i+1} < \cdots < r_k = d $$
be the ranks of the subbundles $ (E_i)_{i=1}^k$. 
Let $e_1,\dots, e_d$ be a local trivialization on $\mathcal U \subset M $ of the tangent bundle $ TM \to M$ such that for every $i=1, \dots, k$ the family
 $ e_{1}, \cdots, e_{r_i} $ is a trivialization of $ E_i$ on $ \mathcal U$. Then the family of vector fields on $\mathbb R \times \mathcal U \subset  \mathbb R \times M $
  $$  \left(t e_1, \cdots, t e_{r_1} , t^2 e_{r_1+1}, \cdots, t^2 e_{r_2}, t^3 e_{r_2+1}, \cdots, t^{k-1}e_{r_{k-1}} ,t^k e_{r_{k-1}+1}, \cdots, t^k e_{d} \right) $$
  (where $t$ is the parameter on $\mathbb R $) generates $ \mathcal E [t] $. Since $e_1,
  \ldots, e_d $ is a trivialization of $TM$, there is of course no non-trivial relations between these generators, hence the result.

  Notice that the isotropy Lie algebras (see Section \ref{sec:isotropy}) of $\mathcal E[t] $ at $ t=0$ are nilpotent Lie algebras of dimension $d$.
 \end{example}

\subsection{Miscellaneous examples}

\begin{enumerate}
\item An important class of examples of singular foliations are the so-called \emph{Riemannian singular foliations}. Those are smooth singular foliations on a Riemannian manifold such that any geodesic orthogonal to a leaf is orthogonal to all leaf it crosses. Leaves given by actions of compact groups (for instance, concentric spheres as in Example \ref{exa:concentric}) are singular Riemannian foliations.
In \cite{nahari2022singularriemannianfoliationsmathcalipoisson}, Hadi Nahari and Thomas Strobl gave an interpretation of those in terms of Poisson structures. 
See also Oleksii Kotov and Thomas Strobl's \cite{zbMATH07119821}.
\item We said very little about dimension $1$ and codimension  $1$ singular foliations on varieties, sometimes of small dimension. There is a wide literature on the matter, see e.g., the book \cite{zbMATH07372950}.
\end{enumerate}

Let us now describe a type of singular foliation that appeared in \cite{fischer2024classification}.
Given 
\begin{enumerate}
\item $G$ be a connected\footnote{The idea of \cite{fischer2024classification} is that any singular foliation is of this type near a leaf, at least at formal level, but $G$ has to be taken infinite dimensional.} Lie group with Lie algebra $ \mathfrak g$,
\item $P \stackrel{\pi}{\to} L $ be a principal $G$-bundle,
\item and $ S$ be a manifold on which the Lie group $G$ acts,
\end{enumerate}
a singular foliation can be constructed as follows.
To start with, $G$ acts freely and properly on $ P \times S$, so that the quotient space
 \begin{equation} \label{eq:MisPS} M = \frac{P \times S}{G}  \end{equation}
 is a manifold whose elements are denoted by $ \overline{(p,s)}$ with $p \in P$ and $s \in S$, while the horizontal bar stands for the class modulo the diagonal $G$-action. 
 The map $\Pi \colon \overline{(p,s)} \mapsto \pi (p) $ is a surjective submersion, 
 that turns the manifold $M$ into a fiber bundle over $L$, with fibers diffeomorphic to $S$ (i.e., each fiber of $\Pi$ is diffeomorphic to $S$, and every point $\ell $ in $L$ has a neighborhood $ \mathcal U$ such that $ P^{-1}(\mathcal U) \simeq \mathcal U \times S$, the projection $ \Pi$ being then given by the projection on the first component). 
 Now, let us equip $M$ with a singular foliation. To start with:
 \begin{enumerate}
 \item the sheaf $ \mathfrak X(P)$ of  vector fields on $P$ 
is a singular foliation on $P$,
\item and  the  infinitesimal $ \mathfrak g$-action induces a singular foliation  $ \mathcal F_S$ on $S$ (as in Section \ref{sec:LieAlgebroidsAreSingFoliation}),
\item hence $ P \times S$ is equipped with the direct product singular foliation   (as in \ref{sec:directproduct}) that we denote by~$ \mathcal G$.
\end{enumerate}
  Since  the projection $ \Pi$ satisfies condition \emph{(ii)} in Proposition \ref{prop:iffforfibers},  there exists a unique singular foliation $ \mathcal F$ on the quotient space $M$, whose pull back singular foliation $ \Pi^{-1}( \mathcal F)$  (see Section \ref{sec:pull-back}) coincides with $\mathcal G $. Moreover, every point $ \ell$ in $L$ has a neighborhood $ \mathcal U $ such that the restriction of $ \mathcal F$ to $ \Pi^{-1}(\mathcal U) \simeq \mathcal U \times S$ is the direct product of the singular foliation of all vector field on $\mathcal U $ with the singular foliation $ \mathcal F_S$ on $S$. 

  \begin{exo}
  Let $ E \to M$ be a vector bundle or rank $N$, and $ P \to M$ be its frame bundle. 
  Let $ G $ be the Lie group of $ N \times N$ invertible matrices. Show that the previous construction applied to the singular foliation of vector fields on $S:=\mathbb R^N $ that are zero at the origin gives the singular foliation of vector fields on $E$ tangent to the zero section.
  \end{exo}

\begin{exo}
This exercise requires the notion of  leaves, and of transverse singular foliation to a leaf. For the singular foliation $ \mathcal F$  on $ M$ as in Equation \eqref{eq:MisPS} in the discussion above, establish the following points.
\begin{enumerate}
\item The classes $ \overline{(p_1,s_1)},  \overline{(p_2,s_2)}$ of two points $ (p_1,s_1),  (p_2,s_2) \in P \times S$ in $M$ belong to the same leaf of $\mathcal F $ if and only if $ s_1,s_2$ belong to the same $G$-orbit.
\item Any leaf of $ \mathcal F$ is a fiber bundle over $ L$.
\item For every fixed point $ x$ of the $G$-action,  $\Pi( P \times \{x\})$ is a leaf of $ \mathcal F$ diffeomorphic to $L$.
\item Show that for a leaf in the previous item, a representative of the transverse singular foliation is given by a neighborhood of $x$ in $S$.
\end{enumerate}
\end{exo}

We list as exercises several instances of singular foliations that do not enter any of the previous categories.

\begin{exo}
\label{exo:deg2form}
Let $ \omega$ be a closed $2$-form. Show that
$$  \left\{ X \, \in \mathfrak X (M) \, \middle| \, \mathfrak i_X \omega =0 \, \right\} $$ 
is a singular foliation on $M$, provided that it is locally finitely generated.
\end{exo}

\begin{exo}
Does the conclusion of Exercise \ref{exo:deg2form} still holds for $\omega$ a closed $n$-form?
\end{exo}

\begin{exo}
Yahya Turki \cite{zbMATH06783221} introduced the following notion: we say that a bivector field $\pi\in \Gamma(\wedge^2 TM )$ is \emph{foliated} if $\pi^\sharp (\Omega^1(M ))$ is closed under the Lie bracket, i.e., if is a singular foliation.
\begin{enumerate}
    \item Show that for any twisted Poisson structure $(\pi, \Omega) $ (also, called ``Poisson structures with background''  or ``WSW-structures'', see \cite{zbMATH02196301}-\cite{zbMATH01886781} for a definition) on a manifold $M$, $ \pi$ is a foliated bi-vector field.
    \item 
Show that in the neighborhood of a regular point of $\pi $, there exists a closed $3$-form $\Omega$ such that the pair $(\pi, \Omega) $ is a twisted Poisson structure. \item  Give an example of a foliated bivector field  $\pi$ that admits a singular point $m$ such that  there is no closed $3$-form $ \Omega$, defined on a neighborhood of $m$, such that $ (\pi,\Omega)$ is a  twisted Poisson structure. 
\end{enumerate}
(\emph{Hint:} this is done in \cite{zbMATH06783221}).
\end{exo}

\section{New constructions from old ones}
\label{sec:newforold}

In the present section, we work indifferently in the context of smooth, complex or real analytic geometry. Most arguments presented here, however, make no sense in algebraic geometry, and have to be adapted. Conversely, some of them only make sense in algebraic geometry. We will be more precise in due time.

Here is a first exercise to train on these notions.

\begin{exo}
Let $\mathcal F $ be a smooth singular foliation on $M$ and $\varphi \in \mathcal O_M $ be a function. Show that
$$  \varphi \mathcal F := \{ \varphi X, X \in \mathcal F\} $$
is a singular foliation again. 
State and show the corresponding result in the real analytic, complex and algebraic settings.
\end{exo}

\subsection{Direct products of singular foliations}
\label{sec:directproduct}

For $X_1, X_2$ vector fields on $M_1,M_2$ respectively, we shall denote by  $(X_1,X_2)$ the vector field on $ M_1 \times M_2$ whose value at $(m_1,m_2) \in M_1 \times M_2$ is $ (X_1|_{m_1} , X_2|_{m_2} ) \simeq T_{(m_1,m_2)} M_1 \times M_2$.

For $(M_1, \mathcal F_1) $
and 
$(M_2, \mathcal F_2)  $
foliated manifolds, the product manifold $M_1 \times M_2 $ can be equipped with the direct product of both foliations. 

\vspace{.5cm}

\begin{definitions}{Direct product of singular foliations}{ref:directproduct}
The \emph{direct product} of two foliated manifolds $(M_1,\mathcal F_1)$ and $(M_2,\mathcal F_2$  is the singular foliation $\mathcal F_1 \times \mathcal F_2 $ on $ M_1 \times M_2$ such that, for every open subset $\mathcal U_1 \subset M_1, \mathcal U_2 \subset M_2$, $\mathcal F_1 \times \mathcal F_2 $ is the $\mathcal O_{\mathcal U_1 \times \mathcal U_2} $-module  generated by vector fields of the form $(X_1,X_2)$ with $ X_1 \in \mathcal F_1$ and $ X_2 \in \mathcal F_2$.

It is denoted by  $(M_1 \times M_2 ,\mathcal F_1 \times \mathcal F_2)$. 
\end{definitions}

\vspace{.5cm}

 \begin{exo}
 \label{exo:finitelygenerateddirectproduct}
 Show that the direct product of finitely generated singular foliations is a finitely generated singular foliation. Compare their ranks.
 \end{exo}

\vspace{.3cm}

\subsection{Pull-back (through surjective submersions)}

\label{sec:pull-back}

Let us give the easiest version of the pull-back of a singular foliation: the pull-back through surjective submersion. We will come back to this notion later on, using a more general definition due to Androulidakis and Skandalis \cite{AS}.

\vspace{.5cm}
We work indifferently in the smooth, holomophic or real analytic settings. We have two manifolds $P$ and $M$, with respective sheaves of sections $\mathcal O^P $ and $ \mathcal O^M$. Given $\psi \colon P \to M$ a map in the relevant category.

\begin{definition}\label{def:related}
 A vector field $X\in\mathfrak{X}(P)$ said to be \emph{$\psi$-related} to a vector field $\widetilde{X}$ on $M$ if for all $p \in P$,
\begin{equation}
  (T \psi)_p (X_{|_p} )= \widetilde{X}_{|_{\psi(p)}}  
\end{equation}
Equivalently, for any  $f\in \mathcal O^M$, 
$  X[f\circ \psi ] =\widetilde{X}[f]\circ \psi$, or, equivalently, such that the following diagram commutes:  

$$\xymatrix{ \mathcal O^M\ar[r]^{\tilde X}\ar[d]_{\psi^*}  &\ar[d]^{\psi^*}  \mathcal O^P\\ \mathcal O^M \ar[r]_{X} & \mathcal O^M }$$
\end{definition}

 A vector field on $M$ is said to be a \emph{vertical vector field} if it takes values in $\ker (\dd_p \psi) \subset T_p P M $ for all $p \in M$. Equivalently,  vertical vector fields are vector fields $\psi$-related to $0 \in \mathfrak X(M)$.

\vspace{0.5cm}

\begin{definitions}{Pull-back of singular foliations}{def:PullBackSubmersion}
Let $\mathcal F$ be a singular foliation on a manifold $M $ and let $\psi: P \to M$ be a surjective submersion. We call \emph{pull-pack of $\mathcal F $ by $\psi$} and denote by $\psi^{-1}(\mathcal F) $ the singular foliation on $P$  generated, as an $ \mathcal O^P$-module, by vector fields $ \psi$-related to a vector field in $\mathcal F$.
\end{definitions}

\vspace{0.5cm}
The definition needs to be justified. First, one has to check that this definition indeed gives a sheaf of $ \mathcal O^P$-modules on . We then have to check that all three conditions $(\alpha),(\beta),(\gamma) $ in Definition \ref{def:consensus2} or \ref{def:consensus2alg} are satisfied.
\begin{enumerate}
\item[$( \alpha)$]  is an obvious consequence of the fact that if $ Y_1,Y_2$ are $ \psi$-related to $ X_1,X_2$ respectively, then $ [y_1,Y_2]$ is $ \psi$-related to $ [X_1,X_2]$.
\item[ $ (\beta)$] holds by definition.
\item[ $ (\gamma)$]  holds since, as for every submersion, every point $ p \in P$ admits a neighborhood on which isomorphic to the product of two open balls $\mathcal U_1,\mathcal U_2 $ such that $ \psi$ is the projection onto $ \mathcal U_1$, which is seen as an open subset of $M$. Under this isomorphism, $ \psi^{-1}(\mathcal F)_\mathcal U$ is the direct product of $ \mathcal F_{\mathcal U_1}$ with $ \mathfrak X_{\mathcal U_2}$. For $ \mathcal U_1,\mathcal U_2$ small enough, it is therefore finitely generated (see Exercise \ref{exo:finitelygenerateddirectproduct}).  
\end{enumerate}

\begin{remark}
In particular, all vertical vector fields, i.e., vector fields tangent to the fibers of $\psi $, are contained in $ \psi^{-1}(\mathcal F)$.
\end{remark}

\begin{remark}
In the smooth case, one can equivalently define singular foliations using compactly supported vector fields as in Definition \ref{def:consensus}. 
Assume $ \mathcal F_c$ is a singular foliation on $M$ as in Definition \ref{def:consensus}.
Defining pull-back causes then a technical difficulty: if the fibers of $\psi \colon P \to M$ are not connected, then
the pull-back singular foliation $\psi^{-1}(\mathcal F_c)$ can not be defined as being the $\mathcal C^\infty(P) $-module generated by compactly supported vector fields on $P$ which are $\psi$-related to a vector field in $\mathcal F_c$.   Indeed, if the fibers of $ \psi$ are not compact, there is no compactly supported vector field on $P$ which is $\psi$-related to a non-zero vector field on $M$.
The singular foliation $\psi^{-1}(\mathcal F) $ has to be then defined as the $ \mathcal C^\infty_c(P)$-module generated by vector fields on $P$ which are $ \psi$-compatible to a vector field in $ \mathcal F$.
In equation, if one denotes by $\mathfrak X (P)_\psi$ the Lie algebra of vector fields on $P$ which are $ \psi$-compatible to a vector field on $M$, and $\psi_* \colon \mathfrak X (P)_\psi \to \mathfrak X(M) $ the natural Lie algebra morphism, we have: 
$$ \psi^{-1}(\mathcal F_c)  := \mathcal C^\infty_c(P) \, (\psi_*)^{-1} ( \mathcal F_c) .$$
\end{remark}

\begin{exo}
Let us work in the setting of smooth differential geometry.
A \emph{horizontal distribution on the surjective submersion $\psi \colon P \to N $}, is a regular distribution $ p \mapsto \mathcal H_p $ on $P$ such that\footnote{Those are also called Ehresmann connection. They exist for any fiber bundles.} 
 $$ \mathcal H_p  \oplus {\mathrm{ker}}(T_p \psi) = T_p P \hbox{ for all $ p \in P$ } .$$
 We call \emph{horizontal lift } of $ X \in \mathfrak X(M) $ and denote by $ \mathcal H(X)$ the unique section of $\mathcal H $ such that $T\psi ( \mathcal H(X)_{|_p}) =X_{|_{\psi(p)}}$ for all $ p \in P$. Show that
 $ \psi^{-1}(\mathcal F)$ is generated, as a sheaf of $\mathcal C^\infty(P)$-modules, by horizontal lifts of vector fields in  $ \mathcal F$ and vertical vector fields (= vector fields tangent to the fibers of $\psi$).
\end{exo}

\begin{exo}
\label{exo:discretegroup}
Let $G$ be a discrete group acting freely and properly on a smooth manifold $ M$ by smooth diffeomorphism.
Recall that this implies that $ M/G$ is a manifold and that the natural projection $ \Pi \colon M \to M/G$ is a surjective local diffeomorphism.
Let $ \mathcal F$ be a singular foliation on $M$. 
\begin{enumerate}
\item 
Show that if $G$ acts on symmetries of $\mathcal F $, then  there exists a singular foliation $ \mathcal F^G$ on $ M/G$ whose pull-back to $M$ is $ \mathcal F$. 
\item Is the converse true?
\end{enumerate}
\end{exo}

\begin{exo}
Show that the $p$-vertical vector fields are contained in $p^{-1}( \mathcal F_B)$.
\end{exo}

\begin{exo}
Let $ (M,\mathcal F)$ be a foliated manifold.
Let $ \psi \colon P \to M$ be a surjective submersion
\begin{enumerate}
\item Show that the leaves of $\psi^{-1}(\mathcal F)$ are the connected components of the inverse images through $\psi$ of the leaves of $\mathcal F $.
\item 
Show that the isotropy Lie algebra of   $p^{-1}(\mathcal F_B) $ at a point $m$ is canonically isomorphic to the isotropy Lie algebra of $\mathcal F_B $ at $p(m) $.
\item 
Show that the transverse singular foliation of $ \psi^{-1}(\mathcal F)) $ of a leaf through a point $p \in P$ is canonically isomorphic to the transverse singular foliation of $ \mathcal F$ at the leaf through $\psi (p )$.
\end{enumerate}
\end{exo}

We conclude this section with a statement which we insist is not an obvious statement, for it will use the assumption ``locally finitely generated''. It is wrong  for general involutive distributions.
For instance, for the ``infinite comb'' of Exercise \ref{exo:infinitecomb}, the projection $ (x,y) \mapsto y$ onto the horizontal axis satisfies Conditions \emph{(ii)} and \emph{(iii)} but does not satisfy condition \emph{(i)}.

\begin{proposition} \label{prop:iffforfibers}
Let $\psi \colon P \to M$ be a surjective submersion with connected fibers, and $\mathcal F_P$ a singular foliation on $P$. Then the following are equivalent:
\begin{enumerate}
    \item[(i)] There  exists a singular foliation $ \mathcal F_M$ on $M$ such that $\mathcal F_P = \psi^{-1}( \mathcal F_M) $.
    \item[(ii)] Each fiber of $\psi$ is contained in a leaf of $L $.
    \item[(iii)] For every $m \in M $,  we have ${\ker}(T_m \psi) \subset T_m \mathcal F_P  $.
\end{enumerate}
\end{proposition}

\subsection{Restriction of a singular foliation to a transverse submanifold}

\label{sec:transverse}

Let $\mathcal F $ be a singular foliation on a smooth manifold $M$, and let $S \subset M$ be a sub-manifold.
We would like to restrict the singular foliation $\mathcal F $ to $S$.

The next exercise presents a naive idea – which works, but has to be made more precise.

\begin{exo}
\label{exo:restriction}
Let $ \mathcal F$ be a singular foliation on a singular foliation $M$, seen as a subspace of $ \mathfrak X_c(M)$ as in Definition \ref{def:consensus}.
Let $S \subset M $ be a closed embedded submanifold.
Consider $\mathfrak i_S^* \mathcal F_{naive} \subset \mathfrak X_c(S) $ to be the sub-space of all vector fields on $S $  obtained by restricting to $S$ vector fields in $\mathcal F $ that happen to be tangent to $S$. Show that $\mathfrak i_S^* \mathcal F_{naive} $
\begin{enumerate}
    \item is a sub-Lie algebra of $ \mathfrak X(S)$,
    \item  is a sub-$\mathcal C^\infty(S) $-module of $ \mathfrak X(S)$,
\end{enumerate}
It is therefore a Lie-Rinehart subalgebra of vector fields on $S$ (i.e., it satisfies $ (\alpha)$ and $ (\beta)$ in Definition \ref{def:consensus}). However,

\begin{enumerate}
    \item [3.] Show that if the submanifold  $S$ is embedded but not closed in $M$, the restriction to $S$ of a vector field compactly supported on $M$ may even not be compactly supported on $S$.
    (\emph{Hint}: take $  M =\mathbb R^2,\mathcal F= \mathfrak X_c(M)$ and $ \Sigma$ the spiral obtained as an  integral curve of $X=x \partial_y - y \partial x - x \partial _x - y \partial_y $. Then multiply $X$ by  a compactly supported function which  is $1$ at $ (0,0)$.). 
    \item[4.] Show that if $S$ is only immersed and not embedded,  $\mathfrak i_S^* \mathcal F_{naive}$ many even not be a $\mathcal C^\infty(S)$-module.
    {\emph{Hint}:} the previous counter-example will do as well…
\end{enumerate}
In order to have an induced singular foliation on $S$, we have to be more sophisticated. Let $S \subset M$ be an immersed 
submanifold of $M$: we now work in the smooth, real analytic or complex settings altogether, and we consider that $ \mathcal F$ is a sheaf as in Definition \ref{def:consensus2}. We
  denote by $\mathfrak i  \colon S  \hookrightarrow M $ the canonical inclusion. We define a sheaf $\mathfrak i_S^* \mathcal F  \subset \mathfrak X(S)$ as follows. To every $\mathcal U \subset S$, we associate the space of all vector fields $Y \in \mathfrak X(S)_{\mathcal U}$ such that for every $s \in \mathcal U$, there exists $X \in \mathcal F_{\mathcal W}  $ (for some open subset $\mathcal W \subset M $ containing $ \mathfrak i(s)$) such that 
 $$ T_{s'}\mathfrak i(Y_{|_{s'}}) = X_{|_{\mathfrak i(s')}}  $$
 for every $s' $ is a neighborhood of $s$ in $ S$.
 Check that the previously defined object:
\begin{enumerate}
    \item[5.] is a sub-sheaf of the sheaf of vector fields on $S$,
    \item[6.] is closed under Lie bracket,
    \item[7.] and is a module over the relevant sheaf of functions on $S$.
\end{enumerate}

Consider the foliation of $\mathbb R^2$ by horizontal lines, i.e., $\mathcal F$ is generated by $\frac{\partial}{\partial x}$. Let $f$ be a function which has support $[0,\infty)$. Then the graph of $f$, namely,
$$S = \{(x,f(x) | x \in \mathbb R\} $$ is an embedded submanifold of $M$.
\begin{enumerate}
\item[8.] Show that, in this case, $\mathfrak i_S^* \mathcal F $ is not locally finitely generated. (\emph{Hint:} $\mathfrak i^*_S \mathcal F$ is exactly the space of vector fields which are supported in $(-\infty,0]\subset \mathbb R$.) 
\end{enumerate}
\end{exo}

However,  there is a situation where the sheaf $ \mathfrak i_S^* \mathcal F $ defined in Exercise \ref{exo:restriction} is a locally finitely generated module, and is therefore an induced singular foliation on $S$. 

\begin{definition}
\label{def:beingtransverse}
We say that a submanifold $S$ of a foliated manifold $(M,\mathcal F )$  \emph{intersects cleanly} $\mathcal F $ if $T_s S + T_s\mathcal F =T_s M $ for all $s \in S $. We also say that $S$ is  \emph{transverse to $\mathcal{F}$},
or is a transverse sub-manifold. 
\end{definition}

The condition about a clean intersection is enough to  guarantee that $\mathfrak i^*_S \mathcal F $  is locally finitely generated.
\vspace{0.5cm}

\begin{propositions}{Submanifolds intersecting $ \mathcal F$ cleanly}{prop:transverse}
Let $S \subset M$ be a submanifold that intersects cleanly\footnote{Also called transverse. See Definition \ref{def:beingtransverse}.} a smooth singular foliation $\mathcal F$. Then $\mathfrak i^*_S \mathcal F $  is a singular foliation on $S$.

It is called the \emph{restriction of the singular foliation to $S$}.
\end{propositions}
\vspace{0.5cm}

The following exercises describe this structure more precisely.

\begin{exo}Let $S$ be a submanifold that cleanly intersects $ (M,\mathcal F)$, and let $\mathfrak i^*_S \mathcal F $ be its induced singular foliation.
\begin{enumerate}
    \item Show that the rank of $ \mathfrak i^*_S \mathcal F$ at a point $s$ is $ {\mathrm{rk}}_s \mathcal F - {\mathrm{codim}}(S)$.
    \item Show that $T_s \mathfrak i^*_S \mathcal F = T_s \mathcal F \cap T_s S $ for all $s \in S$.
    \item (Supposes that the notion of leaf is known, see Section \ref{sec:leavesExist}.) Show that the leaf of $\mathfrak i^*_S \mathcal F$ through a point $s \in S$ is the connected component containing $s$ of $S$ with the leaf through $s$ of $\mathcal F $. 
    \item (Supposes that the notion of isotropy Lie algebra is known, see Section \ref{sec:isotropy}.) Show that the isotropy Lie algebra of $\mathcal F $ and $\mathfrak i^*_S \mathcal F $ coincide  at any point $s \in S$.
\end{enumerate}
\end{exo}

\begin{exo}
The goal of this exercise is to show that there is a neighborhood of a transverse submanifold $S$ in a foliated manifold $(M,\mathcal F)$ on which $\mathcal F $ coincides with a neighborhood of the zero section in the normal bundle $N_S := TL/TS \stackrel{p}{\to} S $, equipped with the pull-pack singular foliation $ p^* \mathfrak i^*_S \mathcal F$. We work in the smooth setting (it is not true in the holomorphic or real analytic settings).
\begin{enumerate}
    \item Show the ``tubular neighborhood theorem'', i.e., that there is a neighborhood $\mathcal U $ of $S$ in $L$ diffeomorphic to a neighborhood $ \mathcal U$ of the zero section in the normal bundle $N_S := TL/TS \stackrel{p}{\to} S $, through a diffeomorphism which is the identity on $S$.
    \item Show that the tubular neighborhood $ \mathcal U$ in the previous item can be chosen such that  vector fields tangent to the fibers of $p \colon  \mathcal U  \to S $ are included in $\mathcal F $.
    \item Conclude that the restriction to $ \mathcal U$ of the singular foliation 
 $\mathcal F $ is isomorphic to the pull-back singular foliation $ p^{-1}( \mathfrak i^*_S \mathcal F)$ (\emph{Hint:} $ \mathcal U$ 
 see Section \ref{sec:pull-back}, and Proposition \ref{prop:iffforfibers}.)
\end{enumerate}
\end{exo}

\begin{exo}
This exercise requires the notion of anchored bundle over a singular foliation, see Section \ref{sec:AnchoredBundle}.
Let $(A, \rho)$ be an anchored bundle over a singular foliation $\mathcal F $. Let $S \subset M$ be a submanifold that intersects $ \mathcal F$ cleanly. Show that $\rho^{-1}(TS) \subset \mathfrak i_S^* A $ is an anchored bundle over $\mathcal F_S $, when equipped with the restriction of the anchor map.
(Here, $\mathfrak i_S^* A$ stands for the restriction of the vector bundle $A$ to $S$.)
\end{exo}

\begin{exo}
\label{exo:restrictionAndSymmetry}
This section requires the notion of symmetry of a singular foliation (see Section \ref{sec:symmetry}).
Let $ \Phi\colon M \to M$ be a symmetry of a singular foliation $ \mathcal F$ on a manifold $M$.
Let $S \subset M$ be a submanifold. 
\begin{enumerate}
\item Show that $S $ intersects $ \mathcal F$ cleanly if and only if $ \Phi(S)$ does.
\item Show that, in that case, the restriction of $\Phi$  to $S$ is an isomorphism of foliated manifolds $(S,\mathfrak i_S^* \mathcal F) \longrightarrow (\Phi(S),\mathfrak i_{\Phi(S)}^* \mathcal F) $.
\end{enumerate}
\end{exo}

\subsection{Pull-back of singular foliations (beyond immersions and submersions)}

\label{sec:pullback2}

We have already defined pull-back through surjective submersions, but also the restriction to some submanifolds.
Let us unify these constructions, following an idea of Androulidakis and Skandalis \cite{AS}. 

In this section, we restrict ourselves to the case of smooth manifolds: the complex or real analytic cases are similar.

Let $L,M$ be manifolds together with a smooth map $\phi\colon L\to B$. Let $\phi^*TM $ be the pull-back through $\phi$ of the tangent bundle $TM$:
 $$ \phi^*TM:=\{(\ell,u) \in L \times TM \mid u \in T_{\phi(\ell)} M\} .$$
 
There are two natural maps:
$$ \xymatrix{ \mathfrak X (L) \ar[dr]^{T\phi} &  & \mathfrak X (M) \ar[dl]_{\phi^*} \\ & \Gamma(\phi^* TM) & } $$
defined as follows.
\begin{enumerate}
    \item 
Any vector field $X$ on $M$ gives a section $\phi^* X $ of $\phi^* TM $ defined by
 $$ \ell \mapsto (\ell,X_{\phi(\ell)}). $$
 called the pull-back of $X$.
\item 
There is a natural vector bundle morphism defined for all $ u \in T_\ell L$ by
 $$  \begin{array}{llll}T\phi :& TL &\to& \phi^* TM \\  & u &\mapsto & (\ell, T_\ell \phi \, (u))\end{array} $$
  At the level of sections, it induces a map $\mathfrak X(L) \to \Gamma(\phi^*TM)  $. 
\end{enumerate}

Let $\mathcal F_M $ be a singular foliation on $M$, seen as in Definition \ref{def:consensus} as a subspace of compactly supported vector fields.
We denote by  $\phi^* \mathcal F_M  $ the $ \mathcal C^\infty_c(L)$-submodule\footnote{Recall that the index $ {}_c$ means "compactly supported".} of $\Gamma_c(\phi^*TM)$  generated by $\{\phi^* X | X \in \mathcal F_M\}\subset \Gamma(\phi^*TM)$. 
We now present the construction of a singular foliation on $L$ as an exercise.

 \begin{exo}
 Consider the submodule of $\mathfrak X_c(L)$ defined by\footnote{Equivalently, a vector field  $X\in \mathfrak X_c (L)$ belongs to $\phi^{-1}(\mathcal{F}_M)$ if and only if  there exists smooth functions $g_1, \ldots, g_k\in \mathcal C_c^\infty(L)$ and  $Y_1,\ldots, Y_k$  in $\mathcal{F}_M$ such that $T_\ell \phi(X|_\ell)=\sum_{i=1}^k g_i(\ell) \, Y_i|_{\phi(\ell)}$ for all $ \ell \in L$. }
$$\phi^{-1}(\mathcal F_M):=\left\lbrace X\in\mathfrak{X}_c(L)\mid T\phi(X)\in \phi^*\mathcal{F}_M\right\rbrace .$$ 
\begin{enumerate}
\item 
Show that $\phi^{-1}(\mathcal{F}_M)$ is involutive, i.e., closed under Lie bracket.
\item  We say that $ \phi \colon L \to M$ is \emph{transverse to $ \mathcal F_M$} if for all $\ell \in L$, we have $ T_{\phi(\ell)} \mathcal F_M + T_\ell \phi \, (T_\ell L) = T_{\phi(\ell )}M$. Show that for $ \phi$ an immersion, this definition matches the transversality condition given in Equation \eqref{eq:transverse}. Show that submersions are  transverse to any singular foliation on $ M$.
\item Show that if $\phi $ is transverse to $\mathcal F_M $, then  $\phi^{-1}(\mathcal F_M)$ is a singular foliation on $M$. 
\end{enumerate}
\end{exo}

This exercise justifies the following definition.

\vspace{.5cm}

\begin{definitions}{Pull-back w.r.t. a transverse map}{transversemaps}
Let $(M,\mathcal F_M)$ be a foliated manifold, and $ \phi: L\to M$  a smooth map
 transverse to $\mathcal F_M$. We call the singular foliation $\phi^{-1}(\mathcal F_M)$ the \emph{pull-back of $\mathcal{F}_M$ through $\phi$}. 
\end{definitions}

\vspace{.5cm}

\begin{exo}
Explain why this notion ``unifies'' (= i.e., admits as particular cases) pull-back with respect to surjective submersions seen in Section \ref{sec:pull-back} and 
restrictions to transverse submanifolds seen in Section \ref{sec:transverse}.
\end{exo}

\subsection{The suspension of a singular foliation}
\label{sec:suspense}

We work here in the setting of smooth differential geometry. Most results could be adapted to real analytic or complex settings.

In this section, we denote by $ \mathfrak X(N)$ the singular foliation of all vector fields on a manifold $N$.
 
\subsubsection{Suspension in dimension $1$}

We call \emph{suspension} of the manifold $M$ with respect to a diffeomorphism $ \phi \colon M \to M$ the quotient of $M \times \mathbb R $ by the action of the additive group $\mathbb Z$ by: 
 \begin{equation}  \label{eq:Zaction} k \cdot (m,t) \, := \, (\phi^k (m), t+k) \end{equation}
 for all $k \in \mathbb Z, m \in M, t \in \mathbb R $. 
 Since the action of $\mathbb Z $ is free and proper, the quotient is a manifold - that we call \emph{suspension of $M$ by $\phi $} and denote by  $ M_\phi := \frac{M \times \mathbb R}{{\mathbb Z}_\phi}$.
 We also denote by $ \pi \colon M \times \mathbb R \to \frac{M \times \mathbb R}{{\mathbb Z}_\phi} $ the natural projection.
 
Let us assume now that $M$ comes equipped with a singular foliation $\mathcal F $ (seen as a sheaf as in Definition \ref{def:consensus2}) and that $\phi\colon M \to M $ is a symmetry of $\mathcal F $, i.e., that $\phi_* (\mathcal F) =\mathcal F $. 
Then $ M \times \mathbb R$ comes with the direct product singular foliation $(M \times \mathbb R, \mathcal F \times \mathfrak X(\mathbb R)) $ (where $\mathfrak X(\mathbb R)$ stands for the sheaf of all vector fields on $ \mathbb R$). 
For all $k \in \mathbb Z   $:
 $$  (m,t) \mapsto (\phi^k (m), t+k) $$
 is a symmetry of 
  the direct product singular foliation $(M,\mathcal F) \times (\mathbb R,\mathfrak X(\mathbb R)) $. By exercise \ref{exo:discretegroup},  the following results hold true.

\begin{proposition}
\label{prop:suspension1}
Let $\phi: M \to M$ be a symmetry for a singular foliation $\mathcal F $. There exists an unique singular foliation on the suspension $  M_\phi := \frac{M \times \mathbb R}{\mathbb Z_\phi} $ of $M$ by $ \phi$ whose pull-back on $ M \times \mathbb R  $ is the direct product singular foliation $(M \times \mathbb R, \mathcal F \times \mathfrak X(\mathbb R)) $.  
\end{proposition}

We call the  singular foliation in Proposition \ref{prop:suspension1} the \emph{suspension of $ \mathcal F$ by the symmetry $\phi $} and denote it by $
\mathcal F_\phi$. Before describing this singular foliation in more details, let
 us recall a classical result of differential geometry about suspensions of diffeomorphisms:

\begin{lemma}
\label{lem:suspension}
If a diffeomorphism $\phi$ of a manifold $M$ is the time $1$ flow of a complete vector field $X \in \mathfrak X(M)$, then the suspension $ M_\phi := \frac{M \times \mathbb R}{\mathbb Z_\phi}$ of $M$ by $ \phi$ is diffeomorphic to the direct product\footnote{i.e., the suspension $\frac{M \times \mathbb R}{\mathbb Z_{id_M}}$ associated to the identity map of $M$.} $M \times S^1 $. 
\end{lemma}
\begin{proof}
The vector field on $ M \times \mathbb R$ whose value at $ (m,t)$ is $ t X_{|_m}$ has a flow at time $t $ that intertwines the $\mathbb Z $-action as in Equation (\ref{eq:Zaction}) with the $ \mathbb Z$-action:
 $$  k \cdot (m,t)= (m,t+k) .$$
 Since the quotient of $ M \times \mathbb R$ through this action is $ M \times S^1$, this completes easily the proof.
\end{proof}

When $ X$ belongs to $ \mathcal F$, then the vector field that appears in the proof of
Lemma \ref{lem:suspension} belongs to $ $. 
Hence (by the highly non-trivial Corollary \ref{coro:flowissymmetry} that will be proven later on),  the next statement holds true.
\vspace{0.5cm}

\begin{propositions}{Inner symmetries have trivial suspensions}{Prop:nosuspension} 
If a symmetry $\phi $ of a singular foliation $\mathcal F$ on $M $ is the time $1$ flow of a complete vector field\footnote{With a little more work, this result can be extended to the case where $ \phi$ is an inner symmetry of $ \mathcal F$.} in $\mathcal F $, then its suspension $ (M_\phi,\mathcal F_\phi)$ is isomorphic to the direct product singular foliation 
$(M \times S^1,\mathcal F \times \mathfrak X(S^1))$.
\end{propositions}
\vspace{0.5cm}

\begin{example}
\label{ex:snake}
An important example of such a singular foliation is the so-called \emph{self-eating snake}, which is defined by:
\begin{enumerate}
\item $M= \mathbb R^2$ with coordinates $ (x,y)$,
\item $ \mathcal F$ is the singular foliation defined by the vector field $ x \partial_y - y \partial_x$ (whose leaves are concentric circles),
\item $ \phi \colon M \to M$ is the division by $2$, namely $ \phi (x,y) = (x/2,y/2).$
\end{enumerate}
The leaves of the suspension of $ (M,\mathcal F)$ by $ \phi$ are then as follows. To start with, $ M_\phi = M \times S^1$, and the singular foliation is as follows: 
\begin{enumerate}
\item  the circle $ \{0\} \times S^1$ is a leaf,
\item  all the other leaves are diffeomorphic to a cylinder. These cylinders wrap around this circle.
\end{enumerate}

\begin{center}
\includegraphics[width=5cm]{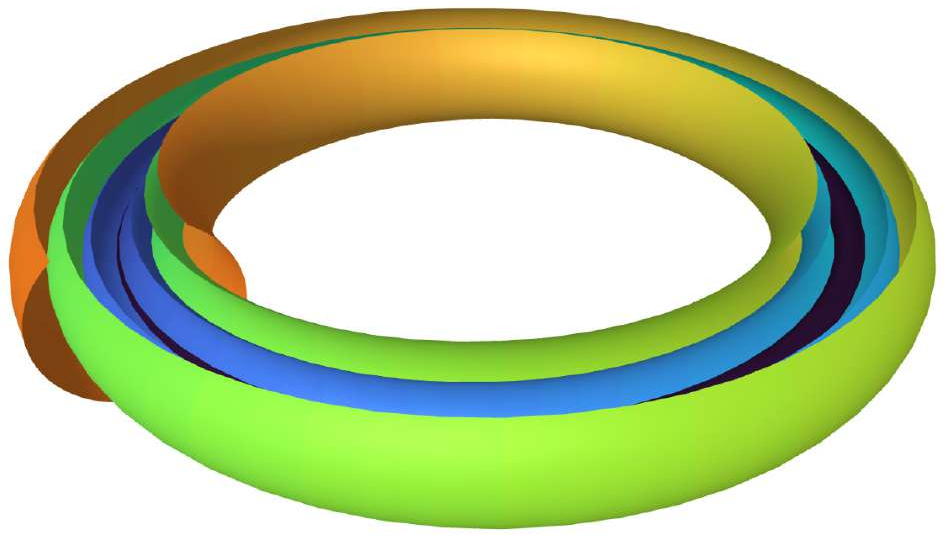}
\end{center}

\end{example}

\subsubsection{Suspension in dimension $\geq 2$}
\noindent
The suspension can be defined in a more general context. 
 Let $ \Sigma$ be any connected manifold, $ \sigma \in \Sigma$ a point. 
 Let $ \pi_1(\Sigma,\sigma)$ be the fundamental group\footnote{i.e.,  $ \pi_1(\Sigma,\sigma)$ is the group of homotopy classes of loops based at $ \sigma$} based at $ \sigma$. Let $\tilde{\Sigma}_\sigma $ be the universal cover of $ \Sigma$, again computed\footnote{i.e., $\tilde{\Sigma}_\sigma $ is the set of all homotopy classes $ [\gamma]$ of paths $\gamma \colon [0,1]  \to \Sigma$ with $\gamma(0)= \sigma$.} with respect to $ \sigma$. The fibers of the natural projection
 $$  \begin{array}{rcl}\tilde{\Sigma}_\sigma&  \to &\Sigma \\ \left[\gamma\right]& \mapsto &\gamma(1) \end{array} $$
are equipped with a natural $ \pi_1(\Sigma, \sigma)$-action by 
 $$  [g] \cdot [\gamma] := [ g \star \gamma]$$
 where $\star$ refers to concatenation of paths.
This action turns $\tilde{\Sigma}_\sigma$  it into a $ \pi_1(\Sigma,\sigma)$-principal bundle. 

Let $M$ be a manifold. For any group morphism:
 $$ \varphi \colon \pi_1(\Sigma,\sigma)  \longrightarrow \mathrm{Diff}(M) ,$$
we call \emph{suspension of $ \varphi$} the quotient manifold: 
 $$  \frac{M \times \tilde \Sigma_\sigma}{\pi_1(L,\ell)_\varphi}   $$
where the quotient is with respect to the diagonal action  
 \begin{equation}\label{eq:diagonalAction} [g] \cdot (m,[\gamma])  =  (\psi([g]) (m) , [g \star \gamma] )\end{equation} for all $ [g] \in \pi_1(\Sigma,\sigma)$, $ m \in M$, $ [\gamma] \in \tilde \Sigma_\sigma$.
 Let $ \mathcal F$ be a singular foliation on $M$, and assume that $ \varphi$ is valued in symmetries of $ \mathcal F$. Consider the direct product singular foliation 
\begin{equation}   \label{eq:directproduct}
(M,\mathcal F )\times ( \Sigma , \mathfrak X(\tilde \Sigma_\sigma)) .
\end{equation}
The group action defined in Equation \eqref{eq:diagonalAction} is valued in symmetries of the direct product singular foliation \eqref{eq:directproduct}. 
In particular, the singular foliation \eqref{eq:directproduct} descends to the quotient to define a singular foliation $  \frac{\mathcal F\times \mathfrak X (\tilde \Sigma_\sigma)}{\pi_1(L,\ell)_\varphi}   $  on 
$   \frac{M \times \tilde \Sigma_\sigma}{\pi_1(L,\ell)_\varphi}  $.
We call \emph{suspension of the singular foliation $(M,\mathcal F) $ with respect to $\varphi \colon \pi_1 (\Sigma,\sigma) \to \mathrm{Sym}(M,\mathcal F)$} the pair
$$  \left(  \frac{M \times \tilde \Sigma_\sigma}{\pi_1(\Sigma,\sigma)_\varphi}
, 
\frac{\mathcal F\times \mathfrak X (\tilde \Sigma_\sigma)}{\pi_1(\Sigma,\sigma)_\varphi} \right) 
$$

\begin{example}
    For $ \Sigma = S^1$ one recovers the previous construction. 
  
\end{example}

    \begin{remark}
It is tempting to generalize Proposition \ref{thm:Prop:nosuspension}   to an arbitrary $ \Sigma$ by stating that if
$ \varphi$ is valued in inner symmetries of $ \mathcal F$, then the suspension should be isomorphic to the direct product of $(M,\mathcal F) \times ( \Sigma, \mathfrak X (\Sigma) )$. But such a statement is wrong, see the  discussion about torus in
\cite{fischer2024classification}.
 \end{remark}
 
\begin{exo}
\label{exo:leafSigma}
This exercise requires the notion of leaves.
Show that if $m \in M$ is a point where $ T_m \mathcal F=\{0\}$, and if $\varphi ([g])(m) =m $ for every $[g] \in \pi_1(\Sigma, \sigma) $, then the suspension of $ (M,\mathcal F)$ with respect to $ \varphi$ has a leaf $L$ diffeomorphic to $ \Sigma$, whose inverse image is $\{m\} \times \tilde{\Sigma}_\sigma $.
\end{exo}

\begin{exo}
This exercise requires the notion of flat leaves of Exercise \ref{exo:flatleaves}. Show that
the leaf $ L \simeq \Sigma$  in Exercise \ref{exo:leafSigma} is flat.

{\emph{Hint}}: Show that the quotient of the direct product singular foliation $ (M,0) \times (\tilde \Sigma_\sigma,\mathfrak X(\tilde \Sigma_\sigma))$ under the $\pi_1(\Sigma, \sigma) $-action is a regular foliation admitting $ L \simeq \Sigma$ as a leaf. \end{exo}

\subsection{New constructions from old ones in algebraic geometry}

In this section, we work with algebraic singular foliations as in Definition \ref{def:algebraic}.
Let us repeat the context. In this section, $\mathcal O $ is an Abelian algebra and ${\mathrm{Der}}(\mathcal O) $ stands for the $ \mathcal O$-module of derivations of $ \mathcal O$ (which is a Lie algebra\footnote{It is even a Lie-Rinehart algebra} when equipped with the commutator). We define algebraic singular foliation over $ \mathcal F $ with respect to $ \mathcal O$ to be sub-$\mathcal O$-modules of $ {\mathrm{Der}}(\mathcal O) $ which are finitely generated\footnote{The assumption “finitely generated” is never used in this section and could be erased. However, we decided to keep it, since it is inherited to the new sub-modules  that we will construct.} and stable under the Lie bracket of derivations.
The purpose of the section is to explain how purely algebraic methods, allow defining new algebraic singular foliations out of this one.

We invite the reader to think that $ \mathcal O$ is the algebra of polynomial functions on $\mathbb K^d $,
or is the algebra of functions on some affine variety\footnote{It is here implicit that our affine varieties are over the field $ \mathbb C$.} $ W$, i.e., a quotient of the algebra of polynomial functions on $\mathbb C^d $ by a prime ideal $ \mathcal I$. Recall that since $ \mathcal O$ is Noetherian, any sub-module of the module of derivations is finitely generated. In this particular case, our constructions have a geometric meaning that we will detail.

\vspace{0.5cm}

Throughout this section, we choose an algebraic singular foliation $\mathcal F$ over the algebra $\mathcal O$.

\subsubsection{Restriction}
An ideal $\mathcal I $ is said to be a \emph{foliated ideal} if  $$\mathcal F [\mathcal I] \subset \mathcal I.$$
The quotient space $\mathcal F / \mathcal I \mathcal F $ then inherits a natural algebraic singular foliation structure over $\mathcal O / \mathcal I$. We call the latter algebraic singular foliation the \emph{restriction w.r.t the ideal $ \mathcal I$}. 

\begin{example}
When $ \mathcal O$ is the algebra of functions on an affine variety $W$, and $\mathcal I$ is the ideal of functions vanishing on an affine subvariety $W' \subset W$, then $ \mathcal I$ is a foliated ideal if and only if all vector fields in $\mathcal F $ are tangent to $W$, and the previous construction corresponds to the restriction of the singular foliation $ \mathcal F$ to $W $.
\end{example}

\subsubsection{Algebra Extension}

Assume that the algebra $\mathcal O $ has no zero divisor, and let $\mathbb O $ be its field of fractions.
Any derivation of $ \mathcal O$ extends to a derivation of $ \mathbb O$, so that we have a natural inclusion\footnote{This inclusion is even an equality for $\mathcal O$ a polynomial algebra over $ \mathbb K$.} $ \mathbb O \otimes_{\mathcal O} {\mathrm{Der}}(\mathcal O) \subset  {\mathrm{Der}}(\mathbb O)$.
For any subalgebra $\tilde{\mathcal O}$ with $\mathcal O \subset \tilde{\mathcal O} \subset \mathbb O$ such that every derivation $X\in \mathcal F$ (extended to a derivation of $ {\mathbb O}$) preserves $ \tilde{\mathcal O} $, there is natural algebraic singular foliation  over  $\tilde{\mathcal O} $ given by $ \tilde{\mathcal O}\otimes_{\mathcal O}\mathcal F \subset \mathrm{Der} (\tilde{ \mathcal O})$. 

\begin{example}
We will not try to give a complete geometric description of this construction in the context of affine varieties. However, let us mention that the blow-up of a singular foliation at a point of Section \ref{sec:blowup1} is a construction of that type on any affine chart.
\end{example}

\subsubsection{Localization}

Let us recall the definition of localization.

\begin{definition}
 A subset $S\subset\mathcal O$ is called \emph{multiplicatively closed} if $1\in S$ and if $S$ is stable under multiplication. For a multiplicative set $S\subset \mathcal O$, the \emph{localization} of $\mathcal O $ at $S$ is the algebra defined as follows:
\begin{enumerate}
\item If $\mathcal O$ has no zero divisor, then it is the subalgebra $ S^{-1}\mathcal O$ of its field $ \mathbb O$ of fractions given by: 
\begin{equation}
     S^{-1}\mathcal O:=\left\lbrace \frac{f}{s} \middle| f\in\mathcal O,\, s\in S\right\rbrace
 \end{equation}
\item If $\mathcal O $ has zero divisor, the previous definition can be enlarged by considering the quotient of $\mathcal O\times S$ by the equivalence relation $(f, s) \sim (g, t)$ if there is an element $u\in S$ such that $u(ft-gs)=0$. Addition and multiplication are then defining by checking that the following operations go to the quotient:
\begin{equation}(f,s) +(g,t):=(ft+ gs, st)\quad \text{and}\quad (f,s) \times (g,t)  := (fg,st). \end{equation}
\end{enumerate}
\end{definition}
 
\begin{remark}
The algebra $\mathcal O$ is a subalgebra of $S^{-1}\mathcal O$ via the homomorphism $\mathcal O \hookrightarrow S^{-1}\mathcal O,\, f\mapsto \frac{f}{1}.$
\end{remark}

\begin{example}
    For an affine variety $W$, localization can be interpreted as meaning that we restrict ourselves to the Zarisky open subset $ \mathcal U$ of $W$. Then $S$ is the multiplicative subset of all polynomial functions $P$ whose zeros are outside $\mathcal U $ (i.e., such that $ P(w) \neq 0$ for all $ w \in \mathcal U$). This interpretation explains the name.
\end{example}

Any derivation $ X \in {\mathrm{Der}}(\mathcal O)$ extends to a derivation of the localization $ S^{-1}\mathcal O$. When  $ \mathcal O $ has no zero divisor, the extension is defined by:
\begin{align*}
   X : S^{-1}\mathcal O&\longrightarrow S^{-1}\mathcal O\\\frac{f}{s}&\longmapsto \frac{X[f]s-fX[s]}{s^2}. 
\end{align*}
When zero divisors exist, then one has to check that the following map goes to the quotient:
 $$ X \colon (f,s) \mapsto ( X[f]s-fX[s],s^2  )  $$
with respect to the equivalence relation defined above. 

In both cases, the $ S^{-1} \mathcal O$-module generated by extensions to  $S^{-1}\mathcal O $ of derivations in an algebraic singular foliation $ \mathcal F$  over $ S^{-1}\mathcal O$  that we denote by $S^{-1}\mathcal{F}$ and call \emph{localization of $\mathcal{F}$ at $S$}.
When no zero divisor exists on $ \mathcal O$,  elements in $S^{-1}\mathcal{F}$ can be thought of  as quotients $ \frac{X}{s} $ with $X \in \mathcal F$ and $s \in S$. 
The Lie bracket restricted to $S^{-1}\mathcal F$ is given as follows:  \begin{equation}
    \forall\; X,Y\in\mathcal F,\, \forall\;(s,t)\in S^2,\,\;\left[\frac{1}{s}X,\frac{1}{t} Y\right] =\frac{1}{st}\left[X,Y\right] + \frac{Y [s]}{s^2t} \, X-\frac{X [t]}{st^2}Y.
\end{equation}

\begin{example}
Geometrically, localization corresponds to restriction to a Zariski open subset.
\end{example}

\subsection{Blowup of a singular foliation}
\label{sec:blowup1}

In this section, we work in the realm of complex algebraic geometry over  $\mathbb C $.
Most construction extend to the  smooth setting: indeed, this is the context in which Debord and Skandalis \cite{DS:Boutet} introduced the notion of blow-up of a singular foliation.

\subsubsection{Blow-up at a point}

Recall that for any $d\in \mathbb{N}$, the set $\mathbb{P}^{d}$ 

of all straight lines through the origin of  $\mathbb{C}^{d+1}$ is a complex manifold of dimension $d$ over $\mathbb{C}$, called the \emph{$d$-dimensional projective space}. Formally, it is defined as the equivalence classes of relation on the quotient $\mathbb{C}^{d+1}\setminus \{(0,\ldots,0)\}$
under the  equivalence relation:
$$u=(u_0 , u_1,\ldots, u_d)\sim v=(v_0, v_1,\ldots, v_d) \Longleftrightarrow  \exists\lambda\in\mathbb{C}\setminus\{0\}\;\text{such that}\; u=\lambda v.$$
Equivalently, it can be defined as the quotient manifold $$\mathbb{P}^{d}:=\mathbb{C}^{d+1}\setminus \{(0,\ldots,0)\}/\mathbb{C}\setminus\{0\},$$
where the group $\mathbb{C}\setminus\{0\}$ acts by diagonal multiplication on $\mathbb{C}^{d+1}$.
In particular, elements in $\mathbb{P}^{d}$ shall be denoted as $ d+1$-tuples of elements not all equal to zero and defined up to a non-zero constant, and denoted by 
 $[x_1 \colon \dots \colon  x_{d+1}] $.

\begin{lemma}
The projective space
$\mathbb{P}^{d}$ is a complex manifold of dimension $d$. It is given by the $d+1 $ following charts:
 $$ \psi_i \colon (x_1, \dots,\hat{x}_i,\ldots, x_{d+1}) \mapsto [ x_1 \colon \dots \colon x_{i-1}\colon \underbrace{1}_{ \hbox{\tiny{$i^{th}$ term}}} \colon  x_{i+1}\colon \dots \colon x_{d+1}] .$$
\end{lemma}

The idea of the blow-up at the origin consists in replacing $\mathbb C^{d+1} $, by pairs made of straight lines through the origin (=elements of $\mathbb P^d $) and a point on that straight line.

\begin{definition}
The blow-up $\mathrm{Bl}_0(\mathbb C^{d+1})$ of $\mathbb C^{d+1}$ at the origin  consists of all pairs
$(D, z) \in \mathbb{P}^d \times \mathbb{C}^{d+1}$ such that $z\in D$.
\end{definition}

Given coordinates $[x_0 \colon \dots \, \colon  x_{d}]$ and $(z_0, \dots, z_d) $ on $\mathbb{P}^d $ and $ \mathbb{C}^{d+1}$ respectively,
we can describe 
$Bl_0(\mathbb C^{d+1})$ 
in terms of coordinates:
$$\mathrm{Bl}_0(\mathbb C^{d+1})=\{(x,z)\in\mathbb{P}^d \times \mathbb{C}^{d+1}\mid z_ix_j = z_jx_i,\; i,j = 0 ,\ldots, d\}.$$
These equations make sense, because multiplying all the $x_i$'s by a non-zero factor leave them invariant.

\begin{lemma}\label{lem:chartsBloxup}
$\mathrm{Bl}_0(\mathbb C^{d+1})$ is a complex manifold of dimension $d+1$. It is given by the $d+1 $ following charts:
 $$ \phi_i \colon ( x_1, \dots, x_{d}) \mapsto ([ x_1 \colon  \dots \colon  x_{i-1} \colon  \underbrace{1}_{ \hbox{\tiny{$i^{th}$ term}}} \colon  x_{i+1} \colon  \dots \colon  x_{d}],( x_i x_1 , \dots x_i x_{i-1},\underbrace{x_i}_{ \hbox{\tiny{$i^{th}$ term}}}, x_i x_{i+1} ,\dots, x_i x_{d} )) .$$
\end{lemma}

In these charts, the natural projection $ \sigma\colon \mathrm{Bl}_0(\mathbb C^{d+1})\longrightarrow \mathbb C^{d+1}$ mapping the pair $(D,z) $ to $z$  is given by the projection onto the second factor.

For $z\neq 0$ the pre-image $\sigma^{-1}(z)$ is pair $(D,z) $ with $D$ being the unique line $D\in \mathbb{P}^d$ passing through $z\in \mathbb{C}^{d+1}$. But $\sigma^{-1}(0) \simeq \mathbb{P}^d$.
Last:
\begin{equation}
\label{eq:ourtsidesingular}
\sigma \colon \mathrm{Bl}_0 (\mathbb{C}^{d+1}) \backslash \sigma^{-1}(\{0\})  \longrightarrow  \mathbb C^{d+1} \backslash \{0\} 
\end{equation}
is a biholomorphism. In particular,  any vector field on $ \mathbb C^{d+1}$ can be transported to a vector field on  $ \mathrm{Bl}_0 (\mathbb{C}^d) \backslash \sigma^{-1}(\{0\}) $. It is natural to ask whether  this vector field can be extended to $ \mathrm{Bl}_0 (\mathbb{C}^{d+1})$ or not. In general  the answer is “no”, but it is “yes” if the vector field vanishes at $0$, as we now see.

\begin{proposition}\label{prop:blowUpExists}
For a holomorphic vector field $X$ of $\mathbb C^{d+1} $, the following two points are equivalent:
\begin{enumerate}
    \item[(i)] $X$ vanishes at $0$
    \item[(ii)] there exists a vector field $\tilde{X} $ on $ \mathrm{Bl}_0 (\mathbb{C}^{d+1}) $ $ \sigma$-related\footnote{i.e., $ T_{(D,z)} \sigma (\tilde{X}_{_{(D,z)}})= X_{|_z}$. Equivalently $\sigma_* \tilde{X}=X $} to $ X$.
\end{enumerate}
If it exists, then the vector field in item (ii) is unique.
\end{proposition}
\begin{proof}
On the $i$-th of the $d+1 $ charts of Lemma \ref{lem:chartsBloxup}, $\sigma $ reads:
 $$ ( x_1, \dots, x_{d+1}) \mapsto ( x_i x_1 , \dots x_i x_{i-1},x_i,x_i  x_{i+1} ,\dots, x_i x_{d+1}) .$$
The pull-back of the coordinate functions $(z_1, \dots, z_{d+1}) $ of $\mathbb C^{d+1} $ are therefore given by 
 $$
 \sigma^* (z_j) = \left\{ \begin{array}{ll} x_i x_j & j \neq i \\  x_i & j=i \end{array} \right. 
 $$
 This implies  that the unique vector field $X_j $ on that chart such that 
 $$ \sigma_* ( X_j) = \frac{\partial}{ \partial z_j} $$ 
 is given in the coordinates $ (x_1, \dots, x_d)$ by
 $$
X_j = \left\{ \begin{array}{ll} \frac{1}{x_i} \frac{\partial}{ \partial x_j}  & j \neq i \\   \frac{\partial}{ \partial x_i} - \sum_{j \neq i} \frac{x_j}{x_i} \frac{\partial}{ \partial x_j}  & j=i \end{array} \right. 
$$
In turn, this implies that for every vector field $ X  = \sum_{i=1}^{d+1} P_i(z_1, \dots, z_d)  \frac{\partial}{ \partial z_j}   $, the unique vector field on the $i$-th chart 
such that $\sigma_* (Z)=X $ is 
$$ Z = \sum_{j \neq i} \left( \frac{P_j(x_1 x_i, \dots, x_i, x_i x_d)}{x_i} - \frac{x_j P_i(x_1 x_i, \dots, x_i, x_i x_d)}{x_i} \right) \frac{\partial}{ \partial x_j} + \frac{P_i(x_1 x_i, \dots, x_i, x_i x_d)}{x_i}  \frac{\partial}{ \partial x_i} $$
This vector field is well-defined on the whole chart if and only if the functions $P_1, \dots, P_{d+1}$ vanish at the origin.
This proves the claim.
\end{proof}

\begin{propositions}{Blow-up of a singular foliation at the origin}{prop:blowUp}
Let $\mathcal F $ be a holomorphic or algebraic singular foliation on  $\mathbb C^{d+1} $.
Assume all vector fields on $\mathbb C^{d+1}$ vanish at $0$. Then there exists a unique singular foliation $\widetilde{\mathcal F} $ on
$\mathrm{Bl}_0 (\mathbb{C}^d)$ such that \eqref{eq:ourtsidesingular}
is an isomorphism of foliated manifolds.

We call $\widetilde{\mathcal F} $ the blow-up  of $\mathcal F $ at the origin.
\end{propositions}

\begin{exo}
In this exercise, we call $\sigma^{-1}(0) $ the \emph{exceptional divisor} of the blow-up: its points are canonically identified with straight lines through the origin.
Let $\mathcal F $ be a singular foliation on $\mathbb C^{d+1} $ made of vector fields vanishing at $0$, and let $\tilde{\mathcal F} $ be its blow up.
\begin{enumerate}
    \item Let $X \in \mathfrak X (\mathbb C^d) $ be a vector field vanishing at $0  $ and $\tilde{X} \in \mathfrak X (\mathrm{Bl}_0 (\mathbb{C}^d))  $ such that $\sigma_* (\tilde{X})= X $. Show that $\tilde{X} $ vanishes at every point of the exceptional divisor
    
    if and only if 
     $$X = \lambda \sum_{i=1}^{d+1} z_i \frac{\partial}{\partial z_i} + \hbox{  quadratic terms } $$
     for some $\lambda \in \mathbb C $ where ``quadratic terms'' means vector fields vanishing at least quadratically at zero. 
    \item Show that some point $ D $ in the exceptional divisor is a point-leaf\footnote{I.e, a leaf reduced to a point - equivalently a point where the tangent space of the foliation is zero.} if and only if $D$ (seen now as a straight line) is an eigenvector for all the linearizations of all vector fields in $\mathcal F $. 
\end{enumerate}
\end{exo}

\subsubsection{Blow-up along a smooth submanifold}

This construction of the blow-up of a singular foliation at a point can be extended considerably to a blow up along a submanifold to which the singular foliation is tangent. 
Since we presented the previous construction in the complex setting, we decided to present this construction in the smooth setting: it can of course be extended to the complex or real setting, but this requires using a  procedure more algebraic than the one described below, that uses the tubular neighborhood theorem (which holds in the smooth case only).

Let $N$ be a submanifold of $ M$.
Denote by $\mathcal{N}_{N/M}$ the normal bundle $TM_{|_N}/TN $ of a submanifold $N \subset M$. The fiber bundle $\mathbb{P}(\mathcal{N}_{N/M}) \to N$ can be interpreted as the projective space of directions normal to $N$ in $M$. The blow-up of $M$ along $N$ is a manifold obtained by gluing $ M \backslash N$ and\footnote{
For a vector bundle $E \to N$, we denote by
$\mathbb{P}(E)$ the projective bundle over $N$ namely, the complex manifold obtained by taking the projective space of all the fibers of $E$.} $ \mathbb{P}(\mathcal{N}_{N/M}) $. The construction goes as follows.
According to the tubular neighborhood theorem, there exists a diffeomorphism $ \Phi$ from a neighborhood $ \mathcal U$ of the zero section in $\mathcal{N}_{N/M} $ to a neighborhood $ \mathcal V$ of $N$ in $M$: it can be obtained, for instance, by considering $\mathcal{N}_{N/M}$ as a sub-bundle of $ TM_{|_N}$, then using the geodesic flow: this requires to choose a metric on $M$. This diffeomorphism $ \Phi$ is the identity when restricted to $ N$.

Now, let $\mathrm{Bl}_N(\mathcal{N}_{N/M}) $ be the fiber bundle over $N$ obtained by applying to each fiber of $\mathcal{N}_{N/M}$ the construction of the blow-up at a point: an element in $ 
\mathrm{Bl}_N(\mathcal{N}_{N/M})$ is a pair $ (P,u)$ with $P$ a straight line through $0$ in some fiber of $\mathcal{N}_{N/M}$ and $u$ a point in that straight line. We denote by $ \tau$ the projection $ \mathrm{Bl}_N(\mathcal{N}_{N/M}) \to \mathcal{N}_{N/M}$. By construction, $ \tau$ is a diffeomorphism  from $\tau^{-1}(\mathcal{N}_{N/M}\backslash N )$ to $ \mathcal{N}_{N/M}\backslash N  $.
The idea now consists in gluing
\begin{enumerate}
\item the manifold $ M \backslash N$
\item with the manifold $\mathrm{Bl}_N(\mathcal{N}_{N/M})  $
\end{enumerate}
by identifying 
\begin{enumerate}
\item
the open subset $ \mathcal V \backslash N $
\item with the open subset 
$ \tau^{-1}(\mathcal U \backslash N )$ 
\end{enumerate}
with the help of the diffeomorphism $$ \Phi \circ \tau \colon  \tau^{-1}(\mathcal U \backslash N ) \longrightarrow \mathcal V \backslash N $$
The result is a manifold $\mathrm{Bl}_N(M)$, called \emph{blow-up of $M$ along $N$}, equipped with a natural projection $ \sigma $ onto $M$, which is a diffeomorphism from $  \sigma^{-1}(M \backslash N) $ to $M \backslash N$, and whose fiber over a point  $n \in N$ is the projective space of the normal bundle $\mathcal{N}_{N/M}$ at $n$.
The set $\sigma^{-1}(N)\simeq \mathbb{P}(\mathcal{N}_{N/M})\subset \mathrm{Bl}_N(M)$ is a hypersurface called the \emph{exceptional divisor} of the blow-up $\sigma\colon \mathrm{Bl}_N(M)\longrightarrow M$.

\begin{example}
Let $N=\mathbb{R}^m $ be a linear subspace of $ \mathbb R^d$ be a linear subspace. Then 
$\mathrm{Bl}_N(M) \simeq  \mathrm{Bl}_0(\mathbb R^{m})   \times \mathbb R^{d-m}$.

\end{example}

The next statement has been established by Debord and Skandalis \cite{DS:Boutet}.

\vspace{.5cm}

\begin{propositions}{Blow-up of a singular foliation, general case}{prop:Blow-up}
Let $\mathcal{F}$ be a singular foliation on $M$ and $L\subset M$ a leaf.

There exists a unique singular foliation $\tilde{\mathcal F} $ on the blow-up $\mathrm{Bl}_L(M)$ of $M$ along $L$ such that $\sigma $ is an isomorphism from $\mathrm{Bl}_L(M) \backslash \sigma^{-1}(L)$ to $M \backslash L $. 

It is called the \emph{the blow-up of $\mathcal{F}$ along $L$}.
\end{propositions}

\vspace{.5cm}

In fact, we not need to take $L$ to be a leaf: the construction would work for any submanifold to which all vector fields in $\mathcal F $ are tangent. The proof is based on a lemma: a vector field  $X \in \mathfrak X (M)$ reads $ X = \sigma(\tilde{X})$ for $\tilde{X} $ a vector field in $\mathrm{Bl}_L(M) $ if and only if $X$ is tangent to $L$.

\subsubsection{Monoidal transformation}

The constructions of blow-ups above can be enlarged to any monoidal transformation with respect to any foliated ideal. 

Let us work in the algebraic setting\footnote{Sertöz does the construction on any complex manifold, see \cite{zbMATH04118703}.}. 
Let $M:= \mathbb C^n $ and $ \mathcal O := \mathbb C[z_1, \dots, z_d]$.
Consider an ideal $ \mathcal I$, and choose $ \varphi_1, \dots, \varphi_k$ generators of this ideal.
The monoidal transformation of $M$ with respect to $ \mathcal I$ is the sub-variety $M_\mathcal I $ of  $\mathbb P^{k-1} \times \mathbb  C^d  $ defined by the equations
 \begin{equation}\label{eq:monoidal}   x_i \phi_j (z_1, \dots, z_d) =  x_j \phi_i(z_1, \dots, z_d) \hbox{ for $ 1 \leq i <  j \leq k$}, \end{equation}
where $ z_1, \dots, z_d$ are the coordinates on $ \mathbb C^d$ and $ [x_1 \colon \cdots \colon x_k ]$ are the homogeneous coordinates on the projective space $\mathbb P^{k-1} $, one then chooses the component of this subvariety that projects onto $M$.
Alternatively, $ M_\mathcal I  $ can be seen as the Zariski closure of the graph of the map
 $$   \begin{array}{rcl}  M \backslash V_{\mathcal I} & \to &  \mathbb P^{k-1}\\ (z_1, \dots, z_d)&\mapsto  &  [ \phi_1(z_1, \dots, z_d ) \colon \dots \colon\phi_k(z_1, \dots, z_d ) ]. \end{array}  $$
 Here $ V_{\mathcal I} \subset M$ is the zero locus of $ \mathcal I$, i.e., the sub-variety of all points where all functions in $ \mathcal I$ are zero.
It is a classical result that $ M_\mathcal I$ is a quasi-projective variety, and that the natural projection on the first component:
\begin{equation}\label{eq:projection} 
\phi \colon M_\mathcal I \to M 
\end{equation}
is a proper map, which is a biholomorphism\footnote{In fact, it is even biregular.} when restricted to 
\begin{equation}\label{eq:projection2} \phi^{-1}( M \backslash V_{\mathcal I} ) \to  M \backslash V_{\mathcal I} .\end{equation}

Now, consider $ \mathcal F$ an algebraic singular foliation on $ \mathbb C^d$.

If $ \mathcal I$ is a foliated ideal, i.e., if\footnote{Geometrically, this should be interpreted as meaning that a leaf of $ \mathcal F$ is either contained in $V_\mathcal I $ or has empty intersection in $V_\mathcal I $, i.e that $V_\mathcal I $ is a union of leaves of $ \mathcal F$.}  $ \mathcal F [\mathcal I] \subset \mathcal I$, then $ M_\mathcal I$ comes induced with a natural singular foliation $ \mathcal F_\mathcal I$ such that the map defined in \eqref{eq:projection2} is an isomorphism of foliated manifolds. In particular, it means that for every $X \in \mathcal F$, the vector field on $ \phi^{-1}( M \backslash V_{\mathcal I})$ defined by transporting $X$ through Equation \eqref{eq:projection2}  extends to a vector field $ \tilde{X}$ to the whole variety $ M_\mathcal I$. 
We will not make a more precise statement, simply because doing it would require to make sense of singular foliations on schemes, or at least quasi-projective varieties - something that we have not done in the present manuscript. We prefer to refer the reader to the literature on the subject, namely \cite{zbMATH04118703}, or the more recent \cite{Ruben3,louis2024nash}.

\begin{exo}
Let us explain how to extend a vector field $X$ on $ \mathbb C^d$ such that $ X[\mathcal I] \subset \mathcal I$ to a vector field on $ M_\mathcal I$.  Let $\varphi_1, \dots,\varphi_k$ be generators of $ \mathcal I$. 
\begin{enumerate}
\item Show that there exists $ (\lambda_{ij})_{i,j=1, \dots, k} \in \mathcal O$ such that $$ X [\varphi_i] = \sum_{j=1}^k \lambda_{ij} \varphi_j .$$
\item Show that the vector field defined on $ \mathbb P^{k-1} \times \mathbb C^d $ by\footnote{Notice the inversion of the indices $a$ and $b$.}
 $$  \left( \sum_{a,b=1}^k \lambda_{ba} \,  x_a \frac{\partial}{\partial x_b} ,X\right)    $$
 is tangent to the subset defined by Equations \eqref{eq:monoidal}.
 \item Conclude.
\end{enumerate}

\end{exo}

\subsection{Nash-blowup of a singular foliation}

\label{sec:Nash}

We now introduce another construction of a ``blow-up'' of a singular foliation, which is due to Omar Mohsen \cite{OmarMohsen} who used it, with his collaborators, to prove several hard results about PDEs \cite{AMY}. In fact, this blowup belongs to the large class of Nash blowups, and we will therefore call it the Nash blow-up of a singular foliation. A particular feature of this construction is that it produces a \underline{Debord} singular foliation, which comes therefore from a Lie algebroid. To be more precise, it might be that the Nash blowup of a foliated manifold is not a smooth manifold, i.e., it may have singular points. But 1) whenever the Nash blowup is a smooth manifold, it becomes a foliated manifold with a Debord singular foliation, and 2) the conclusion is still valid in general, but we have to make sense of the notion of projective foliation on a singular quasi-projective variety, see \cite{louis2024nash}. The original article of Mohsen insists more on the topological groupoid that lies on the top of this Lie algebroid, showing that it has several desirable topological properties that we do not discuss here, see \cite{OmarMohsen}.

\begin{bclogo}[  arrondi = 0.1, logo = \bcdz]{Warning !}
The Nash blowup method uses several concepts that are only introduced much later in the text. We recommend the reader unfamiliar with these notions to look first at Section \ref{sec:almost}. The ``level'' of this section is, in general, higher than the level of the neighboring sections.

As usual, we try to deal with the smooth, real analytic and complex cases altogether. Last, since many exercises and 
 remarks divide the text, we decided to write them using smaller font, to make the main road of the text easier to follow.
\end{bclogo}
\subsubsection{Generalities about Grassmann bundles}

We start by recalling the notion of $k$-Grassmannian. We work in the complex setting, but we will mention at some point what happens in the real one.
For any $N, k \in  \mathbb N$ such that $0\leq k\leq N$, the set  $\mathrm{Grass}_{k}(\mathbb C^N)$  of all $k$-dimensional vector subspaces of $\mathbb C^N$ is a complex compact manifold of dimension $k(N-k)$, called  the \emph{Grassmanian of $k$-planes in  $\mathbb C^N$} (see e.g., \cite{sontz, Bredon,Nakahara}). 
For $k=1$, one recovers the notion of projective space.

Let us spell out the topology and the manifold structure of $\mathrm{Grass}_{k}(\mathbb C^N)$. The groups $\mathrm{GL}(\mathbb C)$ and $\mathrm{U}(N)$ of invertible $\mathbb C$-linear maps, and of unitary linear transformations acts transitively on $\mathrm{Grass}_{k}(\mathbb C^N)$ by $ g \cdot V = g(V) $
for all $V \in \mathrm{Grass}_{k}(\mathbb C^N)$ and $g \in \mathrm{GL}(\mathbb C)$ or $\mathrm{U}(N) $.
This action is well-defined since $g$ is invertible. Therefore,  \begin{eqnarray*}
   \mathrm{Grass}_{k}(\mathbb C^N) & \stackrel{set}{\simeq}& \frac{\mathrm{ U}(N)}{\mathrm U(k)\times \mathrm U(N-k)} \\ &\stackrel{set}{\simeq}   & \frac{\mathrm{GL}(N)}{\mathrm{GL}(k,N)} 
\end{eqnarray*} 

Above the subgroups $\mathrm U(k)\times \mathrm U(N-k)$ and $\mathrm{GL}(k,N)$  are the respective stabilizers of $\mathbb C^k \times \{0\}^{N-k} \subset \mathbb C^N $ in their corresponding groups.  Since both groups in the first line above are compact Lie groups, the first description equips $\mathrm{Grass}_{k}(\mathbb C^N)$ with a structure of compact manifold. Since both groups in the second line are complex Lie groups, the second description equips the quotient with a structure of complex manifold. If one replaces $\mathbb C$ by $\mathbb R $, then this second structure does not make any sense, but the compact smooth manifold structure still exists.

{\small{
\begin{remark}
The manifold $\mathrm{Grass}_{k}(\mathbb C^N) $ has natural holomorphic coordinates, called \emph{canonical affine coordinates}. One such chart $\phi := \mathbb C^{k(N-k)} \to \mathrm{Grass}_{k}(\mathbb C^N)  $ maps a $ k \times (N-k)$-matrix $S$ to the subspace generated by the column vectors:
$$
 \begin{pmatrix} 
  {\mathrm{id}}
  \\ 
  S
     \end{pmatrix} 
  $$   The image of this chart is the set of $k$-dimensional subspace that do not intersect $ \{0\}^k \times \{\mathbb C\}^{N-k}$.
The others natural charts are obtained by composition with the $ \mathrm U(N) $-action, i.e., are of the form $ g \cdot \phi $ for some $g \in \mathrm U(N) $.
\end{remark}
}}

\vspace{0.5cm}

\begin{con}
For our current purpose, it will be convenient to consider the Grassmannian of all sub-spaces of \underline{co}-dimension $r$ in $\mathbb C^N $ rather than the set of subspaces of a given dimension. It is convenient to denote this manifold by 
$\mathrm{Grass}_{-r}(\mathbb C^N)$ (notice the use of a minus sign). In other words, we set $\mathrm{Grass}_{-r}(\mathbb C^N) := \mathrm{Grass}_{N-r}(\mathbb C^N)$.
\end{con}

Now, given a vector bundle over a  manifold, one can the Grassmannian at each fiber. The next definition makes sense in the smooth or complex cases without adaptation.

\begin{definition}
Let $A\rightarrow M$ be a vector bundle of rank $d$ over a manifold $M$. Let $k\leq \mathrm{rk}(A)=d$. The disjoint union:
  $$\mathrm{Grass}_{k}({A}):= \coprod_{x\in M} \mathrm{Grass}_{k}(A_{x})  $$
  comes equipped with a natural manifold structure. Also
 \begin{equation}\label{eq:Pi}
     \Pi\colon \mathrm{Grass}_{k}(A)\longrightarrow M
 \end{equation}
 is a fibration with  fiber $ \simeq \mathrm{Grass}_{k}(\mathbb K^d) $ ($\mathbb K=\mathbb R$ or $\mathbb C $). It is called \emph{$k$-th Grassmann bundle of $A$}. The same holds upon replacing $k$ by $ -k$.
 
\end{definition}

The following exercise uses the notion of linear vector fields on an vector bundle (see Section \ref{sec:almost} for details). Again, it is valid as it is in both smooth or complex cases.

\begin{exo}\label{exo:grass+vf}
Let $A\longrightarrow M$ be a vector bundle over $M$.  Show that a linear vector field on $A$ induces a vector field on $\Pi\colon \mathrm{Grass}_{-r}(A)\longrightarrow M$ that is $\Pi$-projectable\footnote{I.e. $ \Pi$-related (see Definition \ref{def:related}) to a vector field in $\mathfrak X(M) $} on $M$.  

{\small{\emph{Hint:} the flow of a linear vector fields on $A$ is a pseudo-group of diffeomorphism of $A$, hence it induces a pseudo-group of vector bundle isomorphisms, which then induces a pseudo-group of fiber bundle diffeomorphisms on $\mathrm{Grass}_{-r}(A)$.}}
\end{exo}

\subsubsection{Nash-blowup I: the space}

Throughout this subsection $(M, \mathcal{F})$ is  a foliated  manifold with $M$ connected, and  $(A, \rho)$ is an anchored bundle over $\mathcal F$.\\

We denote by $M_{\mathrm{reg}}\subset M$ the regular part of $(M,\mathcal F) $ (see Section \ref{sec:openReg}). By construction, it is the open dense subset such that the maps  $x\mapsto\mathrm{im}(\rho_x)$ and  $x\mapsto \ker (\rho_x)$ are locally constant, i.e., the set of regular points of $\mathcal{F}$. 

Here we have to make a distinction between the smooth or complex settings.
\begin{enumerate}
\item[$ \bullet$]
In the complex case, the codimension $r$ of $\mathrm{im}(\rho_x)=T_x \mathcal F \subseteq T_xM$ is the same for all points  $x\in M_{\mathrm{reg}}$ since $M$ is connected. This is also true in the real analytic setting. 
\item[$ \bullet$] In the smooth case, we will have from now on to \underline{assume} that it is the case\footnote{Equivalently, we have to assume that all the regular leaves have the same dimension. It is not the case, for instance, for the singular foliation on  $M=\mathbb R $ generated by the vector field $ \chi(t) \partial_t $ with $ \chi(t)= e^{-1/t^2}$ for $ t \geq 0$ and $ \chi(t)=0$ for $ t \leq 0$, since  then $M_{\mathrm{reg}}=  \mathbb R_-^* \cup  \mathbb R_+^* $ but regular leaves have dimension $0$ on $ \mathbb R_-^*$ and one on  $ \mathbb R_+^*$. See Section \ref{sec:openReg} for an overview of these notions}.
\end{enumerate}

Notice that for every point $x\in M_{\mathrm{reg}}$ in the regular part, 
\begin{enumerate}
\item $\mathrm{im}(\rho_x)$ is an element of the Grassmannian  $\mathrm{Grass}_{r}(T_xM)$, and
\item $\ker(\rho_x)$ is an element of $\mathrm{Grass}_{-r}(A_x)$.
\end{enumerate}
Consider now:
\begin{enumerate}
\item  the Grassmann bundle  $\mathrm{Grass}_{r}(TM)$, and
\item  the Grassmann bundle $\mathrm{Grass}_{-r}(A)$.
\end{enumerate}
As before, denote by $ \Pi$, in both cases, their natural projections onto $M$. Consider the two natural sections of $\Pi$  defined on the regular part $M_{\mathrm{reg}}$ by: 
\begin{enumerate}
\item  $\sigma_{\mathrm{im}}\colon M_{\mathrm{reg}}\longrightarrow \mathrm{Grass}_{r}(TM),\, x\longmapsto \mathrm{im}(\rho_x)$, and
\item  $\sigma_{\mathrm{ker}}\colon M_{\mathrm{reg}}\longrightarrow \mathrm{Grass}_{-r}(A),\, x\longmapsto \mathrm{ker}(\rho_x)$.
\end{enumerate}
Then we define 
\begin{enumerate}
\item  $\mathrm{Bl}_{tgt} (M,\mathcal{F}):= \overline{\sigma_{\mathrm{im}}(M_{\mathrm{reg}})}$ to be the closure of the image of the section $\sigma_{\mathrm{im}}$ in $\mathrm{Grass}_{r}(TM)$, and
\item $\mathrm{Bl}_{A}(M,\mathcal{F}):= \overline{\sigma_{\mathrm{ker}}(M_{\mathrm{reg}})}$ to be the closure of the image of the section $\sigma_{\mathrm{ker}}$ in $\mathrm{Grass}_{-r}(A)$.
\end{enumerate}
We shall denote by $ \pi$ the restriction of $\Pi$ to both $\mathrm{Bl}_{\mathrm{im}}(M,\mathcal{F})$ and $\mathrm{Bl}_{\mathrm{ker}}(M,\mathcal{F})$.

\begin{bclogo}[  arrondi = 0.1, logo = \bcdz]{Comparison with existing literature.}
We claim the definition of $\mathrm{Bl}_{A}(M,\mathcal{F})$  above matches the definition of the blowup of a singular foliation given by \cite{OmarMohsen}, and used in \cite{AMY}.
This requires a careful line by line check,  but the difference is mainly a difference of presentation.
It also matches Ali Sinan Sertöz's construction \cite{zbMATH04118703} of a Nash blowup coherent of a coherent sheaf when applied to the coherent sheaf $ \mathcal F$ (see Proposition \ref{prop:coherent}) up to a difference of context: smooth setting in Omar Mohsen, complex setting in Sertöz\footnote{ 
To be precise, Ali Sinan Sertöz 's  construction consists in applying the method described above to a vector bundle morphism $B \stackrel{d}{\to} A$ whose co-image is $\mathcal F $: the existence of such a vector bundle morphism is the very definition of coherent sheaves.}.

Also, in \cite{Ruben3}, a sequence depending on $ n \in \mathbb N$ of ``blow-ups'' is constructed, which for $n=0$ gives back $\mathrm{Bl}_{\mathrm{tgt}}(M,\mathcal{F})$, and for $n=1$ gives  $\mathrm{Bl}_{A}(M,\mathcal{F})$. 
\end{bclogo}

{\small{
\begin{exo}
Spell out using adapted coordinates the definition of $\mathrm{Bl}_{A}(M,\mathcal{F})$. Check that if there exists coordinates in which the anchored bundle has polynomial expressions (i.e., the anchor map of constant sections is a polynomial vector field), then $\mathrm{Bl}_{A}(M,\mathcal{F})$ is an affine variety in every chart given by adapted coordinates.  
(\emph{Hint:} A solution to this problem is given in \cite{Ruben3}-\cite{louis2024nash},
where it is proven that it is even a monoidal transformation.)
\end{exo}
}}

\begin{remark}
\label{rmk:limitesGrassmann}
 Intuitively, for $x\in M$, $\pi^{-1}(x)= \mathrm{Bl}_{A}(M,\mathcal{F})\cap \Pi^{-1}(x)$ is the set of  all possible limits of the subspaces $\mathrm{ker} \rho_{y}$ when $y\in M_{\mathrm{reg}}$ converges to $x$.
More precisely, for any $x\in M$, there is a   an open neighborhood $\mathcal U\subset M$ of $x$ such that $\mathrm{Grass}_{-r}(A)\simeq \mathcal{U}\times \mathrm{Grass}_{-r}(\mathbb{K}^{\mathrm{rk}(A)})$. By construction,   \begin{equation*}
        \pi^{-1}(x)=\left\{\text{Codim. $r$ subspaces } V\subset A_{_x}\;\middle|\; \exists\, (x_n)\in M_{\mathrm{reg}},\, \text{such that},\,\;\ker \rho_{x_n} \underset{n \to +\infty}{\longrightarrow} V \text{ as }  x_n\underset{n \to +\infty}{\longrightarrow}x\right\}.
        \end{equation*}
One can make  a similar construction with $\mathrm{Bl}_{tgt}(M,\mathcal{F})$:
  \begin{equation*}
        \pi^{-1}(x)=\left\{ \text{Dim. $r$ subspaces }  W\subset T_{_x}M\;\middle|\; \exists\, (x_n)\in M_{\mathrm{reg}},\, \text{such that},\,{\mathrm{im}} \rho_{x_n} \underset{n \to +\infty}{\longrightarrow} W \text{ as }  x_n\underset{n \to +\infty}{\longrightarrow}x\right\}.
        \end{equation*}
\end{remark}

{\small{
\begin{remark}
Since ${\mathrm{im}}(\rho_x)$ coincides with the vector subspace of $T_x M $ that we denoted by $T_x \mathcal F $, which also coincides with the tangent space of the leaf through $x$, the set $\mathrm{Bl}_{tgt}(M,\mathcal{F})$ does not depend on the choice of an anchored bundle. It is also the case for $\mathrm{Bl}_{A}(M,\mathcal{F})$ as we  now see.
\end{remark}
}}

\begin{lemma}\label{lem:indep} For any two anchored bundle $ (A,\rho)$ and $ (A',\rho')$, there is an unique homeomorphism\footnote{A close look at the proof shows that it is in fact more than a simple homeomorphism: the pull-back of a function on $\mathrm{Bl}_{A'}(M,\mathcal{F}) $ which is, locally, the restriction of a  locally defined smooth, real analytic or holomorphic function on $ {\mathrm{Grass}}_{-r}(A)$ is a function on $\mathrm{Bl}_{A}(M,\mathcal{F}) $ of the same type.} making the following diagram commutative
$$\xymatrix{\mathrm{Bl}_{A}(M,\mathcal{F}) \ar[dr]|-{\, \pi \, } \ar@{-}[rr]^{\simeq}& & \mathrm{Bl}_{A'}(M,\mathcal{F}) \ar[dl]|-{\, \pi' \, } \\ &M & }.$$
\end{lemma}
\begin{proof}[Proof\footnote{
We acknowledge for a discussion with Cédric Rigaud for this proof.}]
 This homeomorphism, if it exists, is unique, by density of the regular part. 
We saw in Proposition \ref{thm:prop:UniqueAnchored}, Section \ref{sec:AnchoredBundle}, that there exists anchored bundle morphisms (see Section \ref{sec:AnchoredBundle}) \begin{equation}
    \xymatrix{(A,  \rho)\ar@<2pt>[r]^\Phi&(A', \rho')\ar@<2pt>[l]^\Psi}.
\end{equation}
Now, $(A \oplus A', \rho+\rho') $  is also an anchored bundle over $\mathcal F $, and there is a commutative diagram as follows:
$$ \xymatrix{  &&& A \oplus A' \ar@{->>}[dlll]|-{ a + \Psi(a')
 \longleftarrow (a,a') }\ar[dd]|-{\begin{array}{c}  \hspace{1mm}\hbox{ \tiny{$ \rho+\rho'\,$}} \\ \end{array} }   \ar@{->>}[drrr]|-{(a,a') \longrightarrow \Phi(a)+a'}  &&& \\ A \ar[drrr]|-{\, \rho \, } &&&&&& A' \ar[dlll]|-{\, \rho' \, } \\ &&&TM &&&} .$$
Since the vector bundle morphisms $ A \oplus A' \to A$ and $ A \oplus A' \to A'$ above are surjective anchored bundle morphisms, it suffices to show that Lemma \ref{lem:indep} holds under the additional assumption that the anchored bundles be related by a surjective anchored bundle morphism. 

Consider two anchored bundles $(A,\rho) $ and, say, $(C,\rho_C) $ over $ \mathcal F$ such that there exists a surjective anchored bundle morphism $ \phi : C \to A$. Let $ K \to M$ be the kernel of $ \phi$. Consider the subset  $$\mathrm{Grass}_{-r}(C,K) \subset \mathrm{Grass}_{-r}(C)$$
of all vector sub-spaces of a fiber $C_m $ which satisfy two conditions: their codimension is $r$ and they contain $K_m$.  Since $K$ is a sub vector bundle of $C$, $\mathrm{Grass}_{-r}(C,K)$ is a \underline{closed} submanifold. Since $ K$ is in the kernel of $ \phi$, this submanifold is moreover canonically diffeomorphisc\footnote{Replace by biholomorphic in the complex case in the present discussion, and similarly for the real analytic case.} to $ \mathrm{Grass}_{-r}(A)$. We denote by $ \underline{\phi}$ this diffeomorphism.
Since, the anchor of $C$ is zero on $K$, the section 
$\sigma_{\mathrm{ker}}^C \colon M_{\mathrm{reg}}\longrightarrow \mathrm{Grass}_{-r}(C)$ used to construct $\mathrm{Bl}_{C}(M,\mathcal{F})$ is valued in  $\mathrm{Grass}_{-r}(C,K)$.
Moreover, it is related with the section 
$\sigma_{\mathrm{ker}}^A \colon M_{\mathrm{reg}}\longrightarrow \mathrm{Grass}_{-r}(A)$ used to construct $\mathrm{Bl}_{A}(M,\mathcal{F})$ by
$$\underline{\phi} \circ \sigma_{\mathrm{ker}}^C  =\sigma_{\mathrm{ker}}^A  .$$  
The diffeomorphism  $\underline{\phi}$ therefore also intertwines the closures of the graphs of  $\sigma_{\mathrm{ker}}^A$ and $ \sigma_{\mathrm{ker}}^C$, hence induces an homeomorphism $\mathrm{Bl}_{A}(M,\mathcal{F}) \simeq \mathrm{Bl}_{C}(M,\mathcal{F})$ that has all desired properties.
\end{proof}

 Lemma \ref{lem:indep}  implies\footnote{For singular foliations $ \mathcal F$ which are not finitely generated - which is the case of most of them, especially in the complex case, see Section \ref{sec:globally}} that $\mathrm{Bl}_A(M,\mathcal{F})$ is a \underline{global} object, defined over the whole manifold $M$. Since anchored bundles exist near every point,   $\mathrm{Bl}_A(M,\mathcal{F})$  can be defined in a neighborhood of every point. Since it does not depend on the choice of an anchored bundle, these locally defined sets glue in a natural manner to define a bundle $\mathrm{Bl}(M,\mathcal{F})$ over the whole complex manifold $M$. 

The next exercise\footnote{We refer the reader to Section \ref{sec:isotropy} for definitions of the isotropy Lie algebra and the definition of strong kernel at a point. We also invite the reader to use remark
\ref{rmk:limitesGrassmann}.
The questions asked in this exercise are solved in \cite{Ruben3}, although in  a slightly different context.} describes some important features of that construction.

\begin{exo}\label{exo:def-blup}
Let $\mathcal F$ be a singular foliation on a connected manifold $M$. In the smooth case, we assume that all regular leaves have the same dimension $r$.
Assume that an anchored bundle $(A,\rho) $ exists.
Prove that $\pi\colon \mathrm{Bl}_A(M,\mathcal{F})\rightarrow M$ satisfies the following properties:
 \begin{enumerate}
\item 
$\pi$ is proper and surjective. In particular, for each point $x\in M$, the fiber $\pi^{-1}(x)$ is a non-empty compact set. {\small{(\emph{Hint:} Use the fact that the projection $\Pi$ admits compact fibers.)}}
\item {For every $x\in M$ and $V\in \pi^{-1}(x)$, one has $\mathrm{Sker}(\rho,x)\subseteq V\subseteq \mathrm{ker}(\rho_x)$.} 
\item For every $x\in M$ and $V\in \pi^{-1}(x)$, the image $ [V]$ of $V$ in the isotropy Lie algebra at $x$, i.e., $$ \mathfrak g_x (\mathcal F) =
 \frac{\mathrm{ker}(\rho_{x})}{\mathrm{Sker}(\rho,x)}$$ is a sub-space codimension $r-\dim(L_x)$, where $\dim (L_x)$ is the dimension the leaf through $x$.
 Also, show that $V \to [V ] $ is an injective map from $\pi^{-1}(x) $ to ${\mathrm{Grass}}_{-r + \dim(L_x)} (\mathfrak g_x (\mathcal F)  )  $.
\item For every $x\in M$ and $V\in \pi^{-1}(x)$, $[V]$ is a Lie sub-algebra of $ \mathfrak g_x (\mathcal F) $ of codimension $r-\dim(L_x)$, where $\dim (L_x)$ is the dimension of the leaf through $x$.
\emph{Hint:} equip the anchored bundle with an almost Lie algebroid bracket and show that the almost Lie algebroid bracket restricts to a bilinear map $\wedge^2 V \to V$.
\item  For every $x\in M_{\mathrm{reg}}$,  $\pi^{-1}(x)=\mathrm{ker}(\rho_x)$ is reduced to a point in $\mathrm{Grass}_{-r}(A)$. Also, $\pi^{-1}(M_{\mathrm{reg}})$
is a smooth manifold and the restriction $\pi\colon \pi^{-1}(M_{\mathrm{reg}})\longrightarrow M_{\mathrm{reg}}$ is invertible\footnote{Invertible here means: diffeomorphism, in the smooth case, biholomorphism, in the complex case.}.
\end{enumerate}
\end{exo}

The conclusion of the previous exercise is the following. Let $x \in M$. Within the Grassmannian ${\mathrm{Grass}}_{r - \dim (L_x)} \left( \mathfrak g_x (\mathcal F)  \right) $ of all sub-spaces of codimension $r - \dim (L_x)$ in the isotropy Lie algebra  $ \mathfrak g_x (\mathcal F) $ at $x$ lies a subset, denoted by 
 $$ {\mathrm{Grass}}_{-r + \dim (L_x)}^{Lie} \left( \mathfrak g_x (\mathcal F)  \right) $$
made of all sub-Lie algebras  of codimension $r - \dim (L_x)$. This subset is compact. What the previous exercise gives is an injective inclusion of $\pi^{-1}(x) \subset  {\mathrm{Grass}}_{r}(A_{|_x})$ as a compact subset inside $ {\mathrm{Grass}}_{-r + \dim (L_x)}^{Lie} \left( \mathfrak g_x (\mathcal F)  \right)$.

{\small{
\begin{remark}
A remarkable result by Omar Mohsen \cite{OmarMohsen} is that ${\pi}^{-1}(x) \subset \mathrm{Grass}_{-r+\dim(L_x)}^{Lie}\left( \mathfrak g_x(\mathcal F)\right)$ lies in fact inside the set of sub-Lie algebras that integrate to a \underline{closed} sub-Lie group of the simply connected Lie group\footnote{In fact, \cite{OmarMohsen} established an even stronger result: it integrates to a closed subgroup of the isotropy of Androulidakis-Skandalis holonomy groupoid, which is a quotient of the universal one.} integrating $ \mathfrak g_x (\mathcal F) $. This is highly non-trivial: such Lie sub-algebras do not form a compact subset of $ {\mathrm{Grass}}_{-r + \dim (L_x)}^{Lie} \left( \mathfrak g_x (\mathcal F)  \right)$. 
\end{remark}
}}

Overall, we have therefore constructed an inclusion
\begin{equation}
\label{eq:insideholonmies}
\mathrm{Bl}_A(M,\mathcal{F}) \hookrightarrow \coprod_{x\in M}\mathrm{Grass}^{Lie}_{-r+\dim(L_x)}\left( \mathfrak g_x(\mathcal F)\right)
\end{equation}
The image of this inclusion does not depend on the choice of an anchored bundle.
Lemma \ref{lem:indep} and the previous discussion lead to the following conclusion:
 
\vspace{.5cm}

\begin{propositions}{$\mathrm{Bl}(M,\mathcal{F})$ is well-defined}{prop:NashWellDefined}
$\mathrm{Bl}_A(M,\mathcal{F})$ does not depend on the choice of  an anchored bundle $ (A,\rho)$ over  $\mathcal{F}$, or precisely:
\begin{enumerate}
\item the image of the injective map \eqref{eq:insideholonmies} does not depend on $ (A,\rho)$,
\item and for any other anchored bundle $(A',\rho') $
 over $ \mathcal F$, there is canonical homeomorphism 
  $\mathrm{Bl}_A(M,\mathcal{F}) \simeq \mathrm{Bl}_{A'}(M,\mathcal{F})$. 
  \end{enumerate}
\end{propositions}

\vspace{.5cm}

This independence allows the next definition.

\vspace{.5cm}
\begin{definitions}{Nash blow-up of a singular foliation}{def:NashWellDefined}
In view of Proposition \ref{thm:prop:NashWellDefined}, it makes sense to  denote by $\mathrm{Bl}(M,\mathcal{F})$ (with no reference to the chosen anchored bundle) the blowup of a foliated manifold $(M,\mathcal F)$ whose regular leaves all have the same dimension. We call it\footnote{Following \cite{zbMATH07128538}-\cite{zbMATH04118703}.} the \emph{Nash blowup of $ (M,\mathcal F)$}. 
\end{definitions}

\vspace{0.5cm}
The following definition now makes sense.

\begin{definition}
\label{def:limitLie}
    Let $ (M,\mathcal F)$ be a singular foliation and $r$ the dimension of the regular leaves.
    For every point $x \in M$, we call \emph{limit Lie subalgebras} the subalgebras of codimension  $r -\dim(L_x) $ of the holonomy Lie algebra $\mathfrak g_x(\mathcal F) $ that appear in exercise \ref{exo:def-blup}, i.e., which are the projections on $\mathfrak g_x(\mathcal F) $ of the limits of kernels of the anchor maps.
\end{definition}

Let us recapitulate what we have established on $\mathrm{Bl}(M,\mathcal{F})\subset \mathrm{Grass}_{-r}(A)$ and its natural projection on $M$. Firstly, it satisfies two properties resolutions have to satisfy, namely 
\begin{enumerate}
\item the map $\pi$ is proper - in particular, it is surjective, 
\item and $\pi^{-1}(M_{\mathrm{reg}}) $ is a manifold which is isomorphic\footnote{i.e., diffeomorphic, diffeomorphic through real analytic maps, or biholomorphic depending on the context.} to $M_{\mathrm{reg}} $
\end{enumerate}
 Secondly, despite its extremely pleasant two properties, the closed subset $\mathrm{Bl}(M,\mathcal{F})\subset \mathrm{Grass}_{-r}(A)$ has a major problem: {\textbf{it is not a sub-manifold in general}}. 

 {\small{
 \begin{remark}
 In \cite{Ruben3}, it is shown that $\mathrm{Bl}(M,\mathcal{F}) $ is obtained by a locally monoidal transformation, whose center lies within the closed subset of points in $M$ that admit no neighborhood on which $ \mathcal F$ is Debord (a subset of $ M_{sing}$). Also, smoothness is addressed.
\end{remark}
}}

\subsubsection{Some natural vector bundles}

Proposition \ref{thm:prop:NashWellDefined} has given the space on which the Nash blowup will be defined. We now have to equip it with a singular foliation. Before doing so, we have to describe two exact sequences of vector bundles over $\mathrm{Gr}_{-r}(A)$, for $A \to M$ a vector bundle and $r$ an integer $ \leq {\mathrm{rk}}(A)$.

\begin{enumerate}
\item Let $ \Pi^! A$ be the pull-back of $A \to M$ to $ \mathrm{Grass}_{-r}(A)$.
\item We call \emph{tautological subbundle} the vector subbundle of $\Pi^! A $ whose fiber at a point  $V_x\in\mathrm{Grass}_{-r}(A)$ is precisely the vector space $V_x$, seen a subspace of $\left(\Pi^! A\right)_{V_x}=A_x$. We denote it by $\tau_A^{-r} $ (the $ -r$ being a reminder that it is of corank $r$ in $\Pi^! A$).
\item  We denote by $A_{Bl}^r$  and call \emph{tautological quotient bundle} the bundle over $\mathrm{Grass}_{-r}(A)$ obtained by taking the quotient of the first bundle by the second one, i.e., $ A_{Bl}^r := \Pi^! A / \tau_A^{-r}$. (The ``$r$'' reminds that it has rank $r$).
\end{enumerate}
These three vector bundles fit into the exact sequence
$$\xymatrix{
      \tau_A^r\ar@{^(->}[r] \ar[dr]&\ar[d] \Pi^! A \ar@{->>}[r] &  A_{Bl}^r \ar[dl] \\  & \mathrm{Grass}_{-r}(A) & 
 }$$

{\small{
\begin{remark}
\label{eq:dualexact}
There is of course a dual canonical isomorphism and inclusion of vector bundles:
 \begin{equation}
 \label{eq:annihilator}
 \xymatrix{
\left(A_{Bl}^r\right)^* \ar[dr] \ar@{-}[r]^{\sim}&   \ar[d]\left( \tau_A^r \right)^\perp \ar@{^(->}[r] & \Pi^! A^*\ar[dl] \\ & \mathrm{Grass}_{-r}(A) & }
 \end{equation}
 where the symbol $ {}^\perp$ stands for the annihilator.
\end{remark}
}}
\subsubsection{Nash-blowup II: the Debord foliation}

We now explain how to lift the singular foliation $\mathcal F $ to a singular foliation on $\mathrm{Bl}(M,\mathcal{F})$. Let $ (A,\rho)$ be an anchored bundle over $\mathcal F$. Lemma \ref{lem:FlowOfLinear}, in Section \ref{sec:almost} implies that for every $X \in \mathcal F$ there exists a linear vector field $\xi_X$ on $A$ (by choosing an almost Lie algebroid bracket on 
$(A, \rho)$) fulfilling the following two properties
 \begin{enumerate}
\item 
 $p$-projectable\footnote{I.e. $p$-related with $X$.} to $X$ 
 \item and whose flow $\phi^{\xi_X}_t$ is an isomorphism of anchored bundles whenever defined 
\begin{equation}\label{eq:linearvf+anchor}
\xymatrix{A\ar[r]^{\phi_t^{\xi_X}}\ar[d]_{\rho}&A\ar[d]^{\rho}\\TM\ar[r]^{T\phi_t^X} &TM.}
\end{equation}
\end{enumerate}
Using Exercise \ref{exo:grass+vf}, we see that $\xi_X$ induces a vector field  on $\mathrm{Grass}_{-r}(A)$ which is $\Pi$-projectable{I.e. $\Pi$-related (see Definition \ref{def:related}) to $X$.} on $X$. 

We call $\hat{\xi}_X $ a \emph{lift of $X$ on $ \mathrm{Grass}_{-r}(A)$}.

\begin{lemma}(Lift of vector fields of $\mathcal{F}$)
\label{grass_pullback}
Consider the lift $ \hat{\xi}_X$ of some $ X \in \mathcal F$ to the Grassmann bundle $\mathrm{Grass}_{-r}(A)$ of some anchored bundle $(A,\rho) $. 
\begin{enumerate}
\item the flow of $ \hat{\xi}_X$ of $X$, whenever it is defined, preserves the subset\footnote{Recall that we erased the index "$A$" from the notation $\mathrm{Bl}_A(M,\mathcal{F})$, since it does not depend on the anchored bundle $(A,\rho)$ by Proposition \ref{thm:prop:NashWellDefined}, but, when an anchored bundle $ (A,\rho)$ is given, we still see it as a subset of  $\mathrm{Grass}_{-r}(A)$.} $\mathrm{Bl}(M,\mathcal{F}) \subset \mathrm{Grass}_{-r}(A)$, i.e., 
    \begin{equation}
    \xymatrix{\mathrm{Bl}(M,\mathcal{F})\ar[r]^{\phi_t^{\hat{\xi}_X}}\ar[d]_{\pi}&\mathrm{Bl}(M,\mathcal{F})\ar[d]^{\pi}\\M\ar[r]^{\phi_t^X}&M}
\end{equation}
\item In particular, $\hat{\xi}_X$ is tangent to $\mathrm{Bl}(M,\mathcal{F})$ in the neighborhood of any point where $\mathrm{Bl}(M,\mathcal{F})$  is a sub-manifold. 
\end{enumerate}
 \end{lemma}
 \begin{proof}
The second item in the lemma is a direct consequence of the first one. Let us prove the first one. For any $x $ is the regular part of $\mathcal F $, it follows from \eqref{eq:linearvf+anchor}
$$\phi_t^{\xi_X}|_{x}\left(\mathrm{\ker}{\rho}_{x}\right)=\mathrm{ker}{\rho}_{\phi^X_t(x)}.$$
Equivalently, in terms of the section $\sigma_{ker}: M_{\mathrm{reg}} \to   \mathrm{Grass}_{-r}(A)$, it means that
 $$ \phi^{\hat{\xi}_X}_t \circ \sigma_{\mathrm{ker}} = \sigma_{\mathrm{ker}} \circ \phi^X_t    .$$
 In particular, since $\phi^X_t(M_{\mathrm{reg}})\subseteq M_{\mathrm{reg}}$, this implies that $  \phi^{\hat{\xi}_X}_t $ preserves the closure $ \overline{\sigma_{\mathrm{ker}}(M_{\mathrm{reg}})}$, i.e., it preserves $\mathrm{Bl}(M,\mathcal{F})$.
 Let us spell out this argument: for any element $V_x\in\mathrm{Bl}(M,\mathcal{F})$, there exists a sequence $ (x_n)_{n \in \mathbb N}$ in $M_{\mathrm{reg}}$ such that $x_n\underset{n \to +\infty}{\longrightarrow}x$ and such that $\mathrm{ker}{\rho}_{x_n}\underset{n \to +\infty}{\longrightarrow} V_x$.  In view of Equation \eqref{eq:linearvf+anchor},  one has  $$\phi_t^{\xi_X}|_{x_n}\left(\mathrm{\ker}{\rho}_{x_n}\right)=\mathrm{ker}{\rho}_{\phi^X_t(x_n)},\quad \text{for every $n\in\mathbb{N}_0$}.$$ 
The limit  $ \phi_t^{\hat \xi_X} (V_x) $ of the sequence $$ n \mapsto
\phi_t^{\hat \xi_X} \left(\mathrm{\ker}\rho_{x_n}\right) =
 \phi_t^{\xi_X}|_{x_n} \left(\mathrm{\ker}\rho_{x_n}\right) $$
 therefore belongs to $ \mathrm{Bl}(M,\mathcal{F})$  by construction. This completes the proof.
 \end{proof}

\vspace{0.5cm}
From now on, we assume for the sake of simplicity that $\mathrm{Bl}(M,\mathcal{F})$ is a smooth submanifold. However, we claim that this assumption is not required for most  statement below, provided that a correct notion of singular foliation on a singular subset is given: we refer to \cite{louis2024nash} for a detailed study of the matter.

\vspace{.5cm}

\begin{proposition}
\label{prop:NashBlowUp}
Let $(M,\mathcal F) $ be a foliated manifold such that 
\begin{enumerate}
\item[a)]
all regular leaves have the same dimension\footnote{Again, this is always true if $M$ is connected in the complex and real analytic settings.} and \item[b)] $\mathrm{Bl}(M,\mathcal{F})$ is a sub-manifold\footnote{In \cite{louis2024nash}, a framework is given to make sense of the three statements of the theorem even when this second assumption is not satisfied.} of  $\mathrm{Grass}_{-r}(A)$ for some anchored bundle $(A,\rho)$. 
\end{enumerate}
Then: 
\begin{enumerate}
    \item Every vector field $X\in \mathcal{F}$, there exists a unique vector field $\pi^! (X) $ on $\mathrm{Bl}(M,\mathcal{F})$ which is $\pi $-related to $X$;
    \item The vector fields  $\{ \pi^! (X) , X \in \mathcal F\} $ generate a singular foliation. The latter shall be denoted by $\pi^!(\mathcal F)$.

    \item The restriction of $ \pi$ to $\pi^{-1}(M_{\mathrm{reg}}) $ is an isomorphism of singular foliation from $( \pi^{-1}(M_{\mathrm{reg}}),\pi^! (\mathcal F))$ to $(M_{\mathrm{reg}}, \mathcal F)$. 
    
\end{enumerate} 
\end{proposition}  
\begin{proof}
Uniqueness in item 1 is a direct consequence of the fact that $ \pi$ is an isomorphism when restricted to the dense open subset $\pi^{-1}(M_{\mathrm{reg}}) $. Existence is an immediate consequence of the second item Lemma \ref{grass_pullback}: the restriction of $\hat{\xi}_X $ to ${\mathrm{Bl}}(M,\mathcal F) $ is the required vector field.
Now, since $ \pi$ is an isomorphism when restricted to the dense open subset $\pi^{-1}(M_{\mathrm{reg}}) $, it is obvious that:
 $$ [\pi^! (X), \pi^{!}(Y)] = \pi^{!}\left([X,Y]\right) \, \, \hbox{ for all $X,Y \in \mathcal F$} $$
 The sheaf considered in item 2 is therefore stable under Lie bracket: it is also obviously locally finitely generated. Item 3 is obvious.
\end{proof}

This statement allows to make sense of the following definition. In the complex case, it does correspond to Nash blowup of coherent sheaves\footnote{Holomorphic singular foliations are instances of coherent sheaves, see Proposition \ref{prop:coherent}.}  in \cite{zbMATH04118703}, hence the chosen name. We will explain later on how it relates with Mohsen's groupoid in the smooth setting.

\vspace{0.5cm}

\begin{definitions}{The Nash-blowup}{nashblow}
The singular foliation $(\mathrm{Bl}(M,\mathcal{F}), \pi^!({\mathcal{F}}))$ shall be referred to as the \emph{Nash blowup of $(M, \mathcal{F})$}.
\end{definitions}

\vspace{0.5cm}

Here is now the most surprising point about the Nash blowup: the singular foliation that we obtain is Debord, i.e., it is isomorphic, as a sheaf, to the sheaf of sections of a vector bundle, which is of course $ A^r_{Bl}$ restricted to $\mathrm{Bl}(M,\mathcal{F})$. Equivalently, it is given by a Lie algebroid whose anchor map is injective on a dense open subset (see Section \ref{examples-of-SF}).
This is established in items 1 and 2 of the next theorem.

\vspace{0.5cm}

\begin{theorems}{The Nash-blow up is a Debord foliation}{NashIsProjective}
Let $(M,\mathcal F) $ be a foliated manifold such that 
\begin{enumerate}
    \item[a)] all regular leaves have the same dimension\footnote{Again, the first assumption is always valid in the complex and real analytic settings if $M$ is connected.} and 
    \item[b)]  $\mathrm{Bl}(M,\mathcal{F})$ is a manifold\footnote{Again, \cite{Ruben3, louis2024nash} give a framework that allows to make sense of the conclusion of the theorem even when this second assumption is not satisfied.}
\end{enumerate}
Then  
\begin{enumerate}
\item
The Nash blowup $ \pi^! (\mathcal F)$ is a Debord singular foliation on $\mathrm{Bl}(M,\mathcal{F})$.
\item
Its associated Lie algebroid is the restriction to $ \mathrm{Bl}(M,\mathcal{F})$ of the canonical quotient bundle $ A_{Bl}^r$. In equation:
 $$\pi^!(\mathcal{F})\simeq \Gamma({A^r_{Bl}}_{|_{\mathrm{Bl}(M,\mathcal{F})}}).$$
     \item In the smooth setting, it is the Lie algebroid of Mohsen's Lie groupoid\footnote{Item 3 in the  theorem does not make much sense yet, since Mohsen's Lie groupoid of a singular foliation has not been defined. It is a topological groupoid defined in \cite{OmarMohsen}, which is Lie when $ \mathrm{Bl}(M,\mathcal{F})$ is a manifold. } of $ \mathcal F$. In the complex setting, as a vector bundle, it matches the vector bundle described in \cite{zbMATH04118703}, section IV.
 \end{enumerate}

\end{theorems}

\vspace{0.5cm}

In view of the second item, we can make sense of the following definition:

\begin{definition}
    \label{def:Nashblowup}
    We call the Lie algebroid in item 2 of Theorem \ref{thm:NashIsProjective} the Nash blow up Lie algebroid, and we denote it by $ A_{Bl}$.  
\end{definition}

We will only prove items 1 and 2 of this theorem in the smooth setting only, and leave it to the reader to adapt to real analytic or complex settings. We refer to \cite{Ruben3,louis2024nash} for item 3. We start with a lemma.

\begin{lemma}
\label{lem:anchorGrass}
There exists a vector bundle morphism $$\hat{\rho}\colon\Pi^!A\longrightarrow T(\mathrm{Grass}_{-r}(A)) $$ 
over $\mathrm{Grass}_{-r}(A)$ such that
 \begin{enumerate}
 \item[a.] the following diagram commutes:
$$ \xymatrix{\Pi^!A \ar[r]^<<<<{\hat\rho} \ar[d]^{}& \ar[d]^{T\Pi}   T \mathrm{Grass}_{-r}(A)) \\A \ar[r]^{\rho} & TM    } . $$
  \item[b.] and for any $a \in \Gamma(A)$, $\hat{\rho}(\Pi^!a) \in \mathfrak X\left(\mathrm{Grass}_{-r}(A)\right)$ is a lift of $X=\rho(a)\in \mathcal{F}$ (induced by a linear vector field $\xi_X \in \mathfrak X(A)$ as in Lemma \ref{grass_pullback}). 
 \end{enumerate} 
\end{lemma}
\begin{proof}
Let $V_x \in  {\mathrm{Grass}}_{-r}(A)$ be an element in the fiber of $x\in M$. Consider $a_1, \dots, a_n$ a local trivialization of $A$ in a neighborhood $\mathcal U $ of $x$. The value at $V_x\in {\mathrm{Grass}}_{-r}(A)$ of the lifts $\hat{\xi}_{a_1}, \dots,\hat{\xi}_{a_i}:= \hat{\xi}_{{\rho}(a_i)}, \cdots \hat{\xi}_{a_n}$ constructed as in lemma \ref{grass_pullback} define a vector bundle morphisms $\hat{\rho}_{\mathcal U} $ which satisfies both properties above by construction. Using a partition of unity $ (U_i,\chi_i,\rho_i)_{i \in I}$, this construction can be made global: it suffices to check that the gluing $ \sum_i \Pi^* \chi_i \,  \hat{\rho}_i $  still satisfies the required properties if each $\hat{\rho}_i $ defined as before does.
\end{proof}

We can now prove the two first items of the theorem.

 \begin{proof}
By the second item in Lemma \ref{lem:anchorGrass},
for every $X \in \mathcal F$, and for every $a \in \Gamma(A)$  such that $\rho(a)=X $, the vector field $\hat{\rho}(\Pi^! a)$ comes from a linear vector field as in Lemma \ref{grass_pullback} (here $\Pi^!a $ is the pull-back section of $a$). As a consequence, by construction
\begin{enumerate}
\item the vector field $\pi^! X$ 
\item  the restriction of  $\hat{\rho}(\Pi^! a) $ to $\mathrm{Bl}(M,\mathcal{F})$
\end{enumerate}
coincide on ${\mathrm{Bl}(M,\mathcal{F})}$. 
In other words, 
$$\hat{\rho}_{|_{\mathrm{Bl}(M,\mathcal{F})}} (\Pi^! a) = \pi^! X .$$
This implies that the pair
 $$  \left(\Pi^!A|_{\mathrm{Bl}(M,\mathcal{F})},\; \hat{\rho}_{|_{\mathrm{Bl}(M,\mathcal{F})}}  \right)   $$
is an anchored bundle over $\pi^! (\mathcal F)  $.
Let us check that the tautological sub vector bundle $\tau_A^{-r}$ lies in the kernel of the anchor map $\hat\rho$, restricted to ${\mathrm{Bl}(M, \mathcal{F})}$. 
This follows from the fact that $\ker\hat{\rho}_{V_x}\subseteq \ker\rho_x$ for all $x\in M$ by its very definition, and that both spaces coincide on $M_{\mathrm{reg}}$ since $\dim (\mathrm{im}\hat{\rho}_{V_x})=\dim(\mathrm{im}{\rho}_{x})$ for all $x\in M_{\mathrm{reg}}$.

Hence, the anchor map $\hat{\rho}_{|_{\mathrm{Bl}(M,\mathcal{F})}}$ goes to the quotient to define a map as follows:
\begin{equation}
 \xymatrix{0\ar[r] &\tau_A^{-r}\ar[r]&\Pi^!A\ar[r]\ar[d]^{} &A_{Bl}^r\ar[r]\ar@{-->}[ld]&0\\& &T(\mathrm{Bl}(M,\mathcal{F})) & &.}
 \end{equation} This map is injective on $M_{\mathrm{reg}} $ with image $T\mathcal F $.    Therefore, $\pi^!(\mathcal{F})$ is a projective singular foliation and is the image of a Lie algebroid bracket on ${A_{Bl}^r}_{|_{\mathrm{Bl}(M,\mathcal{F})}}$. This finishes the proof of the two first items.
\end{proof}

We finish this section by presenting several exercises, describing examples and properties.

\begin{exo}
Show that if $(M,\mathcal{F})$ is a Debord singular foliation (in particular a regular foliation) then  $\mathrm{Bl}(M,\mathcal{F})\simeq M$  and  $\pi^! ({\mathcal{F}})\simeq\mathcal{F}$.
\end{exo}

\begin{exo}
    Show that 
    \begin{enumerate}
        \item $\pi^!\mathcal{F}$ is a regular foliation of rank $r$ if and only if $\tau^r_A=\ker(\hat{\rho})$, i.e., $V_x=\ker(\hat{\rho}_{V_x})$ for all $x\in M$. In particular, if $\mathrm{Sker}(\hat{\rho},x)\neq 0$ for all $x\in M$ and  $\tau^r_A$ is a line bundle, then $\pi^!\mathcal{F}$ is regular.
        \item  $\mathrm{Sker}(\rho,x)=V_x$ for all $V_x\in \pi^{-1}(x)$ with $x\in M$, if only if $\mathcal{F}$ is a regular foliation.
    \end{enumerate}
\end{exo}

\begin{exo} Show that for $ \mathcal F$ the singular foliation on $\mathbb R^n $ of all vector fields vanishing at the origin, the Nash blowup coincides with the blowup at $0$ as in Section \ref{sec:blowup1}.
\end{exo}

\begin{exo}
Let $M$ be a manifold and a $N \subset M$ submanifold. Consider $\mathcal{F}_N\subset \mathfrak X(M)$ the singular foliation of vector fields that are tangent to $N$. Show that the Nash blowup coincides with the blowup along the submanifold $N$ in this case.

\end{exo}

\begin{exo}
Let $(M=\mathbb R^3, \mathcal{F})$  the singular foliation given by the transformation Lie algebroid of the action of ${\mathfrak{so}}(3) $ on $\mathbb R^3 $. It is generated by the vector fields 
 $$X=y \frac{\partial}{\partial x} - x\frac{\partial}{\partial y},\;\; Y=z\frac{\partial}{\partial x} - x\frac{\partial}{\partial z},\;\; Z=z \frac{\partial}{\partial y}- y\frac{\partial}{\partial z}.$$

 \begin{enumerate}
     \item Show that $\mathrm{Bl}(M,\mathcal{F})$ is the usual blowup of $\mathbb R^3$ at $0$. In particular, $ \mathrm{Bl}(M, {\mathcal{F}})$ is smooth.
     \item Show that  $\pi^! (\mathcal{F})$ is  generated in the $x$-chart by $$\pi^!(X)=xy\frac{\partial}{\partial x}-(y^2+1)\frac{\partial}{\partial y}- yz\frac{\partial}{\partial z}\;\;\text{and}\;\; \pi^!(Y)=xz\frac{\partial}{\partial x}-yz\frac{\partial}{\partial y}- (z^2+1)\frac{\partial}{\partial z}$$
     while $\pi^!(Z)$ still has the same expression.
     Compute the generators on the other charts.
     
     \item  Deduce that $\pi^!(\mathcal{F})$  is a regular foliation on $\mathrm{Bl}(M,\mathcal{F})$ of rank $2$. 
 \end{enumerate}
 \end{exo}

The next exercises will describe more precisely the leaves of the Nash blowup. 
We claim that we never really need the simplifying assumption that $ \mathrm{Bl}(M,\mathcal{F})$ is a submanifold, but we invite the reader to make it for the sake of simplicity.
 
\begin{exo}
Let $ (\mathrm{Bl}(M,\mathcal{F}), \pi^! (\mathcal F)) $ be the Nash blowup of a singular foliation $\mathcal F $, and let $\pi\colon  \mathrm{Bl}(M,\mathcal{F}) \to M$ be the projection.
\begin{enumerate}
\item Show that $\pi \colon  (\mathrm{Bl}(M,\mathcal{F}), \pi^! (\mathcal F))\longrightarrow (M,\mathcal F)$ maps leaves into leaves.
\emph{Hint: Use $\Pi\circ \phi^{\hat{X}_k}_{t_1}\circ\cdots\circ\phi^{\hat{X}_k}_{t_k}= \phi^{{X}_1}_{t_1}\circ\cdots\circ \phi^{{X}_k}_{t_k}\circ \Pi$. }
\item Let us fix $x \in M$, and consider $\pi^{-1}(x) $ to be a subset of $\mathrm{Grass}_{-r+\dim(L_x)}^{Lie}\left(\mathfrak g_x(\mathcal F)\right)  $ as in exercise \ref{exo:def-blup}.
Let $\mathrm G_x(\mathcal F)$ be the simply connected Lie group integrating $\mathfrak g_x(\mathcal F) $. It naturally acts on 
$$ \mathrm{Grass}_{-r+\dim(L_x)}^{Lie}\left(\mathfrak g_x(\mathcal F)\right) .$$
There is therefore a Lie algebra morphism
 $$ \phi\colon \mathfrak g_x(\mathcal F) \longrightarrow \mathfrak X\left(\mathrm{Grass}_{-r+\dim(L_x)}^{Lie}\left(\mathfrak g_x(\mathcal F)\right) \right) ,$$
which in turn yields an action of $ G_x(\mathcal F)$ on 
 $\mathrm{Grass}_{-r+\dim(L_x)}^{Lie}\left(\mathfrak g_x(\mathcal F)\right) $.
Show the following points.
\begin{enumerate}

\item The fiber $ \pi^{-1}(x)$ of $ (\mathrm{Bl}(M,\mathcal{F}) \to M $  is stable under the $\mathrm G_x(\mathcal F)$-action. 
\item We assume that $T_x\mathcal{F}=\{0\}$. Show that the leaf of $\pi^! (\mathcal F) $ through some point in $  \pi^{-1}(x)$ is precisely the $\mathrm G_x(\mathcal F)$-orbit.
\item  Show that a leaf of  $\pi^! (\mathcal F) $ is reduced to a point\footnote{Seen as a codimension $ r$ Lie algebra of $ \mathfrak g_x(\mathcal F) $.} $V \in \pi^{-1}(x)$  if and only if $ T_x \mathcal F=\{0\}$ and $V$ is a Lie ideal of the isotropy Lie algebra $\mathfrak g_x(\mathcal F) $.
\item  Show that a point\footnote{Seen as a codimension $ r-{\mathrm{dim}}(L_x)$ Lie algebra of $ \mathfrak g_x(\mathcal F) $.} $V\in \pi^{-1}(x)$ is a regular point if and only if it is equal to its own stabilizer.
\end{enumerate}

\end{enumerate}
\end{exo}

\begin{exo}
    Assume that $A\rightarrow M$ is a Lie algebroid over $M$.  Show that $A\rightarrow M$ acts  naturally on $\mathrm{Bl}(M,\mathcal{F})\rightarrow M $.
\end{exo}

\subsection{Push-forward}

Let $\phi \colon P \to M$ be a smooth, complex or real analytic  map, depending on the context. We will assume that $\phi$ is a surjective submersion.

The \emph{push-forward} $T_p \phi:T_p P \to T_{\phi(p)}M$ does not extend in general to vector fields: for $X$ 
a vector field on $P$ and $n= \phi (p) = \phi (p')$ with $p \neq  p'\in P$, then $X_p$ and $X_{p'}$ are both pushed forward to tangent vectors at $n\in M$, but in general $T_p \phi (X_p ) \neq T_{p'}\phi(X_{p'})$.  When this happens, we denote this vector field by $\phi_*(X)$  and we call it the \emph{push-forward of $X$ through $\phi$}.\\

Let us introduce a notation: for $\phi \colon P \to M $ a surjective submersion, we denote by $\mathfrak X (P)_\phi $ the space of vector fields $X $ on $P$ which are $\phi$-related to a vector field on $M$, that we denote by $\phi_*(X)$. 

Assume that we are now given a singular foliation $\mathcal F^P $ on $P$. 
Then $\mathcal F^P  \cap \mathfrak X (P)_\phi  $ is both a $\mathcal C^\infty(M) $-module and stable under Lie bracket, and so is 
 $$  \phi_*( \mathcal F^P  \cap \mathfrak X (P)_\phi )  \subset \mathfrak X (M) .$$
 
 When the latter is finitely generated, it is a singular foliation that we call \emph{push-forward} singular foliation and denote by $\phi_* (\mathcal F^P) $.

\begin{exo}
Let $ \mathcal F^M$ be a singular foliation on $M$, and $ \phi \colon P \to M$ a surjective submersion.
Show that the push-forward of the pull-back $ \phi^{-1}(\mathcal F^M)$ (see Section \ref{sec:pull-back}) is the singular foliation  $ \mathcal F^M$.
\end{exo}
 
 \begin{exo}
 Here are examples where 
 $\phi_*( \mathcal F^P  \cap \mathfrak X (M)_\phi) $ is not finitely generated.
 
 Let $ \mathcal F^P$ be  Androulidakis-Zambon's ``non-finitely-many-generators'' singular foliation of exercise \ref{exo:infinitestability}. In this case, we have $ P := \mathbb R^2$. Let $ M:=\mathbb R$, and $\phi \colon (x,y) \mapsto y $ the projection onto the horizontal axis. 
 
\begin{enumerate}
    \item Show that for every $\phi$-projectable\footnote{I.e., $ \phi$-related to a vector field on $ \mathbb R$, see Definition \ref{def:related}.} vector field $X$ on $\mathbb R^2 $ whose derivatives vanish at order $n$ at the point of coordinates $(n,0)$, its projection $\phi_* (X) \in \mathfrak X(\mathbb R)$ is a vector field that vanishes at order $n $ at~$0$. 
    \item Show that $\phi_*( \mathcal F^P  \cap \mathfrak X (P)_p) $  coincides with the space of vector fields on $\mathbb R $ vanishing at $0$ with all their derivatives.
    \item Conclude (\emph{Hint}: use the second item in Exercise \ref{exo:finitelygeneratedex}).
\end{enumerate} 
 \end{exo}

\section{Morphisms of singular foliations} 

\emph{Isomorphisms} of singular foliations are easily defined; they are diffeomorphisms (biholomorphisms in the complex case) of the underlying varieties that intertwine their respective singular foliations: see Section \ref{sec:symmetry}. But defining general morphisms of singular foliations is more involved: we even dare to say that finding a definition that makes consensus is still an open problem. 

\vspace{.2cm}
There is a case, however, for which an easy definition of morphism exists: surjective submersions.

\begin{definition}
\label{def:morphSubm}
Let $(P,\mathcal G) $ and $ (M,\mathcal F)$ be foliated manifolds. A submersion $ \phi \colon P \to M $ is said to be a morphism of singular foliation if $\mathcal G \subset \varphi^{-1}(\mathcal F) $.
\end{definition}

This definition satisfies the following two interesting properties.
\begin{enumerate}
\item If two points are in the same leaf of $\mathcal G$, then their images are in the same leaf of $ \mathcal F$. 
\item For every $ p \in P$, the inclusion $ T_p \phi \left( T_p \mathcal G \right) \subset T_{\phi(p)} \mathcal F$ holds.
\item If $ P=M$ and $ \Phi$ is the identity map, then morphisms are simply inclusions of singular foliations. More precisely, the identity map is a morphism 
 $(M,\mathcal G) $ and $ (M,\mathcal F)$ if and only if $\mathcal G \subset \mathcal F $.
\end{enumerate}

However, it is strange that with this definition, an immersion $S \hookrightarrow M $ could not be a morphism of singular foliation.  For instance,  for sub-manifolds intersecting $S \subset M$ cleanly\footnote{See Section \ref{sec:transverse}.} a singular foliation $\mathcal F_M $ on $M$, we would like the inclusion map to be a morphism of singular foliation from $(S,\mathfrak i^*_S \mathcal F_M )$ to $ (M,\mathcal F_M)$. More generally, for the Androulidakis-Skandalis pull-back  $(L,p^{-1}(\mathcal F_M))  $ of $ (M,\mathcal F_M)$ by  a map $ p \colon L \to M$ transverse to $ \mathcal F$ as in Section \ref{sec:pullback2}, we would like $p$ to be a morphism.
Recall from that section that we say that a smooth map $ p \colon  L \to M$ is transverse to $\mathcal F_M $ if for all $ \ell \in L$
  $$ T_{\phi(\ell)} \mathcal F_M + T_\ell \phi (T_\ell P) = T_{\phi(\ell)} M .$$
  This is enough to define the pull-back $\varphi^{-1}(\mathcal F_M) $. The latter is a singular foliation on $L$, see Section \ref{sec:pullback2}.
This clears the way to our next definition, which is more general than Definition \ref{def:morphSubm}, since any submersion $P \to M $ is transverse to any singular foliation on $M$.

\vspace{.5cm}

\begin{definitions}{Morphisms of singular foliations: the transverse case}{morphtrans}
Let $(L,\mathcal F_L)$ and $(M,\mathcal F_M) $ be foliated manifolds. A map $ \phi \colon L \to M $ is said to be a morphism of singular foliation if
\begin{enumerate}
    \item $\phi$ is a transverse to $\mathcal F_M $,
    \item $\mathcal F_L \subset \phi^{-1}(\mathcal F_M) $. 
\end{enumerate}
\end{definitions}
\vspace{.5cm}

The above notion is constructed exactly such that every map $L \to M $ transverse to $ \mathcal F_M$ is a morphism of singular foliation  $ (L,\phi^{-1}(\mathcal F_M))  $ to $ (M,\mathcal F_M)$.
Again, this definition satisfies several interesting properties.
\begin{enumerate}
\item If two points are in the same leaf of $\mathcal F_L$, then their images through $\phi$ are in the same leaf of $ \mathcal F_M$. 
\item For every $ \ell \in L$, the inclusion $ T_\ell \phi \left( T_\ell \mathcal F_L \right) \subset T_{p(\ell)}\mathcal F_M$ holds.
\end{enumerate}

\begin{exo}
Show that the inclusion of a transverse submanifold $S$ in a foliated manifold $(M,\mathcal F) $ is a morphism in the previous sense for every sub-singular foliation included in the restriction $\mathfrak i_S^* \mathcal F $ (see Section \ref{sec:transverse}).
\end{exo}

\begin{exo}
Show that for the direct product singular foliation $\mathcal F_1\times \mathcal F_2$ on $M=M_1\times M_2$, the projections are morphisms. Moreover, show that the direct product singular foliation is the largest singular foliation on $M$ with this property.     
\end{exo}

\begin{exo}
This exercise requires the notion of Lie algebroid morphism \cite{MR2157566}.
Assume that the base map of a Lie algebroid morphism is a submersion: is it a morphism of singular foliations?
Does a Poisson submersion induce morphism of their respective symplectic foliations? 
\end{exo}

While the above definition is sufficient to treat many important cases, it has weaknesses. For instance, it is strange that for $L$ a leaf in a singular foliation $ \mathcal F$ on $M$ (see Section \ref{sec:leavesExist} for a definition of leaves), the inclusion $ \mathfrak i \colon L \hookrightarrow M$ is never a morphism, unless $L$ is an open subset of $M$. 
Hence a natural question:

\begin{question}
Is Definition \ref{def:morphtrans} the definitive answer to the question of defining "morphisms of singular foliations"?
\end{question}

For the coming lines, we are indebted to a discussion with Hadi Nahari and his advisor, Thomas Strobl. Let $(L,\mathcal F_L)$ and  $(M,\mathcal F_M)$ be singular foliations and $\phi:L\to M$ be a smooth map.
 
\begin{enumerate}
\item The most naive condition to impose on $\phi$ in order to say that $ \phi$ is a morphism would be to ask that $T_m \phi \, (T_m\mathcal F_1)\subset T_{\phi(m)} \mathcal F_2$ for any $m\in M_1$. However, with this notion, the map $t\mapsto t^3$ would be a morphism of singular foliations from $(\mathbb R,\mathfrak X_c(\mathbb R))$ to $(\mathbb R,t\cdot \mathfrak X_c(\mathbb R))$ (= compactly supported vector fields on $ \mathbb R$ that vanish at $0$). This latter map, however, does not preserve leaves and hence should not be called morphism. This is therefore \emph{not} a good definition.
\item 
A more clever idea is to use the vector bundle $ \phi^* TM \to L$ as in Section \ref{sec:pullback2}.
Recall that both $ T\phi(\mathcal F_L)$ and the $ \mathcal C^\infty_c(L) $-module generated by $p^* X $ with $ X \in \mathcal F_M$ belong to the compactly supported sections of that vector bundle. We denote the second one by $ C^\infty_c(M) p^* \mathcal F_M $. A possible definition of morphism of singular foliations would be to say that $ \phi$ is a morphism  if  $  T\phi(\mathcal F_L)\subset \mathcal C^\infty_c(M) p^* \mathcal F_M  $.
\end{enumerate}
Assume for a moment that  we define morphisms using the second item above. Then the condition of the first item would also be satisfied. Also, 
a morphism as in Definition \ref{def:morphtrans} would then still be a morphism in the new sense.
Last, with such a definition, it would still be true that two points in the same leaf $\mathcal F_L $ are mapped to the same leaf of $ \mathcal F_M$.

This notion of morphism is however too weak for certain purposes\footnote{For instance, it does not induce a morphism to the level of isotropy Lie algebras (cf. section \ref{sec:isotropy}) or fundamental groupoids (cf. section \ref{sec:fundamental}).} An alternative definition can be found in Garmendia and Villatoro's \cite{GV}, where morphisms of foliated manifolds are defined as sheaf comorphisms compatible with the Lie bracket.
Hence, we prefer not to take the risk of making this definition into a formal one, and to leave room for more discussion on the matter. We only insist that, whatever definition one chooses, maps as in Definition \ref{def:morphtrans} should be morphisms.

\section{Leaves of a singular foliations} 

\label{sec:leavesExist}

We show in this section that to any singular foliation  is attached a smooth partitionifold.

\subsection{What is a leaf?}

 Let $\mathcal F $ be a singular foliation on a manifold $M$. For the present discussion, we will place ourselves in the context of smooth differential geometry, and consider $  \mathcal F$, as in Definition \ref{def:consensus}, as an involutive $\mathcal C^\infty(M)  $-submodule of compactly supported smooth vector fields on~$M$. The discussion can be easily adapted to the complex and real analytic settings\footnote{But can \emph{not} be adapted to the algebraic setting.}.

\vspace{.5cm}
 
\begin{questions}{What is a leaf?}{ques:whatLeaf}
What are the leaves of $\mathcal F $? And do they exist?
\end{questions}

\vspace{.5cm}

 There are two natural notions of leaves, two different notions that deserve to be called ``leaves''.

\begin{enumerate}
    \item The first idea is that leaves are ``reachable points''. That is, we will define an equivalence relation on $M$ by pairing two points in $M$ such that one can be reached one to the other by following the flow of vector fields in $\mathcal F$. 
    \item But we  may also use the tangent space of $\mathcal F $. A leaf should be  submanifold (by definition !) whose tangent space at a point $m$ is the tangent distribution $T_m\mathcal F $ at that point. 
\end{enumerate}

Here is a formal definition.

\begin{definition}
\label{def:reachable}
Let $\mathcal F \subset \mathfrak X_c(M) $ be a singular foliation.
We say that a point $y \in M$ is \textit{reachable from a point $x \in M$} if there exists:
\begin{enumerate}
    \item a finite sequence $x_0, \dots,x_N $ of points in $M$ with $x_0=x$ and $ x_N =y$
    \item  time-dependent vector fields\footnote{I.e a linear combination of the type $ \sum_{i=1}^r f_i X_i$ with $X_i \in \mathcal F$ and $f_i \in \mathcal C^\infty(M \times \mathbb R) $. These vector fields being zero for all $t$ outside some compact subset of $M$,  their flow at time $t$ are  defined on the whole manifold $M$.}  $(X_t^{(i)})_{t \in \mathbb R} \in \mathcal F$ for $i =0, \dots,N-1 $, with $ X_i$ being defined in a neighborhood of $x_i$ and $x_{i+1}$,
\end{enumerate}
such that for all indices $i=0, \dots, N-1$, the integral curve starting at $x_i$ at time $t=0$ of $X_t^{(i)} $ reaches $x_{i+1}$ at time $ t=1$.
\end{definition}

\begin{exo}
Show that could suppress “time-dependent” in the second item Definition \ref{def:reachable}, and that the $R$-leaves would stay the same.
\end{exo}

If one uses the sheaf definition \ref{def:consensus2}, the definition is easily adapted: one has to assume that the time $1$-flow of $ X_t^{(i)}$ is defined near $ x_i$ for all $ i=1, \dots, N-1$.
The following statement is obvious.

\begin{proposition}
The relation on $M$ defined by $ x \sim y$ if $y$ is reachable from $x$ is an equivalence relation.
\end{proposition}

We call {\emph{reachable leaves}} or \emph{R-leaves} for short the equivalence classes of the previous relation.

But there is a second natural definition of what a leaf should be.

\begin{definition}
A \emph{tangent-leaf}, or \emph{T-leaf} for short, is a connected submanifold $L \subset M$ such that for every $ \ell \in L$, 
 $$T_\ell L = T_\ell \mathcal F ,$$
 and which is maximal among connected sub-submanifolds that satisfy the same property\footnote{I.e., it cannot be strictly included in a submanifold that satisfy the same property. }.
\end{definition}

While defining singular foliation in Definition \ref{def:consensus}, we assumed “locally finitely generated”. 
$ R$-leaves and $ T$-leaves could be defined for any sub-module $ \mathcal F \subset \mathfrak X_c(M)$ stable under Lie bracket, even if they do not satisfy the “locally finitely generated” assumption. The next exercise shows that without this assumption, $R$-leaves and $T$-leaves are different concepts. 

\begin{exo}
\label{exo:comblast}
"\emph{The infinite comb (revisited) (after \ref{exo:comb} and \ref{exo:infinitecomb})}"
Let $ M:= \mathbb R^2$ be the Cartesian plane with coordinates $(x,y)$. 
Let $\mathcal I_{-} \subset \mathcal C^\infty(\mathbb R^2)$ be the ideal of functions vanishing identically on $\mathbb R_- \times \mathbb R $. 

Consider all vector fields of the form
 $$  \mathcal F_{comb} = \left\{  f(x,y) \frac{\partial}{\partial x} + g(x,y)  \frac{\partial}{\partial y}  \middle| g \in \mathcal C^\infty(\mathbb R^2),  f \in \mathcal I_{-}  \right\}  $$ 

\begin{enumerate}
    \item Show that $\mathcal F_{comb} $:
     \begin{enumerate}
         \item is stable under multiplication by $\mathcal C^\infty(\mathbb R^2)$,
         \item is involutive, i.e., is closed under the Lie bracket of vector fields:
          $$ [\mathcal F_{comb}, \mathcal F_{comb}] \subset \mathcal F_{comb}.$$
     \end{enumerate}

    \item Draw what vector fields in $\mathcal F_{comb} $ look like.
    \item Determine $T _{x,y} \mathcal F_{comb} $ for all $ (x,y)$ (\emph{Hint}: it depends on the sign of $x$).
    \item Show that any point in $M=\mathbb R^2 $ is reachable from any point in $M=\mathbb R^2 $. How many $R$-leaves exists? 
    \item Does $\mathcal F_{comb} $ admit $T$-leaves?

\end{enumerate}
\end{exo}

\vspace{0.5cm}

\begin{definitions}{Definition of leaves}{leaves}
A \emph{leaf} of an involutive distribution $\mathcal F \subset \mathfrak X_c(M)$ is a submanifold $L \subset M$ which:
\begin{enumerate}
    \item[(i)] is a T-leaf,
    \item[(ii)] and a R-leaf.
\end{enumerate}
\end{definitions}

\vspace{.5cm}

Here is the main result of this section, which is attributed to Hermann \cite{Hermann}.

\vspace{.5cm}

\begin{theorems}{Hermann: Singular foliations do admit leaves!}{Hermann}
Every singular foliation on a smooth manifold $M$ partitions $M$ into leaves.
\end{theorems}

\vspace{.5cm}

Here is an immediate consequence of this theorem. Any leaf of $ \mathcal F$ is an immersed submanifold, since so are $T$-leaves by definition.  
Leaves partition $M$, since $R$-leaves partition $M$ par definition. Leaves form therefore a partitionifold. It is moreover a smooth partitionifold, since the tangent space of the leaf through a point is $ T_m \mathcal F$ by the definition of $T$-leaves. 
Hence, every singular foliation induces a smooth partitionifold on $M$ that we will denote by $L_\bullet $. By definition, for every $m \in M$, $L_m $ is the set of reachable points from $m$, and it satisfies that $T_m \mathcal F = T_m L_m$.

Here is an even more precise statement than  Theorem \ref{thm:Hermann}. It is the one that we will indeed prove, and it immediately implies Theorem \ref{thm:Hermann}. 

\vspace{.5cm}

\begin{theorems}{Second version}{Hermann2} 
Let $\mathcal F $ be a singular foliation on a smooth manifold $M$.
Every R-leaf $L$ 
        is a (maybe immersed) submanifold of $M$,
           whose tangent space $T_\ell L$ coincides with $T_\ell \mathcal F  $ at every $ \ell \in R$. Lastly, leaves form a smooth partitionifold of $M$.
\end{theorems}

We rest of the present section is mainly dedicated to the proof of this statement.
In Section \ref{self:symmetry}, we will prove that the flow of a vector field in $ \mathcal F$ is a symmetry of $ \mathcal F$, provided that it exists. In Section
\ref{sec:splitting}, we will prove an important splitting theorem explaining the local structure of a singular foliation. Only then, we will be able to address Theorem \ref{thm:Hermann2} in Section \ref{proofHermann}.

\subsection{A singular foliation is a symmetry of itself}  
\label{self:symmetry} 

The first step to prove Theorems \ref{thm:Hermann} and \ref{thm:Hermann2} is to prove that vector fields in a singular foliation have flow which are infinitesimal symmetries of themselves.
The arguments presented in this section are elementary, but quite complicated. Much better conceptual arguments proving the same results will be given using the notion of anchored bundle and almost Lie algebroids.

This is actually a particular instance of the following more general statement, see e.g. \cite{zbMATH07106151}.

\vspace{.5cm}

\begin{propositions}
{The flow of an infinitesimal symmetry is a symmetry}{prop:symInt}
Choose $ t \in \mathbb R$.
Let $Y \in \mathfrak X(M)$ be a vector field whose time $t$-flow $\varphi^Y_t \colon M \to M $ exists. If\footnote{That is, if $Y$ is an ``infinitesimal symmetry'' of $\mathcal F $.} $[Y,\mathcal F]\subset \mathcal F$, then~$\varphi^Y_t$ is a symmetry of $\mathcal F$. 
\end{propositions}

\vspace{.5cm}

In fact, we are going to prove a more general result.

\begin{proposition}
\label{prop:symInt}
Let $Y$ be a vector field such that $ [Y,\mathcal F] \subset \mathcal F$. For every open neighborhood $ \mathcal U$ on which $\mathcal F $ is generated by vector fields $X_1,\dots, X_r $, and any $\mathcal V \subset \mathcal U$ an open subset such that $\phi_t^Y (x) $ exists and takes values in $\mathcal U $ for all $x \in \mathcal V $ and $ |t|\leq \epsilon $, there exists a matrix $\mathbf A(t,x) $, whose coefficients are functions on $\mathcal V $ depending on $t $ such that for all $ i=1, \dots, r$:
$$ (\phi_t^Y)_* \left(  \begin{array}{ccc}X_1 \\ \vdots \\ X_r \end{array} \right) =   \left(  \begin{array}{ccc} && \\ &\mathbf A(t,x) & \\ & &  \end{array} \right) \left(  \begin{array}{ccc}X_1 \\ \vdots \\ X_r \end{array} \right)  $$
Moreover, we can assume that
\begin{equation}\label{eq:chainA}  \mathbf A(s,\phi_t^Y(x)) \circ \mathbf A(t,x) = \mathbf A(t+s,x)  \end{equation}
for all $s,t \in \mathbb R,  x \in \mathcal V $ for which $\phi_t^Y(x) \in \mathcal V$, and $|t|,|t+s| \leq \epsilon $.
\end{proposition}

\begin{proof}
Consider an open neighborhood $ \mathcal U$ of a point $ m \in M$ on which $\mathcal F$ is generated by $X_1,...,X_n$.
Let us chose $\epsilon>0 $ and a smaller neighborhood $\mathcal V \subset U $ such that if $ |t| \leq \epsilon $, $ \phi_t^Y (\mathcal V) \subset \mathcal U$.  
By definition of a symmetry of a singular foliation, there exist  smooth functions $b_i^j \in \mathcal C^\infty(\mathcal U)$, such that $[Y,X_i]=\sum_{j=1}^r b_i^j \, X_j$.
Let us write this expression as a matrix:
 \begin{equation}\label{eq:defad} {\mathrm{ad}}_Y \left( \begin{array}{c} X_1  \\ \vdots   \\ X_r   \end{array}\right)  =  \left( \begin{array}{ccc} & &  \\ & {\mathbf{ad_Y}}(x)&  \\ & &   \end{array}\right) \,  \left( \begin{array}{c} X_1  \\ \vdots   \\ X_r   \end{array}\right) \end{equation}
 with ${\mathbf{ad_Y}} $ being a shorthand for the matrix of functions on $\mathcal U $ whose $i$-th line and $j$-th column is $ b_j^i $.

 For any diffeomorphism $ \phi \colon \mathcal V \to \phi(\mathcal V)$,  the push-forward map
  $\phi_* \colon \mathfrak X( \phi(\mathcal V)) \simeq \mathfrak X(\mathcal V)$ is
defined by $\phi_* (X)|_m = T_{\phi(m)}\phi^{-1} (X_{\phi(m)})  $. It satisfies for all $F \in \mathcal C^\infty(\phi(\mathcal V))$ and 
 $X \in \mathfrak X (\phi(\mathcal V)) $ the relation:
 \begin{equation}\label{eq:pushforward} \phi_* (FX)= \phi^* F \, \, \phi_*(X)\end{equation} 
 Also, if $\phi=\phi_t^Y $ is the flow of $Y$ at time $t$:
  \begin{equation}\label{eq:derPhiY}  \frac{\partial}{\partial t}(\phi_t^Y)_* X =(\phi_t^Y)_* [Y,X] = [ Y ,  (\phi_t^Y)_* X]\end{equation}
 
We want to show that there exist time-dependent functions  $A_i^j(t,x)$ on $\mathcal V $ such that 
\begin{align}
    (\phi^Y_{t})_*(X_i)=\sum_{j=1}^r A_i^j(t) \, X_j
\end{align}
where $(\phi^Y_{t})_*(X_i)$ is to be understood as the image through the push-forward map of the restriction of $ X_i$ to $\phi_t^Y(\mathcal V) $. We also want the matrix of functions $(A_i^j(t,x)) $ to be invertible for all $t,x$.

Again, let us write the expression we wish to obtain in a matrix form. Below, both sides are column vectors of vector fields on $\mathcal V $:
 \begin{align}\label{eq:flowansatz} (\phi_t^Y)_*  \left( \begin{array}{c} X_1  \\ \vdots   \\ X_r   \end{array}\right)  =  \left( \begin{array}{ccc} & &  \\ & \mathbf A(t,x)&  \\ & &   \end{array}\right) \,  \left( \begin{array}{c} X_1  \\ \vdots   \\ X_r   \end{array}\right) \end{align}
 with $\mathbf A (t,x)$ being a shorthand for the matrix of functions on $\mathcal V $ whose $i$-th line and $j$-th column is $ A_j^i(t,x) $. Consider the initial value problem with parameters $x \in \mathcal V $:
 $$ \frac{ \partial A_i^j(t,x)}{\partial t}= \sum_{k= 1}^r b_i^k(\phi^X_{t}(x)) \, A_k^j(t,x)$$  with initial conditions $A_i^j(0,x)=\delta_{i,j}$. Or, equivalently, consider the initial value problem on the vector space $r \times r$ matrices:
  \begin{align}\label{eq:IVB}\frac{\partial}{ \partial t} \left( \begin{array}{ccc} & &  \\ & \mathbf A(t,x)&  \\ & &   \end{array}\right)  = \left( \begin{array}{ccc} & &  \\ & {\mathbf{ad_Y}}\left(\phi_{t}^Y (x)\right)&  \\ & &   \end{array}\right)  \left( \begin{array}{ccc} & &  \\ & \mathbf A(t,x)&  \\ & &   \end{array}\right)  \end{align}
  with initial condition $\mathbf  A(0,x)= {\mathrm{id}} $. The initial value problem have solutions for all $ x \in \mathcal V$ and $|t|\leq \epsilon$, upon changing $ \mathcal V$ for a smaller neighborhood that we still call $\mathcal V $ if necessary. Those solutions depend smoothly on the parameters $ x \in \mathcal V$. Also, the matrix $ \mathbf A(t,x)$ is invertible for all $|t| \leq \epsilon$ and $x \in \mathcal V$.
Last, as any differential equation, it satisfies  \eqref{eq:chainA}.
  
  We claim that Equation \eqref{eq:flowansatz} holds. To show it, let us introduce the column vector whose components are vector fields on $\mathcal V$:
   $$ R(t,x) = \left( \begin{array}{ccc} & &  \\ & \mathbf A(t,x)&  \\ & &   \end{array}\right)^{-1} \hspace{.3cm} \circ \hspace{.3cm}   (\phi_t^Y)_* \hspace{.3cm}   \left( \begin{array}{c} X_1  \\  \vdots   \\ X_r   \end{array}\right)   $$
 An easy computation gives (we now abbreviate the matrix notations, also $ (X_\bullet)$ stands for the column vector $X_1, \dots, X_r $):
  \begin{align*} \frac{\partial R(t,x)}{\partial t} = - \mathbf A^{-1} \circ \frac{\partial \mathbf A}{ \partial t}  \circ \mathbf A^{-1}  \circ (\phi_t^Y)_* (X_\bullet) + \mathbf A^{-1} \circ \underbrace{(\phi_t^Y)_* \circ {\mathrm{ad}}_Y }_{\hbox{by Eq. \eqref{eq:derPhiY}}}(X_\bullet)  \\=
   - \mathbf A^{-1}  \circ \frac{\partial \mathbf A}{ \partial t}  \circ\mathbf A^{-1} \circ (\phi_t^Y)_* (X_\bullet) + \mathbf A^{-1} \circ (\phi_t^Y)_* \circ \underbrace{{\mathbf{ad_Y}}(x)}_{\hbox{by Eq. (\ref{eq:defad})}} (X_\bullet) \\=
    - \mathbf A^{-1} \circ \frac{\partial \mathbf A}{ \partial t} \circ \mathbf A^{-1}  \circ (\phi_t^Y)_* (X_\bullet) + \mathbf A^{-1} \circ \underbrace{{\mathbf{ad_Y}}(\phi_{t}^Y (x) ) \circ  (\phi_t^Y)_*  }_{\hbox{by Eq. \eqref{eq:pushforward}}} (X_\bullet)  \\  =
   \mathbf A^{-1} \circ \underbrace{\left( - \frac{\partial \mathbf A}{ \partial t} \circ \mathbf A^{-1}+ {\mathbf{ad_Y}}(\phi_{t}^Y (x) ) \right) }_{ \hbox{$=0$ by Eq. \eqref{eq:IVB}} } \circ (\phi_t^Y)_*   (X_\bullet) \\
    =0 
  \end{align*}
  Since $ R(0,x) = (X_\bullet)$, we have $R(t,x) = (X_\bullet) $ for all $t \leq \epsilon $ and the \eqref{eq:flowansatz} follows. This implies that the push-forward of any vector field in $ \mathcal F$
 under the flow of $Y$ is a vector field in $ \mathcal F$ at least for $t$ small enough. Composing such push-forward maps, we obtain that it is still  true for all $t$ such that the flow of $Y$ is well-defined.  

\end{proof}

Let us restate Proposition \ref{prop:symInt} differently. We call \emph{infinitesimal symmetry of $ \mathcal F$} a vector field $Y$
 such that $ [Y,\mathcal F] \subset \mathcal F$ (in contrast with \emph{symmetry of $ \mathcal F$} which are diffeomorphisms such that $ \phi_*(\mathcal F)=\mathcal F$).
 
 \vspace{.5cm}
 
 \begin{propositions}{Symmetries and infinitesimal symmetries}{flowofinner}
 When the flow of an infinitesimal symmetry of $\mathcal F $ exists, it is a symmetry of $ \mathcal F$.
 \end{propositions}

\vspace{.5cm}
 
 Again, a much better proof will be given using the notion of anchored bundle and almost Lie algebroids.
Proposition \ref{thm:flowofinner} has several immediate and very important corollaries.

\begin{corollary}
\label{coro:flowissymmetry}
Let $X \in \mathcal F$ be a vector field whose time $t$-flow $\varphi^X_t \colon M \to M $ exists. Then~$\varphi^X_t$ is a symmetry of $\mathcal F$. 
\end{corollary}

\begin{remark}
It deserves to be noticed that the conclusion of the corollary is \emph{not} true for the infinite comb (See Exercises \ref{exo:comb}-\ref{exo:infinitecomb}-\ref{exo:comblast}). Its proof indeed made an intense use of the assumptions ``locally finitely generated''. 
\end{remark}

Here is a second corollary, which is not totally trivial. Inner symmetries were defined in Section \ref{sec:symmetry}.

\begin{corollary}
\label{coro:Inner_is_normal}
Let $(M,\mathcal F) $ be a foliated manifold. 
    Any inner symmetry\footnote{See definition \ref{def:InnerSym}.} of $\mathcal F $  is a symmetry of $\mathcal F $.

    Moreover, the group $\mathrm{Inner}(\mathcal F) $
of inner symmetries is a normal subgroup in the group  $\mathrm{Sym}(\mathcal F) $ of symmetries of $ \mathcal F$.\end{corollary}
\begin{proof}
This is not totally obvious, since Corollary \ref{coro:flowissymmetry} only deals with vector fields in $ \mathcal F$ that do not depend on the time. However, it can be deduced from \ref{coro:flowissymmetry} as follows.

Let $I$ be an open interval of $ \mathbb R$ containing $ [0,1]$. We denote the real parameter by $t$. Let $ (X_t)_{t \in I}$ be a smooth time-dependent vector field\footnote{See Definition \ref{innersymmetries}} in $\mathcal F $, whose time-$1$ flow is an inner symmetry that we denote by $ \phi$ and let $\tilde{\mathcal F} $ be the singular foliation on $M \times I $ which is the direct product of $(M,\mathcal F) $ with $(I,\mathfrak X(I)) $.

By the equivalence between item (i) and (iii) in the first question of  Exercise \ref{exo:smoothtimedepDirectProduct}, the vector field 
$$ Y := X_t + \frac{\partial}{\partial t} $$
belongs to $\tilde{\mathcal F} $. By Proposition \ref{thm:flowofinner}, its time $1$ flow $ \phi_1^Y$ is a symmetry of $\tilde{\mathcal F} $, at least in the open subset where it is well-defined.
Since $ \phi_1^Y$ maps $ M \times \{0\}$ to $M \times \{1\} $, the symmetry $ \phi_1^Y$ induces a singular foliation isomorphism from the restriction of $\tilde{\mathcal F} $ to $ M \times \{0\}$ to the restriction of $\tilde{\mathcal F} $ to $ M \times \{0\}$. Since these restrictions coincide with $\mathcal F $ (upon identifying $M \times \{0\}$ and $M \times \{1\} $ with $M$), $ \left.\phi_1^Y\right|_{M \times \{0\}}$ induces a symmetry of $ (M,\mathcal F)$.
Now, the time-$1$ flow $\phi^Y_1 $ of $Y$ and the time $1$-flow $ \phi$ of  $ (X_t)_{t \in I}$ are related by
 $$  \phi^Y_1 (m,0) = \left(\phi(m),1\right) \hbox{ for all $ m \in M$.} $$
 This completes the proof of the claim.
\end{proof}

 \begin{remark}
 The results of this section can be extended to the real analytic or complex settings.
 \end{remark}

\subsection{The local splitting theorem}
\label{sec:splitting}

The second step in the proof of Theorems \ref{thm:Hermann} and  
\ref{thm:Hermann2} is an equivalent of Weinstein's splitting theorem in Poisson geometry. 
It was apparently rediscovered many times, see e.g., Paul Baum and Raoul Bott's Theorem 0.30 in \cite{zbMATH03423310}, or Proposition 1.12 in \cite{AS}.

Let us state this theorem first.

\subsubsection{The statements}

The results of this section are valid in the smooth, complex (upon replacing $\mathbb R$ by $ \mathbb C$ in the statements below), or real analytic cases. They are \emph{not} true in algebraic geometry.

\vspace{.5cm}

\begin{theorems}{Local splitting, version 1}{split}
Consider $\mathcal F $ a singular foliation on a manifold $M$ of dimension $d$. 

Any $m \in M$ a point admits a neighborhood on which $\mathcal F $ is isomorphic to the direct product of
\begin{enumerate}
    \item the singular foliation of all vector fields on an open ball in $\mathbb R^l $, with
$l= {\mathrm{dim}} (T_m \mathcal F) $,
    \item with a singular foliation $\mathcal T $ on an open ball in $\mathbb R^{d-l} $,
    contained in the space of vector fields vanishing at the center of the open ball. 
\end{enumerate}
Moreover, the rank of $\mathcal T$ at the center of the open ball is  $ r= {\mathrm{rk}}_m(\mathcal F)-l$.
\end{theorems}

\vspace{.5cm}

Alternatively, it can be practical to state this result in local coordinates.

\vspace{.5cm}

\begin{theorems}{Local splitting, version 2}{localsplitting3}
Let $\mathcal F $ be a singular foliation on a smooth, complex or real analytic singular foliation $M$.
Let  $ m \in M$ be a point and let $l= {\mathrm{dim}} (T_m \mathcal F)$. Every point $m$ admits a chart $\mathcal U_m$ with local coordinates $ (x_1, \dots, x_{l}, y_1, \dots, y_{d - l}) $, centered at $m$, on which the restriction of  $ \mathcal F$  admits the following generators:
\begin{enumerate}
    \item[a)] the $l$ vector fields $\frac{\partial}{ \partial x_1} , \dots , \frac{\partial}{ \partial x_l}  $,
    \item[b)] and $k$ vector fields of the form $$ f_1(y_1, \dots, y_{d-l}) \frac{\partial}{ \partial y_1} + \dots +   f_{n-l}(y_1, \dots, y_{d-l}) \frac{\partial}{ \partial y_{d-l}}  $$
    with $ f_1(0, \dots, 0)= \dots =  f_{n-l}(0, \dots, 0)=0$.
\end{enumerate}

Moreover, one can assume that $ k+l$ is equal to ${\mathrm{rk}}_m(\mathcal F)$. 
\end{theorems}

\vspace{.5cm}

Here is a third version of the local splitting theorem (inspired by \cite{AS}). Notice that there is no equivalent statement for the Weinstein splitting theorem in Poisson geometry.

\vspace{.5cm}

\begin{theorems}{Local splitting, version 3}{thm:localsplitting3}
Let $\mathcal F $ be a singular foliation on a smooth, complex or real analytic singular foliation $M$ of dimension $d$.
For every $m \in M$, 
there exists 
\begin{enumerate}
    \item an open neighborhood $ \mathcal U$ of $m$ in $M$
    \item a singular foliation $\mathcal T $ of rank ${\mathrm{rk}}_{m}(\mathcal F) - l$ on an open neighborhood $ \mathcal V$ of $0$ in $ \mathbb K^{d-l}$, with $ l = {\mathrm{dim}}(T_m \mathcal F)$, admitting $\{0\}$ as a leaf, 
    \item a surjective submersion $ \phi \colon \mathcal U \to \mathcal V $, 
\end{enumerate}
such that the restriction of $\mathcal F$ to $\mathcal U $ coincides with the pull-back singular foliation $\phi^{-1}( \mathcal T) $.
\end{theorems}

\vspace{0.5cm}
Before proving these theorems, let us recall the following lemma:

\begin{lemma} \label{lem:hadamard}
If a vector field $X$ is not zero at some point $m \in M$, then there exists a local chart $\mathcal U $ with coordinates $ (x, y_1, \dots, y_{d-1})$, centered at $m$,
such that, on $\mathcal U $, we have $ X = \frac{\partial}{\partial x} $.
\end{lemma}

\begin{proof}[Proof of Theorems \ref{thm:split}, \ref{thm:localsplitting3}, \ref{thm:thm:localsplitting3}]
We leave it to the reader to verify that all three versions of the local splitting theorem are equivalent. 
We will prove Theorem \ref{thm:localsplitting3}.
Our proof is by recursion: since the statement is local by nature, it suffices to consider the following recursion assumption 

\begin{center}
$ \mathcal H_l$ = ``The statement is proved at $m=0 $ for any singular foliation $\mathcal F$ on an open ball in a finite dimensional vector space such that $ {\mathrm{dim}}(T_0 \mathcal F)  \leq l$''.
\end{center}

For $ l=0$, $ \mathcal H_0$ is automatically true and there is nothing to prove. Assume now  $\mathcal H_l$ is valid, and let us prove  $\mathcal H_{l+1}$.

Let $X^1, \dots, X^r$ be generators of a singular foliation $\mathcal F $ defined in an open neighborhood of $0 \in \mathbb K^d$. Without any loss of generality, one can assume $ X^r|_{m} \neq 0 $. By the Hadamard lemma \ref{lem:hadamard}, there exists local coordinates  $ (x, y_1, \dots, y_{d-1})$ centered at $0$, such that, on $\mathcal U $, in which $ X^r = \frac{\partial}{\partial x} $. In these coordinates, the remaining generators read as:
$$ X^j = \sum_{i=1}^{r-1} F^j_{i}(x, y_1, \dots, y_{d-1}) \frac{\partial}{\partial y_i} + g^j (x, y_1, \dots, y_{d-1}) \frac{\partial}{\partial x}. $$
Since $ X^r=\frac{\partial}{\partial x} $ belongs to $\mathcal F $, there is a second family of generators of  $\mathcal F $ giver by $X^r$ together with the $r-1$ vector fields: 
$$ \widehat{X}^{j}:= X^j -  g^j (x, y_1, \dots, y_{d-1}) \frac{\partial}{\partial x} = \sum_{i=1}^{r-1} F^j_{i}(x, y_1, \dots, y_{d-1}) \frac{\partial}{\partial y_i} . $$ 
Let $\mathcal G $ be the module generated by $\widehat{X}^{1}, \dots, \widehat{X}^{r-1}$. 
This module has the following description:
$\mathcal G $ is the intersection of $\mathcal F $ with vector fields on the fiber of the map $$\Pi \colon (x, y_1, \dots, y_{d-1}) \mapsto x .$$
In equation:
$$ \mathcal G = \mathcal F \cap \{ \Pi-vertical\} $$
Since both $\mathcal F $ and $ \Pi$-vertical vector fields are closed under Lie bracket, it defines, in particular, a singular foliation of rank $ d-1$ on some neighborhood of $0 $.

Now, $ [X^r , \mathcal G] \subset \mathcal F$, since $X^r \in \mathcal F $ and $ \mathcal G \subset \mathcal F$. Also, $\mathcal G $ being vertical with respect to  $\Pi $ while $X^r $ is $\Pi $-related to the vector field
$\frac{\partial}{\partial x} $ on $\mathbb R $, the Lie bracket $ [X^r , \mathcal G]  $ is valued in $\Pi$-vertical vector fields, so that  
 $$ [X^r , \mathcal G] \subset \mathcal F \cap \{ \Pi-vertical\} = \mathcal G.   $$
Said differently, $X^r$ is an infinitesimal symmetry of $\mathcal G $. By Theorem \ref{prop:symInt}, its flow is a symmetry of  $\mathcal G $.
Concretely, it means that for all $ (x,y_1, \dots,y_{d-1})$ and all $t \in \mathbb R$ such that  $ (x+t, y_1, \dots,y_{d-1}) $ is still within the considered open subset,
$$ \left(\phi_{t}^{X^r}\right)_*  \left( \begin{array}{c} \widehat{X}_1(x,y) \\ \vdots \\ \widehat{X}_{r-1}(x,y) \end{array} \right) = \left( \begin{array}{ccc} && \\ & \mathbf{A}(t,x,y) & \\ && \end{array} \right) \left( \begin{array}{c} \widehat{X}_1(x,y) \\ \vdots \\ \widehat{X}_{r-1}(x,y) \end{array} \right)  , $$
where $ \mathbf A(t,x,y) $ is an invertible matrix that satisfies:
 $$  \mathbf A \left(s,\phi_t^{X^r} (x,y)\right) \circ \mathbf A (t,x,y)  = \mathbf A (t+s,x,y) $$
Since the flow at time $t$ of $ X^r$ reads 
 $$   \phi_t^{X^r} \colon (x, y_1, \dots, y_{d-1}) \longrightarrow  
  (x+t, y_1, \dots, y_{d-1}) $$
it means that there exists an invertible matrix $\mathbf A(t,x,y) $ such that:
 $$ \left( \begin{array}{c} \widehat{X}_1(x+t,y) \\ \vdots \\ \widehat{X}_{r-1}(x+t,y) \end{array} \right) = \left( \begin{array}{ccc} && \\ & \mathbf{A}(t,x,y) & \\ && \end{array} \right) \left( \begin{array}{c} \widehat{X}_1(x,y) \\ \vdots \\ \widehat{X}_{r-1}(x,y) \end{array} \right)  ,$$
 where the invertible matrix $\mathbf A (t,x,y) $ satisfies:
 $$  \mathbf A \left(s,\phi_t^{X^r} (x,y)\right) \circ \mathbf A (t,x,y)  = \mathbf A (t+s,x,y) $$
 In particular, the vector fields
  $$     \left( \begin{array}{c} Z_1 (x,y)\\ \vdots \\ Z_{r-1}(x,y) \end{array} \right) = \left( \begin{array}{ccc} && \\ & \mathbf{A}(x,0,y) & \\ && \end{array} \right)^{-1} \hspace{.3cm} \left( \begin{array}{c} \widehat{X}_1(0,y) \\ \vdots \\ \widehat{X}_{r-1}(0,y) \end{array} \right) . $$
  are well-defined in a neighborhood of $0$ and satisfy the following two properties:
  \begin{enumerate}
      \item  they are local generators of $\mathcal G $ (since the matrix $\mathbf{A}(t,x,y)$ is invertible for $t,x,y$ small enough),
      \item they are invariant under the flow of $ X^r$, since
      
        $$    (\Phi^{X^r}_t)_* \left( \begin{array}{c} Z_1 (x,y)\\ \vdots \\ Z_{r-1}(x,y) \end{array} \right) = (\Phi^{X^r}_t)_*\left(\left( \begin{array}{ccc} && \\ & \mathbf{A}(x,0,y) & \\ && \end{array} \right)^{-1} \hspace{.3cm} \left( \begin{array}{c} \widehat{X}_1(0,y) \\ \vdots \\ \widehat{X}_{r-1}(0,y) \end{array} \right)\right) . $$
$$
= \left( \begin{array}{ccc} && \\ & \mathbf{A}(x+t,0,y) & \\ && \end{array} \right)^{-1}
\left( \begin{array}{ccc} && \\ & \mathbf{A}(t,0,y) & \\ && \end{array} \right)
\hspace{.3cm} \left( \begin{array}{c} \widehat{X}_1(0,y) \\ \vdots \\ \widehat{X}_{r-1}(0,y) \end{array} \right)
$$

$$
= \left( \begin{array}{ccc} && \\ & \mathbf{A}(x,0,y) & \\ && \end{array} \right)^{-1}
\hspace{.3cm} \left( \begin{array}{c} \widehat{X}_1(0,y) \\ \vdots \\ \widehat{X}_{r-1}(0,y) \end{array} \right)
=\left( \begin{array}{c} Z_1 (x,y)\\ \vdots \\ Z_{r-1}(x,y) \end{array} \right)
$$

\end{enumerate}
In coordinates, it means that they are of the form:
 $$ Z^i = \sum_{i=1}^{d-1} f^i_j (y_1, \dots,  y_{d-1}) \frac{\partial}{\partial y_j}  .$$
They therefore define a singular foliation $\mathcal G $ of rank $ r-1$ on $\mathbb K^{d-1} $. 
By construction, the dimension of $ T_0\mathcal G$ is $ l-1$. We can then apply the recursion hypothesis, and we obtain the existence of coordinates $(x_1, \dots, x_{l}, y_1', \dots, y_{d-l-1}') $ on which $\mathcal G $ is of the form described in Theorem \ref{thm:localsplitting3}. These variables, together with $x_{l+1} := x $ form a system of coordinates on which $\mathcal F $ is also of the form decreed by Theorem \ref{thm:localsplitting3}. 
\end{proof}

\subsection{Proof of Theorems \ref{thm:Hermann} and  
\ref{thm:Hermann2}}
\label{proofHermann}

We will use the following property of immersed submanifolds:

\begin{proposition}
\label{prop:immersed}
If a connected subset $ L \subset M$ satisfies that every $ m \in L$ has a neighborhood $ \mathcal U$ such that the connected component of $ m $ in $ L \cap \mathcal U $ is a submanifold of dimension $k$, then it is a (maybe immersed) submanifold of dimension $k$.
\end{proposition}

Now we can prove Theorem \ref{thm:Hermann2} as follows.  Choose a $R$-leaf $L$. 
An immediate consequence of the local splitting theorem is that every point $ m\in L$ admits a neighborhood $\mathcal U \subset M $ admitting the following property: For the restriction $\mathfrak i_{\mathcal U}^* \mathcal F$ the set of reachable points $L^{\mathcal U}_m $ is the submanifold $y_1 = \dots = y_{d-\ell}=0 $ in some local coordinates $(x_1, \dots, x_d, y_1, \dots = y_{d-\ell}) $ on which $ m=(0, \dots, 0)$. Said otherwise, the connected component of $m$ in $L \cap \mathcal U $ is a submanifold. It therefore satisfies the assumptions of Proposition \ref{prop:immersed} and is an immersed submanifold. It is therefore also a $T$-leaf. This concludes the proof of Theorem 
\ref{thm:Hermann2}. Theorem  \ref{thm:Hermann} is an immediate consequence.
 
 \begin{remark}
 Notice that the functions $ x_1, \dots, x_{d-l}$ that appear in the local splitting theorem define a diffeomorphism   $\Phi^{\mathcal U}_m $ from the submanifold $ L^{\mathcal U}_m  $ to an open neighborhood of $\mathbb K^{d-l}$. The families $(L^{\mathcal U}_m ,\Phi^{\mathcal U}_m)_{m \in M} $ form an atlas for the leaf $L_m$.
 \end{remark}

\section{Near a leaf: the transverse singular foliation and $\mathcal F$-connections}

\vspace{1mm}
\noindent
In Section \ref{sec:leavesExist}, we saw that a singular foliation indeed decomposes the underlying manifold into submanifolds called leaves. In this subsection, we will discuss the ``shape'' of singular foliations near a fixed leaf. 
We refer to \cite{fischer2024classification}
to a general classification, which is way beyond the purpose of the present section, and uses the notion of Yang-Mills bundle as in \cite{PhDSimon,fischer2023integratingcurvedyangmillsgauge,kotov2015curvingyangmillshiggsgaugetheories}.

\subsection{Traveling along a leaf} 
\label{sec:traveling}

\vspace{1mm}
\noindent
Our first result means that ``if you travel along a leaf, the landscape you will see is always the same'', i.e., it means that any two points in the same leaf of a singular foliation $\mathcal F $ have arbitrarily small open neighborhoods on which the restrictions of $\mathcal F $ are isomorphic.
\vspace{.5cm}

\begin{propositions}{Along a leaf, the landscape is always the same}{leavesAreBoring} 
Let $L$ be a leaf of a singular foliation.
\begin{enumerate}
\item In the smooth case, for any two points $\ell_0,\ell_1 $ of the leaf $L $ of a singular foliation $\mathcal F $, there exists an inner symmetry of $\mathcal F $ mapping $ \ell_0$ to $ \ell_1$. 
\item In the holomorphic or real analytic cases, the same results hold, but the inner symmetry is in general only defined in a neighborhood of $\ell_0 $.
\end{enumerate}
\end{propositions}
\vspace{.5cm}
\noindent
We start with a lemma. Recall that we say that two points  $\ell_0,\ell_1 \in L $ are \emph{$ \mathcal F$-reachable one from the other} if there exists vector fields $ X_1, \dots, X_s \in \mathcal F$ 
and $ t_1, \dots, t_s \in \mathbb R$ such that 
 \begin{equation}\label{eq:beingreachable}  \Phi_{t_1}^{X_1} \dots \Phi_{t_s}^{X_s} (\ell_0) = \ell_1  .\end{equation}
 In the smooth setting, we assume $ X_1, \dots, X_s$ to be complete. In the complex or real-analytic settings, we can only impose that the composition of flows $\Phi_{t_1}^{X_1} \dots \Phi_{t_s}^{X_s}$ that appears in Equation \eqref{eq:beingreachable} is well-defined in a neighborhood of $ \ell_0$. The previous lemma was in fact already proven in the course of Section \ref{sec:leavesExist}, but we reprove it so that the chapter can be read independently.

\begin{lemma}
Any two points in the same leaf of $ \mathcal F$ are $ \mathcal F$-reachable one from the other.
\end{lemma}
\begin{proof}

The notion of ``being $ \mathcal F$-reachable one from the other'' defines an equivalence relation on $L$ that we call the \emph{$ \mathcal F$-reachability relation}.
Since the tangent space of the leaf $T_\ell L $ at a point $ \ell \in \mathcal F$ coincides with $ T_\ell \mathcal F$, every $u \in T_\ell L$ is of the form $ X_{|_\ell}$ for some $X \in \mathcal F$ that one can without any loss of generality assume to be of compact support, hence complete. 
Applying this reasoning to a basis $ e_1, \dots, e_k$ of $ T_\ell L$, we find complete vector fields $ Y_1, \dots, Y_k \in \mathcal F$ through $ e_1, \dots, e_k$. Now the differential of the map 
$$ \begin{array}{rcl} \mathbb R^k& \to &   L  \\ (t_1, \dots, t_k )& \mapsto&   \Phi_{t_1}^{Y_1} \dots \Phi_{t_k}^{Y_k} (\ell_0) \end{array}$$
is invertible at $ t_1= \dots = t_k =0$, since $ Y_1|_{\ell}=e_1, \dots,  Y_k|_{\ell}=e_k$ are in the image. This implies that the image of this map contains a neighborhood of $ \ell$, so that the set of points of $L$ which are $ \mathcal F$-reachable from any point  $\ell \in L$ contains a neighborhood of $ \ell$. In particular, each equivalence class of the $\mathcal F $-reachability relation is an open subset of $ L$. Since the leaf $L$ is a connected manifold by definition, this implies that there is only one equivalence class, namely $L$ itself.
\end{proof}

\begin{proof}[Proof of Proposition \ref{thm:leavesAreBoring}]
The result is now in an immediate consequence of Proposition \ref{thm:flowofinner}, which states that each flow in Equation \eqref{eq:beingreachable} is a symmetry of $\mathcal F $ (local in the complex or real analytic settings and global in the smooth setting).
\end{proof}

Proposition \ref{thm:leavesAreBoring} has several natural consequences. 
\vspace{0.5cm}

\begin{definitions}{}{}
    Let $L$ be a leaf. A pointed submanifold $(S,\ell)$  of $M$ with $ \ell \in L \cap S$ is said to be a \emph{$\mathcal F$-cut} of the leaf $L$ if
\begin{enumerate}
\item[(i)] $ S$ is transverse to $L$ at $ \ell$, i.e., $T_\ell L \oplus T_\ell S = T_\ell M $, and  
\item[(ii)] $S$ cuts $ \mathcal F$ transversally, i.e., $ T_s \mathcal F + T_s S = T_s M$ for every $s \in S$.
\end{enumerate}
\end{definitions}

\begin{remark}
Notice that any pointed submanifold $(S,\ell)$ that satisfies condition (i) admits a neighborhood $ S' \subset S$ of  $ \ell$ such that condition (ii) is also satisfied.
Notice also that any neighborhood $S'$ of a $ \mathcal F$-cut $S$ is a $\mathcal F$-cut again.
\end{remark}

\noindent
Consider a $ \mathcal F$-cut $(S,\ell) $. Then consider the restriction $\mathcal T_S :=\mathfrak i_S^* \mathcal F $ as in Section \ref{sec:transverse}, i.e., the restriction to $ S$ of vector fields in $\mathcal F$ that are tangent to $S$. The following lemma holds true.

\begin{lemma}
   For every $ \mathcal F$-cut $ (S,\ell)$ of a leaf $L$ of a singular foliation $ \mathcal F$, $\mathcal T_S :=\mathfrak i_S^* \mathcal F $ is a singular foliation on $S$. Also, the point $ \{\ell\}$ is a leaf of $ \mathcal T_S$.
\end{lemma}
 \begin{proof}
The first statement is a consequence of Proposition \ref{thm:prop:transverse}, since item (ii) in the definition of a $ \mathcal F$-cut means that $S$ intersects $ \mathcal F$ cleanly (see Definition \ref{def:beingtransverse}).
The second statement follows from the fact that $ T_\ell \mathcal T_S = T_\ell S \cap T_\ell \mathcal F$ (see Question 2 in \ref{exo:rankof restriction}). Since $L$ in the leaf, $T_\ell \mathcal F=T_\ell L$. Item (i) in the definition of a $ S$-cut then implies that $ T_\ell \mathcal T_S = T_\ell S \cap T_\ell L=\{0\}$.  Hence, the leaf through $L$ of $ \mathcal T_S$ reduces to $ \{\ell\}$.
 \end{proof}

\vspace{.5cm}

\begin{theorems}{Any two $ \mathcal F$-cuts to $L$ have isomorphic germs}{transverseFol}
Let $L$ be a leaf of a smooth, real analytic or complex singular foliation.
For any two $ \mathcal F$-cuts $ (S_0,\ell_0)$ and $ (S_1, \ell_1)$,  there exists neighborhoods $ S_0'$ of $ \ell_0 \in S_0$  and $ S_1'$ of $ \ell_1 \in S_1$ and an isomorphism of singular foliations: 
$$ \xymatrix{  (S_0',\ell_0, \mathcal T_{S_0'})\ar[rr]^{\simeq }&&
 (S_1',\ell_1,\mathcal T_{S_1'}) }$$
where $\mathcal T_{S_0'},\mathcal T_{S_1'}$ are the induced singular foliations on $S_0' $ and $ S_1'$.
\end{theorems}
\begin{proof}
We prove it in two steps:
\begin{enumerate}
    \item[\emph{Step 1}] We prove that there exists a $\mathcal F$-cut $ (\tilde{S}_0,\ell_1)$ through the point $ \ell_1$ whose restricted singular foliation $ \mathcal T_{\tilde{S_0}}$ is isomorphic\footnote{In the holomorphic or real analytic settings, one has to replace $S_0$ by a neighborhood of $ \ell_0$ in $ S_0$. This does not affect the argument} to $ (S_0,\ell_0, \mathcal T_{S_0})$. 
    \item[\emph{Step 2}] We then prove Theorem \ref{thm:transverseFol} for the case $ \ell_0=\ell_1$.
\end{enumerate}
By the first step, one can assume without any loss of generality that $ \ell_0=\ell_1$ in  Theorem \ref{thm:transverseFol}. The second step then provides a proof of the result.

\vspace{0.3cm}
\noindent
{\emph{Step 1}}.
Let $ \Phi$ be an inner symmetry as in Proposition \ref{thm:leavesAreBoring} such that $ \Phi(\ell_0)=\ell_1$. Since $ \Phi$ is a symmetry of $\mathcal F $, and since it restricts to a diffeomorphism of $L$ such that $ \Phi(\ell_0)= \ell_1$, $ (\Phi(S_0), \ell_1) $ is a $\mathcal F $-cut. Also, the restriction of $ \Phi$ to a diffeomorphism $ S_0 \to \Phi(S_0)$ is an isomorphism between their respective restricted singular foliations. This proves the first point.

\vspace{0.3cm}
\noindent

{\emph{Step 2}}.
Let $ \ell \in L$ be a point. After \ref{thm:localsplitting3}, we can without loss of generality assume $\ell$ to be the origin of $U\subset \mathbb R^{n}$ with coordinates $(x_1,...,x_l,y_1,...,y_{n-l})$ such that the foliation is generated by $\frac{\partial}{\partial x_1},...,\frac{\partial}{\partial x_l}$ and $Y_1,...,Y_{r-l}$ where $Y_i$ only depend on the $y$ variables and only contain $\frac{\partial}{\partial y_j}$ components. In particular $L=\{y=0\}$. By construction, $S_0=\{x=0\}$ is a $\mathcal F$-cut, so it suffices to show that any other $\mathcal F$-cut $S_1$ through the origin is (locally) equivalent to it. The space $T_\ell S_1$ being transverse to $T_\ell L$ translates to $T_\ell S_1\to T_\ell S_0$ being surjective, i.e., shrinking $S_0$ to $S_0'\ni \ell$, there is a local section $\sigma:S_0'\to S_1$. i.e., near $\ell$ we have $S_1=\{(f(y),y)\}$ for some function $f:\mathbb R^{n-l}\to \mathbb R^{l}$ (with $f(0)=0$). The family of diffeomorphisms $\Phi_t(x,y)=(x+t\cdot f(y),y)$ map $S_0$ to $S_1$ and are the flow of a vector field in $\mathcal F$ (fixing $\ell$), hence induces an isomorphism of transverse foliations.   This completes the proof of the Theorem.

\end{proof}

\vspace{.5cm}
\noindent
We call \emph{germ at $0 \in \mathbb K^s$} of a singular foliation on $ \mathbb K^s$ an equivalence class of pairs $ (\mathcal U, \mathcal T_\mathcal U)$ where 
\begin{enumerate}
\item $ \mathcal U$ is a neighborhood of $ 0$ in $ \mathbb K^s$, and
\item $\mathcal T_\mathcal U$ is a singular foliation on $\mathcal U $,
\end{enumerate}
under the equivalence relation that consists in identifying  $ (\mathcal U, \mathcal T_\mathcal U)$ 
and $   (\mathcal V, \mathcal T_\mathcal V)$ if there exists a neighborhood  $\mathcal W $ of $ 0$ contained in $\mathcal U \cap \mathcal V $
on which the restrictions of $ \mathcal T_\mathcal U$ and $ \mathcal T_\mathcal V$ coincide.
Lastly, we call \emph{isomorphism classes of germs at $0 \in \mathbb K^s$} of singular foliations the equivalence  classes for the equivalence relation on germs at $0$ of singular foliations on $\mathbb K^s $ that identifies two germs if they have representatives $ (\mathcal U,\mathcal F_\mathcal U)$ and $ (\mathcal V,\mathcal F_\mathcal V)$ which are isomorphic as foliated manifolds, through an isomorphism that maps $0$ to $0$.

\begin{remark}
\label{rmk:meaningOfEquivalnceForIsoOfGerms} 
Alternatively, an isomorphism class of germs of singular foliations at $0\in \mathbb K^s$ is an equivalence class for the equivalence relation on pairs as in items 1. and 2. above under the equivalence relation that identifies  $ (\mathcal U, \mathcal T_\mathcal U)$   and  $ (\mathcal V, \mathcal T_\mathcal V)$ if and only if there exists open neighborhoods $\mathcal U' \subset \mathcal U$ and $ \mathcal V' \subset \mathcal V$ of $ 0$ and an isomorphism of singular foliation $$ \Phi \colon  (\mathcal U', \mathfrak i^*_{\mathcal U'} \mathcal T_{\mathcal U}) \longrightarrow (\mathcal V', \mathfrak i^*_{\mathcal V'} \mathcal T_{\mathcal V})  $$
that maps $0 $ to $0$.
\end{remark}

Let $L$ be a leaf of dimension $k$ of a singular foliation $\mathcal F $ on a manifold $M$ of dimension $d$. A pair  $ (\mathcal U,\mathcal T_\mathcal U)$, with $ \mathcal U$ a neighborhood of $0$ in $\mathbb K^{d-k} $, is called a \emph{representative of the transverse singular foliation of $L$} if there exists a $ \mathcal F$-cut $(S,\ell)$ whose restricted singular foliation is isomorphic to  $ (\mathcal U,\mathcal T_\mathcal U)$ through an isomorphism that maps $\ell \in S$ to $0 \in \mathbb K^{d-k}$. Theorem \ref{thm:transverseFol} implies that any two representatives of the transverse singular foliation of $L$ are in the same class for the equivalence relation defining isomorphism classes of germs at $0$ of singular foliations.
The next definition therefore makes sense.

\vspace{0.5cm}

\begin{definitions}{``The'' transverse singular foliation of a leaf}{transverse}
Let $L$ be a dimension $k$ leaf of a singular foliation $ \mathcal F$ on a manifold of dimension $d$.    We call \emph{transverse singular foliation of $L$} \underline{the} class in isomorphism classes of germs at $0\in \mathbb K^{d-k}$ of representatives of the transverse singular foliation of $L$.
\end{definitions}
 
 \subsection{Tubular neighborhoods and $ \mathcal F$-connections}

We now introduce a type of Ehresmann connection in a neighborhood of a leaf that appeared in \cite{LGR}, Section 2.2., and is generalized and used in \cite{Ryvkin2} and \cite{fischer2024classification} to classify neighborhoods of leaves.
In this section, we work in the smooth setting. 
In the real analytic or holomorphic setting, the objects introduced here (tubular neighborhood, $ \mathcal F$-connections) make sense and can be defined mutatis mutandis, but the issue is that they may not exist, because their construction goes through partitions of unity.
\vspace{.5cm}

Let $L$ be an embedded submanifold of a manifold $M$. There always exist a pair $(\mathcal U_L,p)$ where:
\begin{enumerate}
    \item $\mathcal U_L $ is an open neighborhood of $L $ in $M $,
    \item $ p \colon \mathcal U_L \longrightarrow L$ is a surjective submersion (whose restriction to $L$ is the identity)
\end{enumerate}
Moreover, one can assume $ \mathcal U_L$ is diffeomorphic to a neighborhood of the zero section in the normal bundle $ N_L:=TM_{|_L}/TL$ through a diffeomorphism that intertwines $p$ and the natural projection of $ N_L$ on its base. Here, we will never need such an isomorphism, and we will simply call \emph{tubular neighborhoods} pairs $(\mathcal U_L,p)$ as above.

In the smooth context, tubular neighborhood always exist, and any two tubular neighborhoods $(\mathcal U_L,p)$ and $(\mathcal U_L',p')$ have restrictions near $L$ which are isomorphic through an isomorphism that intertwines $p$ and $p'$. In the real analytic or complex setting, tubular neighborhoods do not exist in general.

In the smooth setting moreover, any tubular neighborhood admits an \emph{Ehresmann distribution}, i.e., there exists a smooth distribution $ H \subset T \mathcal U_L$ of constant rank $k= {\mathrm{dim}}(L)$ such that for every point $m \in\mathcal U_L $, we have 
  $$  {\mathrm{Ker}}(T_m p) \oplus H_m = T_m \mathcal U_L. $$
Moreover, there exists an Ehresmann distribution $H$ such that $ H_\ell = T_\ell L$ for every $\ell \in L$. 
Given a vector field $X \in \mathfrak X(L)$, one defines a smooth vector field  $H(X)$ on $ \mathcal U_L$ by imposing that $ H(X)$  be  section of $ H$ such that $ T_m p (H_m(X)) =X_{|_m} $ for all $ m \in \mathcal U_L$. Equivalently, $ H(X)$  is the unique vector field on $\mathcal U_L $ which is valued in the distribution $H$ at every point and is $p$-related to $X$. One call \emph{horizontal lift} of an Ehresmann distribution $H$ the henceforth induced map:
   \begin{equation}\label{eq:defhorizintallift}  \begin{array}{rrcl} H \colon &\mathfrak X(L) &\longrightarrow &\mathfrak X (\mathcal U_L) \\ & X &\mapsto & H(X) \end{array}    \end{equation}
   The horizontal lift satisfies several properties that we list below:
\begin{enumerate}
\item for every $X \in \mathfrak X(L)$ and $f \in \mathcal C^\infty(L)$, we have $ H(f X) =p^* f \, H(X)$, 
\item for every $X \in \mathfrak X(L)$, we have $ p_* (H(X))=X$,
\item for every $X,Y \in \mathfrak X(L)$,  $ \kappa_H(X,Y) := H([X,Y]) - [H(X), H(Y)]$ is a vector field tangent to the fiber of $ p$.
\item $ \kappa_H(X,Y)$ is a skew-symmetric and $ \mathcal C^\infty(L)$-bilinear map. It can therefore be seen as a $2$-form on $L$ valued in $p$-vertical vector fields. It is called the \emph{curvature} of the Ehresmann distribution.
\end{enumerate}
Now, in the context of singular foliations, when $L$ is a leaf, one will impose one more constraint on the Ehresmann connection.
In Definition \ref{def:Fconnection} below, items 1. and 2. mean that $(\mathcal U_L,p)$ is a tubular neighborhood.

\vspace{.5cm}

\begin{definitions}{$\mathcal F$-connection}{Fconnection}
Let $L$ be an embedded leaf of a singular foliation $\mathcal F $. We say that a triple $(\mathcal U_L,p,H)$ where:
\begin{enumerate}
    \item $\mathcal U_L $ is an open neighborhood of $L $ in $M $,
    \item $ p \colon \mathcal U_L \longrightarrow L$ is a surjective submersion (whose restriction to $L$ is the identity), and
    \item $H$ is an Ehresmann distribution with respect to $p$,
\end{enumerate}
is a \emph{$\mathcal F$-connection} if sections of $H$ are included in $\mathcal F $.
\end{definitions}

\vspace{.5cm}
\begin{remark}
    \cite{fischer2024classification} explains why $ \mathcal F$-connections are a particular case of the so-called Yang-Mills connections \cite{PhDSimon,fischer2023integratingcurvedyangmillsgauge,kotov2015curvingyangmillshiggsgaugetheories}, at least at formal level.
\end{remark}

\begin{exo}
\label{exo:justtangentenough}
Could the phrase ``sections of $H$ are included in $\mathcal F $'' in Definition \ref{def:Fconnection} above could be equivalently replaced by 
\begin{enumerate}
\item ``$ H \subset T\mathcal F$'',
\item or by ``$H_m$ is tangent to the leaf through $m$ for every $m \in M$''.
\end{enumerate}
\emph{Hint}: Consider the singular foliation on $\mathbb R^2 $ generated by $ \frac{\partial}{\partial x} $ and $ y^2 \frac{\partial}{\partial y}$, the leaf $L=\{(x,0), x \in \mathbb R\}$, the projection $p \colon (x,y) \mapsto (x,0) $, and $ H =\left\langle \frac{\partial}{\partial x}+ y  \frac{\partial}{\partial y}   \right\rangle$. 
\end{exo}

\begin{remark}
\label{rem:fibersofFconnection}
If on a tubular neighborhood $ (\mathcal U_L,p)$ there exists an Ehresmann connection which is an $ \mathcal F$-connection, then 
for every $ \ell \in L$, the fiber $p^{-1}(\ell)$ of $p\colon \mathcal U_L \longrightarrow L$ is a submanifold that cleanly intersects $\mathcal F$. Each fiber is therefore a $ \mathcal F$-cut of the leaf $L$.
In particular, for every $\ell \in M$, $p^{-1}(\ell) $ admits a restricted singular foliation that we will denote by $\mathcal T_\ell$ when needed.
\end{remark}

\begin{remark}
For any $\mathcal F $-connection $(\mathcal U_L,p,H) $, we have $ H_m \subset T_m \mathcal F$ for all $ m \in M$. This condition is however not sufficient to guarantee that $H$ defines an $\mathcal F$-connection, see Exercise
\ref{exo:justtangentenough}.
\end{remark}

The phrase  ``sections of $H$ are included in $\mathcal F $'' in Definition \ref{def:Fconnection} above could be equivalently replaced by the following condition ``the horizontal lift $H$ (see Equation \eqref{eq:defhorizintallift}) is valued in $\mathcal F $''.

\begin{lemma} A $\mathcal F$-connection $(U,p, H)$ for an embedded leaf $L$ is equivalent to the data\footnote{In \cite{fischer2024classification}, which works in the transversally formal setting, this lemma becomes the definition of an $ \mathcal F$-connection.} given by: 
\begin{enumerate}
\item[1] a neighborhood $\mathcal U$ of $L$ in $M$ and a projection $p$ as above,
\item[2] and a $\mathcal C^\infty(L)$-linear section 
of the natural projection 
$\mathcal F_{\mathcal U}^{proj} \to \mathfrak X(L)$, where $\mathcal F_{\mathcal U}^{proj} $ stands for vector fields in $ \mathcal F_{\mathcal U}$ which are $p$-related with a vector field in $\mathcal F$.
\end{enumerate}
\end{lemma}
\begin{proof}
Given a $\mathcal F$-connection as in Definition \ref{def:Fconnection}, the morphism $s$ of item 2 is the horizontal lift defined  as in Equation \eqref{eq:defhorizintallift}. Conversely, given a section $s$ as in item 2, the horizontal distribution  $H$ is the distribution generated by the vector fields $s(\mathfrak X(L))$. 
\end{proof}

Using this characterization one can verify that, at least in the smooth setting, $\mathcal F$-connections exist:
\vspace{0.5cm}

\begin{propositions}
    {$ \mathcal F$-connections exist.}{Fconnexionsexist}
Any embedded leaf of a smooth singular foliation admits a $\mathcal F$-connection.
\end{propositions}

\begin{proof}   
Fix a tubular neighborhood $(\mathcal U,p)$ of an embedded $l$-dimensional leaf $L$ of the singular foliation $ \mathcal F$. The local splitting Theorem \ref{thm:split} obviously implies that $\mathcal F$-connections exist locally, i.e., every point $ \ell \in L$ admits a neighborhood $\mathcal V $ in $ M$ which an $ \mathcal F$-connection exists: it suffices to consider (we use notations of Theorem \ref{thm:localsplitting3}) the vector fields 
$$ H := \left\langle\frac{\partial}{ \partial x_1} , \dots , \frac{\partial}{ \partial x_l} \right\rangle, $$
which define a distribution $H$ whose sections are in $ \mathcal F$ and is in direct sum with $ {\mathrm{Ker}}(Tp)$ in a neighborhood of $ y=0$. As a consequence, upon shrinking $\mathcal U$ if necessary, we can cover a neighborhood of a leaf $L$ by small open subsets $(\mathcal U_i)_{i \in I}$ such that $p^{-1}(\mathcal U_i)\cap \mathcal U$ admits a $\mathcal F$-connections  for every index $i$.  Consider
$s_i:\mathfrak X(\mathcal U_i)\to \mathcal F_{p^{-1}(\mathcal U_i)}^{proj}$ their horizontal lifts.
One then glues these local connections by the use of a partition of unity $(\chi_i)_{i \in I}$ for the open cover $(\mathcal U_i)_{i \in I}$ of $ L$, and defines $s(X)=\sum_{i}p^*\chi_i\cdot s_i(X|_{U_i})$ for every vector field $ X \in \mathfrak X(L)$. This completes the construction.
\end{proof}

For $L$ a leaf a singular foliation, here is a natural question: does there exist a regular foliation $ \mathcal R$ contained in $\mathcal F$ and admitting $L$ as a leaf? The answer is that such a regular foliation exists near a given embedded leaf $L$ if and only there exists a flat $ \mathcal F$-connection, i.e., a $\mathcal F$-connection for which the curvature (defined above) is zero (i.e., if $H$ is an integrable distribution).

\begin{exo}
\label{exo:flatleaves}
Let $L$ be an embedded leaf of a smooth singular foliation $\mathcal F$. Show that the following items are equivalent:
\begin{itemize}
\item[(i)] There exists a neighborhood $U$ of $L$ and a regular foliation $\mathcal R\subset \mathcal F_U$ of rank $dim(L)$.
\item[(ii)] $L$ admits a flat $\mathcal F$-connection.
 \end{itemize}
\end{exo}

 We say that a leaf is \emph{flat} if one of the equivalent conditions in Exercise \ref{exo:flatleaves} is satisfied.

\subsection{Fibered products along a leaf}
\label{sec:fibered}

We now present a construction of ``pull-back'' of a singular foliation near a leaf, when an $ \mathcal F$-connection is given. It does not depend on the choice of an $ \mathcal F$-connection, but the existence of  an $ \mathcal F$-connection  is however needed to guarantee its existence of the pullback. 
Let $ L$ be an embedded leaf of a singular foliation $ \mathcal F$.
Assume that one is given:
\begin{enumerate}
\item An $\mathcal F $-connection $(\mathcal U_L,p,H) $
\item Any manifold $\Sigma$ and any map $ \varphi \colon \Sigma \longrightarrow L$. Notice that we do not require $ \varphi$ to be a submersion, nor an immersion: it just has to be in the relevant category of maps (smooth, real analytic or holomorphic).
\end{enumerate}
Since $ p$ is a surjective submersion,  the fibered product 
$$  \varphi^! \mathcal U_L := \mathcal U_L \times_{p,L,\varphi} \Sigma := \{(m,\sigma) \in M \times \Sigma  \, | \, p(m)=\varphi(\sigma) \} $$
is a submanifold of $\mathcal U_L \times \Sigma $. Now,  $\mathcal U_L \times \Sigma $ comes with the direct product singular foliation $\mathcal F \times  \mathfrak X(\Sigma)  $ , i.e., the direct product of the foliation $ \mathcal F$ (restricted to $ \mathcal U_L$) and of the foliation $ \mathfrak X(\Sigma)$ of all vector fields on $ \Sigma$.   

\begin{lemma} The submanifold $\mathcal \varphi^! \mathcal U_L := U_L \times_{p,L,\varphi} \Sigma$ intersects cleanly the singular foliation
$\mathcal F \times  \mathfrak X(\Sigma)  $.
\end{lemma}
\begin{proof}
The proof relies on the existence of an Ehresmann connection $H$. Since $H$ is in direct sum with $ {\mathrm{Ker}}(Tp)$, it is easily checked that any vector of  $T (\mathcal U_L \times \Sigma )$ reads as a sum of an element in $u+v+w $ with $u \in H$, $ v \in {\mathrm{Ker}}(Tp)$ and $w \in T\Sigma$. Now by definition, $v$ is tangent to the submanifold  $T (\mathcal U_L \times_{p,L,\varphi} \Sigma )$, $w$ belongs to $ T(\mathfrak X(\Sigma)) \subset T(\mathcal F \times \mathfrak X(\Sigma))$. Moreover, $u$ belongs to $T(\mathcal F) \subset T(\mathcal F \times \mathfrak X(\Sigma))$ if $H$ is an $ \mathcal F$-connection.
The conditions in Definition \ref{def:beingtransverse} are therefore satisfied.
\end{proof}

We denote by $ \varphi^!  \mathcal F  $
the restriction of the direct product singular foliation to $ \varphi^! \mathcal U_L := \mathcal U_L \times_{p,L,\varphi} \Sigma$.

\begin{remark}
The existence of a distribution $H$ whose sections are in $\mathcal F $ was used to claim that it is a singular foliation, $\varphi^!  \mathcal F$ was defined without any reference to $H$. It therefore does not depend on the latter.
\end{remark}

We have two lists of comments. The first one is a list of generalities about fibered products, and the second one  relates the items of the first list to the properties of the foliated manifold $ (\varphi^! \mathcal U_L,\varphi^!  \mathcal F)$. 
\begin{enumerate}
 \item The natural projection $ \varphi^!p : \varphi^! \mathcal U_L    \to \Sigma$ is a surjective submersion, whose fiber over $ \sigma \in \Sigma$ is diffeomorphic to the fiber of $p $ over $ \varphi(\sigma)$.
 \item  The submanifold $L \times_{p,L ,\varphi} \Sigma  $ is canonically diffeomorphic to $ \Sigma$, making $ \Sigma $ a submanifold of $ \mathcal U_L \times_{p,L,\varphi} \Sigma$. 
 This inclusion is a right inverse of the above projection.
  \item Any Ehresmann connection $H$ on $(\mathcal U_L, p) $ induces an Ehresmann connection $ \varphi^! H$ on $   \mathcal U_L \times_{p,L,\varphi} \Sigma$. By construction,
 $$  \varphi^! H_{(m,\sigma)} := \{ (h,v)\in H_m \times T_\sigma \Sigma \, | \, T_m p (h) = T_\sigma \varphi (u) \} .$$
\end{enumerate}
We leave it to the reader to check the following list of points: 
\begin{enumerate}
\item The fibers of $\varphi^! p  $ and of $p$ intersect cleanly  $\varphi^! \mathcal F  $ and $\mathcal F $ respectively, and the diffeomorphism described in the first item above is an isomorphism of their respective restrictions.

\item The manifold $\Sigma\in \phi^!\Sigma$ is a leaf of $\varphi^! \mathcal F$.
\item If the Ehresmann connection $H$ is an $ \mathcal F$-connection for $ \mathcal F$, then $ \varphi^! H$ is a  $ \mathcal F$-connection for $ \varphi^! \mathcal F$.
\end{enumerate}
As a conclusion, $ ( \varphi^! \mathcal U_L , \varphi^! p, \varphi^! H)  $ is an $ \mathcal F$-connection for the leaf $\Sigma $ of $\varphi^! \mathcal F $. We call it the \emph{fibered product of $(  \mathcal U_L , p,  H)$ with respect to $ \varphi\colon \Sigma \to L$}. 

\label{sec:fiberleaf}

\subsection{Parallel transportation and $ \mathcal F$-connections}

For this section, in addition to the own works of the authors \cite{LGR,Ryvkin2}, we acknowledge ideas coming from the PhD \cite{PhDSimon,fischer2023integratingcurvedyangmillsgauge} of Simon Raphael Fischer, see also \cite{fischer2024classification}.  

Let $ (\mathcal U_L, p, H)$ be an $ \mathcal F$-connection for an embedded leaf $L$ of a singular foliation $ \mathcal F$.
Let $ H: \mathfrak X(L) \to \mathcal F_{proj}$ be the horizontal lift as in Equation \eqref{eq:defhorizintallift}.
Recall that the notion of smooth time-dependent vector field valued in $ \mathcal F$ has some subtleties: they are dealt with in Definition \ref{innersymmetries}. Let $I$ be an open interval of $ \mathbb R$.

\begin{lemma}
 \label{lem:liftisisomorphism}
 For any smooth time dependent vector field $(X_t)_{t \in I} $ on the leaf $ L$, $H(X_t)$ is a smooth time dependent vector field in $\mathcal F $.
\end{lemma}
\begin{proof}
For every $ m \in \mathcal U_L$, let $X_1, \dots, X_k \in \mathfrak X(L) $ be a local trivialization of $ TL$ on a neighborhood $ \mathcal W$ of $ p(m)$. There exists functions $ f_1(\ell,t), \dots, f_k (\ell,t)$ on $ \mathcal W \times I$ such that 
 $$ \left.X_t\right|_m = \sum_{i=1}^k   f_i(m,t) \left.X_i\right|_{m}.   $$
This implies that
 $$  H(X_t)_{|_m}= \sum_{i=1}^k f_i(p(m),t) H (X_i) $$
 is a smooth time-dependent vector field on $ \mathfrak X (L)$.
\end{proof}

\noindent Here is an immediate consequence of this lemma, together with the fact that $ H(X)$ is $p$-related with~$X$.

\begin{proposition}
    \label{prop:liftissymmetry}
    For every time dependent vector field $ (X_t)_{t \in I}$ on $L$ whose time $t_0$ flow is well-defined, the flow $ \Phi_{t_0}^{H(X_t)}$ at time $t_0$ of the horizontal lift $ H(X_t)$
    \begin{enumerate}
    \item is well-defined in a neighborhood $ \mathcal U_L'$ of $L$ in $ \mathcal U_L$,
    \item is a symmetry of $ \mathcal F$,
    \item commutes with $p$, more precisely 
     $$ \xymatrix{ \ar[d]^p \mathcal U_L' \ar[r]^{\Phi_{t_0}^{H(X_t)}}& \mathcal U_L \ar[d]_p \\ L \ar[r]^{\Phi_{t_0}^{X_t}}&  L}   .$$
    \end{enumerate}
\end{proposition}

Recall from Section \ref{sec:fiberleaf} that for every $ \ell \in L$, $p^{-1}(\ell) $ is an $ \mathcal F$-cut for the leaf $L$, so that the restriction of $\mathcal F $ to $p^{-1}(\ell) $ is a singular foliation that we denote by $ \mathcal T_\ell$. Recall from Theorem \ref{thm:transverseFol} that the foliated manifolds $ (p^{-1}(\ell_0), \mathcal T_{\ell_0})$ and $ (p^{-1}(\ell_1), \mathcal T_{\ell_1})$ have restrictions to neighborhoods of $ \ell_0$ and $\ell_1 $ respectively which are isomorphic, through an isomorphism of foliated manifolds that maps $ \ell_0$ to  $ \ell_1$.  

We now explain how parallel transportation with respect to an $ \mathcal F$-connection realizes such an isomorphism of foliated manifolds.
Let us briefly recall the notion of parallel transportation.
Given an Ehresmann connection $H$ on a tubular neighborhood $(\mathcal U_L,p) $ of an embedded manifold $L $, one says that a path  $ \tilde{\gamma} \colon [0,1] \to \mathcal U_L $  in $ \mathcal U_L$ 
is \emph{parallel} if for all $ t\in [0,1]$:
$$\frac{d\gamma}{dt}(t) \in H_{\gamma(t)}  .$$
Given a path $ \gamma \colon [0,1] \to L $ on $L$ such that $ \gamma(0)=\ell_0 $, and given a point $ m \in p^{-1}(\ell_0)$, there exists at most one parallel path $ \gamma^H_m\colon [0,1] \to \mathcal U_L $ such that 
 $$\tilde{\gamma}_m (0) = m \hbox{ and }  p \left(\gamma^H_m (t) \right) = \gamma(t) \hbox{ for all $ t\in [0,1]$} .$$
 We call this path the \emph{horizontal lift} of $ \gamma(t)$ starting from $ m$.
If, moreover, $ H_\ell= T_\ell L$ for every $ \ell \in L$, then the path $ \gamma_m^H (t)$ is well defined at time $1$ for every $m $ in a neighborhood $\mathcal U_{p^{-1}(\ell_0)} $ of $ \ell_0$ in $ p^{-1}(\ell_0)$. 
We call \emph{parallel transportation} over $\gamma(t) $ the map:
$$ \begin{array}{rrcl} P_\gamma^H \colon & \mathcal U_{p^{-1}(\ell_0)} & \to & p^{-1}(\ell_1)\\ &m &\mapsto  & \gamma_m^H(1) \end{array}   $$
 The following results are classical, cf. e.g., \cite{kolarmichor}:
\begin{enumerate}
\item The parallel transportation $P_\gamma^H$ is a diffeomorphism onto its image. 
\item Its inverse is the parallel transportation over the path $ t \mapsto \gamma(1-t)$.
\item The parallel transportation does not depend on a parametrization of $ \gamma (t) $. To be more precise $ P_{\gamma}^H= P_{\gamma \circ \psi }^H$ for every $ \psi \circ [0,1] \to [0,1]$ such that $ \psi(0)=0$ and $ \psi (1)=1$.
\item The previous item allows making sense of the following statement: given  two paths $ \gamma_1,\gamma_2  \colon [0,1] \to L $ such that $ \gamma_1(1)=\gamma_2(0)$, one has:
 $$ P_{\gamma_2}^H \circ P^H_{\gamma_1} = P^H_{\gamma_1 * \gamma_2} $$ 
 where $ \gamma_1 * \gamma_2 \colon [0,1] \to L$ is any path that merges $ \gamma_1$ and $ \gamma_2$, for instance 
  $$  \gamma_1 * \gamma_2 (t) =\left\{ \begin{array}{ll} \gamma_1(2t)& \hbox{ for $t \in [0,1/2]$} \\ \gamma_2(2t-1) & \hbox{ for $t \in [1/2,1]$} \end{array} \right. . $$
\end{enumerate}
Recall that an Ehresmann connection $H$ is said to be \emph{complete}\footnote{For instance, the central circle of the restriction to a relatively compact neighborhood of it of the so-called self-eating snake does not admit a complete $ \mathcal F$-connection, see Example \ref{ex:snake}. One can avoid completeness by working at the formal level (as in \cite{fischer2024classification}), or at the level of germs. } if its associated lift $s_H$ (defined as in Equation \eqref{eq:defhorizintallift}) maps complete vector fields on $L $ to complete vector fields on $\mathcal U_L $. In this case, it can be shown that $ P_\gamma^H$ is defined on the whole fiber $ p^{-1}(\ell_0)$, and is a diffeomorphism onto the fiber $ p^{-1}(\ell_1)$.
 \emph{We will from now on assume that the $H$-connection is complete in order to simplify the proofs and the statements}.
 
 For an arbitrary $\mathcal F$-connection we would get the corresponding statement on neighborhoods of $L$, but we will leave the reader to make the generalization.

Now, assume the complete Ehresmann connection is a $\mathcal F$-connection on some tubular neighborhood of an embedded leaf. 

\begin{lemma}\label{lem:elementaryhol0}
Let $\gamma:[0,1]\to L$ be a smooth path from $\ell_0$ to $\ell_1$.
The parallel transportation $P^H_\gamma:p^{-1}(\ell_0)\to p^{-1}(\ell_1)$ is an isomorphism of singular foliations from $(p^{-1}(\ell_0),\mathcal T_{\ell_0})$ to $(p^{-1}(\ell_1),\mathcal T_{\ell_1})$.   We denote it by $ P_{\ell_0,\gamma, \ell_1}^H$
\end{lemma}
\begin{proof} 
Let $I = [0,1]$.
For any path $ \gamma \colon I \to L$, there exists a compactly supported time dependent vector field $(X_t)_{t \in [0,1]}$ on $L$ such that 
$$ X_t |_{\gamma(t)} = \frac{d \gamma(t)}{dt}$$
for all $t \in I$. The integral curve of $(X_t)_{t \in I} $ starting from $ \ell_0$ coincides with the path $ \gamma \colon I \to L$ by construction. 
Since $X_t$ and $H$ are complete, Proposition \ref{prop:liftissymmetry} applies with  $ \mathcal U_L' =\mathcal U_L $, and yields 
the time $1$-flow $ \Phi_{1}^{H(X_t)}$ in $ H(X_t) $. By construction, the restriction to the fiber $ p^{-1}(\ell_0)$ of  $\Phi_{1}^{H(X_t)}$ coincides with $ P_{\ell_0,\gamma,\ell_1}$.
Since $  \Phi_{1}^{H(X_t)}$ is a symmetry of $ \mathcal F$, Exercise \ref{exo:restrictionAndSymmetry} applies and yields the desired isomorphism of the restricted singular foliations.
\end{proof}

Let us recall some vocabulary from Section \ref{sec:symmetry}. Given  two isomorphisms of singular foliations $\Phi, \Psi : (M_1,\mathcal F_1) $ and $(M_2,\mathcal F_2) $, the following are equivalent:
\begin{enumerate}
\item[(i)]  $ \Psi^{-1} \circ \Phi$ is an inner symmetry of $(M_1,\mathcal F_1) $.
\item[(ii)] $ \Phi \circ \Psi^{-1}$ is an inner symmetry of $(M_2,\mathcal F_2) $,
\end{enumerate}
We say that $ \Phi,\Psi$ \emph{differ by an inner symmetry} of one of these equivalent conditions are satisfied.
In the next lemma, again, the assumption ``complete'' could be deleted, at the expense of restricting ourselves to a smaller neighborhood of $ L$

\begin{lemma}\label{lem:elementaryhol}
Let $H,H'$ be two complete $ \mathcal F$-connections on the same tubular neighborhood $(\mathcal U_L,p) $.
Let $\gamma:[0,1]\to L$ be a smooth path from $\ell_0$ to $\ell_1$.
The singular foliations isomorphisms  $ P_{\ell_0,\gamma, \ell_1}^H$ and $ P_{\ell_0,\gamma, \ell_1}^{H'}$ obtained by parallel transportation over $ \gamma(t)$ with respect to $H $  and $H'$ differ by an inner symmetry. 

\end{lemma}
\begin{proof}
Observe that $H(X)-H'(X)$ is for every $ X \in \mathfrak X(L)$ valued in the singular foliation of vector fields in $\mathcal F$, which are $p$-vertical, which we denote by $\mathcal T_\bullet $.
Let $(X_t)_{t \in I}$ be as in the proof of Lemma \ref{lem:elementaryhol0}.
Then $Y_t:=H(X_t)-H'(X_t)$
is a smooth time dependent vector field in $\mathcal T_\bullet $.
We will use the following classical formula (see, e.g., \cite{PosilicanoFlow}) for time-dependent vector fields $A_t,B_t$, where $\phi_{t,s}$ denotes the flow at time $t$ with starting from time $s$:
$$
\phi^{A,B}_{t,t_0}=
\phi^A_{t,0}\circ 
\phi^{(\phi_{0,t}^A)_*B_t}_{t,t_0}
\circ \phi^A_{0,t_0}
$$
Evaluating it at $t_0=0$, $A=s(X)$ and $B=Y$ yields:
$$
\phi^{s(X)+Y}_{t,0}=
\phi^{s(X)}_{t,0}\circ 
\phi^{(\phi_{0,t}^{s(X)})_*Y_t}_{t,0}
$$
Since the flow $\phi^{s(X)}$ preserves both $\mathcal F$ and verticality, $Z={(\phi_{0,1}^{s(X)})_*Y_1}|_{p^{-1}(x)}\in \mathcal T_x$, i.e., $
\phi^{s(X)+Y}_{1,0}=
\phi^{s(X)}_{1,0}\circ \psi
$
for an inner symmetry $\psi$. A similar technique has been used to prove Proposition 2.3 in \cite{AZ2}.
\end{proof}

\noindent
Now we can finally show the following result. The setting is the one of Lemmas \ref{lem:elementaryhol}-\ref{lem:elementaryhol0}.

\begin{lemma}
\label{lem:elementaryhol2}
Let $H$ be a complete $\mathcal F $-connection. The isomorphisms of singular foliations obtained by parallel transportation with respect to $H$ over homotopic paths differ by an inner symmetry\footnote{See Definition \ref{def:InnerSym}}. 
\end{lemma}
Let us restate Lemma \ref{lem:elementaryhol2} more precisely: it says that for $ \gamma_0$ and $ \gamma_1$ two homotopic paths in $L$ from $ \ell_0$ and $ \ell_1$,  the parallel transportations 
$P_{\ell_0,\gamma_0,\ell_1}^H$ and $P_{\ell_0,\gamma_1,\ell_1}^H $, which are isomorphisms of singular foliations by Lemma \ref{lem:elementaryhol0}, differ by an inner symmetry. 
\begin{proof}
In this proof we essentially follow part of the proof of \cite{kolarmichor} of the (generalized) Ambrose-Singer theorem, adapted to our setting. Given two homotopic paths, $\gamma,\tilde \gamma$ we can reparametrize them such that they are constant near the boundary, so that composing them yields a smooth null homotopic loop. Hence, it suffices to show the statement for a (smoothly) null homotopic loop $\gamma$. Let $H$ be a homotopy (with fixed endpoint $x$) such that $H_1=\gamma$ and $H_0\equiv x$. We can consider $f_t= P^s(H_t)$ and want to show that $\frac{\partial f_t}{\partial t}\circ f_t^{-1}$ is a smooth time-dependent\footnote{Defined as Definition \ref{innersymmetries}.} vector field $Z_t$ on $(p^{-1}(x),\mathcal T_x)$, then $f_t$ is its flow and the claim follows. In order to show the claim, one can pull back the whole situation along $H$ and obtain a tubular neighborhood with foliation and connection over $[0,1]^2$. In \cite[9.11, Claim 2]{kolarmichor} it is shown that for $X=\frac{\partial}{\partial x}$, $Y=\frac{\partial}{\partial y}$ on the unit square we have:

\begin{align*}
Z_t=\frac{\partial f_t}{\partial t}\circ f_t^{-1}=\int_0^1
&-(\phi^{s(X)}_\tau)^* [s(X),s(Y)] +(\phi^{s(X)}_\tau)^*[s(X),(\phi^{s(Y)}_t)^*(\phi^{s(X)}_{-\tau})^*s(Y)]\\&
-(\phi^{s(X)}_\tau)^*(\phi^{s(Y)}_t)^*(\phi^{s(X)}_{-\tau})^*[s(X),s(Y)]
d\tau
\end{align*}

Applying flows, commutators and integrals to (time-dependent) elements in the foliation yields elements in the foliation, so $Z_t$ is a smooth time-dependent vector field in the transverse foliation, i.e., its flow $f_t$ is an inner symmetry.
\end{proof}

Let $(\mathcal U_L,p,H) $ be an $ \mathcal F$-connection with the Ehresmann connection $H$ being complete.
Altogether, Lemmas \ref{lem:elementaryhol0}-\ref{lem:elementaryhol2} above imply that for every $ \ell \in L$, there exists a group morphism
from $ \pi_1(L,\ell)$ (the fundamental group of the leaf $L$) to the group of outer symmetries of the singular foliation $(\pi^{-1}(\ell), \mathcal T_\ell) $, i.e., the quotient of the group of symmetries of the singular foliation $(\pi^{-1}(\ell), \mathcal T_\ell) $ by the group of inner symmetries of $(\pi^{-1}(\ell), \mathcal T_\ell) $, see Section \ref{sec:symmetry}.
In equation:
 \begin{equation}\label{eq:holonomymorphism} \Xi \colon  \pi_1(L,\ell) \longrightarrow \mathrm{Out}(\pi^{-1}(\ell), \mathcal T_\ell) .\end{equation}
Lemma \ref{lem:elementaryhol} implies that $\Xi$ does not depend on the choice of an Ehresmann connection $H$, provided it is complete. Moreover, if one drops this completeness assumption, the Lemmas above can be extended, but the group morphism above is only defined at the level of germs\footnote{Let us be more precise on this point. Let $ (\mathcal V,\mathcal T) $ be a representative of the transverse singular foliation of $ \mathcal F$. By its definition (Definition \ref{def:transverse}), there exists an isomorphism of singular foliations from a neighborhood of $0$ in  $(\mathcal V,\mathcal T) $ to a neighborhood of $ \ell$ in $ (\pi^{-1}(\ell), \mathcal T_\ell) $ mapping $ 0$ to $ \ell$. Using this local isomorphism, the group morphism $ \Xi$ becomes a group morphism valued in outer symmetries of the germs of singular foliations at $0$ represented by $(\mathcal V,\mathcal T) $.}. 

This group morphism appeared first in C.L.G. and L.R.'s \cite{LGR} (last line of Definition\footnote{$\mathrm{Out}(\pi^{-1}(\ell), \mathcal T_\ell)$ is denoted by ${\mathrm{Diff}}(p^{-1}(\ell)/T_\ell)$ in the referenced work.} 2.29). It extends an older construction by Dazord \cite{Dazord}. This morphism\footnote{\cite{fischer2024classification} works in the slightly different context of formal singular foliations, the construction is however similar} is called  \emph{outer holonomy} by Simon Raphael  Fischer and C.L.G. in \cite{fischer2024classification}. It also appears for some important particular singular foliations associated to submanifolds in the works of Francis \cite{PhD_Francis,francis2023singular} and of Bischoff, del Pino, and Witte \cite{BDW}. There are analogous statement by  Rui Loja Fernandes \cite{zbMATH01836813}, and by Rui Loja Fernandes and Yvan Struchiner \cite{zbMATH07438668}, for Lie algebroids.

\begin{bclogo}[arrondi = 0.1, logo = \bcdz]{An easy mistake !}
For a regular foliation, the outer holonomy is just the usual holonomy. It is therefore tempting to believe that if it is a trivial group morphism, then there is a neighborhood of $L$ in $M$
where the singular foliation is trivial, i.e., is isomorphic to a neighborhood of $ (L,\{\ell\})$ in the direct product of $(L,\mathfrak X_c (L) $ with $  (\pi^{-1}(\ell), \mathcal T_\ell)$. For regular foliation, hence for regular leaves, this is true. 

But for a leaf of a singular foliation, this is completely wrong: let $E$ be a non-trivial orientable vector bundle, and consider the singular foliation of all vector fields on $E$ tangent to the zero section. Then the outer holonomy is trivial (because orientation preserving symmetries are inner symmetries) but the singular foliation can not be a direct product otherwise $E$ itself would be a direct product.  

However, for transversally quadratic singular foliation, it is true that outer holonomy determines $ \mathcal F$ in a formal neighborhood of $L$, see \cite{fischer2024classification} for details.
\end{bclogo}

The outer holonomy depends only on the choice of the isomorphism $ \Psi$. Its class modulo conjugation by an outer isomorphism of singular foliation on the arrival space by is however canonical.
In particular, its kernel $K_\ell$ does only depend on the choice of the point $ \ell$. Finally, this notion of holonomy has been generalized in  \cite{LGR} as a sequence of group morphisms from $ \pi_n(L)$, or all $n \geq 2 $. This construction relies on the notion of universal Lie-$\infty$ algebroid of a singular foliation, which we will encounter later in this text.

\vspace{0.5cm}

There is a similar but however slightly different approach of this ``holonomy of a leaf'' which consists of seeing it as a groupoid morphism. 

More precisely, there are two transitive groupoids over $L$. 
\begin{enumerate}
\item The fundamental groupoid $\Pi_1(L)$ is the set of all homotopy classes of paths with fixed end points \cite{Mackenzie}. It is a transitive groupoid over $L$, that we denote by $\Pi_1(L) \toto L $. By construction, arrows between $ \ell_0,\ell_1$ are the homotopy classes of paths from $ \ell_0$ to $\ell_1 $. Composition and inversion of paths amount to a Lie groupoid structure.
\item Let $(\mathcal U_L,p,H) $ be an $\mathcal F$-connection. Consider the groupoid over $L$ for which the arrows with source $\ell_0 \in L$ and target $ \ell_1 \in L$ is the set of all isomorphisms of singular foliations from  $(p^{-1}(\ell_0),\mathcal T_{\ell_0})$ to  $(p^{-1}(\ell_1),\mathcal T_{\ell_1})$. 
This groupoid admits a natural quotient if one identifies two arrows that differ by an inner symmetry. We denote by  ${\mathrm{Out}}(\mathcal U_L,p) \toto L$ this quotient groupoid.
\item[] In the construction of item 2, one could also consider germs of isomorphisms of singular foliations from  $(p^{-1}(\ell_0),\mathcal T_{\ell_0})$ to  $(p^{-1}(\ell_1),\mathcal T_{\ell_1})$, defined from a neighborhood of $ \ell_0$ to a neighborhood of $ \ell_1$. We then denote the henceforth obtained groupoid by $ {\mathrm{Out}}_{L}(\mathcal U_L,p)$.
\end{enumerate}
\noindent
We can now define the holonomy $ Hol$ as a groupoid morphism.

\vspace{0.5cm}

\begin{propositions}{The holonomy of a singular foliation \cite{LGR,fischer2024classification,Dazord}}{holonomy}  Let $L$ be an embedded leaf of a smooth singular foliation $\mathcal F$ on a manifold $M$. Let $(\mathcal U_L,p,H)$ be an $\mathcal F$-connection, with $H$ a complete Ehresmann connection.
There is a groupoid morphism, called (first) holonomy: 
$$\xymatrix{  \Pi_1(L)\ar[rr]^{Hol} \ar@<1ex>[d]\ar@<-1ex>[d]& & \ar@<1ex>[d]\ar@<-1ex>[d]
{\mathrm{Out}}(\mathcal U_L,p) \\
L \ar[rr]^{=} & &  L}$$

This morphism does not depend on the choice of $H$. Moreover, 
if no complete $ \mathcal F$-connection $H$ exists, then $Hol$ still exists, provided that $\mathrm{Out}(\mathcal U_L,p)$ is replaced by the groupoid ${\mathrm{Out}}_{L}(\mathcal U_L,p) \toto L$. 
\end{propositions}
\vspace{0.5cm}
\noindent

\begin{bclogo}[arrondi = 0.1, logo = \bcdz]{Warning !}
We warn the reader not to confuse the groupoid
 ${\mathrm{Out}}(\mathcal U_L,p) \toto L$ 
with the groupoid ${\mathrm{OutSym}}_\mathcal F \toto M$ of Section \ref{sec:holonomygroupoid}.
They seem similar, but there is an important difference. The first groupoid is over a leaf $L$ and the second one is over $M$. But even the restriction of the second one to a leaf is \emph{not} the same as the first one. In short, ${\mathrm{Out}}(\mathcal U_L,p)$ is obtained by dividing by all inner symmetries of the transversal, while ${\mathrm{OutSym}}_\mathcal F$ is obtained by dividing by inner symmetries having a very-fixed point.   
 
\end{bclogo}

Let us conclude this section by recalling that
Androulidakis and Zambon \cite{AZ2} have defined a more sophisticated holonomy, using the holonomy groupoid of the leaf instead of the fundamental groupoid of the leaf.  The group morphism\footnote{They work at the level of germs, and therefore do not need the completeness assumption.} above can be seen as a quotient of that morphism.

\chapter{Canonical geometric and algebraic structures hidden behind a singular foliation}

Throughout this chapter, we take many sophisticated ideas and break them down into simpler parts to explain the hidden structures of a singular foliation. We begin in Section \ref{sec:AlmostLieAlgebroid}  with the concept of “anchored bundles” in the context of a singular foliation, and proceed to present their morphisms and equivalences. This part concentrates exclusively on the $\mathcal{C}^\infty(M)$-module structure of the singular foliation $\mathcal{F}$. 
In Section \ref{sec:almost}, we go  further by adding a bracket  to an anchored bundle. This brings us to a concept known as the “almost Lie algebroid” associated to a singular foliation. 
This part now makes use of the Lie bracket.
Subsequently, in Section \ref{sec:isotropy}, we discuss the notion of “isotropy Lie algebra and  and holonomy Lie algebroid” of Androulidakis-Skandalis. In Section \ref{sec:Bisubmersion}, we discuss the concept of “bisubmersions”, also introduced by Androulidakis and Skandalis. These ideas help to explain how to define the “holonomy groupoid” of a singular foliation in Section \ref{sec:holonomygroupoid}. In Section \ref{sec:Geometric-resolutions}, we discuss the notion of geometric resolution of a singular foliation (again, this uses only the structure of module over functions of a singular foliation), while Section \ref{sec:Universal}  expends the notion of almost Lie algebroid over a singular foliation to something more general called the  “universal Lie $\infty$-algebroid” (or "universal $Q$-manifold) of a singular foliation.

\section{Anchored bundles over a singular foliation} 
\label{sec:AlmostLieAlgebroid}

Throughout this section $M$ is a smooth, real analytic or complex manifold\footnote{It could also mostly be an affine variety or a Zarisky open subset of $\mathbb C^d$, but we will not detail these cases}. Also, $\mathcal O$ stands for the corresponding sheaf of functions on $M$. 

\subsection{Anchored bundles}\label{sec:AnchoredBundle}

As we will see, the smooth setting is considerably simpler, and has much better properties. However, we aim to address all possible all settings, as much as we can. The reader not interested in sheaves can, in the discussion below, simply ignore the sheaf vocabulary, and consider compactly supported sections, smooth functions and vector fields instead.

We choose $M$ a manifold in the relevant category.
We denote by $ \mathcal O_\bullet$ or simply by $\mathcal O$ the sheaf of functions and by $ \mathfrak X_\bullet$ the sheaf of vector fields on a manifold $M$. 

\vspace{.5cm}

\begin{definitions}{Anchored bundle}{anchored}

An \emph{anchored vector bundle} is a pair $(A,\rho)$ made of a vector bundle $A\to M$, and a vector bundle morphism  called its \emph{anchor map}. $$\xymatrix{ A \ar[d]\ar[r]^\rho& T M\ar[d]\\ M\ar@{=}[r]&M}$$ 
\end{definitions}
\vspace{.5cm}

Readers who are exclusively interested in the smooth setting are welcome to bypass the following lines.
As we saw about Lie algebroids in Section \ref{sec:LieAlgebroidsAreSingFoliation},
$$\mathcal U \mapsto \rho(\Gamma_\mathcal U(A))$$
is a pre-sheaf in the complex and real analytic cases, but it can be sheafified (in the smooth case, it is always a sheaf, so sheafification is useless). We  denote by 
$\rho(\Gamma(A))$ this sheaf and call it the \emph{image} of $\Gamma(A) $ through~$\rho$. 
By construction, 
$\rho(\Gamma(A)) \subseteq \mathfrak{X}_\bullet$ is a sub-sheaf of  $\mathcal{O}$-module which is locally finitely generated. Furthermore, it is generated, locally, by a maximum of $\mathrm{rk}(A)$ generators.
\vspace{.5cm}

\begin{definitions}{Anchored bundle over $\mathcal F $}{anchoredoverF}

Let $\mathcal F $ be a singular foliation on $M$. We say that an anchored bundle $(A,\rho) $
\begin{enumerate}
    \item \emph{terminates within $ \mathcal F$} if $\rho(\Gamma(A)) \subseteq \mathcal F $
    \item \emph{is over $ \mathcal F$} if $\rho(\Gamma(A)) = \mathcal F $.
\end{enumerate} 
\end{definitions}
\vspace{.5cm}

Notice that anchored bundle over $\mathcal F $ could be defined for any  locally finitely generated sub-sheaf  of $\mathfrak X_\bullet$. We have not used $ [ \mathcal F, \mathcal F] \subset \mathcal F$ at this point. 

\vspace{.5cm}

\begin{questions}{Behind a singular foliation?}{ques:anchoredbundle}
Let $\mathcal F $ be  a singular foliation on $ M$. 
\begin{enumerate}
    \item Does there always exist an anchored bundle $(A, \rho) $ over $  \mathcal F $?
    \item If yes, how unique (= canonical) are they?
    \item If yes, what properties and additional structures do they have?
\end{enumerate} 
\end{questions}

\begin{remark}
For Debord  foliations (see Section \ref{sec:Debord}), an anchored bundle exists on the whole manifold $M$, by Serre-Swann theorem.
\end{remark}

\vspace{.5cm}

\begin{propositions}{Answer to the first part of Question \ref{thm:ques:anchoredbundle}}{prop:AnchoredBundle}
Let $\mathcal F$ be a singular foliation on $M$. 
\begin{enumerate}
    \item If $\mathcal F $ is finitely generated, then there exists an anchored bundle $ (A,\rho)$ over $ \mathcal F$, and $A$ can be chosen to be a trivial vector bundle.  
   In particular, an anchored bundle exists in a neighborhood of any point.
   \item
In the smooth setting, the following points are equivalent\footnote{Notice that a statement equivalent to this one was already established in Proposition \ref{prop:localglobalnumberofgenerators}. We give here an alternative proof.}.
\begin{enumerate}
    \item[(i)] $\mathcal F $ is finitely generated\footnote{See Section \ref{sec:globally}}.
    \item[(ii)] There exists an anchored bundle $(A,\rho) $ over $ \mathcal F  $.
\end{enumerate}
\end{enumerate}
\end{propositions}
\begin{proof}
Assume that $\mathcal F $ is finitely generated, and  $ X_1, \dots, X_r$ are generators. Let $A $ be the trivial vector bundle of rank $r$, i.e., 
 $$ A = M\times \mathbb K^r \longrightarrow M.$$
 Denote the canonical trivialization of $A$ by $ a_1 \dots, a_r$
and define the anchor map by 
$\rho(a_i)= X_i $ for all $ i=1, \dots, r$. We have $\rho(\Gamma(A)) = \mathcal F  $ by construction. This proves the first item of the statement. It also proves the implication $(i)  \implies  (ii)$. Let us show that  $(ii) \implies (i)$. Let $ (A,\rho)$ be as in $1$.
It is a classical theorem in smooth differential geometry that there exists a vector bundle $B \to M$ such that $A \oplus B  $ is a trivial vector bundle $E \to M$. Define a vector bundle morphism on that trivial vector bundle by 
$$ \rho_E : \, \xymatrix{ E \ar[r]^{\mathrm{pr}_A} & A  \ar[r]^{\rho} &  TM} $$
where $\mathrm{pr}_A$ is the projection onto $A$ with respect to $B$. The pair $ (E,\rho_E)$ is a trivial vector bundle such that $ \rho_E(\Gamma(E)) = \mathcal F$. In particular, $\mathcal F $ has ${\mathrm{rk}}(E)$ generators. This concludes the proof.

\end{proof}

\subsubsection{Are two anchored bundles over $\mathcal F $ really different?}

Let us define morphisms of anchored bundles – and add an equivalence class of them. 
\emph{Until the end of the present section, we work in the smooth setting, and let the complex or real analytic contexts in remarks or footnotes.}

\vspace{.5cm}

\begin{definitions}{Morphisms and Equivalences}{anchored:morphisms}
Let $(A_1 \to M_1,\rho_1)$ and $(A_2 \to M_2,\rho_2)$ be anchored bundles on a smooth manifold $M$.
\begin{enumerate}
    \item We call \emph{morphism of anchored bundles} 
    any vector bundle morphism  $\Phi\colon A_1\longrightarrow A_2$ over a map $ \phi \colon M_1 \to M_2$ making the following diagram commutative:
    \begin{equation}\label{eq:mophismAnchored} \xymatrix{ A_1 \ar[rrr]^{\Phi}\ar[rd]\ar[dd]_{\rho_2} & & &A_2\ar[dd] ^{\rho_1} \ar[ld] \\ &M_1 \ar[r]^{\phi} &M_2 & \\ TM_1\ar[rrr]^{T\phi}\ar[ru] & & & TM_2 \ar[lu] }
    \end{equation}
\item[] We speak of an \emph{isomorphism of anchored bundle} when $\Phi$ is an isomorphism of vector bundles.
 \item Two morphisms of anchored bundles $\Phi,\Phi'$ as in item 1 are said to be \emph{equivalent} if $\rho \circ (\Phi-\Phi')=0 $.
\item An \emph{equivalence of anchored bundles} is a pair of anchored bundle morphisms\footnote{In the complex or real analytic settings, an equivalence of vector bundle morphisms shall be a covering $ (\mathcal U_i)_{i \in I}$ of $M$ and an equivalence $(\Phi_i,\Psi_i) $ on each one of the open sets $\mathcal U_i $. We also assume $\Phi_i,\Phi_j $ and $\Psi_i,\Psi_j $ to be equivalent on $\mathcal U_i \cap \mathcal U_j $.} \begin{equation}
    \xymatrix{A_1\ar@<2pt>[r]^\Phi&A_2\ar@<2pt>[l]^\Psi}
\end{equation}
such that $\Psi\circ\Phi $ and $\Phi \circ \Psi $ are equivalent to the identities of $ A_1 $ and $ A_2$.

\item[] It is easily checked that both equivalences above are indeed equivalence relations on the class of anchored bundles and their sets of morphisms.
\end{enumerate}
\end{definitions}

\vspace{.5cm}
For a good understanding of the next theorem, recall 
that an anchored bundle terminates within a singular foliation $\mathcal F $ if $\rho(\Gamma(A)) \subset \mathcal F $ and is over $\mathcal F $ if $\rho(\Gamma(A)) = \mathcal F $.

\vspace{.5cm}

\begin{propositions}{The unique (up to equivalence) anchored bundle}{prop:UniqueAnchored}

Any two anchored bundles over the same singular foliation are equivalent.
\end{propositions}
\begin{proof}
Let\footnote{The reader used to an algebraic point of view can prove this statement in one sentence: "$\Gamma(A_2)$ is a projective module over functions".} $(A_1 \to M_1,\rho_1)$ and $(A_2 \to M_2,\rho_2)$ be anchored bundles over a singular foliation $\mathcal{F}$. Let $U\subset M$ be an open subset of $M$ and fix a local trivialization $ e_1, \dots,e_r$ of $A_1$. To define a $\mathcal{O}(U)$-linear map $$\varphi_U\colon \Gamma_U(A_1)\longrightarrow\Gamma_U(A_2)$$ such that $\rho_1(a)=\rho_2(\varphi(a))$ for every $a\in \Gamma_U(A_1)$, 
it suffices to map $ e_i$ to any $ f_i in \Gamma_U(A_2)$ such that $ \rho_2(f_i) =e_i$ for all $ i=1, \dots, r$, then extend by linearity. 
Likewise, we have a map $$\psi_U\colon \Gamma_U(A_2)\longrightarrow\Gamma_U(A_1)$$  In the smooth case, we use partition of unity to glue these local maps to a global one\footnote{The gluing cannot be completed in the complex setting, but it is not needed in view of the definition suggested in the footnote of Definition \ref{def:anchored:morphisms}.}.
It is straightforward that those define an equivalence.
\end{proof}

\begin{exo}
 Let $L$ be an embedded leaf of a singular foliation $L$ on a smooth manifold $M $. Let $A_L \to L$ be the holonomy Lie algebroid\footnote{See Section \ref{sec:holonomy-Lie-algebroid}.}. Show that there exists a neighborhood $ \mathcal U$ of $L$ in $M$, equipped with a submersion $p \colon \mathcal U \to L$, on which there exists an anchored bundle of the form $ (A,\rho)$ with  $A =p^* A_L$ while $ \rho$ is an anchor whose restriction  to $L$ coincides with the anchor bundle of the Lie algebroid $ A_L$.
\end{exo}

\subsubsection{Leaves of an anchored bundle}

Up to this point, we have mainly relied on the fact that $\mathcal{F}$ is a locally finitely generated module over functions. The property of stability under the Lie bracket has not played a role yet. Here is, however, a first result that makes use of leaves.

\begin{proposition}
    \label{prop:anchoredBundleLeaf}
Let $(A, \rho) $ be an anchored bundle over a singular foliation $\mathcal F $.
Any two points in the same leaf have neighborhoods on which the restrictions of $(A,\rho)$ are isomorphic.
\end{proposition}
\begin{proof}
An even stronger statement will in fact be proven, namely Proposition \ref{prop:crucial}, which immediately implies this one.
\end{proof}

\section{Almost Lie algebroids: definition and existence}\label{sec:almost}

The existence and (up to equivalence) uniqueness of an anchored bundle over a singular foliation now clarified. Now comes the third part of  Question \ref{thm:ques:anchoredbundle}: What kind of structure does this bundle possess? Here, we propose a potential candidate.

\vspace{.5cm}

\begin{definitions}{Almost Lie algebroids}{almost}
\cite{zbMATH06668113} Let  $(A,\rho)$ be an anchored vector bundle over a smooth, real analytic or complex manifold  $M$. We call \emph{almost-Lie algebroid structure} a skew-symmetric bilinear (over $\mathbb K $) map $$[\cdot\,,\cdot]_A : \Gamma(A) \wedge \Gamma(A)  \longrightarrow \Gamma(A) $$
that satisfies the \emph{Leibniz identity},
\begin{equation}\label{eq:Leibniz}
   [x, f y]_A = \rho(x)[f]y +f[x, y]_A, \quad\text{for all $x, y \in \Gamma(A), f \in \mathcal{O}(M )$}
\end{equation}
 and the \emph{anchor condition}:

\begin{equation}\label{eq:bracket-morphism}
   \rho([x, y]_A) = [\rho(x), \rho(y)], \quad\text{for all $x, y \in \Gamma(A)$}.
\end{equation}

\end{definitions}

\begin{remark}\label{rmk:Jacobi-kernel}
In the definition of an almost Lie algebroid, we do not assume  $[\cdot\,,\cdot]_A$ to satisfy the Jacobi identity, i.e., for all $x, y, z\in \Gamma(A)$, the \emph{Jacobiator} $$J(x,y,z)=[x,[y,z]_A]_A +[y,[z,x]_A]_A+[z,[x,y]_A]_A$$does not vanish. When it does, it turns to a Lie algebroid whose image through the anchor map is $\mathcal F$. However, it satisfies for any sections $x,y,z \in \Gamma(A)$:
\begin{equation}\label{AlmostJacobi}
     \rho\left(J(x,y,z)\right)=0.
\end{equation}
 \end{remark} 

The following Lemma makes almost Lie algebroids a good candidate to answer item 3 in Question \ref{thm:ques:anchoredbundle}.

\begin{lemma}\label{lemma:almost&foliation}
For every almost-Lie algebroid on $(A\to M,\rho, [\cdot\,, \cdot]_A)$, the image of the anchor map $\rho (\Gamma(A)) \subseteq \mathfrak X_\bullet $ is a singular foliation on $M$.
\end{lemma}
\begin{proof}
It is an immediate consequence of the anchor condition.
\end{proof}

We can now answer the third point of Question \ref{thm:ques:anchoredbundle}. We learned from Marco Zambon the following result (the earliest written form we found is Proposition 2.1.4 of \cite{zbMATH07105909}):

\vspace{.5cm}

\begin{propositions}{Almost Lie algebroids}{Almost} \label{prop:Almost-existence} 
Every finitely generated foliation on $M$ is the image under the anchor map of an almost-Lie algebroid.

In the smooth case, moreover,
\begin{enumerate}
    \item 
Every anchored vector bundle $(A, \rho)$ over $M$ such that $\rho(\Gamma(A)) =  \mathcal F$ can be endowed with an almost-Lie algebroid bracket. 
\item 
A singular foliation is the image under the anchor map of an almost-Lie algebroid if and only if it is finitely generated.
\end{enumerate}

\end{propositions}
\begin{proof}

Let $ X_1, \dots, X_r$ be generators of $\mathcal F $. By Lemma \ref{lem:christoff} and Exercise \ref{exo:christo1}, there exist functions\footnote{called Christoffel symbols in  Lemma \ref{lem:christoff}, Exercise \ref{exo:christo1} and in the discussion around Definition \ref{def:globallyfinitelygenerated}. See also  Proposition \ref{prop:nonminimal}. } $c_{ij}^k$ such that
$$ [X_i,X_j]= \sum_{k=1}^n c_{ij}^k X_k  \hbox{ and } c_{ji}^k = - c_{ij}^k.$$
Let $A $ be the trivial vector bundle over $M$ with fibers $\mathbb K^r $, and let $ e_1, \dots, e_r$ be the canonical trivialization of this bundle. We define the almost Lie algebroid anchor and brackets on generators by
$$ \left\{ \begin{array}{rcl} 
\rho(e_i)
&=&
X_i 
\\
\left[e_i,e_j\right]_A
& = 
& \sum_{k=1}^r c_{ij}^k e_k
\end{array} \right. $$
 and extend them using linearity (for the anchor) or Leibniz identity (for the bracket).
  This is easily checked to be an almost Lie algebroid.

The second part of the statement (i.e., the smooth case) comes from the observation that almost Lie algebroid brackets on a given anchored bundle $ (A,\rho)$ can be glued using a partition of unity. More precisely,  given an anchored bundle $ (A \to M,\rho)$, a partition of unity $(\chi_i)_{i \in I} $ relative to an open cover $(\mathcal U_i)_{i\in I} $, and almost Lie algebroid brackets $ [\cdot\,, \cdot ]_i$ (relative to $ \rho$) on  $\mathcal U_i$ for all $i \in I$, the following expression: 
 $$ [\cdot\,, \cdot ] = \sum_{i\in I} \chi_i \, [\cdot\,, \cdot ]_i   $$
 is an almost Lie algebroid bracket\footnote{This would be totally wrong with Lie algebroids: this comes from the fact that Jacobi identity is quadratic in the bracket, while the Leibniz identity and the anchor condition are linear with respect to the bracket, once the anchor is fixed.} on $\cup_{i \in I} \mathcal U_i $ - relative to the anchor $ \rho$.
\end{proof}

We call an almost Lie algebroid that terminates, as an anchored bundle, within a given singular foliation $\mathcal F $ on $M$ an \emph{almost Lie algebroid that terminates within $\mathcal F $}. Let us turn it into a category by defining morphisms. In fact, we will only deal with morphisms over the identity of $M$, which are much simpler. The subtlety is that we do not assume morphisms of almost Lie algebroid structures to be compatible with the bracket, but only to be compatible with the anchor! This is absolutely counter-intuitive, but makes perfect sense having Lie $ \infty$-algebroids\footnote{Lie $ \infty$-algebroid is introduced in Section \ref{sec:hiherstructures}} in mind.

\begin{definition}

Let $M$ be a manifold.
\begin{enumerate}
    \item We call \emph{morphism of almost Lie algebroids}  morphisms of anchored bundles – forgetting the almost Lie algebroid bracket.
    \item Two such morphisms are equivalent if and only if they are equivalent as anchored bundle morphisms.
    \item In particular, an equivalence between almost Lie algebroids is simply an equivalence of their underlying anchored bundles.
\end{enumerate}
\end{definition}

This deserves justification: why did we not require that ``morphisms''  respect the almost Lie algebroid brackets?
The answer comes from the following proposition that says that they automatically do, up to an element in the kernel of the anchor.

\begin{proposition}
\label{prop:morphismuptoker}
Let $(A_1,[\cdot\,,\cdot]_{A_1},\rho_1)$ and $(A_2,[\cdot\,,\cdot]_{A_2},\rho_2)$ be almost Lie algebroids that terminate within the same singular foliation $\mathcal F $. 
For any morphism $\Phi $ from the first one to the second one:
 $$  [ \Phi(a), \Phi(b)]_{A_2}- \Phi([a,b]_{A_1} ) \in \ker (\rho_2)$$
\end{proposition}
\begin{proof}
By definition of an almost Lie algebroid: 
 \begin{align*} \rho_2( [ \Phi(a), \Phi(b)]_{A_2}- \Phi([a,b] )) =   [\rho_2 \circ  \Phi(a), \rho_2 \circ \Phi(b)]_{A_2}- \rho_2 \circ \Phi([a,b] )) \\ 
  =   [\rho_1(a), \rho_1 (b)]_{A_2}- \rho_1 ([a,b] )) \\ 
  =0
 \end{align*}
 This proves the claim.
\end{proof}

Let us conclude this section by a theorem that follows from Propositions \ref{thm:Almost} and \ref{prop:EquivAlmost}.

\vspace{0.5cm}

\begin{proposition}
\label{prop:EquivAlmost}
Any two almost Lie algebroids over a finitely generated singular foliation $ \mathcal F$ are equivalent.    
Moreover, any almost Lie algebroids that terminates within $ \mathcal F$ admits a morphism to any almost Lie algebroids over $ \mathcal F$, and this morphism is unique up to equivalence. 
\end{proposition}

This means that, given a singular foliation $ \mathcal F$,  in the category where 
\begin{enumerate}
\item objects are almost Lie algebroids that terminates within $ \mathcal F$ and \item arrows are equivalence classes of morphisms of almost Lie algebroids,
\end{enumerate}
the terminal\footnote{As we will see, these properties are a sort of "toy-model" for the properties of the universal Lie $ \infty$-algebroid of Section \ref{sec:Universal}.}  objects are almost Lie algebroids over $ \mathcal F$. .

\begin{exo}
Let $ \varphi$ be a function on a smooth manifold $M$.
\begin{enumerate}
\item 
Show that the module $ \mathcal F_{\dd \phi}$ generated by vector fields of the form:
 $$\{   X[\phi] Y - Y[\phi] X \, | \, X,Y \in \mathfrak X_c(M)\} $$
form a singular foliation.
\item  Show that $(\wedge ^2 TM, \rho=\mathfrak i_{\dd \phi}) $ (i.e., the bivectors, equipped with the contraction by the $1$-form $\dd \phi $) is an anchored bundle over  $ \mathcal F_{\dd \phi}$  
\item Is it true\footnote{$ \mathcal L$ stands for the Lie derivative} that $ [P,Q]:= \mathcal L_{ \mathfrak \rho( P)} Q $, with $P,Q $ bivector fields,
satisfies $ \rho([P,Q])=[\rho(P),\rho(Q)]$. Is it  is an almost Lie algebroid bracket on $(\wedge ^2 TM, \rho=\mathfrak i_{\dd \phi}) $?
\item For $M=\mathbb R^n$, show that there exists an almost Lie algebroid bracket on $(\wedge ^2 TM, \rho=\mathfrak i_{\dd \phi}) $ whose restriction to constant bivector fields is given by $[\cdot,\cdot] $.
\end{enumerate}
{\emph{Hint:}} A related example is dealt with in Example 3.13 in \cite{LLS} and in Section 3.2.1 in \cite{CamilleLouis}.
\end{exo}

In \cite{AMY}, Androulidakis, Mohsen, and Yuncken introduced the Helffer-Nourrigat cone. This is a very important object, that can be described easily out of the notion of an anchored bundle, as we do in the first question of the next exercise (the subsequent questions involve notions that will be seen only later on in the text). 

\begin{exo}
\label{exo:Helffer}
Let $(M,\mathcal F)$ be a smooth singular foliation such that all regular leaves have the same dimension $r$. We denote by $ M_{\mathrm{reg}} \subset M$ the open dense subset of all regular points of $\mathcal{F}$.
\begin{enumerate}
\item Let $ (A,\rho)$ be an anchored bundle over $ \mathcal F$.
Let $ \rho^* \colon  T^* M \to A^*$ be the dual of the anchor map.
We call \emph{Helffer-Nourrigat cone, computed with respect to $ (A,\rho)$}, the closed subset of $ A^*$ given by 
 $$  {\mathrm{NH}}_{A}(\mathcal F) :=  \overline{{\mathrm{Im}}(\rho^*)|_{M_{\mathrm{reg}}} } = \overline{ \coprod_{m\in M_{\mathrm{reg}}} \rho_m^* (T_m^*M) }.$$
 The horizontal bar refers to the closure in the usual topology. 
 We denote by $ \pi :  {\mathrm{NH}}_{(A,\rho)}(\mathcal F) \to M $ the restriction to the Helffer-Nourrigat cone of the projection $ A^* \to M$.
 Show that
\begin{enumerate}
\item Show that for every $m \in M_{\mathrm{reg}}$, the fiber of $ \pi$ over $M$ coincides\footnote{We use the symbol $ {}^\perp$ for the annihilator.} with $ {\mathrm{Im}}(\rho^*_m) = ({\mathrm{ker}}(\rho_m) )^{{\perp}}$.
\item Show that for every $m \in M$, the fiber of $ \pi$ is non-empty, and is contained the annihilator of the strong kernel\footnote{See Section \ref{sec:isotropy}.} of $ \rho$ at $m$.
\item Show that for every $m \in M$, the fiber of $ \pi$ is a union of sub-vector spaces of $ A_m^*$, all of dimension $r$.

\item (\textit{For Poisson geometers}) Show that if $(A,\rho)$ admits a Lie algebroid bracket, the Helffer-Nourrigat cone is a union of symplectic leaves of the  Poisson structure on $A^*$ associated to the Lie  algebroid bracket.
\end{enumerate}
\item[] We will now explain how the Helffer-Nourrigat cone can be seen as a subset of $\coprod_{L \in {\mathrm{Leaves\, of\,}} \mathcal F } A_L^* $.
Here $A_L$ is the holonomy Lie algebroid of a leaf $L $, defined in \ref{sec:holonomy-Lie-algebroid}.
\item Use question 1.b  to show that for every $m \in M$, the fiber of $ \pi$ over $m$ can be included  into $ \left. A_L^* \right|_{m}$. {\emph{Hint}}: We recall that  $\left. A_L \right|_{m} $ is the quotient of $A_m$ by the strong kernel of $ \rho$ at $m$.
\item Show that the image of the inclusion in question 2. is made of the union, for all $m \in M$, of the annihilator of all limit subalgebras\footnote{See Definition \ref{def:limitLie}. The limit Lie algebras being vector subspaces of the isotropy Lie algebra $\mathfrak g_m(\mathcal F) $, which is  included into $ \left. A_L \right|_{m}$,  can be seen as vector subspaces of $ \left. A_L \right|_{m}$. Since the  limit Lie algebras are precisely the points of the Nash resolution (which is a sub-set of the Grassmannian of  $\mathfrak g_m(\mathcal F) $, see Section \ref{sec:Nash}), we can also say the fiber of the Helffer-Nourrigat cone over a point $m$ is the union of all annihilators of the “points” in fiber over $m$ of the Nash resolution.} at $m$.
\item[] We denote this set by $ {\mathrm{HN}}(\mathcal F)$, without reference to a particular anchored bundle, since by the previous two questions, it does depend on the choice of an anchored bundle. We call it the \emph{Helffer-Nourrigat cone} of $ \mathcal F$.

\end{enumerate}
The presentation above seems different from the presentation done in \cite{AMY} but the difference is only a difference of presentation, see \cite{Ruben3}. We acknowledge discussions with Fani Petalidou and Mohsen Masmoudi when writing this exercise.
\end{exo}

\begin{exo}
Show that the Helffer-Nourrigat cone for the singular foliation on $ M=\mathbb R^n$ made of all vector fields vanishing at the origin is given as follows:
\begin{enumerate}
    \item if $m$ is not the origin, the fiber of the Helffer-Nourrigat cone over $m$ is $\simeq \mathbb R^n $, 
    \item if $m$ is the origin, the fiber is made of square $n \times n $ matrices of rank $\leq 1$.
\end{enumerate}
We thank Cédric Rigaud for this exercise.
\end{exo}

\subsection{An alternative proof of Proposition \ref{prop:symInt}}

We use the notions of anchored bundle and the almost Lie algebroids to give a much simpler proof of a result that was crucial to establish the existence of leaves:
Proposition \ref{prop:symInt}.
This proof is inspired from a proof that appeared in \cite{zbMATH07105909} by Henrique Bursztyn, Hudson Lima, and Eckhard Meinrenken for almost-Lie algebroids and internal symmetries. The result that we will prove, and which immediately implies Proposition \ref{prop:symInt}, is the following proposition\footnote{We work in the smooth setting in the present section: it can be adapted to the complex and algebraic settings, but not to the algebraic one.}.

\begin{proposition}
\label{prop:crucial}
Let $\mathcal F $ be a singular foliation on $M$. Choose $ Y \in \mathfrak X (M)$ be a vector field such that $ [Y,\mathcal F] \subset \mathcal F$.
For any open subset $ \mathcal U \subset M$ on which the time $1$-flow $\phi_1^Y $ of $Y$ is well-defined, the restrictions of any anchored bundle $ (A,\rho)$ over $M$ to $\mathcal U $  and $\phi_t^Y(\mathcal U )=  $ are isomorphic\footnote{We insist on “isomorphic” and not "equivalent".  Of course, it is part of the statement that the base map of the isomorphism is $\phi_1^X $}.
\end{proposition}

These isomorphisms can be seen, when  $A$ is a trivial bundle, as families indexed by $m \in M$ of invertible matrices as in Proposition \ref{prop:symInt}. Proposition  \ref{prop:crucial} is therefore a generalization of Proposition \ref{prop:symInt}.

The proof is based on the notion of linear vector field. A vector field $Y $ on a vector bundle $E \stackrel{p}{\to} M$ is said to be \emph{linear} if one of the following equivalent conditions holds:
\begin{enumerate}
    \item[(i)] For any function $f$ on $E $ whose restriction to any fiber of $p \colon E \to M$ is a polynomial of degree $\leq k$, $Y[f] $ is a polynomial of degree $\leq k $. 
    \item[(ii)]
    $ Y [p^* \mathcal O] \subset p^* \mathcal O$ and $Y [ \Gamma(A^*) ] \subset  \Gamma(E^*)$, with the understanding that $ \Gamma(A^*)$ must be considered as a smooth function on $E$  linear on each fiber of $ p \colon E \to M$.
    \item[(iii)] In any local coordinates $ (x_1, \dots, x_n, y_1, \dots, y_r)$ on $E$, with $(x_i)_{i=1, \dots, n}$ being local coordinates on the base manifold $M$, and $(y_j)_{j=1, \dots, r}$ linear coordinates on the fibers, the vector field $Y$ is of the form:
$$ Y = \sum_{i=1}^n A_i(x_1, \dots, x_n ) \frac{\partial}{\partial x_i} + \sum_{i,j=1}^n   B_{ij}(x_1, \dots, x_n ) y_i \frac{\partial }{ \partial y_j} .$$
\end{enumerate}
Linear vector fields on $E\to M$ are stable under Lie bracket.
Moreover, for any linear vector fields $Y$, there is a unique vector field\footnote{i.e., a unique vector field $p$-related to $Y$, if one uses the terminology of Section \ref{sec:pull-back}.} $ p_* {Y} \in \mathfrak X(M)$ such that
 $$  Y[p^* F] = p^* (p_* {Y}) [F]  \hbox{ for every $ F \in \mathcal O_M$}, $$
 and the assignment $ Y \mapsto p_* {Y}$ is a Lie algebra morphism that we will call projection.

\begin{lemma} \label{linearlemma}
Let $E \stackrel{p}{\to} M$ be a vector bundle and $X \in \mathfrak X(M) $ be a vector field. 
For any linear map:
 $$ \delta_X \colon \Gamma (E^*) \longrightarrow \Gamma(E^*)  $$
 such that for every function $f\in \mathcal O_M $ and every $\epsilon \in \Gamma(E^*) $:
  \begin{equation}\label{condition:lnear_vector_field} \delta_X (f \epsilon) = f \delta_X ( \epsilon) + X[f] \, \epsilon  \end{equation}
 there exists a unique linear vector field on $E$ that projects on $X$ and whose restriction to fiberwise linear functions on $E$ is $ \delta_X$. 
\end{lemma}

The Lie algebra of linear vector fields on $ E\to M$ can be seen as the Lie algebra of the group of vector bundle isomorphism of $ E\to M$. Below is a more precise statement, that we leave to the reader:

\begin{lemma}
\label{lem:FlowOfLinear}
Let $Y$ be a linear vector field on a vector bundle $ E \stackrel{p}{\to} M$,
and $p_* Y$ the vector field on $M$ to which it is $p$-related.
The flow $\phi_t^Y $ at time $t$ of a linear vector field is defined if and only the flow at time $t$ of its projection $p_*Y$ on $M$ is defined. In that case, it is vector bundle isomorphism
$$ \xymatrix{ E \ar[d] \ar[rr]^{\phi_t^Y}&& E\ar[d] \\M \ar[rr]^{\phi_t^{p_*  Y}} &&M } .$$
\end{lemma}

\begin{proof}[Proof of Proposition \ref{prop:crucial}]
Let $Y$ be a vector field satisfying $[Y,\mathcal F]\subset \mathcal F$. Let us first observe that the classical Lie derivative:
$$ \begin{array}{rrcr} \delta_Y^{TM} \colon   &\Gamma(T^* M) & \to  & \Gamma(T^* M) \\ &\alpha &\mapsto& L_Y \alpha \end{array} $$
satisfies condition \eqref{condition:lnear_vector_field}  
and therefore defines a linear vector field $ \widehat{Y}^{TM} $ on $ TM \stackrel{p}{\to} M$ which is $p$-related to $Y$. This vector field is called the tangent prolongation of $Y$. Our goal is to construct a linear vector field on $A$ which is $\rho$-related with $\hat Y^{TM}$.

Let us first assume that $A=M\times \mathbb R^r$ is a trivial bundle. 
\begin{itemize}
\item In this case, there is a matrix-valued smooth function ${\mathfrak Y}$ on $M$ such that $[Y,\rho(e_i)]=\sum_{j=1}^r {\mathfrak Y}_{i,j}\rho(e_j)$ for every $ i\in \{1, \dots, r\}$. This allows to define  $\gamma_Y^A:\Gamma(E)\to \Gamma(E)$ by $$\gamma_Y^A\left(\sum_{i=1}^r f_i e_i\right)=\sum_{i=1}^r Y\left[f_i\right] \, e_i+\sum_{i,j=1}^r f_i {\mathfrak Y}_{i,j}\, e_j$$
for every $r$-tuple of functions $ f_1, \dots, f_r$.
This construction imposes that $\gamma_Y^A$ adheres to Equation~\eqref{condition:lnear_vector_field}.
    
\item Next, we construct a “dual” of $\gamma_Y^A$, that we denote by $\delta_Y^{A}:\Gamma(A^*)\to \Gamma(A^*)$, by imposing the following duality relation: 
\begin{align}\label{eq:dualizederivation}
   \left\langle  \delta_Y^{A}(\alpha) \middle| e\right\rangle :=\gamma_Y^A \left[ \langle \alpha,e \rangle \right]- \left\langle \alpha , \gamma_Y^A(e) \right\rangle.
\end{align}
for every $e\in \Gamma(A), \alpha \in \Gamma(A^*)$. This definition requires a justification: the right-hand term in Equation \eqref{eq:dualizederivation} is $ \mathcal C^\infty(M)$-linear in $e$, so that there exists a section $ \delta_Y^{A}(\alpha)$ that satisfies  Equation \eqref{eq:dualizederivation}. Moreover, $ \alpha \to  \delta_Y^{A}(\alpha)$ satisfies for any $f \in \mathcal C^\infty(M)$ the relation  

$\delta_Y^{A}(f\alpha)=Y(f)\alpha+f\delta_Y^{A}(\alpha)$.
\item In view of Lemma \ref{linearlemma} therefore, $\delta_Y^{A}$ defines a linear vector field $\hat Y^A$ $p$-related to $Y$. Let us check that $\hat Y^A$ is $\rho$-related to $\hat Y^{TM}$, i.e.,
   $$   \hat Y^A \circ \rho^*  (F) = \rho^*  \circ \hat Y^{TM} (F)$$
for every $F \in \mathcal C^\infty(T M)$. 
For $F \in  p^* \mathcal C^\infty(M) $, the condition holds automatically true. 
It suffices therefore to check  that the condition holds true when applied on a fiberwise linear function, i.e., a $1$-form $ \eta \in \Omega^1 (M)$. 
In view of the definition of  $\hat Y^A$  and $ \hat Y^{TM}$ on such function, it therefore suffices to check that
$$ \delta_Y^A \circ \rho^*(\eta)=\rho^*\circ L_Y (\eta) $$
 Let us check this relation: 
for every $\eta\in\Gamma(T^*M),a \in \Gamma(A)$,  we have by construction
\begin{align*}
    \left\langle \delta_Y^A(\rho^*\eta),a \right\rangle &=Y\left[\langle \rho^*\eta,a\rangle \right]-\langle \rho^*\eta,\gamma_Y^A (a)  \rangle \\
    &=Y\left[ \langle \eta, \rho(a)\rangle \right] - \langle \eta, \rho(\gamma_Y^A (a))\rangle \\
    &= L_Y  \circ \iota_{\rho(a)} \, \eta -\iota_{\rho(\gamma_Y^A(a))} \, \eta \\
     &=  \iota_{\rho(a)} \circ L_Y   \, \eta + \iota_{[Y,\rho(a)] -\rho(\gamma_Y^A(a))} \, \eta 
\end{align*}
We therefore have
 $$    \left\langle \left(\delta_Y^A \circ \rho^*- \rho^* \circ L_Y \right) \eta ,a  \right\rangle  =  \left\langle \eta , [Y,\rho(a)] -\rho\left(\gamma_Y^A(a)\right)  \right\rangle $$
By definition of $\gamma_Y^A$, the previous relation holds true if  $a=e_i $ belongs to the canonical trivialization of the trivial bundle  $A \to M $.
By $ \mathcal C^\infty(M)$-linearity, it holds true for any $a \in \Gamma(A)$.
This proves that $\hat Y^A $ is $ \rho$-related with $ \hat Y^{TM}$.

\end{itemize}

Let us now turn to the case of a non-trivial anchored bundle $(A\to TM, \rho)$. By choosing a second vector bundle $E$ such that $A'=A\oplus E$ is a trivial vector bundle and setting $\rho_{A'}=\rho_A \circ \pi_A$, with $\pi_A $ the projection onto $A$, we obtain a new anchored bundle. Since it is trivial, we can apply the previous considerations and obtain $\delta^{A'}_Y$. We now define $\delta_A^Y:\Gamma(E^*)\to \Gamma(E^*)$ by:
$$
\delta^A_Y(\alpha)(e)=\delta^{A'}_Y(\pi_A^*\alpha)(i_{A,*}e),
$$
where $\pi_A$ and $i_A$ are the projection and inclusion of $A\subset A'$. This map can also be obtained from $\pi_A\circ \gamma_Y^{A'}\circ i_A$ with the dualization procedure, as in Equation \eqref{eq:dualizederivation}.
We can now verify:

\begin{align*}
    \delta^A_Y\rho_A^*\beta(e)&=\delta^{A'}_Y(\pi_A^*\rho_A^*\beta)(i_{A,*}e)
    =\delta^{A'}_Y(\rho_{A'}^*\beta)(i_{A,*}e)\\
    &=(\rho_{A'}^*L_Y\beta)(i_{A,*}e)
    =(i_{A}^*\rho_{A'}^*L_Y\beta)(e)\\
    &=(\rho_{A}^*L_Y\beta)(e)
\end{align*}
i.e., the corresponding linear vector field $\hat Y^{A}$ is $\rho$-related to $\hat Y^{TM}$. In particular, where defined, the flow of $\hat Y^{A}$ induces an isomorphism of anchored bundles.
\end{proof}

\section{Isotropy Lie algebra and holonomy Lie algebroids}\label{sec:isotropy}

Let us use the almost Lie algebroids associated to a singular foliation in the previous section to associate a Lie algebra, called isotropy Lie algebra, to any point of a singular foliation. We then relate it with the original definition of Androulidakis and Skandalis.

We work here in the smooth setting. All arguments can easily be adapted to the complex or real analytic ones.

\subsection{Kernel and Strong-kernel of a morphism of vector bundles}

Consider a vector bundle morphism over  the identity of $M$:
$$  \xymatrix{
B \ar[d] \ar[r]^{\Phi}
& C \ar[d] \\
M \ar@{=}[r]
& M } 
$$
Choose a point $m \in M$.
There are two subspaces in $B_m  $ that deserve to be called  ``kernels''.
\begin{enumerate}
    \item the usual kernel $ {\mathrm{ker}}(\Phi_m)  $, i.e., 
     $$  \{ u \in B_m  \, | \, \Phi_{|_m}(u)=0\},$$
    \item and there is the \emph{strong kernel}, i.e., the subspace $\mathrm{Sker}(\Phi,m)\subset B_m $ of all elements through which there is a neighborhood $ \mathcal U$ of $m$ in $M$ and a local section in the kernel of 
     $  \Phi\colon \Gamma_{\mathcal U} (B) \to \Gamma_{\mathcal U} (C)  $. In equation:
     $$ \mathrm{Sker}(\Phi,m):= \{ u \in B_m \hbox{ s.t. } \exists U \in \Gamma(B) \hbox{ with } \Phi(U)=0 \hbox{ and }  U_{|_m} = u \}. $$ 
\end{enumerate}
Of course, there is an inclusion:
 $$\mathrm{Sker}(\Phi,m) \subset \mathrm{ker}(\Phi|_{m}). $$
Moreover, the dimensions of the distributions have opposite behavior:
\begin{enumerate}
    \item the map $m \mapsto {\mathrm{dim}}(\mathrm{ker}(\Phi_{|_m}) ) $ is upper semi-continuous, i.e., if a sequence $(x_n) $ in $M$ has limit $x$, then 
     $$ {\mathrm{dim}}(\mathrm{ker}(\Phi_{|_{x}})) \geq \hbox{ upper limit of } {\mathrm{dim}}(\mathrm{ker}(\Phi_{|_{x_n}})) $$
    \item the map $m \mapsto {\mathrm{dim}}(\mathrm{Sker}(\Phi,m ) $ is lower semi-continuous, i.e., if a sequence $(x_n) $ in $M$ has limit $x$, then 
     $$ {\mathrm{dim}}(\mathrm{Sker}(\Phi,x)) \leq \hbox{ lower limit of } {\mathrm{dim}}(\mathrm{Sker}(\Phi,{x_n})).$$
\end{enumerate}

In particular, if kernel and strong kernel coincide at a point, they coincide in a neighborhood of that point. In particular, their dimensions are constant at all points in a neighborhood, so that they form a vector sub-bundle. Let us state this conclusion for future reference.

\begin{proposition}
Let $\Phi \colon B \to C $ be a vector bundle morphism over the identity of $M$. For any $m \in M$, the following two assertions are equivalent:
\begin{enumerate}
  \item the kernel and the strong kernel coincide at $m$. 
    \item There is a neighborhood $ \mathcal U$ of $m$ in $M$ on which  the kernel and the strong kernel coincide at all points. 
\end{enumerate}
In this case, moreover, these coinciding kernels form a sub-vector bundle of  the restriction to $\mathcal U $ of $B$. 
\end{proposition}

\subsection{The isotropy Lie algebra (I: the space)}

Let $\mathcal F $ be a singular foliation on a manifold $M$.
Let $\mathcal U $  be an open neighborhood of $m$ on which $\mathcal F $ is finitely generated. 
In view of Proposition \ref{thm:Almost}, there exists an anchored bundle $(A \to \mathcal U, \rho) $ over $\mathcal F $. 
 We call \emph{isotropy vector space at $m$} the quotient space:
 $$ \mathfrak g_m (\mathcal F) =
 \frac{\mathrm{ker}(\rho_{m})}{\mathrm{Sker}{(\rho,m)}} . $$
Notice that the notation $\mathfrak g_m (\mathcal F)$ makes no reference to the chosen anchored bundle. This is justified by the following proposition:

\vspace{.5cm}

\begin{propositions}{The isotropy vector space at $m$ makes sense}{prop:canonicalSpace}

Let $\mathcal F $ be a singular foliation. The isotropy vector spaces associated to any two anchored bundles are canonically isomorphic\footnote{It makes sense, therefore, to denote it by $ \mathfrak g_m (\mathcal F)$.}.
\end{propositions}
\begin{proof}
This is an immediate consequence of Theorem \ref{thm:prop:UniqueAnchored}, since a morphism of anchored bundle morphisms over $M$ maps kernel to kernel and Strong kernel to Strong kernel, and since homotopic morphisms  induce the same map  at the quotient level. 
\end{proof}

Here is an important theorem, due to \cite{AS}.

\vspace{.5cm}

\begin{theorems}{Ranks and dimensions}{RankDimension}
Let $(M,\mathcal F)$ be a singular foliation.
For every $ m \in M$, 
\begin{enumerate}
    \item the rank ${\mathrm{rk}}_m(\mathcal F)$ of $\mathcal F $ at $ m$ (i.e., the minimal number of local generators),
    \item the dimension $ \dim (L_m)$ of the leaf through $m$, 
    \item and the dimension $\dim \left( \mathfrak g_m ({\mathcal F}) \right)
$ of the holonomy vector space\footnote{that will be soon equipped with a Lie bracket making it Androulidakis-Skandalis isotropy Lie algebra.} at $m$,
\end{enumerate}
are related by the relation
 $$ {\mathrm{rk}}_m(\mathcal F) = \dim \left( \mathfrak g_m ({\mathcal F}) \right)+ \dim (L_m) .$$
\end{theorems}

\begin{proof}
For simply notations, we set $ r= {\mathrm{rk}}_m(\mathcal F)$, $ g:= \dim \left( \mathfrak g_m ({\mathcal F}) \right) $, and $ \ell =  \dim (L_m) $.
Let $Y_1, \dots, Y_r  $ be a minimal family local generators of $ \mathcal F$ near $m$. Without any loss of generality, one can assume that $ Y_1|_{r-\ell+1}, \dots, Y_{r}|_m$ form a basis of $ T_m L_m$ while $Y_{1}, \dots, Y_{r}-\ell $ are all zero at the point $m$.

Since these vectors vanish at $m$, 
the classes $[Y_1], \dots, [Y_{r-\ell}]$ of these vectors in $\mathfrak g_m (\mathcal F) $ are well defined, and, altogether, these classes form a system of generators of $\mathfrak g_m(\mathcal T) $. Hence $ g \leq r-\ell$.
 It remains to show that the classes $[Y_1], \dots, [Y_{r-\ell}] \in \mathfrak g_m(\mathcal T)$ are linearly independent. Assume that they are not, i.e., (without loss of generality):
$[Y_1]=\sum_{i=2}^{r}\alpha_i[Y_i]$ for some $\alpha_i\in\mathbb K$. 
This means that 
$ Y_1=\sum_{i=2}^{r}\alpha_iY_i + Y$ for some $Z$ in $\mathcal I_m\mathcal F$. 
This element $Z$ reads as a sum $Z=\sum_{i=1}^r g_i Y_i$ for some functions $g_i$ that vanish at $m$, so that we obtain 
$$(1-g_1)Y_1=\sum_{i=2}^{r-\ell}(\alpha_i+g_i)Y_i+ \sum_{i=1}^\ell g_i Y_{r-\ell+i}$$
Now, there exists a neighborhood of $m$ where we can invert $(1-g_1)$. On that neighborhood,  $Y_1$ is  a linear combination of the other $Y$'s which contradicts the minimality of the family $Y_1,...,Y_r$. This concludes the proof. 
\end{proof}

\subsection{The isotropy Lie algebra II: the bracket} 

Now, let $\mathcal F $ be a singular foliation, $m$ is  point, $\mathcal U$ an open subset containing $m$ and $ (A \to \mathcal U , \rho)$ an anchored bundle of $\mathcal F $ restricted to $ \mathcal U$.

According to Proposition \ref{thm:Almost}, $(A \to \mathcal U , \rho)$ can be equipped with an almost Lie algebroid bracket\footnote{At least in the smooth setting. In complex or real analytic setting, one may have to restrict to a smaller open neighborhood.} that we denote by $$[\cdot\,, \cdot]_A  \colon \Gamma(A) \times \Gamma(A) \longrightarrow \Gamma(A). $$
 Let us fix a point $ m \in M$. Consider two elements $a,b \in A_m $. For any two sections $ \tilde a, \tilde b $ of $A \to \mathcal U $ through $a,b$. In view of the Leibniz identity, the value at $m$ of the almost Lie algebroid bracket $ [\tilde a, \tilde b]_A $ depends on the $1$-jet at $m$ of the sections $ \tilde a, \tilde b $.
 However, if $ \rho(a)=0$, then  \begin{align*} [\tilde a,f \tilde b]|_{m}= f(m) [\tilde a,\tilde b]|_{m} + \rho(a)[f] b \\  =f(m) [\tilde a,\tilde b]|_{m} \end{align*}
 for any local function $f $. This implies that $  [\tilde a,\tilde b]|_{m}$ depends only on the value of the section $\tilde b$ at $m$, i.e., depends only on $b$. As a consequence, if $ a,b \in {\mathrm{Ker}}(\rho_m)$, then  $  [\tilde a,\tilde b]|_{m}$ depends only on $ a$ and $b$ so that the bracket $ [\cdot\,, \cdot]_A$  induces a bilinear map
  $$ [\cdot\,,\cdot ]_{A,m} \colon  \wedge^2 \mathrm{ker}(\rho_m) \longrightarrow A_m  $$ 
 given for all $a,b \in \mathrm{Ker}(\rho_m) $ by
 $$[ a,b ]_{A,m} = [\tilde{a},\tilde{b}]_{A}(m)$$
for any sections $\tilde{a},\tilde{b}$ through $a$ and $b$.
Moreover, $[ a,b ]_{A,m}$ is in fact valued in $\mathrm{ker}(\rho_m)$: this follows easily from the anchor condition.
Lastly, the anchor condition implies that the strong kernel at $m$ is an ``ideal'' of that bracket, i.e.,
 $ [{\mathrm{Sker}}(\rho,m),{\mathrm{ker}}(\rho_m) ]_{A,m} \subset {\mathrm{Sker}} (\rho,m)$
 so that the skew-symmetric bilinear map $ [ \cdot\,, \cdot  ]_{A,m} $ goes to the quotient to a bilinear map
  \begin{equation}\label{eq:isLie} [\cdot\,, \cdot]_m \colon  
  \wedge^2 \mathfrak g_m(\mathcal F) \longrightarrow
   \mathfrak g_m(\mathcal F)
  \end{equation}

\vspace{0.3cm}

\begin{propositions}{The Lie bracket}{braisotropy}
The bilinear map \eqref{eq:isLie}:
\begin{enumerate}
    \item is a Lie bracket on the holonomy vector space $\mathfrak g_m(\mathcal F) $,
    \item is canonically defined, i.e., does not depend on the choice of an anchored bundle and of an almost Lie algebroid bracket.
\end{enumerate}
\end{propositions}
\begin{proof}
The Jacobi identity follows from the fact the Jacobiator (see Remark \ref{rmk:Jacobi-kernel}) of three sections of the bracket $[\cdot\,,\cdot]_A$ lies in the kernel of $\rho\colon \Gamma(A)\to \mathfrak X(M)$, see Equation \eqref{AlmostJacobi}. Now, for any two almost Lie algebroid brackets $[\cdot,\cdot]  $ and $[\cdot,\cdot]' $ on $A $, we have, for any two sections $\tilde{a}, \tilde{b} $
 $$  \rho \left( [\tilde{a}, \tilde{b}]_A' - [\tilde{a}, \tilde{b}]_A \right) =0$$
 so that $[\tilde{a}, \tilde{b}]_A' - [\tilde{a}, \tilde{b}]_A$ is valued in the Strong kernel of $ \rho$. This implies that the induced bracket \eqref{eq:isLie}
 does not depend on the choice of an almost Lie algebroid bracket on a given anchored bundle $(A,\rho)$.
 More generally, given two anchored bundles $ (A, \rho)$ and $(A',\rho') $, Proposition \ref{thm:prop:UniqueAnchored} implies the existence of anchored bundle morphisms $\Phi : (A,\rho) \longrightarrow (A',\rho') $.
 For any two almost Lie algebroid structures on $A$ and $A'$, we have by Proposition \ref{prop:morphismuptoker}
  $$ \rho' \left( \Phi([\tilde{a}, \tilde{b}]_A) -[\Phi(\tilde{a}), \Phi(\tilde{b})]_{A'}\right)=0  $$
  so that $\Phi([\tilde{a}, \tilde{b}]_A) -[\Phi(\tilde{a}), \Phi(\tilde{b})]_A$ is valued in the Strong kernel of $ \rho'$ at $m$. This implies that $ \Phi$  induces a Lie algebra morphism $$ 
(\Phi,{m}) \colon \frac{\mathrm{ker}(\rho_{m})}{\mathrm{Sker}(\rho,m)}  \longrightarrow
 \frac{\mathrm{ker}(\rho_{m}')}{\mathrm{Sker}(\rho',m)} .$$
 Equivalent morphism would induce the same Lie algebra morphism.
 The same construction, applied to an “inverse” $\Psi \colon A' \to A $ as in Proposition \ref{thm:prop:UniqueAnchored}, gives an inverse map to that Lie morphism. The result follows
 \end{proof}

\begin{definitions}{Isotropy Lie algebra: definition}{istropyLieA}
We call \emph{isotropy Lie algebra of $\mathcal F $ at $m \in M$} the pair $ (\mathfrak g_m(\mathcal F), [\cdot, \cdot])$ with $\mathfrak g_m(\mathcal F)$ the vector space defined as in Proposition \ref{thm:prop:canonicalSpace} $ [\cdot, \cdot]$ as in Proposition \ref{thm:braisotropy}.
\end{definitions}
We end this subsection with the following lemma.

\begin{lemma}
\label{lem:almostminimalrank}
    Let $(M, \mathcal{F})$ be a singular foliation and $x\in M$ such that $T_x\mathcal{F}=\{0\}$. Then, $\mathrm{rk}_x(\mathcal{F})=\dim \mathfrak g_x(\mathcal{F})$ and $\mathcal{F}$ admits on an open neighborhood $\mathcal{U}$ of $x$ an almost Lie algebroid structure on the trivial bundle $A=\mathcal{U}\times \mathfrak g_x(\mathcal{F})$ whose anchor is $0$ at the point $x$ and whose bracket is given by the Lie bracket of $\mathfrak g_x(\mathcal{F})$ at this point.
\end{lemma}
\begin{proof}
 If $T_x\mathcal{F}=\{0\}$, then  the leaf through $x$ is reduced to a point. By Theorem \ref{thm:RankDimension}
(3), $\mathrm{rk}_x(\mathcal{F})=\dim \mathfrak g_x(\mathcal{F})$. It follows that $\mathcal{F}$ admits a minimal set of generators $X_1,\ldots, X_k$ in a neighborhood $\mathcal{U}$ of $m$ with $k=\dim \mathfrak g_x(\mathcal{F})$ elements. that induce a basis $e_1, \dots, e_k $ of $\mathfrak g_x(\mathcal{F})$. Hence, $\mathcal{U}\times \mathfrak g_x(\mathcal{F})$ admits an almost Lie algebroid whose image is $\mathcal{F}|_{\mathcal{U}}$.
The bracket is constructed as follows. The functions $ c_{ij}^k$ such that
 $$ [X_i,X_j]=\sum_{k} c_{ij}^k X_k  $$
are not unique but their values $ c_{ij}^k(x)$ at $x$ have to coincide with the Christoffel symbol of the Lie algebra $\mathfrak g_x(\mathcal{F})$ in the basis $ e_1, \dots,e_k$. The functions $ c_{ij}^k$ then define an almost Lie algebroid structure as in Section \ref{sec:AlmostLieAlgebroid} which satisfies the required properties. 
\end{proof}

\begin{exo}
Let  $m$ be a point in a foliated manifold $ (M,\mathcal F)$
Show that the following points are equivalent:
\begin{enumerate}
    \item[(i)] $m$ is a regular point,
    \item[(ii)] $ \mathfrak g_m(\mathcal F) =0$. 
\end{enumerate}
\end{exo}

\begin{exo}Let $(M,\mathcal F)$ be a foliated manifold and $m\in M$ a point. Show the following points:
\begin{enumerate}
    \item Let $U$ be an open subset of $M$ containing $m$, then $\mathfrak{g}_m(\mathcal F)= \mathfrak{g}_m(\mathcal F|_U)$.
    \item Let $(\tilde M, \tilde {\mathcal F})$ be another foliated manifold and $\tilde m\in \tilde M$. Then $\mathfrak{g}_{(m,\tilde m)}(\mathcal F\times \tilde{ \mathcal F})= \mathfrak{g}_{m}(\mathcal F)\oplus \mathfrak{g}_{(\tilde m)}(\tilde{\mathcal F})$.
    (The direct sum should be understood as a Lie algebra direct sum).
    \item For $ \varphi \colon N \mapsto M $ a surjective submersion, $\mathfrak {g}_{n}(\varphi^{-1}(\mathcal F)) = \mathfrak{g}_{m}(\mathcal F)  $ for every $ n \in \varphi^{-1}(m)$.
\end{enumerate}
\end{exo}

\begin{exo}
Several exercises in the Section \ref{examples-of-SF} consisted into computing explicitly the isotropy Lie algebra at a given singular point.  We invite the reader to look for the word “isotropy” is that section.
\end{exo}

\begin{exo}
Let $m$ be a point of a foliated manifold $(M,\mathcal F) $.  
\begin{enumerate}
\item Show that the isotropy Lie algebra $\mathfrak g_m(\mathcal F)$ at $m$ is canonically isomorphic to the singular foliation defined by the vectors fields as in item b) of the splitting theorem \ref{thm:localsplitting3}. 
\item  Show that the isotropy Lie algebra at any point  of $L$ is isomorphic to the isotropy Lie algebra at the origin of any representative of the transverse singular foliation (see Definition \ref{def:transverse}). \emph{Hint:} This in fact the same question as the previous one, with different wording.
\item Let $L$ be the leaf through $m$. Show that the isotropy Lie algebras at two points in $ L$ are isomorphic (\emph{Hint:} Use Proposition \ref{thm:leavesAreBoring}). 
\item Is this isomorphism canonical?
\item \emph{Hard\footnote{Except for the reader that knows 1° the holonomy Lie algebroid, see Section \ref{sec:holonomy-Lie-algebroid}, and 2) that there is a similar phenomenon for the kernels of two points on any transitive Lie algebroid.}} Show that if the leaf $L$ is simply connected, then this isomorphism is canonical up to an inner Lie algebra automorphism.
\item \emph{Hard.} Construct a group morphism $ \pi_1(L,\ell) \to {\mathrm{Out}}( \mathfrak g_m(\mathcal F))$, where ${\mathrm{Out}}$ is the Lie group of outer automorphism. \emph{Hint}: this of course related to the outer holonomy constructed in Equation \eqref{eq:holonomymorphism}. Can be also seen as a general result on transitive Lie algebroids, then one applies the results of Section \ref{sec:holonomy-Lie-algebroid}.
\end{enumerate}

the isotropy Lie algebras of  $\mathcal F$ at $m$ and the  isotropy Lie algebra of any representative of the transverse singular foliation (as in Theorem \ref{thm:transverseFol}) at $\ell$ are canonically isomorphic.  (\emph{Hint:} use the splitting theorem of Section \ref{sec:splitting})
\end{exo}

\subsection{Androulidakis-Skandalis construction of the isotropy Lie algebra}

 Let us now present the original definition of the isotropy Lie algebra of a singular foliation by Androulidakis and Skandalis. 
 
Let $\mathcal F $ be a singular foliation, defined as Definition\footnote{Again, the interested reader will easily adapt these arguments to the real analytic and complex settings.} \ref{def:consensus}. For every $m \in M$, we consider $\mathcal F_m  \subset \mathcal F$ the sub-Lie algebra of vector fields in $\mathcal F $ vanishing at $m$.   
Let $ \mathcal I_m$ be the ideal of functions vanishing at $m$. 
There is an inclusion $ \mathcal I_m \mathcal F \subset \mathcal F_m$, where $ \mathcal I_m \mathcal F$ stands for the space of vector fields on the form $\sum_{i=1}^s f_i X_i $ with $ f_1, \dots, f_s \in \mathcal I_m$ and $X_1, \dots, X_s \in \mathcal F$.
Moreover,  $ \mathcal I_m \mathcal F$ is a Lie ideal of $\mathcal F_m$, since for all $ X \in \mathcal F, Y \in \mathcal F_m $ and $F \in \mathcal I_m $:
 $$ [FX, Y] = \underbrace{F}_{\in \mathcal I_m} \underbrace{[X,Y]}_{\in \mathcal F} - \underbrace{Y[F]}_{\in \mathcal I_m} \, \underbrace{X}_{\in \mathcal F} . $$
 
\vspace{.5cm}

\begin{propositions}{The isotropy Lie algebra: original construction}{holonomyLieAlgebra}
Let $(M,\mathcal F)$ be a singular foliation. For every $ m \in M$, the isotropy Lie algebra at $m$ is canonically isomorphic to the quotient Lie algebra:
$$
\xymatrix{\mathfrak g_m(\mathcal F)\ar[r]^{\simeq} & \frac{\mathcal F_m}{\mathcal I_m\mathcal F}.}
$$
\end{propositions}
\begin{proof}
Let $ (A,\rho)$ be an anchored bundle over $\mathcal F $, defined in a neighborhood of $m$. For any $a \in \ker(\rho|_m)$, consider $\tilde{a} $ is any section of $A$ through $a$. By construction, $\rho(\tilde{a})$ is in $\mathcal F_m$. Since any two such sections  $ \widetilde{a}_1, \widetilde{a}_2 $ through $a $ satisfy
  $$\widetilde{a}_1- \widetilde{a}_2= \sum_{i=1}^{\mathrm{rk(A)}} F_i b_i  $$
  with $ F_i \in \mathcal I_m $ and $b_i\in\Gamma(A)$ for $i=1, \ldots, \mathrm{rk}(A)$,
  we have
   $$\rho \left(\widetilde{a}_1\right) -\rho \left(  \widetilde{a}_2 \right)=\sum_{i=1}^{\mathrm{rk}(A)} F_i\, \rho(b_i).$$
   As a consequence,
   $\rho \left(\widetilde{a}_1\right) -\rho \left(  \widetilde{a}_2 \right) \in \mathcal I_m \mathcal F $, and the map
    $$\xymatrix{\ker(\rho_m) \ar[r]&\frac{\mathcal F_m}{\mathcal I_m\mathcal F}}$$
    is therefore well-defined. It is also surjective by construction. 
   It has the strong kernel in its kernel, and therefore goes down to a morphism:
$$ \xymatrix{\mathfrak g_m(\mathcal F) = \frac{\ker(\rho_m)}{\mathrm{Sker}(\rho,m) } \ar[r]& \frac{\mathcal F_m}{\mathcal I_m\mathcal F} }.$$
The anchor condition implies that it is a Lie algebra morphism from
 $\mathfrak g_m(\mathcal F)$ onto $ \frac{\mathcal F_m}{\mathcal I_m\mathcal F}$.
Let us show that the kernel of this map is zero. An element $a \in {\mathrm{Ker}}(\rho_m)$ has a class in $ {\mathfrak g}_m(\mathcal F)$ mapped to zero by the previous Lie algebra morphism if  there exists a section $\tilde  a \in \Gamma(A)$ through $a$ be such that $  \rho(\tilde a) \in \mathcal I_m \mathcal F $. This implies that $ \rho( \tilde{a}) = \sum_{i=1}^k f_i \rho(\tilde{a}_i) $ for some sections $ \tilde{a}_1, \dots, \tilde{a}_k$ and some functions $ f_1, \dots, f_k$ vanishing at $ m$.
In particular, the section $ \tilde{a}- \sum_{i=1}^k f_i \tilde{a}_i$ belongs to the strong kernel of $\rho$ at $m$. Since its value at $m$ is $a$, this completes the proof.
\end{proof}

\vspace{.5cm}

\begin{exo}
   Let $(M, \mathcal{F})$ be a singular foliation. For $m\in M$, consider the evaluation map $\mathrm{ev}_m\colon \mathcal{F}\rightarrow T_mM$, $X\mapsto X(m)$.  Show that
\begin{enumerate}
       \item the image of $\mathrm{ev}_m$ is $T_m\mathcal{F}$ and the kernel is $\mathcal{F}_m$. Also, $\mathrm{ev}_m$ goes to quotient to a map  $\frac{\mathcal{F}}{\mathcal{I}_m\mathcal{F}}\rightarrow T_mM$ that we denote by $\underline{\mathrm{ev}}_m$.
       \item $\mathfrak g_m(\mathcal{F})\simeq \ker(\underline{\mathrm{ev}}_m)=\frac{\mathcal{F}_m}{\mathcal{I}_m\mathcal{F}}$.
       \item $\mathrm{rk}_{m}(\mathcal F)=\dim (\frac{\mathcal{F}}{\mathcal{I}_m\mathcal{F}})$.  
   \end{enumerate}
\end{exo}

\begin{exo}
\label{exo:action}
Let $m$ be a point where all vector fields in a singular foliation $ \mathcal F$ vanish. 
Show that the tangent bundle of $ M$ at $ m$ is the dual of $\mathcal I_m/\mathcal I_m^2 $, where $ \mathcal I_m$ is as above. Then show that the Lie algebra $ \mathcal F$ acts on the vector space $\mathcal I_m/\mathcal I_m^2 $, and that this action is zero for any element in $\mathcal I_m \mathcal F $.
Conclude that $ \mathfrak g_m(\mathcal F)$ acts on $ T_m^* M$ and therefore on $ T_m M$.
\end{exo}

\subsection{The linear isotropy Lie algebra}

\subsubsection{The linear part of a vector field vanishing at a point}

Let $ m $ be a point in a manifold $M$, and $ \mathcal I_m $ be the ideal of functions vanishing at $m \in M$. Denote by $  \mathfrak X_m(M)$ the Lie algebra of vector fields vanishing at $m$.
The purpose of this preliminary section is to show that there exists a natural Lie algebra morphism:
 $$ \mathfrak X_m (M) \longrightarrow {\mathrm{gl}}(T_m M) .$$
There are several equivalent manners to see this Lie algebra homomorphism, that we now detail.
\begin{enumerate}
	\item One manner is simply to take local coordinates $(x_1, \dots, x_n) $ in which $m$ has coordinates $ (0, \dots, 0)$. The vector fields  
	$$ \frac{\partial}{\partial x_1}, \dots,\frac{\partial}{\partial x_n} ,$$
	restricted to $ T_mM$ form a basis  of that vector space that we shall denote by $ \delta_1, \dots, \delta_n $. We then map a vector field: 
	 $$  \sum_{i=1}^n X_i{\scriptstyle {(x_1, \dots, x_n)}}  \frac{\partial}{\partial x_i} $$ to
	 the linear endomorphism of $T_m M$ whose matrix in the basis  $ \delta_1, \dots, \delta_n $ is 
	  $$   \begin{pmatrix}
	  & &  \\
	  &  \frac{\partial X_i}{\partial x_j} {\scriptstyle {(0, \dots, 0)}}  &  \\ & &  \end{pmatrix}$$   
	  \item[] We leave it to the reader to check that this is indeed a Lie algebra morphism.
	  \item[] Although very explicit, this method has a drawback: we have to check it does not depend on the choice of local coordinates. It is therefore better to use the coming two descriptions then show that, in local coordinates, they take the previous form.
	\item The second manner is to use the flow $\phi_t^X $ of a vector field $X \in \mathfrak X_m (M)$. Since $X$ vanishes at $m $, for every $ \eta>0$, there is a neighborhood $ \mathcal U_m$ of $m$ on which $ \phi_t^X$ is well-defined for all $ t \in -]\eta, \eta[$. Also, $\phi_t^X(m)=m $, so that the differential of  $\phi_t^X$ at $m$ is a family depending on $t \in -]\eta, \eta[$ of invertible  linear endomorphisms
	 $$ T_m \phi_t^X \colon T_m M \longleftrightarrow T_m M $$
	We then define a linear endomorphism of $ T_mM$ by 
	 $$  X \mapsto \left.\frac{\partial}{\partial t} \right|_{t=0}  T_m \phi_t^X$$
	 \item[] The previous map is well-defined, but it is not clear that it is a Lie algebra morphism. Also, defining it required the notion of flow, which does not make sense in algebraic geometry.
	 \item The third manner is to look, for any vector field $X$ vanishing at $m$, at the adjoint action:
	  $$ Y \mapsto [X,Y]  $$
	  and to check that $ [X,Y]_{ |_{m}} $ only depends on $Y_{|_m} $, so that the adjoint action induces a linear endomorphism of $T_m M $. The Jacobi identity implies that this map is a Lie algebra morphism. 
	 \item A fourth manner is to use the canonical identification $$T_m ^* M \simeq \frac{\mathcal I_m}{\mathcal I_m^2} $$
	with $ \mathcal I_m$
 the ideal\footnote{
 In real analytic or complex geometry setting, ``ideal'' must be understood as ``sheaf of ideals''
 } 
  of functions vanishing on $M$ (which, in the algebraic geometry setting, is in fact a definition of the cotangent space). 
Consider vector fields as derivations of the sheaf of functions: a vector field $X$ vanishes at $m$ if and only if $X[\mathcal I_m] \subset \mathcal I_m $. By derivation properties, this implies
$X[\mathcal I_m^2] \subset  \mathcal I_m^2 $, so that $ X$ induces a linear endomorphism of
 of $T_m ^* M \simeq \frac{\mathcal I_m}{\mathcal I_m^2}$. Since the bracket of vector fields is their commutator, when seen as a derivation, it follows that the map above is a Lie algebra morphism: 
  $$ \mathfrak X_m (M) \longrightarrow {\mathrm{gl}}(T_m^* M) .$$
  The desired Lie algebra morphism is obtained by composing the latter morphism with the canonical dualization Lie algebra isomorphism $ {\mathrm{gl}}(T_m^* M) \simeq  {\mathrm{gl}}(T_m M)$. 
  \item Last, one can recognize that for the singular foliation $ \mathfrak X_m(M)$ of all vector fields vanishing at $m$, the isotropy Lie algebra at $m$ is $\mathrm{gl}(T_mM)$. The map is then just the Lie algebra morphism 
   $  \mathfrak X_m(M) \longrightarrow   \frac{\mathfrak X_m(M)}{\mathcal I_m \mathfrak X_m(M)}= \mathrm{gl}(T_mM)$.

 \end{enumerate}

 \vspace{0.5cm}

\begin{notations}{Linear Part of a vector field}
{}
We denote by ${\mathfrak{Lin}} $ the Lie algebra morphism
 $$\{ \hbox{Vector fields vanishing at $m$}\} \longrightarrow {\mathrm{gl}}(T_m M). $$
described in the lines above.

\end{notations}

\subsubsection{The linear isotropy Lie algebra of a foliation vanishing at a point}

Let us consider a foliation $\mathcal F$ on a manifold $M$ made of vector fields vanishing at a point $m$ (equivalently, such that $ \{m\}$ is a leaf).
To our knowledge, Dominique Cerveau \cite{Cerveau} is the first to have understood the importance and studied the following Lie algebra.

\vspace{0.5cm}

\begin{definitions}{Linear isotropy Lie algebra}{linearpart}
Let  $\mathcal F $ be a singular foliation, and $m$ a point.
We call \emph{linear isotropy Lie algebra of $\mathcal F $} and denote by $\mathfrak g_m^{lin}(\mathcal F) $ the image of $\mathcal F $ through the linear part morphism
 $ {\mathfrak{Lin}} $.  In equation:
 $$\mathfrak g_m^{lin}(\mathcal F)  := {\mathfrak{Lin}} (\mathcal F).  $$
\end{definitions}

\vspace{0.5cm}

\begin{remark}
Upon choosing local coordinates, and therefore a basis of $T_m M $, the linear isotropy Lie algebra of $ \mathcal F$ at $m$ at the origin is the sub-Lie-algebra of all matrices $(a_{ij})  $ such that there exists $ X \in \mathcal F$ whose Taylor expansion at the origin reads:
$$
X = \sum_{i,j}a_{ij}x_i\frac{\partial}{\partial x_j} + \mathrm{higher~order~terms}
$$
\end{remark}

\begin{example}Let $\mathcal F$ be the singular foliation induced by a Lie algebra action of $\mathfrak g\subset \mathfrak{gl}(\mathbb R^d)$ on $\mathbb R^d $. Then the linear holonomy of $\mathcal F$  at $0$ is $\mathfrak g/\mathfrak k$ with $\mathfrak k \subset \mathfrak g$ being the Lie algebra of all elements that act trivially.
\end{example}

\begin{example}
Let $\mathcal F\subset\mathcal  I_0^2\mathfrak X(\mathbb R^d)$, i.e., a foliation made of vector fields vanishing quadratically in the origin (See section \ref{ex:singFolVanish}). Then the linear isotropy Lie algebra of $\mathcal F$ at $0$ is $\{0\}$. 
\end{example}

By construction, 
\begin{equation}\label{eq:sujectiveMorphism} { {\mathfrak{Lin}} } \colon \mathcal F \longrightarrow \mathfrak g_m^{lin}(\mathcal F) \end{equation}
is a surjective morphism of Lie algebras. 
Here is therefore a natural question, that one could ask for any Lie algebra morphism: 
Does \eqref{eq:sujectiveMorphism} admit a section which is a Lie algebra morphism? If yes, it means, geometrically, that $\mathcal F $ contains a sub-singular foliation associated to the Lie algebra action of  $\mathrm g_m^{lin} (\mathcal F) $ on~$M$.

\begin{question}
Does ${\mathfrak{Lin}} $ admit sections? I.e, does $\mathcal F $ contain, in a neighborhood $\mathcal U $ of $m$, a sub-singular foliation given by a Lie algebra action of ${\mathfrak g}_m^{lin} (\mathcal F) $ on $\mathcal U $?
\end{question}

In general, the answer to this kind of question tends to be ``no, unless the image is semi-simple''. And in the infinite dimensional case, the answer tends to be ``no,  unless the image is compact and semi-simple. If the image is semi-simple, then there are only formal sections''. There are several results in that vein, by Conn for Poisson structures and Zung for Lie algebroids. To our knowledge, the singular foliation case is widely open: we will discuss this in Section \ref{sec:linearization}.

Here is an important result by Dominique Cerveau for the semi-simple case. A more recent proof can also be found in \cite{Ryvkin2}, Theorem 2.8.

\vspace{.5cm}

\begin{theorems}{A linearization theorem by Dominique Cerveau}{LinearCerveau}

If the linear isotropy Lie algebra of $\mathcal F$ at $m$ is a semi-simple Lie algebra, then the map:
 $$  {\mathfrak{Lin}} \colon \mathcal F \longrightarrow \mathfrak g_m^{lin} (\mathcal F) $$
 admits a formal section\footnote{See discussion around Equation \eqref{eq:formalsection}.}  which is a  Lie algebra morphism.
\end{theorems}

\vspace{.5cm}

We will not prove this theorem, but we will at least us say a word about its meaning\footnote{Also, a more general result by Dominique Cerveau will be stated later, see Proposition \ref{prop:cerveau}.}. Formal functions\footnote{In the smooth setting, it is the quotient of $ \mathcal C^\infty (M)$ by the ideal of functions vanishing with all their derivatives. In the other settings, it is a formal completion, i.e., the ring of formal power series in $d$ variables near $m$.} at a point $m \in M$  form an algebra that we should denote by $\hat{\mathcal O}_m $. 
Formal functions $\hat{\mathcal O}_m $ are a module over
over germs of smooth, complex, polynomials, or real analytic functions (that we should denote by $\mathcal O$). 
As a consequence, the tensor product 
$$\hat{\mathcal O}_m \otimes_{\mathcal O} \mathcal F$$
is a finitely generated $\hat{\mathcal O}_m $ module stable under Lie bracket\footnote{It is  an algebraic singular foliation in the sense of Definition \ref{def:algebraic} for the ring $\hat{\mathcal O}_m $.}, and ${\mathfrak{Lin}}$ extends by linearity to a Lie algebra morphism: 
\begin{equation}\label{eq:formalsection} {\mathfrak{Lin}} \colon \hat{\mathcal O}_m \otimes_{\mathcal O} \mathcal F \longrightarrow \mathfrak g_m^{lin} (\mathcal F) \end{equation}
The result of Dominique Cerveau states that this Lie algebra morphism admits a section which is a Lie algebra morphism.

\begin{exo}
\label{exo:rankof restriction}
Let $(M,\mathcal F)$ be a singular foliation. Assume all vector fields in $ \mathcal F$ are zero at  a given point  $m\in M$. Show that the quotient space\footnote{where $\mathcal I_m\subset \mathcal C^\infty(M)$ is the ideal of all functions vanishing in $m$}
$$
\frac{\mathcal F }{\mathcal F\cap \mathcal I_m^2\mathfrak X(M)},
$$
is a Lie algebra isomorphic to $\mathfrak g^{lin}_m(\mathcal F)$.
\end{exo}

The previous exercise extends as follows. We leave the proof to the reader.

\begin{proposition}
Let $(M,\mathcal F)$ be a foliated manifold and $m\in M$. The linear isotropy Lie algebra of $\mathcal F$ at $m$ is canonically isomorphic to the Lie algebra
$$
\mathfrak g^{lin}_m(\mathcal F)=\frac{\mathcal F\cap \mathcal I_m\mathfrak X(M)}{\mathcal F\cap \mathcal  I_m^2\mathfrak X(M)},
$$
where $\mathcal I_m\subset \mathcal C^\infty(M)$ is the ideal of all functions vanishing in $m$, equipped with the Lie bracket induced by the Lie bracket of vector fields. 
\end{proposition}

The linear isotropy Lie algebra captures the ``linear approximation'' of the foliation at a given point. As a consequence, for foliations vanishing quadratically at a point, this Lie algebra is trivial.

\begin{example}
Let $\mathcal F=\mathcal I_0^n\mathbb{R}^d$. Then ${\mathfrak g}_0(\mathcal F)$ will have dimension $$d\times {{n+d-1}\choose n},$$ while the linear one will be trivial for $n\geq 2$.  
\end{example}

When $\mathcal F$ admits real analytic generators, one can prove that the linear holonomy contains all the semi-simplicity of $\mathfrak g_m(\mathcal F)$, i.e.:

\vspace{0.5cm}

\begin{propositions}{The semi-simple part is linear}{prop:ArtinRees} 
\label{prop:kernel}
For a real analytic foliation $\mathcal F$, the kernel of the linearization map $\mathfrak g_m(\mathcal F)\to \mathfrak g^{lin}_m(\mathcal F)$
is a nilpotent Lie algebra.
\end{propositions}

\begin{proof}
We refer to Theorem 1.10 in \cite{Ryvkin2} for the complete proof. The main ingredient of the proof is the Artin-Rees Lemma, which is valid for Noetherian rings (i.e., the ring of germs of analytic functions but not the ring of smooth functions). 
\end{proof}

Whether the theorem holds also in the smooth category is still an open problem:

\begin{question}
\label{ques:smoothcase}
Is it possible to omit the assumption ``locally real analytic'' in Proposition \ref{prop:kernel}? 
\end{question}

\vspace{0.5cm}
Let us finish with a Levi-Malcev type of theorem, in the smooth setting.

Consider the following series of Lie algebras\footnote{$\mathcal F(m)$ stands here for vector fields in $\mathcal F$ that vanish at $m$, and are defined in a neighborhood of $m$.}

\begin{align*}
\mathcal F(m) \to {\mathfrak j}^\infty_m(\mathcal F)\to ...\to {\mathfrak j}^N_m(\mathcal F)\to {\mathfrak j}^{N-1}_m(\mathcal F) \to ...\to  {\mathfrak j}^{1}_m(\mathcal F)=\mathfrak g^{lin}_m(\mathcal F),
\end{align*} Here ${\mathfrak j}^N_m(\mathcal F)=\frac{\mathcal F(m)}{I_m^{N+1} \mathfrak{X}(M)}$ are $N$-jets of vector fields on $\mathcal F(m)$ and ${\mathfrak j}^\infty_m(\mathcal F)$ their projective limit, i.e., the space of Taylor expansions of elements in $\mathcal F(m)$.
For  $N \in \mathbb N$, the kernel of ${\mathfrak j}^N_m(\mathcal F)\to \mathfrak g^{lin}_m(\mathcal F)$ is a nilpotent Lie ideal. The linear isotropy Lie algebra $\mathfrak g^{lin}_m(\mathcal F)$ might still contain a solvable ideal. By dividing out the maximal solvable ideal $\mathfrak r$, we obtain a semisimple Lie algebra  $\mathfrak  g^{lin}_m(\mathcal F)^{ss}=\frac{g^{lin}_m(\mathcal F)}{\mathfrak r}$, which can be added on the right in the above filtration to obtain a series of surjections:
\begin{align*}
\mathcal F\to {\mathfrak j}^\infty_m(\mathcal F)\to ...\to {\mathfrak j}^N_m(\mathcal F)\to {\mathfrak j}^{N-1}_m(\mathcal F) \to ...\to  \mathfrak g^{lin}_m(\mathcal F)\to \mathfrak  g^{lin}_m(\mathcal F)^{ss}. 
\end{align*}
Using classical techniques, and using nilpotency of the kernels as above, one can prove that singular  foliations satisfies a sort of Levi-Malcev-theorem:

\begin{proposition}[Cerveau, 1977]
\label{prop:cerveau} \cite{Cerveau}-\cite{Ryvkin2}
There is a formal Lie algebra section of ${\mathfrak j}^\infty_m(\mathcal F)\to  g^{lin}_m(\mathcal F)^{ss}. $, i.e., a Lie algebra homomorphism from the semi-simple part of the linear isotropy Lie algebra at $m $ to the Lie algebra ${\mathfrak j}^\infty_m(\mathcal F) $.  
\end{proposition}

\subsection{The holonomy Lie algebroid of a leaf}\label{sec:holonomy-Lie-algebroid}

So far, we have attached several Lie algebras to a point $m$ of a foliated manifold $ (M,\mathcal F)$. 
There is a more general object that captures the dynamics along a leaf, namely the holonomy Lie algebroid. 

As for the isotropy Lie algebra, the simplest way to define it  is purely algebraic, it consists in  seeing it  as a quotient. But this is only valid upon the assumption that $L$ is embedded. We will start by assuming that it is embedded, then extend the construction to the general case using almost Lie algebroid. Also, we leave it to the reader to see how to adapt the result of this section to the complex or real analytic settings.

\vspace{0.5cm}

Let $L$ be a leaf of a singular foliation $\mathcal F $ on a manifold $M$. To start with, let us assume that $ L$ is an embedded submanifold.  Consider the quotient space 
  \begin{equation}
  \label{eq:quotientalgebroid}
       \frac{\mathcal F}{\mathcal I_L \mathcal F}  
  \end{equation}
Since $ L$ is an embedded leaf, $ \mathcal I_L$ is what we used to call a foliated ideal, i.e., 
  $$ \mathcal F [\mathcal I_L]  \subset \mathcal I_L $$ 
  so that the quotient described in Equation \eqref{eq:quotientalgebroid} inherits a Lie algebra bracket. There is also a $ C^\infty(M)$-module structure on the quotient  \eqref{eq:quotientalgebroid}, but since $ \mathcal I_L$ acts by zero, it is in fact a $ C^\infty(M) / \mathcal I_L \simeq \mathcal C^\infty(L)$-module. 
Altogether, these two structures turn the quotient \eqref{eq:quotientalgebroid}  into a Lie-Rinehart algebra over $ \mathcal C^\infty(L)$.   If one can show that it is in fact a \emph{projective} $ \mathcal C^\infty(L)$-module, then by the Serre-Swan theorem, there exists a vector bundle $ A_L \to L$ such that
 \begin{equation}\label{SerreSwann}  \Gamma_c(A_L) \, \simeq \,  \frac{\mathcal F}{\mathcal I_L\mathcal F}  ,\end{equation}
as $ \mathcal C^\infty(L)$-module. The Lie bracket then equips $ A_L$ with a Lie algebroid structure, whose anchor is given for all $ \tilde{a} \in \Gamma(A_L)$  by 
 $$ \rho(\tilde{a}) = \hat{a}_{|_L} $$
where $\hat{a} \in \mathcal F$ is any element whose class modulo $ \mathcal I_L \mathcal F $ corresponds to $ \tilde{a}$ in Equation \eqref{SerreSwann}. 

But how do we prove that $ \mathcal F/\mathcal I_L\mathcal F$ is a projective $ \mathcal C^\infty(L)$-module? And can we deal with the non-embedded case? To address both issues, we will, again, use almost Lie algebroids.
From now on, we make no assumption on $L$.
Let $ (A,\rho)$ be an anchored bundle defined in a neighborhood $ \mathcal U$ of a point  $\ell \in L$.
We denote by $ (\mathcal U \cap L)_\ell $ the connected component of $ \ell$ in $\mathcal U \cap L$. Without any loss of generality, we can assume this submanifold is now embedded.
Theorem \ref{prop:anchoredBundleLeaf} implies that the Strong kernel of $ \rho$ has the same dimension at every point in $ L \cap \mathcal U$. Since  any element of the strong kernel admits a section through it which is valued in the Strong kernel at every other point, this guarantees that, altogether, the strong kernels of points in  to $  (L \cap \mathcal U)_\ell$ assemble to a sub-vector bundle of $ A_{|_L}$.
We denote by $ A_L$ the quotient. 
We leave it to the reader to check that the isomorphism \eqref{SerreSwann} holds true in the following form: every point $ \ell \in L$ has a neighborhood $ \mathcal V$ in $M$ such that
\begin{equation}
    \label{serreswannlocal} \Gamma_c(A_L|_ {\mathcal U}) \, \simeq \,  \frac{\mathcal F|_\mathcal U}{\mathcal I_{(L \cap \mathcal U)_\ell}\mathcal F|_\mathcal U} \end{equation}
where $\mathcal F|_\mathcal U$ is the restriction of the singular foliation to $\mathcal U $, and $\mathcal I_{(L \cap \mathcal U)_\ell}$ the ideal of functions vanishing on $(L \cap \mathcal U)_\ell$. Since the right-hand side of the previous equation does not make any reference to the anchored bundle $(A,\rho) $, it means that $A_L$ does not depend on the choice of a particular anchored bundle. As a consequence, $ A_L \to L$, together with its anchor, are defined on the whole leaf $L$. Moreover, the Lie bracket that the right-hand side of the previous equation equips $ A_L$ with a Lie algebroid structure, whose anchor map $ A_L \to TL$ is surjective at every point.

\vspace{0.5cm}

\begin{definitions}{Holonomy Lie algebroid of a leaf}{holonomy}
Let $ L$ be a leaf of a foliated manifold $ (M,\mathcal F)$. We call
\emph{holonomy Lie algebroid of $\mathcal F$ along $L$} the unique Lie algebroid $(A_L,\rho,[\cdot\,, \cdot]) $ such that every point $ \ell \in L$ admits a neighborhood $ \mathcal U$ on which there is an isomorphism of both Lie algebra and $ \mathcal C^\infty(L)$-module structures as in Equation \eqref{serreswannlocal}.
\end{definitions}
\vspace{0.5cm}

Also, there is a short exact sequence:
 $$  \xymatrix{  \mathfrak g_L(\mathcal F)\ar@{^(->}[r]\ar[d] & A_L \ar^{\rho}@{->>}[r] \ar[d]& TL\ar[d]\\ L \ar^{=}[r]& L\ar^{=}[r] &L}  ,$$
where $\mathfrak{g}_L(\mathcal F)=\bigsqcup_{\ell\in L}\mathfrak{g}_\ell(\mathcal F)$ is a bundle of Lie algebras over $L$ whose fiber at $ \ell \in L$ is the isotropy Lie algebra $\mathfrak g_\ell (\mathcal F) $.

\begin{exo}
Show that the Lie algebroid $A_L$ acts on the normal bundle $ TM_{|L}/TL$ of $L$ in $ M$.
This action is exploited by Androulidakis and Zambon in \cite{AZ2}. Show that if a leaf is a point,  we recover the action described in Exercise \ref{exo:action}\footnote{\emph{Hint:}
If the leaf $L$ is embedded, Exercise \ref{exo:action} can be imitated as follows: Show that the normal bundle of $ M$ at $ m$ is the dual bundle of  a bundle whose sections are $\mathcal I_L/\mathcal I_L^2 $, where $ \mathcal I_L$ functions vanishing on $L$. Then show that the Lie algebra $ \mathcal F$ acts on the vector space $\mathcal I_L/\mathcal I_L^2 $, and that this action is zero for any element in $\mathcal I_L \mathcal F $. Conclude that $A_L$ acts on dual of the normal bundle and therefore on the normal bundle. Use sheaves for the general case.}.
\end{exo}

\begin{exo}
Find an example where there is no section $ TL \to A_L$ of the anchor map which is a Lie algebra morphism. Show that if the leaf is flat\footnote{See Exercise \ref{exo:flatleaves} and the discussion following it.}, then such a section exists.
\end{exo}

\section{Bisubmersions over a singular foliation}\label{sec:Bisubmersion}

 \subsection{Definitions} 

\noindent
In order to understand singular foliations, the most crucial and intriguing object is certainly the so-called holonomy groupoid. This groupoid is constructed from an object called bisubmersion. Both concepts were introduced by Androulidakis and Skandalis in \cite{AS}.

\vspace{1mm}
\noindent
A word of caution. The notion of bisubmersion seems to be extremely basic. At first look, it seems to be only a ``Lie-groupoid-like but without a product'' object, that is, it seems to be a very poor structure. However, as we will see in the coming lines, although its definition is extremely short, it is indeed a very subtle and rich notion. In particular, it is easy to make wrong statements about them, and many true statements are hard to prove.

\vspace{1mm}
\noindent
There is an analogy, that has its limits, but may sound familiar to some readers: bisubmersions (more precisely those called atlases) are to the holonomy groupoid what Lie groupoids representing a stack are to the  differential stack in question.   Androulidakis and Zambon also pointed us out that bisubmersions can also be thought of as the plots for the holonomy groupoid, and that equivalence is like equivalence of plots in diffeology - an enriching perspective.

\vspace{2mm}
\noindent
Let us give the definition. From now on, we work in the \underline{smooth} setting, and leave it to the reader to adapt to real analytic or complex ones.

\vspace{.5cm}

\begin{definitions}{Androulidakis-Skandalis' bisubmersions}{bisub}
Let $M$ be a manifold equipped with a singular foliation $\mathcal F $.
A bisubmersion over  $(M,\mathcal F) $ is a triple $(W,s,t) $ where:
\begin{enumerate}
    \item $W$ is a manifold,
    \item $s,t\colon W \to M$ are  submersions\footnote{We do not assume $s$ and $t$ to be surjective. See Section \ref{sec:pull-back} for a definition of the pull-back singular foliation.}, respectively called \emph{source} and \emph{target},
\end{enumerate}
such that 
\begin{enumerate}
    \item the pull-back singular foliations $ s^{-1}(\mathcal F)$ and $t^{-1}( \mathcal F)$ are equal\footnote{Notice that this implies that $ \Gamma(\ker (Ts))+\Gamma(\ker (Tt)) \subset s^{-1} (\mathcal F) = t^{-1} (\mathcal F)$},
    \item and any vector field $Z \in s^{-1} (\mathcal F) = t^{-1} (\mathcal F )$ is of the form\footnote{In general, this decomposition is not unique at all. This second condition can be equivalently stated as $ \Gamma(\ker (Ts))+\Gamma(\ker (Tt)) = s^{-1}(\mathcal F) = t^{-1} (\mathcal F)  $.} $Z=X+Y $ with $X\in \Gamma(\ker (Ts)) $ and $Y \in \Gamma(\ker (Tt)) $.
\end{enumerate}
\vspace{1mm}
\noindent
Also, we will use the name \emph{bi-pull-back singular foliation} and the notation $ {}^{s}\mathcal F{\, }^t $ for the singular foliation on $W$ given by:
\begin{equation} \label{eq:singFolonW}  {}^{s}\mathcal F{\,}^t := s^{-1}(\mathcal F)=t^{-1}(\mathcal F)= \Gamma({\mathrm{ker}}(Ts))+ \Gamma({\mathrm{ker}}(Tt)). \end{equation}

\vspace{1mm}
\noindent
Last, vector fields in $\Gamma({\mathrm{ker}}(Ts))\cap \Gamma({\mathrm{ker}}(Tt))$ shall be said to be is \emph{bi-vertical}.
\end{definitions}

\vspace{.5cm}
\noindent
When we will need to insist on the names of all structures, such a bisubmersion shall be denoted by $$ M \stackrel{s}{\leftarrow} W \stackrel{t}{\rightarrow} M .$$ Quite often, we will simply say ``a bisubmersion $W$ over a foliated manifold $(M,\mathcal F) $''.

\begin{bclogo}[arrondi = 0.1, logo = \bcdz]{Warning !}
Notice that it does not make much sense to say  "Let
$ M \stackrel{s}{\leftarrow} W \stackrel{t}{\rightarrow} M $ be a bisubmersion",  without mentioning over which singular foliation $ \mathcal F$ it is a bisubmersion. We will always say "Let
$ M \stackrel{s}{\leftarrow} W \stackrel{t}{\rightarrow} M $ be a bisubmersion over $ \mathcal F$ " or "over $ (M,\mathcal F)$".
\end{bclogo}

\begin{exo}
Let $ W$ be a bisubmersion over $ (M,\mathcal F)$.
Show that if ${\mathrm{dim}}(W)= {\mathrm{dim}}(M) $, then $\mathcal F=0 $.
\end{exo}

\begin{exo}
Show that, for any bisubmersion, compactly supported\footnote{Again, one could avoid the use of compactly supported vector fields by using sheaves.} bivertical vector fields form  a $\mathcal C^\infty(W) $-module which is closed under Lie bracket. 
Notice that it may not be a singular foliation, for it may not be finitely generated, see Exercise \ref{exo:bivertical_nonSF}.
\end{exo}

\begin{exo}
\label{exo:bivertical_nonSF}
Let $X$ be a compactly supported vector field on a manifold $M$, and $ \mathcal F_X = \{f X | f \in \mathcal C^\infty(M) \}$ the singular foliation it generates as in Exercise \ref{exo:rank1}. Show that   $ W:=\mathbb R  \times M$ equipped with the following source and target maps
 $$  s (u,m) := m \hbox{ and } t(u,m) := \phi_u^X(m) \hbox{ for all $ m \in M,u\in \mathbb R$}$$
 is a bisubmersion for $ \mathcal F_X$. Show that if $ M = \mathbb R$ and $ X$ is a vector field with support $ [-1,1]$, bivertical vector fields do not form a locally finitely generated module. 
\end{exo}

\begin{exo}
Let $M$ be  a manifold, and $ \mathcal F=\mathfrak X_c (M)$ be the singular foliation of all compacted supported vector fields. Show that $  W:= M \times M$, equipped with source and targets defined to be the projections on the first and second components, respectively, is a bisubmersion over  $\mathcal F$
with no non-trivial bi-vertical bivector fields.
\end{exo}

\noindent
Here is a technical lemma: the reader may jump to Definition \ref{def:unitsbisub} directly. The importance of this lemma will appear later on, and may not be obvious at first sight. 

\begin{lemma} 
\label{lem:localstructure}
Consider a bisubmersion $ M \stackrel{s}{\leftarrow} W \stackrel{t}{\rightarrow} M $ over a singular foliation $ \mathcal F$, and let  ${}^{s}\mathcal F{ }^t $ be the bi-pull-back singular foliation as in Equation \eqref{eq:singFolonW}.
Choose $x \in W$. Let $ \underline{X}_1, \dots, \underline{X}_{r_s}$ be local generators of $ \mathcal F$ in a neighborhood of $ s(x)$
and $ \underline{Y}_1, \dots, \underline{Y}_{r_t}$ be local generators of $ \mathcal F$ in a neighborhood of $ t(x)$. Let $k=2\, \mathrm{rk}(\ker Ts)=2\, \mathrm{rk}(\ker Tt)$.
Then $x$ admits a neighborhood $\mathcal U $ on which $ {}^{s}\mathcal F{\,}^t $ is generated by vector fields $ X_1, \dots, X_{r_s},Y_1, \dots,Y_{r_t},Z_1, \dots,Z_{k}$ satisfying the following properties:
\begin{enumerate}
\item  $ X_1, \dots, X_{r_s} $ belong to $ \Gamma({\mathrm{ker}}(Tt))$ and are $s$-related to $ \underline{X}_1, \dots, \underline{X}_{r_s}$, 
\item $ Y_1, \dots, Y_{r_t} $ belong to $ \Gamma({\mathrm{ker}}(Ts))$ and are $t$-related to $ \underline{Y}_1, \dots, \underline{Y}_{r_t}$,
\item  $Z_1, \dots,Z_{k}\in\Gamma({\mathrm{ker}}(Ts)\cap \Gamma({\mathrm{ker}}(Tt))   $ are bi-vertical vector fields, 
\item the vector fields $ X_{1}, \dots, X_{r_s}, Z_1, \dots, Z_{k}$ generate ${\mathrm{ker}}(Tt) $ at every point of $ \mathcal U$, and the vector fields $  Y_{1}, \dots, Y_{r_t}, Z_1, \dots, Z_{k}$ generate ${\mathrm{ker}}(Ts) $ at every point of $ \mathcal U$.
 \end{enumerate}
\end{lemma}
\begin{proof}
By definition of the pull-back singular foliation $ s^{-1}(\mathcal F)$, there exists for every $ i \in \{1, \dots, r_s\}$ a vector field $ X_i^! \in s^{-1}(\mathcal F)$ which is $s$-related to $\underline X_i$.
By definition of a bi-submersion, there exists $ X_i^s\in\Gamma({\mathrm{ker}}(Ts))$ and $X_i^t\in\Gamma({\mathrm{ker}}(Tt))$ such that
 $X_i^! = X_i^s + X_i^t$. By construction, $ X_i:=X^t_i\in \Gamma({\mathrm{ker}}(Tt))$ is $s$-related to $\underline{X}_i$ and $t$-related to $0$. 
 \emph{Item 1 is therefore satisfied}.

 \vspace{1mm}
 \noindent
 The same argument, with $s$ replaced by $t$, yields vector fields $Y_1, \dots, Y_{r_t} \in \Gamma({\mathrm{ker}}(Ts))$ that are $t$-related to $\underline{Y}_1, \dots, \underline{Y}_{r_t}$ and $s$-related to $0$. 
 \emph{Item 2 is therefore satisfied}.

  \vspace{1mm}
 \noindent
 Now, it is an elementary property of $s^{-1}(\mathcal F)$ that for any vector field $ Z \in  s^{-1}(\mathcal F)$, there exists  functions $ (g^i)_{i=1}^{r_s}$, defined in a neighborhood $\mathcal U$ of $x$ in $W$, such that $Z-\sum_{i=1}^{r_s} g^i X_i \in \Gamma({\mathrm{ker}}(Ts))$. In particular, for $Z\in\Gamma({\mathrm{ker}}(Tt))$, it means that 
$$Z-\sum_{i=1}^{r_s} g^i X_i $$
belongs to $\Gamma({\mathrm{ker}}(Ts))  \cap  \Gamma({\mathrm{ker}}(Tt))$, i.e., is bi-vertical.

Consider now a local trivialization $ Z_1', \dots, Z_{k/2}'$ of  ${\mathrm{ker}}(Tt)$ on $ \mathcal U' \subset \mathcal U$. By the above discussion, there exists local functions $$ \left(g_i^j\right)_{i=1,\dots,r_s}^{j=1, \dots, k/2},$$
defined on some neighborhood $ \mathcal U'' \subset \mathcal U'$ of $x$, such that 
$$ Z_j := Z_j'-\sum_{i=1}^{r_s} g^j_i X_i \in  \Gamma({\mathrm{ker}}(Ts)) \cap \Gamma({\mathrm{ker}}(Tt)),\quad j\in \{1,\ldots,k/2\}, $$
 \emph{as claimed in Item 3}.
The family of vector fields defined on $\mathcal U'' $:
$$
\left(  X_1, \dots, X_{r_s},Z_1, \dots, Z_{k/2}  \right)
$$
generates $\Gamma({\mathrm{ker}}(Tt))$, \emph{as required in item 4}.
The same argument, with $s$ replaced by $t$, yields vector fields  $ Z_{k/2+1} \dots, Z_{k} \in  \Gamma({\mathrm{ker}}(Ts))\cap \Gamma({\mathrm{ker}}(Tt))$,
\emph{therefore satisfying Item 3},
such that the family 
$$
\left(Y_1, \dots, Y_{r_t},Z_{k/2+1}, \dots, Z_{k}  \right)
$$
generates $\Gamma({\mathrm{ker}}(Ts))$, \emph{as required in item 4}.

\vspace{1mm}
Altogether, the families $X_1, \dots, X_{r_s}$, $Y_1, \dots, Y_{r_t}$ and, $Z_1, \dots, Z_k$ satisfy all the requirements of Lemma~\ref{lem:localstructure}. 
\end{proof}

\begin{exo}
\label{exo:sprojtproj}
Consider a bisubmersion over $\mathcal F $ as before.
\begin{enumerate}
\item 
Show that $\Gamma({\mathrm{Ker}}(Ts)) $ is generated by vector fields $t$-related to vector fields in $ \mathcal F$. 
\item
Show also that any element of $\mathcal F$ is $t$-related to a section of $\Gamma({\mathrm{Ker}}(Ts))$. 
\end{enumerate} 
\noindent
 \emph{Hint:} Use Lemma \ref{lem:localstructure}.
\end{exo}

\begin{remark}
The vector fields $ Z_1, \dots, Z_k$ in
Lemma   \ref{lem:localstructure} are \emph{not}, in general, generators of bi-vertical vector fields. Otherwise,  bi-vertical vector fields would be finitely generated, which is not true in view of Exercise \ref{exo:bivertical_nonSF}.
\noindent
$\square$
\end{remark}

\vspace{.5cm}
\noindent
We claimed that bisubmersions may be thought of as an equivalent of a Lie groupoid, as a sort of ``Lie-groupoid-without-a-product''. To justify the analogy, as for Lie groupoids, we now define units and bisections.
\vspace{.5cm}

\begin{definitions}{Some important notions: units and bisections}{unitsbisub}

  A \emph{bisection} of a bisubmersion 
$ M \stackrel{s}{\leftarrow} W \stackrel{t}{\rightarrow} M $ over $ \mathcal F$
is a submanifold $\Sigma \subset W$ such that
\begin{enumerate}
\item $ s(\Sigma) \subset M$ and $ t(\Sigma) \subset M$ are open subsets\footnote{We will speak of \emph{global bisections} when $s(\Sigma)=t(\Sigma)=M $, and \emph{local bisections} otherwise.},
\item the restriction of $s$ to $\Sigma$ is a diffeomorphism from $\Sigma $ to $s(\Sigma)$, and the restriction of $t$ to $ \Sigma$ is a diffeomorphism from $\Sigma $ to $t(\Sigma)$.
\end{enumerate}

Let $\mathcal U \subset M $ be an open subset.
A map $\epsilon \colon \mathcal U \longrightarrow W $ is said to be a \emph{unit map}\footnote{Again, we will speak of \emph{global units} when $\mathcal U=M $, and \emph{local units} otherwise.} if it is a section of both $s$ and $t$, i.e., $s\circ \epsilon =t\circ\epsilon=\mathrm{id}_{\mathcal
 U}$. 
\end{definitions}

\vspace{.5cm}
\begin{remark}
The image $\epsilon(\mathcal U)$ of a (local/global) unit map $\epsilon$ of a bisubmersion is a (local/global) bisection. 
Notice that, unlike for groupoids, the unit map may not exist. Moreover, even if it exists, it may not be unique. Examples will be given below.

\noindent
$\square$
\end{remark}

\noindent
The next exercise explains how an anchored bundle can be recovered out of the data of a bisubmersion with unit. An even more general result will be proven in Lemma \ref{lem:BisectionAnchoredBundle}.

\begin{exo}
\label{exo:normalbundleidentity}

Consider a bisubmersion $ M \stackrel{s}{\leftarrow} W \stackrel{t}{\rightarrow} M $ over $ \mathcal F$. Assume a unit map $\epsilon \colon M \to W $ exists.
         \begin{enumerate}
             \item Show that the normal bundle $N = TW/T(\epsilon(M))$ of $\epsilon(M)$ into $W$  is canonically isomorphic to $\ker(Ts) \subset TW|_{ \epsilon(M)} $ and to $\ker(Tt) \subset TW|_{ \epsilon(M)}$
             \item Show that the vector bundle morphism $Ts - Tt: TW|_{\epsilon(M)}\rightarrow TM $ goes to the quotient to give a vector bundle morphism $\rho_W \colon N \to TM $ (over the natural diffeomorphism  $\epsilon(M) \to M$).
             \item[] From now on, we  use the diffeomorphism $ \epsilon \colon M \to \epsilon(M)$ in order to consider $N$ as a vector bundle over $M$ (and not as a vector bundle over $ \epsilon (M)$) that we denote by $A$. Also, we consider $ \rho_W$ as a vector bundle morphism from $A$ to $TM$ (over the identity of $M$).
              \item (\emph{Not easy!}) Show that the pair $ ( A , \rho_W) $ is an anchored bundle over $\mathcal F $ (i.e., $\rho_W \left( \Gamma(A)\right)=\mathcal F$).  
              {\emph{Hint:}} Use Lemma \ref{lem:localstructure} and more precisely its consequence, Exercise \ref{exo:sprojtproj}. 
         \end{enumerate}
         The results of this exercise will soon be generalized, see Lemma \ref{lem:BisectionAnchoredBundle} below.
\end{exo}

  \begin{lemma} 
  \label{lem:bisectionsExist}
 Consider $ M \stackrel{s}{\leftarrow} W \stackrel{t}{\rightarrow} M $ a bisubmersion over $ \mathcal F$.
 For every $x \in W$, there exists at least one bisection through $x$.
  \end{lemma}
  \begin{proof}
  There exists vector subspaces $V \subset T_x W$ of dimension equal to ${\mathrm{dim}}(M)$ whose intersection with $ \ker(T_x s)$ and with  $ \ker(T_x t)$ are both reduced to zero. In particular, both $T_x s\colon V \to T_{s(x)}M$ and $ T_x t\colon V \to T_{t(x)}M$ are linear invertible maps. Any submanifold of $W$ through $x$ admitting $V$ as its tangent space admits a restriction to a neighborhood of $x$ which is a bisection.
  \end{proof}
\begin{exo}
Prove the following statements.
\begin{enumerate}
\item Show that there exists infinitely many local bisections through a point $x \in W$, except maybe if $\mathcal F=0 $ in a neighborhood of $s(x)$.
\item Show that for any two local bisections $\Sigma_0,\Sigma_1 $ through $x \in W$, there exists a neighborhood $\mathcal U $ of $s(x) $ and a smooth family  $(\Sigma_\epsilon)_{\epsilon \in [0,1]}$ of bisections through $ x$ such that $ s(\Sigma_\epsilon)=\mathcal U$ for all $\epsilon \in [0,1] $.
\end{enumerate}
\end{exo}

\vspace{0.4cm}
\noindent
As for Lie groupoids, every bisection $\Sigma $ induces a diffeomorphism:
 $$ \underline{\Sigma} \colon s(\Sigma) \longrightarrow t(\Sigma),$$
 that makes the following diagram commutative:
 $$
 \xymatrix{ & \Sigma \ar[dl]_{s_{|_\Sigma}} \ar[dr]^{t_{|_\Sigma}}& \\ s(\Sigma) \ar[rr]^{\underline{\Sigma}}& &t(\Sigma)}
$$

\begin{exo}
\label{exo:units}
Show that units of a bisubmersion are in one-to-one correspondence with bisections $\Sigma$ such that $\underline{\Sigma} $ is the identity map of $M$.
\end{exo}

\noindent
Here is an important result, very similar to what happens for Lie groupoids\footnote{See Proposition 2.8 in \cite{AS}. All concepts and results of this section come from \cite{AS}, we will not give the exact correspondence from now on.}.
\vspace{0.5cm}

\begin{propositions}{Bisections induce Symmetries}{bisectionsandsymmetries}
Let $\mathcal F $ be a singular foliation on a manifold $M$ and let
$ M \stackrel{s}{\leftarrow} W \stackrel{t}{\rightarrow} M $ be a bisubmersion over $ \mathcal F$. For every local bisection $\Sigma$, the induced (local)  diffeomorphism $$\underline{\Sigma} \colon s(\Sigma) \longrightarrow t(\Sigma)  $$ is an isomorphism of singular foliations from\footnote{Also called “local symmetry of $\mathcal F $" in Section \ref{sec:symmetry}. 
Recall that $ \mathcal F_{s(\Sigma)}$ and  $ \mathcal F_{s(\Sigma)}$ stand for the restrictions of $ \mathcal F$ to $ s(\Sigma)$ and $ t(\Sigma)$, respectively.} $ \left(s(\Sigma), \mathcal F_{s(\Sigma)}\right)$ to $\left(t(\Sigma), \mathcal F_{t(\Sigma)}\right) $. In particular, when $ \Sigma$ is  a global bisection, $ \underline{\Sigma}$ is a symmetry of $(M, \mathcal F)$.
\end{propositions}

\vspace{.5cm}
\noindent
We prove Proposition \ref{thm:bisectionsandsymmetries} for a local bisection. The global case follows as a particular case. Let ${}^s\mathcal F{\,}^t $ be the bi-pull-back singular foliation on $W$ as in Equation \eqref{eq:singFolonW}.
Denote by $\left({}^s\mathcal F{\,}^t\right)_{\Sigma}$ its restriction to $\Sigma$ (see Section \ref{sec:transverse}). Recall that $\left({}^s\mathcal F{\,}^t\right)_{\Sigma}$, i.e., the subspace of $\mathfrak X(\Sigma) $ obtained by considering the restriction to $ \Sigma$ of all vector fields in ${}^s\mathcal F{\,}^t $ that are tangent to $ \Sigma$, need not be  a singular foliation in general.  However, it is so in the present situation:

\begin{lemma}
\label{lem:SigmaRestriction}
Let $\Sigma $ be a bisection  of $W$.
\begin{enumerate}
    \item 
The restriction $ \left({}^s\mathcal F{\,}^t\right)_{\Sigma}$ of ${}^s\mathcal F{\,}^t  $ to $\Sigma$ is a singular foliation on $ \Sigma $. 
\item Moreover, both $s_{|_\Sigma}$ and $t_{|_\Sigma} $ are diffeomorphisms of singular foliations from 
$ (\Sigma,{}^s\mathcal F_\Sigma{\,}^t ) $ to $(s(\Sigma),\mathcal F_{s(\Sigma)}) $ and $ (\Sigma,{}^s\mathcal F_\Sigma{\,}^t ) $ to $(t(\Sigma),\mathcal F_{t(\Sigma)})$, respectively.
\end{enumerate}
 \end{lemma}

\begin{proof}

The following decomposition holds $\forall \sigma \in \Sigma$:
 $$ T_\sigma \Sigma \oplus {\mathrm{ker}}(T_\sigma s)= T_\sigma W .
 $$ 
 Since $ {\mathrm{ker}}(T_\sigma s) \subset T_\sigma {}^s\mathcal F{\,}^t$, we therefore have
 $$   T_\sigma \Sigma+ T_\sigma {}^s\mathcal F{\,}^t= T_\sigma W .
$$
In terms of the concepts introduced in Section \ref{sec:transverse}, it means that $\Sigma$ intersects cleanly the bi-pull-back singular foliation ${}^s\mathcal F_\Sigma{\,}^t$. Proposition \ref{thm:prop:transverse} therefore applies:  $\left({}^s\mathcal F{\,}^t\right)_{\Sigma}\subset \mathfrak X(\Sigma)$ is a singular foliation on~$\Sigma$. {\emph{This proves item 1.}}

\vspace{1mm}
\noindent
Since the restriction $s_{|_\Sigma} $ of $s$ to $ \Sigma$ is a diffeomorphism, 
one can consider the push-forward singular foliation $$(s_{|_\Sigma})_*  \left(\left({}^s\mathcal F{\,}^t\right)_{\Sigma} \right) .$$ It is by construction a singular foliation on $s(\Sigma)$.
We first show the inclusion $\mathcal F_{s(\Sigma)}\subseteq ({s_{\Sigma}})_{*}(\left({}^{s}\mathcal F{\,}^t\right)_\Sigma)$. Let $u \in \mathcal F_{s(\Sigma)}$. There exists a unique vector field $u^{\Sigma}\in\mathfrak X(\Sigma)$ such that $ Ts_\Sigma(u^\Sigma)=u$. 
 Let $ v$ be a vector field in $s^{-1}(\mathcal F) $ that $s$-related to $u$. 
For every $\sigma\in \Sigma$, the difference $  u_{\sigma}^\Sigma-v_\sigma  $  is valued in ${\mathrm{ker}}(T_\sigma s)$. 
In view of the decomposition
 $$  T_\sigma W = T_\sigma \Sigma \oplus {\mathrm{ker}}(T_\sigma s) ,$$
there exists  a vector field in $Z \in \Gamma( {\mathrm{ker}}(T s))$ such that $u_{\sigma}^\Sigma -v_\sigma= Z_\sigma$ for all $ \sigma \in \Sigma$. Consider $u^! := v +Z $. The vector field $ u^!$ belongs to $ s^{-1}(\mathcal F)$ by definition. Also, by construction, its restriction to $\Sigma$ is tangent to $\Sigma$. Last, it coincides with $u^\Sigma$ on $ \Sigma$. This proves that $u^{\Sigma} \in \left({}^{s}\mathcal F{\,}^t\right)_\Sigma $. Therefore,  the desired inclusion holds.

\vspace{1mm}
\noindent
Let us show the opposite inclusion $(s_{|_{\Sigma}})_* \left( {}^{s}\mathcal F{\,}^t\right)_\Sigma  \subseteq \mathcal F_{s(\Sigma)} $ : let $ v$ be a vector field in ${}^{s}\mathcal F{\,}^t$ that happens to be tangent to $\Sigma$. We show that there exists another vector field $ \tilde{v}$ which coincides with $v$ on $ \Sigma$ and is $s$-related to a vector field in $\mathcal F$. By construction, $v= \sum_i g_i v_i$ where the $g_i$ are smooth functions on 
$W$ and $v_i\in {}^{s}\mathcal F{\,}^t $ are $s$-related to elements  $u_i \in \mathcal F$. Let $ \tilde{g}_i$ be functions on $W$ that coincide with $g_i$ on $\Sigma$ and are constant along $s$-fibers, i.e., $\tilde{g}_i=s^*f_i$ for some smooth function $f_i$ on $M$. The vector field $ \tilde{v} = \sum_i \tilde{g}_i v_i$ is tangent to $\Sigma $ (since it coincides with $v$ on $ \Sigma$) and is $s$-related with $\sum_i f_i u_i \in \mathcal F $. Hence $ s_*(v|_{\Sigma}) = \sum_i f_i u_i |_{s(\Sigma)}$
and in particular $(s_\Sigma)_*(v|_{\Sigma})\in \mathcal{F}_{s(\Sigma)}$. 
In turn, this implies $(s_{\Sigma})_* \left( {}^{s}\mathcal F_\Sigma{\,}^t \right) \subset \mathcal F_{s(\Sigma)} $. The same argument holds for $(t_{\Sigma})_* \left( {}^{s}\mathcal F_\Sigma{\,}^t \right) \subset \mathcal F_{t(\Sigma)} $. This completes the \emph{proof of item 2}.

\end{proof}

\begin{proof}[Proof of Proposition \ref{thm:bisectionsandsymmetries}]
From lemma \ref{lem:SigmaRestriction}, it follows that both $s_{|_\Sigma}$ and $t_{|_\Sigma}$ below are diffeomorphisms of singular foliations: 
$$  \xymatrix{ & \left(\Sigma,\left({}^s\mathcal{F}{\,}^t\right)_{\Sigma}\right) \ar[dl]_{s_{|_\Sigma}} \ar[dr]^{t_{|_\Sigma}} &  \\ (s(\Sigma),\mathcal F_{s(\Sigma)}) \ar@{..>}[rr]^{\underline{\Sigma}}& &  (t(\Sigma),\mathcal F_{t(\Sigma)}) }  $$
Hence the horizontal line $\underline{\Sigma}$ is a diffeomorphisms of singular foliations. This proves the claim. 
\end{proof}

\noindent
Let us now associate an anchored bundle over $\mathcal F $ to any  bisection $\Sigma $, that we will assume to be global for the sake of simplicity.
 The vector bundle morphism:
 $$  \begin{array}{rcl}TW_{|_\Sigma} &\longrightarrow &TM \\ u &\mapsto  & Tt (u) - T \underline{\Sigma} \circ Ts (u) \end{array}  $$
 admits $T\Sigma$ in its kernel. 
Therefore, it goes to the quotient to yield a vector bundle morphism $ \rho_\Sigma$
$$
\xymatrix{ TW_{|_\Sigma}/T\Sigma = N_\Sigma \ar[rr]^{\rho_\Sigma} \ar[d] & &TM  \ar[d]\\ \Sigma \ar[rr]^t& & M } $$
where $ N_\Sigma$ stands for the normal bundle of $ \Sigma$ in $W$. 
In view of the decomposition:
 $$ T_\sigma W = T_\sigma \Sigma \oplus {\mathrm{ker}}(T_\sigma s)  \hbox{ for all $ \sigma \in \Sigma$} ,$$
 the  restriction ${\mathrm{ker}}(Ts)_{|_\Sigma} $  of ${\mathrm{ker}}(Ts) $ to $\Sigma$ is, as a vector bundle, canonically isomorphic to the normal bundle $N_\Sigma$.
Under this isomorphism, the following diagram is commutative
  \begin{equation}\label{eq:anchor}  \xymatrix{ N_\Sigma \ar[rr]^{\simeq} \ar@{..>}[rrd]_{\rho_\Sigma} & & {\mathrm{ker}}(Ts)_{|_\Sigma}  \ar[d]^{Tt}\\ &  & TM }    .  \end{equation}

\vspace{2mm}
\noindent
Now, since $t : \Sigma \to M$ is a diffeomorphism on its image, 
$N_\Sigma$ can be considered as a vector bundle $t^{-1}N_\Sigma $ over the open subset $t(\Sigma) \subset M$ (rather than over $\Sigma $). Under this identification, $\rho_\Sigma$ becomes a vector bundle morphism $ t^{-1}N_\Sigma \to TM $  over the identity map of $t(\Sigma) \subset M$ that we denote by the same letter $ \rho_\Sigma$ and call the \emph{anchor of $\Sigma $}.

 \begin{lemma}
 \label{lem:BisectionAnchoredBundle}
 Let $\Sigma$ be a  bisection of a bisubmersion $W$ over a foliated manifold $(M,\mathcal F) $. The pair $(t^{-1}N_\Sigma, \rho_\Sigma)$ is an anchored bundle for the restriction of $\mathcal F $ to $t(\Sigma) $.
 \end{lemma}
\begin{proof}
Choose any point $ \sigma \in \Sigma.$
In some neighborhoods, there exists vector fields as in Lemma \ref{lem:localstructure}. By the fourth item, $Y_1, \dots, Y_{r_t}$ and $ Z_1, \dots, Z_k $ generate the kernel of $ Tt$, and therefore its restriction to $\Sigma $. They identify therefore to sections of the normal bundle that generate it. By Equation \eqref{eq:anchor},
their images under $ \rho_\Sigma$ are the local generators $ {\underline{Y}}_1, \dots, {\underline{Y}}_{r_t} $ of $ \mathcal F$. This completes the proof.
\end{proof}

\vspace{.5cm}
\begin{center} 
\emph{The end of this subsection is required for what follows, but is more technical/involved. We recommend the reader unfamiliar with the notion to look directly at Section \ref{sec:Example}}.
\end{center}
\vspace{.5cm}

\noindent
Consider two local bisections  $\Sigma_0,\Sigma_1$ of a bisubmersion $W$ that contain the same point $w \in W$. In general, the induced local symmetries $\underline{\Sigma_0}, \underline{\Sigma_1} $  of $(M,\mathcal F) $ do not coincide, even locally. However, we will see that they differ by an inner symmetry admitting $s(w)$ as a very-fixed point\footnote{See Definition \ref{def:veryfixed}.},
(equivalently, they differ by an inner symmetry admitting $t(w)$ as a very-fixed point).
We refer to Section \ref{sec:symmetry} for the vocabulary about inner symmetries, and very-fixed points of those. We also refer to this section\footnote{The parameter is denoted by $t$ in that section and by $ \epsilon$ here, since $t$ is the target map.} for the subtleties about the precise meaning of ``smoothness'' for a time-dependent vector field $ (X_\epsilon)_{\epsilon \in I}$. 

\begin{proposition}
\label{prop:differSymm}
 Consider a bisubmersion $W$ over a foliated manifold $(M,\mathcal F) $. Let $w \in W$.
 For any two local bisections $\Sigma_0,\Sigma_1$ of $W$ through $w$,
 the local diffeomorphism $\underline{\Sigma_0} \circ \underline{\Sigma_1 }^{-1} $ is an inner symmetry of $ \mathcal F$ with very-fixed point\footnote{See Definition \ref{def:veryfixed}.} $ t(w) $, at least in a neighborhood of $ t(w)$. 

\vspace{1mm}
\noindent
 Said differently, there exists  $ (X_\epsilon)_{\epsilon \in [0,1]}$, a smooth time-dependent vector field in $\mathcal F $, such that 
\begin{equation} \label{eq:correctform}\left. X_\epsilon\right|_{x} \, = \, \sum_{i=1}^r f_i(x,\epsilon) \left. X_i \right|_{x},\end{equation}
with $ X_1, \dots, X_r$ being local generators of $ \mathcal F$, and $ f_1, \dots,f_r$ being smooth functions that vanish on $ \{m\} \times [0,1]$, whose time-$1$ flow $\Phi_1^{X_\epsilon} $ coincides with 
 $ \underline{\Sigma_0} \circ \underline{\Sigma_1 }^{-1}   $ in a neighborhood of $ t(x)$.

\end{proposition}

\vspace{1mm}
\noindent
The proof starts with a Lemma.
Let $(\Sigma_\epsilon)_{\epsilon \in I} $ be a family of global bisections. We say that such a family is \emph{smooth} if the sections $ \Sigma^{-1}_\epsilon \colon M \to W  $ that inverts $s_{|_{\Sigma_\epsilon}} \colon  \Sigma_\epsilon \to W$ depends smoothly on the parameter~$ \epsilon$.  

\begin{lemma}
\label{lem:familySmooth}
If there exists a smooth interpolation by bisections between two given bisections $\Sigma_0 $ and $\Sigma_1 $, then $\underline{\Sigma_1 }^{-1} \circ  \underline{\Sigma_0} $ is a local inner symmetry of $\mathcal F $.
\end{lemma}
\begin{proof}
\noindent
For any small enough $ \epsilon \in \mathbb R$ and any $w \in \Sigma_\epsilon $, the path
$$ \eta \mapsto s_{|_{\Sigma_\eta}} (  s(m))$$
is a path in a fixed $s$-fiber and goes through $m$ at $\epsilon =\eta $. Its derivative:
 $$  u_\epsilon :=  \left. \frac{\partial s_{|_{\Sigma_\eta}} (  s(m)) }{\partial \eta } \right|_{\eta=\epsilon} $$ 
 is therefore a section of ${\mathrm{ker}}(Ts)_{|_{\Sigma_\epsilon}} $, that naturally identifies to a section of the anchored bundle $(A_{\Sigma_\epsilon},\rho_{\Sigma_\epsilon}) $ of Lemma \ref{lem:BisectionAnchoredBundle}.   Its image through $\rho_{\Sigma_\epsilon} $ is a vector field that we denote by $ X_\epsilon$. In view of Lemma \ref{lem:BisectionAnchoredBundle}, $ X_\epsilon $ belongs to $ \mathcal F$. It is by construction a smooth time-dependent vector field in $\mathcal F $ as defined in Section \ref{sec:symmetry}.  
The relation $Tt (u_\epsilon) = X_\epsilon $ holds by construction, and implies that $ \eta \mapsto t \circ s_{|_{\Sigma_\eta}} (  s(w))   $ is an integral curve of the vector field $X_\epsilon $, that relates $ \underline{\Sigma_0}(s(w)) $ to $\underline{\Sigma_1}(s(w))  $ by construction. This proves the claim.
\end{proof}

\begin{proof}[Proof of Proposition \ref{prop:differSymm}]
Upon replacing $ \Sigma_0,\Sigma_1$ by neighborhoods of $ w$  in   $ \Sigma_0,\Sigma_1$  that we denote by the same symbols, one can assume that a smooth family of local bisections interpolating between $ \Sigma_0$  and $\Sigma_1 $ exists. Moreover, those can be chosen such that \emph{(i)} we have $w \in \Sigma_\epsilon $ for all $ \epsilon \in [0,1]$, and \emph{(ii)} there exists a neighborhood $ \mathcal U$ of  $ t(w)$ contained in  
$t(\Sigma_\epsilon) $  for all value of $\epsilon \in [0,1] $. This last point allows to use Lemma \ref{lem:familySmooth} (adapted to local bisections). It therefore shows that $ \phi=\underline{\Sigma_0} \circ \underline{\Sigma_1}^{-1}$ in an inner symmetry of $ \mathcal F$. Since all the bisections $\Sigma_\epsilon $ are through $w$, the section $u_\epsilon $ that appears in Lemma \ref{lem:familySmooth} vanishes at $t(w) $. As a consequence (see Exercise \ref{exo:veryfixed}), the smooth time-dependent vector field $(X_\epsilon)_{\epsilon \in [0,1]}$, whose time-$1$ flow is $ \phi$, is of the type given in Equation \eqref{eq:correctform}. 
The construction in Lemma \ref{lem:familySmooth} yields therefore a smooth time-dependent vector field in $\mathcal F $ whose time $1$-flow satisfies the required assumptions that makes $ t(w)$ a very fixed point.

\end{proof}

 \vspace{1mm}
 \noindent
For any two points in a foliated manifold $(M,\mathcal{F})$.
Let us denote by $ \mathrm{Sym}_{\mathcal F}(m,m')$ the set of local symmetries of $ \mathcal F$, defined from a neighborhood of $m$ to a neighborhood of $m'$. We denote by $\mathrm{Sym}_{\mathcal F}(m,m')$  and call germs of such local symmetries the quotient of $\mathrm{Sym}_{\mathcal F}(m,m')$ through the equivalence relation that identifies two diffeomorphisms that coincide in a neighborhood of $m$.
Of course, in general, 
$ \mathrm{Sym}_{\mathcal F}(m,m') $ may be empty : it is however, non-empty if $m$ and $m'$ are in the same leaves: see Section \ref{sec:leavesExist}.

 \vspace{1mm}
 \noindent
Let us introduce an equivalence relation on 
$ \mathrm{Sym}_{\mathcal F}(m,m') $.
We identify two elements $\bar{\Phi}, \bar{\Psi} $ if one of the equivalent conditions below holds:
\begin{enumerate}
\item[(i)] $  \Phi^{-1} \circ \Psi $ coincides  in a neighborhood of $m$ with an inner symmetry of $ \mathcal F$ with very-fixed point $m$. 
\item[(ii)] $  \Psi^{-1} \circ \Phi $ coincides  in a neighborhood of $m'$ with an inner symmetry with very-fixed point $m'$.
 \end{enumerate} 

\begin{exo}
Show that both conditions (i) and (ii) above are indeed equivalent, and that the latter is an equivalence relation.
\emph{Hint:} Look at exercise \ref{exo:InnerAndSymmetry2}.
\end{exo}

\noindent
We now denote  the quotient of $\mathrm{Sym}_{\mathcal F}(m,m')$ by the previous equivalence relation by $\mathrm{OutSym}_{\mathcal F}(m,m') $. We  call this set the \emph{germs of outer symmetries from $m$ to $m'$}.

\vspace{1mm}
\noindent
Let $W$ be a bisubmersion for $(M,\mathcal F) $.
It follows from Proposition \ref{thm:bisectionsandsymmetries} that to any pair $ (x,\Sigma)$, with $\Sigma$ a bisection of $W$ through $x$, one can associate the germ of $\underline{\Sigma} \in {\mathrm{Sym}}_\mathcal F(s(x),t(x)) $.
By Proposition \ref{prop:differSymm}, the class of $\underline{\Sigma} $ in ${\mathrm{OutSym}}_\mathcal F(s(x),t(x)) $ does not depend on the choice of $\Sigma $.   In view of Lemma \ref{lem:bisectionsExist}, we have therefore defined a map  that makes the following diagram commute: 
$$  \xymatrix{W \ar[rr]^{{\mathrm{OutSym}}_W } \ar@<-.5ex>[d]_s \ar@<.5ex>[d]^t &&  \mathrm{Out}_{\mathcal F}\ar@<-.5ex>[d]_s \ar@<.5ex>[d]^t \\M\ar@{=}[rr]&&M }  $$
where $\mathrm{OutSym}_{\mathcal F} = \coprod_{m,m' \in M} \mathrm{OutSym}_{\mathcal F} (m,m')$ and where the $s$ and $t$ are defined for any element in $  \mathrm{OutSym}_{\mathcal F} (m,m') $ to be $m$ and $m'$, respectively.

\vspace{0.5cm}

\begin{definitions}{Outer-germ map}{outergerm}
We call \emph{Outer-germ map} and denote by ${\mathrm{Out}}_W $ the map 
$$  W \longrightarrow  \mathrm{OutSym}_{\mathcal F} $$
defined by associating to any $w \in W$ the class, modulo inner symmetries with very-fixed points\footnote{See Definition \ref{def:veryfixed}.}, of the germ of symmetries\footnote{See Section \ref{sec:symmetry}.} of $ \mathcal F$ induced by an arbitrary\footnote{The discussion before the statement explains why any two bisections through $w$ induce the same class in the quotient space.} local bisection through $ w$.
\end{definitions}

\vspace{0.cm}

\begin{exo}
\label{exo:OutConstAlongBiverti}
We intend to prove that the map $ \mathrm{Out}_{W}$ is constant along integral curves of bi-vertical vector fields.
Let $X \in \mathfrak X (W)$ be a bi-vertical vector field on a bisubmersion $W$ over $(M,\mathcal F) $. Let $\Phi_t^X $ be the flow of $X$ and let $\epsilon \mapsto \gamma(\epsilon) = \Phi_\epsilon^X(m) $ be an integral curve of $ X$ starting from $x \in W$, defined for all $\epsilon $ is some interval $I$ containing $0$.
\begin{enumerate}
\item Let $\Sigma $ be a bisection of $W$ through $\gamma(0) $. Show that for all $\epsilon \in I $, there exists a neighborhood $\Sigma' $ of $x$ in $\Sigma $ on which $\Phi_\eta^X$ is well-defined for all $\eta \in [0,\epsilon] $.
\item Show that $\Sigma_\eta := \Phi_\epsilon^X (\Sigma') $ is a bisection for all $\eta  \in [0,\epsilon]$, and that
 $\underline{\Sigma_\eta} $ is constant for all $\eta \in [0,\epsilon] $.
\item Conclude that the map $ \mathrm{Out}_{W}$ is constant along integral curves of bi-vertical vector fields.
\end{enumerate}
\end{exo}

\begin{exo}
Let $(M,\mathcal F)$ be a foliated manifold, $W$ and $W'$ two bisubmersions over $ \mathcal F$. Show that if $x \in W$ and $x' \in W'$ are related by $ {\mathrm{Out}}_W(x)={\mathrm{Out}}_W(x')  $, then there exists local bisections $\Sigma $ and $\Sigma' $ through $ x$ and $x'$ such that $\underline{\Sigma} = \underline{\Sigma'} $.
(\emph{Hint}: We will have more tools to solve this exercise after right and left actions are defined. This exercise will be then repeated as Exercise \ref{exo:Outiff}.)
\end{exo}

 \subsection{Examples of bisubmersions}

\label{sec:Example}

\subsubsection{Basic examples and non-examples}

\noindent
The first exercise shows that Lie groupoids are bisubmersions - a non-trivial fact.
 
 \begin{exo}
 \label{exo:groupasbosubmersions}

 Let $\Gamma $ be a Lie groupoid over $M$. Let $A$ be its Lie algebroid and $\rho \colon A \to TM $ its anchor map. The basic singular foliation is (by definition) the singular foliation $\mathcal F := \rho\left(\Gamma(A)\right) $ (see Section \ref{sec:LieAlgebroidsAreSingFoliation}). 
The purpose of this exercise is to show that any  Lie groupoid $\Gamma $ is a bisubmersion for its basic singular foliation $\mathcal F $.

\vspace{1mm}
\noindent
 Let us denote by $\overrightarrow{a} $ and $\overleftarrow{a} $ the right and left invariant vector fields on $\Gamma $ associated to $a \in \Gamma(A)$ as in \cite{Mackenzie}. There are miscellaneous conventions: we chose them such that
  $ \overrightarrow{a}$ is $s$-related to $0$ and $t$-related to $\rho(a) $, while $ \overleftarrow{a})$ is $t$-related to $0$ and $s$-related to $\rho(a) $.
\begin{enumerate}
\item  Show that the $ \mathcal C^\infty(\Gamma)$-module $\mathcal G$ generated by vector fields of the form 
$ \overrightarrow{a} + \overleftarrow{b} $ for some $ a,b \in \Gamma(A)$ is a singular foliation on $\Gamma$. 
\item In passing: describe its leaves. We may assume that $\Gamma$ is source-connected for simplicity.
\item 
Use the Lie groupoid axioms to show that 
$ \Gamma({\mathrm{ker}}(Tt))  $ (resp. $ \Gamma({\mathrm{ker}}(Ts))  $) coincides with the $ \mathcal C^\infty(\Gamma)$-module generated by $ \{ \overrightarrow{a} \mid a \in \Gamma(A)\} $ (resp.  $ \{ \overleftarrow{a} \mid a \in \Gamma(A)\} $).
\item Show that $\mathcal G = \Gamma({\mathrm{ker}}(Ts)) + \Gamma({\mathrm{ker}}(Tt))$.
\item Show that $ \mathcal G \subset  s^{-1}(\mathcal F)  $ and  $ \mathcal G \subset  t^{-1}(\mathcal F)  $.
\emph{Hint}: prove that $ \mathcal G $ is generated by vector fields $s $-related (resp. $t$-related) to vector fields in $ \mathcal F$.
\item Show that $s^{-1}(\mathcal F)  \subset \mathcal G  $ and $ t^{-1}(\mathcal F)) \subset \mathcal G$.
{\emph{Hint:}} Show that every $X \in \mathcal F$ is $s$-related  (resp. $t$-related) to a vector field in  $\Gamma({\mathrm{ker}}(Tt))$, equivalently $\Gamma({\mathrm{ker}}(Ts))$.
\item Conclude that $\Gamma $ is a bisubmersion for its basic singular foliation $\mathcal F $.
\item Show that ${\mathrm{Out}}_{\Gamma} \colon \Gamma \to {\mathrm{OutSym}}_{\mathcal F}$ is a groupoid morphism (the groupoid structure on ${\mathrm{OutSym}}_{\mathcal F}$ is defined in the next section \ref{sec:holonomygroupoid}). 
\end{enumerate}

 \end{exo}

\vspace{0.5cm}
\noindent
It is also very important to have non-examples in mind.

 \begin{exo}
 Here are non-examples of bisubmersions. 
 Let $\mathcal F$ be a singular foliation on $M$.
 \begin{enumerate}
     \item Show that $W:=M \times M $, equipped with the projections  onto the first and second components as source and target, is \emph{not} a bisubmersion over $\mathcal F$, unless $\mathcal F = \mathfrak X_c(M) $ is the singular foliation of all compactly supported vector fields on $M$.
     \item Show that $W:=M  $ equipped with the identity map as source and target is \emph{not} a bisubmersion over $\mathcal F$, unless $\mathcal F = 0 $.
     \item Give an example of a manifold $W$, equipped with two surjective submersions $s,t\colon W \to M$, that do satisfy $s^{-1}(\mathcal F) =t^{-1}(\mathcal F)$, and is still not a bisubmersion over $\mathcal F $.
 \end{enumerate}
 \end{exo}

\subsubsection{The crucial example: A bisubmersion for every finitely generated singular foliation}

\noindent
There is a very natural bisubmersion over any finitely generated singular foliation $\mathcal F $ on $M$. Let $ X_1 , \dots, X_r$ be generators of $ \mathcal F$. For the sake of simplicity, we assume that they are complete vector fields, and leave it to the reader to generalize. Consider the following triple:
\begin{enumerate}
    \item The manifold $\mathbb R^r \times M$.
    \item The map $s \colon \mathbb R^r \times M \to M $ given by the projection on the second factor.
    \begin{equation}\label{eq:projfirstfactor} \begin{array}{rrcl} s\colon   & \mathbb R^r \times M &\to & M \\ &((\lambda_1 \dots, \lambda_r) , m)& \to & m \end{array}.\end{equation}
    \item The map $t \colon \mathbb R^r \times M \to M $ given by:
      \begin{equation}\label{eq:projsecondfactor} \begin{array}{rrcl} t\colon   & \mathbb R^r \times M &\to & M \\  & ((\lambda_1 \dots, \lambda_r) , m) &\mapsto & 
      \phi^{X_1 }_{\lambda_1} \circ \cdots \circ \phi^{X_r}_{\lambda_r} (m) 
   \end{array}.
     \end{equation}

\end{enumerate}
\noindent  If one do not assume that the vector fields $ X_1, \dots,X_r$ are complete vector fields, then the previous map $t$ still makes sense, but the initial manifold $\mathbb R^r \times M$ has to be replaced by a neighborhood of $ (0, \dots, 0) \times M $ (i.e., of the zero section of the trivial bundle $ \mathbb R^r \times M \to M$). The map $s$ is always a surjective submersion. The map $t$ is also a surjective submersion at least, again, in a neighborhood of the zero section.

\vspace{1mm}
\noindent
The following proposition is very important and very non-trivial. Here, our inspiration is not \cite{AS} but we adapted a proof by Claire Debord \cite{debord2000local}. We  will soon state and prove a slightly different but more general statement (see Proposition \ref{thm:Prop:FundamentalExampleAbstract} below).

\vspace{.5cm}

\begin{propositions}{A crucial example}{prop:crucialExample}
Let $\mathcal F $ be a finitely generated singular foliation.
There is a neighborhood  $ \mathcal V$ of the zero section in $\mathbb R^r \times M \to M $ such that $M \stackrel{s}{\leftarrow} \mathcal V \stackrel{t}{\rightarrow} M $, with $s,t$ as in Equations \eqref{eq:projfirstfactor}-\eqref{eq:projsecondfactor}, is a bisubmersion over $\mathcal F $.
\end{propositions}

\begin{remark}
\normalfont
 If $\mathcal F $ is generated, as a module over smooth functions, by vector fields which are integrable by quadrature ($ \simeq $ such that the flow can be computed  ``by hand''), then Proposition \ref{thm:prop:crucialExample} means that there is an \emph{explicit} bisubmersion for $ \mathcal F$.
\end{remark}

\begin{proof} Recall from \eqref{eq:projsecondfactor} that the target map $t$ is given by $ ((\lambda_1 \dots, \lambda_r) , m) \mapsto \phi^{X_1 }_{\lambda_1} \circ \cdots \circ \phi^{X_r}_{\lambda_r} (m) $ and from Equation \eqref{eq:projfirstfactor} that the source map $s$ is given by $ ((\lambda_1 \dots, \lambda_r) , m) \mapsto m$. The vector field on $ \mathbb R^r \times M $ given by $$Z_i\colon ((\lambda_1 \dots, \lambda_r) ,m)\mapsto (0, \left.X_i\right._{|_m})$$ is $s$-related to $X_i\in \mathcal{F}$ for every $i=1, \ldots,r$.
Hence, $ Z_i \in s^{-1}(\mathcal F)$. As a consequence, $\lambda_i Z_i  \in s^{-1}(\mathcal F) $ so that its flow $\phi_{1}^{\lambda_i Z_i} = \phi_{\lambda_i}^{Z_i} $ is a symmetry\footnote{(at least whenever it is defined - hence a local symmetry)} of $s^{-1}(\mathcal F) $ by Corollary \ref{coro:flowissymmetry}. Hence,
 $ \phi^{Z_1 }_{\lambda_1} \circ \cdots \circ \phi^{Z_r}_{\lambda_r} $
is a symmetry of $s^{-1}(\mathcal F) $ - at least on the open set where it makes sense.
Now, notice that $$s\circ \phi^{Z_1 }_{\lambda_1} \circ \cdots \circ \phi^{Z_r}_{\lambda_r} (\lambda_\bullet,m)=\phi^{X_1 }_{\lambda_1} \circ \cdots \circ \phi^{X_r}_{\lambda_r} (m)=t(\lambda_\bullet ,m)$$ at least on the open set where it makes sense\footnote{Here and below, we use the shorthand $\lambda_\bullet=(\lambda_1, \dots,\lambda_r)$. Also, we will not recall at every line that expressions are only defined on an open subset containing the zero section $M$.}.  
Since there is a symmetry $\Phi$ of $s^{-1}(\mathcal F) $ such that $s \circ \Phi = t $, we have $ s^{-1}(\mathcal F) = t^{-1}(\mathcal F)$ at least in a neighborhood of  $\{0\}  \times M$.

 \vspace{1mm}
 \noindent
 The rest of the proof uses the same techniques as in \cite{debord2000local}.
 We have to prove
that any vector field in $s^{-1}(\mathcal F) = t^{-1}(\mathcal F) $ 
is the sum of a vector field in ${\mathrm{Ker}}(Ts) $ and a vector field in  ${\mathrm{Ker}}(Tt) $
at least in a neighborhood of  $\{0\}  \times M$. 
 Since we have already proven that $s^{-1}(\mathcal{F})=t^{-1}(\mathcal{F})$, and since $\Gamma(\ker Ts) \subset t^{-1}(\mathcal F)$ while $\Gamma(\ker Tt) \subset s^{-1}(\mathcal F) $ , we already have the inclusion $$t^{-1}(\mathcal{F})\supset \Gamma(\ker Ts)+\Gamma(\ker Tt).$$ We need to show the converse inclusion. For every $i\in \{1, \ldots,r\}$, let  ${\mathbf A}_i(\lambda_i,m):=(a_{ij}^l)_{jl}(\lambda_i,m)$ be a $r\times r$ matrix with coefficients in $\mathcal C^\infty(M)$ such that ${\mathbf A}_i(0,m)=\mathrm{Id}$ for all $m\in M$ and
\begin{equation}
    \label{eq:sym-matrix}\begin{pmatrix} \left(\phi^{X_i}_{\lambda_i}\right)_*(X_1) \\\vdots\\
\left(\phi^{X_i}_{\lambda_i}\right)_*(X_r) \end{pmatrix}(m)= {\mathbf A}_i(\lambda_i, m)\begin{pmatrix}
    X_1\\\vdots\\X_r
\end{pmatrix}(m).
\end{equation}

Such a matrix exists by Proposition \ref{prop:symInt} Section \ref{sec:leavesExist}. A direct computation gives for every $i=1 \dots,r$:
\begin{align*}
    T_{(\lambda, m)}t\left(\frac{\partial}{\partial \lambda_i}\right)&=\left.\frac{d}{d\epsilon}\right|_{\epsilon =0} \left(\phi^{X_1 }_{\lambda_1} \circ \cdots \circ \phi^{X_i}_{\lambda_i+\epsilon} \left(\phi^{X_{i+1} }_{\lambda_{i+1}} \circ \cdots \circ \phi^{X_r}_{\lambda_r}  (m)\right) \right)\\&=T\left(\phi^{X_1 }_{\lambda_1}\circ \cdots \circ  \phi^{X_{i-1}}_{\lambda_{i-1}}\right) \left(\left. X_i\right|_{\phi^{X_{i}}_{\lambda_{i}} \circ\phi^{X_{i+1}}_{\lambda_{i+1}} \circ \cdots \circ \phi^{X_r}_{\lambda_r}(m)}\right)\\&=T \left(\phi^{X_1 }_{\lambda_1} \circ \cdots \circ \phi^{X_{i-1}}_{\lambda_{i-1}}\right) \left(\left.X_i\right|_{\phi^{X_{i-1}}_{-\lambda_{i-1}} \circ \cdots \circ \phi^{X_{1}}_{-\lambda_1}(t(\lambda_\bullet,m))}\right)\\
    &= \left. \left(\phi^{X_1 }_{\lambda_1} \circ \cdots \circ \phi^{X_{i-1}}_{\lambda_{i-1}}\right)_* X_i \, \, \, \right|_{t(\lambda_\bullet,m)} \\
    &= \left. \left(\phi^{X_1 }_{\lambda_1} \right)_* \circ \cdots \circ \left(\phi^{X_{i-1}}_{\lambda_{i-1}}\right)_* X_i \, \, \, \right|_{t(\lambda_\bullet,m)}
\end{align*}
which, in view of Equation \eqref{eq:sym-matrix} gives
$$T_{(\lambda_\bullet, m)}t\begin{pmatrix}
  \frac{\partial}{\partial \lambda_1}\\\vdots\\ \frac{\partial}{\partial \lambda_r}
\end{pmatrix}= \mathfrak M(\lambda_\bullet, m)\begin{pmatrix}
        X_1\\\vdots \\X_r
    \end{pmatrix}(t(\lambda_\bullet,m))$$

where $\mathfrak M(\lambda_\bullet, m)$ is the $r\times r$ matrix  whose $i$-th line is the $i$-th line of the following matrix:
 
$$ \begin{array}{rcl}\mathfrak M(\lambda_\bullet, m) &=& {\mathbf A}_{i-1}\left(\lambda_{i-1},\phi^{X_{i-2}}_{-\lambda_{i-2}} \circ \cdots \circ \phi^{X_{1}}_{-\lambda_1}\circ t(\lambda_\bullet ,m)\right)  \\
&& \times \,  {\mathbf A}_{i-2}\left(\lambda_{i-2},\phi^{X_{i-3}}_{-\lambda_{i-3}} \circ \cdots \circ \phi^{X_{1}}_{-\lambda_1}\circ t(\lambda_\bullet,m)\right)  \\ & &\times \cdots  \\ &&\cdots \times {\mathbf A}_1(\lambda_1,t(\lambda,m))\end{array}$$
for all $1\leq i\leq r$. Here $ \times$ stands for the product of $ r \times r$ matrices
.
Despite its extraordinary description, the only point that matters in that this matrix exists and that,  
since $\mathfrak M(0, \ldots, 0, m)= \mathrm{Id}$ for all $m\in M$,  it is invertible on an open subset $\mathcal{V}$ of $\mathbb R^r\times M$ that contains $\{0\}\times M$. In turn, this invertibility implies there exists for every $i \in 1, \dots, r$ a linear combination of the vector fields $\frac{\partial}{\partial \lambda_1}, \dots, \frac{\partial}{\partial \lambda_r}$, depending smoothly on the point in $ \mathcal V$,  which is $t$-related to $ X_i$. Said otherwise, there exists for every $i \in 1, \dots, r$ a vector field $Z_i \in \mathfrak{X}(\mathcal{V})$ which lies inside ${\mathrm{ker}}(Ts) $ and is $t$-related with $ X_i$. By construction,  $ t^{-1}(\mathcal F) $ is generated by $Z_1, \dots, Z_r $ and vector fields in ${\mathrm{ker}}(Tt) $. This implies that,
after restriction to $\mathcal{V}$:

$$s^{-1}(\mathcal{F})=t^{-1}(\mathcal{F})\subseteq \Gamma(\ker Ts)+\Gamma(\ker Tt).$$
This completes the proof.
\end{proof}

\subsubsection{The previous example made more abstract: anchored bundles as bisubmersions}

\vspace{1mm}
\noindent
Proposition \ref{thm:prop:crucialExample} can be made more abstract by using an  anchored bundle\footnote{We refer to Section \ref{sec:AnchoredBundle} for results are the existence of anchored bundles. In short: in the smooth case, it exists if and only if $\mathcal F $ is finitely generated, and it always exists locally.} over $ \mathcal F$. 
The idea is to mimic the construction of  the so-called parallel $A$-paths in integration of Lie algebroids as in  \cite{MR1973056}.
Let $ \mathcal F$ be a singular foliation on $M$. Let  $(A,\rho) $ be an anchored bundle such that $\rho(\Gamma(A)) = \mathcal F$. We denote by $\pi \colon A \to M $ the projection onto the base manifold.

\vspace{0.5cm}

\begin{definitions}{Anchored paths}{anchoredpaths}
Let $ (A,\rho)$ be an anchored bundle over $ \mathcal F$.
We say that a path $a\colon I \to A $ is \emph{anchored} if 
 $$ \frac{d \gamma(t)}{dt} = \rho_{\gamma(t)} (a(t))  $$
 where $\gamma = \pi \circ a : I \to M $ is
the projection of $a(t) $ onto $M$.
\end{definitions}
\vspace{0.3cm}
\noindent
Notice that the path $ \gamma (t) = \pi \circ a (t)$ in Definition \ref{def:anchoredpaths} can not ``jump'' from one leaf of $ \mathcal F$ to another leaf.

\begin{exo}
    Show that for any anchored path $ t \mapsto a(t)$, there exists a smooth time-depending section $ (a_t)_{t \in I}$ of $A $ such that $ a_t(\gamma(t))=a(t)$. Show that the $ \gamma(t)$ is an integral curve of a smooth time-dependent vector field in $\mathcal F $.
\end{exo}

\vspace{0.3cm}
\noindent
Let us choose an affine  connection\footnote{They always exist in the smooth setting, and always exist locally in the complex or real-analytic settings.} on $A$:
$$\begin{array}{rrcl} \nabla :& \mathfrak X(M) \times \Gamma(A) & \to  & \Gamma(A) \\ &(X,a) & \mapsto & \nabla_X a \end{array} .$$

\vspace{.5cm}

\begin{definitions}{Parallel anchored paths}{Parallelanchoredpath2}
Let $A \stackrel{\pi}{\to}M$ be an anchored bundle and $\nabla $ be a connection on $A$.
 We say that an anchored path $a(t)$ is  $\nabla $-parallel if it satisfies:
  $$ \nabla_{\dot{\gamma}(t)} a(t) = 0 $$
  where  $\gamma = \pi \circ a : I \to M $ is
the projection of $a(t) $ onto $M$.
\end{definitions}

\vspace{.5cm}
\noindent
Here is a result which is purely a differential geometry result. 

\begin{lemma}
\label{lem:geodesics}
Every element $a\in A $ is the starting point of a $ \nabla$-parallel anchored path $t \mapsto a^{\nabla}(t) $, defined on some open interval containing $0 \in \mathbb R$, that we call the \emph{geodesic starting at $a \in A$}.
Moreover, there is a unique linear vector field $\Xi$ on $A$ whose integral curves are these geodesics, i.e.,
 $$\begin{array}{rcll} a^\nabla(t) &=& \Phi_t^\Xi (a) &\hbox{ for all $a\in A,$}\\ & & &\hbox{ and all $t \in \mathbb R$ for which it is defined.}  \end{array}$$
 We call $\Xi$ the \emph{geodesic vector field}.
\end{lemma}
\begin{proof}
The affine connection $ \nabla$ can be seen a sub-vector bundle $ H^\nabla \subset TA $ in direct sum with the kernel of $T\pi $. In particular, for any $a\in A$, $ T_a\pi \colon  H^\nabla \longrightarrow T_{\pi(a)} M $ is one-to-one. We denote its inverse map by $u \mapsto H^\nabla_a (u)$ for any $u \in T_{\pi(a)}M$. We define $\Xi$ by 
$$ \Xi_{|_a} :=   H^\nabla_a  \left( \rho_{\pi(a)}(a)\right)  .$$ 
Equivalently, $ \Xi$ is the only section of $H^\nabla $ such that $ T_a \pi (\Xi_{|_a}) = \rho(a) $ at every point $a \in A$. It is routine to check that the integral curves of $ \Xi$ are precisely the above defined geodesics.
\end{proof}

\begin{exo}
\label{exo:constantpath}
Show that the geodesic starting from a point of the zero section is a constant path.
\end{exo}

\begin{exo}
Show that for $A=TM$ and $\rho={\mathrm{id}} $, the geodesics defined above are the usual geodesics of an affine connection on $M$.
\end{exo}

\begin{exo}
Let $(A, \rho)$ be an anchored bundle for a singular foliation $\mathcal F $. Let $ \nabla$ be a connection on $A$. Choose $m \in M$. 
\begin{enumerate}
\item Show that for any $a \in A_m$, the geodesic path starting from $a$ remains inside the leaf $L_m$ through $M$, i.e., $ t \mapsto \pi \circ \Phi^t_\Xi (a) \in L_m$ for all $t$ for which it is defined.  
\item Show that there is a neighborhood $\mathcal V $ of $ m$ in $L_m$ such that for any $m' \in \mathcal V$, $m$ and $m'$ are the starting and ending point of the base path of some geodesic.
\item Here is an open question: are any two points in the same leaf the starting and end point of the base path of a geodesic?
\end{enumerate}
\end{exo}

\noindent
For a given $a \in A$, the geodesic path $ t \mapsto a^{\nabla}(t) = \Phi_t^\Xi(a)$ may not be defined for all $t$, but it follows from Exercise \ref{exo:constantpath} there is a neighborhood $ \mathcal U_A \subset A$ of the zero section where it is defined for all $t \in [0,1] $.
We denote by $\Phi^\Xi_1 $ the map from $ \mathcal U_A $ to $A$ that sends $a$ to $ \Phi^\xi_1 (a) = a^\nabla(1) $.
 Now, let us consider the triple made of 
\begin{enumerate}
    \item the neighborhood $ \mathcal U_A$ of the zero section in $A$,
    \item the projection $ A \stackrel{\pi}{\to} M $ that we rename $s$,
    \item the composition $t$ of $a \mapsto \Phi_1^\Xi (a) $ with the projection $  \colon A \stackrel{\pi}{\to} M $.
\end{enumerate}
Notice that, by construction:
\begin{equation}
\label{eq:stPhi}
t = s \circ \Phi_1^\Xi
\end{equation}
\vspace{.5cm}

\begin{propositions}{The same as Proposition \ref{thm:prop:crucialExample} but more abstract: anchored bundles are "fundamental" bisubmersions.}{Prop:FundamentalExampleAbstract}

Let $(A,\rho)$ be any anchored bundle such that $ \mathcal F = \rho(\Gamma(A))$, and $ \nabla$ be a connection on $A$. There is a neighborhood $\mathcal A $ of the zero section in $A$ on which $ M \stackrel{s}{\leftarrow} \mathcal A \stackrel{t}{\rightarrow} M $, with $s,t$ as above, is a bisubmersion of $\mathcal F $.

\end{propositions}
\vspace{.5cm}

We call \emph{fundamental bisubmersion of $ \mathcal F$} a bisubmersion $M \stackrel{s}{\leftarrow}\mathcal A \stackrel{t}{\rightarrow} M$
of the form given in Proposition \ref{thm:Prop:FundamentalExampleAbstract}, associated to some anchored bundle $(A,\rho)$ and some connection~$ \nabla$.
We start with a lemma about the geodesic vector field $\Xi$ of Lemma \ref{lem:geodesics}.

\begin{lemma}
\label{lem:XiinF}
The geodesic vector field $ \Xi $ belongs to $s^{-1}(\mathcal F)$.
\end{lemma}
\begin{proof}
By construction (see section \ref{sec:pull-back}), $ s^{-1}(\mathcal F)$ is generated by:
\begin{enumerate}
\item vertical vector fields, i.e., vector fields tangent to the fiber of $s=\pi \colon A \to M $,
\item  the horizontal lifts $H^\nabla (X) $ of vectors $X \in \mathcal F$.  
\end{enumerate}
As a consequence, a vector field which is locally of the form  $\sum_{i=1}^r  f_i H^\nabla\left(\rho(e_i)\right) $, where   $e_1, \dots, e_r $ is a local trivialization of $A$ and $ f_1, \dots, f_r$ smooth real-valued functions on $A$, belongs to $ s^{-1}(\mathcal F)$.  
By construction,  for any local trivialization  $e_1, \dots, e_r $ of $A$, and any $a \in A$:
$$ \Xi_{|_a} = \sum_{i=1}^r \langle e_i^*, a  \rangle \, H^\nabla(\rho(e_i)) $$ 
where $e_1^*, \dots, e_r^* $ is the dual trivialization of $A^*$. This proves the claim.
\end{proof}

\noindent
Here is an immediate consequence of Lemma \ref{lem:XiinF}.

\begin{lemma}
\label{lem:sIst}
There is a neighborhood of the zero section in $A$ where $ t^{-1}(\mathcal F) = s^{-1}(\mathcal F)$.
\end{lemma}
\begin{proof}
Lemma \ref{lem:XiinF} and Corollary \ref{coro:flowissymmetry} imply that the flow of $ \Phi^\Xi_t $ is a symmetry of $s^{-1}(\mathcal F) $. In particular, $\Phi^\Xi_1 $ is a symmetry\footnote{And even an inner-symmetry, but this is not important here.} of $s^{-1}(\mathcal F) $ on some neighborhood $\mathcal U_A $ where this time $1$-flow is well-defined:
$$   \left(\Phi_{1}^\Xi\right)^{-1} \left( s^{-1}(\mathcal F)\right) = s^{-1}(\mathcal F). $$
Since the target map consists precisely in composing $\Phi_{1}^\Xi$ with $s$ (see \eqref{eq:stPhi}), we have:
$$  t^{-1}(\mathcal F) =  (s \circ \Phi_{1}^\Xi)^{-1}(\mathcal F) =\left(\Phi_{1}^\Xi\right)^{-1} \left(  s^{-1}(\mathcal F)\right) =s^{-1}(\mathcal F).$$
This proves therefore that $s^{-1}(\mathcal F) = t^{-1}(\mathcal F) $.
\end{proof}

\begin{proof}[Proof of Proposition \ref{thm:Prop:FundamentalExampleAbstract}]
In view of Lemma \ref{lem:sIst},
we have to prove that any vector field $Z \in s^{-1}(\mathcal F)=t^{-1}(\mathcal F)$ decomposes as $X+Y$ with $X$ a vector field tangent to the fiber of $s$, i.e., the fiber of the canonical projection  $A \to M $, and the fiber of $t$. We are allowed to replace $\mathcal A $ by a smaller neighborhood of the zero section.

To start with, there is a natural anchored bundle structure over the singular foliation $ s^{-1}(\mathcal F)$ given by the vector bundle  $s^! A \oplus s^! A \to A $ equipped with the anchor\footnote{For $ \phi: M \to N$ a map and $ E\to N$ a vector bundle, we denote by $ \phi^! E$ the pull-back vector bundle $\phi^! E_m \simeq E_{\phi(m)}   $ and by $\phi^! e \in \Gamma(\phi^! E)$ the pull-back of a section $ e \in \Gamma(E)$.}
 $$  \rho_{A \oplus A} : (s^!a_1,s^!a_2) \mapsto  H(\rho(a_1))  + a_2^v$$
 where $H: \mathfrak X(M) \to \mathfrak X(A)$ is the Ehresmann connection on $ A$ associated to $ \nabla$, and where $ a^v$ is the vertical vector field\footnote{For very point $a$, there is an injection $A_{s(a)} \hookrightarrow T_a A$ due to the vector bundle structure that the fibers of $ A \to M$ possesses. This allows to see sections of $A$ as vector fields on $A$ valued on $ {\mathrm{ker}}(Ts)$.} on $A$ associated to $ a \in \Gamma(A)$.
 Since $ \Xi \in s^{-1}(\mathcal F)$, Proposition \ref{prop:crucial} allows lifting $\Phi_1^\Xi $ to an isomorphism of anchored bundles
  $$    \xymatrix{ s^! A \oplus s^! A\ar[r]^{\Psi_1} \ar[d]& \ar[d] s^! A \oplus s^! A\\ \mathcal A\ar[r]^{\Phi_1^\Xi}& A } $$
Since $ 0 \oplus s^! A $ has by construction an image through the anchor map $\rho_{A \oplus A}$ which lies in ${\mathrm{ker}}(Ts)$, and since $ \Phi_1^\Xi$ intertwines ${\mathrm{ker}}(Ts)$ and ${\mathrm{ker}}(Tt)$, 
the sub-vector bundle $ \Psi_* ( s^!A \oplus 0) \subset s^! A \oplus s^! A$ has an image through the anchor map that belongs to ${\mathrm{ker}}(Tt)$.
If we can show that for some open neighborhood $ \mathcal A' \subset \mathcal A$ of the zero section,  the vector bundles $\Psi_* ( 0 \oplus s^!A )$ and  $ 0 \oplus s^!A $ are in direct sum inside $ s^! A \oplus s^! A $, then the result follows, since for any vector field $Z$ in $s^{\mathcal F} $, it suffices to decompose a section  $\alpha \in \Gamma( s^! A \oplus s^! A) $ such that 
 $$ \rho_{A \oplus A} (\alpha)= Z  $$
 under the form $ \alpha =\alpha_1 + \alpha_2$
with $\alpha_1, \alpha_2  $ in $ 0 \oplus s^!A $ and  $\Psi_* ( 0 \oplus s^!A )$  respectively. By construction
 $ Z=X_1 +X_2 $ 
 with $ X_1 = \rho_{A \oplus A} (\alpha_1)$ and  $ X_2 = \rho_{A \oplus A} (\alpha_2)$ to have a decomposition of $Z $ as the sum of an element in $ \Gamma({\mathrm{ker}}(Ts))$
and an element in $ \Gamma({\mathrm{ker}}(Tt))$.

Let us show this point. The vector field $\Xi$ vanishes on the zero section, so that each point of the zero section is a fixed point of $\Phi_1^\Xi $. Also, it satisfies, on any point $m$  of the zero section, and for any $a \in \Gamma(A)$, the relation
 \begin{equation} \label{eq:inirdertoleift} [\Xi,a^v]_{|_m} = \rho(a_m) \end{equation}
 where the right-hand side, which belongs to $ T_m M$, is to be seen as an element in $ T_m A$ with the help of the zero section $T_m A \hookrightarrow T_m A $.
Now, for proving Proposition \ref{prop:crucial}, we  lifted $\Xi$ to a fiberwise linear vector field on $ s^! A \oplus s^! A $. Since $m$ is a fixed point, upon identifying linear vector fields with linear endomorphism of $ A_m \oplus A_m$, 
Equation \eqref{eq:inirdertoleift} implies that this lift can be chosen to be
given by the matrix:
 $$  \begin{pmatrix}
     0&  {\mathrm{id}}_{A_m  } \\ 0 &  0
 \end{pmatrix} .$$
 Its flow is therefore given by the matrix
  $$  \begin{pmatrix}
     {\mathrm{id}}_{A_m }&  {\mathrm{id}}_{A_m  } \\ 0 &  {\mathrm{id}}_{A_m } 
 \end{pmatrix} . $$
 Hence, the restriction of $ \Phi_1 (0 \oplus s^!A )$ to $M$ is the diagonal sub-bundle of $ A \oplus A$.
 Since the vector sub-bundle $ 0 \oplus A $ is in  direct sum with the diagonal sub-bundle of $ A \oplus A$, the desired property holds true in a neighborhood $ \mathcal A'$ of the zero section. This concludes the proof.
\end{proof}

\begin{remark}
\normalfont
When the anchored bundle is the trivial bundle $A = \mathbb R^r \times M \to M $ with anchor $$ \rho\left( (t_1, \dots, t_r),m \right) = \sum_{i=1}^r t_i X_i(m)$$ as in the proof of Proposition \ref{thm:prop:AnchoredBundle}, and the connection is the trivial connection:
 $$ \nabla_X (f_1, \dots, f_r) = \left( X[f_1], \dots,X[f_r]\right) \hbox{ for any $ f_1, \dots, f_r \in \mathcal C^\infty(M)$}$$
 then the  parallel path stating from $((t_1, \dots, t_r),m) \in A $ is the path:
 $$ t \mapsto   \left((t_1, \dots, t_r),\phi_t^{\sum_{i=1}t_i X_i}(m)\right)  $$
In this case,  the bisubmersion in Proposition \ref{thm:Prop:FundamentalExampleAbstract} is similar (but different) to the bisubmersion described in Proposition \ref{thm:prop:crucialExample}.
\end{remark}

\begin{remark}
    \normalfont
   The fundamental bisubmersion of Proposition \ref{thm:Prop:FundamentalExampleAbstract} could more be equivalently defined, in a more symmetric and therefore pleasant way, by using the maps $ s=\pi \circ  \Phi^\Xi_{-1/2}$ and $t = \pi \circ \Phi^\Xi_{+1/2} $ with $ \pi \colon A \to M$ the projection onto the base manifold. 
   Proposition \ref{thm:Prop:FundamentalExampleAbstract} still stands with this definition.
\end{remark}

 \subsubsection{Discussion on the notion of bisubmersion}

 \vspace{1mm}
 \noindent
As already mentioned, the word ``bisubmersion'' used alone does not make sense if used alone: only the expression ``bisubmersion over the singular foliation $\mathcal F$'' makes sense. To clarify this point,  we introduce the following notion (maybe not interesting by itself, but  practical for pedagogical purposes). 
 
 4\begin{definition}
We call \textit{twin-submersion} the data of 
 \begin{enumerate}
     \item two manifolds $X,M$,
     \item two surjective submersions $s,t \colon X \to M $.
 \end{enumerate}
 A \textit{unit map for twin-submersions} is a smooth map $\epsilon\colon M \to X $ which is a section of both $s$ and $t$.
 \end{definition}

  \vspace{1mm}
 \noindent
Here is a natural question: 

\begin{question}
Given a  twin-submersion $(X,M,s,t) $, when is it a bisubmersion for some singular foliation $\mathcal F $?
\end{question}

  \vspace{1mm}
 \noindent
 We answer this question through the following exercises. 

\begin{exo}
   Show that a twin-submersion can not be a bisubmersion for two different singular foliations on $M$.

 \vspace{1mm}
 \noindent
   {\emph{Hint:}} Show that $\mathcal F \neq \mathcal F' $ implies $ s^{-1}(\mathcal F) \neq s^{-1}(\mathcal F') $.
\end{exo}

\begin{exo}
Let $(X,M,s,t)$ be a  twin-submersion. 
\begin{enumerate}
\item 
Show that the following two conditions are equivalent. 
\begin{enumerate}
    \item[(i)] $\mathcal G := \Gamma(\ker(Ts)) + \Gamma(\ker(Tt)) $ is stable under Lie bracket. 
    \item[(ii)]  $ \left[ \Gamma(\ker(Ts)) , \Gamma(\ker(Tt)) \right] \subset   \Gamma(\ker(Ts)) + \Gamma(\ker(Tt))$
\end{enumerate}
Show that the previous equivalent two conditions are in turn equivalent to: 
\begin{enumerate}

    \item[(iii)] $\ker(Tt) $ is generated by $s$-projectable\footnote{I.e. $s$-related to a vector field on $M$, see Definition \ref{def:related}.} vector fields and $\ker(Ts) $ is generated by $t$-projectable\footnote{I.e. $t$-related to a vector field on $M$, see Definition \ref{def:related}.} vector fields.
\end{enumerate}
{\emph{Hint:} Show that $\mathcal G $ is generated by vector fields $X$ which are $s$-projectable. Then write $X=X^s+X^t$ with $X^s$, $ X^t$ in the kernels of $Ts$ and $Tt$, respectively. Show that $X^t $ is $s$-projectable and belongs to the kernel of $Tt$.}

\item 
We now assume that the fibers of $s$ and $t$ are connected. 
Show that if one of the equivalent conditions above hold, then there exists singular foliations $\mathcal F_s $ and $\mathcal F_t $ on $M$ such that
$$ \mathcal G = s^{-1}(\mathcal F_s)   =   t^{-1}(\mathcal F_t) $$

\item Let us now assume that there exists a unit map $\epsilon \colon M \hookrightarrow X $, i.e a left inverse of both $s$ and $t$. Show that $ \mathcal F_s= \mathcal F_t $.
\item Show that if a twin-submersion admits a unit map and has connected $s$- and  $t$-fibers, then there exists a singular foliation with respect to which it is a bisubmersion if and only if $$ \left[ \Gamma(\ker(Ts)) , \Gamma(\ker(Tt)) \right] \subset   \Gamma(\ker(Ts)) + \Gamma(\ker(Tt)).$$
\end{enumerate}

\end{exo}

\subsubsection{More examples of bisubmersions}

\noindent
By composing the source or the anchor with a symmetry of a bisubmersion over $(M,\mathcal F) $, one still obtains a bisubmersion over $ \mathcal F$. 
The following definition therefore makes sense.

\vspace{.5cm}

\begin{definitions}{Composing bisubmersions with symmetries}{compSym}
Let $ M \stackrel{s}{\leftarrow} W \stackrel{t}{\rightarrow} M$ be a bisubmersion over $\mathcal F $, and  $\phi \colon M \to M$ a symmetry of $\mathcal F$. 
Then $$  M \stackrel{\phi \circ s}{\longleftarrow} W \stackrel{t}{\longrightarrow} M   \hbox{ and } M \stackrel{ s}{\longleftarrow} W\stackrel{{\phi^{-1} \circ t}}{\longrightarrow} M   $$ 
 are bisubmersions of $ \mathcal F$ again. We call them the \emph{right and left composition by the symmetry $ \phi$}, respectively. 	
\end{definitions}

\noindent
An analogous construction can be done for local bisubmersions over $ \mathcal F$.
In fact, any bisubmersion over $ \mathcal F$  is of the type constructed in Proposition  \ref{thm:Prop:FundamentalExampleAbstract}, up to composition by a symmetry, at least locally. This is the topic of the next Exercise.

\begin{exo}
	\label{exo:frombisubToAnchored}
	 Let  $ M \stackrel{s}{\leftarrow} W \stackrel{t}{\rightarrow} M$  be a bisubmersion of $\mathcal F $. Choose $\Sigma $ a bisection. Let $ N_{\Sigma} $ be the normal bundle of $\Sigma $ in $W$ as in Lemma \ref{lem:BisectionAnchoredBundle}.

 Show that there is a neighborhood of any bisubmersion $\Sigma  $ in $W$ isomorphic, as a bisubmersion, to the left composition by $\underline{\Sigma} $ of a fundamental bisubmersion (see Proposition \ref{thm:Prop:FundamentalExampleAbstract}) of the anchored bundle $(A_{\Sigma}, \rho)$. (\emph{Hint}: use Theorem \ref{thm:equiv} below, together with Exercise \ref{exo:samedimension}).

\end{exo}

\subsection{Left and right ``actions'' of anchored bundles} 

\noindent
Let $\mathcal F $ be a singular foliation on $M$. So far, we have seen two classes of objects ``over $\mathcal F$''.
\begin{enumerate}
    \item  bisubmersions $M \stackrel{s}{\leftarrow} W \stackrel{t}{\rightarrow} M$  of $\mathcal F $. 
\item anchored bundles $(A,\rho) $ over $\mathcal F $. 
\end{enumerate}
Let us assume that we are given both. What is the relation between them? 
Proposition \ref{thm:Prop:FundamentalExampleAbstract} was already a partial  to that question, but there is more.
If $W$ was a  Lie groupoid with Lie algebroid $(A,\rho)$, then $W$ would be a bisubmersion and $ (A,\rho)$ would be an anchored bundle for the same singular foliation $\mathcal F =\rho(\Gamma(A))$ (see exercise \ref{exo:groupasbosubmersions}). In that case, left and right Lie algebroids actions of $ A$ on $W $ can be defined with the help of the notion of left and right invariant vector fields. 
Those can also be defined through left and right Lie groupoid action.
Despite the lack of product on a bisubmersion, there is a very similar construction. Or course, it is not an ``action'' \emph{stricto sensu}, but the following object can be defined, following \cite{RubenSymetries}.

\vspace{0.5cm}

\begin{definitions}{Right and Left action}{actions}
Let $ (M,\mathcal F)$ be a foliated manifold. Let $M \stackrel{s}{\leftarrow} W \stackrel{t}{\rightarrow} M$  and $(A,\rho) $ be a bisubmersion and an anchored bundle over $(M,\mathcal F)$, respectively.

\vspace{1mm}
\noindent
We call, respectively,  right and left actions of $(A,\rho)$ on the bisubmersion $ W$ vector bundle morphisms\footnote{$s^!A,t^!A$ are the pull-back vector bundles of the vector bundle $A \to M$ to $W$ through $s$ and $t$, respectively.}:
 $$  L : s^! A \longrightarrow \ker (Tt) \hbox{ and }  R : t^! A \longrightarrow \ker (Ts)  $$ 
 making the following diagrams commutative:
 $$ \xymatrix{ s^! A \ar[r]^{L} \ar[d]^{\rho} & \ker (Tt) \ar[dl]_{Ts}\\ TM &  } \begin{array}{c} \\ \\ \\\ \\ \hbox{ and }\\  \end{array} \xymatrix{ \ker (Ts) \ar[dr]_{Tt}& t^! A \ar[l]_{R} \ar[d]^{\rho} \\ &TM   }  $$

 \end{definitions}

\vspace{.5cm}
\noindent
The notion is interesting for many reasons. To start with, it really exists:
\vspace{.5cm}

\begin{propositions}{Right and left actions exist}{prop:rightleftexist}
Let $ (M,\mathcal F)$ be a foliated manifold. Let $M \stackrel{s}{\leftarrow} W \stackrel{t}{\rightarrow} M$  of $\mathcal F $ and $(A,\rho) $ be a bisubmersion and an anchored bundle over $(M,\mathcal F)$ respectively. 
Right and left-actions  of $(A,\rho)$ on $ W$ exist.
\end{propositions}
\begin{proof}
We will prove it for the right action $R$.
Let $x \in W$, and let $ e_1, \dots, e_r$ be a local trivialization of $A$ near $t(x) $. Then $\rho(e_1), \dots, \rho(e_r) $ are generators of $\mathcal F $. 
In view of Lemma \ref{lem:localstructure}, there exists a neighborhood $\mathcal U $ of $x $ and vector fields $X_1, \dots, X_r $ in $\Gamma_\mathcal U(\ker(Ts)) $ such that $X_i $ is $t$-related to $\rho(e_i) $.
We define a vector bundle morphism $ R_\mathcal U$ from $t^! A_{|_\mathcal U} $ to $ \ker(Ts)_{|_\mathcal U} $ that satisfies the desired property in a neighborhood of $x $ by mapping the sections\footnote{Recall that for $\phi: N \to M $ a map, $ E \to M$ a vector bundle, and $e$ a section of $ E$, we denote by $ \phi^! e$ the section of $ \pi^! E = E \times_M N$ whose value at $n \in M$ is $ (e_{|_{\phi(n)}} , n)$.} $ t^! e_1, \dots, t^! e_r$ to $X_1,\dots, X_r$. 
The proposition then follows from the fact that given a partition of unity $(\mathcal U_i,\chi_i)_{i \in I} $ of $ W$ such that  
such a  vector bundle morphism  $ R_{\mathcal U_i}$ as in the statement exists on $ \mathcal U_i$, the linear combination $$ R := \sum_{i \in I} \chi_i R_{\mathcal U_i} $$
is a  vector bundle morphism\footnote{(now defined on the whole manifold $W$)} from $ t^!A$  to ${\mathrm{ker}}(Ts) $ that satisfies all desired properties.
\end{proof}

\begin{remark}
When the bisubmersion  $W$ is a Lie groupoid and $A$ is its Lie algebroid (see exercise \ref{exo:groupasbosubmersions}), then the usual right and left actions of  the Lie algebroid are instances of the previous left and right actions. 
They are not the unique ones in general.
\end{remark} 

\vspace{1mm}
\noindent
For any left or right actions and any $a \in \Gamma(A)$, a vector field on $ W$ is defined by $ L(s^!a)  $ or $ R(t^!a)$ respectively. We denote these vector fields by $\overrightarrow{a} $ and $\overleftarrow{a}  $ respectively.

The maps:
 $$ \begin{array}{rcl} \Gamma(A) & \to & \mathfrak X(W) \\ a & \mapsto & \overrightarrow{a} \\ a & \mapsto  & \overleftarrow{a} \end{array}  $$
satisfy for any $ a \in \Gamma(A)$ the following conditions by construction:
\begin{enumerate}
\item the vector field $\overrightarrow{a}$
(resp. $\overleftarrow{a}$) $s$-projects (resp. $t$-projects) to $\rho(a) \in \mathfrak X(M) $
    \item the vector field $\overrightarrow{a}$
(resp. $\overleftarrow{a}$) is tangent to the fibers of $t$ (resp. $s$).
\item For any $F \in \mathcal C^\infty(M)$, we have $ \overrightarrow{Fa}=t^* F \, \overrightarrow{a} $ and $ \overleftarrow{Fa}=s^* F \, \overleftarrow{a} $. 
\item For any $a,b \in \Gamma(A) $, the vector field $ [\overrightarrow{a},\overleftarrow{b}] $ is bi-vertical.
 \end{enumerate}

\begin{remark}
\normalfont
\label{rmk:rightaction-modulobivertical}
For a different choice for the right action, the corresponding vector fields $\overrightarrow{a}$ will differ by bivertical vector fields. Hence, the class of $\overleftarrow{a}$, $\overrightarrow{a}$ modulo bivertical vector fields does not depend on the choice of a right action. The same holds for the left action. 
\end{remark}
 
\begin{lemma}
\label{lem:sometimesIso}
Let $ W $ be a bisubmersion over $(M,\mathcal F) $. Any $w \in W$ admits a neighborhood $ \mathcal V$ such that there exists \emph{(i)} anchored bundles $(A_R,\rho_R) $ and $(A_L,\rho_L) $ over $ \mathcal F$, defined on $ s(\mathcal V)$ and $ t(\mathcal V)$, and \emph{(ii)} a choice of right and left actions $R$ and $L$ near $w$ such that both $R$ and $L$
$$  L : s^! A_L \stackrel{\simeq}{\longrightarrow} \ker (Tt)|_{\mathcal V} \hbox{ and }  R : t^! A_R \stackrel{\simeq}{\longrightarrow} \ker (Ts)|_{\mathcal V}  $$
are vector bundle isomorphisms over ${\mathcal V}$.
Moreover, if $ s(w)=t(w)$, then one can assume $ (A_R,\rho_R)$ and $ (A_L,\rho_L)$ coincide near $s(w)=t(w)$. 
\end{lemma}
\begin{proof}
Let $\Sigma $ be a local bisection through $w $ and $ N_\Sigma$ be the normal bundle of $ \Sigma$ in $W $.
In view of Lemma \ref{lem:BisectionAnchoredBundle}, this normal bundle can be seen as an anchored bundle for $M$ on $ s(\Sigma)$ or $ t(\Sigma)$: it suffices to identify the henceforth obtained vector bundle with $ {\mathrm{ker}}(Tt)$ or to 
 $  {\mathrm{ker}}(Ts)$  through $s$ and $t$ and to equip it with $ Ts$ and $Tt$, respectively. The proof of Proposition \ref{thm:prop:rightleftexist} can then be adapted so that $R$ and $L$ are the identity maps on  $ {\mathrm{ker}}(Ts)_{|_\Sigma}$ or to 
 $  {\mathrm{ker}}(Tt)_{|_\Sigma}$. Hence, both $R$ and $L$ have to be invertible on a neighborhood of $\Sigma $ in $W$. This completes the proof. 
\end{proof}

\vspace{1mm}
\noindent
Now, assume that the anchored bundle $ (A,\rho)$ is equipped with a bracket making it an almost Lie algebroid: such a bracket exists by Proposition \ref{prop:Almost-existence}.  
It is interesting to notice that for  every choice of such an almost Lie algebroid bracket on $(A,\rho)$, the vector fields that measures the default of the left and right actions to preserve the brackets, i.e., the vector fields  
 $$ \overleftarrow{[a,b]_A}- [\overleftarrow{a},\overleftarrow{b}] \hbox{ and }  \overrightarrow{[a,b]_A}- [\overrightarrow{a},\overrightarrow{b}] \hbox{ with $a,b \in \Gamma(A)$},$$ are bi-vertical vector fields\footnote{We also saw that the vector field $[\overrightarrow{a},\overleftarrow{b}] $ is always bi-vertical.}.
 Of course, it does not make sense to hope for a Lie algebra morphism, since $\Gamma(A)$ is not a Lie algebra. However, the following exercise describes a natural and canonical Lie algebra morphism.

 \begin{exo}
Let $ W $ be a bisubmersion over $(M,\mathcal F) $.
Show that:
\begin{enumerate}
\item $t$-projectable\footnote{I.e. $ t$-related (see Definition \ref{def:related}) to a vector field in $\mathfrak X(M) $.} sections of $ \ker(Ts)$ form a Lie algebra, denoted by $ \Gamma(\ker(Ts))^t$.
\item Bivertical vector fields form a Lie ideal (denoted ${\mathrm{Bivert}} $) of the previous Lie algebra.

\end{enumerate}
Let $ (A,\rho) $ be an anchored bundle over $\mathcal F $.  Let ${\mathrm{Ker}}(\rho) \subset \Gamma(A)$ be the subspace of all $a \in \Gamma(A) $ such that $\rho(a)=0 $.
\begin{enumerate}
\item[3] Show that ${\mathrm{Ker}}(\rho)  \subset \Gamma(A) $ is a left and right ideal for the bracket of $\Gamma(A) $ and that the quotient space $\Gamma(A) /{\mathrm{Ker}}(\rho)$ is a Lie algebra isomorphic to $ \mathcal F$.
\item Show that for any choice of a right-action, $\overrightarrow{a} $ is bivertical if $a \in {\mathrm{Ker}}(\rho) \subset \Gamma(A) $.
\item Show that the induced map:
 $$ \mathcal F \simeq \frac{\Gamma(A)}{{\mathrm{Ker}}(\rho) } \longrightarrow \frac{\Gamma(\ker(Ts))^t}{{\mathrm{Bivert}} }  $$
 is a Lie algebra morphism. 
 \item Show that this Lie algebra morphism does not depend on the choice of the right action $R$ (see Remark \ref{rmk:rightaction-modulobivertical}).
 \item Show that this Lie algebra morphism is in fact an isomorphism.
\end{enumerate}

 \end{exo}

\vspace{0.3cm}
\noindent
We now use the  notion of left and right action  to show the following interesting point about $ {\mathrm{Out}}_W$. 
Let $(M,\mathcal F)$ be a foliated manifold, $W$ and $W'$ two bisubmersions over $ \mathcal F$. We intend in Exercise \ref{exo:Outiff} below to show that if $x \in W$ and $x' \in W'$ are related by $ {\mathrm{Out}}_W(x)={\mathrm{Out}}_W(x')  $, then there exists local bisections $\Sigma $ and $\Sigma' $ through $ x$ and $x'$ respectively such that $\underline{\Sigma} = \underline{\Sigma'} $:

\begin{exo}
\label{exo:Outiff}
Let $\Sigma_0 $ and $ \Sigma_0'$ be bisections of $W$ and $W'$ through $x$ and $x'$. By assumption $\Phi := \underline{\Sigma_0} \circ (\underline{\Sigma_0'})^{-1}$ is an inner symmetry admitting $ m=t(x)=t'(x')$ as a very fixed point, i.e., there exists a smooth time-dependent vector field $(X_\epsilon)_{\epsilon \in [0,1]} $ in $\mathcal F $  whose time $\epsilon$-flow $ \Phi_\epsilon$ is defined for all $ \epsilon\in [0,1]$ and whose time $1$-flow is $\Phi$ (at least in a neighborhood of $X$). Moreover, by Exercise \ref{exo:veryfixed}, 
$ X_\epsilon = \rho(a_\epsilon)$ where $ (a_\epsilon)_{\epsilon \in [0,1]}$ is a smooth\footnote{I.e., smoothly depending on the parameter $ \epsilon \in I$.} time dependent section of $ \Gamma(A)$  that we can assume to satisfy $ a_\epsilon(m)=0$ for every $ \epsilon \in [0,1]$.
\begin{enumerate}
\item Choose a right action on $W$. Let $\Psi_\epsilon $ be the flow of the time dependent vector field $( \overleftarrow{a_\epsilon} )_{\epsilon \in [0,1]} $. Show  that there exists a neighborhood $\Sigma $ of $x$  in $\Sigma_0 $ such that $\Sigma_\epsilon := \Psi_\epsilon(\Sigma) $ is well-defined for every $ \epsilon \in [0,1]$, and is a bisection of $ W$ for all $\epsilon \in [0,1]$.
\item Show that the bisection  $\Sigma_\epsilon$ contains $x$ for every $ \epsilon \in [0,1]$.
\item Show that $ \underline{\Sigma_\epsilon} = \Phi_\epsilon \circ \underline{\Sigma}$ for every $ \epsilon \in [0,1]$.
\item Conclude that the bisections $ \Sigma_1:=\Psi_1 (\Sigma)$ and $ \Sigma_0'$ induce diffeomorphisms that coincide in a neighborhood of $ s(x)=s'(x')$.
\item Conclude.
\end{enumerate} 
\end{exo}

 Last, right and left actions are extremely practical to check that bisubmersions satisfy the following property.

 \begin{definition}
 \label{def:smallinner}
 Let $I$ be an interval containing zero.
  Let $ (M,\mathcal F)$ be a foliated manifold, and $ \mathcal U \subset M$ an open subset.
     We say that a bisubmersion $ M \stackrel{s}{\leftarrow} W \stackrel{t}{\rightarrow} M $ \emph{represents all small inner symmetries on $ \mathcal U$} if, for every smooth time-dependent vector field $ (X_\epsilon)_{\epsilon \in I}$ whose time $t$-flow is defined on $ \mathcal U$ for all $\epsilon \in I $, there exists for all $ m \in M$ and all $\epsilon$ small enough a local bisection $ \Sigma_\epsilon$ of $W$ such that $ \underline{\Sigma_t}$ and the flow $ \phi_\epsilon^{X_\bullet}$ of   $ (X_\epsilon)_{\epsilon \in I}$ coincide in a neighborhood of $m$.
 \end{definition}

Here is a simple condition to have this condition satisfied, and we will use right actions to prove it.

 \begin{proposition}
 \label{prop:smallinner}
     If a  bisubmersion  $ M \stackrel{s}{\leftarrow} W \stackrel{t}{\rightarrow} M $  over $\mathcal F $ admits a unit map $ \epsilon$, defined on $ \mathcal U \subset M$, then it represents all small inner symmetries on $ \mathcal U$.
 \end{proposition}
 \begin{proof}
 Let $ (X_\epsilon)_{\epsilon \in I}$ be as in Definition \ref{def:smallinner}.
 Let $m \in \mathcal U$, and let $ (A,\rho)$ be an anchored bundle over $\mathcal F $ near $m$. Let $(a_\epsilon)_{\epsilon \in A}$ be a smooth time-dependent section of $A$ such that $ \rho(a_\epsilon)=X_\epsilon$. 
 Choose a right action, and consider the smooth time-dependent vector field $\overrightarrow{a_\epsilon} $ on $W$. Since this vector field tangent to the fibers of $s$ and is $t$-related to $ X_\epsilon$, its flow $ \Phi_t^{a_\bullet}$ and the flow $ \phi_t^{X_\bullet}$ are related by $ t \circ \Phi_t^{a_\bullet} \circ \epsilon = \phi_t^{X_\bullet}$ for every $\epsilon$ for which it is defined. The submanifold $ \Sigma_\epsilon =  \Phi_\epsilon^{a_\bullet} \circ \epsilon (\mathcal U) $ is a bisubmersion that satisfies the desired property.
 \end{proof}.

\begin{example}
\label{ex:smallinner}
The bisubmersion of Proposition \ref{thm:prop:crucialExample} represents all small inner symmetries on the $ \mathcal U$ where it is defined.
    The fundamental bisubmersions in Proposition \ref{thm:Prop:FundamentalExampleAbstract}, associated to an anchored bundle $ (A,\rho)$ defined on $ \mathcal U$, represents all small inner symmetries on $ \mathcal U$.
\end{example}
 
\subsection{Products and inverse of bisubmersions}
\label{sec:productinverse}

\noindent
Bisubmersions for a given singular foliation behave like Lie groupoids, but so far there is still no inverse and no product.    
The following notions present an analogy of those\footnote{See Proposition 2.4 in \cite{AS}.}:

\vspace{.5cm}

\begin{definitions}{Product and Inverse}{ProdInv}
Let $\mathcal F$ be a singular foliation on a manifold $M$. 
\begin{enumerate}
    \item The \emph{inverse} $W^{-1}$ of a bisubmersion $ M \stackrel{s}{\leftarrow} W \stackrel{t}{\rightarrow} M $ over $ \mathcal F$ is simply the bisubmersion 
    $ M \stackrel{t}{\leftarrow} W \stackrel{s}{\rightarrow} M $.
    \item The \emph{composition} $W*W'$ of two bisubmersions $ M \stackrel{s}{\leftarrow} W \stackrel{t}{\rightarrow} M $ and $ M \stackrel{s'}{\leftarrow} W' \stackrel{t'}{\rightarrow} M $ over $ \mathcal F$ is the fibered product 
     $$W \times_{s,M,t'} W' = \{(x,x') \in W \times W' \hbox{ s.t. } t(x)=s'(x')\} $$ equipped with the source $ (x,x') \mapsto s(x) $ and target $(x,x') \mapsto t'(x') $.
\end{enumerate}
\end{definitions}
\begin{exo}
We leave it as an exercise to check that the product of
bisubmersions for $\mathcal F $, defined as above, is a bisubmersion  for $\mathcal F $ again.
\end{exo}

We finish this discussion  with some explanation of the names ``inverse'' and ``compositions'' of bisubmersions.

\begin{proposition}
\label{prop:productAndInverse}
Let $(M,\mathcal F) $ be a foliated manifold.
\begin{enumerate}
\item 
Consider a bisubmersion $ M \stackrel{s}{\leftarrow} W \stackrel{t}{\rightarrow} M $ over a singular foliation $ \mathcal F$.  Then the outer-germ maps of $W$ and $ W^{-1}$ are related as follows
 $$  {\mathrm{Out}}_{W^{-1}}(x) 
 = 
 \left({\mathrm{Out}}_{W}(x)\right)^{-1} \hbox{ for all $x\in W$}.  $$
\item Consider a second bisubmersion $ M \stackrel{s}{\leftarrow} W' \stackrel{t}{\rightarrow} M $ over a singular foliation $ \mathcal F$, the outer-germ map of $ W * W'$, $W$, and $W'$ are related by:
$${\mathrm{Out}}_{W * W'}(x,x') = {\mathrm{Out}}_{W'}(x') \circ {\mathrm{Out}}_{W}(x)  $$
for all $ x \in W, x' \in W'$ such that  $t(x)=s(x') $.
\end{enumerate}
\end{proposition}
\begin{proof}
Let us prove item $1$, for every local bisection $\Sigma \subset W$ through $x$, $ \Sigma$ is also a bisection for $W^{-1}$. But the corresponding isomorphism of $\mathcal F $ being the inverse of the source, restricted to $\Sigma $, composed with the target, they give isomorphisms which are inverse one to the over. 

Let us prove Item $2$. For two local bisections $\Sigma, \Sigma' $ through $x$
and $x'$, the fibered product $ \Sigma \times_{s,M,t'} \Sigma'$ is a bisection of the product bisubmersion $W*W'$. Moreover,  we have the following property:
 \begin{equation} \label{eq:inversionbisub} \underline{\Sigma \times_{s,M,t'}\Sigma'} \, = \, \underline{\Sigma'} \circ \underline{\Sigma} \end{equation}
This completes the proof of the claim.
\end{proof}

\subsection{Equivalence of bisubmersions (and their compositions)}

There is\footnote{We follow here Section 2 in \cite{AS}.} a ``Morita equivalence-like'' equivalence relation on the set of all bisubmersions over $\mathcal F$. 
Its definition is very natural for the reader used to Morita equivalence of Lie groupoid. There is also a natural notion of morphism, that, surprisingly, will be (more or less) the same as an equivalence.
  
  \vspace{.5cm}
  
\begin{definitions}{Equivalence of bisubmersions}{EquvAndMorph}
Consider two bisubmersions $ M \stackrel{s}{\leftarrow} W \stackrel{t}{\rightarrow} M $ and $ M \stackrel{s'}{\leftarrow} W' \stackrel{t'}{\rightarrow} M $ over $ \mathcal F$.
\begin{enumerate}
    \item A \emph{morphism} from the first one to the second one is a map $ W \to W'$ making the following diagram commutative:
    $$\xymatrix{M \ar@{=}[d] &W \ar[d] \ar[l]_{s}\ar[r]^{t}& M\ar@{=}[d]\\ M & W' \ar[l]^{s'}\ar[r]_{t'}& M} $$
 \item  An \emph{equivalence between them} is a third bisubmersion $ M \stackrel{s''}{\leftarrow} P \stackrel{t'' \, }{\rightarrow} M $
of $ \mathcal F$ equipped with two surjective submersions $P \to W $ and $P \to W'$ making the following diagram commutative:
 \begin{equation}
     \label{eq:def:bisubequiv}
 \xymatrix{  &W \ar[dl]_{s} \ar[dr]^{t} & \\ M&P\ar@{->>}[u]_{\stackrel{}{\pi'}} \ar@{->>}[d]_{ \stackrel{\pi}{}} \ar[l]_{ \, \, \, \, s''} \ar[r]^{t'' \, \, \,}  &M\\ &W' \ar[ul]^{s'} \ar[ur]_{t'}   & } 
 \end{equation} 
\end{enumerate}
\end{definitions}  

\vspace{0.5cm}
 
\begin{exo}
Show that the so-called equivalence defined above is indeed an equivalence relation on bisubmersions over $\mathcal F$.
\end{exo} 

\begin{exo}
Let $ \Gamma \toto M $ be a Lie groupoid with Lie algebroid $ (A,\rho)$.
For any connection $\nabla $ on $A$, its exponential map ${\mathrm{exp}} \colon \mathcal U_A \to \Gamma $ is a (iso)morphism of bisubmersions from the fundamental bisubmersion $\mathcal U_A \subset A $ to an open neighborhood of the unit manifold $ M$ in $\Gamma $. 
The following exercise details this point (and recalls the definition of the exponential map in the context of Lie algebroids). 
\begin{enumerate} 
\item For every $m \in M$ and every $ a \in A_m$, there exists a unique curve $\epsilon \mapsto \gamma(\epsilon) $ such that for every $ \epsilon $ for which it is defined:
 $$ s \circ \gamma(\epsilon) = m  \hbox{ and } a(\epsilon) = \gamma^{-1}(\epsilon) \frac{d \gamma(\epsilon)}{d \epsilon}  $$
 is the geodesic of $A$ starting from $a$. We define ${\mathrm{exp}}(a)  \in \Gamma $ to be $ \gamma(1) $. 
\item  Show that $t \circ \Phi_\epsilon^\Xi (a)= t (\gamma(\epsilon)) $. 
\item Conclude. 
\end{enumerate} 
\end{exo}

Here is an important and surprising theorem.

\vspace{.5cm}

\begin{theorems}{Equivalence of bisubmersions}{thm:equiv}\label{thm:equiv}
Consider two bisubmersions $ M \stackrel{s}{\leftarrow} W \stackrel{t}{\rightarrow} M $ and $ M \stackrel{s'}{\leftarrow} W' \stackrel{t'}{\rightarrow} M $ for $ \mathcal F$. 
The following statements are equivalent:
\begin{enumerate}
    \item[(i)] Both bisubmersions are equivalent. 
    \item[(ii)] The following two conditions hold:
    \begin{enumerate}
        \item any $x \in W$ admits a neighborhood $\mathcal U $ on which a morphism $ \mathcal U \to W'$ exists,
        \item and any $x' \in W'$ admits a neighborhood $\mathcal U' $ on which a morphism $ \mathcal U' \to W$ exists.
    \end{enumerate}
    \item[(iii)] Both bisubmersions induce the same outer-germs\footnote{Equivalently: for every $x \in W$ there exists $x' \in W'$ such that $ \mathrm{Out}_{W} (x)=\mathrm{Out}_{W'} (x') $ and conversely for every $x' \in W'$ there exists $x\in W$ such that $ \mathrm{Out}_{W'} (x')=\mathrm{Out}_{W} (x)$.}, i.e.
     $$ \mathrm{Out}_{W} (W)= \mathrm{Out}_{W'} (W').$$
\end{enumerate}
\end{theorems}

\begin{proof}

    Let us prove $(i)\Rightarrow (ii)$. Consider an equivalence as in Equation \eqref{eq:def:bisubequiv}.
    For every $x \in W$, a local section $ \sigma $ of $\pi\colon W'' \to W $ can be defined in an open neighborhood of $x$.   The composition of $ \sigma$ with the projection $\pi'\colon P \to W' $ is a morphism as in $(ii)$ item $(a)$.
    Similarly, for every $x \in W'$ the composition of the projection $\pi\colon P \to W $ with a local section of $\pi'\colon P \to W'$ is a morphism as in $(ii)$ item $(b)$. This shows that $(i)$ implies $(ii)$. 

   Let us show $(ii)\Rightarrow (iii)$. The image of any local bisection $\Sigma $ of $W$ through a morphism $ \phi \colon \mathcal U \to W' $ (with  $\mathcal U $ an open subset of $W$) is a bisection of $ W'$, and 
     $ \underline{\phi(\Sigma)} = \underline{\Sigma}  $.
     In particular, for any $x,x'$ as in item $(ii)$, we have ${\mathrm{Out}}_W(x)= {\mathrm{Out}}_{W'}(x') $.

    Let us prove that $(iii)$ implies $(i)$.
    First, let us use the convenient notations  $ \mathcal G$ and $\mathcal G' $ for the bi-pull-back singular foliations on  the bisubmersions $W$ and $W' $, respectively.
    
    Let $x \in W$ and $x' $ be such that ${\mathrm{Out}}_W(x)= {\mathrm{Out}}_{W'}(x') $. By exercise \ref{exo:Outiff}, there exists local bisections $ \Sigma$ and $\Sigma' $ through $x$ and $x'$ such that $ \underline{\Sigma}=\underline{\Sigma'} $. In particular, the diffeomorphism $$
    \begin{array}{rrcl}\phi & \Sigma &\to& \Sigma' \\ & x& \mapsto &  \left(s'_{|_{\Sigma'}}\right)^{-1}\circ s_{|_{\Sigma}}\colon\Sigma (x) \end{array} $$ 
makes the following diagram a commutative diagram:
    $$ \xymatrix{ & M& \\  \Sigma \ar[rr]^{\phi} \ar[ru]^{s} \ar[rd]_{t}& & \Sigma' \ar[lu]_{s'} \ar[ld]^{t'} \\ & M &  }  .$$

     Let $\{\underline{X}_i, i= 1,\ldots, r\}$ be a set of generators of $\mathcal F$ in a neighborhood of $ s(x)$ and let  $\{\underline{Y}_i, i= 1,\ldots, r\}$ be a set of generators of $\mathcal F$ in a neighborhood of $ t(x)$. Let  
      $$(X;Y;Z)=(X_1, \dots, X_{r_s};Y_1,\dots, Y_{r_t}; Z_1, \dots,Z_k) $$
      and
      $$(X';Y';Z')=(X_1', \dots, X_{r_s}';Y_1',\dots, Y_{r_t}'; Z_1', \dots,Z_{k'}') $$
      be generators of $ \mathcal G$ and $ \mathcal G'$ as in Lemma \ref{lem:localstructure}.
     Without any loss of generality, one can assume $k=k'$: it suffices for instance to add $k'-k $ times the vector field $0$ if $ k' >k$.
     Consider the manifold $$ P_x'' :=\mathrm{Gr}_\phi(\Sigma) \times \mathbb R^{r_s} \times \mathbb R^{r_t} \times \mathbb R^k ,$$
     where $\mathrm{Gr}_\phi(\Sigma) \subset \Sigma \times \Sigma'$ is the graph of $ \phi$.

     From now on, we use for any $n$-tuple of vector fields $ \xi_\bullet = (\xi_1, \dots, \xi_n)$  the notation $$\exp{(\lambda_\bullet \xi_\bullet )}:= \Phi^{\xi_1}_{\lambda_1}\circ\cdots\circ \Phi^{\xi_n}_{\lambda_n}$$ 
 for any    $\lambda_\bullet =(\lambda_1,\cdots, \lambda_n)\in \mathbb{R}^n$ such that the flows make sense. 
   There is an open subset of $  P_x'' $ containing $ \mathrm{Gr}_\phi(\Sigma) \times (0, \dots, 0) \times (0, \dots, 0) \times (0, \dots, 0) $, that we still denote by $  P_x'' $ with a slight abuse of notations, on which the following two maps are well-defined:
     $$ \begin{array}{rrcl}
        \pi \colon & P_x''& \to  &  W \\ 
      & \left((\sigma,\sigma'),(\lambda_1, \dots,\lambda_{r_s}), (\mu_1, \dots, \mu_{{r}_t}), (\nu_1, \dots, \nu_k)\right)  &\mapsto& \exp{(\lambda_\bullet X_\bullet )}\circ \exp{(\mu_\bullet Y_\bullet )} \circ \exp{(\nu_\bullet Z_\bullet)}(\sigma)\\ & & & \\ 
      \pi'\colon & P_x'' & \to & W'   \\  
  &\left((\sigma,\sigma'),(\lambda_1, \dots,\lambda_{r_s}), (\mu_1, \dots, \mu_{{r}_t}), (\nu_1, \dots, \nu_k)\right)   &\mapsto &\exp{(\lambda_\bullet X_\bullet')}\circ \exp{(\mu_\bullet Y_\bullet')} \circ \exp{(\nu_\bullet Z_\bullet')}(\sigma')  \end{array}
  $$
     Above $ (\sigma,\sigma')$ is an element of the graph of $ \phi$, i.e., $ \sigma'=\phi(\sigma)$.
     
Consider the point inside $ P_x''$  given by
$$ O_\sigma := ((\sigma,\sigma'),\underbrace{0,\dots,0}_{r_s+r_t+k\,\hbox{ times}})$$ It is not hard to see that the image of the differential at $ O_\sigma$ of $ \pi$ (resp. $ \pi'$) is generated by $ T_\sigma\Sigma$ and $ T_\sigma \mathcal G$ (resp. $ T_{\sigma'}\Sigma'$ and $ T_{\sigma'}\mathcal G'$). Both $ \pi,\pi'$ 
are therefore surjective submersions in a neighborhood  of $ O_\sigma  $ that we denote by $ P_x''  $ again.    
Moreover, we claim that the commutativity of diagram \eqref{eq:def:bisubequiv} holds, i.e.,
\begin{equation}\label{diag:weak-equi-proof}
    \scalebox{0.7}{ \xymatrix{&& M&& \\W\ar[urr]^{s}\ar[drr]_{t}&&\ar[ll]_<<<<<\pi  P_x'' \subset \mathrm{Gr}_\phi(\Sigma)\times \mathbb R^{r_s+r_t+k}\ar@{..>}[u]\ar@{..>}[d]\ar[rr]^<<<<<<{\pi'} && W'\ar[llu]_{s'}\ar[dll]^{t'}\\ && M&& }}
 \end{equation}
This is an easy consequence of the following facts, valid for every index
\begin{enumerate} 
\item the vector fields $X_i$ and $ X_i'$ are $s$-related to the same vector field in $ \mathcal F$ (namely $ \underline{X_i}$) and $t$-related with the same vector field on $ \mathcal F$ (namely, zero).
\item the vector fields $Y_i$ and $ Y_i'$ are $t$-related to the same vector field in $ \mathcal F$  (namely $ \underline{Y_i}$) and $s$-related with the same vector field on $ \mathcal F$ (namely, zero).
\item the vector fields $Z_i$ and $ Z_i'$ are $s$- and $t$-related to the same vector field in $ \mathcal F$  (namely zero).
\item the commutativity holds true on $\mathrm{Gr}_\phi(\Sigma)\times 0^{r_s+r_t+k}  $ by definition of a morphism of bisubmersions (which means that $s(\sigma)=s'(\sigma')$ and $t(\sigma)=t'(\sigma')$ for all $(\sigma,\sigma')\in \mathrm{Gr}_\phi(\Sigma)$.
\end{enumerate}
This shows that every point in $x \in W$ admits a neighborhood $ \mathcal U_x$ on which an equivalence of bisubmersion exists between that neighborhood and some open subset of $W'$.
There exists a countable family $ (x_i)_{i \in \mathbb N}$ such that the open subsets  $\mathcal U_{x_i}$ cover $W$. The disjoint union $P_1:= \coprod_{i \in \mathbb N} P_{x_i} $ is an equivalence between $W$ and an open subset of $W'$. Now, by the same reasoning, an equivalence $P_2$ between $ W'$ and an open subset of $W$ can be constructed. Their disjoint union $ P_1 \coprod P_2$ is an equivalence between the bisubmersions $W$ and $W'$.

 This completes the proof of the claim.
\end{proof}

\begin{exo}
\label{exo:ifsurjective}
Consider two bisubmersions $ M \stackrel{s}{\leftarrow} W \stackrel{t}{\rightarrow} M $ and $ M \stackrel{s'}{\leftarrow} W' \stackrel{t'}{\rightarrow} M $ for $ \mathcal F$. Using Theorem \ref{thm:equiv}, show that if there exists a surjective morphism of bisubmersions from $W$ to $W'$, then $W$ and $ W'$ are equivalent.
\end{exo}

\begin{exo}
\label{exo:samedimension}
Consider two bisubmersions $ M \stackrel{s}{\leftarrow} W \stackrel{t}{\rightarrow} M $ and $ M \stackrel{s'}{\leftarrow} W' \stackrel{t'}{\rightarrow} M $ for $ \mathcal F$ of the same dimension.
Show that if $x \in W$ and $x'\in W'$ are equivalent, then there exist neighborhoods $\mathcal U $, $ \mathcal U'$ of $x$ and $x'$ and an isomorphism of bisubmersions $ \psi \colon \mathcal U \to \mathcal U'$, i.e., a morphism which is a diffeomorphism onto its image. 
\emph{Hint:} By the proof of Theorem \ref{thm:equiv}, there exists an equivalence $y \in W''$ such that $ \pi(y)=x$ and $\pi'(y) =x' $. Consider a local bisection $ \Sigma$ though $y$, i.e., a submanifold of $W''$ on which $ \pi,\pi'$ are diffeomorphisms onto their images. Such local bisections exist. The induced map $\underline{\Sigma} $ does the job. 
\end{exo}

\vspace{.5cm}

Theorem \ref{thm:equiv} can be restated as in the following manner, which is interesting by itself:

\begin{corollary}
\label{cor:restatement}
Consider two bisubmersions $ M \stackrel{s}{\leftarrow} W \stackrel{t}{\rightarrow} M $ and $ M \stackrel{s'}{\leftarrow} W' \stackrel{t'}{\rightarrow} M $ for $ \mathcal F$. 
For any two points $x \in W$ and $x' \in W' $, the following statements\footnote{Notice that all these statements imply $ s(x)=s(x')$ and $ t(x)=t(x')$.} are equivalent:
\begin{enumerate}
    \item[(i)] A neighborhood $\mathcal U$ of $x$ is $W$ is equipped with a morphism of bisubmersions $ \phi\colon \mathcal U \to W'$ mapping $x$ to $x'$.
    \item[(ii)] A neighborhood $\mathcal U'$ of $x'$ is $W'$ is equipped with a morphism of bisubmersions $ \phi' \colon \mathcal U' \to W$ mapping $x'$ to $x$.
    \item[(iii)] There exist local bisections $\Sigma $ through $x$ and $ \Sigma'$ through $x'$ that induce the same germ of isomorphisms of $ \mathcal F$ (i.e., $\underline{\Sigma}=\underline{\Sigma}' $ near $s(x)$).
    \item[(iv)] $\mathrm{Out}_W(x) =\mathrm{Out}_{W'}(x')  $.
    \end{enumerate}
\end{corollary}

\section{Holonomy groupoid of a singular foliation}

We now introduce the holonomy groupoid of Androulidakis and Skandalis, defined in the  \cite{AS}, Section 3. We shall use a presentation which may seem quite different from the original one, but the difference is only a difference of presentation.

\subsection{The holonomy groupoid without its topology}
\label{sec:holonomygroupoid}

An important point about the holonomy groupoid is that it is a topological groupoid, and even a diffeological groupoid (although we will not develop this point here). To start with, we present a construction strongly inspired by a construction given by Garmendia and Villatoro \cite{GV} of Androulidakis-Skandalis' holonomy groupoid. But this construction will be incomplete,  because it only describes it as a set - without a topology\footnote{Villatoro and Garmendia do equip it with a topology, but using an infinite dimensional manifold structure that we do not wish to introduce here.}. We will call it for the moment the ``flowing groupoid'' of a singular foliation, and later on prove that it is isomorphic, as a groupoid (in the category of sets), to the properly defined holonomy groupoid (which will be a topological groupoid, i.e., a groupoid in the category of topological spaces - and even more). 

Let $ (M,\mathcal F)$ be a foliated manifold. In Section \ref{sec:symmetry}, we defined two sets, both equipped with two projections onto $M$ called source\footnote{Conventions on source and target of a groupoid are not the same in non-commutative geometry and in Poisson geometry. We use the following convention: for $\Gamma \toto M$ a groupoid over $M$, the product of two elements $ \gamma_1$ and $\gamma_2$ is defined if $t(\gamma_1)=s(\gamma_2) $ and the product $\gamma_1 \cdot \gamma_2 $ admits the source of $ \gamma_1$ as its source and the target of $ \gamma_2$ as its target. } and target and denoted by $s$ and $t$:
\begin{enumerate}
\item The set $  {\mathrm{Sym}}_\mathcal F   $ equipped with two maps $ s: {\mathrm{Sym}}_\mathcal F \to M$ and  $ t: {\mathrm{Sym}}_\mathcal F \to M$ called source and target, respectively, such that $ s^{-1}(m) \cap t^{-1}(n)$ is for all\footnote{For a given $m,n$, it may of course be empty. For $m,n$ on the same leaf, it is never empty in view of Theorem \ref{thm:thm:landscape}. For $m,n$ on the regular part, provided all regular leaves have the same dimension, it is also non-empty.} $m,n \in M $ made of germs of local isomorphisms $ \phi$ of $\mathcal F $ mapping a $m$ to $n $.  
 \item 
The set
${\mathrm{OutSym}}_\mathcal F$
obtained by taking the quotient of the previous set modulo the equivalence relation 
$ \Phi \sim \Psi $ if $ \Phi \circ \Psi^{-1} $ coincides in a neighborhood of $m$ an inner symmetry having $m$ as a very fixed point\footnote{See Definition \ref{def:veryfixed} for a definition. We saw in that section this condition is equivalent to: $\Psi \circ \Phi^{-1}$ coincides in a neighborhood of $n$ with an inner symmetry having $n$ as a very fixed point.}. 
\end{enumerate}
So far, we have not spoken of any structure on these sets. We do it now. 
\begin{enumerate}
\item The set ${\mathrm{Sym}}_\mathcal F $ admits a natural groupoid structure: source and target are already defined, the composition consists in composing germs of local isomorphisms of $\mathcal F $, the inverse is obtained by inverting such an isomorphism, and the unit map consists in mapping $ m\in M$ to the germ of the identity map.
\item  The set
${\mathrm{OutSym}}_\mathcal F$ also admits a groupoid structure, obtained as a quotient of the previous one. Let us describe now.
Consider the bundle of groups  ${\mathrm{K}}_\mathcal F$ over $M$ whose fiber at $m \in M$ are the germs of inner symmetries admitting $m$ as a very fixed point.
This bundle of groups is by construction included into ${\mathrm{Sym}}_\mathcal F $.
It follows from Exercise \ref{exo:InnerAndSymmetry2} that it is a normal bundle or groups, i.e.,  
 $$  [\phi]^{-1}  \left.{\mathrm{K}}_\mathcal F \right|_{m} [\phi] \in \left.{\mathrm{K}}_\mathcal F \right|_{n}  $$
for any  germ $[\phi]$ of a local isomorphism $\phi $ of $\mathcal F $ such that $ \phi(m)=n$.
Since the quotient ${\mathrm{Sym}}_\mathcal F $ of  $ {\mathrm{K}}_\mathcal F $ is precisely ${\mathrm{OutSym}}_\mathcal F $, the latter inherits a groupoid structure. 
\end{enumerate}

By construction, there is a short exact sequence of groupoids:
\begin{equation}
\label{eq:DefOut} \xymatrix{ \ar@{^(->}[r] {\mathrm{K}}_\mathcal F \ar[d] &\ar@<1ex>[d] \ar@<-1ex>[d] {\mathrm{Sym}}_\mathcal F 
\ar@{->>}[r]
&  \ar@<1ex>[d] \ar@<-1ex>[d] {\mathrm{OutSym}}_\mathcal F\\ M\ar@{=}[r]&M\ar@{=}[r] &M  }  
\end{equation}

\vspace{1cm}
Now, there is a second natural groupoid ${\mathrm{InnerSym}}_\mathcal F \subset {\mathrm{Sym}}_\mathcal F $ of ${\mathrm{Sym}}_\mathcal F \toto M $ 
obtained by considering all germs of local inner symmetries of $ (M,\mathcal F)$, see Definition \ref{def:InnerSym}.
As a set, it consists in taking all smooth time dependent vector fields $ (X_\epsilon)_{\epsilon \in [0,1]}$ in $\mathcal F $ (recall that those were defined in Section \ref{sec:symmetry}), then to take germs of their time-$1$ flows at every point where they are well-defined.

\begin{lemma}
Let $ (M,\mathcal F)$ be a foliated manifold.
 $\mathrm{InnerSym}_\mathcal F \toto M $ is a subgroupoid of $ {\mathrm{Sym}}_{\mathcal F} \toto M$.
\end{lemma}
\begin{proof}
 Lemma \ref{lem:innerIs} can be adapted to show that $\mathrm{InnerSym}_\mathcal F$ is indeed a groupoid.
Corollary \ref{coro:flowissymmetry} states that any local inner symmetry is a symmetry of $ \mathcal F$, hence we have a groupoid inclusion ${\mathrm{InnerSym}}_\mathcal F \subset {\mathrm{Sym}}_\mathcal F $.
\end{proof}

Since $\mathrm{InnerSym}_\mathcal F \toto M $  contains $ \mathcal K_\mathcal F $, Equation \eqref{eq:DefOut} allows defining the flowing groupoid\footnote{Our construction of the flowing groupoid is equivalent to the one completed by \cite{GV} using $ \mathcal F$-cuts to the leaves.} 

\begin{definition}
\label{def:Flowinggroupoid}

Let $(M,\mathcal F) $ be a foliated manifold. We call \emph{flowing groupoid\footnote{The flowing groupoid will be soon isomorphic (as a set) to the holonomy groupoid, see Proposition \ref{prop:flow=holonomy}! So we do not insist too much on the notion.} of $\mathcal F$ } and denote by 
$$  {\mathrm{Flow}}_\mathcal F \toto M   $$
the image of $\mathrm{InnerSym}_\mathcal F $ through the projection $\mathrm{Sym}_\mathcal F \longrightarrow \mathrm{OutSym}_\mathcal F $. 
\end{definition}
 We will see in Proposition \ref{prop:flow=holonomy} that the flowing groupoid \underline{coincides} with the holonomy groupoid, but it is however important for pedagogical reasons, and also because it is the most practical manner to construct the holonomy groupoid of a given singular foliation, see, e.g., Exercises \ref{exo:flowing1} and \ref{exo:flowing2} below.  
By construction, there is a short exact sequence of groupoids:
$$ \xymatrix{ \ar@{^(->}[r] {\mathrm{K}}_\mathcal F \ar[d] &\ar@<1ex>[d] \ar@<-1ex>[d] {\mathrm{Inner}}_\mathcal F 
\ar@{->>}[r]
&  \ar@<1ex>[d] \ar@<-1ex>[d] {\mathrm{Flow}}_\mathcal F\\ M\ar@{=}[r]&M\ar@{=}[r] &M  }  .$$

\begin{exo}
\label{exo:flowing1}
    Let $ \mathcal F$ be the singular foliation on $M= \mathbb R^n$ with $ n \geq 2$  made of all vector fields vanishing  at $0$  as in Section \ref{ex:singFolVanish}.
    Show that the flowing groupoid is the disjoint union of the pair groupoid of $M \backslash \{0\} $ with the group of invertible matrices of positive determinant (over the origin). 
 \end{exo}

\begin{exo}
\label{exo:flowing2}
    Let $ \mathcal F$ be the singular foliation on $M= \mathbb R^n$ with $ n \geq 2$ made of all vector fields vanishing  at $0$ at order $k$  as in Section \ref{ex:singFolVanish}.
    Show that the flowing groupoid is the disjoint union of the pair groupoid of $M \backslash \{0\} $ with the quotient of the group of formal diffeomorphism admitting the origin as a fixed point, whose derivatives at $0$ vanish up to order $ k-1$ by the group of formal diffeomorphism admitting the origin as a fixed point, whose derivatives at $0$ vanish up to order $ k-1$. Show that this group is in fact isomorphic to a vector space, equipped with the addition as a group product, if $ k \geq 2$.
 \end{exo}

\subsection{Atlases}

Let $ (M,\mathcal F)$ be a smooth singular foliation.
Let us introduce a very particular type of bisubmersions\footnote{Following Section 1.3.2 in \cite{AS}.}, the quotient of which is going to be, by definition, the holonomy groupoid (now equipped with a topology).

\vspace{0.5cm}

\begin{definitions}{Atlases: a groupoid-like bisubmersions}{Atlas}
Let $(M,\mathcal F)$ be a foliated manifold.
We say that a bisubmersion 
$ M \stackrel{s}{\leftarrow} W \stackrel{t}{\rightarrow} M $ over $\mathcal F $ is an \emph{atlas of $\mathcal F $} when 
\begin{enumerate}
    \item $W$ is equivalent\footnote{See Definition \ref{def:EquvAndMorph}} to its inverse\footnote{See Section \ref{sec:productinverse}.} $ W^{-1}$,
    \item the composition\footnote{See Section \ref{sec:productinverse}.}  $W*W $ is equivalent to $W$.
    \item $W$ admits local unit maps\footnote{See Definition \ref{def:unitsbisub}.}
\end{enumerate}
Also, atlases are said to be \emph{equivalent} when they are equivalent as bisubmersions over  $\mathcal F $.
\end{definitions}

\vspace{.5cm}

\begin{exo}
\label{exo:noneedofthird}
We show in this exercise that the third assumption in Definition \ref{def:Atlas} is a consequence of the two first ones, and can therefore be omitted.
Consider $ M \stackrel{s}{\leftarrow} W \stackrel{t}{\rightarrow} M $ an atlas of $\mathcal F $.
Let $m \in M $ be a point. Let $x \in W$ be such that $ t(x)=m $. Let $\Sigma$ be a germ of a bisection of $W$ through $ x$.
\begin{enumerate}
\item Consider the inverse bisection $ W^{-1} := M \stackrel{t}{\leftarrow} W \stackrel{s}{\rightarrow} M $, and $ \Psi : W \to W^{-1} \simeq W$ a morphism of bisubmersions defined in a neighborhood of $x$. Let $ \Sigma' = \Psi(\Sigma)$. Show that $\Sigma'$ is a bisection of $W'$ through $x'=\Psi(x)$ and that $\underline\Sigma'=\underline \Sigma^{-1}$.

\item Show that the graph $ \Sigma'' := \mathrm{Gr}(\Psi|_{\Sigma})=\{(y,\Psi(y)| y \in \Sigma\}$ of the restriction  $\Psi|_{\Sigma} $ of $\Psi$ is a bisection through $ (x,x')$ for the product bisection $W*W$. 
\item  Let $\Xi \colon W*W \to W $ be a morphism of bisubmersions defined in a neighborhood of $ (x,x')$. Show that $\Xi(\Sigma'') $ is a bisection of $W$ through $\Xi(x,x') $ whose induced diffeomorphism is the identity map in a neighborhood of $m$.
\item Use Exercise \ref{exo:units} to conclude the argument.
\end{enumerate}
\end{exo}

Here is a  consequence of Theorem  \ref{thm:equiv}, more precisely of the equivalence of item (i) and (iii) in that Theorem.

\begin{proposition}
A bisubmersion $ M \stackrel{s}{\leftarrow} W \stackrel{t}{\rightarrow} M $ over $\mathcal F $ is an atlas of $\mathcal F $ if and only if the image ${\mathrm{Out}}_W (W)$ of the outer-germ map is a sub-groupoid of ${\mathrm{OutSym}}_\mathcal F \toto M$.

Also, two atlases $W$ and $W'$ are equivalent if and only if the sub-groupoids ${\mathrm{Out}}_W (W) \toto M$ and ${\mathrm{Out}}_{W'}(W')\toto M $ coincide. 
\end{proposition}

The following result, presented as an exercise, can be practical to compute atlases.
\begin{exo}
    Let $ M \stackrel{s}{\leftarrow} W \stackrel{t}{\rightarrow} M $ of $\mathcal F $  be a bisubmersion over a singular foliation $\mathcal{F}$ that admits a (local) unit $\epsilon\colon M\hookrightarrow W$. Show that $W$ is an atlas if and only if $W*W^{-1}$  is equivalent to $W$.
\emph{Hint:} Consider the following commutative diagrams
\begin{align}
     \xymatrix{  &W \ar[dl]_{s} \ar[dr]^{t} & \\ M & P\ar@{->>}[u]^{\pi} \ar@{->>}[d]_{\pi'} \ar[l] \ar[r]  &M\\&W*W^{-1}\ar[ul]^{t\circ\mathrm{pr}_2} \ar[ur]_{t\circ\mathrm{pr}_1}   & } & \hspace{1cm}  \xymatrix{  &W \ar[dl]_{s} \ar[dr]^{t} & \\ M& P\ar@{->>}[u]^<{\mathrm{pr}_1\circ\pi'} \ar@{->}[d]^<<<{\pi\times \mathrm{pr}_2\circ\pi'} \ar[l] \ar[r] &M\\ &W*W\ar[ul]^{s\circ\mathrm{pr}_2} \ar[ur]_{t\circ\mathrm{pr}_1}   & } 
\end{align}
and use Theorem \ref{thm:equiv}.
\end{exo}

\vspace{0.3cm}

We are in particular interested in a particular subclass of atlases, that we call inner atlases.

\vspace{1cm}

\begin{definitions}{An inner atlas}{definitionFundAtlas}
  Let $ (M,\mathcal F)$ be a foliated manifold. 
  We say that an atlas $W$ for $\mathcal F $
 is an inner atlas if  
 $$ {\mathrm{Out}}_{W}(W)   = {\mathrm{Flow}}_\mathcal F .$$  
 \end{definitions}
\vspace{0.5cm}

Definition \ref{def:definitionFundAtlas} can be restated as meaning that  
\begin{enumerate}
\item 
for any local bisection $ \Sigma$ through $ w \in W$, the local symmetry  $\underline{\Sigma}  $ coincides, in a neighborhood of every point where it is defined, with an inner symmetry of $\mathcal F $, and, \item any inner symmetry of $ \mathcal F$ is induced by a bisection of $ W$, at least near any point where it is defined.
\end{enumerate}
Even more explicitly, an atlas is an inner atlas if and only if: 
\begin{enumerate}
\item For any local diffeomorphism of the type $\Phi := \phi_1^{X_1} \circ \dots \circ \phi_n^{X_n} $ with $ X_1, \dots, X_n$ in $\mathcal F $, and any $m \in M$ such that $\Phi(m) $ is well-defined, there exists $ x \in W$
 and a bisection $\Sigma $ through $x$ such that $\underline{\Sigma} =\Phi $ is a neighborhood of $m$. 
 \item Conversely, for any  bisection $\Sigma $ through $x$, there exists  $ X_1, \dots, X_n$ in $\mathcal F $ such that $m=\Phi(s(x)) $ is well-defined and $\underline{\Sigma} =\Phi $ is a neighborhood of $m$. 
\end{enumerate}

Notice that we have not assumed in Definition \ref{def:definitionFundAtlas} that the inner atlas is an atlas: it is a priori simply a bisubmersion. However, it follows from the definition that it has to be an atlas:

\begin{lemma}
Any inner atlas for $ \mathcal F$ is an atlas for $ \mathcal F$.
\end{lemma}
\begin{proof}
By equivalence of item i) and iii) in Corollary \ref{cor:restatement}, a bisubmersion is an atlas if and only if its germs of bisections are stable under inverse and product. 
It is of course the case for all local diffeomorphisms as in the second item of Definition \ref{def:definitionFundAtlas}.
This completes the proof.
\end{proof}

\begin{lemma}\label{lem:onlyone}
Any two inner atlases for a singular foliation $ \mathcal F$ are equivalent.
\end{lemma}
\begin{proof}
This is an immediate consequence of the equivalence of (i) and (iii) in  Theorem \ref{thm:thm:equiv}.
\end{proof}

\begin{exo}
\label{exo:innerleaf0}
Let $W$ be an inner atlas with source $s$ and target $t$ for $ \mathcal F$. Show that two points $\ell_1,\ell_2 \in M$ are in the same leaf of $ \mathcal F$ if and only if there exists $w \in W$ such that $ s(w)=\ell_1$
and $ t(w)=\ell_2$.
The converse is not true, see Exercise \ref{exo:innerleaf}

\end{exo}

\begin{proposition}\label{ex:Atlas-Groupoid}
Any Lie groupoid with source-connected fibers is an inner atlas for its basic singular foliation.
\end{proposition}

\begin{proof}
Let $\Gamma \toto M $ be a Lie groupoid with $s$- and $t$-connected fibers. Let $ (A \to M, \rho, [\cdot\,, \cdot ])$ be its algebroid and $ \mathcal F = \rho(\Gamma(A))$ its basic singular foliation.  
Recall that $\Gamma $ is a bisubmersion of $\mathcal F $ (see Exercise \ref{exo:groupasbosubmersions}).
The space $ \Gamma$ is an atlas of $\mathcal F $, since the inverse map  $ {\mathrm{inv}} \colon \Gamma \to \Gamma$ is an equivalence between 
$$  \xymatrix{ & \Gamma \ar[dd]^{{\mathrm{inv}}} \ar[dr]^{s} \ar[dl]_{t} &  \\M &  & M\\&  \ar[ur]^{t} \ar[ul]_{s} \Gamma& \\}  $$
while the product $\mu $ gives an equivalence of bisubmersions:
$$\xymatrix{ & \Gamma \times_{s,M,t} \Gamma \ar[dd]^{\mu} \ar[dr]^{t\circ\mathrm{pr}_1} \ar[dl]_{s\circ\mathrm{pr}_2} &  \\M &  & M\\&\ar[ur]_{t}\ar[ul]^{s} \Gamma& \\}
$$
This proves that $ \Gamma$ is an atlas. To show that it is an inner atlas, we proceed as follows. Since fibers of $s$ are connected, there exists, for any $\gamma  \in \Gamma$ with source $m$, a time dependent section $ a_t$ of $A$ such that the time $1$-flow $ \Phi_t$ of $ \overrightarrow{a_t} \in \mathfrak X(\Gamma)$ maps $m$ (seen as an element in $\Gamma$ through the unit map $ \epsilon$) to $ \gamma$. Let $ \mathcal U \subset M$ be a neighborhood of $m$ on which $ \Phi_t \circ \epsilon $ is well-defined.  
Then $ \Phi_1 \circ \epsilon (\mathcal U)$ is a bisection of $ \Gamma$ through $ \gamma$. Its induced diffeomorphism $ \underline{\Phi_1 \circ \epsilon (\mathcal U)}$ is the time $1$-flow of the vector field $\rho(a_t) $, which is a vector field on $ \mathcal F$ that depends smoothly on $t$. It is therefore an inner symmetry of $ \mathcal F$.
\end{proof}

Here is a non-example of inner atlas.

\begin{example}
\label{ex:noninner}
Let $\mathcal F $ be a singular foliation on a manifold $M$. 
Let $W$ be an inner atlas with source $s$ and target $t$. Assume that there exists a  finite group $G$ acting on $M $ by symmetries of $\mathcal F $, denoted by $ \phi_g : M \to M$ for a given $g \in G$. Then 
$$  \coprod_{g \in G} \{g\} \times  W  $$
is an atlas for $\mathcal F $, when equipped with the source $(g,w)\mapsto s(w) $ and $ t(g,w) \mapsto \phi_g (t(w))$. If $\phi_g$ is not an inner symmetry for at least one $g \in G $, it is however not an inner atlas, since the symmetry associated to the bisection $m \mapsto (\{g\},\epsilon(m)) $, with $ \epsilon$ a local unit map for $W$, is not an inner symmetry of $\mathcal F$.
\end{example}

\begin{exo}
\label{exo:innerleaf}
Apply the construction of Example \ref{ex:noninner} to vector fields on $\mathbb R^2 $ vanishing quadratically at zero,  and to a finite sub-group of the group of rotations.
Show that the converse of the statement in Exercise \ref{exo:innerleaf0} is not true.
\end{exo}

The most important of all atlases (that will be proven soon to be an inner atlas) was introduced by Androulidakis and Skandalis under the name of  \emph{path holonomy atlas of $\mathcal F $ associated to an anchored bundle and a connection}.
It is obtained by the following procedure:
\begin{enumerate}
    \item Take a family $ (A_i \to  \mathcal U_i,\rho_i)_{i \in I}$ of anchored bundles
such that the open subsets $ (\mathcal U_i)_{i \in I}$ where they are defined cover $ M$. Assume each one of them is equipped with a connection $\nabla_i $. 

    \item For every $i \in I$, there exists a neighborhood $\mathcal A_i $ of the zero section of $ A_i \to \mathcal U_i$ on which there is a bisubmersion over $\mathcal F $ (see Proposition \ref{thm:Prop:FundamentalExampleAbstract}).
    \item Then consider the disjoint union\footnote{These sets are called \emph{longitudinal charts} in \cite{AS}.} for all $n \geq 0$ and all $ i_1, \dots,i_n \in I$ of all direct products 
     $$ \mathcal A_{i_1}^\star \times_M \dots \times_M \mathcal A_{i_n}^\star \hbox{  ($n $ times) } $$
     where $ ^\star$ means that we consider $\mathcal A$ or its inverse $\mathcal A^{-1} $. 
     \item[] Denote by $ W_{path}$ the disjoint union of all manifolds as in the third item.
\end{enumerate}
It deserves to be noticed that every connected component of $W_{path}$ the manifold is finitely dimensional, although the dimension is not bounded. We have the following lemma:

\begin{lemma}
$W_{path}$ is an atlas  for $\mathcal F $.
\end{lemma}
\begin{proof}
There is a natural isomorphism between $W_{path}$ and
$W_{path}^{-1}$. Also $ W_{path}* W_{path}$ injects into $ W_{path}$ by construction.
\end{proof}

Following \cite{AS}, we call \emph{path holonomy atlas  of $\mathcal F $} such an atlas.

\begin{proposition}
    Any path holonomy atlas  of $\mathcal F $ is an inner atlas of $ \mathcal F$.
\end{proposition}
\begin{proof}
This is a direct consequence of Example \ref{ex:smallinner}, which shows that each one of the fundamental bisubmersion $ \mathcal A_i$ contains all small inner symmetries on $ \mathcal U_i$. In turn, since any inner-symmetry $ \Phi$ is, in a neighborhood of any point where it is defined, a composition
$$ \Phi= \phi_1 \circ \dots \circ \phi_n $$
where for every index $k$, we have $\phi_k= \underline{\Sigma_k} \subset \mathcal A_{i_k} $ a bisection as in Definition \ref{def:smallinner} for some $ i_k \in I$. Now, 
$ \Sigma := \Sigma_n \times_M \dots \times_M \Sigma_1 $ is a bisection of the bisubmersion  $  \mathcal A_{i_1} \times_M \dots \times_M \mathcal A_k$. 
By construction, and in view of Equation \eqref{eq:inversionbisub}, we have $$ \underline{\Sigma}= \underline{ \Sigma_n \times_M \dots \times_M \Sigma_1}=\underline{ \Sigma_1} \times \dots \times \underline{ \Sigma_n}  =\phi_1 \dots \phi_n= \Phi.$$ This completes the proof. 
\end{proof}

\subsection{Holonomy groupoid}

We can now at last define the holonomy groupoid of Androulidakis and Skandalis \cite{AS}. We start with a general result.

\vspace{.5cm}

\begin{propositions}{From atlases to groupoids}{EquivRel}
Let $ M \stackrel{s}{\leftarrow} W \stackrel{t}{\rightarrow} M $ of $\mathcal F $ be an atlas of $\mathcal F $. Consider the  equivalence relation on $W$ given by $ x \sim x'$ if and only if $x$ and $x'$ have neighborhoods which are equivalent as bisubmersions of $\mathcal F $.

The equivalence classes of this relation form a topological\footnote{A subset of the groupoid is open if and only if its inverse image in the atlas (which is a manifold) is open.} groupoid over $M$.
\end{propositions}
\begin{proof}
The inverse of an equivalence class represented by $x \in W$ is represented by any point in $W^{-1} $ which  equivalent to $x$. It is easy to check that it is well-defined.
The same applies to product. Given two compatible $x_1,x_2$, there exists $x_3 \in X$ such that $ x_3$ and $(x_1, x_2) $ are equivalent. This construction goes to the quotient w.r.t. the equivalence relation on $X$ and defines the groupoid product. All these maps are continuous.
  \end{proof}

Theorem \ref{thm:thm:equiv}
implies that equivalent atlases  of $ \mathcal F$ induce canonically homeomorphic quotient groupoids.
In particular, Definition \ref{def:Holonomy} makes sense.

\vspace{.5cm}
\begin{definitions}{Holonomy groupoid of a singular foliation.}{Holonomy}
We call holonomy groupoid of $\mathcal F $ the topological\footnote{In \cite{AS}, this groupoid is equipped with much more than just a topology: a $ C^*$-algebra of "smooth functions" is even introduced. We will not go as far here.} groupoid associated to a fundamental atlas of $\mathcal F $. We denote it by $ {\mathrm{Hol}}(\mathcal F) \toto M$.
\end{definitions}

\vspace{.5cm}
To start with, let us state that the flowing groupoid, introduced for pedagogical reason, is a notion we can now get rid of.

\begin{proposition}
\label{prop:flow=holonomy}
For any singular foliation, the flowing groupoid\footnote{See Definition \ref{def:Flowinggroupoid}} and the holonomy groupoid are canonically isomorphic.
\end{proposition}
\begin{proof}
The flowing groupoid can be redefined as being the image through $ {\mathrm{Out}}_\mathcal F$ of any inner atlas. The quotient that defines the holonomy groupoid consists in identifying two points in an inner atlas whose image in ${\mathrm{Sym}}_\mathcal F $ coincide. The result follows.
\end{proof}

\begin{remark}
Let us assume that the singular foliation $(M,\mathcal F)$ comes from a Lie groupoid $\Gamma \toto M $ whose $s$- and $t$-fibers are connected, i.e., let us assume that $\mathcal F =\rho\left(\Gamma(A)\right) $ with $A$ the Lie algebroid of $\Gamma \toto M$. In Proposition  \ref{ex:Atlas-Groupoid}, we saw that $\Gamma $ is an inner atlas for $\mathcal F $. 
There is therefore a surjective groupoid morphism $ \Pi$ as follows: 
$$ \xymatrix{\Gamma \ar[r]^\Pi \ar@<0.5ex>[d] \ar@<-0.5ex>[d]  & {\mathrm{Hol}}({\mathcal F}) \ar@<0.5ex>[d] \ar@<-0.5ex>[d] \\ M \ar[r]^{=}& M.} $$
Let us study the fibers of this map $ \Pi$. For that purpose, let us consider 
$$ \mathcal K := \{ a \in \Gamma(A)   \hbox{ s.t. }   \rho(a) \}  .$$
By construction $\mathcal K$ is the $ \mathcal C^\infty(M)$-module of sections of $A$ which are valued in the kernel of $\rho$ for all point $ m \in M$.
Any section $ a \in \mathcal K$ gives, by right action on $\Gamma $, a vector field $\overrightarrow{a} $ in $\Gamma $.
The $ \mathcal C^\infty(\Gamma)$ module that these vector fields generate is the module of bi-vertical vector fields. Assume for simplicity that it is a singular foliation (i.e., that it is locally finitely generated).
Then two points in $ \Gamma$ that are in the same leaf of the bi-vertical singular foliation are in the same fiber of $ \Pi$. See e.g., Exercise \ref{exo:OutConstAlongBiverti}.
 The quotient $ \Gamma/\sim$ of the Lie groupoid $\Gamma $ by the equivalence class defined by the bi-vertical singular foliation is a topological groupoid, and $ \Pi$ induces a continuous groupoid morphism $ \Gamma/\sim \longrightarrow {\mathrm{Hol}}(\mathcal F)$.
\end{remark}

\subsection{About smoothness of the holonomy groupoid: two theorems by Claire Debord}

Recall that a \emph{Lie groupoid} is a groupoid $ \Gamma \rightrightarrows M$ such that $ \Gamma$ and $M$ are manifolds, the source, and target are smooth surjective submersions, and all structural maps (unit, product, inverse) are smooth\footnote{or real analytic, or holomorphic, depending on the setting}.

The holonomy groupoid is certainly not a smooth groupoid in general. It is a topological groupoid, and the topology may be quite horrible -very far away from a manifold topology.

However, the following theorem was proven by Claire Debord \cite{Debord2}.

\vspace{0.5cm}

\begin{theorems}{Along a leaf}{DebordSmoothness}
\cite{Debord2}
The orbits of the holonomy groupoid $\mathrm{Hol}(\mathcal F) $ of a singular foliation $\mathcal F $ are the leaves of $\mathcal F $.
Moreover, its restriction to any leaf $L$:
\begin{enumerate}
    \item is a Lie groupoid, 
    \item whose Lie algebroid is the holonomy Lie algebroid\footnote{See Section \ref{sec:holonomy-Lie-algebroid}.} of the leaf $L$.
\end{enumerate}
\end{theorems}
\begin{proof}
The first sentence was established in Exercise \ref{exo:innerleaf0}. 
The hard part of the proof relies on a theorem that bounds below the periods of a periodic orbit of a vector field in a neighborhood of a point.
This lower bound forbids a bisubmersion to have ``too many quotients'', it makes the quotient "discrete-like". We refer to \cite{Debord2}
\end{proof}

Smoothness of the holonomy groupoid happens in a second situation. This theorem is also due to Claire Debord \cite{Debord}.

\vspace{.5cm}

\begin{theorems}{Projective case}{DebordCase}
The holonomy groupoid of a Debord\footnote{See Section \ref{sec:Debord}} foliation is a Lie groupoid whose  Lie algebroid is the Lie algebroid in Proposition \ref{thm:Debord}.
\end{theorems}

\subsection{The fundamental groupoid}\label{sec:fundamental}
An alternative approach to the construction of the holonomy groupoid is by anchored paths divided by a certain equivalence relation. This is how the Lie groupoid integrating a Lie algebroid has been constructed in \cite{MR1973056,1707.00265}. In the context of singular foliations, this approach has been taken in \cite{GV}, Section 1.1.2 in \cite{LGR} and Section 4.6 of \cite{LLS}.
\\

Le $\mathcal F$ be a singular foliation on $M$ and $A$ an anchored bundle over $\mathcal F$. We start by introducing the right notion of homotopy between anchored paths (cf. Definition \ref{def:anchoredpaths}). For this, we introduce part of the longitudinal complex\footnote{$(\Lambda^\bullet_{\mathcal C^\infty(M)} \mathcal F)^*$ is \emph{not}, in general, isomorphic to the exterior tensor product $(\Lambda^\bullet_{\mathcal C^\infty(M)} \mathcal F^*)$ of ${\mathrm{Hom}}_{\mathcal C^\infty(M)} (\mathcal F, \mathcal C^\infty(M))$ } $(\Lambda^\bullet_{\mathcal C^\infty(M)} \mathcal F)* $ of $\mathcal F$ (See Section 4.1. in \cite{LLS}). Its cochains are given by skew-symmetric $\mathcal C^\infty(M)$-multilinear maps from $\mathcal F^{\times k}$ to $\mathcal C^\infty(M)$. For us, only the first 
differential $$d^{\mathcal F}:\mathcal F^*=\mathrm{Hom}_{\mathcal C^\infty(M)}(\mathcal F, \mathcal C^\infty(M))\to (\Lambda^2_{\mathcal C^\infty(M)}\mathcal F)^*,$$ defined by the usual formula $$(d^\mathcal F\alpha)(X,Y)=X(\alpha(Y))-Y(\alpha(X))-\alpha([X,Y])$$ will play a role. For an anchored bundle $A$, we will denote by $\rho^*$ the induces map from $\mathcal F^*$ to $\Gamma(A)$ and given an anchored morphism $\phi:A_1\to A_2$, we will denote by $\phi^*$, the induced map $\Gamma(\Lambda A_2^*)\to \Gamma(\Lambda A_1^*)$. With those conventions, we can finally define anchored homotopies:

\vspace{.5cm}

\begin{definition}
Let $A \stackrel{\pi}{\to}M$ be an anchored bundle and $\gamma,\tilde\gamma:I\to A$ anchored paths. We say that $\gamma$ is $A$-homotopic to $\tilde\gamma$, if there exists an anchored bundle morphism:
$$H: T[0,1]^2\to A$$ such that:
\begin{itemize}
    \item $H({t,0})=\gamma(t)$, $H({t,1})=\tilde \gamma(t)$.
    \item $H({0,s})$ and $H({1,s})$ are constant and equal to zero.
    \item The following diagram commutes:

    \begin{center}
    \xymatrix{\Omega^1(T[0,1]^2) \ar[d]^{d} &&\ar[ll]_{H^*\circ \rho^*} \mathcal F^*\ar[d]^{d^{\mathcal F}}\\
    \Omega^2(T[0,1]^2)&&\ar[ll]_{H^*\circ \rho^*}(\Lambda^2_{\mathcal C^\infty(M)}\mathcal F)^*
    } \end{center}
\end{itemize}
Then we call $H$ an anchored homotopy between $\gamma$ and $\tilde \gamma$.
\end{definition}

\vspace{.5cm}

In order to define the fundamental groupoid of a singular foliation, we just need to make paths composable, since in general the concatenation of smooth paths need not be smooth. The lazy solution to this problem is by using lazy paths (cf. e.g., \cite{LazyPeople} also referred to as sitting instants initially \cite{SittingInstants}), i.e., paths which are constant near the boundary.

\begin{definitions}{Fundamental groupoid}{} Let $M$ be a manifold and $\mathcal F$ a singular foliation. Let $A$ be any anchored bundle over $\mathcal F$.
The fundamental groupoid of $\mathcal F$ is the quotient
$$
\Pi_1(M,\mathcal F)=\frac{\mathrm{lazy~anchored~paths~in~}A}{A-\mathrm{Homotopies}}
$$
The composition of the fundamental groupoid is given by the concatenation of lazy paths.
\end{definitions}

\vspace{0.5cm}

As usual for such constructions, the source, and target maps are the projections on the endpoints of a path and the unit elements are given by the constant paths equal to $0$. The orbits of the fundamental groupoid of $\mathcal F $ are exactly the leaves of the foliation $ \mathcal F$. Also, for any regular leaf $L$, $\Pi_1(M,\mathcal F)_L$ is just the fundamental groupoid of the leaf.\\

 $\Pi_1(M,\mathcal F)_L$ can be seen as a universal cover of the holonomy groupoid in the following sense. There is a groupoid morphism:
 $$ \Pi_1(M,\mathcal F) \longrightarrow {\mathrm{Hol}}(\mathcal F)  $$
 such that for every $m \in M$, the restriction of the map above to $s^{-1}(m)$ is a submersion and a local diffeomorphism. 
 Since, by construction, the fundamental groupoid is source-simply-connected, it means that $ \Pi_1(M,\mathcal F)$ is obtained by taking, for every leaf $ L$, the source-simply-connected groupoid integrating the holonomy Lie algebroid $ A_L$. See \cite{GV,LGR,LLS} for several variations on this topic.

\section{Geometric resolutions of a singular foliation}\label{sec:Geometric-resolutions}

\subsection{Definition and universality of existence of geometric resolutions}

\label{sec:geometricresolutionsdef}

\subsubsection{Introduction}

Let us work within the smooth setting for the moment, and use Definition \ref{def:consensus} of singular foliations\footnote{Recall that the index $ {\Gamma}_c(E)$ means "compactly supported sections of $E$" and $ \mathfrak X_c(M)$ means "compactly supported vector fields on $M$"}.
We saw in Section \ref{sec:AnchoredBundle} that for every finitely generated singular foliation $\mathcal F $, there exists an anchored vector bundle\footnote{We can even assume $A$ to be a trivial bundle, see the discussion in Section \ref{sec:AnchoredBundle}.}
  $$  \xymatrix{ A \ar[r]^{\rho} \ar[d] & TM  \ar[d] \\ M \ar@{=}[r]  &  M}$$
  such that
  
  $\rho(\Gamma_c(A))  = \mathcal F$. Now, consider the kernel\footnote{Notice that this kernel consist in compactly supported sections of $ A$ which for every $m \in M$ are valued in the Strong kernel of $ \rho$ at $m$.} $\ker(\rho) $ of  $$ \rho \colon \Gamma_c(A) \longrightarrow \mathfrak X_c(M) .$$
  The space $\ker(\rho) \subset \Gamma_c(A)$ is again a $ \mathcal C^\infty(M)$-submodule of $\Gamma_c(A)$. If it is finitely generated, then there exists a second vector bundle\footnote{We can even assume $B$ to be a trivial bundle, see the discussion in Section \ref{sec:AnchoredBundle}.} $ B \to M$ and a vector bundle morphism $\dd \colon B \longrightarrow A $ such that
   $$ \dd \left(\Gamma_c(B)\right)  = \ker(\rho) $$
   In particular, we have $ \rho \circ \dd = 0$, and 
   $$ \xymatrix{ B \ar[r]^{\dd^{(2)}} \ar[d]&A \ar[r]^{\rho} \ar[d] & TM  \ar[d] \\M \ar@{=}[r]  & M \ar@{=}[r]  &  M}  $$ 
   is a complex of vector bundles which is exact at the level of sections, i.e., 
  that the sequence 
   $$ \xymatrix{ \Gamma_c(B)\ar[r]^{\dd^{(2)}} & \Gamma_c(A) \ar@{->>}[r]^{\rho}&   \mathcal F} 
  $$ is exact. The procedure continues when the kernel of
   $$  \dd^{(2)} \colon \Gamma_c(B) \longrightarrow \Gamma_c(A) $$ is again finitely generated as a $\mathcal C^\infty(M) $-module: there exists a vector bundle $C \to M$ and a vector bundle morphism $\dd \colon C \to B  $ such that $\dd^{(3)} (\Gamma_c(C)) = \ker(\dd^{(2)}) $. By construction, 
   $$ \xymatrix{ C \ar[r]^{\dd^{(3)}} \ar[d]& B \ar[r]^{\dd^{(2)}} \ar[d]&A \ar[r]^{\rho} \ar[d] & TM  \ar[d] \\M \ar@{=}[r] & M \ar@{=}[r]  & M \ar@{=}[r]  &  M}$$ 
   is a complex of vector bundles  and the following is an exact complex:
  $$ \xymatrix{\Gamma_c(C)\ar[r]^{\dd^{(3)}}  & \Gamma_c(B)\ar[r]^{\dd^{(2)}} & \Gamma_c(A) \ar@{->>}[r]^{\rho}&   \mathcal F} .
  $$ 
   Here is a natural question.

\vspace{.5cm}
  
  \begin{question}
  \label{ques:resol}
  
  When can the construction of the complex of vector bundles described above be continued ``up to infinity'' (i.e., can one be certain that the kernels are finitely generated) ?. 
  
  Does it stop at some point? (i.e., can we manage that the kernel of $\dd^{(k)} $ is trivial for $k$ large enough?
  
  Assume it can be constructed, what kind of geometric information is
 encoded in that complex?  \end{question}

\subsubsection{Definitions}
  
Let us start with by precise definition and a precise vocabulary (valid in the smooth, real analytic or complex settings).

\begin{definition}
   An \emph{anchored complex of vector bundles} consists of a triple $(E_{-\bullet},\dd^{(\bullet)},\rho)$, where
   \begin{enumerate}
       \item $ E_{-\bullet} = (E_{-i})_{i \geq 1}$ is a family of vector bundles over $M$, indexed by negative integers.
       \item  for every $i\geq 1$, $\dd^{(i+1)}\in \mathrm{Hom}(E_{-i-1}, E_{-i})$ is a vector bundle morphism over the identity of $M$  called the \emph{differential map}.
       \item $\rho \colon E_{-1} \longrightarrow TM $ is a vector bundle morphism over the identity of $M$ called the \emph{anchor map}. 
   \end{enumerate}
   such that  
    $$   \xymatrix{  \ar[r] & E_{-i-1} \ar[r]^{{\dd^{(i+1)}}} \ar[d] & 
     E_{ -i} \ar[r]^{{\dd^{(i)}}}
     \ar[d] & E_{-i+1} \ar[r] \ar[d] & \ar@{..}[r] & \ar[r]^{{\dd^{(2)}}}& E_{-1} \ar[r]^{\rho} \ar[d]& TM \ar[d] \\ 
      \ar@{=}[r] & \ar@{=}[r] M  &  \ar@{=}[r] M 
      &  \ar@{=}[r] M  &\ar@{..}[r] & \ar@{=}[r]   &  \ar@{=}[r] M  & M}
    $$

 form a complex, i.e., such that 
  $$ \dd^{(i)}\circ\dd^{(i+1)}=0  \hbox{ and }  \rho \circ \dd^{(2)}=0$$  
\end{definition}

Let us fix some vocabulary about such complexes:

\begin{enumerate}
    \item The integer $-i$ is called the \emph{degree} of the vector bundle $E_{-i}$. 
    The choice of negative numbers may seem surprising: it will be justified when introducing Lie $ \infty$-algebroid structures.
\item An anchored complex of vector bundle $(E_{-\bullet},\dd^{(\bullet)}, \rho)$ is said to be of \emph{length} $n\in\mathbb N$ if $E_{-i}=0$ for  $ i \geq n+1 $ and \emph{finite length} if $ E_{-i}=0$ except for finitely many indices.
\item 
 We shall speak of \emph{anchored complex of trivial bundles} when all the vector bundles $(E_{-i})_{ i\geq1}$ are trivial vector bundles. We do not assume $TM$ to be a trivial bundle.
\end{enumerate}

There are two main cohomologies that one can associate to an anchored complex of vector bundles.

\begin{enumerate}
    \item {\textbf{Cohomology at the level of sections}}.
     An anchored complex of vector bundles $(E_{-\bullet},\dd^{(\bullet)},\rho)$ induces a complex of sheaves of modules over functions. More explicitly, for every open subset $\mathcal U \subset M $, there is a complex: $$ \cdots {\longrightarrow} \Gamma_{\mathcal U}({E_{ -i-1}})
     \stackrel{\dd^{(i+1)}}{\longrightarrow} \Gamma_{\mathcal U}({E_{-i}})
     \stackrel{\dd^{(i)}}{\longrightarrow}{\Gamma_{\mathcal U}(E_{-i+1})}{\longrightarrow}\cdots  {\longrightarrow} \Gamma_{\mathcal U}(E_{-1})  
     \stackrel{\rho}{\longrightarrow}\mathfrak X(\mathcal U) .$$  In particular, $\mathrm{Im}\left(\dd^{(i+1)}\right)\subseteq\ker \dd^{(i)}$ for every $i\in\mathbb N$, so that the quotient spaces:
      $$ H^{-i}(E_{\bullet},\mathcal U) = \left\{ \begin{array}{ll} \frac{\ker 
      {\rho}}{\mathrm{Im}\left(\dd^{(2)}\right)} & \hbox{for $i=1 $} 
      \\ &\\\frac{\ker \dd^{(i)}} {\mathrm{Im}\left(\dd^{(i+1)}\right)} &
      \hbox{ if $ i \geq 2$} \end{array}\right. $$
    is a module over functions on $\mathcal U $ that we call \emph{$i$-th cohomology of $(E_{-\bullet}, \dd^{(\bullet)}, \rho) $ at the level of sections}.

    \item {\textbf{Cohomology at an arbitrary point $m \in M$}}. 
    
    An anchored complex of vector bundles $(E_{-\bullet},\dd^{(\bullet)},\rho)$ at an arbitrary point $m \in M$, restricts to a complex of vector spaces $$ \quad  \cdots {\longrightarrow} {E_{ -i-1}}_{|_m}\stackrel{\dd^{(i+1)}_{|_m}}{\longrightarrow} {E_{-i}}_{|_m}  \stackrel{\dd^{(i)}_{|_m}}{\longrightarrow}{E_{-i+1}}_{|_m}{\longrightarrow}\cdots  \stackrel{\dd^{(2)}_{|_m}}{\longrightarrow} E_{-1}  \stackrel{\rho_{{|_m}}}{\longrightarrow} T_m M .$$   In particular, $\mathrm{Im}(\dd^{(i+1)})\subseteq\ker \dd^{(i)}_{|_m}$ for every $i\in\mathbb N$, and we call the quotient vector spaces:
      $$ H^{-i}(E_{\bullet},m) = \left\{ \begin{array}{ll} \frac{\ker 
      {\rho_{|_m}}}{\mathrm{Im}\left(\dd^{(2)}_{|_m}\right)} & \hbox{for $i=1 $} 
      \\&\\ \frac{\ker \dd^{(i)}_{|_m}} {\mathrm{Im}\left(\dd^{(i+1)}_{|_m}\right)} &
      \hbox{ if $ i \geq 2$} \end{array}\right. $$
 the \emph{$i$-th cohomology of $(E_{-\bullet}, \dd^{(\bullet)}, \rho) $ at the point $m$}.  
\end{enumerate}

\begin{bclogo}[arrondi = 0.1, logo = \bcdz]{Warning !}
It is important to notice that 
$H^{-i}(E_{\bullet},m) $ may be non-zero at a point $m$ even if $H^{-i}(E_{\bullet},\mathcal U) $ is  zero in every open neighborhood $\mathcal U $ of $m$. 
For instance, for a Debord singular foliation $ \mathcal F$. Let $A$ stand for its associated Lie algebroid $A$ and $\rho$ the anchor. The complex of vector bundles defined by $ E_{-1}=A$ and $ E_{-i}=0$ for $i \geq 2$, together with a trivial $ \dd$ and the anchor $ \rho$, is exact at the level of sections on any open set. But this is not exact when evaluated at a singular point of $ \mathcal F$. 

The converse, however, is not possible, as we will see in Proposition \ref{prop:atapoontzeroaroundzero} below.
\end{bclogo}

A complex of vector bundles needs not to be exact at a point, even when it is exact at the level of sections. But if it is exact at a point, then it is also exact, locally, at the level of sections, as we now see. 

\begin{proposition}
\label{prop:atapoontzeroaroundzero}
Let $ i \in \mathbb N$.
Every $m \in M $ such that $H^{-i}(E_{\bullet},m) =0$ has an open neighborhood $ \mathcal U$ for which $H^{-i}(E_{\bullet},\mathcal U) =0 $.
\end{proposition}
\begin{proof}
Below, $\mathrm{rk}(E)$ stands for the rank of a vector bundle $E\to M$,
and $\mathrm{rk}(d)$ stands for the rank (= dimension of the image) of a vector bundle morphism $d$. Also, for $i=1$, in this proof, it must be understood that $\dd^{(1)}$ is the anchor $ \rho$. 
Assume that $\ker\dd^{(i)}_{|_{m}}=\mathrm{Im}\left(\dd^{(i+1)}_{|_m}\right)$ for some $m\in M$. This implies that \begin{align}\nonumber\mathrm{rk}\left(\dd^{(i+1)}_{|_m}\right)&=\dim\left(\ker \dd^{(i)}_{|_{m}}\right)\\&=\label{eq:pointwise-exact}\mathrm{rk}(E_{-i})-\mathrm{rk}\left(\dd^{(i)}_{|_m}\right).
    \end{align}
    It is a general fact that the functions from $M$ to $ \mathbb N$ given by
    $$m'\mapsto \mathrm{rk}\left(\dd^{(i)}_{|_{m'}}\right) \hbox{ and }  m' \mapsto \mathrm{rk}\left(\dd^{(i+1)}_{|_{m'}}\right)$$
    are lower semi-continuous functions. Hence, there exists an open neighborhood $\mathcal{U}$ of $m$ such that for all $m'\in \mathcal U$\begin{align*}
        \mathrm{rk}\left(\dd^{(i+1)}_{|_{m'}}\right)&\geq \mathrm{rk}\left(\dd^{(i+1)}_{|_{m}}\right)\\&= \mathrm{rk}(E_{-i})-\mathrm{rk}\left(\dd^{(i)}_{|_{m}}\right),\quad\text{by Equation \eqref{eq:pointwise-exact}}\\&\geq\mathrm{rk}(E_{-i})-\mathrm{rk}\left(\dd^{(i)}_{|_{m'}}\right),\quad\text{by lower semi-continuity of  $m' \mapsto \mathrm{rk}\left(\dd^{(i)}_{|_{m'}}\right)$}\\&=\dim\left(\ker \dd^{(i)}_{|_{m'}}\right).
    \end{align*}
    Since we are considering a complex of vector bundles, we also have $$\mathrm{Im}\left(\dd^{(i+1)}_{|_{m'}}\right)\subseteq\ker \dd^{(i)}_{|_{m'}}$$
    for all $m'\in M$. This contradicts the previous inequality, unless if $\ker\dd^{(i)}_{|_{m'}}$ coincides with $~\mathrm{Im}\left(\dd^{(i+1)}_{|_{m'}}\right)$ for all $m' \in \mathcal U$.  This proves the claim.
\end{proof}
\begin{definition}
Let $(E_{-\bullet},\dd^{(\bullet)}, \rho)$ and $(E_{-\bullet}',(\dd')^{(\bullet)},\rho')$ be anchored complexes of vector bundles.  

\begin{enumerate}
    \item 
An \emph{anchored chain map} or \emph{anchored complex of vector bundle morphisms} between the anchored complexes of vector bundles $(E_{-\bullet},\dd^{(\bullet)}, \rho)$ and $(E_{\bullet}',(\dd')^{(\bullet)},\rho')$ is a collection of vector bundle morphisms (of degree zero) $\varphi_\bullet\colon E_{-\bullet}\longrightarrow E'_{-\bullet}$ such that the following diagram commutes

\begin{equation}\label{eq:exa-bisub}
        \xymatrix{
\cdots\ar[r]&E_{-i} \ar[d]_{\varphi_{i}} \ar[r]^{\dd^{(i)} }&E_{-i+1}\ar[d]^{\varphi_{i-1}} \ar[r]&\cdots \ar[r]^{\dd^{(2)}}& E_{-1}\ar[r]^{\rho}\ar[d]^{\varphi_1}&TM\ar@{=}[d] \\
\cdots\ar[r]&E'_{-i} \ar[r]^{\dd'^{(i)}}&E'_{-i+1}\ar[r]&\cdots\ar[r]^{\dd'^{(2)}} &E_{-1}'\ar[r]^{\rho'}&TM }\end{equation}

that is, $\rho'\circ \varphi_{1}=\rho$\; and\; $\dd'^{(i)}\circ \varphi_{i}= \varphi_{i-1}\circ \dd^{(i)}$ for every $i\geq 2$.

\item A \emph{homotopy} between two anchored complexes of vector bundle morphisms $\varphi_\bullet,\psi_\bullet\colon E_{-\bullet}\longrightarrow E'_{-\bullet}$ is the datum $\{h_i\colon E_{-i}\longrightarrow E_{-i-1}'\}_{i\geq 1}$ of vector bundle morphisms (of degree $-1$) 
that satisfies $\psi_1-\varphi_1=\dd'^{(2)}\circ h_1$ and for each $i\geq 2$, $\psi_i- \varphi_i= \dd'^{(i+1)}\circ h_i +h_{i-1} \circ \dd^{(i)}$

\begin{equation}\label{eq:hot-exa-bisub}
        \xymatrix{
\cdots\ar[r]&E_{-i-1}\ar[d]_{\psi_{i+1}- \varphi_{i+1}}\ar[r]^{\dd^{(i+1)}}&E_{-i}\ar@{.>}[ld]_{h_i} \ar[d]<3pt>_{\psi_i-\varphi_i} \ar[r]^{\dd^{(i)} }&E_{-i+1}\ar[d]^{\psi_{i-1}-\varphi_{i-1}}\ar@{.>}[ld]_{h_{i-1}} \ar[r]&\cdots \\
\cdots\ar[r]&E'_{-i-1}\ar[r]^{\dd'^{(i+1)}}&E'_{-i} \ar[r]^{\dd'^{(i)}}&E'_{-i+1}\ar[r]&\cdots}
\end{equation}

\begin{enumerate}
\item When there exists a homotopy between two anchored complexes of vector bundle morphisms $\varphi_\bullet,\psi_\bullet\colon E_{-\bullet}\longrightarrow E'_{-\bullet}$, we say that $\varphi,\psi$ are \emph{homotopic}  and we write $\varphi \sim \psi$.
    \item $(E,\dd, \rho)$ and $(E',\dd',\rho')$  are said to be \emph{homotopy equivalent}, if there exist anchored chain maps $\varphi_\bullet\colon E_{-\bullet}\longrightarrow E'_{-\bullet}$ and $\psi_\bullet\colon E'_{-\bullet}\longrightarrow E_{-\bullet}$ such that $\varphi\circ \psi \sim \mathrm{id}_{E'_{-\bullet}}$ and $ \psi\circ\varphi \sim \mathrm{id}_{E_{-\bullet}}$.
\end{enumerate}
\end{enumerate}
\end{definition}

 \begin{exo}
     Check that \begin{enumerate}
         \item  $\sim$ is an equivalence relation on the class of complexes of vector bundle morphisms.
         \item ``homotopy equivalence'' on the class of anchored complexes of vector bundles is an equivalence relation.
     \end{enumerate}
 \end{exo}

 \begin{lemma}
 \label{lem:alternatesum}
Let $(E_{-\bullet},\dd^{(\bullet)}, \rho)$ and $(E_{-\bullet}',(\dd')^{(\bullet)},\rho')$ be  homotopy equivalent anchored complexes of vector bundles of finite length $n$ and $n'$ respectively. The alternating sum of the ranks of the vector bundles $(E_{-i})_{i\in\mathbb{N}}$ and $(E'_{-i})_{i\in \mathbb N}$ respectively,
 are equal, i.e., $$\sum_{i=1}^n(-1)^{i}\mathrm{rk}(E_{-i})=\sum_{i=1}^{n'}(-1)^{i}\mathrm{rk}(E'_{-i}).$$Here $\mathrm{rk}(E)$ stands for the rank of a vector bundle $E\to M$.
\end{lemma}

\begin{proof}
Note first that the restriction of both complexes to a point $m \in M$ give two finite length complexes of vector spaces of finite dimension. The result is an immediate consequence of the fact that in every degree the cohomology group of  two equivalent complexes of vector spaces are isomorphic. It follows by taking the alternating sum of their dimensions and using the Rank–nullity theorem.
\end{proof}

For the next definition, we invite the reader to consider singular foliations as defined in Definition \ref{def:consensus2} or \ref{def:consensus2alg}.

\vspace{0.5cm}

\begin{definitions}{Geometric resolution of a singular foliation}{defgeomresol}
  Let $\mathcal{F}\subseteq\mathfrak{X}(M)$ be a singular foliation on a smooth, real analytic or complex  manifold $M$. An anchored complex of vector bundles  $(E_{-\bullet},\dd^{(\bullet)}, \rho):=$
    $$   \xymatrix{  \ar[r] & E_{-i-1} \ar[r]^{{\dd^{(i+1)}}} \ar[d] & 
     E_{ -i} \ar[r]^{{\dd^{(i)}}}
     \ar[d] & E_{-i+1} \ar[r] \ar[d] & \ar@{..}[r] & \ar[r]^{{\dd^{(2)}}}& E_{-1} \ar[r]^{\rho} \ar[d]& TM \ar[d] \\ 
      \ar@{=}[r] & \ar@{=}[r] M  &  \ar@{=}[r] M 
      &  \ar@{=}[r] M  &\ar@{..}[r] & \ar@{=}[r]   &  \ar@{=}[r] M  & M}
      $$
 is said 
 \begin{enumerate}
 \item to \emph{terminate in $\mathcal F $} if $ \rho\left(\Gamma(E_{-1})\right) \subseteq\mathcal F$
     \item to be \emph{over $\mathcal F $} if $ \rho(\Gamma(E_{-1}))=\mathcal F$ 
     \item  to be a \emph{geometric resolution of  $\mathcal F$} if the following complex of sheaves\footnote{See discussion below.} is exact:
      $$ \cdots {\longrightarrow} \Gamma({E_{ -i-1}})
     \stackrel{\dd^{(i+1)}}{\longrightarrow} \Gamma({E_{-i}})
     \stackrel{\dd^{(i)}}{\longrightarrow}{\Gamma(E_{-i+1})}{\longrightarrow}\cdots  {\longrightarrow} \Gamma(E_{-1})  
     \stackrel{\rho}{\longrightarrow} \mathcal F\stackrel{}{\longrightarrow} 0.$$ 
 \end{enumerate}
 A geometric resolution $(E_{-\bullet},\dd^{(\bullet)}, \rho)$ is said to be \emph{minimal at a point}  $m\in M$ if for each $i\geq 2$ the linear map $\dd^{(i)}_{|_m}\colon {E_{-i}}_{|_m}\longrightarrow {E_{-i+1}}_{|_m}$ vanishes.
 \end{definitions}

\vspace{0.5cm}

 Let us recall what we mean precisely by the sheaf condition above, and explain how this condition simplifies in the smooth case. 
``Exact as sheaves'' means that for any $i \in \mathbb N $ and any $m \in M $, there is a neighborhood $\mathcal V $ of $m$ such that for any $\mathcal U $  in $\mathcal V $ the complex: 
 $$ \cdots {\longrightarrow} \Gamma_{\mathcal U}({E_{ -i-1}})
     \stackrel{\dd^{(i+1)}}{\longrightarrow} \Gamma_{\mathcal U}({E_{-i}})
     \stackrel{\dd^{(i)}}{\longrightarrow}{\Gamma_{\mathcal U}(E_{-i+1})}{\longrightarrow}\cdots  .$$
 is exact.
  In the smooth setting, this is equivalent to demand that the complex be exact at the level of global compactly supported sections. 
In the smooth setting, therefore, the notion of geometric resolution is much easier, and Definition \ref{def:defgeomresol} can be rewritten as follows. 

\vspace{.5cm}

\begin{definitions}{Geometric resolution of a smooth singular foliation, an equivalent definition.}{}
  Let $\mathcal{F}\subseteq \mathfrak{X}_c(M)$ be a singular foliation on a smooth manifold $M$. A \emph{geometric resolution} of the singular foliation $\mathcal{F}$ is a complex of vector bundles  $(E_{-\bullet},\dd^{(\bullet)}, \rho):=$
    $$   \xymatrix{  \ar[r] & E_{-i-1} \ar[r]^{{\dd^{(i+1)}}} \ar[d] & 
     E_{ -i} \ar[r]^{{\dd^{(i)}}}
     \ar[d] & E_{-i+1} \ar[r] \ar[d] & \ar@{..}[r] & \ar[r]^{{\dd^{(2)}}}& E_{-1} \ar[r]^{\rho} \ar[d]& TM \ar[d] \\ 
      \ar@{=}[r] & \ar@{=}[r] M  &  \ar@{=}[r] M 
      &  \ar@{=}[r] M  &\ar@{..}[r] & \ar@{=}[r]   &  \ar@{=}[r] M  & M}
    $$

 such that  the following complex  is exact:
      $$ \cdots {\longrightarrow} \Gamma_c({E_{ -i-1}})
     \stackrel{\dd^{(i+1)}}{\longrightarrow} \Gamma_c({E_{-i}})
     \stackrel{\dd^{(i)}}{\longrightarrow}{\Gamma_c(E_{-i+1})}{\longrightarrow}\cdots  {\longrightarrow} \Gamma_c(E_{-1})  
     \stackrel{\rho}{\longrightarrow} \mathcal F  .$$ 
 \end{definitions}

\vspace{.5cm}

  \begin{remark}
  \label{rmk:minimalcohom}
 When a geometric resolution $(E_{-\bullet},\dd^{(\bullet)}, \rho)$ is minimal at a point  $m\in M$ then one has, $H^{-i}(E_{\bullet},m)={E_{-i}}_{|_m}$  for all $i\geq 2$.
 \end{remark}

\begin{remark}
A  singular foliation $\mathcal{F}$ is Debord if and only if there exists a geometric resolution of length $1$.
\end{remark}

\subsubsection{Universality}

We conclude this section with a theorem which says that, given a foliated manifold  $ (M,\mathcal F)$, in the category where
\begin{enumerate}
\item objects are anchored complexes of vector bundles that terminate within $\mathcal F $.
\item arrows are homotopy classes of morphisms,
\end{enumerate}
geometric resolutions are terminal\footnote{Also called “final”: it is an object $\gamma$ such that for any other object $ \alpha$ there is one and only one arrow from $ \alpha$ to $\gamma $. It is a general property of category theory that for any two terminal (a.k.a. final, a.k.a. universal) objects  $\gamma,\gamma' $ there  exist a unique arrow $\gamma \to \gamma' $ and $ \gamma' \to \gamma$ which are inverse one to the other.} (a.k.a. “universal”) objects.  For the reader not familiar with categories, this gives the following result, that extends Theorems \ref{thm:prop:UniqueAnchored}.

\vspace{.5cm}

\begin{theorems}{Geometric resolutions are terminal objects}{uniquegeomresol}
Let $\mathcal{F}\subseteq \mathfrak{X}(M)$ be a singular foliation on a smooth manifold $M$ that admits a geometric resolution $(E_{-\bullet},\dd^{(\bullet)}, \rho)$. 
\begin{enumerate}
    \item For any anchored complex of vector bundles $(E_{-\bullet}', \dd^{(\bullet)}, \rho')$ that terminates within $\mathcal F$,  there exists a chain map of anchored vector
bundles 
$$(E_{-\bullet}', \dd^{(\bullet)}, \rho') \longrightarrow (E_{-\bullet
}, \dd^{(\bullet)}, \rho)$$ and any two such chain maps are homotopy equivalent.
    \item In particular, 
two geometric resolutions of the same singular foliations are homotopy equivalent. 
\end{enumerate}

The same results hold in the complex and real analytic setting, but in a neighborhood of a point only.

\end{theorems}

\begin{proof}
These are reinterpretations of classical results of algebraic topology, see Section 3.2 in \cite{LLS} for more explanations: for any commutative algebra $ \mathcal O$, resolutions  of an $ \mathcal O$-module, say $ \mathcal F$, by projective $ \mathcal O$-modules are  universal objects in the category where objects are complexes of $ \mathcal O$-modules ending in  $ \mathcal F$ and arrows are homotopy classes of chain maps.
\end{proof}

The fact that two geometric resolutions of $\mathcal{F}$, when they exist, are homotopy equivalent, has many consequences on the topic ``whatever is canonically invariant under homotopy equivalence is canonically attached to the singular foliation - provided it admits geometric resolutions.'' 
It is the case, for instance, of the alternating sums of the ranks by Lemma \ref{lem:alternatesum}. 

\begin{corollary}\label{prop:alternatesum} [{\cite{LLS}, Proposition 2.5}]
Let $\mathcal{F}\subseteq \mathfrak{X}(M)$ be a singular foliation on a smooth, real analytic or complex connected manifold $M$ that admits a geometric resolution $(E_{-\bullet},\dd^{(\bullet)}, \rho) $ of finite length. Then   
    \begin{enumerate}
        \item all the regular leaves have the same\footnote{This point is in fact automatic in the complex and real analytic settings, it is only non-trivial in the smooth case, see Remark \ref{rmk:lowercontinuous}.} dimension $r$,
        \item the alternating sum of the ranks of $ E_{-\bullet}$ is equal to the dimension of the regular leaves, i.e.
        $$r=\sum_{i\geq 1}(-1)^{i+1}\mathrm{rk}(E_{-i}).$$
    \end{enumerate}
    \end{corollary}

If two geometric resolutions are homotopy equivalent,  their restrictions to a point $m \in M$ are also homotopy equivalent. In consequence, the complexes obtained by evaluation at this point have isomorphic cohomologies.
This proves the first part of the following corollary, we leave the rest to the reader.

\begin{corollary}
\label{cor:canonicalpoint}

Let $\mathcal{F}\subseteq \mathfrak{X}(M)$ be a singular foliation on a smooth, real analytic or complex manifold $M$ that admits a geometric resolution. Then, for every $m \in M$, the cohomologies $ H^{-i}(E_{-\bullet}, m)  \simeq H^{-i}(E_{-\bullet}', m)  $ are canonically isomorphic.  
\end{corollary}

It therefore makes sense to denote these spaces of cohomologies by $ H^{-i}(\mathcal F,m)$, erasing the dependency\footnote{Some readers may recognize that $H^{-i}(E_{-\bullet}, m)$ is  the so-called Tor-functor (see \cite{zbMATH00595200}) in the category of $ \mathcal O$-modules, with $\mathcal O $ the sheaf of functions on $M$, applied to the $\mathcal O $-modules $ \mathcal F$ and  $\mathbb K$. On the second one, the action is given by $ F \cdot \lambda = F(m) \lambda$ for all $F \in \mathcal O $ and $\lambda \in \mathbb K$.} on the geometric resolution. We call them \emph{isotropy spaces of $ \mathcal F$ at $m$}. In particular, the dimensions $ d_1, \dots, d_i, \dots$ 
of these spaces are canonically attached to $\mathcal F $. Also, Proposition 
\ref{prop:atapoontzeroaroundzero}  implies that the following items are equivalent:
\begin{enumerate}
    \item $m$ is a regular point,
    \item $H^{-1}(\mathcal F,m)=0 $
    \item $H^{-i}(\mathcal F,m')=0 $ for every $i \geq 1$ and every $m'$ in a neighborhood of $m$.
\end{enumerate}

\begin{remark}
The integers $ d_1, \dots, d_i, \dots$  were constructed without making any use of the Lie bracket of vector fields, so that they are, as a matter of fact, attached to $\mathcal F $ seen as a module over functions, and not to $\mathcal F $ seen as a singular foliation. 
We suggest interpreting them as follows:
\begin{enumerate}
    \item $ d_1$ is the rank of $\mathcal F $  at $m$ minus the dimension of the leaf through $m$, 
    \item $ d_2$ is the rank at $m$ of the  module of relations between the previous generators near $m$,
    
    \item $ d_3$ is the rank of the module over relations between the generators of the module of relations between the generators of $ \mathcal F$,
    
    \item … and so on.
\end{enumerate}
\end{remark}

\subsection{Existence of geometric resolutions (Noetherianity, Syzygies and a flatness theorem by Tougeron-Malgrange)} 

\label{sec:existencegeometric}

We defined and studied geometric resolutions in Section \ref{sec:geometricresolutionsdef}, but do they exist?

Here are some cases where geometric resolutions of a singular foliation  always exist at least locally, and are of
finite length.

\vspace{0.5cm}

\begin{propositions}{A particular case of Syzygy theorem}{prop:Syzygy}
    Let $ \mathcal F$ be a real analytic or holomorphic singular foliation on a manifold $M$ of dimension $d$.
    Any point in $M$ admits a neighborhood on which a geometric resolution by trivial vector bundles exists.
    
    Moreover, its length can be chosen to be less than $ d+1 $.
\end{propositions}

\begin{proof}
This is in fact a general theorem for coherent sheaves on a complex or real analytic manifold, and has nothing to do with singular foliation. It is called Syzygy theorem, and is related to Oka's coherence theorem. We invite the reader to refer to the classical literature on the matter, e.g.,
Griffiths and Harris's \emph{Principles of Algebraic Geometry} \cite{zbMATH03634395} page 696.

\end{proof}

\begin{remark} \label{rmk:algebraic}
For an algebraic singular foliation on an affine variety $W$ as in Section \ref{sec:affinevarieties} over $ \mathbb C$ or $ \mathbb R$, since the ring of functions $\mathcal O_W $ is Noetherian, a geometric resolution of $ \mathcal F$ by trivial vector bundles has to exist, for the following reason:
\begin{enumerate}
    \item $ \mathcal F$ being finitely generated over $ \mathcal O_W$ by definition, there exists an integer $ d_1$ and a surjective $ \mathcal O_W$-linear map $ \rho : \mathcal O_W^{d_1} \to \mathcal F$.
    \item The kernel of $ \rho$ being an $ \mathcal O_W$-module, by Noetherianity, there exists an integer $ d_2$ and a surjective $ \mathcal O_W$-linear map $ \dd^{(2)} : \mathcal O_W^{d_2} \to \mathrm{Ker}(\rho)$.
    \item The kernel of $\dd^{(2)}$ being an $ \mathcal O_W$-module, by Noetherianity, there exists an integer $ d_3$ and a surjective $ \mathcal O_W$-linear map $ \dd^{(3)} : \mathcal O_W^{d_2} \to \mathrm{Ker}(\dd^{(2)})$
    \item ...and so on up to infinity...
\end{enumerate}
Now, the module $\mathcal O_W^{d_i}$ can be considered as a module of sections of a trivial vector bundle  $ E_{-i}$, and the maps $ \rho$ and $ \dd^{(\bullet)}$ can be considered a vector bundle morphisms. This completes the construction.
In general, the procedure described above never stops. But for $ W =\mathbb C^d$, we can make it stop after at most $ d+1$ steps. This is a deep result in algebraic geometry called Syzygy theorem (See, e.g., Theorem 1.1. in \cite{zbMATH02136428}). 
The same results hold for $ W =\mathbb R^N$, and for a polynomial singular foliation on $W$. See \cite{LLS}, Section 3.2, for more details.
\end{remark}

\vspace{0.5cm}

Here is a technical but important question. 
Assume that we are given a real analytic foliated manifold  $ (M,\mathcal F_{ra})$
a geometric resolution $(E_{\bullet}, \dd^{(\bullet)},\rho)$ of $ \mathcal F_{ra}$.
The real analytic manifold can be considered as a smooth manifold. The singular foliation can be considered as a smooth singular foliation $ \mathcal F_{sm}$ (it suffices to consider the sheaf of $ \mathcal C^\infty_c(M)$-module  generated\footnote{If the reader prefers to think in terms of global section compactly supported vector fields as in Definition \ref{def:consensus}, $\mathcal F_{sm} $ must be seen as compactly supported vector fields $X$ such that every $m\in M$ admits a neighborhood $\mathcal U  $ on which $ X =\sum_{i=1}^k f_i X_i$ with $ f_i \in \mathcal C^\infty(\mathcal U)$ and $ X_i$ a local section of $ \mathcal F_{ra}$.} by $ \mathcal F$): 
 $$  \mathcal F_{sm} := \mathcal C^{\infty}(M) \, \mathcal F_{ra}  $$
Moreover, a real analytic geometric resolution $ (E_{-\bullet},\dd^{(\bullet)}, \rho) $, being a family real analytic vector bundles equipped with a family real analytic vector bundle morphisms,  can be considered as a smooth complex of vector bundles over $\mathcal F_{sm} $:
But it is \underline{not} obvious that it is still a geometric resolution of $ \mathcal F_{sm}$. 
It is not obvious that exactness of  the complex of sheaves\footnote{$\Gamma_{ra} $ stands for the sheaf of real analytic sections.}
$$ \cdots {\longrightarrow} \Gamma_{ra}({E_{ -i-1}})
     \stackrel{\dd^{(i+1)}}{\longrightarrow} \Gamma_{ra}({E_{-i}})
     \stackrel{\dd^{(i)}}{\longrightarrow}{\Gamma_{ra}(E_{-i+1})}{\longrightarrow}\cdots  {\longrightarrow} \Gamma_{ra}(E_{-1})  
     \stackrel{\rho}{\longrightarrow} \mathcal F_{ra}  $$
     implies the exactness of the complex of sheaves\footnote{$\Gamma_{sm} $ stands for the sheaf of smooth sections.}
$$ \cdots {\longrightarrow} \Gamma_{sm}({E_{ -i-1}})
     \stackrel{\dd^{(i+1)}}{\longrightarrow} \Gamma_{sm}({E_{-i}})
     \stackrel{\dd^{(i)}}{\longrightarrow}{\Gamma_{sm}(E_{-i+1})}{\longrightarrow}\cdots  {\longrightarrow} \Gamma_{sm}(E_{-1})  
     \stackrel{\rho}{\longrightarrow} \mathcal F_{sm}  .$$
Notice that this is absolutely not true if we deal with global sections. If true, it can only be true at the level of sheaves, i.e., in a neighborhood of a point.
     Such results are called flatness. What we are asking is: are smooth functions flat over real analytic functions? The answer is “no”, and counter-examples are easily found. However, it is true\footnote{See
     Theorem 4 in \cite{MR240826} by Tougeron - a result generally attributed to Malgrange.} that \emph{germs} of smooth functions at a point form a flat module over converging real analytic functions.
This is enough for our present purpose: because it means that given a smooth section $a \in \Gamma_{sm}(E_{-i-1}) $ which is $\dd^{(i)}$-closed, defined on some open subset $ \mathcal U$ of $M$, each point   $m \in \mathcal U$ admits a neighborhood $\mathcal V_m $ where it is exact $ a_{|_{\mathcal V_m}} = \dd^{(i+1)} b_{\mathcal V_m}$. Using an open cover and partitions of unity\footnote{See the proof of Proposition 2.3 in \cite{LLS}.}, we can glue the  $b_{\mathcal V_m}$'s to obtain some $b \in \Gamma_{sm}( E_{-i-2}) $ such that $\dd^{(i+1)} b =a$.  Hence:

\begin{proposition}[\cite{LLS}, Proposition 2.3]
\label{prop:flatness}

A real analytic geometric resolution of a real analytic singular foliation $\mathcal F_{ra} $ is also a smooth geometric resolution of its induced smooth singular foliation $ \mathcal F_{sm}$. 

\end{proposition}

\begin{remark}\label{rmk:fromplynomialtosmooth}
In a similar fashion, an algebraic singular foliation on $ \mathbb R^d$ can be seen as a smooth singular foliation\footnote{See Proposition \ref{prop:fromalgtoholo}}. 
Any geometric resolution as in Remark \ref{rmk:algebraic} can be seen as a smooth geometric resolution.
 This is follows from the consecutive use of
 \begin{enumerate}
 \item flatness\footnote{Polynomial functions form a valuation ring. A module over a valuation ring is flat if and only if it is torsion free (see, e.g., stack project, section 15.22), which is the case of real analytic functions on an open subset of $ \mathbb R^d$.} of  real analytic functions over the ring of polynomials, which implies that a polynomial geometric resolution is a real analytic geometric resolution,
 \item then Proposition \ref{prop:flatness}, which says that any real analytic geometric resolution is also a smooth geometric resolution.
 \end{enumerate}
For similar reasons,   an algebraic singular foliation on $ \mathbb C^d$ can be seen as a holomorphic singular foliation\footnote{See the previous footnote: the argument extends to holomorphic function on any open connected subset of $ \mathbb C^d$.}. 
\end{remark}

We say that a smooth singular foliation is \emph{locally real analytic} if every point admits a local chart on which
the generators of the singular foliation are real analytic, i.e., have real analytic coefficients. We \emph{do not} assume that these charts are patched together by real analytic transition function: it is fine if those are smooth only. This is enough to guarantee the existence of real analytic geometric resolution near every point, by using the Syzygy-type argument of Proposition \ref{thm:prop:Syzygy}. Then these real analytic geometric resolutions can be considered as smooth ones, see Proposition \ref{prop:flatness}.  Now, in the smooth context, and only in this context, smooth geometric resolutions can be “glued” - in some sense, and after modification to homotopy equivalent ones: we refer to the discussion of Proposition 3.24 in \cite{LLS}.
Eventually, one obtains the following theorem (Theorem 2.4 in \cite{LLS}).

\vspace{.5cm}

\begin{theorems}{Existence of geometric resolutions}{th:existence}
A locally real analytic singular foliation on a manifold of dimension $d$ admits a geometric resolution of length $\leq d+1$ on any relatively compact open subset of $M$.
\end{theorems}

\vspace{.5cm}

Below, we describe some examples of geometric resolutions of singular foliations.

\begin{exo}
Show that the singular foliation on  $M=\mathbb R $ generated by the vector field $ \chi(t) \partial_t $ with $ \chi(t)= e^{-1/t^2}$ for $ t > 0$ and $ \chi(t)=0$ for $ t \leq 0$ admits no geometric resolution. (\emph{Hint}: Use Corollary~\ref{prop:alternatesum}).
\end{exo}

\begin{example}
\label{ex:fromideal_to_geometricresol}
Let $M$ be a smooth, real analytic or complex manifold, with sheaf of functions $ \mathcal O$.
Let $ \mathcal I \subset \mathcal O$ be a sheaf of ideals. 
We say that $ \mathcal I$ admits a geometric resolution when there exists a complex of vector bundles 
$$  \xymatrix{\dots \ar[r]^{\dd} & I_{-2} \ar[d]\ar[r]^{\dd} & I_{-1}   \ar[d]\ar[r]^{\epsilon \hspace{0.2cm}}& \mathbb K \times M \ar[d]\\ \dots \ar[r]^{=} & M \ar[r]^{=}& M \ar[r]^{=}& M    } $$
which are exact at the level of sections and such that the image of $ \epsilon $ is $ \mathcal I$.
Now, $ \mathcal I \mathfrak X$, i.e., the sheaf of vector fields whose coefficients are in  $ \mathcal I$, form a singular foliation. A geometric resolution of it is given by taking the tensor product\footnote{Notice the use of $(\mathbb K \times M) \otimes TM\simeq TM$ in the last column.} with $ TM$, as follows: 
$$  \xymatrix{\dots \ar[r]^{\dd \otimes {\mathrm{id}}\hspace{0.5cm}} & I_{-2} \otimes TM\ar[d]\ar[r]^{\hspace{0.1cm} \dd \otimes {\mathrm{id}}} & I_{-1}  \otimes TM \ar[d]\ar[r]^{\hspace{0.5cm}\epsilon \otimes {\mathrm{id}}}& TM \ar[d]\\ \dots \ar[r]^{=} & M \ar[r]^{=}& M \ar[r]^{=}& M    } $$
This type of geometric resolution is in fact quite common, as we now see in \ref{ex:nullatzero} and \ref{ex:completeintersection}.
\end{example}

\begin{example}
\label{ex:nullatzero}
Let $\mathcal F_0=\{X\in\mathfrak X(V)\mid X(0)=0\}$ be the singular foliation made of all vector fields vanishing at the origin of a vector space $V$ of dimension $N$ over $\mathbb  C$ or $\mathbb R$). 
We are precisely in the situation of Example \ref{ex:fromideal_to_geometricresol}, with $ \mathcal I$ the ideal of polynomial functions vanishing at $0$.

In view of the previous example, therefore, we have to find a geometric resolution of that ideal, then take the tensor product with the tangent bundle.
This can be done as follows. The contraction by the Euler vector field $$\displaystyle{\overrightarrow{E}=\sum_{i=1}^N x_i\frac{\partial}{\partial x_i}}$$
gives rise to a complex of trivial vector bundles 
\begin{equation}\label{eq:Koszul1}
 \quad  \cdots {\longrightarrow} \wedge^3T^*V  \stackrel{\iota_{\overrightarrow{E}}}{\longrightarrow}\wedge^2T^*V   \stackrel{\iota_{\overrightarrow{E}}}{\longrightarrow}T^*V\stackrel{\iota_{\overrightarrow{E}}}{\longrightarrow} \mathbb{C}\times V=:\underline{\mathbb C},\end{equation} whose complex on the level of sections is $(\Omega^\bullet(V), \iota_{\overrightarrow{E}})$. Here $(x_1,\ldots, x_N )$ are the
canonical coordinates on $V$.  Since the $x_i$'s form a regular sequence, it is well known that $(\Omega^\bullet(V), \iota_{\overrightarrow{E}})$ is exact, and the previous complex is therefore a geometric resolution of the ideal of polynomial functions that vanish at $0$.
A geometric resolution of $ \mathcal F_0$ is therefore given by the trivial bundles with fiber
 $E_{-i} :=  (\wedge^{i} V^* \otimes V ) \times V \to V $, the differential 
  $ \iota_{\overrightarrow{E}}\times {\mathrm{id}}$, and the anchor given at a point $v \in V$:
  $$ \begin{array}{rcl}  E_{-1}|_{v} \simeq V^* \otimes V & \to & V \simeq T_v V\\ \alpha \otimes u & \to &(\iota_{\overrightarrow{E}_{|_v}} \alpha) u = \alpha(v) \, u.\end{array}$$
\end{example}

\begin{example}
\label{ex:completeintersection}
Example \ref{ex:nullatzero} can be further extended.
 Consider the singular foliation $\mathcal I \mathfrak X(V) $ of all  vector fields on a finite dimensional vector space $V$ over $\mathbb K=\mathbb R $ or $ \mathbb C$ whose coefficients\footnote{For $ \mathbb K=\mathbb C$, when $ \mathcal I$ is the ideal of functions on an affine variety $W$, this is precisely the singular foliation of vector fields that vanish on $W$, see Section \ref{sec:algebraiccase_subvariety}} are in $ \mathcal I$. 
 When the ideal $ \mathcal I$ is generated by a regular sequence $ \phi_1, \dots, \phi_a$ of polynomial functions, then it admits a Koszul resolution\footnote{See Section 4.5 in \cite{zbMATH00595200}, or the original 
 (in French) \cite{zbMATH03058924}}. This is constructed as follows. Let $B $ be a 
vector space of dimension $a$, and $\alpha $   be the section of the trivial vector bundle $ B^* \times V \to V$ given by 
$$ \alpha = \sum_{i=1}^a \phi_i \, e_i^* $$
with $ e_1^*, \dots, e_a^*$ a basis of $V^*$.
Consider the trivial vector bundle $ \wedge^\bullet  B \times V \to V$ with fiber $ \wedge^\bullet  B $. The pair  $ \left( \Gamma(\wedge^\bullet  B \times V ) ,\iota_\alpha\right)$ is a complex of vector bundles\footnote{In degree $ -i$, the vector bundle is $\wedge^i B \times V \to V $.}. Koszul theorem states that this construction is a geometric resolution of the ideal $ \mathcal I$, which in turn gives a geometric resolution of $ \mathcal I \mathfrak X(V)$ by taking the tensor product with $TV$ as in Example \ref{ex:fromideal_to_geometricresol}.
\end{example}

\
\begin{example}

Here is a case of singular foliation on $\mathbb K^d $ which is of the form $ \mathcal I \mathfrak X(V)$ for some vector space $V$ and some ideal $ \mathcal I$, but the construction of Example \ref{ex:completeintersection} does not apply because $ \mathcal I$ is not generated by a regular sequence.
Let $\mathcal{F}_2$ be the  algebraic singular foliation made of vector fields vanishing at order at least $2$ at the origin of $V=\mathbb{K}^2$ (see Section \ref{ex:singFolVanish}). Here  $\mathcal{I}$ is the ideal generated by the monomials $x^2, xy, y^2$. It is \emph{not} generated by a regular sequence.
However, it admits a geometric resolution  \begin{equation}\label{eq:ideal-reso1}
    0\longrightarrow \mathcal{O}(\mathbb{K}^2)\oplus\mathcal{O}(\mathbb{K}^2)\stackrel{d}{\longrightarrow} \mathcal{O}(\mathbb{K}^2)\oplus\mathcal{O}(\mathbb{K}^2)\oplus\mathcal{O}(\mathbb{K}^2)\stackrel{\epsilon}{\longrightarrow}\mathcal{I}\longrightarrow 0,
\end{equation}where for all $f,g,h\in\mathcal{O}(\mathbb{K}^2)$, $$\epsilon (f,g,h)=x^2f+xyg+y^2h\,\;\text{and}\,\; d(f,g)=(-y f,xf-yg,xg).$$
Here $ \mathcal O(\mathbb K^2)$ stands for polynomial functions on $ \mathbb K^2$. This corresponds to a complex of trivial vector bundles 
\begin{equation}
    0\longrightarrow \mathbb K^2 \times V \stackrel{d}{\longrightarrow} \mathbb K^3 \times V \stackrel{\epsilon}{\longrightarrow} \mathbb K \times V.
\end{equation}
Taking the tensor product with $TV$, which is a trivial bundle of rank $2$, we see that a geometric resolution is therefore given by the following trivial bundles:
\begin{equation}
    0\longrightarrow \mathbb K^2 \otimes \mathbb K^2 \times V \stackrel{d \otimes \mathrm{id}}{\longrightarrow} \mathbb K^3 \otimes \mathbb K^2 \times V  \stackrel{\epsilon \otimes \mathrm{id}}{\longrightarrow} \mathbb K \otimes \mathbb K^2 \simeq TV.
\end{equation}
 This geometric resolution can possibly be seen as a real analytic or a smooth singular foliation, as in Remark \ref{rmk:fromplynomialtosmooth}, when $ \mathbb K=\mathbb  R $ or a holomorphic one when $ \mathbb K=\mathbb  C$.

\end{example}

\begin{example}\label{exo:phi-Koszul}
Let $\varphi \in \mathcal O(V)$ be a polynomial function on a finite dimensional vector space $V$ over $ \mathbb R$ or $\mathbb C$. Consider the algebraic singular foliation on $ V$ of all polynomial vector fields that “kill”\footnote{This algebraic singular foliation was considered in Section \ref{sec:algebraiccase_subvariety}}. $ \varphi$:
$$  
\mathcal{F}_\varphi:=\{X\in\mathfrak{X}(V)\mid X[\varphi]=0\}.$$
The contraction by the  exact $1$-form $\dd\varphi \in \Omega^1(V)$ turns  sections\footnote{I.e., polynomial polyvector fields,  which here identify to polynomial maps from $V$ to $\wedge^\bullet V $.} of $ \wedge^\bullet TV$,
over $V$ into a complex of vector bundles 

\begin{equation}
    \label{eq:KoszulComplex1}
\cdots\stackrel{\iota_{\dd\varphi}}{\longrightarrow}\wedge^3 TV \stackrel{\iota_{\dd\varphi}}{\longrightarrow}\wedge^2 TV\stackrel{\iota_{\dd\varphi}}{\longrightarrow} TV \stackrel{\iota_{\dd\varphi}}{\longrightarrow}\mathbb C\times V.\end{equation}

Consider the complex on the level of sections:
\begin{equation} \label{eq:KoszulComplex2}\ldots\xrightarrow{\iota_{\dd\varphi}}\mathfrak X^3 (V) \xrightarrow{\iota_{\dd\varphi}}\mathfrak X^2 (V)\xrightarrow{\iota_{\dd\varphi}}\mathfrak X(V) \xrightarrow{\iota_{\dd\varphi}}\mathcal{O}(V)\end{equation}
(where  $\mathfrak X^i(V):= \Gamma(\wedge^i T V )$ stands for  $i$-vector fields on $V$).
Koszul theorem\footnote{See Section 4.5 in \cite{zbMATH00595200}, or the original 
 article \cite{zbMATH03058924}} states that this complex 
is exact in all degree, except in degree $0$, 
if $\left(\frac{\partial\varphi}{\partial x_1} ,\cdots,\frac{\partial\varphi}{\partial x_N}\right) $ is a regular sequence. The theorem of Cohen-Macaulay\footnote{See e.g., \cite{zbMATH00464113} for detailed proofs or \cite{zbMATH05043646}, Section 3, for a quick overview.} then says that $\left(\frac{\partial\varphi}{\partial x_1} ,\cdots,\frac{\partial\varphi}{\partial x_N}\right) $ is a regular sequence if $ \phi$ is weight homogeneous of non-zero degree with one isolated singularity at the origin. In all these cases, the truncated complex: 
\begin{equation} \label{eq:KoszulComplex3}\cdots\xrightarrow{\iota_{\dd\varphi}}\mathfrak X^3 (V) \xrightarrow{\iota_{\dd\varphi}}\mathfrak X^2 (V)\xrightarrow{\iota_{\dd\varphi}}\ker(\iota_{\dd{\varphi}})=\mathcal{F}_\varphi\end{equation}
is  an algebraic geometric resolution of $\mathcal{F}_\varphi$. It can therefore be seen as a real analytic or a smooth singular foliation, as in Remark \ref{rmk:fromplynomialtosmooth}.

In the case $ \mathbb K = \mathbb C$, it can be seen as a holomorphic one, in view of the same remark. 
\end{example}

\subsection{Geometric resolutions of length $\leq 2$ and singular foliations}
\label{sec:length2}

As we saw in Section \ref{sec:Debord}, a singular foliation $\mathcal F$ admits a geometric resolution $ 0 \to E_{-1} \to TM$ of length $1$ if and only if it is Debord.  $ E_{-1}$ then always admits an almost Lie algebroid bracket $ [\cdot\, , \cdot]_{E_{-1}}$ which is automatically a Lie algebroid bracket. In conclusion, if a geometric resolution of length $1$ exists, then it admits a Lie algebroid structure.
  
In this section, we discuss the case where a singular foliation $\mathcal{F}$ admits a geometric resolution of length $2$. In this case, we claim that there are Lie algebra-like structures on it. This will be generalized in Section \ref{sec:Universal}.

Let $(M,\mathcal{F})$ be a singular foliation that admits a geometric resolution of length $2$, namely 
 \begin{equation}\label{eq:length2}
(E_{-\bullet},\dd^{(\bullet)}, \rho): \quad    0{\longrightarrow} E_{-2}  \stackrel{\dd^{(2)}}{\longrightarrow}E_{-1}   \stackrel{\rho}{\longrightarrow}TM.
\end{equation}
Since Equation \eqref{eq:length2} is a geometric resolution of $\mathcal{F}$, the pair $(E_{-1}\rightarrow M, \rho)$ is an anchored bundle over~$\mathcal F$.
Therefore, by Proposition~\ref{prop:Almost-existence}, item 1, in the smooth case,  $ E_{-1}$  can be endowed with an almost Lie algebroid structure $(E_{-1},[\cdot\,,\cdot]_{E_{-1}},\rho)$, and in the real analytic or complex cases, the almost Lie algebroid structure exists in a neighborhood of any point. 

Let us assume from now that we are given on an almost Lie algebroid structure on $E_{-1} $. It is quite judicious to ask whether we can extend this bracket to sections of  $E_{-2}$. If yes, what kind of structures will we have?

Since the complex \eqref{eq:length2} is a geometric resolution of $\mathcal F$, the complex 
\begin{equation}\label{complex:sections-level}
    0{\longrightarrow} \Gamma(E_{ -2})_U \stackrel{\dd^{(2)}}{\longrightarrow}\Gamma(E_{ -1})_U\stackrel{\rho}{\longrightarrow}\mathcal{F}_U\longrightarrow 0
\end{equation}
is exact for all open subsets $U\subset M.$

\begin{enumerate}
    \item For all $ a \in \Gamma_U(E_{-1}), b \in \Gamma_U({E_{-2}})$:
\begin{equation*}
        \rho([a,\dd^{(2)} b ]_{E_{-1}})=[\rho(a),\rho\circ\dd^{(2)}b]=0,\; \hbox{(by the anchor condition\footnote{I.d. Equation \eqref{eq:bracket-morphism}}, and since   $\rho\circ\dd^{(2)}\equiv 0$).}
    \end{equation*}
In other words, the bracket $[a,\dd^{(2)} b]_{E_{-1}}$ is in the kernel of the anchor map $ \rho$. By exactness of the complex \eqref{complex:sections-level}, there exists a unique local section denoted by $\nabla_{a}b\in \Gamma(E_{ -2})_U$ such that \begin{equation}\label{eq:nabla}
    \dd^{(2)}(\nabla_{a}b)=[a,\dd^{(2)}b]_{E_{-1}}.
\end{equation}
The  bilinear map:
     $$   \begin{array}{lll} \Gamma(E_{-1})_U \otimes \Gamma( E_{-2})_U &\to & \Gamma(E_{-2})_U \\ \hspace{0.8cm}(a,b) & \mapsto & \nabla_a b \end{array}$$
     does not depend on the chosen open subset $ U$, i.d. is globally defined on $M$, and satisfies:
     \begin{enumerate}
         \item $\dd^{(2)}\nabla_{a}b=[\dd^{(2)}a,b]_{E_{-1}}$,\;$\forall a\in \Gamma(E_{ -2})_U,\, b\in  \Gamma(E_{ -1})_U$, by construction, and
         \item for all function $f\in \mathcal{O}(U)$:
      $ \nabla_a (f b) =  f \nabla_a b + \rho(a)[f] \, b$ and $ \nabla_{fa} b = f  \nabla_{a}b$, \text{for  all}\, $a \in  \Gamma(E_{-1}), b \in \Gamma(E_{-2})$,
\end{enumerate}
\item For all $a,b,c \in \Gamma(E_{-1})$, we have $$\mathrm{Jac} \,(a,b,c):=[a,[b,c]_{-1}]_{E_{-1}} +[b,[c,a]_{E_{-1}}]_{E_{-1}}+[c,[a,b]_{E_{-1}}]_{E_{-1}}\in \ker \rho.$$ By using again exactness of the complex \eqref{complex:sections-level} there is a unique local section that we denote by $[a,b,c]_{E_{-1}}\in\Gamma(E_{-2})_U$ that satisfies \begin{equation}
    \dd^{(2)}[a,b,c]_{E_{-1}}=\mathrm{Jac}(a,b,c).
\end{equation} 
It is easily checked that the map $ (a,b,c) \mapsto [a,b,c]_{E_{-1}} $ is linear over functions, and therefore comes from a pointwise linear vector bundle morphism $ \wedge^3 E_{-1} \to E_{-2}$. In particular, it is globally defined on $M$.
\end{enumerate}

The following lemma recapitulates the discussion above.

\begin{lemma}\label{lem:2-Lie-algebroid}
Let $(M,\mathcal{F})$ be a singular foliation that admits a geometric resolution of length $2$ as in \eqref{eq:length2} such that $(E_{-1},[\cdot\,,\cdot]_{E_{-1}},\rho)$ admits an almost Lie algebroid bracket. 
\begin{enumerate}
    \item There is 
     a bilinear map\footnote{Here, $ \Gamma(E)$ refers to local sections of the vector bundle $E$. In the smooth case, it suffices of course to work with global sections.}:
     $$   \begin{array}{lll} \Gamma(E_{-1}) \otimes \Gamma( E_{-2} ) &\to & \Gamma(E_{-2}) \\ \hspace{0.8cm}(a,b) & \mapsto & \nabla_a b \end{array}    $$ 
     \item[] and a skew-symmetric trilinear map:
     $$[\cdot\,,\cdot\,,\cdot]_{E_{-1}}\colon\Gamma(E_{-1})\wedge\Gamma(E_{-1})\wedge\Gamma(E_{-1})\longrightarrow \Gamma(E_{-2})$$ 
     \item such that  for every function $f$:
     \begin{enumerate}
         \item  $ \nabla_a f b =  f \nabla_a b + \rho(a)[f] \, b$ and $ \nabla_{fa} b = f  \nabla_{a}b$, \text{for  all}\, $a \in  \Gamma(E_{-1}), b \in \Gamma(E_{-2})$,
         \item $[fa,b,c]_{E_{-1}} = f [a,b,c]_{E_{-1}}$ for  all\, $a,b,c \in  \Gamma(E_{-1})$,
     \end{enumerate}
     \end{enumerate}
     \end{lemma}

We can now state our main result, which will be soon enlarged (see Theorem \ref{thm:existenceuniversal}).

\begin{proposition}
\label{prop:2-Lie-algebroid}
Let $(M,\mathcal{F})$ be a singular foliation that admits a geometric resolution of length $2$ as in \eqref{eq:length2} such that $(E_{-1},[\cdot\,,\cdot]_{E_{-1}},\rho)$ admits an almost Lie algebroid bracket. 
Let $ \nabla $ and $ [\cdot, \cdot,\cdot]_{E_{-1}}$ be as in Lemma \ref{lem:2-Lie-algebroid}.
     The \emph{$2$-ary bracket} on $\Gamma(E_{-1} \oplus E_{-2}) $ defined by: 
     $$ [a,b]_2 = \left\{  \begin{array}{ll} [a,b]_{E_{-1}}  & \hbox{ for}\quad a,b \in \Gamma(E_{-1}) \\ 
     \nabla_a b  & \hbox{ for}\quad a \in \Gamma(E_{-1}),\, b \in \Gamma(E_{-2}) \\ 
     -\nabla_b a  & \hbox{ for}\quad a \in \Gamma(E_{-2}),\, b \in \Gamma(E_{-1}) \\ 
     0  & \hbox{ for}\quad a,\,b \in \Gamma(E_{-2}) \\ 
     \end{array} \right. $$
     together with the \emph{$3$-ary bracket} on $\Gamma(E_{-1} \oplus E_{-2}) $ defined by $[a,b,c]_3=[a,b,c]_{E_{-1}}  $ if $a,b,c \in \Gamma(E_{-1})$ and zero otherwise,
    satisfies
    \begin{enumerate}
    \item $ [\cdot,\cdot,\cdot]_3$ is linear over functions, while for every function $f$ and for all $a\in\Gamma(E_{-1}), b\in\Gamma(E_{-1})$ or $ \Gamma(E_{-2})$,
     $$  [a,fb] = f  [a,b]_2 + \rho(a)[f] \, b $$
    \item  for all $a\in\Gamma(E_{-2}), b\in\Gamma(E_{-1})$, \begin{equation}
        \dd^{(2)}[a,b]_2=[\dd^{(2)}a,b]_2,
        \end{equation}
        \item for all $a,b,c \in\Gamma(E_{-1})$
     $$ \dd^{(2)}[a,b,c]_3+  [a,[b,c]_2]_2 +[b,[c,a]_2]_2+[c,[a,b]_2]_2=0
    $$ 
    \item for all $a,b\in\Gamma(E_{-1})$ and $ c \in \Gamma(E_{-2})$
   $$
      [a,b,\dd^{(2)}c]_3 + [a,[b,c]_2]_2 +[b,[c,a]_2]_2+[c,[a,b]_2]_2=0.
   $$
    \end{enumerate}
\end{proposition}

\vspace{0.3cm}

\begin{definition}
The structure $(E_{-1}\oplus E_{-2},\dd, \rho,[\cdot\,,\cdot]_2, [\cdot\,,\cdot, \cdot]_3)$ described in the third item of Proposition \ref{prop:2-Lie-algebroid} is known as a \emph{$2$-Lie algebroid}.
\end{definition}
A generalization to arbitrary singular foliation admitting a geometric resolution will be discussed in the Section \ref{sec:Universal}.

\section{Universal $Q$-manifolds (a.k.a. universal Lie $ \infty$-algebroids)}\label{sec:Universal}
Beyond the world of manifolds is the universe of ``manifolds up to homotopy'', which are known under various names: some or more or less equivalent, and some are mostly dual notions:
\begin{enumerate}
\item[$\clubsuit$] Lie $\infty$-algebroids, also called Lie algebroids up to homotopy. 
    \item[$\diamondsuit$] $Q$-manifolds, also called dg-manifolds (dg= differential graded) of positive degree.
\end{enumerate}
The notions $\clubsuit$ and $\diamondsuit$ are in fact ``equivalent'' in the sense that they are \underline{dual}\footnote{The reader familiar to these concepts knows that it is a bit more subtle: one needs to choose a splitting.} one to the other. 
To explain where the notion of $Q$-manifold comes from, and the duality with Lie $ \infty$-algebroids, let us start with some basic points about Lie algebras.
We invite the reader familiar with Lie $ \infty$-algebroids and $NQ$-manifolds to go directly to Theorem \ref{thm:existenceuniversal}.

\subsection{Two dual point of views on Lie algebras}
\label{sec:twodual}

For $V$ a vector space, the exterior algebra $\wedge^\bullet V = \oplus_{k=0}^\infty \wedge^k V$ is a graded algebra
with respect to a product that we denote by~$ \wedge$.

\begin{definition}
A \emph{co-Lie} algebra is a vector space $V$ equipped with a degree $+1 $ derivation 
$$ \delta \colon \wedge^{\bullet} V \longrightarrow  \wedge^{\bullet+1} V   $$
such that $\delta^2=0 $.
\end{definition}

Before explaining this definition, let us start with a few comments.
\begin{enumerate}
\item We write $\delta \colon \wedge^{\bullet} V \mapsto  \wedge^{\bullet+1} V $ to mean that for every $k \geq 0$, $ \delta$ maps $\wedge^{k} V$ to $\wedge^{k+1} V$, i.e., it is of degree~$ +1$.
    \item By a degree  $+1$ derivation, we mean that
     $$ \delta (\alpha \wedge \beta) = \delta(\alpha) \wedge \beta + (-1)^k \alpha \wedge \delta(\beta) $$
     for all $\alpha \in \wedge^{k} V  $ and $ \beta \in \wedge^{\bullet} V $. The signs are exactly those of the de Rham differential (which is also a degree $+1 $ derivation). 
    \item  For any degree $+1$ derivation, $\delta^2 $ is easily seen to be a degree $+2 $ derivation\footnote{Of course, this is not true for degree $0$ derivation, otherwise the formula $(fg)``=f''g + f g"$ would be true.}.
    \item A degree $+1$ derivation of $\wedge^\bullet V$ is entirely determined by its restriction to $V $, which is a map 
     $ \mu \colon V \longrightarrow \wedge^2  V$ that we call the \emph{co-Lie-bracket}.
     This comes from the derivation property:
      \begin{equation}\label{eq:deriv} \delta (v_1 \wedge  \dots \wedge v_k ) = \sum_{i=1}^{k}  (-1)^{i+1} v_1 \wedge  \dots \wedge v_{i-1} \wedge \mu (v_i) \wedge v_{i+1} \wedge \dots \wedge v_k \end{equation}
      for every $ \in \mathbb N$ and $ v_1, \dots, v_k \in V$.
     \item  Conversely, any linear map $ \mu \colon  V \longrightarrow \wedge^2 V$ extends to a unique degree $ +1$ derivation by using~\eqref{eq:deriv}. 
\end{enumerate}
\vspace{0.3cm}

\begin{propositions}{Lie algebras are dual to co-Lie algebra}{prop:duality}
There is a one-to-one correspondence\footnote{More precisely, to a Lie algebra structure on $ \mathfrak g$ corresponds a co-Lie algebra structure on $ \mathfrak g^*$, and to a co-Lie algebra structure on $ V$ corresponds a co-Lie algebra structure on $ V^*$.} between finite dimensional Lie algebras and finite dimensional co-Lie algebras.
\end{propositions}

\vspace{0.3cm}
\begin{proof}
The correspondence goes as follows. 
\begin{enumerate}
    \item Given a Lie algebra $(\mathfrak g, [\cdot\,, \cdot ])$, the dual of the Lie algebra bracket $ [\cdot\,, \cdot ] \colon \wedge^2 \mathfrak g \to \mathfrak g$ is a map 
     $\mu \colon  \mathfrak g^* \to (\wedge^2 \mathfrak g)^* $. Since the dimensions are finite, there is a canonical isomorphism $ (\wedge^2 \mathfrak g)^*\simeq \wedge^2 \mathfrak g^*$, and we still  denote by $\mu $ the map  $\mu \colon  \mathfrak g^* \to \wedge^2 \mathfrak g^* $. Using the derivation property in Equation \eqref{eq:deriv}, one  extends $ \mu$ to a degree $+1 $ derivation $ \delta$ of $\wedge^\bullet  \mathfrak g^*$. It is a routine to check that if the Jacobi identity holds for $[\cdot\,, \cdot] $, then it implies that $ \delta^2 =0$.
    \item Conversely, given a co-Lie algebra, the dual of  the co-Lie bracket $\mu \colon  V \to \wedge^2 V$ is a linear map $[\cdot\,, \cdot]  \colon \wedge^2 V^* \to V^* $.   It is a routine to check that  $ \delta^2 =0$ implies that the Jacobi identity holds for the bracket $[\cdot\,, \cdot] $. 
\end{enumerate}
Moreover, the two maps above are inverse one to the other.
\end{proof}

\begin{remark}
The degree $ +1$ derivation corresponding to a finite-dimensional Lie algebra over the field $  \mathbb K$ is the Chevalley-Eilenberg differential computing the Lie algebra cohomology (valued in $ \mathbb K$).
\end{remark}

To put it all in a nutshell:
\vspace{.5cm}
\begin{center}
\begin{tabular}{|l|l|
} \hline
 \multicolumn{2}{|c|}{} \\
\multicolumn{2}{|c|}{What is a finite dimensional Lie algebra?} \\ 
\multicolumn{2}{|c|}{Two dual answers.} \\ 
 \multicolumn{2}{|c|}{} \\ \hline & \\ Direct notion: &  Dual notion:\\ & \\ 
 A vector space $\mathfrak g $    &  A vector space $V$ \\ \textit{and} a linear map  &\textit{and} a degree $ +1$ derivation \\ $[\cdot\,,\cdot]\colon \wedge^2\mathfrak g \longrightarrow \mathfrak g $ & $\delta \colon \wedge^\bullet V \longrightarrow \wedge^{\bullet+1} V $\\ \hline 
 \multicolumn{2}{|c|}{} \\
 \multicolumn{2}{|c|}{such that} \\
 \multicolumn{2}{|c|}{} \\ \hline
 & \\
 The Jacobi identity holds & $ \delta^2 =0 $ \\
 & \\ \hline
\end{tabular}
\end{center}

\vspace{0.3cm}

\subsection{Graded symmetric algebras}

Throughout of this section we are  working on a field $\mathbb{K}\in\{\mathbb{R},\mathbb{C}\}$.

\vspace{1mm}
Let us introduce some terminology. 
Consider a graded vector space: 
 $$  Z_\bullet := \oplus_{i \in \mathbb Z} V_i .$$
 An element in $V_\bullet$ is said to be \emph{homogeneous} when there exists $i \in \mathbb Z$ such that it belongs to $Z_i$.
\vspace{0.5cm}

\begin{definitions}{Graded Symmetric algebras}{}
Let  $Z_\bullet := \oplus_{i \in \mathbb Z} Z_i$ be a graded vector space.
 We call  \emph{graded symmetric algebra of $Z_\bullet $} the quotient of the tensor algebra $ \oplus_{k \geq 0} Z^{\otimes k} $
 by the ideal generated by all elements of the form
 $$  \{ x \otimes y - (-1)^{ij} y \otimes x  \, | \,  i,j \in \mathbb Z,  x \in V_i, y \in V_j \}.$$
 We denote it by $ S(Z_\bullet)$. We denote its product by $\odot$.
\end{definitions}

\vspace{0.3cm}
Let us state a few basic points about graded symmetric algebras and their terminology. To start with, the graded symmetric algebra comes with two different ``degrees'' that we have to distinguish: linear combinations of elements of the form
 $$  z_{1} \odot \cdots \odot z_{k} $$
with $ z_{1} \in Z_{i_1}, \dots, z_{k} \in Z_{i_k}$ shall be said of
\begin{enumerate}
    \item \emph{polynomial degree} $k$,  because they are products of $k$ elements,
    \item and of \emph{degree} $i_1 + \cdots +i_k $, because the individual degrees of the terms in the product add up to this integer\footnote{The “degree” is sometimes called “ghost degree” by theoretical physicists.}. 
\end{enumerate}
With respect to both the polynomial degree and the degree,  $ S(Z_\bullet)$ is ``graded'' in the sense that the degree of the product of two terms is the sum of their degrees. But with respect to the degree, it is also graded commutative, i.e., 
 $$  P \odot Q = (-1)^{ij} Q \odot P$$
 for any $P,Q\in S(Z_\bullet)$ of degrees $i$ and $j$, respectively.

For homogeneous elements $z_1 ,\ldots, z_k\in Z_\bullet$ and permutation $ \sigma \in \mathfrak S_k$, the \emph{Koszul sign}, denoted by $\epsilon(\sigma,z_1,\ldots,z_k)$ or simply by $\epsilon(\sigma)$ when there is no ambiguity, is the sign  induced by the permutation of the $v_i$ in the graded symmetric algebra:\begin{equation}
    z_{\sigma(1)}\odot\cdots\odot z_{\sigma(k)}\,=\,\epsilon(\sigma,z_1,\ldots,z_k)\, z_{1}\odot\cdots\odot z_{k}.
\end{equation}
For $Z_\bullet =\oplus_{i \in \mathbb Z} Z_{i}$ a graded vector space,
we call  \emph{graded dual} denote by $Z^*_\bullet$ the graded vector space:
$$Z^*_\bullet =\oplus_{i=-\infty}^\infty Z_{-i}^* ,$$
with the understanding that elements in $Z_{-i}^*$ are of degree $+i$.
Notice that $Z^*_\bullet$ is strictly contained in the dual $ Z^*$.
According to our conventions $S(Z^*_\bullet)$ stands for the graded commutative symmetric algebra generated by $Z^*=\oplus_{i=1}^\infty Z_{-i}^*$. Elements in $Z^*_{-i_1}\odot\cdots\odot Z_{-i_k}^* $ are therefore of polynomial degree $k$ and degree $i_1+\cdots+ i_k$. We shall define elements of polynomial degree zero to be elements in $\mathbb{K}$. 
In the present section, we will be interested in two kinds of symmetric algebras:
\begin{enumerate}
    \item[$\clubsuit$] those of the form $ S(E_\bullet)$ with $ E_\bullet := \oplus_{i \geq 1} E_{-i}$ a negatively graded vector space. Except $ S^0(E_\bullet)=\mathbb K$, all its components are of negative degrees,
    \item[$\diamondsuit$] and those of the form $ S(V_\bullet)$ with $ V_\bullet := \oplus_{i \geq 1} V_{i}$ a positively graded vector space. Except $ S^0(V_\bullet)=\mathbb K$, all its components are of positive\footnote{Since it might be confusing for speakers of Latin languages, let us be precise: here, positive means $ \geq0$, and negative means $ \leq -1$.} degrees.
\end{enumerate}
If each one of the vector spaces $E_{-i}$ and $V_i$ are of finite dimension, and  if $ E_{-i}$ is the dual of $V_i $ for all $ i\geq 1$, then 
there is a natural duality between:
\begin{enumerate}
    \item[$\clubsuit$] elements of polynomial degree $k$ and degree $ -i$ in $ S(\oplus_{i \geq 1} E_{-i} )$, 
    \item[$\diamondsuit$] elements of polynomial degree $k$ and degree $ +i$ in $ S(\oplus_{i \geq 1} V_{i} )$, 
\end{enumerate}
given for all homogeneous $ v_1, \dots, v_k \in V_\bullet $ and $e_k, \dots, e_1 \in E_\bullet $ by the pairing\footnote{Some authors would divide by $ k!$}:
 $$  \left\langle \, v_1 \odot \dots \odot 
 v_k  \, \middle| \, e_k \odot \dots \odot e_k\, \right\rangle = \sum_{\sigma \in \mathfrak S_k} \epsilon(\sigma,v_1,\dots,v_k ) 
   \prod_{i=1}^k  \langle v_{\sigma (i)} , e_i \rangle .$$

 \vspace{0.5cm}

\subsubsection{Lie $ \infty$-algebras}

We are now ready for the following question\footnote{It is complicated to make a history of the notion presented in this section: let us just say that no idea here is ours, and that the notion of Lie $\infty $-algebra originated in \cite{zbMATH00065837}.}.

\vspace{0.5cm}

\begin{questions}{Towards Lie $\infty $-algebras}{}

Let $V_\bullet = \oplus_{i \geq 1} V_i $ be a positively graded vector space with each $ V_i$ of finite dimension. Let $ E_{-i}$ be the dual of $V_i$.

Assume $ S(V_\bullet)$  comes equipped with a degree $+1 $ derivation $\delta$ such that $\delta^2 =0 $. What kind of structures do we obtain on the dual spaces $ E_\bullet := \oplus_{i \geq 1} E_{-i} $? 
\end{questions}

\vspace{.5cm}

We saw in the Section \ref{sec:twodual} that if $ V_i=0$ for $i \geq 2$, such  derivations are in one-to-one correspondence  with Lie algebras structures on $ E_{-1}:=V_1^*$. 
Let us go back to the general case.
 The derivation $ \delta$ is entirely determined by its restriction to $V_\bullet$. Decomposing according to the polynomial degree, we see that $\delta = \sum_{k \geq 1 } \delta^{(k)} $, with $\delta^{(k)} \colon V_\bullet \mapsto S^k V_\bullet $ a degree $ +1$ map. 
By duality, there is a one-to-one correspondence between:
\begin{enumerate}
    \item[$\clubsuit$] the datum $\left(E_\bullet = \sum_{i=1}^\infty E_{-i} ,(\ell_k)_{k\geq 1}\right)$ made of a collection of vector spaces $E_\bullet =(E_{-i})_{i\geq 1}$ of finite dimension together with a family of degree $+1$ linear maps $\left(\ell_k\colon S^\bullet(E)\longrightarrow E \right)_{k\geq 1}$ called $k$-ary brackets,

    \item[$\diamondsuit$] a sequence $\delta^{(k)} $ of linear maps $V_\bullet \longrightarrow S^k (V_\bullet) $, with $V_i =E_{-i}^* $ for all $i \geq 1$.
\end{enumerate}

The relation between both is given for all $v \in V_\bullet$, $ e_1, \dots, e_k \in E_\bullet$ by:
 $$ \langle  \delta^{(k)} (v), e_k \odot \cdots \odot e_1 \rangle  =  \langle v , \ell_k( e_k \odot \cdots \odot e_1) \rangle. $$
A direct computation gives that $\delta^2=0$ holds if and only if the $\ell_k $'s equip $ E_{-\bullet}$ with a Lie $ \infty$-algebra structure, the latter being defined as follows\footnote{We prefer to avoid any historical comments. See, e.g. \cite{zbMATH06568059}, Sections 2.1 and 2.2, for a pedagogical introduction.}.

\vspace{0.5cm}

\begin{definitions}{Lie $ \infty$-algebras}{}
A \emph{negatively graded Lie $\infty$-algebra} is the datum $\left(E_\bullet,(\ell_k)_{k\geq 1}\right)$ made of a collection of vector spaces $E_\bullet =(E_{-i})_{i\geq 1}$ together with a family of degree $+1$ linear maps $\left(\ell_k\colon S^\bullet(E_\bullet)\longrightarrow E_\bullet \right)_{k\geq 1}$ called $k$-ary brackets, which fulfill the compatibility conditions the so-called \emph{higher Jacobi identities}: for all homogeneous elements $e_1,\ldots,e_n\in E$\begin{equation} \label{eq:higherjacobi}
    \sum_{i=1}^n\sum_{\sigma\in \mathfrak{S}_{i,n-i+1}}\epsilon(\sigma)\ell_{n-i+1}\left(\ell_i(e_{\sigma(1)},\ldots,e_{\sigma(i)}),v_{\sigma(i+1)},\ldots,e_{\sigma(n)}\right)=0.
\end{equation}
Here $\epsilon(\sigma)$ is the Kozsul sign associated to $ \sigma$ and $e_1,\ldots,e_n$, and
$\mathfrak{S}_{i,n-i+1}$ stands for the set of $(n,i)$-shuffles.
\end{definitions}
\vspace{.5cm}
In particular, for $ n=1$, Equation \eqref{eq:higherjacobi} means $ \ell^2_1=0$, so that $ (E_{\bullet}, \ell_1)$ is a complex.
Let us conclude this section by the following result.
\vspace{0.5cm}

\begin{propositions}{Dual point of views}{prop:dualLieInfinity}
Let $ E_{-\bullet} = \oplus_{i\geq 1} {E_{-i}} $
and $V_\bullet = \oplus_{i\geq 1} {V_{+i}} $
be finite dimensional graded vector spaces in duality\footnote{i.e., $ E_{-i}$ is the dual of $V_i$ for all $i \geq 1$.}.
There is a one-to-one correspondence between
\begin{enumerate}
    \item[$\clubsuit$] Lie $\infty$-algebras brackets $(\ell_k)_{k\geq 1}$ on $ E_{\bullet}$,
\item[$\diamondsuit$] degree $ +1$ derivations squaring to zero of $ S(V_\bullet) $.
\end{enumerate}
\end{propositions}

\vspace{0.5cm}
To put it all in a nutshell, in view of Proposition \ref{thm:prop:dualLieInfinity}, we have
\vspace{0.2cm}
\begin{center}
\begin{tabular}{|l|l|
} \hline
 \multicolumn{2}{|c|}{} \\
\multicolumn{2}{|c|}{What is a finite dimensional Lie $\infty $-algebra?} \\  
\multicolumn{2}{|c|}{Two dual answers} \\ 
 \multicolumn{2}{|c|}{} \\ \hline & \\ Direct notion: &  Dual notion:\\ & \\ 
 A neg. graded vector space $E_\bullet = \oplus_{i=1}^\infty E_{-i}  $    &  A pos. graded vector space $V_\bullet =\oplus_{i=1}^\infty V_{i} $ \\ \textit{and} a family of degree $ +1$ linear maps  &\textit{and} a degree $ +1$ derivation \\ $\ell_k \colon S^k(E_\bullet)  \longrightarrow E_\bullet $ (called “brackets”)& $\delta \colon S( V_\bullet) \longrightarrow S( V_\bullet)  $\\ \hline 
 \multicolumn{2}{|c|}{} \\
 \multicolumn{2}{|c|}{such that} \\
 \multicolumn{2}{|c|}{} \\ \hline
 & \\
  Eq. \eqref{eq:higherjacobi} ("higher Jacobi") holds for all $n$ & $ \delta^2 =0 $ \\
 & \\ \hline
\end{tabular}
\end{center}

\vspace{0.3cm}

\subsection{(Positively graded) $NQ$-manifold}

We will now extend the previous discussion from Lie $ \infty$-algebras to Lie $\infty $-algebroids\footnote{We prefer not to make an historical introduction: we acknowledge being inspired by several works by Pavol Severa \cite{zbMATH05049976}, Theodore Voronov \cite{Voronov}, as well as by Giuseppe  Bonavolont{\`a} and Norbert Poncin's \cite{zbMATH06243829}, see also the first sections of \cite{zbMATH07568450} and Section 3.4 in \cite{LLS} for a similar presentation.}, and from $ S(V_\bullet)$ equipped with a degree $ +1$ differential to the so-called $Q$-manifolds. 

\subsubsection{Graded manifolds}

Let us first define positively graded manifolds. In addition to the purely mathematical definition, we give as footnotes some vocabularies commonly used by mathematical physicists working with this object.

\vspace{.5cm}

\begin{definitions}{Graded manifolds: the objects}{gradedmfd}
Let $M$ be a smooth, real analytic or complex manifold and $\mathbb K=\mathbb R$ or $ \mathbb C$ depending on the chosen context. Let $ \mathcal O$ be the corresponding sheaf of (smooth, real analytic, or holomorphic) functions on $M$.

A \emph{positively graded manifold} over a manifold\footnote{It is customary to  call $M$ the \emph{base} of the sheaf.} $M$ is a sheaf 
$$ \mathcal E \colon \mathcal U \mapsto \mathcal E_{\mathcal U} $$
of graded commutative algebras over $\mathbb K $
such that every $m \in M$ admits an open neighborhood $\mathcal U\subset M$ on which the sheaf\footnote{It is customary to call sections of the sheaf $ \mathcal E$ \emph{graded functions}.} structure takes the form\footnote{It is convenient to denote a graded manifold as a pair $(M,\E)$.} 
$$\E_{\mathcal U}=\mathcal O_{\mathcal U }\otimes_\mathbb{K} S(V_\bullet)$$ for some graded vector space\footnote{It is customary to say that the graded manifolds in \emph{concentrated in degrees} $ [1, \dots,N]$ if all $V_i =0$ when $ i \notin [1, \dots, N]$.} $V=\oplus_{i=1}^\infty V_{+i}$. 
\noindent

\end{definitions}

\vspace{0.5cm}

\begin{remark}
\label{rmk:Qsmooth}
In the smooth setting, it can be proven that there exists a globally and canonically defined graded vector bundle
 $ V_\bullet \to M $ such that $\mathcal E $ is isomorphic to the sheaf of sections of the graded commutative algebra bundle  $ S(V_\bullet) \to M $, obtaining by considering pointwise the symmetric algebra of$ V_\bullet$. Although $V_\bullet$ is canonical, the isomorphism of sheaves 
  $$\xymatrix{\mathcal E \ar@{<->}[rr]^{\simeq \hspace{0.5cm} } & &  \Gamma(S(V)) } $$
 is \emph{not} canonical, and is called a \emph{splitting}. For a statement adapted to the present situation, see \cite{KotovSalnikov}.
Upon choosing such an isomorphism, a function $F\in\E_i$ is a formal sum \begin{equation}\label{eq:polynomialdegreedecomposition}
    F=\sum_{k\geq 0}F^{(k)}
\end{equation}
with $F^{(k)}\in \E$ an element of polynomial degree $k$ and degree $i$. 
For degree reasons, the sum in \eqref{eq:polynomialdegreedecomposition} must be finite.  To be more precise, the only possible values of $k$ are between $ 1 \leq k\leq i$.
\end{remark}

\begin{remark}
Remark \ref{rmk:Qsmooth} has consequences on vocabulary.
Even in the holomorphic or real analytic settings, one calls sections in the sheaf $ \mathcal E$ under the name of \emph{(graded) functions}. But what are these sections “functions” on? In the smooth case, we can answer this question.
Consider the negatively graded vector bundle $ E_{-\bullet} = \oplus_{i\geq 1} {E_{-i}} $ with $ E_{-i}= V_i^*$ for all $i \geq 1 $.
Sections of the sheaf $\E$  are  considered by physicists as \emph{functions on $E_\bullet $}.

\end{remark}

\begin{remark}
\label{rem:rankVi}
If the manifold $M$ is connected, then for every $i \geq 1$, the  dimension of  the component $V_i$ of degree~$i$ in the graded vector space $ V_\bullet$ that appear in Definition \ref{def:gradedmfd} does not depend on the point $m$.  From now on, we will assume that $M$ is connected and denote by $ r_i$ the rank of $ V_i$.
\end{remark}

\subsubsection{Local coordinates of a graded manifold}
Among the somewhat cumbersome vocabulary that makes graded manifolds hard to deal with for the non-used reader is the notion of “local graded coordinates”.
Assume that $M$ is connected and let $r_i$ be as in Remark \ref{rem:rankVi}. For any $\mathcal{U}\subset M$ an open subset such that $({V_i})_\mathcal{U}\overset{\sim}{\longrightarrow}\mathcal{U}\times \mathbb{K}^{r_i}$ for every $i\geq 1$, we call  \emph{graded coordinates} on the graded manifold $(M,\E)$ the data made of
\begin{enumerate}
    \item
  a system of coordinates $(x_1,\ldots,x_n)$ of $M$ on $\mathcal U$.
\item 
and in all degree $i\geq 1$, a local trivialization $\left(\xi_i^{[1]},\ldots, \xi_i^{[r_i]}\right)$ of $V_i$ on $\mathcal U$.
\end{enumerate}

\vspace{1mm}
\noindent
After choosing graded coordinates as above, i.e., a list $$\left(x_1,\ldots,x_n,\xi_1^{[1]},\ldots, \xi_1^{[r_1]},\ldots, \xi_i^{[1]},\ldots, \xi_i^{[r_i]},\ldots\; \right),$$
then any element of $\E(\mathcal U)$  is a  ``polynomial'' in finitely many of the variables  
$$ 
\left(\xi_i^{[j]}\right)_{j=1,\ldots,r_i, \;i\geq 1}
$$ 
with coefficients in the algebra $\mathcal{O}(\mathcal{U})$ of smooth, real analytic or complex functions over $U$. 
For instance, an element in degree $3$ decomposes as a sum of the form
 $$  \sum_{\scalebox{0.5}{ $1 \leq i < j < k \leq r_1 $ }} f_{ijk}(x_\bullet ) \, \xi^{[i]}_1 \odot \xi^{[j]}_1 \odot \xi^{[k]}_1 
 + \sum_{\scalebox{0.5}{$\begin{array}{c}1 \leq i  \leq r_1 \\ 1 \leq j \leq r_2 \end{array}$}}  g_{ij}(x_\bullet )  \, \xi^{[i]}_1 \odot \xi^{[j]}_2 
 +
 \sum_{\scalebox{0.5}{$1 \leq i \leq r_3$}} h_i (x_\bullet )\, \xi^{[i]}_3 $$
 with $ f_{ijk},g_{ij},h_i$ being smooth, real analytic or holomorphic
functions in the variables $ x_\bullet =x_1, \dots, x_n$.

\begin{example}
The sheaf of differential forms $(M, \E=\Omega(M))$ on a manifold $M$ is a graded manifold since for every point $m\in M$, it takes the form $\mathcal{O}_\mathcal{U}\otimes_\mathbb{K}\wedge^\bullet T_m^*M$ on an open neighborhood $\mathcal U$  of $m$. Exterior forms can be seen as functions on the graded vector bundle $E_{-1}=TM$, concentrated in degree $-1$.
\end{example}
\begin{example}
Let $k$ be a positive integer. A finite dimensional vector space $E$ and its dual $V$ can be seen as graded vector bundles of respective degree $-k$ and $k$ over a point.
$E$ is a graded manifold over $M=\{\mathrm{pt}\}$ and the sheaf $ \mathcal E$ is isomorphic (as an algebra) to $ \wedge V$ for $k$ odd and $S(V)$ for $k$ even. 
It must be understood that $ \wedge^i V$ or $ S^i(V)$ are then of degrees $ ki$.
\end{example}

We now define morphisms.

\vspace{0.5cm}

\begin{definitions}{Graded manifolds: the morphisms}{gradedMorphisms}
A \emph{morphism of graded manifolds} from $(M',\mathcal \E')$ to $(M, \mathcal{E})$  is a pair made of a smooth or real analytic or holomorphic map $\phi\colon M'\longrightarrow M$ called the \emph{base
map} and a sheaf morphism over $ \phi^*$, i.e., a family of graded algebra morphisms:
 $$ \mathcal E_{\mathcal U} \to  \E_{\phi^{-1}(\mathcal U)}' ,$$
compatible with the restriction maps, such that\begin{equation}
    \Phi(f\alpha)=\phi^*(f )\Phi(\alpha).
\end{equation}
for all $f\in\mathcal O_{\mathcal U}$ and $\alpha\in\E_{\mathcal U}$.

A \emph{homotopy} between two morphisms of graded manifolds
$ \Phi, \Psi \colon  (M',\mathcal \E') \longrightarrow (M, \mathcal{E})$ is a morphism
of graded manifolds\footnote{The direct product of graded manifolds considered above goes with some subtleties. The henceforth obtained sheaf is a sheaf over $ M' \times [0,1]$. The graded vector bundle as in Remark \ref{rmk:Qsmooth} is $  p_1^!  V_1 \oplus  p_2^! T[0,1]$ in degree $1$ and $ p_1^! V_i$ in degree $ i \geq 2$, with $ p_1,p_2$ the projections on the first and second component respectively. We refer to Section 3.4 in \cite{LLS} for more details.}
 $$ (M,\mathcal E)  \longrightarrow (M',\mathcal E') \times ([0,1], \Omega([0,1])) $$
 whose restrictions to the extremities of the interval $[0,1] $ coincide with $\Phi$ and $\Psi $ respectively.
\end{definitions}

\subsubsection{Vector fields on graded manifolds}

Vector fields on a manifold $M$ are derivations of its sheaf of algebra of functions: this principle is valid in the smooth, real analytic or complex settings. For a graded manifold $(M, \E)$, the equivalent of functions are the sections of the sheaf $\mathcal E$. Since it is not commutative but graded commutative, one has to consider graded derivations. 
Consider $k \in  \mathbb Z$.
A \emph{graded derivation of degree $k$} of $\mathcal E $ is the data, for every $\mathcal U \subset M $ of a linear map
 $$ Q \colon  (\mathcal E_{\mathcal U})_\bullet \longrightarrow (\mathcal E_{ \mathcal U})_{\bullet +k} ,$$
 compatible with all restriction maps, that increases the degree by $+k$ and satisfies: 
  $$  Q[F\,G] = Q[F] \, G +(-1)^{ki}F \, Q[G] $$
  for every $F \in  (\mathcal E_{\mathcal U})_i, G \in \mathcal E_{\mathcal U}$.
    Since we think geometrically, we will simply say \emph{vector fields of degree $k$} instead of graded derivations.

\begin{definition}
We call Let $(M,\E)$ be a smooth graded manifold. For any open subset $\mathcal U \subset M$, we denote by $ \mathfrak X_\bullet (\mathcal E)$ the space of graded derivations of degree $k$
and call its elements \emph{vector fields of degree $k$}

In the complex or real analytic settings, graded derivations of degree $k $ do not form not a sheaf, and one has to define $\mathfrak X_\bullet (\mathcal E)$  using local graded derivations. 

 In all three cases, the sheaf 
 $$\mathcal U\longmapsto \mathfrak{X}_\bullet(\E)_{\mathcal U}=\oplus_{k \in \mathbb Z}  \mathfrak{X}_k (\E)_{\mathcal U}$$ 
 is a sheaf of $\E$-modules that we call \emph{sheaf of graded vector fields\footnote{We will simply say "vector fields on $ \mathcal E$" most of the time.} on $ (M,\mathcal E)$}. 
\end{definition}

Let us list some important facts on vector fields on $(M,\E)$:
\begin{enumerate}
     \item The sheaf $\mathfrak X_\bullet(\mathcal E):=\oplus_{k\mathbb{Z}}\mathfrak X_k(\E)$ of vector fields on $(M,\E)$ is a graded $\E$-module. Also, $\mathfrak X_\bullet(\E)$ of  is a sheaf of graded Lie algebras. The graded Lie bracket \begin{equation}\label{graded:commutator}
         [P,Q]=P\circ Q-(-1)^{kl}Q\circ P\end{equation} of two vector fields $P,Q$ of degree $k,l$ respectively is a vector field of degree $k+l$. It is easily checked that the bracket \eqref{graded:commutator} fulfills
         
         \begin{enumerate}
             \item $[P,Q]=-(-1)^{jk}[Q,P]$\quad (graded skew-symmetry)
             \item $(-1)^{jl}[P, [Q, R]] + (-1)^{jk}[Q, [R, P]]+(-1)^{kl}[R, [P, Q]]=0,$\quad (graded Jacobi identity)
         \end{enumerate}
         for all graded vector fields $P,Q, R\in\mathfrak X(\E)$ of degrees $j,k$ and $l$, respectively.
    \item Let $ E_\bullet$ be a graded vector bundle over $M$ as in Remark \ref{rmk:Qsmooth}, i.e., such that there exists a splitting
     $$ \mathcal E \simeq \Gamma(S(V)) $$
     with $ V_i= E_{-i}^*$ for all $ i \geq 0$.
     Any section $e\in \Gamma(E_{-i})$ corresponds to a vector field\footnote{Such vector fields are called \emph{vertical} because they are linear over functions on the base manifold. In this case, the polynomial degree is $-1$ and its degree is $-i$.} $\iota_e\in\mathfrak X_{-i}(\E)$ defined by contraction with $e$.
\item  Let $(\mathcal U, x_1,\ldots,x_n)$ local coordinates on $M$ and $(\xi_i^{[j]})_{j=1,\ldots, r_i}$ with $i\geq 1$ be local coordinates on $ \E$, defined as above.

Let $(e_i^{[j]})_{j=1,\ldots, r_i},i\geq 1$ be the dual basis of $(\xi_i^{[j]})_{j=1,\ldots, r_i},i\geq 1$. Then for every pair $i,j$, it is customary to consider the contraction $\iota_{e_i^{[j]}}$ of item as being the “partial derivative”:
 $$\frac{\partial}{\partial \xi_i^{[j]}}$$ 
 The notation comes from the computational fact that, in local coordinates, the effect of $\iota_{e_i^{[j]}}$ is similar to the one of the partial derivative with respect to $\xi_i^{[j]}\in \Gamma( V_i)$. 
 
Elements in the following list are vector fields of degree $k$:
$$\displaystyle{\left(\xi_{i_1}^{[j_1]}\odot\cdots\odot\xi^{[j_l]}_{i_l}\frac{\partial}{\partial x_j}\right)_{
\tiny{
\begin{array}{c}

i_1+\cdots +i_l=k\\

\end{array}
}
}\hbox{ and }\left(\xi^{[j_1]}_{i_1}\odot\cdots\odot\xi^{[j_l]}_{i_l}\frac{\partial}{\partial \xi_i^{[j]}}\right)_{\tiny{
\begin{array}{c}

i_1+\cdots+ i_l-i=k

\end{array}}}}$$ 
that form a “basis” for $\mathfrak X_k(\E)(\mathcal{U})$ in the sense that any  vector field $Q\in\mathfrak X_k(\E)(\mathcal U)$ is an infinite\footnote{The sum containing $\frac{\partial}{\partial x_j}$ is finite, but the second one is only finite if finitely many of the $ E_{-i}$'s are non-zero.} sum of the form
$$Q=\sum_{{\tiny{\begin{array}{c}

i_1\cdots i_l=k

\end{array}}}}\frac{1}{l!}\,{}^{j}Q^{j_1\cdots j_l}_{i_1\cdots i_l}\,\xi^{[j_1]}_{i_1}\odot\cdots\odot\xi^{[j_l]}_{i_l}\,\frac{\partial}{\partial x_j}+\sum_{{\tiny{\begin{array}{c}

i_1\cdots i_l-i=k

\end{array}}}}\frac{1}{l!}\,{}^{ij}Q^{j_1\cdots j_l}_{i_1\cdots i_l}\, \xi^{[j_1]}_{i_1}\odot\cdots\odot\xi^{[j_l]}_{i_l}\frac{\partial}{\partial \xi_i^{[j]}}.$$ 
for some smooth, real analytic or holomorphic functions $Q^{j_1\cdots j_l}_{i_1\cdots i_l}$ on $\mathcal U$. These functions are unique, in one chooses an ordering on the indices.
For example, a vector field $Q$ of degree $+1$ can be written in these notations as
$$Q=\sum_{_{\tiny{\begin{array}{c}1\leq u\leq r_1\\j=1,\ldots,n
\end{array}}}}{}^{j}Q^{u}_{1}\, \, \xi^{[u]}_{1}\,\frac{\partial}{\partial x_j}+\sum_{_{\tiny{\begin{array}{c}

i_1+\cdots + i_l-i=1

\end{array}}}}\frac{1}{l!}\,{}^{ij}Q^{j_1\cdots j_l}_{i_1\cdots i_l}\, \, \xi^{[j_1]}_{i_1}\odot\cdots\odot\xi^{[j_l]}_{i_l}\frac{\partial}{\partial \xi_i^{[j]}}.$$ 
\end{enumerate}
We can define $NQ$-manifolds. 

\vspace{0.5cm}

\begin{definitions}{$NQ$-manifolds}{defQmanif}
A \emph{$NQ$-manifold\footnote{Also called in this context positively graded $Q$-manifold or positively graded dg-manifold.}} is a positively graded manifold $(M, \E)$ endowed with a
degree $+1$ homological vector field $Q \in\mathfrak{X}_1(\mathcal E)$  such that $Q^2 = 0$.\\
It shall be denoted as a triple $ (M,\mathcal E, Q)$.
\end{definitions}

\begin{example}
When trying to understand $NQ$-manifolds, the novice reader may be surprised by phrases like "A finite dimensional Lie algebra is a $NQ$-manifold with base manifold a point and concentrated in degree $1$". Let us explain the meaning of such words.
Given a finite dimensional Lie algebra $(\mathfrak g,[\cdot\,,\cdot\,])$ of dimension $d$. 
The pair $(M=\{\mathrm{pt}\},\, \E=\wedge^\bullet\mathfrak{g}^*)$ is a graded manifold over $M=\{\mathrm{pt}\}$, which is generated by a vector space of degree $ -1$ (hence "concentrated in degree $ -1$"). 
When equipped with the Chevalley-Eilenberg differential which is a derivation of degree $ +1$, i.e., a vector field of degree $ +1$, and  squares to zero. 
It is therefore a $NQ$-manifold. Let us write it in coordinates:
Fix a basis $(e_i)_{i=1,\ldots,d}$ of $\mathfrak g$ and let $\xi^i\in \mathfrak g^*$ with $i=1,\ldots,d$ be the dual basis.  Consider the Christoffel symbols
$\lambda_{ij}^k\in\mathbb{K} $ defined by
$$[e_i,e_j]=\sum_{k=1}^d\lambda_{ij}^k e_k.$$
The Chevalley-Eilenberg operator is given by the degree $+1$ “vector field”:
$$ Q=\frac{1}{2}\sum_{i,j,k= 1}^d\lambda_{ij}^k\, \xi^i\wedge\xi^j\frac{\partial}{\partial\xi^k}.$$
\end{example}

\begin{example}
Given a  complex of vector bundles $(E_{-\bullet},\dd^{(\bullet)})$
over $M$. There is a natural dg-manifold
given by its sheaf of sections $(M, \E = \Gamma(S(E^*))$ and whose homological vector field  $Q\in \mathfrak X_{+1}(E)$ is obtained by dualization  $\dd^{*}\colon S^1(E^*)\longrightarrow S^1(E^*)$ of the differential map  $\dd\colon E \longrightarrow E$ map and $\rho^*\colon  T^*M\rightarrow E_{-1}^*$ of the anchor $\rho\colon E_{-1}\rightarrow TM$ 

\begin{equation}
   \Gamma(E_{-1}^*)\stackrel{\dd^{*}}{\longrightarrow}\Gamma(E_{-2}^*)\stackrel{\dd^{*}}{\longrightarrow}\cdots
\end{equation}
that we extend to a  derivation on $\E$ squaring to zero. 
\end{example}

\begin{example}Let $M$ be a manifold.
The sheaf $(M, \E=\Omega(M))$, equipped with De Rham differential, is a $NQ$-manifold. In terms of coordinates, the homological vector field reads
$$\sum_{i=1}^n dx_i\frac{\partial}{\partial x_i}.$$
\end{example}
Let us introduce some vocabulary that will need to use.
\begin{definition}
\label{def:degreedef}
Let $(M,\E', Q')$ and $(M,\E, Q)$ be two $NQ$-manifolds coming with splittings $ \mathcal E \simeq \Gamma(S(V_\bullet))$ and  $ \mathcal E' \simeq \Gamma(S(V_\bullet'))$.
A linear map $\Phi\colon \E \to \E'$ is said to be of \emph{polynomial degree/degree} $j\in\mathbb{Z}$ provided that for all functions $\alpha\in \E$ of polynomial degree/degree $i$, $\Phi(\alpha)$ is of polynomial degree/degree $i+j$.

         \end{definition}

\begin{remark}Assume that we are in the context of Definition \ref{def:degreedef}, then two points are worth noticing. 
\begin{enumerate}
\item Any map $\Phi\colon \E\to \E'$ of degree $i$ decomposes w.r.t the polynomial degree as follows:\begin{equation*}
        \Phi=\sum_{r\in\mathbb{Z}}\Phi^{(r)}
    \end{equation*}
with $\Phi^{(r)}\colon \E \to \E'$ a map of polynomial degree  $r$. 
    \item When  $\Phi\colon \E \to \E'$ is a graded morphism of algebras, its degree is zero, which implies that $\Phi^{(r)} = 0 $  for all $r<0$, i.e., it can not have components of negative polynomial degree. Furthermore, for all $n, r \in\mathbb N$ and all $\xi_1 ,\ldots, \xi_k \in \Gamma (V)$ one has the following decomposition:
\begin{equation}
\Phi^{(r)}(\xi_1\odot\cdots\odot\xi_n)= \sum_{i_1+\cdots+i_n =r}\Phi^{(i_1 )} (\xi_1 )\odot \cdots \odot\Phi^{(i_n)} (\xi_n).\end{equation} In particular, $\Phi$ is uniquely determined by the values on  $\Gamma(V)$.
\end{enumerate}
\end{remark}
In the context of Lie algebroids, the following idea is attributed to Vaintrob \cite{zbMATH01262295}.

\vspace{.5cm}

\begin{definitions}{Morphisms}{morphismsNQ}
Let $(M,\E,Q)$ and $(M',\E',Q')$ be two $NQ$-manifolds  with
sheaves of functions $\E$ and $\E'$ respectively. A \emph{morphism of $NQ$-manifold over $M$} 

is a morphism of graded manifolds
 $\Phi\colon \E \to \E'$  over   which intertwines 
$Q$ and $Q'$, i.e.,
\begin{equation}
 \Phi\circ  Q=Q'\circ \Phi. 
\end{equation}
\end{definitions}

\subsection{Negatively graded Lie $\infty$-algebroids and their morphisms}
\label{sec:hiherstructures}

We denote by $ \mathcal O$ the sheaf of smooth, real analytic or holomorphic functions on $M$, and by $\Gamma(E) $ the sheaf of sections of the graded bundle $ E_\bullet$. 

\vspace{.5cm}

\begin{definitions}{ Lie $ \infty$-algebroid}{NGLA}

\label{NGLA}
A \emph{Lie $\infty$-algebroid\footnote{"Negatively graded Lie $ \infty$-algebroid" would be a more precise name, but since they are all negatively graded, we just say " Lie $ \infty$-algebroid". We do not like the confusing notation $ L_\infty$-algebroid, although it is often used.}} $\left(E_\bullet,(\ell_k)_{k\geq 1}, \rho\right)$ is a collection of vector bundles $E =(E_{-i})_{i\geq 1}$ over $M$ endowed with a sheaf of Lie $\infty$-algebra
structures $(\ell_k)_{k\geq 1}$ over the sheaf of sections of $E$ together with a vector bundle morphism $\rho\colon E_{-1}\to TM$, called the \emph{anchor map}, such that the $k$-ary-brackets $$\ell_k: \underbrace{\Gamma(E_\bullet)\times \cdots \times \Gamma(E_\bullet)}_{k {\text{-times}}}\longrightarrow \Gamma(E_\bullet)$$ are all $\mathcal O$-multilinear  except when $k=2$ and at least one of the arguments is of degree $-1$. The $2$-ary bracket  satisfies the Leibniz identity

\begin{equation}
    \ell_2(x, f y) = \rho(x)[f]y + f\ell_2(x, y),\; f\in \mathcal{O},\; x \in \Gamma(E_{-1}), y\in\Gamma(E_\bullet).
\end{equation}
\end{definitions}

\vspace{0.5cm}

\begin{exo}
    Show that Definition \ref{NGLA} implies that $\rho(\ell_2(x,y) ) = [\rho(x), \rho(y)]$ 
for all $ x, y \in\Gamma(E_{-1})$ and that $\rho\circ \ell_1=0$.
\end{exo}
\begin{remark}
\label{rmk:linearpart}
Definition \ref{NGLA} implies the following facts
\begin{enumerate}
\item The sequence of morphisms of vector bundles
$$\xymatrix{\cdots\ar[r]^{\ell_1}&E_{-2}\ar[r]^{\ell_1}&E_{-1}\ar[r]^{\rho}&TM}$$ is a complex of vector bundles that we call the \emph{linear part}. A Lie $\infty$-algebroid is said to be \emph{acyclic} if its linear part
has no cohomology in degree $\leq -1$.
\item The 2-ary bracket restricts to  an almost-Lie algebroid
structure on $E_{-1}$. Hence, by Lemma \ref{lemma:almost&foliation}, the image\footnote{As usual, one can take global compactly supported sections in the smooth case, otherwise, the image has to be taken in the sense of sheaves.} $\mathcal{F}:=\rho(\Gamma(E_{-1}))\subseteq \mathfrak X(M)$ is a singular foliation on $M$ called the \emph{basic singular foliation} of $\left(E_\bullet,(\ell_k)_{k\geq 1}, \rho\right)$. We say, then, that the Lie $\infty$-algebroid $\left(E_\bullet,(\ell_k)_{k\geq 1}, \rho\right)$ is \emph{over $\mathcal{F}$}.
\end{enumerate}

\end{remark}
Here is the most important theorem of the present section \cite{Voronov,zbMATH06243829}.

\vspace{0.5cm}

\begin{theorems}{Duality $NQ$-manifolds - Lie $ \infty$-algebroids}{theo:correspondence}

Let $M$ be a manifold, $E_\bullet$ and $V_\bullet$ negatively graded and, respectively, positively graded vector bundles over $M$ such that $ E_{-i} = V_i^{*}$ for all $ i \geq 1$. There is a one-to-one correspondence between:
\begin{enumerate}
    \item[(i)] negatively graded Lie $\infty$-algebroids $\left(E_\bullet,(\ell_k)_{k\geq 1}, \rho\right)$,
\item[(ii)] $NQ$-manifolds $(M,\E,Q)$ with $ \mathcal E = \Gamma(S(V_\bullet))$.
\end{enumerate}

\end{theorems}

\vspace{0.5cm}

To make Theorem \ref{thm:theo:correspondence} complete, let us describe the duality relation through which both structures correspond one to the other:

\begin{enumerate}
\item for all $f\in \mathcal{O}$, $e \in\Gamma(E_{-1})$\begin{equation}
    \langle Q(f), e\rangle = \rho(e)[f],
\end{equation}
    \item for all $\xi\in\Gamma(E^*)$ and $e\in\Gamma(E)$:
\begin{equation}
\langle Q^{(0)}(\xi), e\rangle = (-1)^{\lvert\xi\rvert}\langle\xi, \ell_1(e)\rangle\end{equation}.

\item for all homogeneous elements $e_1,e_2\in \Gamma(E)$ and $\xi \in\Gamma(E^*)$
\begin{equation}
  \langle Q^{(1)}(\xi), e_1\odot e_2\rangle = \rho(e_1)[\langle\xi,e_2\rangle]-\rho(e_2)[\langle \xi, e_1\rangle]-\langle \xi, \ell_2(e_1, e_2)\rangle,
\end{equation}
with the understanding that the anchor $\rho$
vanishes on $E_{-i}$ when $i\geq 1$.
\item for every $k \geq 3$, the $k$-ary brackets $\ell_k \colon \Gamma( S_\mathbb{K}^k(E))\to\Gamma (E)$ and the polynomial degree $k-1$ component
$Q^{(n-1)}\colon\Gamma(E^*)~\longrightarrow~\Gamma(S_\mathbb{K}^k( E^*))$ of $Q$ are dual to each other.
\end{enumerate}
Here $\langle\cdot\,,\cdot\rangle$ stands for the duality pairing between sections of a vector bundle and sections of its dual.

\begin{definition}
\label{def:dualNQ}
A Lie $ \infty$-algebroid and a $NQ$-manifold that corresponds one another as in Theorem \ref{thm:theo:correspondence} are said to be \emph{dual} one to the other.
\end{definition}

\begin{example} \cite{zbMATH01262295}
Let $(A,[\cdot\,,\cdot\,]_A,\rho)$ be a Lie algebroid concentrated in degree $-1$. The graded manifold $(M,\E=\Gamma(\wedge A^*))$  carries a dg-manifold structure $Q\in\mathfrak X(A)$ which is given by
\begin{align*}
    \langle Q[f], a\rangle&=\rho(a)[f]\\
    \langle Q[\xi], a\wedge b \rangle&= \rho(a)[ \langle \xi,b\rangle]-\rho(b)[ \langle \xi,a\rangle] -\langle \xi,[a,b]_A\rangle
\end{align*}
for $f\in \mathcal{O},\xi\in \Gamma(A^*)$ and $a,b\in \Gamma(A)$. This is sufficient to extend $Q$ by derivation on $\E$. One can check that $Q^2=0$ because of Jacobi identity.
In particular, the Lie algebroid of vector fields on $M$ correspond to the De Rham complex.
\end{example}

\vspace{0.5cm}

\begin{definitions}{Lie $ \infty$-morphisms are defined through duality}{}Let $\left(E'_\bullet,(\ell'_k)_{k\geq 1}, \rho'\right)$ and $\left(E_\bullet,(\ell_k)_{k\geq 1}, \rho\right)$ be Lie $\infty$-algebroids over $M$ and, respectively, $M'$. Let $(M',\mathcal E',Q')$ and respectively,  $(M,\mathcal E,Q)$ be their dual $NQ$-manifolds\footnote{See Definition \ref{def:dualNQ}.}.
A \emph{Lie $\infty$-algebroids morphism} or \emph{Lie $\infty$-morphism} from $\left(E_\bullet',(\ell'_k)_{k\geq 1},\rho'\right)$ to $\left(E_\bullet,(\ell_k)_{k\geq 1},\rho\right)$ is a morphism\footnote{Notice that it goes backward.} of $NQ$-manifolds from $(M,\mathcal E ,Q)$ to $(M',\mathcal E',Q')$.
\end{definitions}

\subsubsection{Homotopic Lie $\infty$-algebroids}

Having defined objects (Lie $\infty$-algebroids) and arrows (morphisms of Lie $\infty$-algebroids), we still have to define homotopy between morphisms.

\vspace{0.5cm}

\begin{definitions}{Homotopies between morphisms}{}{}
\label{defn:homotopies}
Let $\left(E'_\bullet,(\ell'_k)_{k\geq 1}, \rho'\right)$ and $\left(E_\bullet,(\ell_k)_{k\geq 1}, \rho\right)$ be smooth Lie $\infty$-algebroids over $M$ and, respectively, $M'$. We assume $ E_\bullet$ and $ E_\bullet'$ to be of finite length. Let $(M',\mathcal E',Q')$ and, respectively,  $(M,\mathcal E,Q)$ be their dual $NQ$-manifolds\footnote{See Definition \ref{def:dualNQ}.}
Two Lie $\infty$-morphisms $\Phi,\Psi\colon(M,\mathcal E,Q)\longrightarrow (M',\mathcal E',Q')$ are said to be \emph{homotopic}\footnote{In the  real analytic settings, one can assume forms on $ [0,1]$ are real analytic, and in the complex case, $ [0,1]$ can be replaced by a neighborhood of $ [0,1]$ in $ \mathbb C$. One has then to consider the equivalence relation generated by this relation, which is not necessary in the smooth case.} (and we then write  $\Phi\sim \Psi$) if there is a
$NQ$-manifold morphism\footnote{We leave it to the reader to define the direct product of two $ NQ$-manifolds (for  graded manifolds, see Definition \ref{def:gradedMorphisms}).}
$$  (M,\mathcal E,Q) \longrightarrow \ (M,\mathcal E',Q') \times  ([0,1], \Omega^\bullet([0,1]), \dd^{DR})  ,$$
whose restriction to the extremities of the interval are $ \Phi$ and $ \Psi$, respectively. \end{definitions}
\vspace{0.5cm}

When $M=M'$, Lie $ \infty$-algebroid morphisms can be also defined using Taylor coefficients\footnote{See, e.g., \cite{zbMATH07568450} for a pedagogical introduction.}

\begin{remark} 
Equivalently, an equivalence $\Phi\sim\Psi$, consists of:
\begin{enumerate}
    \item  a  piecewise-smooth path $t\mapsto\Phi_t$ valued in $NQ$-manifold morphisms such that $$\Phi_0 = \Phi\quad
\text{and}\quad \Phi_1 = \Psi,$$
\item a piecewise-smooth path $t\mapsto H_t$ valued in $\Phi_t$-derivations of degree $-1$, such that the following equation:\begin{equation}
    \frac{d\Phi_t}{dt}=Q'\circ H_t + H_t\circ Q 
\end{equation}
holds for every $t\in [0, 1]$ where it is defined.
\end{enumerate}

We refer to Section 3.4.4 in \cite{LLS} or to Section 1.2.5. in \cite{CamilleLouis} for more details - there is a subtlety if the length of $ E_\bullet$ or $E_\bullet'$ is not bounded.
\end{remark}

\begin{exo}
 Show that Definition \ref{defn:homotopies} implies that for every pair of  homotopic Lie $\infty$-morphisms $\Phi, \Psi\colon (M',\E',Q')\longrightarrow (M,\mathcal E,Q)$, there exists an $\mathcal{O}$-linear map $H\colon \E \longrightarrow \E'$
of degree $-1$ such that:
\begin{equation}
    \Psi-\Phi = Q'\circ H+H\circ Q.
\end{equation}
\end{exo}

\begin{exo} Show that, in the smooth setting, homotopy of Lie $\infty$-algebroid morphisms is an equivalence relation and also this equivalence relation is compatible with composition of Lie $\infty$- morphisms.
\end{exo}

\subsection{NQ-manifolds and singular foliations}

For any Lie $\infty$-algebroid over $M$, the image of the anchor map is a singular foliation called its basic foliation. Now we will analyze the opposite direction, i.e., given a singular foliation $\mathcal{F}\subseteq \mathfrak{X}(M)$,  we will try to find a Lie $\infty$-algebroid over $M$ whose basic foliation is $\mathcal{F}$? We will of course construct this structure on a geometric resolution, studied in Section \ref{sec:Geometric-resolutions}. We will show that it is as unique as can be, i.e., unique up to equivalence. This extends a similar discussion that we had about almost-Lie algebroids, see Theorem  \ref{th:universal}.

The next theorem is obtained by proving that every graded almost Lie algebroid  over a geometric resolution can be extended to a (unique up to homotopy) Lie $\infty$-algebroid structure.
It appeared first in an explicit form in the PhD of Sylvain Lavau  \cite{lavauphd}, followed by a referred version by C.L.G., Sylvain Lavau and Thomas Strobl in \cite{LLS}, but the authors acknowledge it was discussed several years earlier by Ralph Mayer and, even more extensively, by Chenchang Zhu. Also, Theodore Voronov and his collaborators discussed this result in a slightly different context.
The construction was then reinterpreted by Yaël Frégier and Rigel Juarez-Ojeda in  \cite{fregier2019homotopytheorysingularfoliations} using semi-models category.
Last, it has been  generalized later for arbitrary Lie-Rinehart algebras over a commutative unitary algebra in \cite{CamilleLouis}.

\vspace{.5cm}

\begin{theorems}{Universal Lie $\infty$-algebroid: existence}{existenceuniversal}
Let $\mathcal{F}$ be a singular foliation on a manifold $M$.
\begin{enumerate}
\item 
In the smooth setting, any  geometric resolution $(E_{-\bullet
}, \dd^{(\bullet)}, \rho)$ comes equipped with a Lie $\infty$-algebroid structure whose linear part\footnote{Defined in Remark \ref{rmk:linearpart}.} is $(E_{-\bullet
}, \dd^{(\bullet)}, \rho)$.
\item 
In the real analytic or holomorphic setting, this still holds true, but in a neighborhood of any point of $M$ only.
\end{enumerate}
\end{theorems}
\vspace{0.3cm}

\begin{proof}
The proof goes by recursion. Section \ref{sec:length2} gave a proof when a geometric resolution of length $2$ exists.
In the general case, the first step of the proof consists in making $ E_{-1}$ an almost Lie algebroid, then in showing that  there exists a family of graded symmetric bilinear maps
 $$ \ell_2 (\Gamma(E_{-i}), \Gamma(E_{-j})) \longrightarrow \Gamma(E_{-i-j+1}) $$
 that give back the almost Lie algebroid bracket if $ i=j=1$, and satisfy for all $ a \in \Gamma(E_{-i}), b \in \Gamma (E_{-j}) $, and for every function $f$,
  $$ \dd^{(i+j-1)} \ell_2(a,b)= \ell_2(\dd^{(i)} (a),b)+(-1)^i  \ell_2(a,\dd^{(j)}(b)) $$
  and $ \ell_2(a,fb) = f \ell_2(a,b)  $ unless $ i=1$, in which case this is replaced by:
  $$ \ell_2(a,fb)= f \ell_2(a,b) + \rho(a)[f] \, b .$$
  Then we look at the graded Jacobiator and continue by recursion.
\end{proof}
\vspace{0.5cm}

\begin{definitions}{Universal Lie $\infty$-algebroid: definition}{universal}
We call \emph{universal Lie $\infty $-algebroid\footnote{Its dual is called a \emph{universal $NQ$-manifold of $ \mathcal F$}} of a singular foliation} a Lie $\infty $-algebroid whose linear part is a geometric resolution of $ \mathcal F$.
\end{definitions}

\vspace{.5cm}
Here is an immediate consequence of Theorem \ref{thm:existenceuniversal} and the existence results of Section \ref{sec:Geometric-resolutions}.

\begin{corollary}
    A real analytic or holomorphic singular foliation admits a universal Lie $ \infty$-algebroid in a neighborhood of any of its point. 
    
    A locally real analytic singular foliation admits a universal Lie $ \infty$-algebroid on any relatively compact open subset.
\end{corollary}

The name ``universal'' is justified: it is indeed a universal object in the category of Lie $\infty $-algebroids whose anchor map is valued in $\mathcal F $. 
The arrows of that category are defined to be homotopy classes of morphisms.
\vspace{.5cm}

\begin{theorems}{Universal Lie $\infty$-algebroids deserve the name ``universal''}{areuniversal} \label{th:universal}
Let $\mathcal{F}$ be a singular foliation over smooth $M$. Given, 
\begin{enumerate}
    \item[a)] a Lie $ \infty$-algebroid  $\left(E'_\bullet,(\ell'_k)_{k\geq 1}, \rho'\right)$ that terminates in $\mathcal F$, i.e, $\rho'(\Gamma(E_{-1}'))\subseteq \mathcal{F}$,
    \item[b)] a universal Lie $\infty $-algebroid $ \left(E_\bullet,(\ell_k)_{k\geq 1}, \rho \right)$ of $\mathcal F $,
\end{enumerate} 
then
\begin{enumerate} 
\item there exists a Lie $\infty $-algebroid morphism from $\left(E'_\bullet,(\ell'_k)_{k\geq 1}, \rho'\right)$ to   $ \left(E_\bullet,(\ell_k)_{k\geq 1}, \rho \right)$.
\item and any two such morphisms  are homotopic.  
\end{enumerate}
In the real analytic or complex settings, the same holds, but in a neighborhood of a point only.
\end{theorems}

\vspace{0.5cm}
As for geometric resolutions (see Theorem \ref{thm:uniquegeomresol}), or as for almost Lie algebroids (see discussion after Proposition \ref{prop:EquivAlmost}), Theorem  \ref{thm:areuniversal} means that given a singular foliation $ \mathcal F$, in the category where
\begin{enumerate}
\item objects are Lie $ \infty$-algebroid morphisms that terminates in $ \mathcal F$,
\item and arrows are homotopy classes of Lie $ \infty$-algebroid morphisms,
\end{enumerate}
universal Lie $ \infty$-algebroids of  $ \mathcal F$ are terminal (a.k.a. final, or universal) objects\footnote{As usual, add "in a neighborhood of a point" in the real analytic and complex cases - or redefine morphisms and homotopies as in Lavau's PhD, (see the last chapter in \cite{lavauphd}) as being local ones that glue up to homotopy on the intersections.}.
Here is an immediate corollary of this result, which is valid for any pair of terminal objects in any category.

\vspace{.5cm}
\begin{corollaries}{The universal Lie $\infty$-algebroid is as unique as it can be}{cor:unique}
Two universal Lie $\infty $-algebroids of a 
 smooth singular foliation are homotopy equivalent.
Moreover, the homotopy equivalence between them is unique up to homotopy.

In the real analytic and complexes cases, the same holds, but in a neighborhood of a point only.
\end{corollaries}
\vspace{.5cm}

Theorems \ref{thm:existenceuniversal}-\ref{thm:areuniversal} are proven in \cite{Lavau18,LLS,CamilleLouis} by a finite and constructive recursion\footnote{Provided  that the geometric resolution $ (E_\bullet,(\ell_k)_{k \geq 1}, \rho)$ and  an almost Lie-algebroid bracket on $ E_{-1}$ are given for Theorem \ref{thm:existenceuniversal}, provided a morphism of geometric resolution is given for Theorem \ref{thm:areuniversal}.}. This does not mean that it is easy to construct it. These theorems are therefore constructive. Here are some examples of universal Lie $\infty$-algebroids of singular foliations.

\begin{example}
For a regular foliation $\mathcal F$ on a manifold $M$, the Lie algebroid $T \mathcal{F}\subset T M$ is a universal Lie $\infty$-algebroid of $\mathcal{F}$.
\end{example}
\begin{example}
    Let $\mathcal F\subseteq \mathfrak X(M)$ be a Debord foliation, i.e., $\rho(\Gamma(A))\simeq\mathcal{F}\subset\mathfrak X(M)$ for some Lie algebroid $(A\to M, [\cdot\,,\cdot]_{A}, \rho)$ whose anchor map $\rho\colon A\to TM$ is injective on an open dense subset of $M$. The latter Lie algebroid is a universal Lie $\infty$-algebroid of $\mathcal{F}$.
\end{example}
\begin{example}
    The Lie $2$-algebroid constructed in Proposition \ref{prop:2-Lie-algebroid} over a singular foliation $\mathcal F$ admitting a geometric resolution of length $2$ is a universal Lie $\infty$-algebroid of $\mathcal{F}$.
\end{example}
\begin{example}
We go back to Example \ref{exo:phi-Koszul}. A universal Lie $\infty $-algebroid of $\mathcal{F}_\varphi\subset\mathfrak{X}(V)$ is given on the free resolution $\left(E_{-\bullet} = \wedge^{\bullet+1}V,\mathrm{\dd}=\iota_{\dd \varphi},\rho= -\iota_{\dd \varphi} \right)$ 
by defining the following $n$-ary brackets:
\begin{equation}
\label{eq:nary}
\left\lbrace\partial_{I_{1}},\cdots, \partial_{I_{n}}\right\rbrace_{n}:=\sum_{i_{1}\in I_{1},\ldots,i_{n}\in I_{n}} \epsilon(i_{1},\ldots,i_{n})\varphi_{i_{1}\cdots i_{n}}\partial_{I_{1}^{i_{1}}\bullet\cdots\bullet I_{n}^{i_{n}}};
\end{equation}and the anchor map given for all $ i,j \in \{1, \dots,n\}$ by\begin{equation}
\rho\left(\frac{\partial}{\partial x_i}\wedge\frac{\partial}{\partial x_j}\right) :=\frac{\partial\varphi}{\partial x_j}\frac{\partial}{\partial x_i}-\frac{\partial\varphi}{\partial x_i}\frac{\partial}{\partial x_j}.
\end{equation}Above, for every multi-index $J=\left\lbrace j_1,\ldots ,j_n\right\rbrace\subseteq\left\lbrace 1,\ldots,d\right\rbrace$ of length $n$, $\partial_J$ stands for the $n$-vector field $\frac{\partial}{\partial x_{j_1}}\wedge\cdots\wedge\frac{\partial}{\partial x_{j_n}}$ and $\varphi_{j_{1}\cdots j_{n}}:=\frac{\partial^{n}\varphi}{\partial x_{j_1}\cdots\partial x_{j_n}}$ . Also, $I_{1}\bullet\cdots\bullet I_{n}$ is a multi-index obtained by concatenation of $n$ multi-indices $I_{1},\ldots,I_{n}$. For every $i_1\in I_1,\ldots,i_n\in I_n$,\;$\epsilon(i_1,\ldots,i_n)$ is the signature of the permutation which brings $i_1,\ldots,i_n$ to the first $n$ slots of $I_{1}\bullet\cdots\bullet I_{n}$. Last, for $i_s\in I_s$, we define $I_{s}^{i_s}:=I_s\backslash i_s$. We refer to \cite{LLS}, Example 3.101 or to Section 3.2.1 in \cite{CamilleLouis} for a proof.
\end{example}

\begin{question}
    For the algebraic\footnote{But seen as holomorphic, as in Remark \ref{rmk:algebraic}.} vector fields on $ \mathbb C^n$ tangent to an affine variety $W \subset \mathbb C^d$, what is the universal Lie $ \infty$-algebroid like?
\end{question}

\subsection{The isotropy Lie $ \infty$-algebra of a singular foliation at a point}

\label{sec:monodromy_infinity}

Let $M$ be a real analytic, smooth, or complex manifold.

Given a point $m \in M$, there is a functor\footnote{It is implicitly assume here that the grading of the Lie $ \infty$-structures go from $ -1$ to $-\infty $.}:
 \begin{equation}
     \label{eq:isotropyFUnctor}{\mathrm{Isotropy}}_m \colon \{ \hbox{ Lie $\infty$-algebroids on $M$ }  \}  \longrightarrow \{ \hbox{ Lie $\infty$-algebras } \} 
     \end{equation}
 that we describe in the next lines.
 Then, we apply this functor to the universal Lie $\infty $-algebroids at an arbitrary  point $m$, and explain why the henceforth obtained Lie $\infty $-algebras deserve to be called isotropy Lie $ \infty$-algebras by relating them to AS-isotropy Lie algebras (see Section \ref{sec:isotropy}). 

\vspace{0.5cm}

\subsubsection{ Isotropy functor: Lie $ \infty$-algebroid + Point  $ \mapsto$   Lie $ \infty$-algebra.  }

Let $\left(E_{-\bullet},(\ell_k)_{k \geq 1}, \rho\right)$ be a Lie $\infty$-algebroid with anchor $\rho$ on a manifold\footnote{It suffices the Lie $\infty$-algebroid structure to be defined in a neighborhood of $m$, it is important for the real analytic and complex settings.} $M$.

Consider the graded vector space $ev(E,m)_\bullet$ given by
$$  ev(E,m)_i = \left\{ \begin{array}{ll} 
 {E_{-i}}_{|_m} & \hbox{ for $i \geq 2$}\\ \ker (\rho_{m}) & \hbox{ for $i =1$}\\ 0 & \hbox{ for $i \leq 0$}\end{array} \right.  $$
\vspace{0.5cm}
 
 \begin{propositions}{Specialization of a Lie $\infty$-algebroid at a point}{isotropyfunctor}
  For every point $m\in M$, the $k$-ary bracket of $ \left(E_{-\bullet},(\ell_k)_{k \geq 1}, \rho\right)$ restrict to  $ev(E,m)$,
and equipped the latter with a Lie $\infty$-algebra structure.

This restriction is functorial\footnote{i.e., a (maybe local near $m$) morphism of Lie $ \infty$-algebroids \underline{over the identity of $M$} induces a morphism of Lie $ \infty$-algebras, and homotopic morphisms are mapped to homotopic morphisms.}.
 \end{propositions}

\begin{proof}
We denote the evaluation of a section $s\in \Gamma(E_\bullet)$ at $m$ by $s(m)  $ or $ s_{|_m}$ depending on the context.
We denote the induced brackets on $ ev(E,m)_\bullet$
For every $k\geq 1$, we set
\begin{equation}
    \{x_1,\ldots,x_k\}_k:=\ell_k(s_1,\ldots,s_k)_{|_m}
\end{equation}
for all $x_1,\ldots,x_k\in ev(E,m) $, where $s_1,\ldots,s_k\in \Gamma(E)$ are sections of $E$ such that ${s_i}(m)=x_i$ with $i=1,\ldots,k$.

      For $k\neq 2$, this is well-defined, since $\ell_k$ is linear over functions.
      But it is not so immediate that a $2$-ary bracket is also well-defined.
Let  $(e^i_1,\ldots e^i_{\mathrm{rk}(E_{-i})})$ be a local trivialization of $E_{-i}$ on a neighborhood $\mathcal U$ of the point $m\in M$. For $x_1\in \ker(\rho_{m})$ and $x_2\in {E_{-i}}_{|_m}$ write $$\displaystyle{x_1=\sum_{k=1}^{{\mathrm{rk}(E_{-1})}}}\lambda_ke^{1}_k(m),\quad\displaystyle{x_2=\sum_{k=1}^{{\mathrm{rk}(E_{-i})}}}\mu_ke^{i}_k(m)$$ for some scalars $(\lambda_i)$ in $\mathbb{K}$. The scalars $(\lambda_k)$, $(\mu_k)$ extend to functions $(f_k)$,  $(g_k)$ on $\mathcal{U}$. Therefore, we have $$\{x_1,x_2\}_2=\ell_2(s_1,s_2)_{|_m}$$ with 
$$\displaystyle{s_1=\sum_{k=1}^{{\mathrm{rk}(E_{-1})}}}f_ke^{1}_k,\quad\displaystyle{s_2=\sum_{k=1}^{{\mathrm{rk}(E_{-i})}}}g_ke^{i}_k.$$ If $\widetilde{s}_2$ is another extension of $x_2$, then $(s_2-\widetilde{s}_2)(m)=0$ and this is equivalent to $(g_k-\widetilde{g}_k)(m)=0$ for  $k=1,\ldots, \mathrm{rk}(E_{-i})$. It follows that
\begin{align*}
    \ell_2(s_1,s_2-\widetilde{s}_2)_{|_m}&=\sum_{k=1}^{\mathrm{rk}(E_{-i})}\ell_2\left(s_1,(f_k-\widetilde{g}_k)e^{i}_k\right)_{|_m}\\&=\sum_{k=1}^{\mathrm{rk}(E_{-i})}\cancel{(f_k-\widetilde{g}_k)(m)}\ell_2\left(s_1,e^{i}_k\right)_{|_m}+\cancel{\rho(s_1)_{|_m}[f_k-\widetilde{g}_k]}e^{i}_k\\&=0.
\end{align*}
This proves the claim.
       
Also, the induces brackets $\{\cdots\}_k$ have values in $ev(E,m)_\bullet$ for degree reasons, except maybe for the $2$-ary bracket when applied to elements of degree $-1$ (i.e., elements of the kernel of $\rho_{m}$): in that case the bracket is still in the kernel of $\rho_{m}$ since \begin{align*}
    \rho_{m}(\{x_1,x_2\}_2(m))&=    \rho_{m}(\ell_2(s_1,s_2)_{|_m})\\&=\rho(\ell_2(s_1,s_2))_{|_m}\\&=[\rho(s_1),\rho(s_2)]_{|_m}=0
\end{align*}
In the last line we have used the fact that the Lie bracket of two  vector fields that vanish at $m$ is a vector field that vanishes again at $m$.
Functoriality is left to the reader.

\end{proof}

\begin{remark}
This proof above is elementary, but let us explain how the proof would work when working with $NQ$-manifolds, exploiting the duality od Definition \ref{def:dualNQ}. Let $ (M,\mathcal E,Q)$ be the $NQ$-manifold dual to the Lie $ \infty$-algebroid structure. It is easy to check that the ideal of $ \mathcal E=\Gamma(S(E_\bullet^*))$ generated by $ \mathcal I_m+Q(\mathcal I_m)$, with $ \mathcal I_m$ the ideal of functions vanishing at $m$, is a $Q$-ideal. Hence, the quotient $ \mathcal E/ ((\mathcal I_m+Q(\mathcal I_m))\mathcal E)$ comes equipped with a derivation of degree $ +1$
squaring to $0$. 
It suffices then to check that this quotient is canonically isomorphic to the graded symmetric algebra of the dual of $ ev(E,m)_\bullet$.
This gives a very aesthetic proof of Theorem \ref{thm:isotropyfunctor}, and also gives functoriality in a trivial manner.
\end{remark}

\subsubsection{  Cohomology functor:  Lie $ \infty$-algebra $ \mapsto$ DGLA  }
 Now, there is a second functor
 \begin{equation}\label{eq:cohom} {\mathrm{Cohom}} \colon   \{ \hbox{ Lie $\infty$-algebras } \}  \longrightarrow \{ \hbox{ Differential Graded Lie algebras } \}  \end{equation} 
that consists in noticing that, for a given Lie $\infty $-algebra on $ F_{-1} \oplus F_{-2} \oplus \cdots $ with brackets $ (\ell_k)_{k \geq 1}$, the high Jacobi identities \eqref{eq:higherjacobi} applied for $ n=1,2,3$ imply that
\begin{enumerate}
\item the $1$-ary bracket squares to zero, and therefore turns $ (E_\bullet,\dd)$ into a complex, whose cohomology we denote by  $ \oplus_{i \geq 1} H_{-i}$,
\item the $2$-ary bracket goes to the quotient to define a degree $+1 $ graded symmetric bilinear map $ H_{-i} \times H_{-j} \to H_{-i-j+1} $ that we denote by $ [\cdot\,, \cdot]$.
\item the Jacobiator of the $2$-ary bracket $\ell_2$ being, the Jacobiator of $ [\cdot\,,\cdot]$ is zero.
\end{enumerate}
Hence, $ (H_{-\bullet},  [\cdot\,,\cdot])$ is a graded Lie algebra\footnote{With a graded symmetric degree $ -1$ Lie bracket, which is not the most  convention may be used to. To get a degree $0$ graded skew-symmetric bracket, one has to define a new degree by declaring $ H_{-i}$ to be of degree $ -i+1$ and to replace $ [\cdot, \cdot]$ by $ (a,b) \mapsto (-1)^{i} [a,b]$ for any $b \in E_{-i}$ .}.
This map behaves well with respect to morphisms and homotopies.

\subsubsection{The isotropy DLGA at a point of a singular foliation}

Now, in the particular case we are interested in, the cohomology of the complex $ (ev(E,m)_m,\{\cdot\}_1)$ is the cohomology obtained by evaluating a geometric resolution at $m$.
We saw these cohomologies in Section \ref{sec:geometricresolutionsdef}, more precisely in Corollary \ref{cor:canonicalpoint}.
We called them isotropy spaces of $ \mathcal F$ at $m$, and denote them by $ H^\bullet(\mathcal F,m)$ and showed that it does not depend on the chosen point. 
For any singular foliation that admits a (local) geometric resolution $(E_{-\bullet}, \dd^{(\bullet)}, \rho) $  therefore, by
\begin{enumerate}
\item constructing a universal Lie $ \infty$-algebroid structure on $ (E_{-\bullet}, \dd^{(\bullet)}, \rho)$,
\item then applying the functor ${\mathrm{Isotropy}}_m $ of  Eq.\eqref{eq:isotropyFUnctor}
\item then applying the functor $ {\mathrm{Cohom}}$ of Eq. \eqref{eq:cohom}
\end{enumerate}
one obtains a  differential graded Lie algebra(=DGLA). Moreover, the whole construction is “canonical” in the sense that Theorem \ref{thm:areuniversal}
states that two different choices of a universal Lie $\infty $-algebroid would lead to equivalent DGLA. Now,  Proposition 4.12 in \cite{LLS} states that the situation is in fact even better: the final DGLA is independent of the choice of a universal Lie $\infty $-algebroid. This justifies the next definition.

\begin{definition}
\label{def:isotropyDGLA}
We call the differential Lie algebra structure on $ H^\bullet(\mathcal F,m):= \oplus_{i \geq 1} H^{-i}(\mathcal F,m)$ described above the \emph{isotropy DGLA} of the singular foliation at the point $m$.
\end{definition}

The name is justified by the coincidence of its first term with  Androulidakis and Skandalis isotropy Lie algebra of Section \ref{sec:isotropy}.

\begin{proposition}\label{prop:isotropy}
For every point $m$ in a foliated manifold that admits a (maybe local) geometric resolution, the restriction of the isotropy differential graded Lie algebra to its component of degree $ -1$ is a Lie algebra that coincides with the Androulidakis and Skandalis isotropy Lie algebra $\mathfrak{g}_m(\mathcal F)$.
\end{proposition}
\begin{proof}
For $m\in M$, we construct a Lie algebra is isomorphism $\zeta\colon\frac{\ker(\rho_m)}{\mathrm{im}(\dd^{(2)}_{m})}\rightarrow \mathfrak{g}_m$ as follows: For an element $u\in\ker(\rho_m)$, let $\widetilde{u}$ be an extension of $u$ to a local section on $E_{-1}$. By construction, one has $\rho(\widetilde{u})\in\mathcal{F}(m)$. Let $\widetilde{\rho}_m$ be the surjective linear map defined by
\begin{align*}
    \widetilde{\rho}_m\colon\ker(\rho_m)\longrightarrow \mathfrak{g}_m,\, u\longmapsto [\rho(\widetilde{u})].
\end{align*}
 Since any other extension $\widetilde u$ for $u$ differs from the first one by a section in $\mathcal{I}_m\Gamma(E_{-1})$, the map $\widetilde{\rho}_m$ is well-defined. Surjectivity is due to the fact that every vector field of $\mathcal{F}$ vanishing at $m\in M$ is of the form $\rho(e)$ with  $e$ a (local) section of $E_{-1}$ whose value at $m$ belongs to $\ker(\rho_m)$. In addition, it is not hard to see that $\widetilde{\rho}_m$ is a morphism of brackets.
 
 It remains to show that $\ker(\widetilde{\rho}_m)=\mathrm{im}(\dd^{(2)}_{m})$: let $u \in \ker(\widetilde{\rho}_m)\subset \ker({\rho}_m)$ and $\widetilde{u}$ be a local section of $E_{-1}$ that extends $u$. By definition of $u$, the class of  $\rho(\widetilde{u})$ is zero in $\mathfrak{g}_m$, therefore, there exists some functions $f_i\in \mathcal{I}_m$ and $X_i\in\mathcal{F}, i=1,\ldots, k$,  local generators such that $$\displaystyle{\rho(\widetilde{u})=\sum_{i=1}^kf_iX_i}.$$ This implies that  $$\rho( \widetilde{u}-\sum_{i=1}^kf_ie_i)= 0.$$
where for $i = 1,\ldots,k$, $e_i$ is a (local) section of $E_{-1}$ whose image through $\rho$ is $X_i$. Since $(E_{-\bullet}, \dd^{\bullet}, \rho)$
is a geometric resolution, there exists a (local) section  $q \in \Gamma(E_{-2})$ such that
\begin{equation}\label{eq:geo-iso}
    \widetilde{u}=\sum_{i=1}^kf_ie_i+ \dd^{(2)}q
\end{equation}
 By evaluating Equation \eqref{eq:geo-iso} at $m$, we find out that $u\in\mathrm{im}(\dd^{(2)}_{m})$. Conversely, for $v\in {E_{-2}}_{|_m}$, choose a (local) section $q$ of $E_{-2}$ through $v$. Therefore, $\dd^{(2)}q\in \ker\rho$, is a (local) extension of $\dd_m^{(2)}v\in \mathrm{im}(\dd^{(2)}_{m})$. The image of $\dd_m^{(2)}v$ through $\widetilde{\rho}_m$ is zero. This proves that $\ker(\widetilde{\rho}_m)=\mathrm{im}(\dd^{(2)}_{m})$.
\end{proof}

\subsubsection{The isotropy Lie $ \infty$-algebra.}

 Now, there are more structures that a “simple” differential graded Lie algebra on 
$H^\bullet(\mathcal F,m)$. The $2$-ary bracket
of that structure is in fact the $2$-ary bracket of a Lie $ \infty$-algebra structure whose $1$-ary bracket is zero. 
To obtain it, we proceed as follows.
For any singular foliation that admits a (local) geometric resolution $(E_{-\bullet}, \dd^{(\bullet)}, \rho) $  therefore, by
\begin{enumerate}
\item We must manage to replace it by another geometric resolution which is minimal\footnote{See Definition \ref{def:defgeomresol}.} at $ m$,
\item then constructing a universal Lie $ \infty$-algebroid structure on $ (E_{-\bullet}, \dd^{(\bullet)}, \rho)$,
\item then applying the functor ${\mathrm{Isotropy}}_m $ of  Eq.\eqref{eq:isotropyFUnctor}
\item and then we stop there (we do \emph{not} apply functor ${\mathrm{Cohom}}$ of Eq. \eqref{eq:cohom}).
\end{enumerate}
If the underlying complex $(E_{-\bullet},\ell_1,\rho)$ of $(E,Q)$ is minimal at $m$ then, for  every $i\geq 2$, the
vector space $H^{-i}(\mathcal{F}, m)$ is canonically isomorphic to ${E_{-i}}_{|_m}$. Also, $H^{-1}(\mathcal{F}, m)$ is canonically isomorphic to $\ker (\rho_{m})$. 
Therefore, the outcome of the construction is a Lie $ \infty$-algebra over the “complex” $H^\bullet(\mathcal F,m) $ whose differential is zero. The $2$-ary bracket is the isotropy differential graded Lie algebra bracket. Proposition 4.12 in
\cite{LLS} shows\footnote{See also \cite{CamilleLouis} Section 2.1.3.} that this Lie $ \infty$-algebra is well-defined up to a strict isomorphism on which there are even more constraints. The following definition then makes sense.

\begin{definitions}{}{}
   For any point $m$ of a foliated manifold $ (M,\mathcal F)$ that admits a geometric resolution (at least, near $m$),  the graded vector space $H^\bullet(\mathcal{F},m)$ carries a Lie $\infty$-algebra structure 
called the \emph{isotropy Lie $\infty$-algebra of the singular foliation $\mathcal F$ at $m$}.
\end{definitions}

This structure is used in Proposition \ref{thm:prop:NMRLAclass}.

\chapter{State of the Art and open questions} 

We list in this section several open problems and questions -some are vague, some are precise, some seem major, some seem mere anecdotes.

\section{Existence of Lie algebroids generating a singular foliation}

Let us present the most intriguing open question regarding singular foliations.

Let us start by making the terminology precise.
So far, it was part of the definition of a ``Lie algebroid $ (A \to M, \rho, [\cdot\,, \cdot])$'' that $A \to M$ had to be a finite rank vector bundle by $M$, i.e.,  that $A \to M$ is a vector bundle modeled over a finite dimensional vector space. 
In this section, however, let us distinguish:
\begin{enumerate}
    \item \emph{finite rank Lie algebroids}, i.e., Lie algebroids as defined so far, with $A \to M$ a finite rank vector bundle,
    \item \emph{infinite rank Lie algebroids}, which have precisely the same definition, except that $A \to M$ is now a vector bundle of infinite rank.
\end{enumerate}

As we saw in Section \ref{sec:LieAlgebroidsAreSingFoliation}, for any finite rank Lie algebroid $  (A \to M, \rho, [\cdot\,, \cdot])$, the image of the anchor map $\mathcal F = \rho \left( \Gamma(A) \right) $ is a singular foliation on $M$.

\begin{exo}
Let $   (A \to \mathcal U, [\cdot\,, \cdot], \rho)$ be an  infinite rank Lie algebroid on a smooth manifold. Check that\footnote{The index $ {}_c$ means compactly supported. } $\mathcal F = \rho \left( \Gamma_c(A) \right) $  is
\begin{enumerate}
    \item a $\mathcal C^\infty(M) $-submodule of $ \mathfrak X_c(M)$,
    \item involutive, i.e., $[\mathcal F, \mathcal F] \subset \mathcal F $.
\end{enumerate} 
Extend the result to the complex or real analytic contexts.
\end{exo}

\begin{example}
Here is an example for which $\mathcal F $ is not locally finitely generated as a  $\mathcal C^\infty(M) $-module but still comes from an infinite rank Lie algebroid:
\begin{enumerate}
    \item[a)] $M =\mathbb R $,
    \item[b)] $A$ is the trivial vector bundle with generators indexed $(e_i)_{i\in \mathbb N}$ indexed by $\mathbb N $,
    \item[c)] the anchor map is $\rho(e_i) = \frac{1}{x^i} e^{-\frac{1}{x^2}}\frac{\partial}{\partial x} $,
     \item[d)] the Lie bracket is defined by $[e_i,e_j]= (i-j)  e^{-\frac{1}{x^2}} e_{i+j+1} $. 
\end{enumerate}
\end{example}

Here is a simple open question, that – as far as we know – first appeared in a printed version in Androulidakis and Zambon's \cite{AZ13}.
We were told by several mathematicians, e.g., Rui Loja Fernandes and Georges Skandalis, that the question was already circulating orally in the early 2000s.

\vspace{.5cm}

\begin{questions}{\cite{AZ13}Lie algebroid?}{ques:HiddenLieOid}

Let $\mathcal F $ be a singular foliation on a manifold $M$. Does every point $m$ admit a neighborhood $\mathcal U $ on which there exists a finite rank Lie algebroid $ (A \to \mathcal U, [\cdot\,, \cdot], \rho)$ such that 
 $ \mathcal F= \rho \left( \Gamma(A)\right) $? 
\end{questions}

\vspace{0.4cm}

Here is a slightly more general formulation of the question:

\vspace{0.4cm}

\begin{questions}{\cite{AZ13} Lie algebroid (version II)?}{ques:HiddenLieOid2}
Is any finitely generated singular foliation the image through the anchor map of a finite rank Lie algebroid?
\end{questions}

\vspace{.5cm}

In addition to the local problem, there is also a global ``gluing'' problem. This one only makes sense on the smooth setting.

\begin{question}
\label{ques:glueing}
If a smooth singular foliation is the image of the Lie algebroid on open subsets $\mathcal U_1, \mathcal U_2 $, is it the image of a Lie algebroid on $ \mathcal U_1 \cup \mathcal
U_2$?
\end{question}

Even if we assume both Lie algebroid structures to be defined on the restrictions to $\mathcal U_1 $ and $ \mathcal U_2 $ of the same vector bundle on $\mathcal U_1 \cup \mathcal U_2 $, Question \ref{ques:glueing} remains non-trivial.

\begin{exemple}
Singular foliations whose number of local generators are not globally bounded can not be, globally, the image through the anchor map of a finite rank Lie algebroid. Hence, the Androulidakis and Zambon's ``non-finitely-many-generators'' singular foliation (see Example \ref{exo:infinitestability}) is not the image through the anchor map of a finite rank Lie algebroid on the whole manifold $M=\mathbb R^2 $.
\end{exemple}

\begin{exo} The purpose of this exercise is to show that any finitely generated singular foliation is the image through the anchor map of an infinite rank Lie algebroid.

\begin{enumerate}
    \item Let $X_1, \dots, X_d$ be vector fields on a manifold $M$, and let $\mathfrak g_{free}^d $ be the free Lie algebra with $d$-generators $e_1, \dots, e_d $. Show that exists a unique Lie algebra morphism $ \rho \colon \mathfrak g_{free}^d \to \mathfrak X (M)$ such that $ \rho(e_i) = X_i$.
    \item Assume now that  $X_1, \dots, X_d$  are generators of a singular foliation $\mathcal F $. Use the previous Lie algebra morphism to construct a Lie algebroid structure on the trivial bundle
    $ \mathfrak g_{free}^d \times M \to M$ such that the image of its anchor map is $\mathcal F $.
\end{enumerate}
\end{exo}

Quite a few singular foliations are the image through the anchor map of a Lie algebroid: symplectic foliations of Poisson structures (see Section \ref{sec:hamiltonian}) for instance, or orbits of a Lie algebra action (see Section \ref{sec:LieAlgebroidsAreSingFoliation}). 
Here is an example of a singular foliation of rank $6$ for which no Lie algebroid is known. Notice that this question is understandable by any master degree student- but is still open!

\vspace{.5cm}

\begin{questions}{A frustrating example}{ques:quadratic}
Is the singular foliation of vector fields on $\mathbb R^2 $ vanishing quadratically\footnote{See section \ref{ex:singFolVanish}} at the origin $ 0 $ the image through the anchor map of a finite rank Lie algebroid? 
\end{questions}

\vspace{0.5cm}

Here are other examples of singular foliations for which no finite rank Lie algebroid is known, except in some particular cases:
\begin{enumerate}
    \item vector fields on $\mathbb C^n $ tangent to a given affine variety $W \subset \mathbb C^n $,
    \item vector fields on $\mathbb C^n $ vanishing at every point of an affine variety $W \subset \mathbb C^n $,
    \item  vector fields  $X \in {\mathfrak X} (\mathbb C^n) $ such that $X[\varphi]=0$ for some polynomial function $\varphi \in \mathbb C[x_1, \dots, x_n] $ (see Example \ref{exo:phi-Koszul}).
\end{enumerate}

\begin{exo}
Show that any singular foliation $\mathcal F$ whose rank can be only $k=1,2$ at a given point comes from a Lie algebroid. (\emph{Hint}: construct an almost Lie algebroid of rank $k$ over $\mathcal F $ and show that its Jacobiator has to be trivial). 
\end{exo}

\begin{exo}
Let $\varphi \in \mathbb C[x,y,z]$ be a polynomial function on $\mathbb C^3$.
Check that the following bivector field:
$$ \{x,y\} =\frac{\partial \varphi}{\partial z} \, , \, \{y,z\} =\frac{\partial \varphi}{\partial x} \, , \, \{z,x\} =\frac{\partial \varphi}{\partial y}  $$
is a Poisson bivector field, and that $\varphi  $ is a Casimir function.
Consider the corresponding Lie algebroid on $A= T^* \mathbb C^3 $. Show that the image of its anchor map is a sub-singular foliation of the singular foliation $\mathcal F_\varphi$ of all vector fields  $X \in {\mathfrak X} (\mathbb C^3) $ such that $X [\varphi]=0$. 
Show that if $\varphi $ is weight homogeneous with an isolated singularity at zero, then 
$\rho(\Gamma(A))=\mathcal F_\varphi$. \emph{Hint:} This is done by Anne Pichereau in \cite{zbMATH05043646}.
\end{exo}

\vspace{1cm}
\noindent
{{\textbf{Discussion}}}
Question \ref{thm:ques:HiddenLieOid} may be misleading, in the sense that ``behind'' a singular foliation is a Lie $\infty $-algebroid\footnote{(= $Q$-manifold = dg-manifold)}.  The Lie algebroid, even if there is one, is certainly not unique (one could take the direct product with any Lie algebra for instance). But the universal Lie $ \infty$-algebroid is unique (up to homotopy, see  Corollary \ref{thm:cor:unique}), so that any homotopy invariant information obtained out of a universal Lie $ \infty$-algebroid is canonically attached to the singular foliation.

Moreover, the universal Lie $ \infty$-algebroid itself gives some hints about a possible Lie algebroid that whose image through the anchor map would be the singular foliation. 

\vspace{0.5cm}

It is shown in \cite{LLS} that some singular foliations of rank $r$ are not the image through the anchor map of a Lie algebroid of rank $r$. 
In fact, the following result is shown in Example 4.32 in \cite{LLS}:

\vspace{.5cm}

\begin{propositions}{No minimal rank Lie algebroid}{prop:NMRLAclass}

The singular foliations of all vector fields $X$ on $\mathbb C^4 $ such that $X[\phi]=0 $ with $\phi(z_1,z_2,z_3,z_4)= z_1^3+z_2^3+z_3^3+z_4^3 $:
\begin{enumerate}
    \item has rank $6$ at the origin,
    \item but can not be  the image through the anchor map of a Lie algebroid of rank $6$ on a given neighborhood of the origin.
\end{enumerate}
\end{propositions}

\vspace{.5cm}
This relatively elementary result uses the universal Lie $ \infty$-algebroid. Proposition 4.29 of \cite{LLS} states that if a Lie algebroid of rank $r$ exists in a neighborhood of a leaf reduced to a point, say $m$, then the isotropy Lie $\infty$-algebra\footnote{See Section \ref{sec:monodromy_infinity}} at $m$ admits a minimal model whose $3$-ary bracket vanishes. Now, there are cohomological obstructions to such a cancellation. Here is the exact statement:
\begin{proposition}[\cite{LLS}, Proposition 4.29]
\label{prop:nonminimal}
A singular foliation, defined in a neighborhood of $ 0 \in \mathbb R^n $ and of rank $r$ at this point, which admits a geometric resolution, and for which the $3$-ary bracket of any minimal model of the Lie $\infty $-isotropy Lie algebra at $0$ is not exact as a Chevalley-Eilenberg cocycle for the isotropy Lie algebra at $0$ can not be the image through the anchor map of a Lie algebroid of rank $r$.
\end{proposition}

Let us state a striking corollary of this statement. Let $ X_1, \dots, X_r $ be generators of a singular foliation $\mathcal F $.  There exists (see Exercise \ref{exo:christo1}) Christoffel coefficients, i.e., functions $ c_{ij}^k$ (with $i,j,k=1, \dots, r $ satisfying
 $$  [X_i,X_j] = \sum_{k=1}^r c_{ij}^k X_k $$
 but those are not unique, since there are relations between the generators. Without any loss of generality, we can assume 
 \begin{equation}
     \label{eq:skew}
     c_{ij}^k = - c_{ij}^k, 
 \end{equation}
 and, since the Jacobi identity holds, we have:
  \begin{align*} 0 = [X_i,[X_j, X_k]] + c.p._{i,j,k} \\ = \sum_{a=1}^r \left( X_i[c_{jk}^a] + \sum_{b=1}^r c_{ij}^b c_{bk}^a +  c.p. \right) X_a  
   \end{align*}
If for every $a \in \{1, \dots,r\}$, 
 \begin{equation}
     \label{eq:JacobiatorTermbyTerm}
   X_i[c_{jk}^a] + \sum_{b=1}^r c_{ij}^b c_{bk}^a \, + \,  c.p.(i,j,k) =0 
  \end{equation}
 then there exists a Lie algebroid of rank $r$ whose image through the anchor map is $\mathcal F $: the Lie algebroid on a trivial bundle of rank $r$ whose bracket is given by
  $$   [e_i,e_j] = \sum_{k=1}^r c_{ij}^k e_k$$
  and whose anchor is $ \rho(e_k)=X_k$ for all $k$.
  Proposition \ref{prop:nonminimal} (i.e., Proposition 4.29 in  \cite{LLS}) explains that, if the isotropy Lie $\infty$-algebra at a point satisfies the cohomological condition linked to its $3$-ary bracket, then there is no way that coefficients $ c_{ij}^k$ could be found that satisfy both Equation \eqref{eq:skew} and Equation \eqref{eq:JacobiatorTermbyTerm}. 
  
  \vspace{.5cm}

\noindent
{\textbf{Another relation between the universal Lie $\infty$-algebroid and a  Lie algebroid over $\mathcal F$ .}}
Here is a known result, which makes more precise the claims made in Theorem \ref{th:universal}.

\begin{proposition}[
\cite{CamilleLouis}, Proposition 2.3]
If a Lie algebroid $A$ over a singular foliation $\mathcal F $ exists, and if $ \mathcal F$ admits a geometric resolution,  then there exists a universal Lie $\infty$-algebroid $(\mathcal E, Q) $ constructed with the help of a geometric resolution $(E_{-i})_{i \geq 1} $ that satisfies the additional conditions:
\begin{enumerate}
\item
$E_{-1}=A$ 
\item the $2$-ary bracket, restricted to $A$, is the bracket of $A$
\item  the restriction to $ E_{-1}=A$ of all the $n$-ary brackets are $0$ for $n \geq 3 $. 
\end{enumerate}
\end{proposition}

This proposition makes the next question a natural one:
\vspace{0.5cm}

\begin{question}

\label{ques:IfLieExists}

If a singular foliation \emph{(i)} admits a geometric resolution and \emph{(ii)} is the image through the anchor map of a Lie algebroid, does it admit a universal Lie $\infty$-algebroid for which all $n$-ary brackets are zero for $n \geq 3$? 
\end{question}

Is the answer to Question \ref{thm:ques:HiddenLieOid} yes or no?  Karandeep Jandu Singh's \cite{Singh} goes in the opposite direction, by showing that the $3$-ary brackets need happen not to be zero for some singular foliation associated to the symplectic Lie algebra. However, this does not answer the question negatively, since other choices could be made in the construction at early steps. 

\vspace{0.5cm}

We have the following conjecture, which rather goes in the direction ``the Lie algebroid seems to exist'' but certainly does not prove it, and leaves room for counter-examples.
If true, then finding negative answers to Question \ref{thm:ques:HiddenLieOid2} will be a hard task.
\vspace{.5cm}

\begin{conjectures}{The isotropy Lie $ \infty$-algebra at a point is formal.}{isotropyInftyAlgebra}
At every point, the isotropy Lie $\infty $-algebra of a singular foliation that admits a geometric resolution is homotopy equivalent to a finite dimensional differential graded Lie algebra.
\end{conjectures}
\vspace{.5cm}
\noindent 
If true, the conjecture implies, for instance, that it is not possible to have a homotopy Lie $\infty $-algebra of a singular foliation at a point isomorphic to the so-called string Lie $2$-algebra. i.e.,  the Lie $2$-algebra of the form 
 $$ \mathfrak g = \mathfrak g_{-2} \oplus \mathfrak g_{-1} $$
 with $ \mathfrak g_{-1}  $ a semi-simple Lie algebra, $\mathfrak g_{-2} = \mathbb R$, the $1$-ary bracket equal to $0$, the $2$-ary bracket equal to the Lie algebra bracket on $ \mathfrak g_{-1}$,  and the $3$-ary bracket:
  $$ \wedge^3 \mathfrak g_{-1} \to \mathfrak g_{-2}$$
  given by the Cartan $3$-form. The statement comes from the fact that the latter is not homotopy equivalent to a finite dimensional differential graded Lie algebra.

\section{About geometric resolutions: when do they exist?}

When do geometric resolutions exist, at least in the neighborhood of every point? 
In the complex or real analytic setting, they exist locally for coherent sheaves, by the classical Hilbert's syzygy theorem.

In the smooth case,
it is relatively simple to construct singular foliations that do not admit geometric resolutions. For instance, we can look at the following function $f:\mathbb R\to \mathbb R$ vanishing on $\{x\leq 0\}$ in $\mathbb R$.
$$
f(x)=\begin{cases}
    e^{\frac{-1}{x^2}} & \text{if } x>0 \\
    0 & \text{if } x\leq 0 
\end{cases}
$$ 
Then $\mathcal F=\mathcal C^\infty_c(M)\cdot f\cdot \partial_x$ defines a singular foliation, which has regular leaves of dimension 0 (=the points $\{x\}$ for $x<0$) and one regular leaf of dimension 1 (=the half-line $ \{x>0\}$). Hence, by Corollary \ref{prop:alternatesum}, this foliation can not admit a geometric resolution.\\

However, we saw in the discussion around theorem \ref{th:universal} that geometric resolutions of a singular foliation exist for locally real analytic singular foliations on a relatively compact open subset. This point was proven in \cite{LLS}, Theorem 2.4. Can one make something more general?
A very natural condition to be imposed on a singular foliation $\mathcal F$ is that it forms a subspace of $\mathfrak X_c(M)$ which is \underline{closed} with respect to the Fréchet topology. Any locally real analytical singular foliation is closed, however not every closed singular foliation is locally real analytic,
hence is a counter-example.

\begin{example} Consider $M=\mathbb R^2$ and the ideal $I=(f)\subset \mathcal C^\infty(M)$ generated by $f(x,y)=y^2-e^{-\frac{1}{x^2}}$. By Example 4.8 in \cite{MR240826} this ideal is closed, hence the foliation $\mathcal F=f\cdot \mathfrak X_c(M)$ is also closed. It can however not be real analytically generated near the origin, since there the 0-dimensional leaf looks like two curves intersecting flatly (i.e., with all derivatives colinear). 
At a consequence, vector fields of the form $  f(x,y) \frac{\partial}{\partial x}$ is not locally real analytic, although it is Fréchet closed.
\end{example}

Since locally real-analytic singular foliations locally admit geometric resolutions and hence, and since Fréchet-closed singular foliations being a generalization of locally real analytic ones, the following question is very natural.

\vspace{0.5cm}
\begin{questions}{Closed singular foliations admit geometric resolutions?}{}
Does any Fréchet-closed singular foliation locally admit a geometric resolution?
\end{questions}
\vspace{0.5cm}

If the answer is “yes”, then Fréchet-closed singular foliation admit a universal Lie $\infty$-algebroids in a neighborhood of every point by Theorem \ref{th:universal}.

\section{Molino-Atiyah classes}

Let us first recall the construction of the Molino class of a regular foliation $\mathcal F $. 
As its name indicates, the Molino class is a class in some cohomology: we first describe the cohomology to which it belongs.

Let $\mathcal F $ be a \underline{regular} foliation on $M$, with tangent bundle $ T \mathcal F \subset TM$. 
Notice that $\mathcal F = \Gamma_c(T\mathcal F) $.
We assume the reader is familiar with the language of Lie algebroids: the presentation is influenced by Zhou Chen, Mathieu Stiénon, Ping Xu's \cite{zbMATH06537581}.
 \begin{enumerate}
    \item The tangent bundle $T \mathcal F $ is a Lie sub-algebroid of the tangent Lie algebroid $TM$, whose anchor map is the inclusion $T\mathcal F \hookrightarrow TM $.
    \item Consider the  \emph{normal bundle} $ N_\mathcal F := TM/T\mathcal F$. Denote by $u \mapsto \overline{u} $ the natural projection
     $TM \longrightarrow N_\mathcal F = TM/T\mathcal F$. The normal bundle comes equipped with a $T\mathcal F$-connection\footnote{For the reader non-familiar with Lie algebroids, items 2) and 3) mean that the restriction of the normal bundle to a leaf is equipped with a natural flat connection.}, called the \emph{Bott connection}, and defined by:
     $$\nabla_X^{Bott} \overline{u} = \overline{[X,u]} $$
     for all $ X \in \mathcal F$ and $u \in \Gamma(TM) $. 
     \item It follows from the Jacobi identity for vector fields on $M$ that the Bott connection is a flat connection. As a consequence $X \mapsto \nabla_X^{Bott}  $ turns $N_\mathcal F$ into a Lie algebroid representation of $T \mathcal F $. 
     \item The dual of a Lie algebroid representation of $ T\mathcal F$, and the tensor or symmetric products of two Lie algebroid representations of $ T\mathcal F$ being  Lie algebroid  representations of $T \mathcal F $ again,  the vector bundle $S^2 N_\mathcal F^* \otimes N_\mathcal F $ (i.e., the vector bundle of symmetric bilinear maps from the normal bundle to itself) is a  Lie algebroid representation of $ T\mathcal F$.
\end{enumerate}
The Molino class is a cohomology class of degree $1$ for the Chevalley-Eilenberg cohomology\footnote{For the reader non-familiar with Lie algebroid: the restriction of $S^2 N^*_\mathcal F \otimes N_\mathcal F $ to a leaf comes with a flat connection.} of $ T\mathcal F$ valued in the module $S^2 N^*_\mathcal F \otimes N_\mathcal F $. By construction, this class has to be represented by a vector bundle morphism: $$ \alpha \colon T\mathcal F \otimes S^2 N_\mathcal F  \longrightarrow N_\mathcal F, $$
which has to satisfy (in-order to be a closed-cocycle):
 \begin{equation}
 \label{eq:tobeclosed}
 \alpha([X,Y],u,v) =  \nabla^{Bott}_X \alpha(Y,u,v) - \alpha(Y,\nabla^{Bott}_X u, v) - \alpha(Y, u, \nabla^{Bott}_X v) - \left( \xymatrix{ X \ar@{<->}[r] & Y }  \right).  
 \end{equation}
Let us now construct the Molino class for a regular foliation.
\begin{enumerate}
    \item Consider a $TM$-connection\footnote{i.e., a linear connection is the usual sense} $\nabla $ on $N_\mathcal F $:
     $$ (X,\overline{u}) \mapsto \nabla_X \overline{u} $$
    whose restriction to $\mathcal F \times \Gamma(N_\mathcal F) $ is the Bott connection, i.e., such that  for all $X \in \mathcal F$:
    $$ \nabla_X \overline{u} = \nabla_X^{Bott} \overline{u} . $$
    \begin{enumerate}
        \item  Such connections always exist.
        \item Without any loss of generality, we can assume that its torsion is zero. The \emph{torsion} is the vector bundle morphism defined by   
         $$ \begin{array}{llll}  T^\nabla  \colon &\wedge^2 TM & \rightarrow & TM \\ & (X,Y) &\mapsto & \overline{ \nabla_X \overline{Y} - \nabla_Y \overline{X} - [X,Y]}\end{array}.$$
    \end{enumerate}
From now on,  we will assume the torsion to be zero.
    \item Consider the curvature $\kappa^\nabla $ of such a connection $\nabla $. By construction, $\kappa^\nabla$ is a vector bundle morphism
     $$ \kappa^\nabla \colon TM \wedge TM \otimes N_\mathcal F \longrightarrow N_\mathcal F .$$
     \item Since the Bott connection is flat, for any $X_1,X_2 \in \mathcal F $,
     $ \kappa^\nabla (X_1,X_2)=0$ if both $ X_1$ and $ X_2$ are in $ \mathcal F$.
 This means that for every $X \in \mathcal F$,  the vector bundle morphism  $\mathfrak i_X \kappa^\nabla \colon  Y \longrightarrow \kappa^\nabla(X, Y) $ vanishes as soon as $Y \in \mathcal F$. It therefore can be seen as a vector bundle morphism $\overline{\mathfrak i_X \kappa^\nabla } \colon N_\mathcal F \otimes N_\mathcal F \longrightarrow N_\mathcal F  $.
     \item The map $X \mapsto \overline{\mathfrak i_X \kappa^\nabla}$ can therefore be seen as a vector bundle morphism from $T\mathcal F $ to $ N_\mathcal F \otimes N_\mathcal F \longrightarrow N_\mathcal F$,
     \item i.e., it can be seen as a vector bundle morphism $T\mathcal F \otimes N_\mathcal F \otimes N_\mathcal F \longrightarrow N_\mathcal F$.
     \item We leave it to the reader to check that the vanishing of the torsion implies that 
     $X \mapsto \overline{\mathfrak i_X \kappa^\nabla}$ is symmetric  in the two last variables $N_\mathcal F  $, and is indeed a vector bundle morphism $$\alpha^\nabla \colon T\mathcal F \otimes S^2 N_\mathcal F \longrightarrow N_\mathcal F.$$
\end{enumerate}
Now, it is a direct computation that the Bianchi identity implies that $\alpha $ satisfies \eqref{eq:tobeclosed} above, and is therefore a cocycle of the Chevalley-Eilenberg cohomology of $T \mathcal F $ is the module $S^2 N^*_\mathcal F \otimes N_\mathcal F$, called the \emph{Molino cocycle of the torsion-free connection $\nabla $}. It can be shown that different choices of connections  $\nabla $ would give the same class in cohomology. This completes the definition of the \emph{Molino class}, which is also called \emph{Atiyah class} (see \cite{zbMATH06537581,zbMATH07369649} for a unification of the Molino class with the Atiyah class in complex geometry -hence the double name. See \cite{zbMATH03334884} for an original construction -in French- by Pierre Molino).
The previous discussion has established the following fact:

\begin{proposition}
\cite{zbMATH03334884}
The Molino class is the obstruction to the existence of an extension of the Bott connection whose curvature $2$-form is zero as soon as one element tangent to the foliation is applied to it. 
\end{proposition}

\vspace{0.5cm}

\begin{questions}{Molino class and meaning?}{ques:Molino}
What is the equivalent of the Molino (also called Atiyah) class for a singular foliation? And what is its geometrical meaning?
\end{questions}

\vspace{0.5cm}

Let us state a few points.
\begin{enumerate}
    \item The Bott connection has a natural extension to the singular case, see \cite{zbMATH07568450,zbMATH07733867}:
    \begin{enumerate}
        \item The formula $(X,\overline{u}) \mapsto \overline{[X,u]}  $ defines a flat Lie-Rinehart connection of $\mathcal F $ on the $\mathcal C^\infty(M) $-module $\mathfrak X (M) / \mathcal F $.
        \item The adjoint representation ``up to homotopy'' of any universal Lie $\infty $-algebroid of $ \mathcal F$ is a flat Lie $\infty $-algebroid connection on a geometric resolution of the $\mathcal C^\infty(M) $-module $\mathfrak X (M) / \mathcal F $ that can also be understood as a generalization of the Bott connection.
        \item The Molino class is an instance of Atiyah classes of Lie algebroid pairs.
    \end{enumerate}
    
    \item Geometrically, the vanishing of the Molino class of a regular foliation has several consequences.
     \begin{enumerate}
         \item For any leaf $L$, and any $x \in L$, the holonomy:
          $$ Hol_{\mathcal F, L, l} : \pi_1 (x,L) \longrightarrow Diff_0(N_\mathcal F|_x) $$
          valued in germs at $0$ of diffeomorphisms of the normal bundle. If the Molino class vanishes, the holonomy is linearizable, i.e., the group morphism $ Hol_{\mathcal F, L, l}$ can be assumed to be valued in linear invertible endomorphisms of $N_\mathcal F|_x $. 
          See, e.g., Theorem 8.5 in \cite{zbMATH07179019}.
         \item We say that two paths $\gamma_1, \gamma_2 \colon [0,1] \to M , \gamma_2 $ are $\mathcal F$-related if there exists $F:[0,1]^2 \to M $ such that $ F(t,0)= \gamma_0(t), F(t,1)= \gamma_1(t) $ and such that for every $t \in [0,1]$, the map $ s \mapsto F(t,s) $ is in a fixed leaf.  
         Notice that parallel transportation for $ \nabla$ along curves of the form $ s \mapsto F(t,s) $  is simply parallel transportation with respect to the Bott-connection.
          If $ \nabla$ is a connection such that the Molino cocycle vanishes, the curvature of $\nabla $ vanishes as soon as a vector tangent to $ \mathcal F$ is applied to it, and in particular on the image of $F$. As the consequence, the following diagram is commutative:
         $$ \xymatrix{ N_\mathcal F |_{\gamma_0(0)}  \ar[rr]^{\Phi^{bott}(F)}\ar[d] && N_\mathcal F |_{\gamma_1(0)} \ar[d] \\ N_\mathcal F |_{\gamma_0(1)} \ar[rr]_{\Phi^{bott}(F)} && N_\mathcal F |_{\gamma_1(1)} }  $$
         where horizontal lines are parallel transportation for the Bott connection, and vertical lines are parallel transportation with respect to $\nabla $ along $\gamma_0 $ and $\gamma_1 $.
      \end{enumerate}

It is not obvious to see what the equivalent of the previous points for a singular leaf is.

\end{enumerate}

    \begin{question}
    Assuming it has been defined, what is the geometrical meaning of the vanishing of the (to be constructed) Molino class for a singular foliation?
    \end{question}

A recent article by Seokbong Seol  \cite{seol2024atiyahclassdgmanifolds} may lead to conjecture that if  the Molino class (to be defined)  vanishes, then the foliation has to be a regular foliation.

\section{Miscellaneous}

Here is a ``potpourri'' of several questions, mostly anecdotal at first sight, but to which we have no immediate answer.

Yahya Turki \cite{zbMATH06783221} suggested the following notion: we say that a bivector field $\pi\in \Gamma(\wedge^2 TM )$ is \emph{foliated} if $\pi^\sharp (\Omega^1(M ))$ is closed under the Lie bracket, i.e.,  is a singular foliation. 

\begin{example}
Poisson bivector fields, but also twisted-Poisson bi-vector fields, are examples.  Yahya Turki \cite{zbMATH06783221} gave examples of foliated bivector fields that are not of this type, but proved that they are twisted Poisson near any one of their regular points (= points in a neighborhood of which $\pi^\# $ has constant rank). 
\end{example}

\begin{question}[Foliated bivector fields]{}
Let $\pi$ be a foliated bivector field. 
Can a Lie algebroid structure with anchor map $\pi^\# $ be constructed on $T^* M $?
\end{question}

It is known that $T^*M$ comes equipped
with a Lie algebroid structure  with anchor $\pi^\sharp\colon T^*M \longrightarrow T M$ when $\pi $ is twisted Poisson \cite{zbMATH01860769}, so the question makes sense.

\vspace{0.5cm}

Sébastien Michéa asked if for any smooth Poisson structure $\pi $ on $\mathbb R^n $, there is another structure $\pi' $ on $\mathbb R^n $ which coincides with $\pi $ in a neighborhood of $0$ and vanishes outside a compact subset. The corresponding question for singular foliations is much easier:

\begin{exo}
Given a smooth singular foliation $\mathcal F $ on $\mathbb R^n $, show that there exist  another singular foliation $\mathcal F' $ on $\mathbb R^n $ which coincides with $\mathcal F $ in a neighborhood of $0$ and vanishes outside a compact subset.
\end{exo}

Here is however, a more delicate question:

\begin{question}
Given a smooth singular foliation $\mathcal F $ on $\mathbb R^n $ such that all regular point have rank $r$, does  there exist  another singular foliation $\mathcal F' $ on $\mathbb R^n $ such that all regular point have rank $r$, which coincides with $\mathcal F $ on the open ball $\sum_{i=1}^n x_i^2 < 1 $, but which is made of vector fields all tangent to the sphere $\sum_{i=1}^n x_i^2 =1 $?
\end{question}

\section{Linearization} 
\label{sec:linearization}

Can we enlarge the classical theorems (i.e., Conn's \cite{Conn1,Conn2} or Zung's \cite{zbMATH05137698},  Crainic-Fernandes' \cite{zbMATH05960681}) about linearizations of Poisson structures or Lie algebroid actions or Lie groupoid actions to the context of singular foliations and of its holonomy groupoid?

The previous linearization theorems all share the same logic. There are first relatively easy results whose patterns are:
$$ \hbox{ Fixed point } + \hbox{ Semi-simple } \Longrightarrow \hbox{ Formally Linearizable }  .$$
For instance, it is not so complicated to show that if a Lie algebroid $(A,\rho,[\cdot\,, \cdot]) $ admits a point $m$ where  $ \rho_m=0$ and the isotropy Lie algebra $\mathfrak g_m = A_m $ is semi-simple, then the Lie algebroid is formally equivalent to the transformation Lie algebroid $ \mathfrak g_m \times T_m M \to T_m M$ for some action of $\mathfrak g_m $ by linear endomorphisms $T_m M $.
 Dominique Cerveau \cite{Cerveau} has a result of this type for singular foliations: it says that if the isotropy Lie algebra at a point is a semi-simple Lie algebra, then the singular foliation is formally equivalent to the one associated to a linear action of this Lie algebra. See \cite{Ryvkin2} Proposition 1.12 for an enlarged version.

Beyond these (relatively easy) results, there are then much more difficult results whose patterns are:
$$ \hbox{ Fixed point } + \hbox{ Compact (and semi-simple) } \Longrightarrow \hbox{ Locally Linear }  $$
The difficulty consists in going from ``formal'' to ``local''.
Recall (see Exercise \ref{exo:action}) that if the leaf through  a point $m$ is reduced to $ \{m\}$, or equivalently if $T_m \mathcal F=0 $, then there is a natural action of the isotropy Lie algebroid $\mathfrak g_m(\mathcal F) $ on $ T_m M$.

\vspace{.5cm}

\begin{questions}{Extend Zung's linearization \cite{zbMATH05137698} to singular foliations}{ques:Zung}
Let $ \mathcal F$ be a singular foliation on a smooth manifold $M$. Le $m$ be  point such that $ T_m \mathcal F=0$.

Assume the isotropy Lie algebra $\mathfrak g_m(\mathcal F) $ of $\mathcal F $ at $m$ is semi-simple of compact type. 
Then, is there a saturated\footnote{i.e., a leaf that has non-empty intersection with that neighborhood is contained in it.} neighborhood of $m$ on which $\mathcal F $ is isomorphic to a saturated neighborhood of $0$ for the singular foliation associated to the natural representation of the isotropy Lie algebra $\mathfrak g_m (\mathcal F)$ on $T_m M $?
\end{questions}

\vspace{0.5cm}
\noindent
In other words, we want to prove that the short exact sequence 
 $$ \xymatrix{ \mathcal I_m \mathcal F \, \, \ar@{^{(}->}[rr] & &\mathcal F \ar@{->>}[rr] && \mathfrak g_m(\mathcal F) } $$
splits with a section $\sigma $: 
$$ \xymatrix{ \mathcal I_m \mathcal F \, \, \ar@{^{(}->}[rr] &&\mathcal F \ar@{->>}[rr] && \mathfrak g_m(\mathcal F) \ar@/_15pt/@{.>}[ll]_{\sigma}} $$
which is a Lie algebra morphism, at least in a neighborhood of $m$. 
 Since any action of a  semi-simple Lie algebra of compact type is linearizable near a fixed point, this section $\sigma $ may be seen as being an action of $\mathfrak g_m(\mathcal F) $ on the vector space $T_m M$.

We could of course enlarge these questions to neighborhood of leaves. 
 Again, the formal case is relatively easy: for  instance it has been proven that \cite{Ryvkin2} that Levi-Malcev style theorems hold: those are formal linearization theorems in a neighborhood of a leaf. Of course, only the semi-simple part of the holonomy Lie algebroid (defined in that article) is formally linearizable.  For Lie algebroids or Poisson structures, several authors e.g., \cite{zbMATH06275355,zbMATH06705120, zbMATH06488204} have proven recently several linearizations or normal form theorems in neighborhood of leaves of Lie algebroids or singular foliations: Pretty much any one of these theorems admit a natural generalization for singular foliations.

\vspace{.5cm}

There are similar questions about the holonomy groupoid \cite{AS}.
Recall that it is a topological groupoid, although it is not a Lie groupoid (but each fiber is a variety \cite{Debord}). The topology is the push-forward topology of any atlas of bisections that define it. It makes sense, therefore, to speak of a singular foliation $\mathcal F$ whose holonomy groupoid  $\mathrm{Hol}(\mathcal F) $ is proper: it is a 
singular foliation for which 
 $$ (s,t) \colon \mathrm{Hol}(\mathcal F)  \longrightarrow M \times M $$
 is a proper map.
 
 \begin{definition}
 We say that a singular foliation $\mathcal F $ is \emph{proper} if $ \mathrm{Hol}(\mathcal F)$ is a proper topological groupoid.
 \end{definition}

\begin{example}
Consider a proper groupoid $\Gamma \rightrightarrows M$, e.g., the action groupoid associated to an action of a compact group on a manifold.
Then the basic singular foliation\footnote{i.e., $\mathcal F = \rho(\Gamma(A)) $
with $(A,\rho, [\cdot\,, \cdot]) $ the Lie algebroid of $G \rightrightarrows M$} is a proper singular foliation (since $G $ itself is an atlas, see Proposition \ref{ex:Atlas-Groupoid}).
\end{example}

Proper groupoids have very strong linearization properties. Here is a theorem by Nguyen Tien Zung\footnote{Recall that for every fixed point $m \in M$ of a Lie groupoid (i.e., any point for which $t (s^{-1}(m))=\{m\}$), the isotropy group at $m$ acts naturally by linear automorphisms of the tangent space $T_m M$}:

\begin{theorem}
\label{th:ZungLinearisation}
\cite{zbMATH05137698}
Consider a proper Lie groupoid $\Gamma \rightrightarrows M $. Every fixed point\footnote{I.e., the $ \Gamma$ orbit through $m$ is $ \{m\}$} $m\in M$ admits a saturated neighborhood $\mathcal U $ on which  the restriction of $ \Gamma$ is isomorphic, as a Lie groupoid, to a transformation groupoid of the action of
the compact isotropy group $G_m$ on the tangent space $V = T_m M$.
\end{theorem}

It is therefore very natural to guess that the following result variation of the previous question should be true.
The action of holonomy Lie algebra $\mathfrak g_m(\mathcal F) $  at $m$ on the tangent space $T_m M$ is defined in Exercise \ref{exo:action}.

\vspace{.5cm}
\begin{questions}{Extend Zung's linearization to SF}{ques:ZunG2}
Consider a  singular foliation $\mathcal F$ on $M$ whose holonomy groupoid $\mathrm{Hol}(\mathcal F) \rightrightarrows M $ is proper.

Is it true that every point $m\in M$ such that\footnote{I.e., the leaf through $m$ is $ \{m\}$.} $ T_m \mathcal F=0$ admits a saturated neighborhood\footnote{i.e., a leaf intersecting this neighborhood entirely belongs to that neighborhood.} $\mathcal U $ on which  the restriction of $ \mathcal F$ is isomorphic, as a singular foliation, to a saturated neighborhood of $0$ for the singular foliation associated to the action of the holonomy Lie algebra $\mathfrak g_m(\mathcal F) $ on the tangent space $T_m M$?

If one replaces $ m$ by a leaf $ L$, then the question is still valid,  upon replacing the action of  $\mathfrak g_m(\mathcal F) $  on $T_m M$ by the natural action\footnote{Considered by Androulidakis and Zambon \cite{AZ13}} of the holonomy Lie algebroid of $ L$ on the normal bundle of $L$ in $M$.
\end{questions}
\vspace{.5cm}

For regular foliations, properness of the holonomy Lie groupoid implies, for instance, that every leaf has a saturated neighborhood on which the holonomy map is by linear automorphisms of a finite group, so that the answer to the question is affirmative in this case.
Also, for singular foliations coming from Poisson manifolds of compact type \cite{compactTypes}, the answer to Question \ref{thm:ques:ZunG2} is ``yes''.
Notice also that it has been proven in \cite{resolGroupoid} that singular foliations arising from a compact Lie groupoid can be made a regular foliation by finitely many blow-up operations of its most singular leaves. It would be interesting to generalize this result to a singular foliation whose holonomy groupoid is compact: a positive answer to the previous question should do it.

\section{Longitudinal differential operators}

For $(M,\mathcal F) $ a singular foliation, we call \emph{longitudinal differential operator} any linear combination of operators of the form
\begin{equation}\label{eq:monomials} 
\begin{array}{rcl} \mathcal C^\infty(M)&\to & \mathcal C^\infty(M)\\F &\mapsto &X_1 \circ \dots \circ X_k [F] \end{array}\end{equation}
with $ X_1, \dots, X_k \in \mathcal F$. 
We denote by $ {\mathrm{Diff}}(\mathcal F)$ the algebra of longitudinal differential operators.

For a regular foliation\footnote{Or, more generally, a Debord foliation},  longitudinal differential operators coincide with the universal enveloping algebra of the Lie algebroid $ T\mathcal F$. This is \underline{not} the case for a generic singular foliation.

 Let us first define the universal enveloping algebra of a singular foliation\footnote{One recognizes here the construction of the universal enveloping algebra of a Lie-Rinehart algebra by, see e.g.,  \cite{zbMATH01287585}, \cite{moerdijk2008universalenvelopingalgebralierinehart} or \cite{zbMATH07453484}. See also \cite{maakestad2022envelopingalgebraliealgebra}.}. 
\begin{enumerate}
\item 
To start with, consider the universal enveloping algebra $U(\mathcal F) $ of the \emph{Lie} algebra $\mathcal F $. 

\item Now,  divide $ U(\mathcal F)$ by the ideal generated  by: 
  \begin{equation}\label{eq:defideals}   X \cdot (f Y) - (fX) \cdot Y - (X[f] Y) \end{equation}
  where $X,Y \in \mathcal F $ and $f \in \mathcal C^\infty(M)$.
  We denote the quotient by $ \mathcal U(\mathcal F) $ and call it the \emph{universal enveloping algebra of the singular foliation $\mathcal F $}. 
  \end{enumerate}
  There is a natural algebra morphism from  $U(\mathcal F) $ to $ {\mathrm{Diff}}(\mathcal F) $ defined by 
\begin{equation}\label{eq:realization} \begin{array}{rcl} U(\mathcal F) & \to & {\mathrm{Diff}}(\mathcal F)\\ X_1 \cdot \dots \cdot X_k &\mapsto& 
X_1 \circ \dots \circ X_k \end{array}
 \end{equation}
  The previously defined map \eqref{eq:realization} being equal to $0$ on the ideal generated by the expressions \eqref{eq:defideals}, it goes to the quotient to define a surjective algebra morphism
   $$   \xymatrix{\mathcal{U}(\mathcal F) \ar@{->>}[r]&  {\mathrm{Diff}}(\mathcal F) \\ P \ar[r]& \underline{P} } .$$
   that we call \emph{realization of $ \mathcal{U}(\mathcal F)$.}

   \begin{exo}
   Show that the realization of $ \mathcal{U}(\mathcal F)$ is injective if $ \mathcal F$ is a Debord singular foliation. Show that it is not injective in general. \emph{Hint:} The singular foliation of all vector fields on X on $\mathbb R^3$ such that $ X[x^2+y^2+ z^2]=0$ provides a counter-example.
\end{exo}

Elements in  $ {\mathrm{Diff}}(\mathcal F) $ 
 of the form \eqref{eq:monomials} are called \emph{monomials of degree $k$}, and we say that a longitudinal differential operator is of degree $\leq k $ if it is a sum of monomials of degree $\leq k $. The degree defines an increasing filtration on the algebra $ {\mathrm{Diff}}(\mathcal F) $, making it a filtered algebra  $  \left(\mathrm{Diff}^{\leq k}(\mathcal F) \right)_{k \geq 0}$ . Similarly, the algebra $\mathcal U(\mathcal F) $ comes equipped with a filtration $  \left(\mathcal U^{\leq k} (\mathcal F)\right)_{k \geq 0}$ defined in a same  manner. The realization \eqref{eq:realization}  is a morphism of filtered algebras.

\vspace{0.5cm}
Let us address the following question (which is not an open question, as we will see).
\vspace{0.5cm}

\begin{questions}{}{symbol}
What is the symbol of a longitudinal differential operator?
\end{questions}

\vspace{0.5cm}

Androulidakis and Skandalis in \cite{zbMATH05843246} gave an answer that is too involved to be dealt with here, using the $C^*$-algebra of half densities of the holonomy Lie groupoid.
Later on, Mohsen, then Androulidakis, Mohsen and Yuncken's \cite{OmarMohsen}-\cite{AMY} gave a second elaborate answer involving representations, more subtle than the  one considered below. In the process, they showed that the Nash blowup invented by Omar Mohsen (see Section \ref{sec:Nash}) provides an answer to Question \ref{thm:symbol} sufficient for several purposes, including dealing with pseudo-differential operator calculus  developed by Androulidakis and Skandalis in \cite{zbMATH05843246}, using the holonomy Lie groupoid and bisubmersions.  We present our interpretation of this answer.

\vspace{0.3cm}

To start with, we will define two ``symbols'':
\begin{enumerate}
\item the symbol of an element in the universal enveloping algebra $ \mathcal U(\mathcal F)$ of $ \mathcal F$, 
\item the symbol  of a longitudinal differential operator, i.e., an element in ${\mathrm{Diff}}(\mathcal F) $. 
\end{enumerate}

\vspace{0.5cm}
Before defining these two symbols,
let us very briefly recall the definition of the symbol in the context of Lie algebroids as in \cite{zbMATH07453484}. Let $B \to M $ be a Lie algebroid. We denote by $\mathcal U(B) $ and call the universal algebra of the Lie algebroid $B$ the quotient of the universal algebra  $U(\Gamma(B)) $ of the Lie algebra of sections of $B$ by the ideal generated by 
  $$   X \cdot (f Y) - (fX) \cdot Y - (\rho(X)[f] \, Y). $$
  with $X,Y \in \Gamma(B) $ and $f \in \mathcal C^\infty(M)$.
  It is again a filtered algebra, with $\mathcal U^{\leq k}(B) $ being the subspace generated by monomials of degree $ \leq k$.
  Now, recall that
\begin{enumerate}
\item upon choosing a Lie algebroid connection, there is a grading preserving coalgebra isomorphism $\mathcal U(B) \simeq \Gamma(S(B)) $ from the universal enveloping algebra of $B$ to the symmetric algebra of $B$. See, e.g., \cite{zbMATH07369649}.
\item since this isomorphism preserves the grading, it induces a $\mathcal C^\infty(M) $-linear map from $\mathcal U^{\leq k}(B) $ to\footnote{We insist that it is $ S^k(B)$, not $ S^{\leq k}(B)$: we project on the top component, i.e., the space generated by monomials of degree $k $.} $ S^k (B)$. 
\item This map does not depend on the choice of a Lie algebroid connection\footnote{Only the top component does not depend on the connection, the low components do. This independence means that we could work with local connections, so that what we say here extends to the real analytic or holomorphic settings.}.
\end{enumerate}
Now, a section of $S^k (B)$ may be seen as a function on $B^*$, which is fiberwise polynomial and homogeneous of degree $ k$. 
For every $ P \in \mathcal U^{\leq k}(B)$, we denote by $\sigma_P $ this function and call it the \emph{symbol} of $P$. 
\vspace{0.5cm}

Let us now define the symbol as an element in $ \mathcal U (\mathcal F)$. For every leaf $ L$ of a singular foliation $\mathcal F $ on a manifold $M$, there exists a natural restriction map
 $$ \mathcal F \longrightarrow \Gamma_L (A_L)  $$
where $A_L $ is the holonomy Lie algebroid of the leaf $L$ (see section \ref{sec:isotropy}). This map is $\mathcal C^\infty(M) $-linear and is a Lie algebra morphism, so that it induces an algebra morphism 
 $$ \begin{array}{ccc} U (\mathcal F) & \to &  U(\Gamma(A_L)) \\ P& \mapsto &P(L) \end{array} $$ from the algebra of longitudinal operators $U (\mathcal F)  $ to the universal enveloping algebra $U(\Gamma(A_L)) $.
 It goes to the quotient to induce an algebra morphism 
 \begin{equation} \label{eq:PofL}\begin{array}{ccc} \mathcal U (\mathcal F) & \to &  \mathcal U(A_L) \\ P& \mapsto &P(L) \end{array} \end{equation}  from the universal enveloping algebra $\mathcal U (\mathcal F)  $ of the singular foliation $ \mathcal F$ to the universal enveloping algebra $\mathcal U(A_L) $ of the Lie algebroid $ A_L$.

 For a given $P \in \mathcal U^{\leq k}(\mathcal F)$ of degree  $\leq k $, consider, for every leaf $L$, the symbol $\sigma_{P(L)} $. It is a fiberwise homogeneous of degree $k$ smooth function on $ A_L^*$.

 \begin{remark}
     For any leaf $L$, the realization of $P$, restricted to $L$, is a differential operator in the usual sense. For a regular leaf, $ A_L=TL$ and $ \sigma_L(P)$ is simply the usual symbol \cite{NistorPingAlan}  of this differential operator.
 \end{remark}
 
 We call \emph{symbol of $P$} the collection $(L, \sigma_{P(L)}) $ indexed by leaves of $\mathcal F $ of the functions $ \sigma_L(P)$. We denote it by $\sigma_P $.

\vspace{0.5cm}
From now on, we will assume that the regular leaves of $(M,\mathcal F) $ are all the same dimension, so that the Nash blowup
 $$ \left( {\mathrm{Bl}}(M,\mathcal F) , \pi^! \mathcal F \right)   $$
makes sense, see Section \ref{sec:Nash}. Recall from Theorem \ref{thm:NashIsProjective} that the Nash blowup is a Debord singular foliation, with associated Lie algebroid (the Nash blowup Lie algebroid of Definition \ref{def:Nashblowup}) the canonical quotient bundle $A_{\mathrm{Bl}} $ of the Grassmann bundle, restricted to $ {\mathrm{Bl}}(M,\mathcal F)$. From now on, we denote by $ ({\mathrm{A_{Bl}}} , \rho_{\mathrm{Bl}}) $ this Lie algebroid.
Consider the lift
 $$ \mathcal F \to \Gamma_{{\mathrm{Bl}}(M,\mathcal F)} ({\mathrm{A_{Bl}}})  $$
 mapping $X \in \mathcal F$ to the unique section of $\Gamma_{{\mathrm{Bl}}(M,\mathcal F)} ({\mathrm{A_{Bl}}}) $ corresponding to $\pi^! X $.
 This lift is 1) a $\mathcal C^\infty(M) $-module morphism and 2) a Lie algebra morphism.
 This implies that it lifts to an algebra morphism
   $$  \mathrm{Diff} (\mathcal F) \longrightarrow \mathcal U\left({\mathrm{A_{Bl}}}\right) , $$
   where $\mathcal U\left({\mathrm{A_{Bl}}}\right)$ stands for the universal algebra 
 of the Nash-blowup Lie algebroid $ {\mathrm{A_{Bl}}}$.
 We denote by $ D \longrightarrow \pi^!{D}$ this algebra morphism.  Since $\pi^!{D}$ is an element of degree $ \leq k$ of the universal Lie algebroid $ {\mathrm{A_{Bl}}}$ , its symbol is a (perfectly well-defined) element of $\Gamma(S^k  {\mathrm{A_{Bl}}}) $, or, equivalently, a function on $A_{\mathrm{Bl}}^*$, which is fiberwise polynomial and homogeneous of degree $k$.
 Since its restriction to the regular part is entirely determined by the image of $D$ through the realization map, and since this symbol depends continuously on the base point by construction, the next definition makes sense. 

 \begin{definition}
 Let $D$ be a longitudinal differential operator of degree $\leq k $.
     We call \emph{symbol} of $ D$ of degree $k$ the symbol of the element $ \pi^! {D}$ in the universal algebra $\mathcal U({\mathrm{A_{Bl}}}) $ of the Nash blowup algebroid ${\mathrm{A_{Bl}}}$.

     By construction, this symbol that we denote $\sigma_{D} $ is a fiberwise homogeneous of degree $k$ polynomial function on the dual ${\mathrm{A_{Bl}}}$ of the Nash blowup Lie algebroid of Definition \ref{def:Nashblowup}. 
 \end{definition}

We now have two symbols: one for the algebra of longitudinal differential operators $ {\mathrm{Diff}}(\mathcal F)$ and one for the universal enveloping algebra $\mathcal U(\mathcal F) $.

To relate these two symbols, we need to give some explanation about the Helffer-Nourrigat cone \cite{AMY} (see also Exercise \ref{exo:Helffer}).
Let us first explain a general phenomenon. 
Given a vector bundle $ E \to M$ and a map $ \Pi \colon N \to M$, consider the pull-back bundle $ \pi^! E \to N$. Now, let $T \subset  \Pi^* E $ be a subvector bundle over $N$ then there is a natural fiberwise injective vector morphism 
 $$ \xymatrix{ \left(\frac{\pi^! E}{T}\right)^*  \ar[d] \ar^{\mathfrak j}[r] & E^* \ar[d]\\ N\ar[r]&M }  .$$
It consists in identifying for every point $n \in N$ the dual of the quotient space $\left(\frac{\pi^! E_n}{T_n}\right)^*$ with the annihilator $ T^\perp_n \subset E^*_{\pi(n)}$ of $ T_n$, and to inject the later in $ E^*$.
Moreover, the image of $\mathfrak j $ can be described as follows. For every $m \in M$,  its intersection with  $ E_m$ is a union of vector spaces of dimension $ {\mathrm{rk}}(E)- {\mathrm{rk}}(T)$. We therefore call \emph{image cone} this image.

We apply this construction to
\begin{enumerate}
\item $E:=A $, with  $ (A,\rho)$ an anchored bundle over $ \mathcal F$,
\item $ N:= {\mathrm{Bl}}(M,\mathcal F)$ the Nash blowup\footnote{Of course, it may not be a manifold, but this defect  has no practical consequence here.} of Section \ref{sec:Nash}, computed with respect to $ (A,\rho)$, 
\item $T$ the canonical bundle $ \tau_A^{-r}$ on the Grassmannian  $ {\mathrm{Grass}}_{-r}(A)$, restricted to ${\mathrm{Bl}}(M,\mathcal F)$.
\end{enumerate}
In this case, Theorem \ref{thm:NashIsProjective} identifies the quotient $ \pi^! E / T $ with the Nash blowup Lie algebroid\footnote{that we recall to be the Lie algebroid of the Mohsen's groupoid, see Theorem \ref{thm:NashIsProjective}} $A_{Bl}$ of Definition \ref{def:Nashblowup} associated to the Nash blow up of $ \mathcal F$.
We therefore obtain a fiberwise injective\footnote{But of course not injective} vector bundle morphism:
 $$  \xymatrix{\mathfrak j\colon {\mathrm{A}_{\mathrm{Bl}}^*} \ar[d] \ar[r] & \ar[d] A^*\\
\mathrm{Bl}(M,\mathcal F) \ar[r]^{\pi}&M}.$$
We leave it to the reader to check\footnote{This can be proven as follows: Since $\pi $ is proper, both sets are closed. It therefore suffices that the image of $ \mathfrak j$ on $ \pi^{-1}(M_{\mathrm{reg}})$ (i.e., the regular part) coincides with the image of $ \rho^*$ (i.e., the annihilator of the kernel of $ \rho$).} that the image cone $\mathfrak j({\mathrm{A}_{\mathrm{Bl}}^*}) $  in this case coincides with the Helffer-Nourrigat cone $  {\mathrm{HN}}(\mathcal F)$ (computed with respect to $(A,\rho)$) defined in Exercise \ref{exo:Helffer}, i.e.,
$$  {\mathrm{HN}}_A(\mathcal F) := \overline{ \cup_{m \in M_{\mathrm{reg}}} (\rho^*(T^* _m M))} \subset A^*.$$
Moreover, for every leaf $L$ of $ \mathcal F$, the image of $ \mathfrak j$ takes values in the annihilator of the strong kernel of $A$,
which can be identified with $ A_L^*$ (the dual of the holonomy Lie algebroid).
The previous vector bundle morphism, therefore,  induces a vector bundle morphism\footnote{Here, we consider $\coprod_{L \in {\mathrm{Leaves}}} A_L^*$ as a “singular” vector bundle over $ M$: the ranks of the fibers vary.} (which is still fiberwise injective)
 $$  \xymatrix{q \colon {\mathrm{A}_{\mathrm{Bl}}^*} \ar[d] \ar[r] & \ar[d] \coprod_{L \in {\mathrm{Leaves}}} A_L^*\\
\mathrm{Bl}(M,\mathcal F) \ar[r]^{\pi}&M}$$
whose image, again,  the Helffer-Nourrigat c\^one  $ {\mathrm{HN}}(\mathcal F)$, described as in question 3 in Exercise \ref{exo:Helffer}.
The concepts of the following proposition are maybe tricky to describe, but its proof it simply a commutative diagram argument. 

\begin{proposition}
\label{prop:pullbackHN}
Let $D \in {\mathrm{Diff}}^{\leq k}(\mathcal F)$  be a longitudinal differential operator of degree $\leq k $ and let $P \in \mathcal U^{\leq k}(\mathcal F)$ be any element whose  realization is $D$. The  symbol $ \sigma_D$ is the pull-back through the  projection $q$ of the restriction of the  symbol $ \sigma_P$ of $P$ to the Helffer-Nourrigat cone bundle of $ \mathcal F$.  
In Equation:
 $$  \sigma_D = q^* \left. \sigma_P \right|_{{\mathrm{HN}}(\mathcal F)}  .$$
\end{proposition}
\begin{proof}
It holds true on the regular part. The result follows by density.
\end{proof}

A remarkable consequence of this statement is that the restriction to  ${\mathrm{HN}}(\mathcal F)$ of the symbol of an element in $P$ depend on its realization $ \underline{P}$ only.

Now, recall that a differential operator is said to be elliptic if its symbol vanishes only at the origin.  For a Lie algebroid \cite{NistorPingAlan}, an element in the universal Lie algebra is said to be elliptic if its symbol vanishes only along the zero section. 

\vspace{.5cm}

\begin{questions}{}{ques:symbol}
Let $\mathcal F $ be a singular foliation. 
What is a longitudinaly elliptic differential operator?
\end{questions}

\vspace{.5cm}

Here is, in our opinion, the correct answer \cite{OmarMohsen}.

\begin{definition}
A longitudinal differential operator $D$ of degree $k$ is said to be \emph{longitudinally elliptic} if its symbol\footnote{Recall that the latter is a function on the dual $ A^*_{\mathrm{{Bl}}} $ of the Nash blowup Lie algebroid.} $ \sigma_{D}$ is a strictly positive function, except on the zero section.
\end{definition}

In particular, this implies that $D$ is longitudinally elliptic if and only if  $\pi^! (D)$ is elliptic for the Nash blowup Lie algebroid $ {\mathrm{A_{Bl}}}$ in the sense of \cite{zbMATH06351314,NistorPingAlan}, which seems to us to be a convincing justification of the notion.
Now, chose some $ P \in \mathcal U(\mathcal F)$ whose realization is $D$.
The symbol $ \sigma_P$ does \emph{not} need to be strictly positive outside the Helffer-Nourrigat cone by Proposition \ref{prop:pullbackHN},  so that our definition does \emph{not} imply that $ P(L) \in \mathcal U(A_L)$ (see Equation \eqref{eq:PofL}) is elliptic for all leaves (although it certainly has to be elliptic on regular leaves).

\section{Cohomologies of a singular foliation}

We already saw that the derived cohomological spaces $\mathrm{Tor}_{\mathcal C^\infty(M)}(\mathcal F,\mathbb K) $
come equipped with a Lie $\infty $-algebra structure, whose cohomology permits to solve some elementary problems. But these are cohomologies associated to points or to leaves.
Our next question is rather vague:

\vspace{0.5cm}

\begin{questions}{Relevant cohomologies?}{}
What are the interesting global cohomology theories for singular foliations?
\end{questions}

Here are several candidates\footnote{We describe them in the smooth context: for the real-analytic or holomorphic settings, one has to add a \v{C}ech-type differential for a good covering – as always in sheaf theory}.
Also, for any $ \mathcal C^\infty(M)$-module $ \mathcal E$, the notation $\mathcal E \wedge_{\mathcal O} \mathcal E$ below stands for the wedge product over $ \mathcal O$, i.e., we allow 
 $$ X \wedge FY = FX \wedge Y \hbox{ for all $X,Y \in \mathcal E $, $F \in \mathcal O$} $$

\begin{enumerate}
    \item \emph{Longitudinal cohomology of a singular foliation} , see
    \cite{LLS}, Section 4.1. Let us describe it. To make the notation easier, we write $ \mathcal O$ instead of $ \mathcal C^\infty(M)$:
\vspace{0.5cm}
    
   \scalebox{0.9}{ \hbox{$ \begin{tabular}{|l|l|} \hline
        Chains in degree $k$ & Differential on chains of degree $k$ \\[3pt] \hline  Skew symmetric and $\mathcal O$-multilinear maps  & \\     $ \underbrace{\mathcal F \wedge_{\mathcal O} \cdots \wedge_{\mathcal O}  \mathcal F}_{k-\hbox{times}} \longrightarrow \mathcal O $  &  $\forall \omega \in {\mathrm{Hom}}_{\mathcal F}^k(\mathcal F, \mathcal O)$\\  & and all $ X_0, \dots, X_k \in \mathcal F$ \\
        i.e., ${\mathrm{Hom}}_{\mathcal F}^\bullet
        (\mathcal F, \mathcal O):={\mathrm{Hom}}_{\mathcal O}(\wedge_{\mathcal O}^\bullet \mathcal F, \mathcal O)$ & $\delta \omega \, (X_0, \dots, X_k) = $\\
        & $\sum_{i=0}^k (-1)^i Xi\left[\omega (X_0, \dots, \widehat{X_i}, \dots, X_k)\right] $ \\ &$ \, \, \, + \sum_{i < j} (-1)^{i+j+1} \omega([X_i,X_j], X_0, \dots, \widehat{X_i},\dots, \widehat{X_j},\dots, X_k)$ \\[53pt]
        \hline  
        Chains in degree $0$ & Differential on chains of degree $0$ \\[3pt] 
        \hline  & \\
        In degree $0$, chains are simply elements of $\mathcal O $  & $\forall F \in \mathcal O$,  \\
        & \\ 
        &  $\begin{array}{llll}\delta (F) \colon  & \mathcal F& \to & \mathcal O\\ & X& \mapsto &X[F]  \end{array} $ \\[3pt] \hline 
    \end{tabular}
    $}}

\vspace{0.5cm}
    
    For a regular foliation, this cohomology is simply the De Rham cohomology along the leaves, \emph{i.e.,} it is the complex $\left(\Gamma(\wedge^\bullet T^* \mathcal F) , d^{dR}_\mathcal F\right)$ with $d^{dR}_{\mathcal F}$ being the De Rham differential restricted to $T\mathcal{F}$, but computed leaf by leaf. 
    \item The \emph{basic cohomology} is the sub-complex of $(\Omega(M) , d^{dR})$ made, in degree $k$, of all $k$-forms that vanish when $k$ vector fields in $\mathcal F $ are applied. Equivalently, these are $k$-forms $\omega $ whose pull-back to any leaf $L$ is zero.
    
    It has been studied by David Miyamoto \cite{Miyamoto}.

    \item The \emph{universal cohomology of $\mathcal F$} is the cohomology of the commutative differential graded algebra of functions\footnote{One can also choose compactly supported functions.} on any universal $Q$-manifold\footnote{See Section \ref{sec:Universal}.} of $\mathcal F $. This is more precisely defined as the cohomology of $\left(\Gamma(S(\oplus_{i \geq 1} E_i^* )), Q \right) $. 
    The definition makes sense: it can be proven that since any two  universal Lie $\infty$-algebroid of $\mathcal F $, say $ (E,Q)$ and $ (E',Q')$ are homotopy equivalent, the differential graded commutative algebras $\left(\Gamma(S(\oplus_{i \geq 1} E_i^* )), Q \right) $ and $\left(\Gamma(S(\oplus_{i \geq 1} (E_i')^* )), Q' \right) $ are homotopy equivalent in a unique up to homotopy manner, see \cite{LLS} Section 3.4.4. In particular, their cohomologies are canonically isomorphic\footnote{It is tempting to believe that because the geometric resolution has no cohomology, the universal cohomology has to be zero as well. This is not at all true.}. 
    
    Universal cohomology is linked to longitudinal cohomology, since there is a map of differential graded commutative algebras:
     $$ \hbox{ Longitudinal cohomology of $ \mathcal F$}  \longrightarrow \hbox{ Universal cohomology of  $\mathcal F$.} $$
    See the discussion on universal and longitudinal cohomologies in \cite{LLS}, Section 4.1.
    \item
     It is also interesting to consider the Chevalley-Eilenberg cohomology for the adjoint representation \cite{zbMATH07568450,Jotz} of any universal Lie $ \infty$-algebroid of $\mathcal F $.
 This coincides with   the cohomology of vector fields on the universal $Q$-manifolds, equipped with $ {\mathrm{ad}}_Q$, and plays a role in deformation theory. 
     
\end{enumerate}
The list will certainly be continued.

\clearpage
\addcontentsline{toc}{chapter}{Bibliography}
\bibliographystyle{alpha}
\bibliography{biblio}

\end{document}